\documentclass[reqno]{memo-l}
\usepackage{hyperref, amssymb, mathtools, url, enumerate}
\usepackage{tikz}
\usepackage{tikz-cd}
\usepackage[scr]{rsfso}
\usepackage{bm}
\usepackage{amssymb}
\usepackage{dsfont}
\usepackage{lettrine}
\usetikzlibrary{matrix,arrows,decorations.pathmorphing, decorations.markings, cd}
\usepackage[alphabetic]{amsrefs}
\usepackage{microtype}
\usepackage[marginal]{footmisc}

\def\twocell[#1]{\arrow[#1, dash, phantom, "\Rightarrow"{scale=1.125, yshift=-.4pt, description, allow upside down, sloped, inner sep=0pt}]}

\newtheoremstyle{mythm}%
{5pt plus 1pt minus 1pt}
{5pt plus 1pt minus 1pt}
{\slshape}
{\parindent}
{\scshape}
{.}
{ }
{}
\theoremstyle{mythm}
\newtheorem{thm}{Theorem}[section]
\newtheorem{lemma}[thm]{Lemma}
\newtheorem{cor}[thm]{Corollary}
\newtheorem{prop}[thm]{Proposition}
\newtheorem*{thm*}{Theorem}
\newtheorem*{prop*}{Proposition}
\newtheorem{introthm}{Theorem}
\newtheorem{introcor}[introthm]{Corollary}

\makeatletter
  \@removefromreset{section}{chapter}
  \@removefromreset{footnote}{chapter}
\makeatother
\AtBeginDocument{\renewcommand{\thesection}{\arabic{section}}}

\theoremstyle{definition}
\newtheorem{rk}[thm]{Remark}

\newtheorem{constr}[thm]{Construction}
\newtheorem{warn}[thm]{Warning}
\newtheorem{ex}[thm]{Example}

\newtheorem{convention}[thm]{Convention}
\newtheorem{notation}[thm]{Notation}
\newtheorem{definition}[thm]{Definition}
\newtheorem{variant}[thm]{Variant}

\newtheorem*{claim*}{Claim}
\AtBeginEnvironment{claim*}{\let\oldqed=\qedsymbol%
\renewcommand\qedsymbol{$\triangle$}}
\AtEndEnvironment{claim*}{\let\qedsymbol=\oldqed}

\numberwithin{equation}{section}

\newcommand{\qednow}{\pushQED{\qed}\qedhere\popQED}

\newcommand{\notehelper}[3]{\textcolor{#3}{$\blacksquare$}\marginpar{\ifodd\thepage\raggedright\else\raggedleft\fi\color{#3}\tiny \textbf{#2:} #1}}

\newcommand{\Aa}{{\mathcal{A}}}
\newcommand{\Bb}{{\mathcal{B}}}
\newcommand{\Cc}{{\mathcal{C}}}
\newcommand{\Dd}{{\mathcal{D}}}
\newcommand{\Ee}{{\mathcal{E}}}
\newcommand{\Ff}{{\mathcal{F}}}
\newcommand{\Mm}{{\mathcal{M}}}

\newcommand{\Ii}{{\mathcal{I}}}
\newcommand{\Jj}{{\mathcal J}}
\newcommand{\Kk}{{\mathcal{K}}}
\newcommand{\Ll}{{\mathcal L}}

\newcommand{\Rr}{{\mathcal R}}

\newcommand{\Uu}{{\mathcal{U}}}
\newcommand{\Vv}{{\mathcal{V}}}
\newcommand{\Ww}{{\mathcal{W}}}
\newcommand{\Xx}{{\mathcal X}}
\newcommand{\Yy}{{\mathcal Y}}

\newcommand{\R} {{\mathbb{R}}}
\newcommand{\C} {{\mathbb{C}}}
\newcommand{\Z} {{\mathbb{Z}}}
\newcommand{\Lhat}{{\widehat{{\smash{\cat{L}}\vphantom{t}}\kern-1.75pt}\kern1.75pt}}

\newcommand{\bbone}{{\mathds{1}}}

\newcommand{\gl}{\textup{gl}}

\DeclareMathOperator{\ev}{ev}

\newcommand{\Sqcof}{{\textup{Sq}^\textup{cof}}}
\newcommand{\exttimes}{\mathbin{\widehat{\smash{\times}}}}

\renewcommand{\phi}{\varphi}
\renewcommand{\epsilon}{\varepsilon}

\newcommand{\sk}{{\textup{sk}}}
\DeclareMathOperator{\Sp}{Sp}
\DeclareMathOperator{\Spc}{Spc}
\DeclareMathOperator{\SPC}{SPC}
\DeclareMathOperator{\Cat}{Cat}
\DeclareMathOperator{\CAT}{CAT}

\DeclareMathOperator{\Tw}{Tw}

\newcommand{\PrL}{\textup{Pr}^{\textup{L}}}
\newcommand{\PrR}{\textup{Pr}^{\textup{R}}}
\newcommand{\PrLT}{\textup{Pr}^{\textup{L}}_{T}}
\newcommand{\PrRT}{\textup{Pr}^{\textup{R}}_{T}}

\renewcommand{\top}{{\textup{top}}}

\newcommand{\ppo}{\mathbin{\raise.25pt\hbox{\scalebox{.67}{$\square$}}}}
\newcommand{\boxsmash}{\mathbin{\raise.5pt\hbox{\rlap{\lower.5pt\hbox{$\square$}}\kern.55pt$\smashp$\kern.55pt}}}
\newcommand{\smallboxsmash}{\mathbin{\raise.35pt\hbox{\rlap{\lower.35pt\hbox{$\scriptstyle\square$}}\kern.375pt$\scriptstyle\smashp$\kern.375pt}}}

\renewcommand{\hom}{\textup{hom}}

\newcommand{\Fun}{{\textup{Fun}}}
\DeclareMathOperator{\PSh}{PSh}

\DeclareMathOperator{\Ar}{Ar}

\DeclareMathOperator{\CMon}{CMon}

\newcommand{\Ntop}{{N_\textup{top}}}

\newcommand{\cat}[1]{\textbf{\textup{#1}}}

\newcommand{\FUN}{\cat{Fun}}
\newcommand{\maps}{\cat{maps}}

\newcommand{\op}{{\textup{op}}}
\newcommand{\triv}{{\textup{triv}}}

\newcommand{\GL}{{\textup{GL}}}
\renewcommand{\O}{{\textup{O}}}
\newcommand{\SO}{{\textup{SO}}}
\newcommand{\Sig}{\smash{\widetilde{\vphantom{t}\smash{\Sigma}}}}
\renewcommand{\hom}{{\textup{Hom}}}
\newcommand{\End}{{\textup{End}}}

\DeclareMathOperator{\core}{\iota}
\DeclareMathOperator{\CAlg}{\textup{CAlg}}

\DeclareMathOperator{\Alg}{\textup{Alg}}

\newcommand{\colim}{\mathop{\textup{colim}}}
\newcommand{\const}{\textup{const}}
\newcommand{\consto}{\ensuremath{c}}
\newcommand{\consts}{\textup{const}}
\newcommand{\id}{\textup{id}}
\newcommand{\pr}{\textup{pr}}
\renewcommand{\Pr}{\textup{Pr}}
\newcommand{\fgt}{\textup{fgt}}
\newcommand{\diag}{\textup{diag}}

\newcommand{\ct}{\textup{ct}}
\newcommand{\BC}{\textup{BC}}
\newcommand{\Un}{\textup{Un}}
\newcommand{\sh}{{\textup{sh}}}
\newcommand{\Fin}{\hbox{$\mathcal F\kern-1.7pt\textit{in}$}}
\newcommand{\bigAll}{\hbox{\hfuzz=20pt\hbox{$\mathcal A\kern-.3pt$}\setbox0=\hbox{$\ell\kern-.75pt\ell$}\hbox to 0pt{\hfuzz=8pt\kern-.075pt\unhcopy0}\hbox to 0pt{\hfuzz=8pt\kern.075pt\unhcopy0}\box0\kern.05pt}}
\newcommand{\smallAll}{\Aa\kern-.3pt\ell\kern-.6pt\ell}
\newcommand{\All}{\mathchoice{\bigAll}{\bigAll}{\smallAll}{\smallAll}}
\newcommand{\dual}{\flat}
\newcommand{\Orb}{{\vphantom{\textup{t}}\smash{\textup{Orb}}}}

\newcommand{\Glo}{\textup{Glo}}
\newcommand{\res}{\textup{res}}

\newcommand{\Ind}{\mathop{\textup{Ind}}\nolimits}

\newcommand{\Top}{\cat{Top}}
\newcommand{\GTop}[1]{\cat{$\bm{#1}$-Top}}

\newcommand{\topCat}{\cat{topCat}}

\newcommand{\BGcat}[1]{\mathbb B{#1}}

\newcommand{\Aut}{\textup{Aut}}

\newcommand{\ulhelper}[2]{\underline{\setbox0=\hbox{$#1#2$}\dp0=0.4pt \box0\relax}}
\newcommand{\ul}[1]{{\mathpalette\ulhelper{#1}}\hbox{\rule[-2pt]{0pt}{0pt}}}

\DeclareMathOperator{\im}{im}

\tikzset{mono/.style={>->},
         epic/.style={->>},
         norm/.style={->, shorten <=3pt,
		              postaction={decorate,
			                      decoration={markings,
                                              mark=at position 2pt with {\node[fill=white,inner sep=-1pt,circle] {$\scriptscriptstyle\otimes$};}}}}}

\newcommand{\oldrightarrownorm}{\mathrel{
\raise1.3pt\hbox{$\scriptscriptstyle\otimes$}\hskip-1.33pt{\to}}}

\let\smashp=\wedge

\newcommand{\incl}{{\mathop{\textup{incl}}}}

\newcommand{\iso}{\xrightarrow{\;\smash{\raisebox{-0.25ex}{\ensuremath{\scriptstyle\sim}}}\;}}

\tikzcdset{pullback/.style = {phantom, "\lrcorner"{very near start}}}

\newcommand\noloc{%
	\nobreak
	\mspace{6mu plus 1mu}
	{:}
	\nonscript\mkern-\thinmuskip
	\mathpunct{}
	\mspace{2mu}
}

\newcommand{\mySp}{{\mathfrak S\kern-.5pt\mathfrak p}}
\newcommand{\myLSp}{{\mathfrak d\kern-.25pt}\mySp}
\newcommand{\mySph}{{\mathfrak R\kern-.5pt\mathfrak e\kern-.5pt\mathfrak p\kern-.5pt\mathfrak S\kern-.5pt\mathfrak p}}

\newcommand{\ulmySplev}{\ul{\myLSp}_\textup{lvl}}
\newcommand{\myS}{{{\mathfrak S\kern-.4pt}}}

\newcommand{\osp}{{\textbf{\textup{OrthSp}}}}
\newcommand{\Gosp}[1]{\cat{$\bm{#1}$-OrthSp}}
\newcommand{\Gospc}[1]{\cat{$\bm{#1}$-$\cat{L}$-Top}}

\newcommand{\Th}{{\textup{Th}}}
\renewcommand{\th}{{\textup{th}}}

\newcommand{\myTh}{{\mathfrak T\kern-1pt\mathfrak h}}
\newcommand{\myth}{{\mathfrak t\kern-.5pt\mathfrak h}}

\newcommand{\pic}{\textup{pic}}
\newcommand{\Rep}{\textup{Rep}}
\newcommand{\VRep}{\textup{VRep}}
\newcommand{\Vect}{\textup{Rep}}

\newcommand{\BOP}{\cat{BOP}}
\newcommand{\BO}{B\textup{O}}
\newcommand{\Gr}{\cat{Gr}}
\newcommand{\eGr}{\cat{EGr}}

\newcommand{\Borcat}[1]{N_{\textup{top}}({#1}^{\dual})}

\newcommand{\grp}{\textup{grp}}

\newcommand{\cont}{\textup{cont}}
\newcommand{\cc}{\text{cc}}

\renewcommand{\thesection}{\thechapter.
\arabic{section}}

\begin{document}
\title{Global $\bm\infty$-categories and global Thom spectra}
\author{\textsc{Emma Brink}}
\author{\textsc{T\kern-.5ptobias Lenz}}
\address{Mathematisches Institut, Rheinische Friedrich-Wilhelms-Universität Bonn, Endenicher Allee 60, 53115 Bonn, Germany}

\begin{abstract}
    We introduce a framework of (Lie-)\emph{global $\infty$-categories}, which formalizes various families of $\infty$-categories indexed by compact Lie groups and equipped with suitable restriction functors along continuous group homomorphisms that occur naturally in equivariant homotopy theory and representation theory. 

    As our main results, we show that in this framework unstable and stable equivariant and global homotopy theory admit universal properties, refining and generalizing the results for finite groups from \cites{CLL_Global,CLL_Clefts}. 
    In particular, we characterize the passage from unstable to stable equivariant and global homotopy theory at the level of global $\infty$-categories as universally inverting the action of representation spheres in an appropriate sense.
    Building on this, we define parametrized equivariant and global Thom spectrum functors and show that they recover classical Thom spectrum constructions defined in terms of pointset models.
\end{abstract}

\thanks{The authors are associate members of the Hausdorff Center for Mathematics at the University of Bonn (DFG GZ 2047/1, project ID 390685813). Parts of this article were written while E.B. was a guest at the Centre for Geometry and Topology (DNRF151) at the University of Copenhagen, she would like to thank the Centre for its hospitality.}

\frontmatter
\maketitle
\setcounter{tocdepth}{1}
\begingroup
\microtypesetup{protrusion=false}
\tableofcontents
\endgroup

\mainmatter
\chapter*{Introduction}
\lettrine[findent=3pt,nindent=.25pt]{E}{quivariant homotopy theory} studies spaces with an action of a suitable group---typically a finite group or more generally a compact Lie group---and their invariants. Throughout the history of the subject, equivariant methods have found spectacular applications also in non-equivariant settings, for instance in Carlsson's proof of the Segal conjecture \cite{carlsson1984SegalConjecture}, Manolescu's disproof of the triangulation conjecture \cite{manolescu}, or the resolution of the Kervaire invariant one problem by Hill, Hopkins, and Ravenel \cite{HHR2016Kervaire}.

While the foundations of equivariant homotopy theory were originally laid in the 1980's using pointset models \cite{LMS}, in recent years there has been a flurry of activity to introduce $\infty$-categorical techniques into the subject. This has led to new descriptions of the fundamental ($\infty$-)categories studied in equivariant homotopy theory \cites{barwick2017spectral,cmnn,CLL_Spans,CHLL_NRings,LLP} as well as to universal properties for them \cites{nardin2016exposeIV,shah2021parametrized, CLL_Global,CLL_Clefts}, formulated using the language of \emph{parametrized higher category theory} \cite{exposeI}. The parametrized perspective
has since been successfully applied to a wide array of equivariant questions, ranging from the construction of multiplicative refinements of equivariant algebraic $K$-theory \cites{CHLL_NRings,equiv-motives} or global versions of TMF and tempered cohomology \cites{GLP24,global-ambidextrous} to more geometric topics like the study of group actions on manifolds \cites{hilman2024equivariantpoincaredualitycyclic,hilman2024parametrisedpoincaredualityequivariant,kirstein2025semifreeisovariantpoincarespaces}; moreover, it has once more also led to new results outside equivariant homotopy theory, for example in chromatic homotopy theory \cite{ben-moshe} or the theory of 6-functor formalisms \cite{CLL_Span2}.  

However, so far many of these $\infty$-categorical results have been confined to the setting of $G$-equivariant homotopy theory for \emph{finite} groups $G$. The purpose of this monograph is to push the theory beyond this case and to in particular establish parametrized universal properties for equivariant homotopy theory in the generality of all compact Lie groups. Before we explain our new results, however, let us spend some time to 
review the situation for finite groups.

\subsection*{Parametrized higher categories} 
A common theme throughout equivariant homotopy theory is that it is often helpful to not restrict one's attention to a fixed group, but to work with a whole family of groups, like all finite groups, all compact Lie groups, or at least all subgroups of a given group. In particular, the $\infty$-categories of $G$-spectra for varying $G$ come with various `change of group' functors (restriction functors, their adjoints, the Hill--Hopkins--Ravenel norms, geometric fixed points,\,\dots), and these are important tools to construct $G$-spectra from simpler data or to conversely reduce questions to smaller groups in an inductive fashion.

It should therefore not be too surprising that in order to obtain a universal property, one should not consider $G$-equivariant homotopy theory for a fixed group $G$, but instead for all suitable groups $G$ at the same time. This kind of structure is actually not unique to equivariant homotopy theory, and instead pops up in many areas of pure mathematics, with other examples being given by categories of representations or (in the finite group case) derived categories of group algebras or the equivariant Kasparov $\infty$-categories of \cite{BEL2023Kasparov}.
If one is willing to restrict to finite groups, these structures were formalized under the name \emph{$\Fin$-global $\infty$-categories} in \cite{CLL_Global}.\footnote{The cited reference uses the name `global ($\infty$-)categories' instead; for us it will however be more convenient to reserve this name for the case where one considers all compact Lie groups.} Such a $\Fin$-global $\infty$-category $\Cc$ consists of an $\infty$-category $\Cc(G)$ for every finite group $G$, together with \emph{restriction} functors $f^*\colon\Cc(G')\to\Cc(G)$ for all homomorphisms $f\colon G\to G'$ as well as higher structure witnessing that restriction along an inner automorphism is the identity; more formally, a $\Fin$-global $\infty$-category is a contravariant functor out of the $(2,1)$-category $\Glo_{\Ff\kern-.4pt\textit{in}}$ of finite connected 1-groupoids,\footnote{There is also a refinement of this notion that additionally encodes structure like the Hill--Hopkins--Ravenel norms and geometric fixed points for equivariant spectra, see e.g.\ \cite{LLP}; however, for the discussion of the universal properties only the restriction functoriality will be relevant.} and a \emph{$\Fin$-global functor} is just a natural transformation. What makes this very basic definition actually useful is that many foundational concepts of higher category theory like adjunctions, limits and colimits, or presentability admit `genuine' versions for $\Fin$-global $\infty$-categories via the general theory of \emph{parametrized higher categories} \cites{exposeI,martiniwolf2021limits}. As one example, \cite{CLL_Clefts} introduces a notion of \emph{equivariant presentability}; while this has a natural formulation internal to parametrized higher category theory, it can also be described in more down-to-earth terms as follows: a $\Fin$-global $\infty$-category $\Cc\colon\Glo^\op\to\CAT_\infty$ is equivariantly presentable if and only it factors through $\PrL$ and restriction along any \emph{injective} homomorphism admits a left adjoint satisfying a basechange condition analogous to the Mackey double coset formula. 

The framework of $\Fin$-global $\infty$-categories is then a natural home for equivariant homotopy theory in the same way ordinary $\infty$-categories provide a natural home for non-equivariant homotopy theory: \cite{CLL_Clefts} constructs $\Fin$-global $\infty$-categories of equivariant spaces and equivariant spectra by localizing classical pointset models and shows that these enjoy universal properties analogous to their non-equivariant counterparts---the former is the free equivariantly presentable $\Fin$-global $\infty$-category, whereas the latter is the free equivariantly presentable $\Fin$-global $\infty$-category that moreover satisfies a genuine version of stability.

Only demanding the existence of left adjoints for \emph{injective} group homomorphisms is of course somewhat arbitrary. If one instead demands the existence of well-behaved left adjoints for \emph{all} group homomorphisms, this naturally leads to \emph{global homotopy theory}, which we will recall next.

\subsection*{Global homotopy theory}
Equivariant enhancements of important cohomology theories like singular cohomology, topological $K$-theory, or various versions of bordism have a tendency to not crop up individually and isolated from each other, but to come to us as a `compatible family' of $G$-spectra as the (compact Lie) group $G$ varies. Just like it can be convenient to work with the $\infty$-categories of $G$-spectra for various $G$ simultaneously, \emph{global homotopy theory} takes the point of view that in the above situation one should consider the whole compatible family as the fundamental object. In Schwede's approach to the subject \cite{schwede2018global}, this idea is implemented using the notion of a \emph{global spectrum}, again defined via a model; the relation to the heuristic notion of a `compatible family of equivariant spectra' was then later made precise in \cite{LNP}.
Let us highlight three examples where the global perspective is useful even if one is only interested in results for a fixed group:
\begin{enumerate}
    \item\label{item:global-group-law} In \cite{hausmann2022global}, Hausmann proves that the equivariant homotopy ring of tom Dieck's $A$-equivariant \emph{homotopical complex bordism spectrum} $\cat{MU}_A$ carries the universal $A$-equivariant formal group law in the sense of Cole--Greenlees--Kriz for any compact abelian Lie group $A$, generalizing a famous non-equivariant result of Quillen's. While this result can be stated purely in terms of equivariant homotopy theory, Hausmann's argument starts from the observation that the $A$-spectra $\cat{MU}_A$ come from a global spectrum $\cat{MU}$, and one of the key steps in his proof is to show that the universal $A$-equivariant formal group laws for varying $A$ can similarly be organized into a global structure, which he calls a \emph{global group law}. The global structure then allows him to establish various algebraic properties of the universal equivariant formal group laws, from which he deduces his result. In the same paper and using a similar strategy, Hausmann also shows that the homotopy groups of the $A$-equivariant homotopical \emph{real} bordism spectrum $\cat{MO}_A$ carry the universal $A$-equivariant 2-torsion formal group law for any elementary abelian 2-group $A$. 
    \item In \cite{schwede-euler}, Schwede exploits the global structure on $\cat{MU}$ to prove regularity results for certain classes in the equivariant homotopy ring of $\cat{MU}_{\text{U}(n)}$. These results were used by La Vecchia \cite{lavecchia} to prove an Atiyah--Segal type completion theorem for $\cat{MU}_G$ for any compact Lie group $G$ (generalizing earlier work of Greenlees--May \cite{GM-MU}, who considered the case where $G$ has abelian identity component), and by Schwede himself \cite{schwede-chern} to compute the equivariant homotopy ring of $\cat{MU}_{\text{U}(n)}$ after a suitable completion.    
    \item\label{item:global-snaith} In the forthcoming \cite{schwede2025snaith}, Schwede proves equivariant generalizations of Snaith's Theorem \cite{snaith} expressing the complex topological $K$-theory spectrum as a localization of the suspension spectrum of $\C\text{P}^\infty$. Similarly to the previous examples, his proof proceeds by first exhibiting a global refinement $\cat{KU}$ of complex topological $K$-theory as a suitable localization, taking advantage of the global structure, and then deducing the equivariant results from this.
\end{enumerate}

\cite{CLL_Global} shows that if one imposes the stronger, `global' presentability condition (with well-behaved left adjoints for all homomorphisms), this leads to universal descriptions of unstable and stable global homotopy theory with respect to \emph{finite} groups: in particular, if we evaluate the free globally presentable and genuinely stable $\Fin$-global $\infty$-category at the trivial group, this yields the localization of the $\infty$-category of global spectra at the so-called `$\Fin$-global weak equivalences,' i.e.\ where we discard all equivariant information for non-discrete compact Lie groups. More generally, one can describe the universal example completely in terms of the pointset models of \emph{$G$-global spectra} from \cite{g-global}.

\subsection*{Beyond finite groups}
While the $\Fin$-global framework has been successfully applied for example to the study of equivariant algebraic $K$-theory in \cite{swan} and \cites{CHLL_NRings,LLP}, it is not yet the correct home for global and equivariant homotopy theory with respect to general compact Lie groups, as considered for example in the applications $(\ref{item:global-group-law})$--$(\ref{item:global-snaith})$ above. In fact, even if one were ultimately only interested in the case of finite groups, the above global arguments require one to consider non-discrete compact Lie groups: the definition of a global group law in $(\ref{item:global-group-law})$ involves $\text{U}(1)$-equivariant data, and similarly the global localization result from \cite{schwede2025snaith} outlined in $(\ref{item:global-snaith})$ is phrased in terms of $\text{U}(n)$-equivariant homotopy classes.\footnote{Both of these are instances of a common pattern: if we have a global spectrum together with a $G$-equivariant homotopy class for any group $G$ and any (unitary) representation $V$ such that these are suitably compatible with restrictions, then the \emph{universal} such classes, from which all other can be obtained by restriction, come from the tautological $\text{U}(n)$-representations. In settings where one only cares about \emph{abelian} compact Lie groups, one can instead restrict to tori or even just to $\text{U}(1)$ as any representation of an abelian compact Lie group splits into 1-dimensional summands.}

It is therefore natural to look for a framework of \emph{global $\infty$-category theory} refining $\Fin$-global $\infty$-category theory and allowing one to formulate universal properties for equivariant and global homotopy theory with respect to all compact Lie groups. And indeed, the \emph{global orbit category} $\Glo$ of Gepner and Henriques \cite{gepnerhenriques2007orbispaces} (which we review in Definition~\ref{def:Glo}) provides a natural replacement of $\Glo_{\Ff\kern-.4ptin}$, which allows us to simply define the $\infty$-category of global $\infty$-categories as $\CAT_{\Glo}\coloneqq \Fun(\Glo^{\op},\CAT_{\infty})$.  Unfortunately however, generalizing the above $\Fin$-global results to global $\infty$-categories is not at all straightforward: in particular, neither the notion of parametrized stability used in \cite{CLL_Global} nor its generalization from \cite{CLL_Spans} are sufficient to formulate the universal properties of equivariant or global spectra for compact Lie groups. Moreover, while the framework of parametrized higher category theory already allows us to state the expected universal properties of global and equivariant \emph{spaces}, the proof given in the finite group case breaks down in the new setting: namely, it crucially used that $\Glo_{\Ff\kern-.4pt\textit{in}}$ is a $2$-category to reduce the verification of the universal property to a technical, but ultimately doable computation, whereas $\Glo$ is not an $n$-category for any $n<\infty$.

Despite these obstacles, we will demonstrate that there is still a good theory of global $\infty$-categories with respect to all compact Lie groups. In particular, we will prove universal properties of unstable and stable, equivariant and global homotopy theory in this framework, and we will use this to devise an $\infty$-categorical approach to equivariant and global Thom spectra like $\cat{MO}$ and $\cat{MU}$.

\subsection*{Unstable results} Let us now describe our results in more detail, beginning with the unstable situation. The general theory of parametrized presentability from \cites{martiniwolf2022presentable,CLL_Clefts} specializes to notions of \emph{equivariant} and \emph{global presentability} for global $\infty$-categories; as in the finite group case, these differ in that the former only demands the existence of well-behaved left adjoints for injective homomorphisms, while the latter requires them for all homomorphisms. 

Our first main results describe the free equivariantly and globally presentable $\infty$-categories in terms of classical models of unstable equivariant and global homotopy theory, respectively. More precisely, we introduce a general \emph{continuous Borel construction}, which associates to a topologically enriched category $\cat{C}$ a certain global $\infty$-category sending a compact Lie group $G$ to the nerve of the topological category $\cat{$\bm G$-$\cat{C}$}$ of (continuous) $G$-objects in $\cat{C}$, with the evident restriction functoriality in group homomorphisms. Applying this to $\cat{C}=\cat{Top}$ and localizing in each degree at the \emph{$G$-equivariant weak equivalences} classically considered in equivariant homotopy theory yields a global $\infty$-category $\ul\myS$ sending $G$ to the $\infty$-category $\cat{$\bm G$-Top}[\text{$G$-weak equivalences}^{-1}]$ of \hbox{\emph{$G$-spaces}} and sending a homomorphism of compact Lie groups to the derived functor of the pointset level restriction. Similarly, starting with the pointset model of global spaces from \cite{schwede2018global} yields a global $\infty$-category $\ul\myS_\gl$ sending the trivial group to the $\infty$-category of global spaces from \emph{op.\ cit.} and more generally sending a compact Lie group $G$ to the $\infty$-category of $G$-global spaces from \cite{barrero2021}, with restrictions again defined as derived functors of pointset level functors. We then prove:

\begin{introthm}[See Theorem~\ref{thm:unstable-main}]\label{introthm:unstable-main}
    The global $\infty$-category $\ul\myS_\gl$ is the free globally presentable global $\infty$-category. More precisely, for any other globally presentable $\Cc$, evaluation at the terminal object defines an equivalence
    \[
        \ul\Fun^\textup{L}(\ul\myS_\gl,\Cc)\iso\Cc,
    \]
    where the left-hand side denotes the global $\infty$-category of left adjoint global functors (see Remark~\ref{rk:FunL}). In particular, we obtain for any such $\Cc$ a natural equivalence between $\Cc(1)$ and the ordinary $\infty$-category of left adjoint global functors $\ul\myS_\gl\to\Cc$.
\end{introthm}

\begin{introthm}[See Corollary~\ref{cor:equiv-spaces-universal}]\label{introthm:unstable-equiv}
    The global $\infty$-category $\ul\myS$ is the free equivariantly presentable global $\infty$-category.
\end{introthm}

The proofs of both results proceed by comparing the pointset models to certain global $\infty$-categories provided by the general theory of parametrized higher categories, which we already know to satisfy the desired universal properties. Most of the work then goes into proving the global statement, with the equivariant version being a rather easy consequence of it: namely, the general theory from \cite{CLL_Clefts} exhibits the free equivariantly presentable global $\infty$-category as a full subcategory of the free globally presentable one, reducing us to similarly exhibiting the pointset model of $G$-spaces as a full subcategory of $G$-global spaces and checking that these embeddings are compatible with our identification from Theorem~\ref{introthm:unstable-main}.

Let us point out that it would be considerably easier to prove that $\ul\myS_\gl$ and $\ul\myS$ agree \emph{after evaluation at any fixed $G$} with the free globally or equivariantly presentable global $\infty$-category: for example, \cite{CLL_Clefts} describes the free equivariantly presentable global $\infty$-category in terms of $\infty$-categories of presheaves, and the pointwise equivalence to $\ul\myS$ essentially amounts to Elmendorf's Theorem \cite{elmendorf}. In that sense, our new contribution is the description of the whole $\Glo$-functoriality in terms of the model, which requires completely different ideas. While identifying the functoriality including all higher coherences may seem like a mere technicality at first, it is actually crucial if one wants to interpret classical pointset constructions as functors out of or into $\ul\myS_\gl$ and $\ul\myS$ (which then often allows one to use the universal properties to say something about these constructions) or if one wants to relate $\ul\myS$ and $\ul\myS_\gl$ to other global $\infty$-categories defined in terms of models. In particular, both our proof of the universal property of global spectra as well as our results on global Thom spectra require the full strength of Theorem~\ref{introthm:unstable-main} even if one only cares about the case of a fixed group---just like the applications \cites{CHLL_NRings,LLP,swan} of the earlier $\Fin$-global results  mentioned above would not have been possible without a complete model-level description of the functoriality.

\subsection*{Symmetric monoidal structures} Recall that the cartesian symmetric monoidal structure on the $\infty$-category $\Spc$ upgrades it to the initial symmetric monoidal $\infty$-category that is \emph{presentably} symmetric monoidal, meaning that the underlying $\infty$-category is presentable and the tensor product preserves colimits in each variable. There is again a natural parametrized version of this, leading to the notions of \emph{equivariantly presentably symmetric monoidal} and \emph{globally presentably symmetric monoidal} global $\infty$-categories. Building on the results above, we then also prove universal properties for suitable symmetric monoidal enhancements of $\ul\myS$ and $\ul\myS_\gl$ (again defined in terms of the model). For simplicity, we only state the global case in this introduction and refer the reader to Corollary~\ref{cor:S-equiv-times-initial} for the corresponding statement in the equivariant case:

\begin{introcor}[See Corollary~\ref{cor:S-gl-times-initial}]
    The box product from \cites{schwede2018global,barrero2021} derives to the cartesian  symmetric monoidal structure on $\ul\myS_\gl$, and this yields the initial globally presentably symmetric monoidal global $\infty$-category.
\end{introcor}

\subsection*{Stable results} As mentioned above, the notion of genuine stability from \cite{CLL_Global} is not yet sufficient to describe the universal property of global and equivariant spectra in our setting, so we have to take a different approach: 

In Corollary~\ref{cor:equiv-pointed}, we prove a based version of Theorem~\ref{introthm:unstable-equiv}, identifying the free \emph{pointed} equivariantly presentable global $\infty$-category with the analogous localization $\ul\myS_*$ of the continuous Borel category associated to $\cat{Top}_*$. 
The category of equivariantly presentable global $\infty$-categories admits a parametrized version of the Lurie tensor product (Proposition~\ref{prop:lurie-tensor-exists}), and Corollary~\ref{cor:equiv-pointed} implies via the general theory of smashing localizations from \cite{HA} that every pointed equivariantly presentable global $\infty$-category admits a unique $\ul\myS_*$-module structure. 
Adapting an idea from \cite{Linskens2023globalization}, we then define a pointed equivariantly presentable global $\infty$-category $\Cc$ to be \emph{representation stable} if for every compact Lie group $G$ and every finite-dimensional $G$-representation $V$ the tensoring $S^V\otimes{-}\colon\Cc(G)\to\Cc(G)$ is an equivalence, where $S^V$ denotes the 1-point compactification of $V$ (whose underlying non-equivariant space is a $\dim(V)$-dimensional sphere). 

We show in Corollary~\ref{cor:stabilization-exists} that the categories of representation stable equivariantly or globally presentable global $\infty$-categories are smashing localizations of the categories of equivariantly or globally presentable global $\infty$-categories, respectively. This in particular implies the existence of \emph{free} equivariantly and globally presentable representation stable global $\infty$-categories, which we then again describe explicitly via classical models of equivariant and global stable homotopy theory. In more detail, we construct a global $\infty$-category $\ul\mySp$ by localizing the continuous Borel category associated to the topological category of \emph{orthogonal spectra} \cite{MMSS} at the equivariant weak equivalences from \cite{mandell-may}; this then sends a compact Lie group $G$ to one of the classical definitions of the $\infty$-category $\mySp_G$ of genuine $G$-spectra, with functoriality via the (left) derived functors of restrictions. If we instead localize with respect to a finer notion of weak equivalence, this yields a global $\infty$-category $\ul\mySp_\gl$ whose value at the trivial group is now literally the $\infty$-category of global spectra from \cite{schwede2018global} (as opposed to its localization at the $\Fin$-global weak equivalences), while its value at a general compact Lie group $G$ is the $\infty$-category $\mySp_\text{$G$-gl}$ of \emph{$G$-global spectra} in the sense of \cite{schwede-stiefel}. We then prove:

\begin{introthm}[See Theorem~\ref{thm:equiv-main} and Corollary~\ref{cor:equiv-spectra-smash}]\label{introthm:stable-equiv-main}
    The global $\infty$-category $\ul\mySp$ is the free representation stable equivariantly presentable global $\infty$-category. Moreover, the smash product of orthogonal spectra derives to a symmetric monoidal structure on $\ul\mySp$, making it into the initial representation stable equivariantly presentably symmetric monoidal global $\infty$-category.
\end{introthm}

\begin{introthm}[See Theorem~\ref{thm:stable-main} and Corollary~\ref{cor:global-spectra-smash}]\label{introthm:stable-main}
    The global $\infty$-category $\ul\mySp_\gl$ is the free representation stable globally presentable global $\infty$-category. Moreover, the smash product derives to a symmetric monoidal structure on $\ul\mySp_\gl$, making it into the initial representation stable globally presentably symmetric monoidal global $\infty$-category.
\end{introthm}

A similar characterization of the $\infty$-category of global spectra as a suitable stabilization internal to $\myS_\gl$-modules in $\PrL$ had been announced by Gepner and Nikolaus \cite{gepner-nikolaus}. As we learned from them, their argument contained a gap, and it is not clear whether their original claim is true. We comment on the relation to our result in more detail in Remark~\ref{rk:gepner-nikolaus}.

Let us point out that the proofs of Theorems~\ref{introthm:stable-equiv-main} and~\ref{introthm:stable-main} are quite different from each other: In the equivariantly presentable setting, we build on work of Cnossen \cite{twisted-ambidexterity} to show that representation stabilization is essentially a pointwise construction, allowing us to prove the symmetric monoidal part of Theorem~\ref{introthm:stable-equiv-main} by combining Theorem~\ref{introthm:unstable-equiv} with earlier work of Gepner--Meier \cite{gepnermeier2020equivTMF}; we then use the general theory of representation stabilization we develop here to deduce the non-monoidal version from this in a second step. 

In contrast to that, representation stabilization in the \emph{globally} presentable setting is a genuinely parametrized phenomenon; in particular, global spectra are \emph{not} just the na\"ive pointwise stabilization of global spaces (see Remark~\ref{rk:not-naive-stab}), meaning we have to work considerably harder for the proof of Theorem~\ref{introthm:stable-main}.
 
\subsection*{Global and equivariant Thom spectra}
With the above universal properties at hand, one can now try to recast classical pointset level constructions from equivariant and global homotopy theory in a model-independent fashion. We demonstrate the feasibility of this approach in the case of Thom spectra.

In §\ref{subsec:J}, we introduce a symmetric monoidal global $\infty$-groupoid $\ul\VRep^\oplus$ of \emph{virtual representations}, and we give two equivalent constructions of a global refinement 
\[
    \ul{\mathfrak J}_\gl\colon\ul\VRep^\oplus\to\ul\mySp_\gl^\otimes
\]
of the classical $J$-homomorphism $B\O\times\Z\to\pic(\Sp)$: one (Construction~\ref{constr:global-J}) using our continuous Borel construction, and one obtained from the non-equivariant construction by general parametrized techniques (see Theorem~\ref{thm:who-is-afraid-of-the-J-homomorphism} and Remark~\ref{rk:J-with-fewer-models}).
Using general results from parametrized higher category theory, we identify the free globally presentably symmetric monoidal global $\infty$-category under $\ul\VRep^\oplus$ with a certain symmetric monoidal global $\infty$-category $\smash{\ul\Spc_{\Glo/\ul\VRep^\oplus}^\otimes}\vphantom{g_g^\otimes}$ of global $\infty$-groupoids over $\ul\VRep^\oplus$. This allows us to define the \emph{global Thom spectrum functor} as the unique symmetric monoidal left adjoint global functor $\smash{\ul\Th_\gl\colon\ul\Spc_{\Glo/\ul\VRep^\oplus}^\otimes}\to\ul\mySp_\gl^\otimes$ extending the global $J$-homomorphism $\ul{\mathfrak J}_\gl$, in analogy with the construction of the non-equivariant Thom spectrum functor from \cite{thom-oo}, also cf.~\cite{horev-klang-zou}. Evaluating at any compact Lie group $G$, this in particular yields a symmetric monoidal left adjoint Thom spectrum functor 
\[     \Th_\text{$G$-gl}^\otimes\colon\Spc_{\text{$G$-gl}/\triv_G\ul\VRep^{\oplus}}^\otimes\coloneqq\PSh(\Glo_{/\BGcat{G}})_{/\triv_G\ul\VRep^{\oplus}}^\otimes\to\mySp_\text{$G$-gl}^\otimes,\] 
and the fact that these come from a global left adjoint encodes strong compatibilites with respect to restrictions and inductions, see Remark~\ref{rk:incoherent-description-Thom-gl} for a precise statement. 

As our last main result, we then show that the symmetric monoidal global functor $\ul\Th_\gl$ admits a pointset description as a global refinement, originally due to Schwede, of a construction of Sagave and Schlichtkrull \cite{sagave-schlichtkrull-thom}. While we have to refer the reader to §\ref{subsec:model-thom} for the details of this construction and to Theorem~\ref{thm:comparison-Thom} for the precise statement of our comparison, this has the following concrete consequence:

\begin{introthm}[See Theorem~\ref{thm:MO-as-colim}]\label{introthm:mo-as-colim}
    There exists an equivalence
    \[
        \cat{MO}\simeq\colim_{\ul\VRep_{[0]}\!}\ul{\mathfrak J}_\gl
    \]
    of global spectra between the global real homotopical bordism spectrum $\cat{MO}$ from \cite{schwede2018global} and the globally parametrized colimit of the restriction of the global $J$-homomorphism to the subgroupoid of zero-dimensional virtual representations.
\end{introthm}

In fact, we can also describe the universal cocone explicitly in terms of the so-called \emph{equivariant Thom classes} of $\cat{MO}$, see Remark~\ref{rk:univ-cocone-thom-classes}.

As the reader might expect at this point, we moreover discuss an equivariant version of the above results. In particular, we construct an \emph{equivariant Thom spectrum functor}  by similarly extending an equivariant version $\ul{\mathfrak J}\colon\ul\VRep^\oplus\to\ul\mySp^\otimes$ of the $J$-homomorphism. 
At every compact Lie group $G$, this evaluates to a symmetric monoidal left adjoint 
\[  \smash{\Th_{G}\colon\Spc_{G/\ul\VRep_G^\oplus}^\otimes}\to\mySp_G^\otimes,\] and these functors commute with restrictions along arbitrary group homomorphisms and inductions along \emph{injective} group homomorphisms, see Remark~\ref{rk:incoherent-description-Thom-equiv} for a precise statement. 
We will then explain how the equivariant Thom spectrum functor can be recovered from the global Thom spectrum functor in two different ways (see Proposition~\ref{prop:equiv-thom-from-global} and Theorem~\ref{thm:Thom-vs-fgt}), which in particular allows us to describe the $G$-equivariant homotopical real bordism spectrum $\cat{MO}_G$ of Bröcker and Hook (the underlying $G$-spectrum of $\cat{MO}$) as a `$G$-parametrized' colimit in analogy to Theorem~\ref{introthm:mo-as-colim}, see Corollary~\ref{cor:MOG-as-colim}. 

\subsection*{Outline} The first part of this monograph (Sections~\ref{global-cat}--\ref{sec:basepoints}) is devoted to the proof of the universal properties of unstable equivariant and global homotopy theory. We begin in Section~\ref{global-cat} by recalling the necessary background from parametrized higher category theory and by formally defining global $\infty$-categories as well as the continuous Borel construction. In Section~\ref{sec:G-global-spaces}, we recall the relevant pointset models of $G$-equivariant and $G$-global unstable homotopy theory, and we prove some new results on the change-of-group adjunctions for (the pointset models of) $G$-global spaces that are crucial for establishing the universal property. Section~\ref{sec:unstable-univ-prop} is then devoted to the proof of Theorems~\ref{introthm:unstable-main} and~\ref{introthm:unstable-equiv} as well as their symmetric monoidal versions. Finally, we discuss pointed versions of these results in Section~\ref{sec:basepoints}.

The second part (Sections~\ref{sec:rep-stable}--\ref{sec:global-spectra}) is concerned with the universal properties of stable equivariant and global homotopy theory. We begin by defining the notion of representation stability (as well as a generalization to arbitrary parametrized $\infty$-categories) and establishing its basic properties in Section~\ref{sec:rep-stable}. In Section~\ref{sec:equiv-spectra} we then prove Theorem~\ref{introthm:stable-equiv-main} on the universal property of equivariant spectra, while the proof of the global analogue (Theorem~\ref{introthm:stable-main}) is the content of Section~\ref{sec:global-spectra}.
 
In the third part, we study Thom spectra. We introduce the $\infty$-categorical versions of the global and equivariant Thom spectrum functors in Section~\ref{sec:parametrized-thom}. The corresponding pointset construction is then discussed in Section~\ref{sec:model-thom}, where we also show that the approaches agree, in particular yielding the proof of Theorem~\ref{introthm:mo-as-colim} as well as of its equivariant counterpart.

The paper has three appendices. Appendix~\ref{app:thom-Stefan-functor} by Stefan Schwede establishes some homotopical properties of the pointset global Thom spectrum functor that are needed to prove our comparison result. In Appendix~\ref{app:model-cat} we recall several results about model categories and their associated $\infty$-categories for easy reference, while Appendix~\ref{app:here-be-pointsets} similarly collects some facts about the pointset topology of actions of compact Lie groups. 

\subsection*{Notation and conventions} We will write $\cat{Top}$ for the category of compactly generated weak Hausdorff spaces, and we will refer to its objects simply as `(topological) spaces.'

Throughout, we denote generic $\infty$-categories in a calligraphic font ($\mathcal C$, $\mathcal D$,\,\dots) and named $\infty$-categories in an upright font ($\Spc$, $\Sp$,\,\dots) or, in the case of localizations of pointset models, in a blackletter font ($\mathfrak S$, $\mySp$,\,\dots). We denote (enriched) 1-categories in boldface ($\cat{C}$, $\cat{Top}$,\,\dots). 
For an enriched category $\cat{C}$, we denote by $\maps(-,-)$ or $\maps_{\cat{C}}(-,-)$ its enriched mapping objects, and for an $\infty$-category $\Cc$, we denote by $\hom(-,-)$ or $\hom_{\Cc}(-,-)$ its mapping $\infty$-groupoids. 

For a quasi-category $S$ (i.e.\ a simplicial set satsfying the inner horn filling condition), we denote by $\core S\subset S$ its maximal Kan subcomplex. We will also write $\core$ for the functor $\Cat_\infty\to\Spc$ of $\infty$-categories induced by this pointset level construction, i.e.\ the right adjoint to the inclusion. 

\subsection*{Acknowledgements} We would like to thank David Gepner, Markus Hausmann, Christian Kremer, Fabio Neugebauer, Thomas Nikolaus, and Stefan Schwede for helpful feedback on earlier versions of this paper. We are moreover grateful to Stefan Schwede for contributing Appendix~\ref{app:thom-Stefan-functor} and to Thomas Nikolaus for helpful discussions about the relation between our work and \cite{gepner-nikolaus}. Finally, we would like to thank David Gepner, Thomas Nikolaus, and Stefan Schwede for sharing Corollary~\ref{cor:Thom-space-po} with us.

\chapter[The unstable story]{\for{toc}{\phantom{II}}The unstable story}
In this first part we will introduce \emph{global $\infty$-categories} as a rigorous implementation of the idea of  a `compatible' family of $\infty$-categories indexed by all compact Lie groups. The key examples of interest to us here are a certain global $\infty$-category $\ul\myS$ built from the $\infty$-categories $\myS_G$ of $G$-equivariant spaces for all compact Lie groups as well as a global analogue $\ul\myS_\gl$ built from the $\infty$-categories of $G$-global spaces in the sense of Barrero \cite{barrero2021}. 

Using the framework of parametrized higher category theory, we will see that many foundational concepts of higher category theory like colimits or presentability have `genuine' analogues for global $\infty$-categories, and we will show as our main results that the global $\infty$-categories $\ul\myS$ and $\ul\myS_\gl$ admit universal properties analogous to the universal property enjoyed by the $\infty$-category $\Spc$ of spaces in classical higher category theory.

\section[Global $\infty$-categories]{Global \texorpdfstring{\for{toc}{$\infty$}\except{toc}{$\bm\infty$}}{∞}-categories}\label{global-cat}
\subsection{The global orbit category} We will define global $\infty$-categories as $\infty$-category-valued presheaves on a certain $\infty$-category $\Glo$ of compact Lie groups, which we will now construct.

\begin{definition}
    We write $\topCat$ for the 1-category of small topologically enriched $1$-categories. For a compact Lie group $G$, we denote by $\BGcat{G}\in\topCat$ the associated topological groupoid with one object.  
\end{definition}

\begin{rk}
    The category $\topCat$ is cartesian closed \cite{kelly-enriched}*{§2.3}, and the internal homs $\cat{Fun}_\text{cont}(\cat{C},\cat{D})$ can be explicitly described as follows (see §2.2 of \emph{op.\ cit.}): the objects of $\cat{Fun}_\text{cont}(\cat{C},\cat{D})$ are the enriched functors and the space of morphisms $F\to G$ is given by the enriched natural transformations, topologized as a subspace of $\prod_{c\in\cat{C}}\maps(F(c),D(c))$. 

    Using this, we may view $\cat{topCat}$ as a category enriched over itself, which we will again denote by $\cat{topCat}$; if at some point we wish to distinguish this enriched category from its underlying $1$-category, we will denote the latter by $\cat{topCat}_0$. The $1$-functor $\cat{Fun}_\text{cont}(-,-)$ then canonically upgrades to an enriched functor $\topCat^\op\times\topCat\to\topCat$.
\end{rk}

\begin{ex}\label{ex:glo-mapping-space}
    For any compact Lie groups $G,H$ we may describe the enriched category $\cat{Fun}_\text{cont}(\BGcat{G},\BGcat{H})$ explicitly as follows: its objects are the Lie group homomorphisms $f\colon G\to H$, and the space $\maps(f_1,f_2)$ of maps between two homomorphisms $f_1,f_2$ is given by the subspace of $H$ consisting of those $h\in H$ such that $f_2(g)=h\cdot f_1(g)\cdot h^{-1}$ for all $g\in G$. In particular, $\maps(f,f)$ is the centralizer $C(\im f)$ of the image of $f$.
\end{ex}

\begin{constr}
    Given any functor $f\colon\cat{topCat}_0\to\cat{SSet}$ preserving finite products, we can change the enrichment along $f$ to obtain an $\cat{SSet}$-enrichment of $\cat{topCat}_0$. We will apply this to the composite
    \[
        \cat{topCat}_0\xrightarrow{\;\Ntop\;}\cat{QCat}\xrightarrow{\;\core\;}\cat{Kan},
    \]
    where $\Ntop$ denotes the topological nerve (\cite{HTT}*{Definition~1.1.5.5}) and $\core$ is the functor sending a quasi-category to its maximal Kan subcomplex. We write $\cat{topCat}_\Delta$ for the resulting simplicial (in fact, $\cat{Kan}$-enriched) category.
\end{constr}

\begin{definition}\label{def:Glo}
    We write $\cat{Glo}$ for the full simplicial subcategory of $\cat{topCat}_\Delta$ spanned by the objects $\BGcat{G}$ for compact Lie groups $G$. We will moreover denote by $\Glo\coloneqq N_\Delta(\cat{Glo})$ the associated $\infty$-category.
\end{definition}

\begin{definition}\label{def:gl-infty-cat}
    A \emph{global $\infty$-category} is a functor $\Glo^\op\to\CAT_\infty$ into the very large $\infty$-category of large $\infty$-categories. We write $\CAT_\Glo\coloneqq\Fun(\Glo^\op,\CAT_\infty)$ for the very large $\infty$-category of global $\infty$-categories.
\end{definition}

\begin{rk}
    There are various alternative definitions of $\Glo$ appearing throughout the literature. We will review them in Remark~\ref{rk:Glo}, where we will in particular show that they agree with our definition.
\end{rk}

\begin{warn}
    The enriched functor $\Ntop$ descends to a functor $\Glo\to\Spc$ of $\infty$-categories (sending $\BGcat{G}$ to the classifying space $BG$ of the topological group $G$), and by a theorem of Rezk \cite{rezk-classifying}*{Theorem~1.2} this is fully faithful when we restrict to the full subcategory spanned by the compact Lie groups with abelian identity component. It is however \emph{not} fully faithful on all of $\Glo$.
\end{warn}

\begin{rk}
    The $\infty$-category $\Glo$ receives a functor from the $1$-category of compact Lie groups and continuous group homomorphisms, induced by the inclusion $\cat{topCat}_0\hookrightarrow\cat{topCat}_\Delta$. Thus, a global $\infty$-category $\Cc$ in particular encodes a family of $\infty$-categories $\Cc(G)$ indexed by compact Lie groups together with compatible restriction functors for all Lie group homomorphisms. In light of Example~\ref{ex:glo-mapping-space}, we can think of the remaining data of a global $\infty$-category as witnessing that restrictions along inner automorphisms are the identity.
\end{rk}

\begin{rk}
    The $1$-category of compact Lie groups admits a topological enrichment via the compact open topology on mapping spaces, yielding a topologically enriched category $\cat{CpctLie}$. We remark purely for motivational purposes that the functor from the previous remark factors through a functor $\Ntop(\cat{CpctLie})\to\Glo$ (although this is not immediately apparent from our construction!), which in turn lifts uniquely to an equivalence $\Ntop(\cat{CpctLie})\to\Glo_*$ to the $\infty$-category of pointed objects in $\Glo$, see \cite{gepnermeier2020equivTMF}*{Remark 2.14}.\footnote{The cited references uses an a priori different definition of $\Glo$ (and denotes it by $\Orb$), but by Remark~\ref{rk:Glo} this is equivalent to our construction.}
\end{rk}
    
Many interesting examples of global $\infty$-categories arise via the following \emph{continuous Borel construction}:

\begin{constr}\label{constr:continuous-Borel}
    Let $\cat{C}$ be a small topological category. Then we have a $\cat{topCat}$-enriched functor $\cat{Fun}_\text{cont}(-,\cat{C})\colon\cat{topCat}^\op\to\cat{topCat}$, and after changing the enrichment we may view this as a simplicially enriched functor $\cat{topCat}_\Delta^\op\to\cat{topCat}_\Delta$. 

    We will now explain how to also make $\Ntop$ into a simplicially enriched functor $\cat{topCat}_\Delta\to\cat{QCat}$, with the simplicial enrichment on the target given by $\maps(\Cc,\Dd)\coloneqq\core\Fun(\Cc,\Dd)$. For this let us first consider the case where we instead take the enrichment on the target via $\maps(\Cc,\Dd)\coloneqq\Fun(\Cc,\Dd)$, and similarly make $\cat{topCat}$ into a $\cat{QCat}$-enriched category by transferring its self-enrichment along $\Ntop$. Then the unenriched functor $\Ntop\colon\cat{topCat}\to\cat{QCat}$ canonically upgrades to a simplicially enriched functor, given on morphism simplicial sets by the Beck--Chevalley maps\footnote{We refer the reader to \cite{CLL_Adams}*{Appendix C} or \cite{CLL_Clefts}*{Appendix A} for general facts about Beck--Chevalley maps.} 
    \begin{equation}\label{eq:Ntop-enrichment}
        \Ntop\cat{Fun}_\text{cont}(\cat{C},\cat{D})\to\Fun(\Ntop(\cat{C}),\Ntop(\cat{D}))
    \end{equation}
    associated to the canonical isomorphism $\Ntop(-\times-)\cong\Ntop(-)\times\Ntop(-)$. Each of the simplicial maps $(\ref{eq:Ntop-enrichment})$ is a functor between quasi-categories, so it restricts to $\core\Ntop\cat{Fun}_\text{cont}(\cat{C},\cat{D})\to\core\Fun(\Ntop(\cat{C}),\Ntop{\cat{D}})$, giving the desired $\cat{Kan}$-enrichment of $\Ntop\colon\cat{topCat}_\Delta\to\cat{QCat}$.

    In summary, we arrive at a $\cat{Kan}$-enriched functor
    \[
        \cat{Glo}^\op\hookrightarrow\cat{topCat}_\Delta\xrightarrow{\;\cat{Fun}_\text{cont}(-,\cat{C})}\cat{topCat}_\Delta\xrightarrow{\;\Ntop\;}\cat{QCat},
    \]
    which after applying simplicial nerves gives us a functor $\Glo^\op\to\Cat_\infty$ of $\infty$-categories. We denote this by $\Ntop(\cat{C}^\dual)$ and call it the \emph{global $\infty$-category of continuous Borel objects} in $\cat{C}$.
\end{constr}

Performing the above construction in a larger universe then produces a global $\infty$-category $\Ntop(\cat{C}^\dual)$ from any (not necessarily small) topological category $\cat{C}$.

\begin{rk}
    Using the enriched functoriality of $\cat{Fun}_\text{cont}(-,\cat{C})$ in $\cat{C}$, one easily upgrades the assignment $\cat{C}\mapsto\Ntop(\cat{C}^\dual)$ from the previous construction to a simplicially enriched functor $\cat{topCat}_\Delta\to\cat{Fun}_\Delta(\cat{Glo}^\op,\cat{QCat})$, with the evident definition on $0$-simplices (i.e.\ $1$-morphisms). Passing to a larger universe and taking simplicial nerves, we then obtain functor 
    \begin{equation}\label{eq:Borcat}
        \Borcat{(-)}\colon N_\Delta(\cat{TOPCAT}_\Delta)\to\Fun(\Glo^\op,\CAT_\infty).
    \end{equation}
\end{rk}

\begin{ex}
    For $\cat{C}=\BGcat{G}$, the continuous Borel construction just recovers the functor $\Glo^\op\to\Spc\subset\Cat_\infty$ represented by $\BGcat{G}\in\Glo$.
\end{ex}

\begin{ex}
    If we apply the above construction to the topological category $\core\cat{L}$ of finite dimensional inner product spaces and $\R$-linear isometric isomorphisms (topologized as orthogonal groups), we obtain a global $\infty$-category $\ul\Vect\coloneqq\Ntop(\core\cat{L}^\dual)$ sending $\BGcat{G}$ to the $\infty$-groupoid of finite-dimensional $G$-representations and $G$-equivariant isometric isomorphisms.
\end{ex}

\begin{ex}
    A common theme throughout this monograph will be that many \emph{universal} examples of global $\infty$-categories arise as Dwyer--Kan localizations of continuous Borel categories, and that this gives rise to universal properties for categories classically studied in equivariant and global homotopy theory. In particular, we will see in Corollary~\ref{cor:equiv-spaces-universal} that localizing the continuous Borel category $\Ntop(\cat{Top}^\flat)$ in each degree $G$ at the \emph{$G$-equivariant weak equivalences} typically considered in unstable equivariant homotopy theory yields the global $\infty$-category freely generated by the point under certain `genuine' colimits. Similarly, Theorems~\ref{thm:equiv-main} and~\ref{thm:stable-main} will provide universal properties for suitable localizations of the continuous Borel category of {orthogonal spectra}, providing a purely $\infty$-categorical approach to stable equivariant and global homotopy theory for compact Lie groups.
\end{ex}

While we will mostly focus on global $\infty$-categories as defined above (recording information for all compact Lie groups), our results often directly transfer to the setting where one restricts attention to a suitable collection of compact Lie groups:

\begin{variant}
    To any collection $\Ff$ of compact Lie groups closed under isomorphisms we may associate the full subcategory $\Glo_\Ff\subset\Glo$ spanned by the objects of $\Ff$, and we will refer to functors $\Glo_\Ff^\op\to\CAT_\infty$ as \emph{$\Ff$-global $\infty$-categories}. Global $\infty$-categories as introduced in Definition~\ref{def:gl-infty-cat} then correspond to the collection $\All$ of all compact Lie groups. On the other hand, in the case of the collection $\Ff=\Fin$ of \emph{finite groups} (i.e.\ discrete compact Lie groups), $\Glo_{\Ff\!in}$ simply is the usual $(2,1)$-category spanned by the (unenriched) groupoids $\mathbb BG$. The resulting $\Fin$-global $\infty$-categories have been extensively studied in \cites{CLL_Global,CLL_Clefts,CLL_Adams,CLL_Spans,LLP}; note that all of these sources refer to $\Fin$-global $\infty$-categories simply as \emph{global \hbox{($\infty$-)}categories}.
\end{variant}

\subsection{Parametrized higher category theory} The notion of an $\Ff$-global $\infty$-category is only a special case of the general notion of a \emph{$T$-$\infty$-category} as introduced in \cite{exposeI}:

\begin{definition}\label{def:T-oo}
    Let $T$ be any small $\infty$-category. A \emph{$T$-$\infty$-category} is a functor $T^\op\to\CAT_\infty$. We write $\CAT_T\coloneqq\Fun(T^\op,\CAT_\infty)$ for the (very large) $\infty$-category of $T$-$\infty$-categories.
\end{definition}

\begin{rk}
    By the universal property of presheaves, any $T$-$\infty$-category admits a unique extension to a limit-preserving functor $\PSh(T)^\op\to\CAT_\infty$, yielding an equivalence $\CAT_T\simeq\Fun^{\text{lim}}(\PSh(T)^\op,\CAT_\infty)$. In what follows, we will often evaluate $T$-$\infty$-categories at general presheaves, which is to be understood as passing to the limit extension first. In particular, we will freely cite results from \cites{martini2021yoneda,martiniwolf2021limits,martiniwolf2022presentable}, which \emph{defines} $T$-$\infty$-categories as limit-preserving functors $\PSh(T)^\op\to\CAT_\infty$ (and, in fact, works more generally with limit preserving functors $\Bb^{\op}\to\CAT_\infty$ for $\infty$-topoi $\Bb$).
\end{rk}

\begin{ex}
    As already alluded to above, $\Ff$-global $\infty$-categories are precisely $\Glo_\Ff$-$\infty$-categories in the sense of Definition~\ref{def:T-oo}.
\end{ex}

\begin{ex}\label{ex:G-oo}
    For any compact Lie group $G$, we define the $\infty$-category $\Orb_G$ as the full subcategory of $\Ntop(\cat{$\bm G$-Top})$ spanned by the transitive $G$-spaces, i.e.\ those of the form $G/H$ for closed subgroups $H\subset G$. We will refer to $\Orb_G$-$\infty$-categories as \emph{$G$-$\infty$-categories}. For finite $G$ (where $\Orb_G$ is just a 1-category), this notion was already introduced in \cite{exposeI}, building on earlier $1$-categorical work of Hill--Hopkins, and this was the example guiding much of the early development of parametrized higher category theory.
\end{ex}

We now discuss some generic examples of parametrized $\infty$-categories:

\begin{ex}
    The inclusion $\Spc\hookrightarrow\CAT_\infty$ induces a fully faithful functor $\PSh(T)\hookrightarrow\CAT_T$ for any small $T$. We will denote the $T$-$\infty$-category associated to a presheaf $X$ by $\ul X$. If $A\in T$ is arbitrary, then we will typically conflate it with its Yoneda image in $\PSh(T)$, and hence again write $\ul A$ for the corresponding represented $T$-$\infty$-category.
\end{ex}

\begin{ex}
    Let $T$ be any small $\infty$-category. The \emph{$T$-$\infty$-category of $T$-spaces} is defined as the functor $\ul\Spc_T\colon A\mapsto \PSh(T)_{/A}$ with functoriality via pullback, i.e.\ this is the cartesian straightening of the pullback of $\ev_1\colon\Ar(\PSh(T))\to\PSh(T)$ along the Yoneda embedding. 
\end{ex}

\begin{rk}
    Note that since $\PSh(T)$ is a topos, the cartesian unstraightening of $\ev_1$ is a limit preserving functor, i.e.~the limit extension of $\ul\Spc_T$ is given by $X\mapsto\PSh(T)_{/X}$ with functoriality via pullback.
\end{rk}

\begin{ex}\label{ex:T-objects}
    For an $\infty$-category $\Ee$, we define the $T$-$\infty$-category $\ul\Ee_T$ of \emph{$T$-objects in $\Ee$} as the composite
    \[
        T^\op\xrightarrow{\;T_{/-}\;}\Cat_\infty^\op\xrightarrow{\;\Fun(-,\Ee)\;}\CAT_\infty,
    \]
    where $T_{/-}\colon T\to\Cat_\infty$ is the cocartesian unstraightening of $\ev_1\colon\Ar(T)\to T$.

    Note that for $\Ee=\Spc$ this recovers the previous example---more precisely, \cite{CLL_Global}*{Remark~2.1.16} shows that the left adjoints $\PSh(T_{/A})\to\PSh(T)_{/A}$ extending $y_{/A}\colon T_{/A}\to \PSh(T)_{/A}$ for $A\in T$ can be assembled into a global functor and that this is an equivalence; we will therefore typically not distinguish between these two perspectives on $\ul\Spc_T$.

    In addition to the case $\Ee=\Spc$, we will also be interested in the case $\Ee=\CAT_\infty$, where we will abbreviate $\ul\CAT_T\coloneqq\ul{(\CAT_\infty)}_T$.
\end{ex}

\subsubsection{$T$-$\infty$-categories of $T$-functors} Below we want to introduce parametrized notions of adjunctions and colimits. Ultimately, these all rely on the fact that $\CAT_T$ is cartesian closed (since $\CAT_\infty$ is so). 

\begin{definition}
    We denote the internal hom between $\Cc,\Dd\in\CAT_T$ by $\ul\Fun_T(\Cc,\Dd)$. If $T$ is clear from the context, we will often also simply write $\ul\Fun(\Cc,\Dd)$.
\end{definition}

\begin{definition}
    Denote by $1\in\PSh(T)$ the terminal object. 
    For $T$-$\infty$-categories $\Cc,\Dd$, we define $\Fun_T(\Cc,\Dd)\coloneqq\ul\Fun_T(\Cc,\Dd)(1)$ (implicitly passing to the limit extension if $T$ does not have a terminal object), and refer to it as the $\infty$-category of \emph{$T$-functors}.  If $T$ is clear from the context, we will also simply write $\Fun(\Cc,\Dd)$.
\end{definition}

\begin{ex}\label{ex:yoneda}
    For any $T$-$\infty$-category $\Kk$, we write $\ul\PSh_T\hskip0pt minus 1pt(\Kk)\hskip0pt minus 1pt\coloneqq\hskip0pt minus 1pt\ul\Fun_T\hskip0pt minus 1pt(\Kk^\op\hskip0pt minus 1.5pt,\hskip0pt minus 1pt\ul\Spc_T\hskip0pt minus 1pt)$ and call it the \emph{$T$-$\infty$-category of (\kern1pt$T$-)presheaves on $\Kk$}. If $\Kk$ is locally small (meaning that each $\Kk(A)$ is locally small), then \cite{martini2021yoneda}*{discussion before Remark~4.7.7} defines a \emph{parametrized Yoneda embedding} $y\colon\Kk\to\ul\PSh_T(\Kk)$, and by Theorem~4.7.8 of \emph{op.\ cit.} this is fully faithful (meaning that $\Kk(A)\to\ul\PSh_T(\Kk)(A)$ is fully faithful for every $A\in T$ or, equivalently, for every $A\in\PSh(T)$).
\end{ex}

\begin{rk}\label{rk:shift-vs-fun}
    If $\Cc$ is any $T$-$\infty$-category and $X\in\PSh(T)$, then $\ul\Fun(\ul X,\Cc)\simeq\Cc(X\times{-})$ naturally in $\Cc$ and $X$, see \cite{CLL_Global}*{Corollary 2.2.9}.
\end{rk} 

\begin{rk}\label{rk:Yoneda-image}
    Using the equivalence from the previous remark, we can easily describe the parametrized Yoneda embedding for $T$-spaces:

    Specializing \cite{martini2021yoneda}*{Proposition 4.7.20${}^\op$}, an object of $\ul\Spc_T(X\times{-})(A)$ corresponding to a map $(u,v)\colon Z\to X\times A$ in $\PSh(T)$ is contained in the Yoneda image if and only if $v$ is an equivalence, and the proof shows that in this case it is the image of $v^{-1}\circ u\colon A\to X$; in other words, we have $y(f)=\big((f,\id)\colon A\to X\times A\big)$ for any map $f\colon A\to X$ of $T$-spaces. Combining this with the non-parametrized Yoneda lemma, we see that the Yoneda embedding $\ul X\to\ul\Spc_T(X\times{-})$ is the unique $T$-functor sending $\id_X$ to the diagonal $\Delta\colon X\to X\times X$.
\end{rk}

\begin{rk}\label{rk:internal-hom-pw}
    Combining the equivalence from Remark~\ref{rk:shift-vs-fun} with the (internalized) adjunction equivalence, we have specific equivalences  
    \[
        \ul\Fun(\Cc,\Dd)(X)\simeq\Fun(\Cc,\ul\Fun(\ul X,\Dd))\simeq\Fun(\Cc,\Dd(X\times{-}))
    \]
    naturally in $\Cc,\Dd\in\CAT_T$ and $X\in\PSh(T)$, also see \cite{CLL_Global}*{Corollary~2.2.9} again. Throughout, we will usually describe objects of $\ul\Fun(\Cc,\Dd)(X)$ in terms of these (fixed) equivalences, i.e.~by specifying functors $\Cc\to\ul\Fun(\ul X,\Dd)$ or $\Cc\to\Dd(X\times{-})$.
\end{rk}

There is also another description of $\ul\Fun(\Cc,\Dd)$ that we will sometimes need:

\begin{constr}
    Let $T$ be any small $\infty$-category, and let $A\in T$ be arbitrary. Then restriction along the forgetful functor $\fgt\colon T_{/A}\to T$ defines a functor $\fgt^*\colon\CAT_T\to\CAT_{T_{/A}}$. This has a right adjoint $\fgt_*$ given by right Kan extension.
\end{constr}

\begin{rk}\label{rk:fgt*}
    The forgetful functor uniquely extends to a cocontinuous functor $\fgt_!\colon\PSh(T_{/A})\to\PSh(T)$, and we may identify $\fgt^*\colon\CAT_T\to\CAT_{T_{/A}}$ with the functor $\smash{\Fun^\text{R}(\PSh(T),\CAT_\infty)\to\Fun^\text{R}(\PSh(T_{/A}), \CAT_\infty)}$ given by restriction along $\fgt_!$. As the right adjoint $\fgt^*\colon\PSh(T)\to\PSh(T_{/A})$ is again cocontinuous, it then follows from (homotopy) 2-functoriality of $\smash{\Fun^\text{R}(-,\CAT_\infty)}$ that $\fgt_*$ is given by restriction along $\fgt^*$, with unit and counit induced analogously. More explicitly, the unit $\eta\colon\Cc\to\fgt_*\fgt^*\Cc$ can be naturally identified with $\pr^*\colon\Cc\to\Cc(A\times{-})$ for $\Cc\in\CAT_T$, while the counit is given by restriction $\Delta^*\colon\Dd\to\Dd((A\times A)\times_A{-})$ along the diagonal for any $T_{/A}$-$\infty$-category $\Dd$.
\end{rk}

\begin{rk}\label{rk:FunT-slice-description}
    Let $A\in T$ and write $\fgt\colon T\to T_{/A}$ for the forgetful functor again. Combining the previous remark with Remark~\ref{rk:internal-hom-pw}, one can identify 
    \begin{equation}\label{eq:FunT-pointwise-slice}
        \ul\Fun_T(\Cc,\Dd)(A)\simeq\Fun_{T_{/A}}(\fgt^*\Cc,\fgt^*\Dd)
    \end{equation}
    naturally in $\Cc,\Dd$. By  \cite{CLL_Global}*{Corollary~2.2.11} and its proof, the resulting equivalence can also be described more directly as follows: the Beck--Chevalley map
    \begin{multline*}
        \fgt^*\ul\Fun_T(\Cc,\Dd)\xrightarrow{\;\text{coev}\;}
        \ul\Fun_{T_{/A}}\big(\fgt^*\Cc,\fgt^*\Cc\times\fgt^*\ul\Fun_T(\Cc,\Dd)\big)\\
        {}\simeq\ul\Fun_{T_{/A}}\big(\fgt^*\Cc,\fgt^*\big(\Cc\times\ul\Fun_T(\Cc,\Dd)\big)\big)\xrightarrow{\;\ev\;}    
        \ul\Fun_{T_{/A}}(\fgt^*\Cc,\fgt^*\Dd)
    \end{multline*}
    associated to the natural equivalence $\fgt^*\Cc\times\fgt^*(-)\simeq\fgt^*(\Cc\times{-})$ is invertible, and $(\ref{eq:FunT-pointwise-slice})$ is obtained by evaluating this equivalence on global sections.
\end{rk}

\subsubsection{Adjunctions} Since $\CAT_\infty$ is cartesian closed, it is in particular tensored and cotensored over itself. Thus, also $\CAT_T$ is tensored and cotensored over $\CAT_\infty$, by performing these operations levelwise. The tensoring allows us to view $\CAT_T$ as a $\CAT_\infty$-enriched $\infty$-category (that is, an $(\infty,2)$-category), with mapping $\infty$-categories given by the $\infty$-categories $\Fun(\Cc,\Dd)$ introduced above. For our purposes, it will be enough that this allows us to upgrade the homotopy 1-category of $\CAT_T$ to a $(2,2)$-category; we will therefore refrain from using $(\infty,2)$-categorical techniques, and beg the reader who would prefer a more highbrow approach for their forgiveness.

Using the (homotopy) enrichment, we obtain an internal notion of \emph{natural transformations} between functors $\Cc\to\Dd$ as morphisms in $\Fun(\Cc,\Dd)$, or equivalently as $T$-functors $\Cc\times[1]\to\Dd$ or $\Cc\to\Fun([1],\Dd)$. In particular, we obtain an internal notion of \emph{adjunction}, consisting of functors $F\colon\Cc\to\Dd$, $G\colon\Dd\to\Cc$ and natural transformations $\eta\colon\id\to GF$ and $\epsilon\colon FG\to\id$ satisfying the triangle identity up to homotopy. Below we will almost exclusively use the following \emph{pointwise criterion} for adjunctions:

\begin{lemma}[See \cite{martiniwolf2021limits}*{Proposition 3.2.9 and Corollary 3.2.11}]
    \label{lemma:pointwise-criterion}
    Let $F\colon\Cc\to\Dd$ be a $T$-functor. Then $F$ admits a right adjoint if and only if $F_A\colon\Cc(A)\to\Dd(A)$ admits a right adjoint $G_A$ for every $A\in T$ such that these right adjoints satisfy the \emph{Beck--Chevalley condition}: for every $f\colon A\to B$ in $T$ the Beck--Chevalley map $f^*G_B\to G_Af^*$ associated to the naturality equivalence $f^*F_B\simeq F_Af^*$ is an equivalence. 
    
    Moreover, in this case the functor $F_X\colon\Cc(X)\to\Dd(X)$ admits a right adjoint $G_X$ more generally for any $X\in\PSh(T)$, and the Beck--Chevalley condition holds more generally for all maps in $\PSh(T)$.\qed
\end{lemma}

\begin{rk}\label{rk:adjunctions-pw-unit-counit}
    Note that if $F\colon\Cc\to\Dd$ and $G\colon\Dd\to\Cc$ are $T$-functors and $\eta\colon\id\to GF,\epsilon\colon FG\to\id$ are natural transformations exhibiting them as adjoint, then $\eta_X\colon\id\to G_XF_X$ and $\epsilon_X\colon F_XG_X\to\id$ exhibit $F_X$ as left adjoint to $G_X$ for any $X\in\PSh(T)$. It follows formally that in the situation of the previous lemma the right adjoint $G$ is given pointwise by the right adjoints $G_X$, with unit and counit similarly given pointwise.
\end{rk}

\begin{rk}\label{rk:adj-via-colims}
    Assume that for any $f\colon A\to B$ in $T$ the restriction functors $f^*\colon\Cc(B)\to\Cc(A)$ and $f^*\colon\Dd(B)\to\Dd(A)$ admit left adjoints $f_!\colon\Cc(A)\to\Cc(B)$ and $f_!\colon\Dd(A)\to\Dd(B)$. By general nonsense about Beck--Chevalley maps, the Beck--Chevalley map $f^*G_B\to G_Af^*$ from Lemma~\ref{lemma:pointwise-criterion} is an equivalence if and only if the Beck--Chevalley map $f_!F_A\to F_Bf_!$ is so. Below, we will see that (under additional assumptions on $\Cc$ and $\Dd$) the latter Beck--Chevalley condition can be interpreted as $F$ preserving certain \emph{parametrized colimits}.
\end{rk}

\subsubsection{Limits and colimits}\label{subsubsec:colimits} Having introduced adjunctions, we can now talk about colimits and limits in $T$-$\infty$-categories.

\begin{definition}
    Let $\Cc,\Kk$ be $T$-$\infty$-categories. We say that $\Cc$ has \emph{$\Kk$-shaped colimits} if the $T$-functor $\const_{\Kk}\colon\Cc\to\ul\Fun(\Kk,\Cc)$ (induced by $\Kk\to 1$) admits a left adjoint $\colim_\Kk$. If also $\Dd$ has $\Kk$-shaped colimits, then we say that a $T$-functor $\Cc\to\Dd$ \emph{preserves $\Kk$-shaped colimits} if the Beck--Chevalley map ${\colim_\Kk}\circ{\ul\Fun(\Kk,F)}\to F\circ{\colim_\Kk}$ is invertible. 

    Dually, we define what it means to have \emph{$\Kk$-shaped limits} and to preserve them.
\end{definition}

\begin{rk}
    In light of Remark~\ref{rk:adjunctions-pw-unit-counit}, the Beck--Chevalley condition from the previous definition can be rephrased as saying that the ordinary Beck--Chevalley map $\colim_\Kk(A)\circ\ul\Fun(\Kk,F)(A)\to F(A)\circ\colim_\Kk(A)$ is invertible for every $A\in T$.
\end{rk}

\begin{ex}[See \cite{CLL_Global}*{Lemmas~2.3.9 and~2.3.10}]\label{ex:existenceconstantcolimits}
    Let $K$ be an ordinary $\infty$-category, and write $\Kk\coloneqq\mathop\const K$ for the corresponding constant $T$-$\infty$-category. Then a $T$-$\infty$-category $\Cc$ has $\Kk$-shaped colimits if and only if it factors through the non-full subcategory $\CAT_\infty(K)\subset\CAT_\infty$ of $\infty$-categories with $K$-shaped colimits and functors preserving these colimits; moreover, in this case also its limit extension factors through $\CAT_\infty(K)$.

    Likewise, a $T$-functor $F\colon\Cc\to\Dd$ preserves $\Kk$-shaped colimits if and only if $F(A)\colon\Cc(A)\to\Dd(A)$ preserves $K$-shaped colimits for every $A\in T$ or equivalently $A\in\PSh(T)$.
\end{ex}

\begin{definition}
    We say that a $T$-$\infty$-category $\mathcal C$ has an \emph{initial object} if it has $(\const\,{\emptyset})$-colimits. 
    Dually, a $T$-$\infty$-category $\mathcal C$ has a \emph{terminal object} if it has $(\const\,{\emptyset})$-indexed limits. 
\end{definition}

By the above example, a $T$-$\infty$-category $\Cc$ admits a terminal or initial object if and only if $\mathcal C(A)$ admits a terminal or initial object for all $A\in T$ and $f^*\colon\mathcal C(B)\to\mathcal C(A)$ preserves terminal or initial objects, respectively, for all $f\colon A\to B$ in $T$. 

\begin{definition}\label{defi:fiberwise-colims}
    We say that a $T$-$\infty$-category $\Cc$ is \emph{fiberwise cocomplete} if it has $(\const\,K)$-shaped colimits for all small $\infty$-categories $K$, i.e.\ if it factors through the non-full subcategory $\CAT_\infty^\text{cc}\subset\CAT_\infty$. Analogously, we say that $F\colon\Cc\to\Dd$ is \emph{fiberwise cocontinuous} if and only if it preserves $\const\,K$-colimits for all $K\in\Cat_\infty$.
\end{definition}

While it might seem natural to define a \emph{cocomplete} $T$-$\infty$-category $\Cc$ as one having $\Kk$-shaped colimits for all small $T$-$\infty$-categories, it turns out that this is not yet a well-behaved notion, see e.g.\ \cite{martiniwolf2021limits}*{Remark~5.1.2 and Example~5.4.12}. Instead one defines:

\begin{definition}
    Let $\cat{U}\subset\ul\CAT_T$ be any full $T$-subcategory. We say that a $T$-$\infty$-category has \emph{$\cat{U}$-colimits} if the $T_{/A}$-$\infty$-category $\fgt_A^*\Cc$ has $\Kk$-shaped colimits for every $A\in T$ and $\Kk\in\cat{U}(A)\subset\CAT_{T_{/A}}$. If also $\Dd$ has $\cat{U}$-colimits, then we say that a $T$-functor $F\colon\Cc\to\Dd$ is $\cat{U}$-cocontinuous if $\fgt_A^*F$ preserves $\Kk$-shaped colimits for all $A\in T$ and $\Kk\in\cat{U}(A)$.

    Dually, we define what it means to have \emph{$\cat{U}$-limits} and to be \emph{$\cat{U}$-continuous}.
\end{definition}

\begin{definition}\label{defi:T-cc}
    A $T$-$\infty$-category is called \emph{(\kern1pt$T$-)cocomplete} if it has $\ul\Cat_T$-colimits, and a $T$-functor is called \emph{(\kern1pt$T$-)cocontinuous} if it is $\ul\Cat_T$-cocontinuous.
\end{definition}

For $T=\Glo$, we will also speak of \emph{global cocompleteness} and \emph{global cocontinuity}.

\begin{rk}[See \cite{martiniwolf2021limits}*{Propositions~4.3.1 and~4.3.2}]
    Throughout, we will use without further mention that parametrized colimits in functor categories are `pointwise.' More precisely: If $\Kk$ is any $T$-$\infty$-category and $\Cc$ is $\cat{U}$-cocomplete, then $\ul\Fun(\Kk,\Cc)$ is again $\cat{U}$-cocomplete. Moreover, if $\Kk\to\Ll$ is any functor, then the restriction $\ul\Fun(\Ll,\Cc)\to\ul\Fun(\Kk,\Cc)$ is $\cat{U}$-cocontinuous, as is postcomposition with any $\cat{U}$-cocontinuous $\Cc\to\Dd$.
\end{rk}

In practice, one often uses a pointwise description of cocompleteness, similarly to the pointwise criterion of adjunctions. This is particularly easy to state for $\cat{U}\subset\ul\Spc_T$:

\begin{lemma}[See \cite{CLL_Global}*{Lemma 2.3.14 and Remark~2.3.15}]\label{lemma:pointwise-crit-cocomplete}
    Let $\cat{U}\subset\ul\Spc_T\subset\ul\CAT_T$ be any full subcategory. Then a $T$-$\infty$-category $\Cc$ is $\cat{U}$-cocomplete if and only if the following conditions are satisfied:
    \begin{enumerate}
        \item For every $A\in T$ and $(u\colon X\to A)\in\cat{U}(A)\subset\PSh(T)_{/A}$ the restriction $u^*\colon\Cc(A)\to\Cc(X)$ admits a left adjoint $u_!$.
        \item For every $A,B\in T$ and every pullback square
        \[
            \begin{tikzcd}
                Y\arrow[d,"v"']\arrow[r,"g"]\arrow[dr,pullback] & X\arrow[d,"u"]\\
                B\arrow[r,"f"'] & A
            \end{tikzcd}
        \]
        in $\PSh(T)$ such that $u\in\cat{U}(A)$ (whence $v\in\cat{U}(B)$), the Beck--Chevalley map $v_!g^*\to f^*u_!$ is invertible.
    \end{enumerate}
    Moreover, in this case conditions (1) and (2) hold more generally without the representability condition, i.e.\ for all $A,B\in\PSh(T)$.\qed
\end{lemma}

\begin{lemma}[See \cite{CLL_Global}*{Lemma 2.3.16}]\label{lemma:pointwise-criterion-cocontinuous}
    Let $\cat{U}\subset\ul\Spc_T$ be any full subcategory, and let $F\colon\Cc\to\Dd$ be a $T$-functor between $\cat{U}$-cocomplete $T$-$\infty$-categories. Then the following are equivalent:
    \begin{enumerate}
        \item $F$ is $\cat{U}$-cocontinuous.
        \item The Beck--Chevalley map $u_!F\to Fu_!$ is invertible for every $A\in\PSh(T)$ and $u\in\cat{U}(A)$.
        \item The previous condition  holds for all $A\in T$.\qedhere\qed
    \end{enumerate}
\end{lemma}

Together with the notion of fiberwise cocompleteness from Definition~\ref{defi:fiberwise-colims}, this is enough to describe $T$-cocompleteness:

\begin{prop}[See \cite{martiniwolf2021limits}*{Proposition~5.4.1}]
    A $T$-$\infty$-category $\Cc$ is $T$-cocomplete if and only if it is fiberwise cocomplete and $\ul\Spc_T$-cocomplete. Analogously, a functor $F\colon\Cc\to\Dd$ of $T$-cocomplete $\infty$-categories is $T$-cocontinuous if and only if it is fiberwise cocontinuous and $\ul\Spc_T$-cocontinuous.\qed
\end{prop}

\begin{rk}
    As a consequence of Remark~\ref{rk:adj-via-colims}, left adjoint functors of $T$-cocomplete $\infty$-categories are $T$-cocontinuous. More generally, any left adjoint of $\cat{U}$-cocomplete $T$-$\infty$-categories is $\cat{U}$-cocontinuous, see \cite{martiniwolf2021limits}*{Proposition 5.2.5}.
\end{rk}

There is also a notion of \emph{parametrized Kan extensions}:

\begin{thm}[See \cite{martiniwolf2021limits}*{Theorem~6.3.5 and Corollary~6.3.7}]\label{thm:Kan-extension}
   Let $f\colon\Aa\to\Bb$ be a $T$-functor such that $\Aa$ is small and $\Bb$ is locally small, and let $\Cc$ be any $T$-cocomplete $T$-$\infty$-category. Then the restriction $f^*\colon\ul\Fun_T(\Bb,\Cc)\to\ul\Fun_T(\Aa,\Cc)$ admits a left adjoint $f_!$. If $f$ is fully faithful, so is $f_!$. \qed
\end{thm}

Next, we introduce $T$-$\infty$-categories of $T$-cocontinuous functors:

\begin{constr}\label{constr:Fun-T-cc}
    Let $\Cc,\Dd$ be $T$-cocomplete $T$-$\infty$-categories. We write $\ul\Fun^\text{$T$-cc}(\Cc,\Dd)$ for the subfunctor of $\ul\Fun(\Cc,\Dd)\colon\PSh(T)^\op\to\CAT_\infty$ given at an object $X\in\PSh(T)$ by the full subcategory of $\ul\Fun(\Cc,\Dd)(X)$ spanned by the $T$-cocontinuous functors $\Cc\to\ul\Fun(\ul X,\Dd)$. \cite{CLL_Global}*{Proposition 2.3.28} shows that $\ul\Fun^\text{$T$-cc}(\Cc,\Dd)$ is a limit-preserving functor $\PSh(T)\to\CAT_\infty$, i.e.\ a $T$-$\infty$-category.
\end{constr}

\begin{rk}
    The above follows \cite{CLL_Global}, whereas \cite{martiniwolf2021limits} instead defines $\ul\Fun^\text{$T$-cc}(\Cc,\Dd)$ in terms of the description from Remark~\ref{rk:FunT-slice-description}. For the equivalence of these two approaches see \cite{CLL_Global}*{Proposition 2.3.26}.
\end{rk}

Using this, we can now state the universal property of parametrized presheaves:

\begin{thm}[See \cite{martiniwolf2021limits}*{Theorem 7.1.1}]\label{thm:univ-prop-PSh}
    Let $\Kk$ be a small $T$-$\infty$-category and let $\Cc$ be a $T$-cocomplete $T$-$\infty$-category. Then the functor $y^*\colon\ul\Fun(\ul\PSh(\Kk),\Cc)\to\ul\Fun(\Kk,\Cc)$ given by restriction along the Yoneda embedding from Example~\ref{ex:yoneda} induces an equivalence 
    \[
        \ul\Fun^\textup{$T$-cc}(\ul\PSh(\Kk),\Cc)\iso\ul\Fun(\Kk,\Cc).\qednow
    \]
\end{thm}

If $\Kk=1$, then the Yoneda embedding $1\to\ul\PSh(1)\simeq\ul\Spc_T$ is the inclusion of the terminal object as a consequence of Remark~\ref{rk:Yoneda-image}. Thus, we get:

\begin{cor}\label{cor:univ-prop-Spc-T}
    Let $\Cc$ be a $T$-cocomplete $T$-$\infty$-category. Then evaluation at the terminal object defines an equivalence
    \[
        \ul\Fun^\textup{$T$-cc}(\ul\Spc_T,\Cc)\iso\Cc.\qednow
    \]
\end{cor}

\subsubsection{Presentable $T$-$\infty$-categories} We now review the notion of \emph{parametrized presentability}:
\begin{definition}
    A $T$-$\infty$-category $\Cc$ is called \emph{fiberwise presentable} if it factors through the non-full subcategory $\PrL\subset\CAT_\infty$ of presentable $\infty$-categories and left adjoint functors. We say that $\Cc$ is \emph{(\kern1pt$T$)-presentable} if it is fiberwise presentable and $T$-cocomplete (Definition~\ref{defi:T-cc}).
\end{definition}

\begin{ex}
    The $T$-$\infty$-category $\ul\Spc_T$ is $T$-presentable, see~\cite{CLL_Global}*{Example 2.4.4}. More generally, if $\Kk$ is any small $T$-$\infty$-category, then $\ul\PSh_T(\Kk)$ is $T$-presentable by Example~2.4.5 of \emph{op.\ cit.}
\end{ex}

\begin{rk}
    \cite{martiniwolf2022presentable} instead defines $T$-presentable $T$-$\infty$-categories as suitable Bousfield localizations of presheaf $T$-$\infty$-categories; for the equivalence to the above definition see \cite{martiniwolf2022presentable}*{Theorem~2.4.2.5}.
\end{rk}

As a consequence of Remark~\ref{rk:adj-via-colims}, we have the following parametrized version of the Adjoint Functor Theorem, see also \cite{martiniwolf2022presentable}*{Proposition 2.4.3.1~and 2.4.3.3}:

\begin{prop}\label{prop:adj-functor-thm}
    \begin{enumerate}
        \item Let $\Cc$ be $T$-presentable and let $\Dd$ be $T$-cocomplete and locally small. A $T$-functor $F\colon\Cc\to\Dd$ is a left adjoint if and only if it is $T$-cocontinuous.
        \item Let $\Cc,\Dd$ be $T$-presentable. A $T$-functor $G\colon\Cc\to\Dd$ is a right adjoint if and only if it is $T$-continuous and fiberwise accessible (meaning that each $\Cc(A)\to\Dd(A)$ is accessible).\qedhere\qed
    \end{enumerate}
\end{prop}

\begin{definition}
    We denote by $\PrLT\subset\CAT_T$ the non-full subcategory with objects the $T$-presentable $T$-$\infty$-categories and with morphisms the left adjoint (or, equivalently, $T$-cocontinuous) functors. 

    We will moreover write $\PrRT\subset\CAT_T$ for the non-full subcategory with the same objects, but whose morphisms are the \emph{right} adjoint $T$-functors.
\end{definition}

\begin{rk}\label{rk:FunL}
    If $\Dd,\Ee$ are presentable, then the Adjoint Functor Theorem similarly lets us identify $\ul\Fun^\text{$T$-cc}(\Dd,\Ee)(A)\simeq\Fun^\text{$T$-cc}(\Dd,\Ee(A\times{-}))$ with the subcategory of $\Fun(\Dd,\Ee(A\times{-}))$ spanned by the \emph{left adjoint} $T$-functors, and we will therefore typically write $\ul\Fun_T^\text{L}(\Dd,\Ee)\coloneqq\ul\Fun_T^\text{$T$-cc}(\Dd,\Ee)$.
\end{rk}

\begin{prop}\label{prop:FunL-presentable}
    If $\Dd,\Ee$ are $T$-presentable, then so is $\ul\Fun_T^\textup{L}(\Dd,\Ee)$.
    \begin{proof}
        If $\Dd\simeq\ul\PSh_T(\Dd_0)$, then $\ul\Fun_T^\textup{L}(\Dd,\Ee)\simeq\ul\Fun_T(\Dd_0,\Ee)$ by Theorem~\ref{thm:univ-prop-PSh}, and so the claim is an instance of \cite{martiniwolf2022presentable}*{Corollary 2.4.2.7}.

        In the general case, we use \cite{martiniwolf2022presentable}*{Remark 2.4.4.12} to find a pushout
        \[
            \begin{tikzcd}
                \Aa\arrow[r, "f"]\arrow[dr,"\ulcorner"{very near end},phantom]\arrow[d,"p"'] & \Bb\arrow[d,"q"]\\
                \Cc\arrow[r,"g"'] & \Dd
            \end{tikzcd}
        \]
        in $\PrLT$ such that $\Aa,\Bb,\Cc$ are presheaf categories. Taking $\hom\big(-,\ul\Fun(\ul A\times[n],\Ee)\big)$ for varying $A\in T$ and $n\ge 0$ we then arrive at a pullback
        \[
            \begin{tikzcd}
                \ul\Fun_T^\text{L}(\Dd,\Ee)\arrow[r,"g^*"]\arrow[d,"q^*"'] & \ul\Fun^\text{L}(\Cc,\Ee)\arrow[d,"p^*"]\\
                \ul\Fun_T^\text{L}(\Bb,\Ee)\arrow[r,"f^*"'] & \ul\Fun_T^\text{L}(\Aa,\Ee)
            \end{tikzcd}
        \]
        in $\CAT_T$. By the above special case, all objects except possibly the top left one are presentable; moreover, $f^*$ and $p^*$ are cocontinuous since $\Fun^\text{$T$-cc}(\Xx,\Yy)\subset\Fun(\Xx,\Yy)$ is closed under $T$-colimits for any $T$-cocomplete $\Xx,\Yy$ by \cite{martiniwolf2022presentable}*{Lemma 2.6.1.3}. The claim now follows since $\PrLT\hookrightarrow\CAT_T$ preserves small limits by \cite{martiniwolf2022presentable}*{Proposition 2.4.4.10}.
    \end{proof}
\end{prop}

We also want to consider a (typically non-presentable) $T$-$\infty$-category $\ul\Fun^\text{R}(\Cc,\Dd)$ of right adjoint functors for presentable $\Cc,\Dd$:

\begin{constr}
    Let $\Cc,\Dd$ be presentable. We define $\ul\Fun_T^\text{R}(\Cc,\Dd)\subset\ul\Fun_T(\Cc,\Dd)$ as the full $T$-subcategory given in degree $A\in T$ by the right adjoint functors $F\colon\Cc\to\ul\Fun(\ul A,\Dd)$. Note that this is indeed a $T$-subcategory: if $f\colon B\to A$ is any map in $T$, then $f^*F$ is the composite
    \[
        \Cc\xrightarrow{\;F\;}\ul\Fun(\ul B,\Dd)\xrightarrow{\;f^*\;}\ul\Fun(\ul A,\Dd),
    \]
    and the second functor is a right adjoint by Theorem~\ref{thm:Kan-extension}.
\end{constr}

\begin{lemma}\label{lemma:PrR-restrict}
    For any $A\in T$, the adjunction $\fgt^*\colon\CAT_{T}\rightleftarrows\CAT_{T_{/A}}\noloc\fgt_*$ restricts to an adjunction $\smash{\PrR_{T}}\rightleftarrows\PrR_{T_{/A}}$.
\end{lemma}

In particular, we see that $\ul\Fun^\text{R}(\Cc,\Dd)(A)$ corresponds in the description of $\ul\Fun(\Cc,\Dd)(A)$ from Remark~\ref{rk:FunT-slice-description} to the full subcategory spanned by the right adjoints $\fgt^*\Cc\to\fgt^*\Dd$, which is the definition used in \cite{martiniwolf2021limits}*{Remark 3.3.6}.

\begin{proof}
    The functor $\fgt^*$ sends presentable $T$-$\infty$-categories to presentable $T_{/A}$-$\infty$-categories by \cite{martiniwolf2022presentable}*{Remark 2.4.2.10}. Moreover, $\fgt_*$ sends $T_{/A}$-cocomplete $T_{/A}$-$\infty$-categories to $T$-cocomplete $T$-$\infty$-categories by \cite{CLL_Global}*{Proposition 2.3.26}, and it preserves fiberwise presentability by the explicit description given in Remark~\ref{rk:fgt*}. Thus, both adjoints restrict accordingly on the level of objects.
    
    Next, we will show that they also restrict accordingly on the level of morphisms. As $\fgt^*$ commutes with the $\CAT_\infty$-tensors, it preserves adjunctions; in particular, it sends right adjoint $T$-functors to right adjoint $T_{/A}$-functors. Moreover, we may conclude from this that the right adjoint $\fgt_*$ preserves $\CAT_\infty$-\emph{co}tensors, so it also preserves right adjoints. 

    It remains to show that the unit and counit are right adjoints; we will argue for the unit, the argument for the counit being analogous. By Remark~\ref{rk:fgt*}, we may identify the unit $\Cc\to\fgt_*\fgt^*\Cc$ with the map $\pr^*\colon\Cc\to\Cc(A\times{-})$. By $\ul\Spc_T$-cocompleteness, each $\pr^*\colon\Cc(B)\to\Cc(A\times B)$ admits a left adjoint $\pr_!$, and these satisfy the Beck--Chevalley condition as the square
    \[
        \begin{tikzcd}
            A\times B\arrow[d,"\pr"']\arrow[r,"A\times f"] &[.25em] A\times B'\arrow[d,"\pr"]\\
            B\arrow[r,"f"'] & B'
        \end{tikzcd}
    \]
    is a pullback for any $f\colon B\to B'$. Thus, the existence of the adjoint follows from Lemma~\ref{lemma:pointwise-criterion}.
\end{proof}

For (much) later use we record the following consequence:

\begin{cor}\label{cor:PrLT-ambidextrous-fgt}
    Let $A\in T$ arbitrary, let $\Ee$ be a presentable $T$-$\infty$-category, and let $F\colon\Cc\to\Dd$ be a left adjoint functor of presentable $T_{/A}$-$\infty$-categories. Then $\fgt_*F$ is again a left adjoint, and there exists a diagram of spaces
    \begin{equation}\label{diag:wrong-way-adjunction}
        \begin{tikzcd}
            \hom_{\PrLT}(\fgt_*\Dd,\Ee)\arrow[d,"\sim"']\arrow[r, "{-}\circ\fgt_*F"] &[.25em] \hom_{\PrLT}(\fgt_*\Cc,\Ee)\arrow[d,"\sim"]\\
            \hom_{\PrL_{T_{/A}}}(\Dd,\fgt^*\Ee)\arrow[r,"{-}\circ F"'] & \hom_{\PrL_{T_{/A}}}(\Cc,\fgt^*\Ee)
        \end{tikzcd}
    \end{equation}
    commuting up to (unspecified) homotopy.
\end{cor}

A more careful (and more involved) argument would show that $\fgt_*$ and $\fgt^*$ restrict to functors between $\PrLT$ and $\PrL_{T_{/A}}$, with $\fgt_*$ being \emph{left} adjoint to $\fgt^*$.

\begin{proof}
    Write $U\colon\Dd\to\Cc$ for the right adjoint of $F$. By the proof of Lemma~\ref{lemma:PrR-restrict}, the $T$-functor $\fgt_*U$ will then be right adjoint to $\fgt_*F$.

    We now recall from \cite{martiniwolf2022presentable}*{Proposition 2.4.4.7} that there exists an equivalence $\PrLT\simeq(\PrRT)^\op$, given on objects by the identity and sending a left adjoint $T$-functor to its right adjoint. We therefore see that the top arrow in $(\ref{diag:wrong-way-adjunction})$ agrees up to conjugation by equivalences with 
    $\fgt_*U\circ{-}\colon\hom_{\PrRT}(\Ee,\fgt_*\Dd)\to\hom_{\PrRT}(\Ee,\fgt_*\Cc)$.
    By Lemma~\ref{lemma:PrR-restrict} we may further identify this up to equivalence with 
    \[
        U\circ{-}\colon\smash{\hom_{\PrR_{T_{/A}}}(\fgt^*\Ee,\Dd)\to\hom_{{\PrR_{T_{/A}}}}}(\fgt^*\Ee,\Cc),
    \]
    which another application of \cite{martiniwolf2022presentable}*{Proposition 2.4.4.7} then identifies (up to equivalence) with the bottom horizontal map in $(\ref{diag:wrong-way-adjunction})$.
\end{proof}

\section{Model categories of equivariant and global spaces}\label{sec:G-global-spaces}
As the main result of this first part, we want to describe the free globally presentable global $\infty$-category in terms of unstable $G$-global homotopy theory as studied in \cites{schwede2018global,barrero2021}. In this section, we will introduce the relevant pointset models of $G$-global spaces and establish their basic properties; the universal property will then be proven in the next section.

\subsection{Equivariant spaces} We begin with a quick reminder on the (genuine) homotopy theory of $G$-spaces for compact Lie groups $G$.

\begin{definition}
    Let $G$ be a compact Lie group, and let $\Ff$ be any collection of subgroups of $G$. We call a map $f\colon X\to Y$ of topological $G$-spaces an \emph{$\mathcal F$-equivariant weak equivalence} or \emph{$\Ff$-equivariant fibration} if $f^H\colon X^H\to Y^H$ is a weak homotopy equivalence or Serre fibration, respectively, for every $H\in\Ff$.
\end{definition}

\begin{prop}\label{prop:equiv-model-structure}
    The $\Ff$-equivariant weak equivalences and fibrations participate in a model structure on the category $\cat{$\bm G$-Top}$ of topological spaces with an action of $G$. This model structure is proper, topological, and cofibrantly generated with generating cofibrations
    \[
        \{G/H\times\partial D^n\hookrightarrow G/H\times D^n : H\in\Ff,n\ge0\}
    \]
    and generating acyclic cofibrations
    \[
        \{G/H\times D^n\hookrightarrow G/H\times D^n\times[0,1] : H\in\Ff,n\ge0\}.
    \]
    \begin{proof}
        \cite{schwede2018global}*{Proposition B.7} constructs this model structure and proves that it is topological, proper, and cofibrantly generated. The description of the generating (acyclic) cofibrations is pointed out in the proof of \emph{loc.\ cit.}
    \end{proof}
\end{prop}

\begin{ex}\label{ex:G-spaces}
    In the special case where $\Ff=\All$ is the collection of all closed subgroups of $G$, we will refer to the above as the \emph{($G$-)equivariant model structure} and to its weak equivalences as the \emph{($G$-)equivariant weak equivalences}. We write $\myS_G$ for the Dwyer--Kan localization of the $1$-category $\cat{$\bm G$-Top}$ at the $G$-equivariant weak equivalences, and call it the \emph{$\infty$-category of $G$-spaces}.
\end{ex}

\begin{ex}
    If $G,H$ are compact Lie groups, then we write $\mathcal G_{G,H}$ for the collection of \emph{graph subgroups}, i.e.~subgroups of $G\times H$ of the form \[\Gamma_{K,\phi}\coloneqq\{(k,\phi(k)):k\in K\}\] for closed subgroups $K\subset G$ and continuous group homomorphisms $\phi\colon K\to H$. We will typically denote $\Gamma_{K,\phi}$-fixed points as $(-)^\phi$.

    Note that a closed subgroup of $G\times H$ belongs to $\mathcal G_{G,H}$ if and only if it intersects $1\times H$ trivially. In particular, $\mathcal G_{G,H}$ is different from the collection $\mathcal G_{H,G}$ of subgroups intersecting $G\times 1$ trivially, unless $G=H=1$.
\end{ex}

\subsubsection{Functoriality} The following functoriality properties for equivariant spaces are well-known, but we include proofs for the reader's convenience:

\begin{lemma}\label{lemma:restr-equiv-rQ}
    Let $\alpha\colon G\to G'$ be a continuous group homomorphism, and let $\Ff,\Ff'$ be collections of closed subgroups of $G$ and $G'$, respectively, such that $\alpha(H)\in\Ff'$ for every $H\in\Ff$. Then 
    \[
        \alpha_!\colon\cat{$\bm G$-Top}\rightleftarrows\cat{$\bm{G'}$-Top}\noloc\alpha^*
    \]
    is a Quillen adjunction with homotopical right adjoint for the $\Ff$-model structure on the left and the $\Ff'$-model structure on the target.
    \begin{proof}
        Note that we have an equality of functors $(-)^H\circ\alpha^*=(-)^{\alpha(H)}$. Thus, it follows immediately from the definitions that $\alpha^*$ preserves weak equivalences and fibrations.
    \end{proof}
\end{lemma}

\begin{ex}\label{ex:graph-restr-rQ}
    If $H$ is any further Lie group, then the previous lemma shows that
    \[
        \alpha_!\coloneqq(H\times\alpha)_!\colon\cat{$\bm{(H\times G)}$-Top}\rightleftarrows\cat{$\bm{(H\times G')}$-Top}\noloc\alpha^*\coloneqq(H\times\alpha)^*
    \]
    is a Quillen adjunction for the $\mathcal G_{H,G}$-model structure on the left and the $\mathcal G_{H,G'}$-model structure on the target. If $\alpha$ is injective, it is also a Quillen adjunction for the $\mathcal G_{G,H}$- and $\mathcal G_{G',H}$-model structures.
\end{ex}

\begin{lemma}\label{lemma:restr-equiv-lQ}
    Let $\alpha\colon G\to G'$ be a homomorphism of compact Lie groups, and let $\Ff,\Ff'$ be \emph{families} of closed subgroups (i.e.\ non-empty collections of subgroups closed under subconjugates) of $G$ and $G'$, respectively, such that $\Ff$ contains every closed subgroup $H\subset G$ for which $\alpha(H)\in\Ff'$.
    Then
    \[
        \alpha^*\colon\cat{$\bm{G'}$-Top}\rightleftarrows\cat{$\bm G$-Top} \noloc\alpha_*
    \]
    is a Quillen adjunction with homotopical right adjoint.
    \begin{proof}
        It suffices to show that $\alpha^*$ sends generating (acyclic) cofibrations to (acyclic) cofibrations. By the explicit description of the generating cofibrations, and since $\cat{$\bm G$-Top}$ is topological, it then further suffices to show that $\alpha^*(G'/H')$ is $\Ff$-cofibrant for every $H'\in\Ff'$.

        By Illman's triangulation theorem \cite{illman-triangle}*{Theorem~7.1}, $G'$ admits the structure of a $(G\times H')$-CW-complex, where $g\in G$ acts via left multiplication by $\alpha(g)$ and $h\in H'$ acts via right multiplication by $h^{-1}$. Note that $\alpha(g)g'h^{-1}=g'$ if and only if $h=g'^{-1}\alpha(g)g'$, i.e.~the $(G\times H)$-isotropy of any $g'\in G'$ is contained in the collection $\mathcal G$ of all subgroups of the form $\Gamma_{K,\alpha'}$ where $K\subset G$ and $\alpha'\colon K\to G'$ is conjugate to $\alpha$ and has image in $H'$. In particular, every cell in the $(G\times H')$-CW-structure on $G'$ needs to have isotropy contained in $\mathcal G$, i.e.\ the $(G\times H)$-space $G'$ is $\mathcal G$-cofibrant. To finish the proof, it then suffices by Lemma~\ref{lemma:restr-equiv-rQ} to observe that for any $\Gamma_{K,\alpha'}\in\mathcal G$ the image under the quotient map $G\times H'\to G$ is precisely $K$, which satisfies $\alpha'(K)\subset H'$, so that $\alpha(K)\in\Ff'$ and hence $K\in\mathcal F$ by assumption.
    \end{proof}
\end{lemma}

\begin{ex}\label{ex:graph-restr-lQ}
    Dually to Example~\ref{ex:graph-restr-rQ}, the previous lemma shows that
    \[
        \alpha^*\colon\cat{$\bm{(H\times G')}$-Top}\rightleftarrows\cat{$\bm{(H\times G)}$-Top} \noloc\alpha_*
    \]
    is always a Quillen adjunction for the $\mathcal G_{G',H}$- and $\mathcal G_{G,H}$-model structures, and that it is a Quillen adjunction for the $\mathcal G_{H,G'}$- and $\mathcal G_{H,G}$-model structures if $\alpha$ is injective.
\end{ex}

\subsubsection{Homotopy pushouts}\label{subsubsec:homotopy-po-equiv}
A priori, checking whether a commutative square in $\cat{$\bm G$-Top}$ is a homotopy pushout in the $\Ff$-model structure involves replacing one of the legs by a cofibration. Working with this directly is often inconvenient as there are rather few cofibrations, and since cofibrations are not preserved by basic operations like restricting the action. We therefore introduce a larger class of maps such that pushouts along them compute homotopy pushouts:

\begin{definition}
    Let $\cat{C}$ be a topological model category. A map $f\colon A\to X$ in $\cat{C}$ is called an \emph{h-cofibration} if the evident map $X\amalg_{A} \bigl([0,1]\otimes A\bigr)\to [0,1]\otimes X$ admits a retraction.
\end{definition}

\begin{rk}\label{rk:h-cof-closure}
    In \cite{schwede2018global}*{Definition~A.28}, h-cofibrations are instead defined as the maps satisfying the left lifting property against a certain class of maps; the equivalence to the above definition is noted directly afterwards. Note moreover that as a collection of maps defined via a left lifting property, the h-cofibrations are automatically closed under pushouts, retracts, coproducts, and transfinite composition.
\end{rk}

\begin{ex}
    By definition, the h-cofibrations in $\cat{Top}$ are precisely the Hurewicz cofibrations; as we are working in compactly generated weak Hausdorff spaces, these are automatically closed \cite{schwede2018global}*{Proposition~A.31}.
\end{ex}

\begin{ex}
    If $G$ is any compact Lie group and $\Ff$ is any collection of closed subgroups, then any $\Ff$-cofibration is an h-cofibration as a consequence of \cite{schwede2018global}*{Corollary A.30(i)}.
\end{ex}

\begin{rk}\label{rk:fixed-and-restr-preserve-h-cof}
    For any homomorphism $\alpha\colon G\to G'$ of compact Lie groups, the restriction functor $\alpha^*$ preserves h-cofibrations as we may simply restrict a retraction of the cylinder inclusion along $\alpha$. Similarly, it follows from \cite{schwede2018global}*{Proposition~B.1(i)} that also the fixed point functor $(-)^G\colon\cat{$\bm G$-Top}\to\cat{Top}$ preserves h-cofibrations.
\end{rk}

\begin{lemma}\label{lemma:h-cof-compute-po-unbased}
    Let $G$ be a compact Lie group, and let
    \begin{equation}\label{diag:po-along-hcof}
       \begin{tikzcd}
            A\arrow[r,"i"]\arrow[d,"f"']\arrow[dr,phantom,"\ulcorner"{very near end}] & B\arrow[d]\\
            C\arrow[r] & D
        \end{tikzcd}
    \end{equation}
    be any pushout in $\cat{$\bm{G}$-Top}$ such that $i$ is an h-cofibration. Then $(\ref{diag:po-along-hcof})$ is a homotopy pushout in the $\Ff$-model structure for any collection $\Ff$ of subgroups.
    \begin{proof}
        It suffices to factor $f$ into a cofibration $A\to C_0$ followed by a weak equivalence, and to observe that the induced map $B\amalg_AC_0\to B\amalg_AC=D$ is an $\Ff$-weak equivalence by \cite{schwede2018global}*{Proposition~B.6}.
    \end{proof}
\end{lemma}

We also have the following pointwise characterization of homotopy pushouts:

\begin{lemma}\label{lemma:homotopy-po-equiv-fixed-points}
    A commutative square 
    \begin{equation}\label{diag:orig-square}
        \begin{tikzcd}
            A\arrow[r]\arrow[d] & B\arrow[d]\\
            C\arrow[r] & D
        \end{tikzcd}
    \end{equation}
    in $\cat{$\bm G$-Top}$ is a homotopy pushout in the $\Ff$-model structure if and only if for every $H\in\Ff$ the resulting square on $H$-fixed points is a homotopy pushout in $\cat{Top}$.
    \begin{proof}
        Pick a factorization of $A\to B$ into a cofibration $i$ followed by an $\Ff$-weak equivalence $f$ as depicted in the following diagram:
        \[
            \begin{tikzcd}
                A\arrow[r,"i"]\arrow[d]\arrow[dr,phantom,"\ulcorner"{very near end}] & B_0\arrow[r,"f","\sim"']\arrow[d] & B\arrow[d]\\
                C\arrow[r] & C_0\arrow[r,dashed,"g"'] & D\rlap.
            \end{tikzcd}
        \]
        By definition, $(\ref{diag:orig-square})$ is a homotopy pushout if and only if $g$ is an $\Ff$-weak equivalence. As $i$ is in particular an h-cofibration, so is $i^H$ for every $H\in\Ff$ by Remark~\ref{rk:fixed-and-restr-preserve-h-cof}; moreover, \cite{schwede2018global}*{Proposition~B.1(i)} shows that if we apply $H$-fixed points to the left-hand square in the above diagram, this is still an (honest) pushout square. Thus, applying the previous lemma for $G=1$, we see that $g^H$ is a weak homotopy equivalence if and only if $(\ref{diag:orig-square})$ induces a homotopy pushout on $H$-fixed points. Letting $H$ vary, this yields precisely the claim.
    \end{proof}
\end{lemma}

\subsection{Global spaces} In this subsection, we will recall Schwede's model of global spaces in terms of \emph{orthogonal spaces} \cite{schwede2018global} as well as Barrero's $G$-global generalization of this \cite{barrero2021}. 

\begin{constr}
    We view inner product spaces of countable dimension as topological spaces by equipping them with the weak topology with respect to the family of finite subspaces, or equivalently with respect to any chosen exhaustive sequence of finite-dimensional subspaces. If $V$ and $W$ are countably-dimensional inner product spaces, then we write $\cat{L}(V,W)$ for the space of linear isometric embeddings, topologized as subspace of the mapping space $\maps(V,W)$ with the compact-open topology. We write $\Lhat$ for the topological category of countable-dimensional inner product spaces and linear isometric embeddings, and $\cat{L}\subset\Lhat$ for the full subcategory of \emph{finite-dimensional} inner product spaces. By \cite{schwede_orbispaces_2020}*{Proposition A.2}, these are indeed topological categories (i.e., the composition maps are continuous).
\end{constr}

\begin{rk}\label{rk:L-top-explicit}
    We can describe the topologies on the above mapping spaces more concretely as follows, see \cite{schwede_orbispaces_2020}*{Proposition~A.5(iii)}:
    For a finite-dimensional inner product space $U$, denote by $\O(U)$ its ortogonal group. 
    If $V,W$ are finite-dimensional inner product spaces and $\cat{L}(V,W)\not=\emptyset$, then we have for any linear isometric embedding $i\in\cat{L}(V,W)$ a homeomorphism \[\O(W)/\O(W-i(V))\to\cat{L}(V,W),\, [A]\mapsto Ai,\] where $W-i(V)$ denotes the orthogonal complement of $i(V)$ in $W$. In other words, $\cat{L}(V,W)$ is topologized as a Stiefel manifold. If $V$ is finite-dimensional and $W$ is arbitrary, then $\cat{L}(V,W)$ is topologized as the colimit of the mapping spaces $\cat{L}(V,W_0)$ over the finite-dimensional subspaces $W_0\subset W$. Finally, if also $V$ is arbitrary, then $\cat{L}(V,W)$ is topologized as the inverse limit of the spaces $\cat{L}(V_0,W)$ for finite-dimensional subspaces $V_0\subset V$.
\end{rk}
    
\begin{definition}
    The topological category of \emph{orthogonal spaces} is defined as the enriched functor category $\cat{$\cat L$-Top}\coloneqq\cat{Fun}_\cont(\cat{L},\cat{Top})$.
\end{definition}

\begin{rk}
    The references \cites{schwede2018global,barrero2021} use the notation \textit{spc} for the category of orthogonal spaces.
\end{rk}

\subsubsection{The level model structure}
The $G$-global model structure on the category $\cat{$\bm G$-$\cat{L}$-Top}$ of orthogonal spaces with a continuous action by a compact Lie group $G$ is constructed in two steps, beginning with a \emph{level model structure}:

\begin{prop}[See \cite{barrero2021}*{Theorem~A.2 and Lemma~A.7}]\label{prop:orth-lvl}
    Let $G$ be a compact Lie group. There is a unique model structure on $\cat{$\bm G$-$\cat{L}$-Top}$ where a map $f\colon X\to Y$ is a weak equivalence or fibration if and only if $f(V)^\phi\colon X(V)^\phi\to Y(V)^\phi$ is a weak homotopy equivalence or Serre fibration, respectively, for every compact Lie group $H$, every \emph{faithful} $H$-representation $V$, and every continuous group homomorphism $\phi\colon H\to G$.

    We call this the \emph{$\bm G$-global level model structure} and its weak equivalences the \emph{$\bm G$-global level weak equivalences}. This model category is proper, topological, and cofibrantly generated with generating cofibrations
    \[\hskip-0.6pt\hfuzz=0.6pt
        \{\cat{L}(V,-)\times_\phi G\times(\partial D^n\hookrightarrow D^n) : \text{$V$\kern1pt faithful $H$\kern-1pt-representation}, \phi\colon H\to G,n\ge0\}
    \]
    and generating acyclic cofibrations
    \[
        \{\cat{L}(V,-)\times_\phi G\times( D^n\hookrightarrow D^n\times I) : \text{$V$\kern1pt faithful $H$\kern-1pt-representation},\phi\colon H\to G,n\ge0\},
    \]
    where $\cat{L}(V,-)\times_\phi G$ denotes the quotient of $\cat{L}(V,-)\times G$ by the diagonal of the right $H$-action on $\cat{L}(V,-)$ via precomposition and the right $H$-action on $G$ via $\phi$.\qed
\end{prop}

\begin{rk}
    Unravelling definitions, we can equivalently describe the level weak equivalences as those maps $f\colon X\to Y$ such that $f(V)$ is a $\mathcal G_{\O(V),G}$-weak equivalence in $\cat{$\bm{(\O(V)\times G)}$-Top}$, and analogously for the fibrations.
\end{rk}

\begin{rk}\label{rk:level-cof-are-h-cof}
    Any cofibration in the level model structure is in particular an h-cofibration, see \cite{barrero2021}*{Lemma 3.15} or \cite{schwede2018global}*{Corollary~A.30(iii)}.
\end{rk}

\begin{lemma}\label{lemma:orth-spc-lvl-inj-restr}
    Let $\alpha\colon G\to G'$ be an \emph{injective} homomorphism of compact Lie groups. Then $\alpha^*\colon\cat{$\bm{G'}$-$\cat{L}$-Top}\to\cat{$\bm G$-$\cat{L}$-Top}$ is left Quillen with respect to the level model structures.
    \begin{proof}
        We may equivalently show that $\alpha_*$ is right Quillen, which follows from Example~\ref{ex:graph-restr-rQ} applied levelwise.
    \end{proof}
\end{lemma}

\begin{lemma}\label{lemma:orth-spc-ev-lrQ}
    For any $V\in\cat{L}$, the functor $\ev_V\colon\cat{$\bm G$-\cat{L}-Top}\to\cat{$\bm{(\O(V)\times G)}$-Top}$ is both left and right Quillen for the level model structure on the source and the $\mathcal G_{\O(V),G}$-model structure on the target.
    \begin{proof}
        The functor $\ev_V$ has both adjoints (given by enriched Kan extension), and it is clear from the definitions that it preserves weak equivalences and fibrations. It only remains to show that it also preserves cofibrations, for which it again suffices to show that $\cat{L}(U,V)\times_\phi G=\phi_!(\cat{L}(U,V))$ is a cofibrant in the $\mathcal G_{\O(V),G}$-model structure for every homomorphism $\phi\colon H\to G$ and any finite faithful $H$-representation $U$. 
        
        By Example~\ref{ex:graph-restr-rQ}, it will suffice to treat the case $\phi=\id$. If $\cat{L}(U,V)=\emptyset$, there is nothing to prove, otherwise we may identify $U$ with a subspace of $V$. Then the homeomorphism $\O(V)/\O(V-U)\to\cat{L}(U,V)$ from Remark~\ref{rk:L-top-explicit} is just the restriction map, which is easily seen to be $(H\times G)$-equivariant with respect to the action on the source via pre- and postcomposition. By Illman's theorem, $\O(V)$ admits the structure of an $(\O(V-U)\times \O(V)\times H)$-CW-complex; moreover, the isotropy of any cell of this decomposition has to be contained in $\mathcal G_{\O(V-U)\times \O(V),H}$ as $H$ acts freely on $\O(V)$ by faithfulness of $U\subset V$. Thus, $\O(V)$ is $\mathcal G_{\O(V-U)\times \O(V),H}$-cofibrant, and hence $\cat{L}(U,V)\cong\O(V)/\O(V-U)$ is $\mathcal G_{\O(V),H}$-cofibrant by another application of Example~\ref{ex:graph-restr-rQ}
    \end{proof}
\end{lemma}

\subsubsection{The $G$-global model structure} The $G$-global model structure we are actually after will be obtained from the level model structure by Bousfield localization.

\begin{definition}\label{def:g-gl-fib-orth-spc}
    A map $f\colon X\to Y$ in $\cat{$\bm G$-$\cat{L}$-Top}$ is called a \emph{$G$-global fibration} if it is a $G$-global level fibration and for all $H$-representations $V,W$ such that $V$ is faithful the naturality square
    \[
        \begin{tikzcd}
            X(V)\arrow[r]\arrow[d,"f(V)"'] & X(V\oplus W)\arrow[d,"f(V\oplus W)"]\\
            Y(V)\arrow[r] & Y(V\oplus W)
        \end{tikzcd}
    \]
    is a homotopy pullback in the $\mathcal G_{H,G}$-equivariant model structure on $\cat{$\bm{(H\times G)}$-Top}$.

    A $G$-orthogonal space $X$ is called \emph{static} if $X\to1$ is a $G$-global fibration, i.e.\ if the transition map $X(V)\to X(V\oplus W)$ is a $\mathcal G_{H,G}$-weak equivalence for all $V,W$ as above.
\end{definition}

\begin{thm}[See \cite{barrero2021}*{Theorem A.20}]
    The $G$-global level model structure admits a Bousfield localization whose fibrations are the $G$-global fibrations and whose fibrant objects are the static orthogonal $G$-spaces (Definition~\ref{def:g-gl-fib-orth-spc}). We call this the \emph{$\bm G$-global model structure} and its weak equivalences the \emph{$\bm G$-global weak equivalences}. This model structure is again proper, topological, and cofibrantly generated.\qed
\end{thm}

\begin{rk}\label{rk:generating-ayclic-between-cof}
    We refer the reader to \cite{barrero2021}*{Construction~A.14} for an explicit construction of a set of generating acyclic cofibrations. All that we will need later is that these are maps between \emph{cofibrant} objects.
\end{rk}

\begin{definition}
    We write $\myS_\text{$G$-gl}$ for the Dwyer--Kan localization of $\cat{$\bm G$-$\cat{L}$-Top}$ at the $G$-global weak equivalences, and call it the \emph{$\infty$-category of $G$-global spaces}.
\end{definition}

In general, the $G$-global weak equivalences are a bit inconvenient to describe due to some pointset level issues, see \cite{barrero2021}*{Definition~3.2}. For nice enough $G$-orthogonal spaces there is however a simple description, which we will now recall.

\begin{constr}\label{constr:ev-infty}
    Let $X$ be an orthogonal space and let $V\in\Lhat$ be an inner product space of countable dimension. We define
    \[
        \overline{X}(V)\coloneqq\mathop{\textup{colim}}\limits_{{\scriptstyle U\subset V}\atop\scriptstyle\textup{ finite dimensional}} X(U),
    \]
    where the colimit runs over the poset of finite dimensional subspaces $U\subset V$, with the evident transition maps. Note that for finite dimensional $V$, this poset has a terminal object given by $V$ itself, and so we have a canonical isomorphism $\overline{X}(V)\cong X(V)$.
    
    If $f\colon V\to V'$ is any linear isometry, then the maps $U\to f(U)\subset V'$ for $U\subset V$ finite assemble into a map $\overline{X}(V)\to\overline{X}(V')$, and as shown as part of \cite{schwede_orbispaces_2020}*{Construction~3.2}, this makes $\overline{X}\colon\Lhat\to\cat{Top}$ into a topologically enriched functor. Together with the obvious functoriality of $X\mapsto \overline{X}$ in maps of orthogonal spaces, we therefore obtain a functor
    \begin{equation}\label{eq:ext-orth-spc}
        \cat{$\cat L$-Top}=\FUN_\cont(\cat{L},\cat{Top})\to\FUN_\cont(\Lhat,\cat{Top})\eqqcolon\cat{$\Lhat$-Top}
    \end{equation}
    which is a section of the restriction functor. 
\end{constr}

\begin{lemma}
    The functor $(\ref{eq:ext-orth-spc})$ is topologically enriched.
    \begin{proof}
        Recall that mapping spaces in enriched functor categories are topologized as a subspace of the product. By the universal property of the latter, the claim therefore amounts to saying that
        \[
            \maps_{\cat{$\cat{L}$-Top}}(X,Y)\to\maps_{\cat{Top}}(\overline{X}(V),\overline{Y}(V))
        \]
        is continuous for all $X,Y\in\cat{$\cat{L}$-Top}$ and $V\in\Lhat$. Adjoining over and plugging in the definition of $\overline{X}(V)$, this is equivalent to continuity of
        \[X(U)\times\maps_{\cat{$\cat{L}$-Top}}(X,Y)\to \overline{Y}(V)\]
        for all finite-dimensional $U\subset V$. But this map is simply the composite
        \[
            X(U)\times\maps(X,Y)\xrightarrow{\;X(U)\times\ev_U\;}X(U)\times\maps(X(U),Y(U))\xrightarrow{\;\ev\;} Y(U)\xrightarrow{\;\;}\overline{Y}(V),
        \]
        and each of these maps is evidently continuous.
    \end{proof}
\end{lemma}

\begin{constr}
    By the above lemma, $(\ref{eq:ext-orth-spc})$ lifts to a functor $\cat{$\bm G$-$\cat{L}$-Top}\to\cat{$\bm G$-$\Lhat$-Top}$ for every compact Lie group $G$. If now $H$ is a compact Lie group and $\Vv$ is an orthogonal $H$-representation of countable dimension, then postcomposition with evaluation at $\Vv$ yields a functor $\ev_\Vv\colon\cat{$\bm G$-$\cat{L}$-Top}\to\cat{$\bm{(H\times G)}$-Top}, X\mapsto\overline{X}(\Vv)$.
\end{constr}

\begin{convention}
    From now on we will no longer distinguish notationally between an orthogonal space and its extension to $\Lhat$. In particular, we will denote $\ev_\Vv$ on objects by $X\mapsto X(\Vv)$.
\end{convention}

If $X$ is an orthogonal space, then the values of (its extension) at infinite dimensional inner product spaces are in general not homotopically meaningful as sequential colimits in $\cat{Top}$ do not preserve weak homotopy equivalences. We will now single out a class of orthogonal spaces for which this problem goes away:

\begin{definition}\label{def:closed-orth-spc}
    An orthogonal space $X$ (possibly equipped with an action of a compact Lie group) is called \emph{closed} if for every linear isometric embedding $V\to W$ of finite-dimensional inner product spaces $V,W$, $X(V)\to X(W)$ is a closed embedding.
\end{definition}

\begin{lemma}[See \cite{barrero2021}*{Proposition 3.5}]\label{lemma:we-between-closed} 
    A map $f\colon X\to Y$ of \emph{closed} orthogonal $G$-spaces is a $G$-global weak equivalence if and only if the induced map $f(\Uu_H)\colon X(\Uu_H)\to Y(\Uu_H)$ is a $\mathcal G_{H,G}$-weak equivalence for every compact Lie group $G$ and any \emph{complete $\bm H$-universe} $\Uu_H$, i.e.\ any countable $H$-representation into which any finite-dimensional $H$-representation embeds.\qed
\end{lemma}

There is also a sufficient supply of closed orthogonal $G$-spaces:

\begin{lemma}\label{lemma:cofibrant-closed}
    Every cofibrant $G$-orthogonal space is closed.
    \begin{proof}
        By Lemma~\ref{lemma:orth-spc-lvl-inj-restr} it suffices to treat the case $G=1$, which is in turn the content of \cite{schwede2018global}*{Proposition~1.2.11(iii)}.
    \end{proof}
\end{lemma}

\subsubsection{A presheaf model}
As the main result of \cite{schwede_orbispaces_2020}, Schwede proved a \emph{global Elmendorf theorem} in the form of a chain of Quillen equivalences between the global model structure on $\cat{$\cat{L}$-Top}$ and a certain model category of enriched presheaves. For our purposes, the following $\infty$-categorical version will be more convenient:

\begin{prop}\label{prop:global-Elmendorf}
    Let $G$ be a compact Lie group, and write $\mathfrak O_\textup{$G$-gl}\subset\myS_\textup{$G$-gl}$ for the full subcategory spanned by the objects of the form $\cat{L}(V,-)\times_\phi G$ for compact Lie groups $H$, continuous group homomorphisms $\phi\colon H\to G$, and faithful $H$-representations $V$. Then the inclusion $\mathfrak O_\textup{$G$-gl}\hookrightarrow\myS_\textup{$G$-gl}$ extends to an equivalence $\PSh(\mathfrak O_\textup{$G$-gl})\iso\myS_\textup{$G$-gl}$.
\end{prop}

The proof requires some preparations. We write $\myS_\textup{$G$-gl}^\textup{closed}$ for the Dwyer--Kan localization of the full subcategory of $\cat{$\bm G$-$\cat{L}$-Top}$ spanned by the closed orthogonal $G$-spaces, so that the inclusion induces an equivalence $\myS_\text{$G$-gl}^\text{closed}\iso\myS_\text{$G$-gl}$, with inverse given by cofibrant replacement.

\begin{lemma}\label{lemma:fixed-point-corep}
    Let $\phi\colon H\to G$ be a homomorphism of compact Lie groups, let $\Uu_H$ be a complete $H$-universe, and let $V$ be any faithful $H$-representation. Then $\cat{L}(V,-)\times_\phi G\in\myS_\textup{$G$-gl}$ corepresents the composite
    \begin{equation}\label{eq:orthspc-fixed-points-oo}
        \myS_\textup{$G$-gl}\iso\myS_\textup{$G$-gl}^\textup{closed}\xrightarrow{\;(-)^\phi\circ\ev_{\Uu_H}\;}\Spc.
    \end{equation}
    Moreover, this functor is cocontinuous.
    \begin{proof}
        It is clear that $\cat{L}(V,-)\times_\phi G$ is cofibrant and corepresents
        \begin{equation}\label{eq:fixed-points-corep-orthspc}
            \cat{$\bm G$-$\cat{L}$-Top}\xrightarrow{\;\ev_V\;}\cat{$\bm{(H\times G)}$-Top}\xrightarrow{\;(-)^\phi\;}\cat{Top}
        \end{equation}
        in the enriched sense. Corollary~\ref{cor:corep-Top} therefore shows that it corepresents the right derived functor $\myS_\text{$G$-gl}\to\Spc$ of $(\ref{eq:fixed-points-corep-orthspc})$ in the $\infty$-categorical sense.
      
        To see that this right derived functor agrees with $(\ref{eq:orthspc-fixed-points-oo})$, let us fix an $H$-equivariant linear isometric embedding $i\colon V\to\Uu_H$. It will then suffice to show that for any cofibrant-fibrant orthogonal $G$-space $X$ the map $X(i)^\phi\colon X(V)^\phi\to X(\Uu_H)^\phi$ is a weak homotopy equivalence. For this, we pick an exhaustive filtration $i(V)\eqqcolon U_0\subset U_1\subset\cdots\subset\Uu_H$ of $U$ by finite-dimensional subrepresentations: for example, we can first choose a filtration by arbitrary finite subspaces $V_i$ and then inductively apply \cite{schwede_orbispaces_2020}*{Proposition~A.7(i)} to find a finite-dimensional $G$-subrepresentation $U_n$ containing $V_n+U_{n-1}$. By fibrancy of $X$, each transition map $X(U_n)^\phi\to X(U_{n+1})^\phi$ is a weak homotopy equivalence, and it is moreover a closed embedding by cofibrancy. Thus, the induced map $X(V)^\phi\cong X(U_0)^\phi\to \colim_n X(U_n)^\phi$ is again a weak homotopy equivalence. Moreover, \cite{schwede2018global}*{Proposition~B.1(ii)} shows that the natural map from the right-hand side to $(\colim_n X(U_n))^\phi$ is a homeomorphism, finishing the proof that $\cat{L}(V,-)/G$ corepresents $(\ref{eq:orthspc-fixed-points-oo})$.

        To prove that $(\ref{eq:orthspc-fixed-points-oo})$ is cocontinuous, it suffices to show that it preserves coproducts and pushouts. The first statement is clear, while for the second one it suffices by Lemma~\ref{lemma:homotopy-po-equiv-fixed-points} that $\ev_{\Uu_H}\colon\cat{$\bm G$-$\cat{L}$-Top}\to\cat{$\bm{(H\times G)}$-Top}$ sends any pushout square
        \[
            \begin{tikzcd}
                A\arrow[r,"i"]\arrow[d] & B\arrow[d]\\
                C\arrow[r] & D
            \end{tikzcd}
        \]
        where $i$ is a cofibration and $A,B,C,D$ are closed to a homotopy pushout in the $\mathcal G_{H,G}$-model structure. By Remark~\ref{rk:level-cof-are-h-cof}, $i$ is in particular an h-cofibration, hence so is $i(\Uu_H)$ by \cite{schwede2018global}*{Corollary~A.30(ii)}. The claim then follows from Lemma~\ref{lemma:h-cof-compute-po-unbased}.
    \end{proof}
\end{lemma}

\begin{proof}[Proof of Proposition~\ref{prop:global-Elmendorf}]
    As $\mathfrak S_\text{$G$-gl}$ is cocomplete and locally small (being the $\infty$-category associated to a model category), the inclusion $\mathfrak O_\textup{$G$-gl}\hookrightarrow\myS_\textup{$G$-gl}$ extends to a left adjoint functor $\PSh(\mathfrak O_\textup{$G$-gl})\hookrightarrow\myS_\textup{$G$-gl}$, with right adjoint given by the restricted Yoneda embedding. The previous lemma together with Lemma~\ref{lemma:we-between-closed} shows that the right adjoint is conservative, so it only remains to show that the left adjoint $\PSh(\mathfrak O_\text{$G$-gl})\to\myS_\text{$G$-gl}$ is fully faithful. By \cite{HTT}*{Proposition~5.1.6.10} it suffices for this that the functor $\hom(X,-)\colon\myS_\text{$G$-gl}\to\Spc$ is cocontinuous for every $X\in\mathfrak O_\text{$G$-gl}$, which holds by the previous lemma.
\end{proof}

\begin{rk}\label{rk:mathfrak-O-gl-vs-Stefan}
    Let us compare our category $\mathfrak O_\gl$ to (the topological nerve of) the global indexing category $\cat{O}_\gl$ from Schwede's global Elmendorf theorem \cite{schwede_orbispaces_2020}. 
    
    We write $\mathfrak L$ for the topological monoid $\cat{L}(\R^\infty,\R^\infty)$ and $\cat{$\bm{\mathfrak L}$-Top}$ for the topological category of $\mathfrak L$-spaces. The topological category $\cat{O}_\gl$ is defined as the full subcategory of $\cat{$\bm{\mathfrak L}$-Top}$ spanned by the orbits $\mathfrak L/H$ for so-called \emph{universal} subgroups; moreover, every compact Lie group is isomorphic to a universal subgroup of $\mathfrak L$ by \cite{schwede_orbispaces_2020}*{Proposition~1.5}. The objects $\mathfrak L/H$ are cofibrant-fibrant in the projective global model structure of \cite{schwede_orbispaces_2020}*{Proposition~1.11} and hence Lemmas~\ref{lemma:fgt-simpl-structure} and~\ref{lemma:localization-on-homs} show that we may identify $\Ntop(\cat{O}_\gl)$ with the full subcategory $\mathfrak L\kern-.4pt\mathfrak O_\gl$ of the localization $\mathfrak L\kern-.4pt\mathfrak S_\gl$ of $\cat{$\bm{\mathfrak L}$-Top}$ at the weak equivalences  of the global projective model structure.
    
    In \cite{schwede_orbispaces_2020}*{Theorem~3.9}, Schwede then shows that the functor $\ev_{\R^\infty}\colon\cat{$\cat{L}$-Top}\to\cat{$\bm{\mathfrak L}$-Top}$ induces an equivalence $\myS_\gl^\text{cof}\to\mathfrak L\kern-.4pt\mathfrak S_\gl$. If $H$ is any universal subgroup and $V$ is any faithful $H$-representation, then Proposition A.10 of \emph{op.\ cit.} shows that we may identify $\cat{L}(V,\R^\infty)/H\simeq\mathfrak L/H$ (via restriction along any equivariant embedding $V\hookrightarrow\R^\infty$). In particular, $\ev_{\R^\infty}$ induces an equivalence $\mathfrak O_\gl\iso\mathfrak L\kern-.4pt\mathfrak O_\gl\simeq\Ntop(\cat{O}_\gl)$.
\end{rk}

\subsubsection{Change of group}
We will now discuss how the above model structures relate to each other as the group $G$ varies:

\begin{lemma}\label{lemma:restr-right-Quillen}
    Let $\alpha\colon G\to G'$ be a homomorphism of compact Lie groups. Then
    \[
        \alpha_!\colon\cat{$\bm{G}$-$\cat{L}$-Top}\rightleftarrows\cat{$\bm{G'}$-$\cat{L}$-Top}\noloc \alpha^*
    \]
    is a Quillen adjunction with homotopical right adjoint.
    \begin{proof}
        By \cite{barrero2021}*{Lemma~3.7(vii)}, the functor $\alpha^*$ is homotopical. On the other hand, Example~\ref{ex:graph-restr-rQ} shows that if $H$ is any compact Lie group, then $\alpha^*\colon\cat{$\bm{(H\times G')}$-Top}\to \cat{$\bm{(H\times G)}$-Top}$ is right Quillen for the $\mathcal G_{H,G}$- and $\mathcal G_{H,G'}$-model structures; applying this for $H=\O(V)$, we see that $\alpha^*$ preserves \emph{level} fibrations. Similarly, if $f\colon X\to Y$ is even a $G'$-global fibration, then this shows that
        \[
            \begin{tikzcd}
                \alpha^*X(V)\arrow[r]\arrow[d,"\alpha^*f(V)"'] & \alpha^*X(V\oplus W)\arrow[d,"\alpha^*f(V\oplus W)"]\\
                \alpha^*Y(V)\arrow[r] & \alpha^*Y(V\oplus W)
            \end{tikzcd}
        \]
        is a homotopy pullback for all $H$-representations $V,W$ such that $V$ is faithful, i.e.\ $\alpha^*f$ is a $G'$-global fibration.
    \end{proof}
\end{lemma}

\begin{lemma}\label{lemma:restr-inj-left-Quillen}
    Let $\alpha\colon G\to G'$ be an \emph{injective} homomorphism of compact Lie groups. Then the adjunction
    \[
        \alpha^*\colon\cat{$\bm{G'}$-$\cat{L}$-Top}\rightleftarrows\cat{$\bm{G}$-$\cat{L}$-Top} \noloc \alpha_*
    \]
    is a Quillen adjunction.
    \begin{proof}
        By the previous lemma, the left adjoint $\alpha^*$ preserves weak equivalences, while Lemma~\ref{lemma:orth-spc-lvl-inj-restr} shows that it preserves cofibrations.
    \end{proof}
\end{lemma}

Lemma~\ref{lemma:restr-right-Quillen} shows that $\alpha_!$ admits a left derived functor $\cat{L}\alpha_!\colon\myS_\text{$G$-gl}\to\myS_\text{$G'$-gl}$, obtained by precomposition with a cofibrant replacement. Unfortunately, the cofibrant replacement provided by the small object argument is completely intractable in practice; to get some control over the left derived functor, we will therefore exhibit a strictly larger class of objects on which the pointset level functor $\alpha_!$ is still homotopical. We begin with a statement for equivariant spaces.

\begin{lemma}\label{lemma:G-CW-Hausdorff}
    Let $G$ be a compact Lie group and let $\Ff$ be any collection of closed subgroups of $G$ closed under conjugation. Then any $\Ff$-cofibrant $G$-space $X$ is Hausdorff, and every isotropy group of $X$ is contained in $\Ff$.
    \begin{proof}
        By Lemma~\ref{lemma:restr-equiv-lQ}, the underlying non-equivariant space of $X$ is cofibrant, hence Hausdorff. For the second statement, note that the collection of all spaces with isotropy groups contained in $\Ff$ is closed under coproduct, pushout and transfinite composition along closed embeddings, as well as under passage to subspaces. As moreover the source and target of any generating cofibration for the $\Ff$-model structure have isotropy groups contained in $\Ff$ (here we use the closure under conjugation), the lemma follows.
    \end{proof}
\end{lemma}

\begin{lemma}\label{lemma:equiv-quotient}
    Let $\alpha\colon G\to G'$ be a homomorphism of compact Lie groups, and let $H$ be any other compact Lie groups. Then $\alpha_!\colon\cat{$\bm{(H\times G)}$-Top}\to \cat{$\bm{(H\times G')}$-Top}$ sends $\mathcal G_{H,G}$-weak equivalences between Hausdorff $(H\times G)$-spaces with free $\ker(\alpha)$-action to $\mathcal G_{H,G'}$-weak equivalences.
    \begin{proof}
        Consider $G'$ as a $(H\times G\times G')$-space with $G'$ acting via left multiplication, $H$ acting trivially, and with $g\in G$ acting via right multiplication by $\alpha(g)^{-1}$. Then $\alpha_!$ factors as
        \[
            \cat{$\bm{(H\times G)}$-Top}\xrightarrow{\;G'\times{-}\;}\cat{$\bm{(H\times G\times G')}$-Top}\xrightarrow{\;(-)/G\;}\cat{$\bm{(H\times G')}$-Top}.
        \]
        It is clear that the first functor preserves Hausdorffness and sends $\mathcal G_{H,G}$-weak equivalences to $\mathcal G_{H,G'\times G}$-weak equivalences. Moreover, if $X$ is $\ker(\alpha)$-free, then $G'\times X$ is $G$-free, as every point of $G'$ has $G$-isotropy $\ker(\alpha)$. It therefore suffices to show that the second functor sends $\mathcal G_{H,G\times G'}$-weak equivalences between $G$-free Hausdorff spaces to $\mathcal G_{H,G'}$-weak equivalences.

        If $X$ is any $G$-free space, $K\subset H$ is any closed subgroup, and $\phi\colon K\to G'$ is any homomorphism, then \cite{schwede2018global}*{Proposition~B.17} shows that the tautological map
        \begin{equation}\label{eq:splitting-of-fp}
            \begin{aligned}
                \smash{\coprod_{\alpha\colon K\to G}X^{(\alpha,\phi)}/C(\im\alpha)}&\longrightarrow (X/G)^\phi\\
                [x]&\longmapsto[x]
            \end{aligned}
        \end{equation}
        is a homeomorphism, where the coproduct runs over $G$-conjugacy classes of continuous group homomorphisms $K\to G$, and $C(\im\alpha)$ denotes the centralizer of the image of $\alpha$ in $G$. Thus, if $f\colon X\to Y$ is a $\mathcal G_{H,G}$-weak equivalence of $G$-free Hausdorff spaces, then $(f/G)^\phi$ splits as a coproduct of the maps $f^{(\alpha,\phi)}/C(\im\alpha)$. By assumption, each $f^{(\alpha,\phi)}$ is a non-equivariant weak equivalence. As $C(\im\alpha)\subset G$ acts freely on $X$ and $Y$, Lemma~\ref{lemma:quotient-we} therefore shows that also the quotient $f^{(\alpha,\phi)}/C(\im\alpha)$ is a weak homotopy equivalence.
    \end{proof}
\end{lemma}

\begin{prop}\label{prop:quotient-we-global}
    Let $\alpha\colon G\to G'$ be any homomorphism of compact Lie groups, and let $f\colon X\to Y$ be a $G$-global weak equivalence in $\cat{$\bm G$-$\cat{L}$-Top}$ such that $X(V)$ and $V$ are $\ker(\alpha)$-free and Hausdorff for every $V\in\cat{L}$. Then $\alpha_!f$ is a $G'$-global weak equivalence.
    \begin{proof}
        Picking functorial cofibrant replacements in the $G$-global level model structure yields a commutative diagram
        \[
            \begin{tikzcd}
                X'\arrow[r,"f'"]\arrow[d] & Y'\arrow[d]\\
                X\arrow[r,"f"'] & Y
            \end{tikzcd}
        \]
        where $X'$ and $Y'$ are cofibrant and the vertical maps are $G$-global level weak equivalences. By 2-out-of-3, $f'$ is again a $G$-global weak equivalence, and hence $\alpha_!f'$ is a $G'$-global weak equivalence by Lemma~\ref{lemma:restr-right-Quillen}. To complete the proof we will show that $\alpha_!$ sends the vertical maps to $G'$-global level weak equivalences; it will then follow by 2-out-of-3, that $\alpha_!f$ is a $G'$-global weak equivalences.
        
        We argue for the left-hand map, the argument for the right-hand one being analogous. By Lemma~\ref{lemma:orth-spc-ev-lrQ} each $X'(V)$ is $\mathcal G_{\O(V),G}$-cofibrant, and hence in particular $\ker(\alpha)$-free and Hausdorff by Lemma~\ref{lemma:G-CW-Hausdorff}. Applying the previous lemma levelwise therefore shows that $\alpha_!(X'\to X)$ is a $G$-global level weak equivalence, as desired.
    \end{proof}
\end{prop}

In particular, we may compute $\cat{L}f_!(X)$ for an orthogonal  $G$-space $X$ by resolving it by a levelwise Hausdorff $G$-free $X'$.

\medskip
\subsubsection{Multiplicative properties} 
The topological category $\cat{L}$ admits an enriched symmetric monoidal structure given by (orthogonal) direct sum; the associativity, unitality, and symmetry isomorphisms come from the corresponding isomorphisms for the cocartesian symmetric monoidal structure on the category of $\R$-vector spaces. This then induces a topologically enriched symmetric monoidal structure on $\cat{$\cat{L}$-Top}$ via enriched Day convolution; the symmetric monoidal product of this structure is called the \emph{box product} and is denoted by $\boxtimes$.

\begin{lemma}\label{lemma:box-product-bifun-general}
    For any compact Lie groups $G,H$, the box product defines a homotopical left Quillen bifunctor $\cat{$\bm G$-$\cat{L}$-Top}\times\cat{$\bm H$-$\cat{L}$-Top}\to\cat{$\bm{(G\times H)}$-$\cat{L}$-Top}$.
    \begin{proof}
        The pushout-product axiom for cofibrations is verified in \cite{barrero2021}*{Lemma A.6}, while the preservation of weak equivalences is Corollary 3.10 of \emph{op.\ cit.}

        For the pushout product axiom for acyclic cofibrations, we will more generally show that the pushout product of any h-cofibration $i$ (of orthogonal $H$-spaces) with any $G$-global weak equivalence $f$ is a $(G\times H)$-global weak equivalence. As both h-cofibrations and weak equivalences are stable under restriction (see Lemma~\ref{lemma:restr-right-Quillen}), this follows at once by applying \cite{barrero2021}*{Corollary~A.10} to $\triv_Gi$ and $\triv_Hf$.
    \end{proof}
\end{lemma}

\begin{cor}\label{cor:box-product-bifun}
    For any compact Lie group $G$, the functor
    \begin{equation}\label{eq:box-product-single-G}
        {-}\boxtimes{-}\colon\cat{$\bm G$-$\cat{L}$-Top}\times\cat{$\bm G$-$\cat{L}$-Top}\to\cat{$\bm G$-$\cat{L}$-Top}
    \end{equation}
    is homotopical and a left Quillen bifunctor.
    \begin{proof}
        We may factor $(\ref{eq:box-product-single-G})$ as
        \[
            \cat{$\bm G$-$\cat{L}$-Top}\times\cat{$\bm G$-$\cat{L}$-Top}\xrightarrow{\;{-}\boxtimes{-}\;}\cat{$\bm{(G\times G)}$-$\cat{L}$-Top}\xrightarrow{\;\Delta^*\;}\cat{$\bm G$-$\cat{L}$-Top},
        \]
        where $\Delta=(\id,\id)\colon G\to G\times G$ is the diagonal embedding. The first functor is a homotopical left Quillen bifunctor by the previous lemma, while the second one is homotopical and left Quillen by Lemmas~\ref{lemma:restr-right-Quillen} and~\ref{lemma:restr-inj-left-Quillen}.
    \end{proof}
\end{cor}

On the level of model categories, the box product often serves as a better behaved version of the cartesian product---for example, we will see later that the correct global version of the (unreduced) suspension spectrum functor is strong symmetric monoidal on the pointset level for the box product (while it is only \emph{lax} symmetric monoidal for the cartesian product), and there is an interesting homotopy theory of strictly commutative monoids for the box product \cite{schwede2018global}*{Chapter~2}. On the level of $\infty$-categories, however, there is no difference:

\begin{lemma}\label{lemma:box-product-vs-product}
    There exists a natural $G$-global weak equivalence $X\boxtimes Y\to X\times Y$ for any $X,Y\in\cat{$\bm G$-$\cat{L}$-Top}$.
    \begin{proof}
        This follows from \cite{barrero2021}*{Proposition 3.9} by restricting along the homomorphism $\Delta\colon G\to G\times G$ as before.
    \end{proof}
\end{lemma}

\begin{cor}[See also \cite{barrero2021}*{Lemma~3.7(v)}]\label{cor:cart-prod-homotopical}
    If $f\colon X\to Y$ and $g\colon X'\to Y'$ are $G$-global weak equivalences in $\cat{$\bm G$-$\cat{L}$-Top}$, then $f\times g$ is again a $G$-global weak equivalence.\qed
\end{cor}

\subsection{Model categorical aspects of the universal property}\label{subsec:model-cat-prereq}
Our first big goal is to construct an equivalence 
\begin{equation}\label{eq:first-goal}
    \ul\Spc_{\Glo}(\BGcat{G})=\PSh(\Glo)_{/\BGcat{G}}\simeq\myS_\text{$G$-gl}
\end{equation}    
natural in compact Lie groups. The existence of such an equivalence has various consequences for the change-of-group adjunctions discussed above---for example, it implies that for any $\alpha\colon G\to G'$ the endofunctor $\cat{L}\alpha_!\alpha^*$ of $\myS_\text{$G'$-gl}$ is given by the cartesian product with $\cat{L}\alpha_!\alpha^*(1)$, simply because the analogous statement holds for $\ul\Spc_\Glo$. In this subsection, we will establish several such properties, which will turn out to be crucial ingredients in the proof of the comparison $(\ref{eq:first-goal})$ in the next section. We begin with the promised product formula for $\cat{L}\alpha_!\alpha^*$:

\begin{prop}\label{prop:f!f*}
    Let $\alpha\colon G\to G'$ be any homomorphism of compact Lie groups, let $X\in\myS_\textup{$G'$-gl}$ arbitrary, and write $x\colon X\to 1$ for the unique map. Then 
    \[
        \cat{L}\alpha_!\alpha^*X\xrightarrow{(\epsilon,\cat{L}\alpha_!\alpha^*x)} X\times\cat{L}\alpha_!\alpha^*1
    \]
    is an equivalence in $\myS_\textup{$G'$-gl}$.
\end{prop}

The proof will require some preparations:

\begin{lemma}[Projection formula]\label{lemma:proj-formula-underived}
    Let $\alpha\colon G\to G'$ be a homomorphism of compact Lie groups, and let $X\in\cat{$\bm G$-$\cat{L}$-Top}$, $Y\in\cat{$\bm{G'}$-$\cat{L}$-Top}$. Then the map
    \begin{equation}\label{eq:upf}
        (\alpha_!\pr_1, \epsilon\circ \alpha_!\pr_2)\colon \alpha_!(X\times \alpha^*Y)\to \alpha_!X\times Y
    \end{equation}
    is an isomorphism in $\cat{$\bm{G'}$-$\cat{L}$-Top}$.
\end{lemma}

We remark for motivational purposes that $(\ref{eq:upf})$ is just the Beck--Chevalley map
\[
    \alpha_!(X\times \alpha^*Y)\xrightarrow{\;f_!(\eta\times\id)\;}
    \alpha_!(\alpha^*\alpha_!X\times \alpha^*Y)=\alpha_!\alpha^*(\alpha_!X\times \alpha^*Y)\xrightarrow{\;\epsilon\;}\alpha_!X\times \alpha^*Y
\]
associated to the identity $\alpha^*(-)\times \alpha^*Y= \alpha^*(-\times Y)$. Using our usual model for $\alpha_!$, we can also describe it explicitly as the map given pointwise by
\begin{equation}\label{eq:proj-formula}
    \begin{aligned}
        (G'\times X\times \alpha^*Y)_{/G}&\longrightarrow (G'\times X)_{/G}\times Y\\
        [g,x,y]&\longmapsto ([g,x],g.y).
    \end{aligned}
\end{equation}

\begin{proof}[Proof of Lemma~\ref{lemma:proj-formula-underived}]
    As all functors in sight are defined levelwise, it suffices to prove the analogous statement with $\cat{Top}$ in place of $\cat{$\cat{L}$-Top}$. 

    Consider the map $\alpha_!X\times Y\to \alpha_!(X\times \alpha^*Y),([g,x],y)\mapsto [g,x,g^{-1}.y]$. This is well-defined: any other representative of $[g,x]$ is of the form $[g\alpha(h), h^{-1}.x]$ for some $h\in G$, and $[g\alpha(h),h^{-1}.x,\alpha(h)^{-1}g^{-1}.y]=[g,x,g^{-1}.y]$ in $(G'\times X\times \alpha^*Y)/G$ by definition of the $G$-action. It is then clear that this map is inverse to $(\ref{eq:proj-formula})$, so it only remains to establish continuity. As $-\times Y$ preserves quotients, we may equivalently show that the map $G'\times X\times Y\to (G'\times X\times \alpha^*Y)/G, (g,x,y)\mapsto [g,x,g^{-1}.y]$ is continuous, which is clear.
\end{proof}

\begin{lemma}\label{lemma:resolve-pt}
    Let $V$ be a faithful representation of the compact Lie group $G$.\footnote{Such a representation always exists as a consequence of the Peter--Weyl Theorem, see e.g.\ \cite{broecker-tom-dieck}*{Theorem III.4.1}.} Then $\cat{L}(V,-)\to1$ is a cofibrant replacement in the $G$-global model structure.
    \begin{proof}
        It is clear that $\cat{L}(V,-)$ is cofibrant. Moreover, if $\phi\colon H\to G$ is any homomorphism of compact Lie groups, and $\Uu_H$ is any complete $H$-universe, then $\cat{L}(V,-)(\Uu_H)^\phi=\cat{L}(\phi^*V,\Uu_H)^H$ is weakly contractible by \cite{schwede2018global}*{Proposition~1.1.21}, i.e. $\cat{L}(V,-)\to1$ is a $G$-global weak equivalence.
    \end{proof}
\end{lemma}

\begin{proof}[Proof of Proposition~\ref{prop:f!f*}]
    We may assume without loss of generality that $X$ is modelled by a cofibrant object of $\cat{$\bm G$-$\cat{L}$-Top}$ (again denoted by $X$), which is then in particular levelwise $G$-free and $\All$-cofibrant. We now fix a faithful $G$-representation $V$ and pick functorial cofibrant replacements in $\cat{$\bm G$-${\cat L}$-Top}$ to obtain a commutative diagram
    \[
        \begin{tikzcd}
            X'\arrow[d,"x'"']\arrow[r,"\sim"] & \cat{L}(V,-)\times\alpha^
            *X\arrow[d,"\pr_1"]\arrow[r,"\sim"]&X\arrow[d]\\
            Y'\arrow[r,"\sim"']&\cat{L}(V,-)\arrow[r,"\sim"']&1
        \end{tikzcd}
    \]
    where the horizontal maps are weak equivalences and $X',Y'$ are cofibrant. In particular, the horizontal composites still define cofibrant replacements $i\colon X'\to X$ and $j\colon Y'\to 1$ by Lemma~\ref{lemma:resolve-pt} combined with Corollary~\ref{cor:cart-prod-homotopical}. Appealing to the corollary once more, products in $\cat{$\bm G'$-$\cat{L}$-Top}$ are homotopical and hence model products in the $\infty$-categorical localization; the claim thus translates to saying that the left-hand vertical map in the commutative square
    \[
        \begin{tikzcd}
            \alpha_!X'\arrow[d, "{(\alpha_!x', \epsilon\circ \alpha_!i)}"']\arrow[r] & \alpha_!(\cat{L}(V,-)\times \alpha^*X)\arrow[d,"{(\alpha_!\pr_1,\epsilon\circ \alpha_!\pr_2)}\,"]\\
            \alpha_!Y'\times X\arrow[r] & \alpha_!\cat{L}(V,-)\times X
        \end{tikzcd}
    \]
    is a $G'$-global weak equivalence. By Proposition~\ref{prop:quotient-we-global}, the horizontal maps are $G$-global weak equivalences, while Lemma~\ref{lemma:proj-formula-underived} shows that the right-hand map is even an isomorphism. The claim therefore follows by 2-out-of-3.
\end{proof}

If $f\colon X\to Y$ is any map in $\PSh(\Glo)$, then $f_!\colon\PSh(\Glo)_{/X}\to\PSh(\Glo)_{/Y}$ lifts to an equivalence $\ul\Spc_{\Glo}(X)=\PSh(\Glo)_{/X}\iso\big(\PSh(\Glo)_{/Y}\big){}_{/f}=\ul\Spc_{\Glo}(Y)_{/f}$, and more generally it induces equivalences $\ul\Spc_{\Glo}(X)_{/Z}\iso\ul\Spc_{\Glo}(Y)_{/f_!Z}$ for all $Z\in\ul\Spc_{\Glo}(X)$. Our next goal is to prove the analogous statement for the change of group adjunction $\cat{L}\alpha_!\dashv\alpha^*$. We begin with a model categorical version:

\begin{thm}\label{thm:G-global-as-slice}
    Let $\alpha\colon G\to G'$ be a homomorphism of compact Lie groups and let $B\in\cat{$\bm G$-$\cat{L}$-Top}$ cofibrant. Then 
    \begin{equation}\label{eq:f_!-slice}
        \alpha_!\colon\cat{$\bm G$-$\cat{L}$-Top}_{/B}\to \cat{$\bm{G'}$-$\cat{L}$-Top}_{/\alpha_!B}
    \end{equation}
    is the left adjoint in a Quillen equivalence (for the usual slice model structures).
\end{thm}

The proof will require some preparations.

\begin{lemma}
    Let $G,H$ be compact Lie groups, and let $X$ be a $G$-free and Hausdorff $(G\times H)$-space. Then the quotient map $\eta\colon X\to X/G$ (i.e.\ the unit of the adjunction $(-)/G\dashv\triv_G$) is a $\mathcal G_{H,G}$-equivariant fibration.
    \begin{proof}
        Let $K\subset H$ be any closed subgroup and let $\phi\colon K\to G$ be any continuous group homomorphism. Using the homeomorphism $(\ref{eq:splitting-of-fp})$, we may identify $\eta^\phi$ with the quotient map $X^\phi\to X^\phi/C(\im\phi)$ followed by the inclusion of a disjoint summand. The latter is clearly a Serre fibration, while the former is so by Lemma~\ref{lemma:quotient-map-is-fib}.
    \end{proof}
\end{lemma}

\begin{lemma}\label{lemma:unit-fibration}
    Let $G$ be a compact Lie group, and $f\colon G\to 1$ the unique map. Then the unit $\eta\colon X\to\triv_G(X/G)$ is a $G$-global fibration for every levelwise Hausdorff and $G$-free orthogonal $G$-space $X$.
    \begin{proof}
        Applying the previous lemma levelwise shows that $\eta$ is a $G$-global level fibration. It therefore only remains to show that  the naturality square
        \[
            \begin{tikzcd}
                X(V)\arrow[r,"\eta"]\arrow[d] & X(V)/G\arrow[d]\\
                X(V\oplus W)\arrow[r,"\eta"] & X(V\oplus W)/G
            \end{tikzcd}
        \]
        is a pullback in the $\mathcal G_{H,G}$-model structure for every compact Lie group $G$ and all $H$-representations $V,W$ such that $V$ is faithful. As noticed above, the horizontal maps are $\mathcal G_{H,G}$-fibrations, so it will suffice to show that this square is a levelwise pullback. This is an instance of Lemma~\ref{lemma:quotient-pb}.
    \end{proof}
\end{lemma}

\begin{proof}[Proof of Theorem~\ref{thm:G-global-as-slice}]
    The functor $(\ref{eq:f_!-slice})$ is left adjoint to the composite
    \begin{equation}\label{eq:f!-slice-RA}
        \cat{$\bm{G'}$-$\cat{L}$-Top}_{/\alpha_!B}\xrightarrow{\;\alpha^*\;}\cat{$\bm{G}$-$\cat{L}$-Top}_{/\alpha^*\alpha_!B}\xrightarrow{\;\eta^*\;}\cat{$\bm{G}$-$\cat{L}$-Top}_{/B}
    \end{equation}
    where $\eta^*$ denotes pullback along the unit $\eta\colon B\to \alpha^*\alpha_!B$. The unit and counit for this adjunction are as follows: for $X\in\cat{$\bm G$-$\cat{L}$-Top}_{/B}$, the unit is the induced map
    \[
        \begin{tikzcd}
            X\arrow[ddr, bend right=10pt]\arrow[drr, bend left=10pt, "\eta"]\arrow[dr,dashed]&[-1em]\\[-1em]
            & B\times_{\alpha^*\alpha_!B} \alpha^*\alpha_!X\arrow[r]\arrow[d]\arrow[dr,pullback] & \alpha^*\alpha_!X\arrow[d]\\
            & B\arrow[r,"\eta"'] & \alpha^*\alpha_!B\rlap,
        \end{tikzcd}
    \]
    while the counit is given for an $Y\to \alpha_!B$ by the composite
    \[
        \alpha_!(B\times_{\alpha^*\alpha_!B}f^*Y)\xrightarrow{\;\alpha_!(\pr)\;} \alpha_!\alpha^*Y\xrightarrow{\;\epsilon\;} Y.
    \]
    Lemma~\ref{lemma:restr-right-Quillen} immediately implies that $(\ref{eq:f_!-slice})$ is left Quillen. It remains to show that the induced adjunction on homotopy categories is an equivalence.

    Write $\pi\colon G'\to 1$ for the unique homomorphism. Then $\cat{L}\pi_!\cat{L}\alpha_!\simeq\cat{L}(\pi\alpha)_!$, so it will suffice by 2-out-of-3 to prove the theorem for $\pi$ and $\pi\alpha$ in lieu of $\alpha$, i.e.\ we may assume that $G'=1$, so that $\alpha_!=(-)/G$. We will then explicitly verify that the derived unit and counit are isomorphisms in the homotopy category.

    We begin by observing that in the composite $(\ref{eq:f!-slice-RA})$, the first functor is homotopical by Lemma~\ref{lemma:restr-right-Quillen}, while the second one is so by Lemma~\ref{lemma:unit-fibration} and right properness of the model structure. Thus, for any cofibrant $p\colon X\to B$ (i.e.~$X$ is cofibrant in $\cat{$\bm G$-$\cat{L}$-Top}$), the derived unit is already modelled by the \emph{underived} unit; we claim that the latter is even an isomorphism, i.e.\ the square
    \[
        \begin{tikzcd}
            X\arrow[r, "\eta"]\arrow[d] & \triv_G(X/G)\arrow[d]\\
            B\arrow[r,"\eta"'] & \triv_G(B/G)
        \end{tikzcd}
    \]
    is a pullback in $\cat{$\bm G$-$\cat{L}$-Top}$. This can be checked levelwise; as cofibrant orthogonal $G$-spaces are levelwise $G$-free and Hausdorff, this is another instance of Lemma~\ref{lemma:quotient-pb}.

    To show that also the derived \emph{co}unit is an isomorphism, we let $q\colon Y\to B/G$ be any cofibrant-fibrant object; in particular, $q$ is a $G'$-global fibration and $Y$ is Hausdorff. We pick a cofibrant replacement $Z\to B\times_{\triv_G B/G} \triv_G Y$ in the $G$-global model structure, so that the derived counit for $Y$ is modelled by the composite
    \begin{equation}\label{eq:derived-unit-local-equiv}
        Z/G\longrightarrow (B\times_{\triv_G B/G} \triv_GY)/G\xrightarrow{\;\pr/G\;} (\triv_GY)/G\xrightarrow{\;\epsilon\;}Y.
    \end{equation}
    As $B$ is $G$-free, so is $B\times_{\triv_GB/G}\triv_GY$; moreover, the latter is Hausdorff as both $B$ and $\triv_GY$ are so. As also the cofibrant object $Z$ is $G$-free Hausdorff, Lemma~\ref{lemma:equiv-quotient} shows that the leftmost map is a $G'$-global weak equivalence, and we are reduced to showing that the composite of the remaining two maps (i.e.\ the \emph{underived} counit) is a $G'$-global weak equivalence. We will show that it is even an isomorphism. For this, note that Lemma~\ref{lemma:proj-formula-underived} provides an isomorphism
    \begin{equation}\label{eq:projform-before-res}
        \begin{aligned}
            (B \times \triv_GY)/G&\xrightarrow{\kern.4pt\,\cong\,\kern.4pt} B/G\times Y\\
            [b,y]&\longmapsto ([b], y)
        \end{aligned}
    \end{equation}
    so that the map $(B\times\triv_GY)/G\to Y$ restricts to the underived counit in question. The closed embedding $B\times_{\triv_GY/G} \triv_GY\hookrightarrow B\times\triv_GY$  induces a closed embedding on $G$-orbits by \cite{schwede2018global}*{Proposition~B.13(iii)}, and so $(\ref{eq:projform-before-res})$ restricts to a closed embedding $(B\times_{\triv_GB/G}\triv_GY)/G\hookrightarrow B/G\times Y$ with image precisely those $([b],y)$ such that $q(y)=[b]$, which is homeomorphic to $Y$ via the projection. If follows that the derived counit is an isomorphism, finishing the proof of the theorem.
\end{proof}

\begin{cor}\label{cor:slice-adjoint-equiv}
    Let $\alpha\colon G\to G'$ be any homomorphism of compact Lie groups and let $B\in\myS_\textup{$G$-gl}$ arbitrary. Then $\cat{L}\alpha_!\colon\myS_\textup{$G$-gl}\to\myS_\textup{$G'$-gl}$ induces an equivalence $(\myS_\textup{$G$-gl})_{/B}\iso(\myS_\textup{$G'$-gl})_{/\cat{L}\alpha_!B}$.
\end{cor}

Specializing to $G'=1$ and $B=1$, we see that $\myS_\textup{$G$-gl}\iso(\myS_\text{gl})_{/(\cat L(V,-)/G)}$ for any finite faithful $G$-representation $V$.

\begin{proof}
    We may assume without loss of generality that $B$ is the image of a cofibrant object in $\cat{$\bm G$-$\cat{L}$-Top}$. We pick a functorial cofibrant replacement $Q\to\id$ for the $G$-global model structure on $\cat{$\bm G$-$\cat{L}$-Top}$ (so that $\cat{L}\alpha_!\colon \myS_\textup{$G$-gl}\to\myS_\textup{$G'$-gl}$ is induced by $\alpha_!\circ Q$) and note that $Q$ lifts to a functorial cofibrant replacement $Q'$ on $\cat{$\bm G$-$\cat{L}$-Top}_{/B}$ (sending $Y\to B$ to the composite $QY\to Y\to B$). We obtain a commutative triangle
    \[
        \begin{tikzcd}
            \cat{$\bm G$-$\cat{L}$-Top}_{/B}\arrow[r, "(\alpha_!\circ Q)_{/B}"]\arrow[dr, bend right=15pt, "\alpha_!\circ Q'"'] &[1.5em] \cat{$\bm{G'}$-$\cat{L}$-Top}_{/\alpha_!QB}\arrow[d]\\
            & \cat{$\bm{G'}$-$\cat{L}$-Top}_{/\alpha_!B}
        \end{tikzcd}
    \]
    where the right vertical map is given by postcomposition with the map $\alpha_!QB\to \alpha_!B$ induced by the cofibrant replacement $QB\to B$. As the latter is a weak equivalence of cofibrant objects, so is the former, and hence the right-hand vertical map in the above diagram induces an equivalence after Dwyer--Kan localization. The same holds for the diagonal composite by the previous theorem, whence also for the top map by 2-out-of-3. However, by Lemma~\ref{lemma:slice-over-fibrant} the map induced by $(\alpha_!\circ Q)_{/B}$ on localizations can be identified with the map $(\cat{L}\alpha_!)_{/B}\colon(\myS_\textup{$G$-gl})_{/B}\to(\myS_\textup{$G'$-gl})_{/\cat{L}\alpha_!B}$ in question, finishing the proof.
\end{proof}

Finally, we record the following basechange condition:

\begin{prop}\label{prop:basechange-surj}
    Let 
    \[
        \begin{tikzcd}
            A\arrow[r,"g"]\arrow[d,"q"']\arrow[dr,pullback] & B\arrow[d,"p"]\\
            C\arrow[r,"f"'] & D
        \end{tikzcd}
    \]
    be a pullback in the 1-category of compact Lie groups\footnote{Note that the full subcategory of the category of topological groups spanned by the compact Lie groups is closed under finite limits because it is closed under finite products as well as passage to closed subgroups.} such that $p$ is \emph{surjective}. Then the Beck--Chevalley map $\cat{L}q_!g^*\to f^*\cat{L}p_!$ of functors $\myS_\textup{$B$-gl}\to\myS_\textup{$C$-gl}$ is an equivalence.
    \begin{proof}
        If $X\in\cat{$\bm B$-$\cat{L}$-Top}$ is cofibrant, then it is Hausdorff and $B$-free the above; thus, the $A$-orthogonal space $g^*X$ is Hausdorff, and it is $\ker(q)$-free as $g$ is injective when restricted to $\ker(q)$. It thus follows from Proposition~\ref{prop:quotient-we-global} that the (derived) Beck--Chevalley map $\cat{L}q_!g^*X\to f^*\cat{L}p_!X$ is already modelled by the \emph{underived} Beck--Chevalley map $q_!g^*X\to f^*p_!X$ on the pointset level. To complete the proof, we will show that this map is even an isomorphism (without any assumptions on $X$).

        For this note that since $p$ is surjective, $p_!(-)=(-)/\ker(p)$ with unit $\id\to p^*p_!$ given by the quotient map, and similarly for $q$. For this choice of adjunction data, one then directly computes that the Beck--Chevalley map is in fact the identity map of $(g^*X)/\ker(q)=X/g(\ker(q))=X/\ker(p)$. 
    \end{proof}
\end{prop}

\section{The unstable universal properties}\label{sec:unstable-univ-prop}
The continuous Borel construction (Construction~\ref{constr:continuous-Borel}) allows us to assemble the $\infty$-categories $\myS_\text{$G$-gl}$ from the previous section into a global $\infty$-category:

\begin{constr}
    The global $\infty$-category $\Ntop(\cat{$\cat{L}$-Top}^\flat)$ sends the object $\BGcat{G}$ (for a compact Lie group $G$) to $\Ntop(\cat{$\bm G$-$\cat{L}$-Top})$, and a 1-morphism of $\Glo$, corresponding to a homomorphism $f\colon G\to G'$ of compact Lie groups, to the restriction functor $f^*\colon\Ntop(\cat{$\bm {G'}$-$\cat{L}$-Top})\to\Ntop(\cat{$\bm G$-$\cat{L}$-Top})$. If we equip each $\Ntop(\cat{$\bm {G}$-$\cat{L}$-Top})$ with the $G$-global weak equivalences, then all these restriction functors are homotopical (Lemma~\ref{lemma:restr-right-Quillen}), so upon localization we obtain a global $\infty$-category $\ul\myS_\gl$ together with a global functor $\smash{\gamma\colon \Ntop(\cat{$\cat{L}$-Top}^\dual)\to\ul\myS_\gl}$ given in degree $G$ by a localization at the $G$-global weak equivalences.
\end{constr}

\begin{rk}
    If $G$ is any compact Lie group, then Lemma~\ref{lemma:fgt-simpl-structure} shows that we can equivalently define $\ul\myS_\gl(\BGcat{G})$ as the localization $\myS_\text{$G$-gl}$ of the underlying 1-category of $\cat{$\bm{G}$-$\cat{L}$-Top}$, i.e.~after forgetting the enrichment. However, while this suffices to construct the restriction of $\ul\myS_\gl$ to the $1$-category of compact Lie groups, the underlying unenriched categories of our pointset models do \emph{not} yet assemble into a simplicially enriched functor $\cat{Glo}^\op\to\cat{TOPCAT}_\Delta$, and so the enrichment on $\cat{$\bm{G}$-$\cat{L}$-Top}$ is necessary to encode some of the higher structure inherent in the global functoriality of $\ul\myS_\gl$.
\end{rk}

As our first main result we will show:

\begin{thm}\label{thm:unstable-main}
    There exists a unique equivalence $\ul\myS_\gl\simeq\ul\Spc_\Glo$.
\end{thm}

Using the universal property of $\ul\Spc_\Glo$ from Corollary~\ref{cor:univ-prop-Spc-T}, this can be restated as follows, strengthening Theorem~\ref{introthm:unstable-main} from the introduction:

\begin{thm}
    The global $\infty$-category $\ul\myS_\gl$ is globally presentable. For any globally cocomplete global $\infty$-category $\Dd$, evaluation at the terminal object defines an equivalence
    \[
        \ul\Fun^\textup{$\Glo$-cc}(\ul\myS_\gl,\Dd)\iso\Dd.
    \]
\end{thm}

Probably the most natural strategy to prove the theorem would be to first show that $\ul\myS_\gl$ is globally cocomplete, and then use the universal property of $\ul\Spc_\Glo$ to obtain the comparison map---however, we do not know how to prove global cocompleteness of $\ul\myS_\gl$ without invoking the above theorem. Namely, while we have seen that all restriction functors admit left adjoints, the Beck--Chevalley condition for these left adjoints involves the extension of $\ul\myS_\gl$ to $\PSh(\Glo)$, of which we have no concrete description in terms of the model. 

On the other hand, there is also no hope to construct the map directly on the level of simplicial sets: the functoriality of $\ul\Spc_\Glo$ is defined via cartesian straightening, while we have no control over the cartesian \emph{un}straightening of $\Ntop(\cat{$\cat{L}$-Top}^\dual)$.

We will therefore employ a different strategy. In §\ref{subsec:comparison-map}, we will explain how, given a suitable global $\infty$-category $\Cc$ with a terminal object, we can assemble the functors $g_!\colon\Cc(\BGcat{G})\simeq\Cc(\BGcat{G})_{/1}\to\Cc(1)_{/g_!1}$ induced by the left adjoint to restriction along the unique map $g\colon\BGcat{G}\to 1$ into a global functor $F\colon\Cc\to\Cc(1)_{/\Upsilon(-)}$ for some specific $\Upsilon\colon\Glo\to\Cc(1)$ with $\Upsilon(\BGcat{G})= g_!1$. 
Corollary~\ref{cor:slice-adjoint-equiv} shows that for $\Cc=\ul\myS_\gl$ this comparison functor is an equivalence, essentially reducing us to showing that $\Upsilon\colon\Glo\to\myS_\gl$ is fully faithful. Proving this full faithfulness is the hardest part of our argument, and occupies all of §\ref{subsec:Upsilon}.

In §\ref{subsec:proof-unstable-comp} we then combine these results to prove Theorem~\ref{thm:unstable-main} as well as an   $\Ff$-global generalization. Building on this, we will prove a similar universal property for an analogously defined global $\infty$-category $\ul\myS$ of \emph{equivariant spaces} in §\ref{subsec:equivariant-spaces}. Finally, we will discuss symmetric monoidal versions of these results in §\ref{subsec:sym-mon}.

\subsection{The comparison map}\label{subsec:comparison-map} Throughout this section, we let $T$ be a small $\infty$-category with a terminal object $1$. The reader is welcome to think of $T=\Glo$.
\begin{constr}
    Let $\Cc\colon T^\op\to\CAT_\infty$ be a $T$-$\infty$-category, with cartesian unstraighening $p\colon\Un^\ct(\Cc)\to T$. As $1\in T$ is terminal, the constant $T$-$\infty$-category $\mathop\const\Cc(1)$ is left Kan extended from $\{1\}$, so the identity of $\Cc(1)$ extends uniquely to a $T$-functor ${\pi^*}\colon\mathop\const\Cc(1)\to\Cc$; simply by naturality, $\pi^*(A)\colon\Cc(1)\to\Cc(A)$ is given by restriction along the unique map $a\colon A\to 1$ for any fixed $A\in T$. After unstraightening, we can view $\pi^*$ as a map $\Cc(1)\times T\to\Un^\ct(\Cc)$ of cartesian fibrations.
\end{constr}

\begin{lemma}\label{lemma:pi!-rel-la}
    Assume that $\Cc$ factors through the wide subcategory $\CAT^\textup{R}_\infty$ of right adjoint functors, i.e.\ for every $f\colon A\to B$ in $T$, the functor $f^*\colon\Cc(B)\to\Cc(A)$ admits a left adjoint $f_!$. Then $\pi^*\colon \Cc(1)\times T\to\Un^\ct(\Cc)$ admits a relative left adjoint $\pi_!$, given on the fiber over $A$ by $a_!\colon\Cc(A)\to\Cc(1)$.
\end{lemma}

Beware however that $\pi_!$ will typically \emph{not} preserve cartesian edges (so that it does not straighten to a $T$-functor).

\begin{proof}
    For every $A\in T$, the map induced by $\pi^*$ on the fibers over $A$ is $a^*\colon\Cc(1)\to\Cc(A)$, and this has a left adjoint $a_!$ by assumption. The claim now follows from \cite{HA}*{Propositions~7.3.2.5 and~7.3.2.6}.
\end{proof}

\begin{constr}\label{constr:Upsilon}
    As $1\in T$ is terminal, the inclusion $\Cc(1)\hookrightarrow\Cc(1)\times T$ of the fiber over 1 is right adjoint to the projection $\pr\colon\Cc(1)\times T\to\Cc(1)$. Thus, the fully faithful inclusion $i\colon\Cc(1)\hookrightarrow\Un^\ct(\Cc)$ has a left adjoint $\lambda$ given by the composite $\pr\circ\pi_!\colon
        \Un^\ct(\Cc)\to\Cc(1)\times T\to\Cc(1)$;
    in particular, the restriction of $\lambda$ to the fiber $\Cc(A)$ over $A$ is given by $a_!$.

    Assume now that $\Cc$ has a terminal object, and let $s\colon T\to\Un^\ct(\Cc)$ be the corresponding cartesian section, i.e.~$s$ is the unique map of cartesian fibrations such that $s(1)\in\Cc(1)$ is terminal, or equivalently every $s(A)\in\Cc(A)$ is terminal. We then define $\Upsilon\colon T\to\Cc(1)$ as the composite
    \[
        T\xrightarrow{\;s\;}\Un^\ct(\Cc)\xrightarrow{\;\lambda\;}\Cc(1).
    \] 
\end{constr}

\begin{ex}
    The cartesian unstraightening of the $T$-$\infty$-category $\Cc=\ul\Spc_T$ is the pullback $\Ar(\PSh(T))\times_{\PSh(T)}T$ along the target map and the Yoneda embedding $y\colon T\to\PSh(T)$. The composite $\const\circ y\colon T\to\Ar(\PSh(T))$ then lifts to a map $T\to\Un^\ct(\ul\Spc_T)=\Ar(\PSh(T))\times_{\PSh(T)}T$ over $T$, and one directly checks that this is fiberwise terminal, and hence agrees with the functor $s$ from above.
    
    On the other hand, $\ev_0\colon\Ar(\PSh(T))\times_{\PSh(T)}T\to\PSh(T)\simeq\PSh(T)_{/1}$ is easily seen to be left adjoint to the inclusion of the fiber over $1$, and hence agrees with $\lambda$. Thus, for $\Cc=\ul\Spc_T$ the functor $\Upsilon\colon T\to\Cc(1)=\PSh(T)_{/1}$ simply recovers the Yoneda embedding, modulo to the canonical identification $\PSh(T)\simeq\PSh(T)_{/1}$.

    What makes Construction~\ref{constr:Upsilon} useful is that it only refers to the parametrized structure of $\Cc$. In particular, if $\Cc$ is some $T$-$\infty$-category we want to show is equivalent to $\ul\Spc_T$ (like $\Cc=\ul\myS_\gl$), the construction automatically provides a candidate for the Yoneda embedding and hence of a map $T\to\Cc(1)$ that should extend to an equivalence $\ul\Spc_T(1)\simeq\PSh(T)\iso\Cc(1)$.
\end{ex}

\begin{lemma}
    Let $T$ and $\Cc$ be as in the proposition. Then $\lambda$ is a Bousfield localization such that the unit is pointwise given by cocartesian edges.
    \begin{proof}
        The first statement is clear as $\lambda$ is left adjoint to the fully faithful inclusion $\Cc(1)\hookrightarrow\Cc$. For the identification of the units, let $x\in\Cc(A)$ and let $\tilde a\colon x\to x'$ be a cocartesian lift of the unique map $a\colon A\to 1$. Then
        \[
            \begin{tikzcd}
                \hom_{\Un^\ct}(x',y)\arrow[d]\arrow[r, "-\circ\tilde a"] & \hom_{\Un^\ct}(x,y)\arrow[d]\\
                \hom_T(1,1)\arrow[r, "-\circ a"']&\hom_T(A,1)
            \end{tikzcd}
        \]
        is a pullback by definition of cocartesian edges. As both spaces in the bottom row are contractible, we see that $\hom(\tilde a,y)$ is an equivalence. Letting $y$ vary, this is precisely the defining property of a unit for our Bousfield localization.
    \end{proof}
\end{lemma}

\begin{rk}
    Let $\Cc\colon T^\op\to\CAT^\textup{R}_\infty$ admit a terminal object. The previous lemma immediately implies that $\Upsilon(A)=a_!1_A=a_!a^*1$, and an only slightly more elaborate computation shows that $\Upsilon$ sends the homotopy class of a map $f\colon A\to B$ in $T$ to the homotopy class of the composite
    \[
        a_!a^*1\iso b_!f_!f^*b^*1\xrightarrow{\;\epsilon\;} b_!b^*1
    \]
    where the unlabelled equivalence is induced by the unique homotopy $bf\simeq a$ in $T$.
\end{rk}

\begin{ex}\label{ex:Upsilon-for-Sgl}
    Specializing the above to $T=\Glo$ and $\Cc=\ul\myS_\gl$, we obtain a functor $\Upsilon\colon\Glo\to\myS_\gl$ sending $\BGcat{G}$ to $\cat{L}(V,-)/G$ where $V$ is our favourite finite-dimensional faithful $G$-representation. Combining Remark~\ref{rk:mathfrak-O-gl-vs-Stefan} with the computation of mapping spaces in $\cat{O}_\gl$ from \cite{koerschgen}*{Proposition 2.10 and Remark~2.11} already shows that $\hom_{\Glo}(\BGcat{G},\BGcat{H})$ and $\hom_{\myS_\gl}(\Upsilon(\BGcat{G}),\Upsilon(\BGcat{H}))$ are \emph{abstractly} equivalent. It will be considerably harder to show that in fact the map induced by $\Upsilon$ is an equivalence.
\end{ex}

We can now state our general criterion that we will use to produce the equivalence $\ul\Spc_\Glo\simeq\ul\myS_\gl$:

\begin{thm}\label{thm:criterion-spc-T}
    Let $T$ be a small $\infty$-category with a terminal object $1$, and let $\Cc\colon T^\op\to\CAT_\infty$ be a $T$-$\infty$-category. Assume the following:
    \begin{enumerate}
        \item $\Cc$ has a fiberwise finite limits, and every restriction functor $f^*\colon\Cc(B)\to\Cc(A)$ admits a left adjoint $f_!$.
        \item For every $a\colon A\to 1$ in $T$, the functor $\Cc(A)\simeq\Cc(A)_{/1}\to\Cc(1)_{/a_!1}$ induced by $a_!$ is an equivalence.
        \item For every $f\colon A\to B$ in $T$ and every $X\in\Cc(B)$ the diagram
        \[
            f_!f^*1\xleftarrow{\;f_!f^*x\;}f_!f^*X\xrightarrow{\;\epsilon\;}X
        \]
        exhibits $f_!f^*X$ as a product, where $x\colon X\to 1$ denotes the unique map and $\epsilon$ denotes the counit of $f_!\dashv f^*$ as usual.
        \item The functor $\Upsilon\colon T\to\Cc(1)$ extends to an equivalence $\PSh(T)\to\Cc(1)$.
    \end{enumerate}
    Then $\Cc$ is equivalent to $\ul\Spc_T$ (in a necessarily unique way).
\end{thm}

For the proof we will need the following `basechange condition' for $\lambda$:

\begin{lemma}\label{lemma:lambda-basechange}
    Let $T$ and $\Cc$ as above, and let $A\in T$ arbitrary. Then:
    \begin{enumerate}
        \item The unique map $\pi_A^*\colon T_{/A}\times\Cc(A)\to T_{/A}\times_T\Un^\ct(\Cc)$ of cartesian fibrations over $T_{/A}$ extending the inclusion of $\{\id_A\}\times\Cc(A)$ admits a relative left adjoint $\pi_{A!}$.
        \item Writing $\lambda_A$ for the composite $\pr\circ\pi_{A!}$, the diagram
        \[
            \begin{tikzcd}
               T_{/A}\times_T\Un^\ct(\Cc)\arrow[d,"\pr"']\arrow[r,"\lambda_A"]&\Cc(A)\arrow[d,"a_!"]\\
               \Un^\ct(\Cc)\arrow[r,"\lambda"'] & \Cc(1)
            \end{tikzcd}
        \]
        commutes up to homotopy.
    \end{enumerate}
    \begin{proof}
        For the first statement it suffices to apply Lemma~\ref{lemma:pi!-rel-la} to the restriction of $\Cc$ along the forgetful functor $T_{/A}^\op\to T^\op$. 
        
        For the second statement, note that \cite{HA}*{Proposition 7.3.2.6} shows that  $T_{/A}\times_T\pi_!\colon T_{/A}\times_T\Un^\ct(\Cc)\to T_{/A}\times\Cc(1)$ is left adjoint to $T_{/A}\times_T\pi^*$. On the other hand, the composite
        \[
            T_{/A}\times\Cc(1)\xrightarrow{\;\id\times a^*\;}T_{/A}\times\Cc(A)\xrightarrow{\;\pi_A\;}T_{/A}\times_T\Un^\ct(\Cc)
        \]
        agrees with $T_{/A}\times_T\pi^*$ as both are maps of cartesian fibrations over $T_{/A}$ and their straightenings are given at the terminal object $\id_A$ by $a^*\colon\Cc(1)\to\Cc(A)$. Passing to left adjoints, we see that $T_{/A}\times_T\pi_{!}$ is equivalent to $(a_!\times\id)\circ\pi_{A!}$. The claim now follows by postcomposing with the projection $T_{/A}\times\Cc(1)\to\Cc(1)$.
    \end{proof}
\end{lemma}

\begin{proof}[Proof of Theorem~\ref{thm:criterion-spc-T}]
    Consider the cartesian fibration $\ev_1\colon\Ar(\Cc(1))\to\Cc(1)$ for the functor $X\mapsto \Cc(1)_{/X}$. Pulling back along $\Upsilon\colon T\to\Cc(1)$, we obtain the cartesian fibration $q\colon\Ar(\Cc(1))\times_{\Cc(1)}T\to T$. By the final assumption, this is equivalent over $T$ to the cartesian unstraightening $\Ar(\PSh(T))\times_{\PSh(T)}T$ of $\ul\Spc_T$. To complete the proof it will therefore suffice to show that $q$ is also equivalent to the cartesian unstraightening $p\colon\Un^\ct(\Cc)\to T$ of $\Cc$.

    We begin by making the unstraightening of $\Ar(\Cc)\colon A\mapsto\Ar(\Cc(A))$ explicit. Cartesian unstraightening commutes with the $\Cat$-tensoring by direct inspection, hence also with the $\Cat$-cotensoring (as unstraightening is an equivalence). Thus, the cartesian unstraightening in question is given by $\Ar\big(\Un^\ct(\Cc)\big)\times_{\Ar(T)}T$. As $[1]\to[0]$ is a localization (which is just a fancy way to say that the 1-simplex is weakly contractible), the diagonal functor $T\to \Ar(T)$ is fully faithful with essential image the equivalences of $T$. We conclude that the unstraightening can also be described as the full subcategory $\Ar^\text{fw}(\Un^\ct(\Cc))\subset\Ar(\Un^\ct(\Cc))$ spanned by the \emph{fiberwise} arrows, i.e.\ those sent to equivalences in $T$ via $p$.

    Recall the cartesian section $s\colon T\to\Un^\ct(\Cc)$. We form the pullback
    \[
        \begin{tikzcd}
            \mathcal R\arrow[r]\arrow[d, "r"']\arrow[dr,pullback] & \Ar^\text{fw}(\Un^\ct(\Cc))\arrow[d,"\ev_1"]\\
            T\arrow[r,"s"]\arrow[rr,"="', bend right=10pt,yshift=-3pt] & \Un^\ct(\Cc)\arrow[r,"p"] & T\rlap.
        \end{tikzcd}
    \]
    The functor $s$ is fully faithful, with essential image spanned by those $(A,X\in\Cc(A))$ where $X$ is terminal. Thus, the top horizontal arrow is fully faithful, with essential image given by those $(A, X\to Y)$ with $Y\in\Cc(A)$ terminal. Note that these full subcategories assemble into a subfunctor of $\Ar(\Cc)$, so that $r\colon\mathcal R\to T$ is again a cartesian fibration, with cartesian edges given by those squares
    \[
        \begin{tikzcd}
            \cdot\arrow[r,"\alpha"]\arrow[d] & \cdot\arrow[d]\\
            \cdot\arrow[r,"\beta"'] & \cdot
        \end{tikzcd}
    \]
    such that both $\alpha$ and $\beta$ are cartesian. On the other hand, the source map $\ev_0\colon\Ar^\text{fw}(\Un^\ct(\Cc))\to\Un^\ct(\Cc)$ is a map of cartesian fibrations (corresponding to the source map $\Ar(\Cc)\to\Cc$), and it restricts to a (fiberwise) equivalence on the image of $\mathcal R$---namely, it is given fiberwise by $\fgt\colon\Cc(A)_{/1}\iso\Cc(A)$. We are thus reduced to producing an equivalence between $r$ and $q$.

    For this we use that $\pi_!\colon\Un^\ct(\Cc)\to\Cc(1)\times T$ is a map over $T$, so it induces a map $\Ar^\text{fw}(\pi_!)\colon\Ar^\text{fw}(\Un^\ct(\Cc))\to\Ar^\text{fw}(\Cc(1)\times T)$ fitting into
    \[
        \begin{tikzcd}[row sep={3.3em,between origins}, column sep={4.2em,between origins}]
            &&&&& \Ar(\Cc(1))\arrow[ddd,"\ev_1"{description}]\\
            &&&& \Ar^\text{fw}(\Cc(1)\times T)\arrow[ur,"\pr"]\\
            && & \Ar^\text{fw}(\Un^\ct(\Cc))\arrow[ur,"\Ar^\text{fw}(\pi_!)"{description,xshift=-.29in,yshift=.4em}] &&\\
            && T\arrow[rrr,"\Upsilon"{description}] &&&  \Cc(1) \\
            & &&& \Cc(1)\times T\arrow[ur,"\pr"{description}]\arrow[from=uuu,crossing over,"\ev_1"{description}]\\
            T\arrow[uurr,"="{description}]\arrow[rrr,"s"'] &&& \Un^\ct(\Cc)\arrow[ur,"\pi_!"{description}] \arrow[from=uuu,crossing over,"\ev_1"{description}]
        \end{tikzcd}
    \]
    where the bottom face commutes by construction of $\Upsilon$. Passing to pullbacks we then obtain a map $\phi\colon r\to q$ over $T$, which we claim is an equivalence. Unravelling the definitions, $\phi$ is fiberwise given by $a_!\colon\Cc(A)_{/1}\to\Cc(1)_{/a_!1}$, which is an equivalence by our second assumption. It therefore only remains to show that $\phi$ preserves cartesian edges. Plugging in the characterization of cartesian edges of $\mathcal R$ from above, this amounts to saying that given any commutative square
    \begin{equation}\label{diag:should-become-pb}
        \begin{tikzcd}
            X\arrow[r,"\ct"]\arrow[d,"\text{fw}"'] & Y\arrow[d,"\text{fw}"]\\
            s(A)\arrow[r,"\ct"'] & s(B)
        \end{tikzcd}
    \end{equation}
    in $\Un^\ct(\Cc)$ such that the horizontal arrows are cartesian, the vertical maps are fiberwise, and the bottom row is in the essential image of $s$, then the image of $(\ref{diag:should-become-pb})$ under $\lambda=\pr\circ\pi_!$ is a pullback in $\Cc(1)$.
    
    For this we note that the image of $(\ref{diag:should-become-pb})$ under $p\colon\Un^\ct(\Cc)\to T$ lifts uniquely to $T_{/B}$ (as $[1]\times[1]$ has a terminal object), so that we may view $(\ref{diag:should-become-pb})$ as a square in $T_{/B}\times_T\Un^\ct(\Cc)$. By Lemma~\ref{lemma:lambda-basechange} it will then suffice to show that the composite
    \[
        T_{/B}\times_T\Un^\ct(\Cc)\xrightarrow{\;\lambda_B\;}\Cc(B)\xrightarrow{\;b_!\;}\Cc(1)_{/b_!1}\xrightarrow{\;\fgt\;}\Cc(1)
    \]
    sends this lift to a pullback. However, the last of these functors preserves pullbacks, while the functor in the middle is even an equivalence by our assumption, so that it will suffice to show that the image under $\lambda_B$ is a pullback. As $\lambda_Bs(B)$ is terminal, the claim simplifies to $\lambda_B(s(A))\gets\lambda_B(X)\to\lambda_B(Y)$ being a product diagram in $\Cc(B)$. Using the characterization of $\lambda_B$ as the recipient of a transformation consisting of cocartesian edges, we identify this as precisely the diagram from our third assumption, finishing the proof.
\end{proof}

\subsection{Understanding \for{toc}{$\Upsilon$}\except{toc}{\texorpdfstring{$\bm\Upsilon$}{Υ}}}\label{subsec:Upsilon} As mentioned above, we want to obtain the equivalence $\ul\myS_\gl\iso\ul\Spc_\gl$ as an application of Theorem~\ref{thm:criterion-spc-T}. Fortunately, we have already verified the first three assumptions in §\ref{subsec:model-cat-prereq}, while for the final condition it will suffice by the global Elmendorf Theorem (Proposition~\ref{prop:global-Elmendorf}) to show that $\Upsilon$ defines an equivalence onto the full subcategory $\mathfrak O_\gl\subset\myS_\gl$.

Proving full faithfulness of $\Upsilon$ will be by far the hardest part of our argument. The key problem is once again that we do not understand the cartesian unstraightening of $\ul\myS_\gl$, and the definition of $\Upsilon$ crucially uses this unstraightening. We therefore have to find a description of the effect of $\Upsilon$ on morphism spaces that avoids talking about unstraightenings. For this we will use a trick: fixing a compact Lie group $H$ and varying $G$ we obtain a natural transformation
\begin{equation}\label{eq:Upsilon-on-hom}
    \Upsilon\colon\hom_{\Glo}(\BGcat{G},\BGcat{H})\to\hom_{\myS_\gl}(\Upsilon(\BGcat{G}),\Upsilon(\BGcat{H}))
\end{equation}
of functors $\Glo^\op\to\Spc$, and we want to prove this is an equivalence. By the Yoneda lemma, $(\ref{eq:Upsilon-on-hom})$ is uniquely characterized by it sending the identity of $\BGcat{H}$ to the identity of $\Upsilon(\BGcat{H})$. If we can therefore find \emph{some} natural equivalence with this property (for each fixed $H$), this will show that $\Upsilon$ is fully faithful.

We will achieve this in two steps: we will first use $\infty$-categorical arguments to massage the right-hand side so that it does not involve $\Upsilon$ anymore, and then use our model to explicitly write down a natural equivalence with the required properties. 

\medskip
\subsubsection{Understanding the target} The first part will in fact work in much greater generality again, so we return to our general setting of a small $\infty$-category $T$ with a terminal object and a $T$-$\infty$-category $\Cc\colon T^\op\to\CAT^\textup{R}_\infty$ with a terminal object. The basic observation is that we have for fixed $A,B\in T$ an equivalence
\[\hskip-50.3pt\hfuzz=50.3pt
    \hom_{\Cc(1)}(\Upsilon(A),\Upsilon(B))=
    \hom_{\Cc(1)}(\lambda s(A),\Upsilon(B))\iso
    \hom_{\Un^\ct(\Cc)}(s(A),\Upsilon(B))\iso
    \hom_{\Cc(A)}(s(A),a^*\Upsilon(B))
\]
by adjunction and the universal property of cartesian edges. We moreover know that $a^*\Upsilon(B)\simeq a^*b_!1$, so we altogether obtain an equivalence \[\hom_{\Cc(1)}(\Upsilon(A),\Upsilon(B))\iso\hom_{\Cc(A)}(s(A),a^*b_!1).\] The hard part will now be to first explain how to make the right-hand side into a functor $T^\op\to\Spc$ for varying $A$ and then to enhance the above pointwise comparison to a natural equivalence.

\begin{constr}\label{constr:slice-cat}
    Let $Y\in\Cc(1)$ be any global section, and recall that the \emph{slice} $T$-$\infty$-category $\Cc_{/Y}$ is defined by the pullback
    \[
        \begin{tikzcd}
            \Cc_{/Y}\arrow[d]\arrow[r]\arrow[dr,pullback] & \Ar(\Cc)\arrow[d,"\ev_1"]\\
            1\arrow[r, "Y"'] & \Cc\rlap.
        \end{tikzcd}
    \]
    We now write $\Cc_{/Y}^\circ$ for the full subcategory spanned in degree $A$ by those objects $X\to a^*Y$ for which $X$ is terminal in $\Cc(A)$; in other words, this is the pullback
    \[
        \begin{tikzcd}
            \Cc_{/Y}^\circ\arrow[d]\arrow[r]\arrow[dr,pullback] & \Cc_{/Y}\arrow[d,"\fgt\kern.4pt=\kern.4pt\ev_0"]\\
            1\arrow[r, "1"'] & \Cc\rlap.
        \end{tikzcd}
    \]
\end{constr}

\begin{prop}\label{prop:the-horrible-proof-we-dont-talk-about}
    The $T$-$\infty$-category $\Cc_{/Y}^\circ$ is a $T$-$\infty$-groupoid, and there exists an equivalence $\smash{\hom(\Upsilon(-),Y)\iso\Cc_{/Y}^\circ}$ sending a map $a_!a^*1=\Upsilon(A)\to Y$ to the object in the slice given by its adjunct $a^*1\to a^*Y$.
    \begin{proof} 
        Throughout, we will abbreviate $\Uu\coloneq\Un^\ct(\Cc)$ for ease of notation. The adjunction equivalence $\hom_{\Cc(1)}(\lambda s(-),Y)\iso\hom_{\Uu}(s(-),Y)$ is natural, so it will suffice to construct an equivalence $\hom_{\Uu}(s(-),Y)\iso\Cc_{/Y}^\circ$ given on objects by factoring a map $s(A)\to Y$ into a fiberwise map followed by a cartesian edge. We will do so by working on the level of cartesian unstraightenings.

        \smallskip
        \textit{Step 1.} Let us describe the cartesian unstraightening of the target first. In the proof of Theorem~\ref{thm:criterion-spc-T} we explained that the cartesian unstraightening of $\Ar(\Cc)$ is given by $\Ar^\text{fw}(\Uu)$, with evaluation maps corresponding to the obvious maps to $\Uu$. As unstraightening preserves limits, we conclude that the cartesian unstraighening of $\Cc_{/Y}^\circ$ is given by the iterated pullback
        \[
            \begin{tikzcd}[row sep=0ex]
                T\arrow[dr, "s"'] && \arrow[dl,"\ev_0"{description}]\Ar^\text{fw}(\Uu)\arrow[dr,"\ev_1"{description}] && \arrow[dl,"s_Y"]T\rlap,\\
                & \Uu && \Uu
            \end{tikzcd}
        \]
        where $s_Y$ denotes the unique map of cartesian fibration with $s_Y(1)=Y\in\Cc(1)\subset\Uu$.

        \smallskip
        \textit{Step 2.} On the other hand, we can also describe the cartesian unstraightening of the source quite easily: taking the usual description of the unstraightening of the contravariant hom-functor, it is given by the iterated pullback
        \[
            \begin{tikzcd}[row sep=0ex]
                T\arrow[dr, "s"'] && \arrow[dl,"\ev_0"{description}]\Ar(\Uu)\arrow[dr,"\ev_1"{description}] && \arrow[dl,"Y"]1\rlap.\\
                & \Uu && \Uu
            \end{tikzcd}
        \]

        \smallskip
        \textit{Step 3.} We will now massage the description from Step 1. For this recall that $\Uu$ admits a factorization system with epis the fiberwise maps and monos the cartesian morphisms \cite{HTT}*{Example 5.2.8.15${}^\op$}. In other words, if we write $\text{Fact}(\Uu)\subset\Fun([2],\Uu)$ for the full subcategory spanned by those $X\to Y\to Z$ where $X\to Y$ is fiberwise and $Y\to Z$ is cartesian, then $d_1^*\colon\Fun([2],\Uu)\to\Fun([1],\Uu)$ restricts to an equivalence $\text{Fact}(\Uu)\iso\Ar(\Uu)$, see \cite{HTT}*{Proposition 5.2.8.17}. Decomposing $[2]\simeq[1]\amalg_{[0]}[1]$, we then obtain an equivalence $\Ar^\text{fw}(\Uu)\times_\Uu\Ar^\ct(\Uu)\iso\Ar(\Uu)$ over $\Uu\times\Uu$, given on objects by composition, and with inverse given on objects by factorization.

        Write $\Uu^\ct\subset\Uu$ for the wide subcategory spanned by the cartesian edges. Then the restricted map $\Uu^\ct\to T$ is a right fibration, so it induces an equivalence
        \[
           \Ar(\Uu^\ct)\times_{\Uu^\ct}\{Y\}\iso\Ar(T)\times_T\{1\}\iso T,
        \]
        where both pullbacks are via the target map. On the other hand, left cancellability of cartesian edges shows that the source agrees with $\Ar^\ct(\Uu)\times_{\Uu}\{Y\}$. We then claim that the resulting equivalence fits into a commutative diagram
        \[
            \begin{tikzcd}
                \arrow[d,"\pr"']\Ar^\ct(\Uu)\times_\Uu \{Y\}\arrow[r,"\sim"] & T\arrow[d,"s_Y"]\\
                \Ar(\Uu^\ct)\arrow[r,"\ev_0"']\arrow[d,"p\circ\ev_0"'] & \Uu^\ct\arrow[d,"p"]\\
                T\arrow[r,equals] & T\rlap.
            \end{tikzcd}
        \]
        Commutativity of the outer rectangle follows easily from the definitions, and commutativity of the bottom square is tautological; in particular, we can view both composites $\Ar^\ct(\Uu)\times_\Uu\{Y\}\to\Uu^\ct$ as maps over $T$ into the right fibration $\Uu^\ct\to T$. As inclusions of terminal objects are right anodyne, it will then suffice for commutativity of the top square to show that both paths agree on the terminal object $\id_Y$, where one immediately checks that both are given by $Y\in\Cc(1)\subset\Uu^\ct$.

        With this established, we can rewrite the pullback $\Ar^\text{fw}(\Uu)\times_\Uu T$ from Step 1 as
        \[
            \begin{tikzcd}
                \Ar^\text{fw}(\Uu)\arrow[d,"\ev_1"'] &\arrow[l,"d_2^*"'] \text{Fact}(\Uu)\arrow[d]\arrow[dl,phantom,"\llcorner"{very near start}] &\arrow[dl,phantom,"\llcorner"{very near start}] \text{Fact}(\Uu)\times_\Uu\{Y\}\arrow[d]\arrow[l]\\
                \Uu & \arrow[l, "\ev_0"] \Ar^\ct(\Uu) & \arrow[l] \Ar^\ct(\Uu)\times_{\Uu}\{Y\}
            \end{tikzcd}
        \]
        Using the equivalence $\text{Fact}(\Uu)\iso\Ar(\Uu)$, we then see that the unstraightening of $\Cc_{/Y}^\circ$ is equivalent to the limit of
        \[
            \begin{tikzcd}[row sep=0ex]
                T\arrow[dr, "s"'] && \arrow[dl,"\ev_0"{description}]\Ar(\Uu)\arrow[dr,"\ev_1"{description}] && \arrow[dl,"Y"]1\rlap,\\
                & \Uu && \Uu,
            \end{tikzcd}
        \]
        which is precisely the description of the unstraightening of $\hom_\Uu(s(-),Y)$ we obtained in Step 2.

        \smallskip
        \textit{Step 4.} The description of our equivalence on objects follows at once by chasing through definitions, using that $\Fun([2],\Uu)\iso\Ar^\text{fw}(\Uu)\times_\Uu\Ar^\ct(\Uu)$ was given on objects by factorization.
    \end{proof}
\end{prop}

\subsubsection{A model categorical description} 
We will now explain how to describe $(\ul\myS_\gl)_{/Y}^\circ$ in more model categorical terms. Throughout, we let $\cat{T}$ be a small locally fibrant simplicial category with a \emph{strictly terminal} object $1$, i.e.\ such that the Kan complex $\maps(A,1)$ is terminal for every $A\in\cat{T}$ (as opposed to merely being contractible); the reader is welcome to set $\cat{T}=\cat{Glo}$. We write $T\coloneqq N_\textup{top}(\cat{T})$, which will then in particular be an $\infty$-category with terminal object again.

\begin{constr}
    Let $\cat{C}\colon\cat{T}^\op\to\cat{TOPCAT}_\Delta$ be a simplicially enriched functor, and let $Y\in\cat{C}(1)$ be arbitrary. As $1\in\cat{T}$ is \emph{strictly} terminal, $Y$ gives rise to a simplicially enriched transformation $\const(1)\Rightarrow\cat{C}$ of enriched functors $\cat{T}^\op\to\cat{TOPCAT}_\Delta$, and we can therefore define $\cat{C}_{/Y}\colon\cat{T}^\op\to\cat{TOPCAT}_\Delta$ as the (strict) pullback
    \[
        \begin{tikzcd}
            \cat{C}_{/Y}\arrow[d]\arrow[r]\arrow[dr,pullback] & \cat{Ar}(\cat{C})\arrow[d,"\ev_1"]\\
            \const(1)\arrow[r,"Y"'] & \cat{C}\rlap.
        \end{tikzcd}
    \]
\end{constr}

\begin{constr}
    For topological categories $\cat{I},\cat{C}$ we recall the natural map
    \[
        N_\text{top}\big(\cat{Fun}_\text{cont}(\cat{I},\cat{C})\big)\to\Fun\big(N_\text{top}(\cat{I}),N_\text{top}(\cat{C})\big)
    \]
    from Construction~\ref{constr:continuous-Borel}. Specializing to $\cat{I}=[1]$, we obtain a natural comparison map $N_\text{top}(\cat{Ar}(\cat{C}))\to\Ar(N_\text{top}(\cat{C}))$ compatible with the evaluation functors. 
\end{constr}

\begin{constr}
    Let $\cat{C}\colon\cat{T}^\op\to\cat{TOPCAT}_\Delta$ be a simplicially enriched functor. Assume we have equipped $\cat{C}(A)$ for every $A\in T$ with a topological model structure, and assume moreover that all restriction functors $\cat{C}(B)\to\cat{C}(A)$ are homotopical, so that we may view $\Ntop(\cat{C})$ as a functor into marked quasi-categories. Postcomposing with an enriched functorial fibrant replacement in marked simplicial sets \cite{cat-htpy}*{Theorem~13.5.2}, we obtain an enriched transformation $\Ntop(\cat{C})\to\Cc$ of $\cat{Kan}$-enriched functors $\cat{T}^\op\to\cat{QCAT}$ given pointwise by a localization at the given weak equivalences; after passing to homotopy coherent nerves, we may equivalently view this as a natural transformation of ($\infty$-categorical) functors $T^\op\to N_\Delta(\cat{QCAT})=\CAT_\infty$.

    Combining this with the previous construction, we obtain $\Ntop(\cat{Ar}(\cat{C}))\to \Ar(\Cc)$, and hence after passing to strict fibers\footnote{As the Joyal model structure is cartesian, $\ev_1\colon\Ar(\Cc)\to\Cc$ is a Joyal fibration, so that its pointset level fiber in simplicial sets models the homotopy fiber; we do not claim that the same holds for $\Ntop(\cat{Ar}(\cat{C}))\to\Ntop(\cat{C})$.} $\Ntop(\cat{C}_{/Y})\to\Cc_{/Y}$ for any $Y\in\cat{C}(1)$.
\end{constr}

\begin{prop}\label{prop:model-categorical-model}
    In the situation of the previous construction, assume in addition that all model structures are right proper and admit functorial factorizations. Then the map $\Ntop(\cat{C}_{/Y})\to\Cc_{/Y}$ is degreewise a localization at the given weak equivalences. Moreover, if we write $\cat{C}_{/Y}^\circ\subset\cat{C}_{/Y}$ for the full subcategory given in degree $A\in\cat{T}$ by those $X\to Y$ such that $X$ defines a terminal object in the localization of $\cat{C}(A)$, then also the induced map $\Ntop(\cat{C}_{/Y}^\circ)\to\Cc_{/Y}^\circ$ is a localization (at the same weak equivalences).
    \begin{proof}
        We will only prove the second statement, the proof of the first one being analogous. By Lemma~\ref{lemma:fgt-simpl-structure}, the canonical map $N\big((\cat{C}_{/Y}^\circ)(A)_0\big)\to N_\top\big((\cat{C}_{/Y}^\circ)(A)\big)$ from the nerve of the underlying unenriched category induces an equivalence on localizations, so it will suffice to show that the canonical map $\smash{N(\cat{C}_{/Y}^\circ(A)_0)\to\Cc_{/Y}^\circ(A)}$ is a localization. For $N(\cat{C}_{/Y}(A)_0)\to\Cc_{/Y}(A)$ this is an instance of Lemma~\ref{lemma:slice-over-fibrant}; \cite{g-global}*{Proposition~A.1.15} then shows that the induced map from the localization of $(\cat{C}_{/Y}^\circ(A)_0)$ to $\Cc_{/Y}(A)$ is fully faithful, and so it only remains to show that every object of $\Cc_{/Y}^\circ(A)$ is contained in its essential image. Indeed, any object is equivalent to the image of some $(X\to Y)\in\cat{C}_{/Y}(A)$, and since $\cat{C}_{/Y}(A)\to\Cc_{/Y}(A)$ lives over the localization $\cat{C}(A)\to\Cc(A)$, $X$ is weakly contractible.
    \end{proof}
\end{prop}

\subsubsection{Full faithfulness of $\Upsilon$} We now specialize to $\cat{T}=\cat{Glo}$. Proposition~\ref{prop:the-horrible-proof-we-dont-talk-about} provides a natural equivalence $\hom_{\myS_\gl}(\Upsilon(-),\Upsilon(\mathbb BH))\iso({\ul\myS_\gl})_{/\Upsilon(\mathbb BH)}^\circ$ sending the identity of $\Upsilon(\mathbb BH)=h_!1$ to the unit $1\to h^*h_!1$ (where $h\colon\BGcat{H}\to 1$ is the unique map); to see that $\Upsilon$ induces equivalences $\hom_{\Glo}(\BGcat{G},\BGcat{H})\iso\hom_{\Glo}(\Upsilon(\BGcat{G}),\Upsilon(\BGcat{H}))$ for all $G\in\Glo$, it will therefore suffice to give \emph{some} equivalence $\hom_{\Glo}(-,\BGcat{H})\iso ({\ul\myS_\gl})_{/h_!1}^\circ$ sending $\id_{\BGcat{H}}$ to the unit. We will do so by using the concrete model of the right-hand side provided by Proposition~\ref{prop:model-categorical-model}.

\begin{constr}
    Fix a faithful $H$-representation $V$, so that we may take $\Upsilon(\BGcat{H})\coloneqq\cat{L}(V,-)/H$, see Example~\ref{ex:Upsilon-for-Sgl}. We will now define an $\cat{SSet}$-enriched natural transformation $\upsilon\colon\maps(-,\BGcat{H})\to\Ntop\big(\FUN_\text{cont}(-,\cat{$\cat{L}$-Top})_{/\cat{L}(V,-)/H}^\circ\big)$. For this observe that the left-hand side was defined as the topological nerve of $\cat{Fun}_\text{cont}(-,\BGcat{H})$, and that the enriched functoriality of both source and target of $\upsilon$ were induced from the $\cat{topCat}$-enriched functoriality. It will therefore suffice to construct a $\cat{topCat}$-enriched natural transformation \[\upsilon_\top\colon\FUN_\text{cont}(-,\BGcat{H})\to\cat{Fun}_\text{cont}(-,\cat{$\cat{L}$-Top})_{/\cat{L}(V,-)/H}^\circ.\] For this we take the unique one sending $\id_{\BGcat{H}}$ to the quotient map $\cat{L}(V,-)\to\triv_H(\cat{L}(V,-)/H)$. Explicitly, $\upsilon_\top$ is given in degree $G$ as the following topological functor:
    \begin{enumerate}
        \item A Lie group homomorphism $\alpha\colon G\to H$ is sent to the quotient map $\alpha^*\cat{L}(V,-)\to\triv_G(\cat{L}(V,-)/H)$.
        \item An element $h\in H$ defining a natural transformation $\alpha\to\beta$ is sent to the map $\alpha^*\cat{L}(V,-)\to \beta^*\cat{L}(V,-)$ (over $\cat{L}(V,-)/H$) given by acting from the right with $h^{-1}$.
    \end{enumerate}
    The topologically enriched natural transformation $\upsilon$ is then simply given in each degree as the topological nerve of $\upsilon_\top$.
\end{constr}

\begin{prop}\label{prop:lc-upsilon-equiv}
    The composite
    \[
        \smash{\maps(\BGcat{G},\BGcat{H})\hskip0pt minus 1pt\xrightarrow{\,\upsilon\,}\hskip0pt minus 1pt\Ntop\hskip0pt minus .5pt\big(\cat{$\bm G$-$\cat{L}$-Top}_{/\triv_G(\cat{L}(V,-)/H)}^\circ\big)\hskip0pt minus 1pt\xrightarrow{\,\textup{loc}\,}\hskip0pt minus 1pt(\myS_\textup{$G$-gl})_{/\triv_G(\cat{L}(V,-)/H)}^\circ}
    \]
    is an equivalence of $\infty$-groupoids.
\end{prop}
\begin{proof}
    Throughout, we will abbreviate $\cat{L}_V\coloneqq\cat{L}(V,-)$, and we will suppress the functor $\triv_G$ from the notation.

    We will first prove essential surjectivity. As the right-hand map is a localization (Proposition~\ref{prop:model-categorical-model}) and its target is an $\infty$-groupoid, it will suffice to show that any object of $\cat{$\bm G$-$\cat{L}$-Top}_{/(\cat{L}_V/H)}^\circ$ can be connected by a zig-zag of maps to an object in the image of $\upsilon$.        
    Let $f\colon X\to\cat{L}_V/H$ be any map in $\cat{$\bm G$-$\cat{L}$-Top}$ such that $X$ is $G$-globally weakly contractible. After cofibrant replacement, we may assume that $X$ is closed (Definition~\ref{def:closed-orth-spc}). If $\Uu_G$ is any complete $G$-universe, then Lemma~\ref{lemma:we-between-closed} shows that $X(\Uu_G)^{\id_G}$ is contractible, in particular non-empty. Thus, we can find a $G$-representation $W$ together with a point $x\in X(W)^G$. By the Yoneda lemma we obtain a $G$-equivariant map $\cat{L}(W,-)\to X$ sending $\id_W$ to $x$. Note that $\cat{L}(W,-)$ is weakly contractible as an $\O(W)$-orthogonal space by Lemma~\ref{lemma:resolve-pt}, hence also as a $G$-orthogonal space by Lemma~\ref{lemma:restr-right-Quillen}; thus, we have reduced to the case $X=\cat{L}(W,-)$.

    The map $f\colon\cat{L}(W,-)\to \cat{L}_V/H$ corresponds to a point $[u]\in(\cat{L}(V,W)/H)^G$. Splitting this fixed point space via \cite{schwede2018global}*{Proposition~B.18} as before, we see that there exists a continuous homomorphism $\alpha\colon G\to H$ and a representative $u$ that is $G$-equivariant as a map $\alpha^*V\to W$. Then $-\circ u\colon \cat{L}(W,-)\to \alpha^*\cat{L}_V$ is $G$-equivariant, and the diagram
    \[
        \begin{tikzcd}[column sep=0pt]
            \cat{L}(W,-)\arrow[dr,"f"',bend right=15pt] \arrow[rr,"-\circ u"] &&  \alpha^*\cat{L}_V\arrow[dl,bend left=15pt,"\upsilon(\alpha)"]\\
            & \cat{L}_V/H
        \end{tikzcd}
    \]
    commutes; this completes the proof of essential surjectivity.

    For {full faithfulness}, we recall that $\upsilon$ was defined as the topological nerve of the enriched functor $\upsilon_\top$, so it will suffice to show that for all Lie group homomorphisms $\alpha,\beta\colon G\rightrightarrows H$ the composite
    \begin{multline*}
        \hskip-10pt\maps_{\cat{Fun}_\text{cont}(\BGcat{G},\BGcat{H})}(\alpha,\beta)\xrightarrow{\;\upsilon_\top\;}
        \maps_{\cat{$\bm G$-$\cat{L}$-Top}^\circ_{/\cat{L}_V/H}}(\alpha^*\cat{L}_V,\beta^*\cat{L}_V)\\{}\simeq\hom_{\Ntop(\cat{$\bm G$-$\cat{L}$-Top}^\circ_{/\cat{L}_V/H})}(\alpha^*\cat{L}_V, \beta^*\cat{L}_V)       
        \xrightarrow{\;\text{loc}\;}\hom_{(\myS_\text{$G$-gl})_{/\cat{L}_V/H}}(\alpha^*\cat{L}_V, \beta^*\cat{L}_V)\hskip-10pt
    \end{multline*}
    is a weak homotopy equivalence. By Lemma~\ref{lemma:unit-fibration}, $\cat{L}_V\to\cat{L}_V/H$ is an $H$-global fibration, so $\upsilon(\beta)$ is a $G$-global fibration by Lemma~\ref{lemma:restr-right-Quillen}. Moreover, as in the proof of essential surjectivity we can find a $G$-representation $W$ together with a point $u\in\cat{L}(V,W)^\alpha$, and enlarging $W$ if necessary, we may assume that $W$ is faithful; the Yoneda lemma then provides a unique map $f\colon\cat{L}(W,-)\to\alpha^*\cat{L}_V$ sending $\id_W$ to $u$, which is a cofibrant replacement by Lemma~\ref{lemma:resolve-pt}. By the description of $\text{loc}$ on mapping spaces from Lemma~\ref{lemma:localization-on-homs}, it will then suffice to show that
    \[
        \maps(\alpha,\beta)\xrightarrow{\;\upsilon_\top\;}
        \maps(\alpha^*\cat{L}_V,\beta^*\cat{L}_V)
        \xrightarrow{\;{-}\circ f\;}\maps(\cat{L}(W,-),\beta^*\cat{L}_V)
    \]
    is a weak homotopy equivalence. Plugging in the definitions, this composite sends an $h\in H$ with $\beta=c_h\circ\alpha$ to a map $\cat{L}(W,-)\to\beta^*\cat{L}_V$ with $\id_W\mapsto u.h^{-1}$. On the other hand, $\cat{L}(W,-)\to\cat{L}_V/H$ corepresents (via evaluation at $\id_W$) the enriched functor sending $X$ to the fiber of $X(W)^{\id_G}$ over the point $[u]\in\cat{L}_V/H$, so that the claim is equivalent to showing that the map
    \begin{equation}\label{eq:upsilon-top-explicit}
        \begin{aligned}
            \hskip-7pt\{h\in H: \beta=c_h\alpha\} &\rightarrow \{v\in\cat{L}(V,W): [u]=[v]\text{ and } g.v=v.\beta(g)\text{ for all $g\in G$}\}\hskip-1pt\\
            h&\mapsto u.h^{-1}
        \end{aligned}
    \end{equation}
    is a weak homotopy equivalence, where both sides carry the subspace topologies. We will show that this map is even a homeomorphism. The source is compact and the target is Hausdorff, so it will suffice that $(\ref{eq:upsilon-top-explicit})$ is bijective. Injectivity is clear as $H$ acts freely on $\cat{L}(V,W)$ by faithfulness of the action on $V$. For surjectivity we observe that if $[v]=[u]$, then there exists by definition an $h\in H$ with $v=u.h^{-1}$, and it only remains that this $h$ defines an element of the left-hand side whenever $v$ is a point in the target of $(\ref{eq:upsilon-top-explicit})$. But indeed, if $g\in G$ is arbitrary, then
    \[
        v.h\alpha(g)=u.\alpha(g)\stackrel{(*)}{=}g.u=g.v.h\stackrel{\smash{(\dagger)}\vphantom{)}}{=}v.\beta(g)h,
    \]
    where $(*)$ uses that $u$ is an $\alpha$-fixed point, $(\dagger)$ uses that $v$ belongs to the target of $(\ref{eq:upsilon-top-explicit})$, and all other identities follow from $v=u.h^{-1}$. Using once more that the $H$-action is free, we conclude that $h\alpha(g)=\beta(g)h$ as desired.
\end{proof}

\begin{thm}\label{thm:Upsilon-ff}
    The functor $\Upsilon$ restricts to an equivalence $\Glo\iso\mathfrak O_\gl\subset\myS_\gl$.
    \begin{proof}
        We already observed above that $\Upsilon(H)\simeq\cat{L}(V,-)/H$ for any faithful $H$-representation $V$. It remains to show that the induced map $\Upsilon\colon\hom(\BGcat{G},\BGcat{H})\to\hom(\Upsilon(\BGcat{G}),\Upsilon(\BGcat{H}))$ is an equivalence for all compact Lie groups $G,H$. Fixing $H$ but varying $G$, we may equivalently show that the composite
        \begin{equation}\label{eq:Upsilon-hom-composite}
            \hom_{\Glo}(-,\BGcat{H})\xrightarrow{\;\Upsilon\;}\hom_{\myS_\gl}(\Upsilon(-),\Upsilon(\BGcat{H}))\iso(\ul\myS_\gl)_{/(\cat{L}(V,-)/H)}^\circ
        \end{equation}
        with the equivalence from Proposition~\ref{prop:the-horrible-proof-we-dont-talk-about} is invertible. By the explicit description of the latter, $(\ref{eq:Upsilon-hom-composite})$ sends $\id_{\BGcat{H}}$ to the unit $(\cat{L}(V,-)\to\cat{L}(V,-)/H)\in\myS_\text{$H$-gl}$, and by the Yoneda lemma this uniquely characterizes the transformation. As the composite
        \[
            \hom_{\Glo}(-,\BGcat{H})\xrightarrow{\;\upsilon\;}\Ntop\big(\Fun_\textup{cont}(-,\cat{$\cat{L}$-Top})_{/(\cat{L}(V,-)/H)}^\circ\big)\xrightarrow{\;\text{loc}\;}(\ul\myS_\gl)_{/(\cat{L}(V,-)/H)}^\circ
        \]
        has the same property, it therefore necessarily agrees with $(\ref{eq:Upsilon-hom-composite})$. The claim now follows from Proposition~\ref{prop:lc-upsilon-equiv}.
    \end{proof}
\end{thm}

\subsection{The universal property of global spaces}\label{subsec:proof-unstable-comp} Putting all the pieces together, we can now prove the desired equivalence $\ul\Spc_\Glo\simeq\ul\myS_\gl$:

\begin{proof}[Proof of Theorem~\ref{thm:unstable-main}]
    It suffices to verify the assumptions of Theorem~\ref{thm:criterion-spc-T}.
    \begin{enumerate}
        \item Each $\myS_\gl(\BGcat{G})\simeq\myS_\text{$G$-gl}$ is complete as the underlying $\infty$-category of a model category. Moreover, each restriction functor $f^*$ admits a left adjoint $\cat{L}f_!$ by Lemma~\ref{lemma:restr-right-Quillen}.
        \item If $G$ is any compact Lie group and $g\colon G\to 1$ denotes the unique map, then $\cat{L}g_!\colon\myS_\text{$G$-gl}\to(\myS_\gl)_{/\cat{L}g_!1}$ is an equivalence by Corollary~\ref{cor:slice-adjoint-equiv}.
        \item The product decomposition of $\cat{L}f_!f^*X$ is the content of Proposition~\ref{prop:f!f*}.
        \item By Theorem~\ref{thm:Upsilon-ff}, $\Upsilon$ defines an equivalence $\Glo\to\mathfrak O_\gl$, and so it extends to an equivalence $\PSh(\Glo)\iso\myS_\gl$ by Proposition~\ref{prop:global-Elmendorf}.
        \qedhere
    \end{enumerate}
\end{proof}

\subsubsection{Applications} We are now finally ready to explain how our definition of $\Glo$ relates to various other definitions in the literature:

\begin{rk}\label{rk:Glo}
    \begin{enumerate}
        \item Combing Theorem~\ref{thm:Upsilon-ff} and Remark~\ref{rk:mathfrak-O-gl-vs-Stefan}, we have an equivalence $\Glo\iso\Ntop(\cat{O}_\gl)$, where $\cat{O}_\gl$ is the model of the global indexing category from \cite{schwede_orbispaces_2020}. Unravelling definitions, this equivalence can be described on the level of homotopy categories as follows: if $G\subset\mathfrak L$ is a universal subgroup, then $\BGcat{G}$ is sent to $\mathfrak L/G$, and a map $\BGcat{G}\to\BGcat{G'}$ corresponding to a homomorphism $f\colon G\to G'$ is sent to $-\circ u\colon\mathfrak L/G\to\mathfrak L/G'$ where $u\colon\R^\infty\to\R^\infty$ is any isometric embedding that is $G'$-equivariant with respect to the $G'$-action on the target as a subgroup of $\mathfrak L$ and the one on the source via restricting the $G$-action. 
        \item \cites{rezk2014global,LNP} work with the topological category $\cat{Glo}'$ whose objects are compact Lie groups and with space of morphisms from $H$ to $G$ given by the geometric realization of the topological groupoid 
        \begin{align*} \maps(H,G)/\hskip-3pt/G\colon \Delta^{\op}&\longrightarrow \Top,\,\\ [n] &\longmapsto \maps(H,G)\times G^n\end{align*} encoding the conjugation action of $G$ on the space of continuous group homomorphisms $\maps(H,G)$ (with compact open topology). 
        By the main result of \cite{koerschgen}, this is Dwyer--Kan equivalent to $\cat{O}_\gl$, and so the corresponding $\infty$-category is equivalent to $\Glo$. By the explicit construction of the equivalence in \emph{loc.\ cit.} as a zig-zag, we see that the induced map on homotopy categories can be described analogously to the equivalence from the previous item. Thus, the composite equivalence $\Ntop(\cat{Glo}')\iso\Glo$ is given on homotopy categories by sending a compact Lie group $G$ to $\BGcat{G}$ and a homomorphism $f\colon G\to G'$ to $\BGcat{f}\colon\BGcat{G}\to\BGcat{G'}$.
        \item \cite{gepnerhenriques2007orbispaces} uses a variant of preceeding definition where they instead take the `fat' realizations of the action groupoid $\maps(H,G)/\hskip-3pt/G$. 
        This is again equivalent to the previous models by \cite{koerschgen}*{Remark~3.10}. 
    \end{enumerate}
\end{rk}

As another application, we can now describe some pullbacks in $\Glo$:

\begin{prop}\label{prop:pb-surj-glo}
    Let
    \begin{equation}\label{diag:pb-surj-glo}
        \begin{tikzcd}
            A\arrow[r,"g"]\arrow[d,"q"']\arrow[dr,pullback] & B\arrow[d,"p"]\\
            C\arrow[r,"f"'] & D
        \end{tikzcd}
    \end{equation}
    be any pullback square in the 1-category of compact Lie groups such that the homomorphism $p$ is \emph{surjective}. Then the corresponding square in $\Glo$ is a pullback.

    In particular, $\Glo$ has finite products, given by products in the 1-category of compact Lie groups.
\end{prop}

A direct proof of this lemma using results from the representation theory of compact Lie groups can be found in \cite{rezk2014global}*{§6.2}. We will give an alternative homotopy theoretic proof exploiting the above comparison.

\begin{proof}
    We saw in Proposition~\ref{prop:basechange-surj} that the associated Beck--Chevalley map $\cat{L}q_!g^*\to f^*\cat{L}p_!$ for $\ul\myS_\gl$ is invertible, whence so is the Beck--Chevalley map $q_!g^*\to f^*p_!$ for $\ul\Spc_{\Glo}$ by Theorem~\ref{thm:unstable-main}. Plugging in the terminal object $1\in\ul\Spc_{\Glo}(\BGcat{B})=\PSh(\Glo)_{/\BGcat{B}}$ and recalling the definition of the functoriality then precisely says that the square
    \[
        \begin{tikzcd}
            \BGcat{A}\arrow[r]\arrow[d] & \BGcat{B}\arrow[d]\\
            \BGcat{C}\arrow[r] & \BGcat{D}
        \end{tikzcd}
    \]
    induced by $(\ref{diag:pb-surj-glo})$ is a pullback in $\PSh(\Glo)$, and hence also in $\Glo$.
\end{proof}

\subsubsection{$\Ff$-global spaces}
Fix a collection $\Ff$ of compact Lie groups closed under isomorphisms. We will now explain how to generalize Theorem~\ref{thm:unstable-main} to give a description of the universal $\Ff$-globally presentable $\Ff$-global $\infty$-category.

\begin{definition}
    A map $f\colon X\to Y$ of closed orthogonal $G$-spaces is called a \emph{$G$-global $\Ff$-weak equivalence} if $f(\Uu_H)\colon X(\Uu_H)\to Y(\Uu_H)$ is a $\mathcal G_{H,G}$-weak equivalence for every $H\in\Ff$ and some (hence any) complete $H$-universe $\Uu_H$. 

    A map $f\colon X\to Y$ between general orthogonal $G$-spaces is called a $G$-global $\Ff$-weak equivalence if in some (hence any) commutative square
    \[
        \begin{tikzcd}
            X'\arrow[d,"\sim"']\arrow[r,"f'"] & Y'\arrow[d,"\sim"]\\
            X\arrow[r,"f"'] & Y
        \end{tikzcd}
    \]
    where the vertical maps are $G$-global equivalences and $X',Y'$ are closed, the map $f'$ is a {$G$-global $\Ff$-weak equivalence} in the sense defined above.
\end{definition}

\begin{rk}
    If $G=1$, the above are referred to as \emph{$\Ff$-global weak equivalences} in \cite{schwede2018global}*{§1.4}. There, the additional assumption is imposed that $\Ff$ be non-empty and closed under subquotients, in which case $\Ff$ is referred to as a \emph{global family}. While the collections $\Ff=\All$ and $\Ff=\Fin$ we care about are in fact global families, the stability under subquotients is not needed for our purposes.
\end{rk}

\begin{lemma}
    Let $\alpha\colon G\to G'$ be any homomorphism of compact Lie groups. Then $\alpha^*\colon\cat{$\bm{G'}$-$\cat{L}$-Top}\to\cat{$\bm G$-$\cat{L}$-Top}$ sends $G'$-global $\Ff$-weak equivalences to $G$-global $\Ff$-weak equivalences.
    \begin{proof}
        In light of Lemma~\ref{lemma:restr-right-Quillen}, it suffices to consider the case of $G'$-global $\Ff$-weak equivalences $f\colon X\to Y$ of \emph{closed} orthogonal $G'$-spaces. But indeed, if $H\in\Ff$ and $\Uu_H$ is any complete $H$-universe, then $\alpha^*f(\Uu_H)$ is a $\mathcal G_{H,G}$-weak equivalence by Example~\ref{ex:graph-restr-rQ}.
    \end{proof}
\end{lemma}

\begin{constr}
    Localizing $\smash{\Ntop(\cat{${\cat{L}}$-Top}^\dual)}$ in each degree $G$ at the $G$-global $\Ff$-weak equivalences and restricting to $\Glo_\Ff$, we obtain an $\Ff$-global $\infty$-category $\ul\myS_\gl^{\,\Ff}\colon\Glo_\Ff^\op\to\Cat_\infty$ given in degree $G\in\Ff$ by the localization $\myS_\text{$G$-gl}^{\,\Ff}$ of $\cat{$\bm G$-$\cat{L}$-Top}$ at the $G$-global $\Ff$-weak equivalences.
\end{constr}

\begin{thm}\label{thm:univ-prop-F-gl-spaces}
    The $\Ff$-global $\infty$-category $\ul\myS_\gl^{\,\Ff}$ is $\Ff$-globally presentable. For any $\Ff$-globally cocomplete $\Ff$-global $\infty$-category $\Dd$, evaluation at the terminal object defines an equivalence
    \[
        \ul\Fun^\textup{$\Glo_\Ff$-cc}(\ul\myS_\gl^{\,\Ff},\Dd)\iso\Dd.
    \]
\end{thm}

To reduce this to the case $\Ff=\All$ treated above we will need:

\begin{lemma}\label{lemma:id-g-fixed-points-vs-ev}
    Let $\Phi\colon\ul\myS_\gl\iso\ul\Spc_\gl$ be the unique equivalence, and fix any homomorphism $\phi\colon H\to G$ of compact Lie groups. Then the composite
    \[
        \myS_\textup{$G$-gl}\xrightarrow[\raise2pt\hbox{$\scriptstyle\smash{\sim}$}]{\;\Phi_G\;} \PSh(\Glo_{/\BGcat{G}})\xrightarrow{\;\ev_\phi\;}\Spc
    \]
    is equivalent to
    \begin{equation}\label{eq:phixed-points}
        \myS_\textup{$G$-gl}\iso\myS_\textup{$G$-gl}^\textup{closed}\xrightarrow{\;(-)^\phi\circ\ev_{\Uu_H}\;}\Spc
    \end{equation}
    for any complete $H$-universe $\Uu_H$.
    \begin{proof}
        Using that $\Phi$ is a global functor, and comparing the functoriality on both sides, it suffices to consider the case where $H=G$ and $\phi=\id_H$. In this case, combining Lemmas~\ref{lemma:fixed-point-corep} and~\ref{lemma:resolve-pt} shows that $(\ref{eq:phixed-points})$ is corepresented by the terminal object. On the other hand, the terminal object of $\PSh(\Glo_{/\BGcat{H}})$ is represented by $\id_H$, so it corepresents evaluation at $\id_H$ by the Yoneda lemma. Adjoining over, it will therefore suffice that the left adjoint of $\Phi_H\colon\myS_\text{$H$-gl}\to\PSh(\Glo_{/\BGcat{H}})$ preserves the terminal object; but as this left adjoint is again an equivalence, this is clear.
    \end{proof}
\end{lemma}

\begin{proof}[Proof of Theorem~\ref{thm:univ-prop-F-gl-spaces}]
    We may equivalently construct an equivalence $\ul\myS_\gl^{\,\Ff}\iso\ul\Spc_{\Glo_\Ff}$ (which will automatically preserve the terminal object). To do so, consider once more the unique equivalence $\Phi\colon\ul\myS_\gl\iso\ul\Spc_\Glo$ and restrict it to a map of $\Ff$-global $\infty$-categories. By construction, we have a localization $\fgt\colon\ul\myS_\gl|_{\Glo_\Ff}\to\ul\myS_\gl^{\,\Ff}$, induced by the identity of $\cat{$\bm G$-$\cat{L}$-Top}$ for varying $G$. On the other hand, the inclusion $\Glo_{\Ff}\hookrightarrow\Glo$ induces a map $\fgt\colon\ul\Spc_{\Glo}|_{\Glo_\Ff}\to\ul\Spc_{\Glo_\Ff}$; this is again a localization, as each $\PSh(\Glo_{/\BGcat{G}})\to\PSh((\Glo_\Ff)_{/\BGcat{G}})$ admits a fully faithful right adjoint (given by right Kan extension along the fully faithful functor $(\Glo_\Ff)_{/\BGcat{G}}\hookrightarrow\Glo_{/\BGcat{G}}$). In summary, we have the solid part of the square
    \[
        \begin{tikzcd}
            \ul\myS_\gl|_{\Glo_\Ff}\arrow[r,"\Phi","\sim"']\arrow[d,"\fgt"'] & \ul\Spc_{\Glo}|_{\Glo_\Ff}\arrow[d,"\fgt"]\\
            \ul\myS_\gl^{\,\Ff}\arrow[r,dashed,"?"{description}] & \ul\Spc_{\Glo_\Ff}
        \end{tikzcd}
    \]
    where both vertical maps are localizations. It will therefore suffice that a map in the top left corner is inverted by the left-hand vertical functor if and only if its image under $\Phi$ is is inverted by the right-hand vertical map. Plugging in the definition of the $G$-global $\Ff$-weak equivalences, this follows at once from the previous lemma.
\end{proof}

\begin{ex}
    For $\Ff=\Fin$ the global family of finite groups, this in particular provides a pointset model of the free presentable $\Fin$-global $\infty$-category. Another such pointset model (via actions of the `universal finite group') was given in \cite{CLL_Global}*{Theorem~3.3.2}. We conclude that these pointset models define equivalent $\infty$-categories, i.e.~for finite $G$ the localization of $\cat{$\bm G$-${\cat L}$-Top}$ at the $G$-global $\Fin$-weak equivalences is modelled by the various pointset models studied in \cite{g-global}*{Chapter~1}, generalizing the comparison for $G=1$ from Corollary~1.5.59 of \emph{op.\ cit.}
\end{ex}

\subsection{The universal property of equivariant spaces}\label{subsec:equivariant-spaces}
In this subsection, we will introduce the global $\infty$-category $\ul\myS$ of \emph{equivariant spaces} and prove its universal property. Let us jump straight to the construction:

\begin{constr}
    Localizing $\smash{\Ntop(\cat{Top}^\dual)}$ at the equivariant weak equivalences we obtain a global $\infty$-category $\ul\myS$ such that $\ul\myS(\BGcat{G})$ is the $\infty$-category $\myS_G$ of $G$-spaces from Example~\ref{ex:G-spaces}, and with restrictions $f^*\colon\myS_{G'}\to\myS_{G}$ given by the derived functors of restriction at the pointset level. We call $\ul\myS$ the \emph{global $\infty$-category of equivariant spaces}.
\end{constr}

Below we will show that $\ul\myS$ is the free so-called \emph{equivariantly cocomplete} global $\infty$-category. Roughly speaking, the difference between equivariant and global cocompleteness is that for the former we only require the existence of left adjoints for restrictions along \emph{injective} homomorphisms. We will formally define the notion of equivariant cocompleteness in §\ref{subsubsec:clefts}, where we will also explain how the free equivariantly cocomplete global $\infty$-category relates to the free \emph{globally} cocomplete global $\infty$-category $\ul\Spc_{\Glo}$. In §\ref{subsubsec:equivariant-spaces-vs-gl}, we will then similarly study the relationship between $\ul\myS$ to $\ul\myS_\gl$ and use this to prove the universal property of $\ul\myS$.

\medskip
\subsubsection{Clefts}\label{subsubsec:clefts}
The notion of equivariant cocompleteness is an instance of the more general notion of \emph{$S$-cocompleteness} of $T$-$\infty$-categories for so-called \emph{clefts} $S\subset T$. While the reader can safely treat the notion of a cleft as a blackbox, we recall its definition for completeness:

\begin{definition}[See \cite{CLL_Clefts}*{Definition~3.2}]
    Let $T$ be a small $\infty$-category. A \emph{cleft} is a wide subcategory $i\colon S\hookrightarrow T$ satisfying all of the following conditions:
    \begin{enumerate}
        \item $S$ is left-cancellable, i.e.~for any composable maps $f,g$ in $T$ such that $f$ and $fg$ belong to $S$,  also $g$ belongs to $S$.
        \item Given any $f\colon A\to B$ in $S$ and $g\colon B'\to B$ in $T$, there exists a pullback square in $\PSh(T)$ of the form
        \[
            \begin{tikzcd}
                i_!X\arrow[d]\arrow[dr,pullback]\arrow[r,"i_!f'"] & i_!B'\arrow[d,"g"]\\
                i_!A\arrow[r, "i_!f"'] & i_!B\rlap,
            \end{tikzcd}
        \]
        where $i_!\colon\PSh(S)\to\PSh(T)$ denotes left Kan extension along $i$ and we conflate elements of $S$ and $T$ with their Yoneda image as before.
        \item Given maps $f\colon A\to B$ and $g\colon B\to A$ in $T$ such that $gf=\id_A$ and $fg$ belongs to $S$, then both $f$ and $g$ belong to $S$.
    \end{enumerate}
\end{definition}
\begin{rk}
By \cite{CLL_Clefts}*{Theorem~3.9}, a wide subcategory $S\subset T$ is a cleft if and only if left Kan extension along the inclusion exhibits $\PSh(S)$ as a \emph{fracture subcategory} of $\PSh(T)$ in the sense of \cite{SAG}*{Definition 20.1.2.1}. 
\end{rk}

\begin{prop}\label{prop:Orb-Glo-cleft}
    Write $\Orb\subset\Glo$ for the wide subcategory given by the \emph{injective} homomorphisms. Then $\Orb\subset\Glo$ is a cleft.
    \begin{proof}
        By \cite{LNP}*{Proposition 6.13}, $\Orb$ is the right class of a factorization system on $\Glo$. The claim is therefore an instance of \cite{CLL_Clefts}*{Proposition~3.33}.
    \end{proof}
\end{prop}

\begin{prop}[See \cite{CLL_Clefts}*{Lemma~3.17}]
    Let $i\colon S\hookrightarrow T$ be a cleft. Then the unique left adjoint $S$-functor $\ul\Spc_S\to i^*(\ul\Spc_T)$ preserving the terminal object is fully faithful, and its essential image is actually a $T$-subcategory of $\ul\Spc_T$.\qed
\end{prop}

\begin{definition}
    We will denote the $T$-subcategory from the previous proposition by $\ul\Spc_{S\triangleright T}$ and call it the \emph{$T$-$\infty$-category of $S$-spaces}.
\end{definition}

The notation `$S\triangleright T$' is to be read as `$S$ extended to $T$,' reflecting that $\ul\Spc_{S\triangleright T}$ is an extension of $\ul\Spc_S$ to a $T$-$\infty$-category.

\begin{rk}\label{rk:Spc-S-T-as-LKE}
    By the explicit construction of the left adjoint $\ul\Spc_S\to i^*(\ul\Spc_T)$ in \cite{CLL_Clefts}*{Construction 3.15}, $\ul\Spc_{S\triangleright T}$ consists in degree $A\in T$ precisely of those presheaves on $T_{/A}$ that are left Kan extended from $S_{/A}$. \end{rk}

\begin{definition}\label{def:S-cc}
    Let $S\subset T$ be a cleft. A $T$-$\infty$-category is called \emph{$S$-cocomplete} if it is fiberwise cocomplete and $\ul\Spc_{S\triangleright T}$-cocomplete, and it is called \emph{$S$-presentable} if it is in addition fiberwise presentable. A functor $F\colon\Cc\to\Dd$ of $S$-cocomplete $T$-$\infty$-categories is called \emph{$S$-cocontinuous} if it is fiberwise cocontinuous and $\ul\Spc_{S\triangleright T}$-cocontinuous.
    We write $\CAT_T^\text{$S$-cc}\subset\CAT_T$ for the subcategory spanend by the $S$-cocomplete $T$-$\infty$-categories and $S$-cocontinuous functors. We moreover write $\Pr^S_T\subset\CAT_T^\text{$S$-cc}$ for the full subcategory spanned by the $S$-presentable $T$-$\infty$-categories.
\end{definition}

\begin{ex}
    For the minimal cleft $S=\core T$, the above simplifies to the corresponding fiberwise notions.     
    On the other hand, for the maximal cleft $S=T$, this recovers the notions of $T$-cocompleteness and $T$-presentability. We will employ both the notations $\Pr^T_T$ and $\PrLT$, depending on whether we want to emphasize that the morphisms are $T$-cocontinuous functors or left adjoints.
\end{ex}

\begin{ex}
    Specializing Definition~\ref{def:S-cc} to the cleft $\Orb\subset\Glo$ from Proposition~\ref{prop:Orb-Glo-cleft} yields the notions of \emph{equivariant cocompleteness} and \emph{equivariant presentability} for global $\infty$-categories.
\end{ex}

Next, we will describe the free $S$-cocomplete $T$-$\infty$-category, generalizing the case $S=T$ considered in Corollary~\ref{cor:univ-prop-Spc-T}. For this we first observe:

\begin{rk}[See \cite{CLL_Clefts}*{Lemma 4.15}]
    Let $\Cc,\Dd$ be $S$-cocomplete $T$-$\infty$-categories. Generalizing the case $S=T$ considered in §\ref{subsubsec:colimits}, the full subcategories of $\ul\Fun_T(\Cc,\Dd)(X)$ spanned by the $S$-cocontinuous functors $\Cc\to\ul\Fun(\ul X,\Dd)$ for $X\in\PSh(T)$ define a $T$-subcategory of $\ul\Fun_T(\Cc,\Dd)$, which we denote by $\ul\Fun_T^\text{$S$-cc}(\Cc,\Dd)$.
\end{rk}

\begin{prop}[See \cite{CLL_Clefts}*{Corollary~4.27}]\label{prop:Spc-S-T-univ-prop}
    Let $S\subset T$ be a cleft. The $T$-$\infty$-category $\ul\Spc_{S\triangleright T}$ is $S$-presentable. Moreover, for every $S$-cocomplete $T$-$\infty$-category $\Dd$, evaluation at the terminal object yields an equivalence
    \[
        \ul\Fun^\textup{$S$-cc}_T(\ul\Spc_{S\triangleright T},\Dd)\iso\Dd.\qednow
    \]
\end{prop}

\begin{rk}\label{rk:PSh-ST}
    More generally, the free $S$-cocomplete $T$-$\infty$-category generated by a small $T$-$\infty$-category $I$ is $S$-presentable, and it can be analogously described as the essential image of $\ul\PSh_S(i^*I)\hookrightarrow\ul\PSh_T(I)$, see \cite{CLL_Clefts}*{Corollary~4.26}.
\end{rk}

We close this discussion by giving pointwise characterizations of $S$-{\hskip0pt}cocompleteness and $S$-cocontinuity. Note that even for $S=T$ these improve on the characterizations one would obtain by na\"ively applying Lemmas~\ref{lemma:pointwise-crit-cocomplete} and~\ref{lemma:pointwise-criterion-cocontinuous}.

\begin{lemma}[See \cite{CLL_Clefts}*{Lemma~4.9}]\label{lemma:S-cocompleteness-repr}
A $T$-$\infty$-category $\mathcal C$ is $S$-cocomplete if and only if the following three conditions are satisfied: 
\begin{enumerate}
\item The $T$-$\infty$-category $\mathcal C$ is fiberwise cocomplete. 
\item For every morphism $s\colon A\to B\in S$, the restriction $s^*\colon\mathcal C(B)\to \mathcal C(A)$ admits a left adjoint $s_{!}$. 
\item For any pullback
    \[
        \begin{tikzcd}
            X\arrow[d,"s'"']\arrow[dr,pullback]\arrow[r,"f'"] & A\arrow[d,"s"]\\
            C\arrow[r,"f"'] & B
        \end{tikzcd}
    \]
    in $\PSh(T)$ with $s$ in $S$ and $C\in T$, the restriction along $s'$ admits a left adjoint $s'_!$, and the Beck-Chevalley map $s'_!f'^*\to f^*s_!$ is an equivalence. 
\end{enumerate}
\end{lemma}
\begin{rk}\label{rk:cleft-cc-addendum}
By \cite{CLL_Clefts}*{Lemma~4.7}, the existence of the left adjoint $s'_{!}$ in $(3)$ of the above lemma follows from $(1)$ and $(2)$. More generally, \emph{loc.\ cit.} shows that under these assumptions a left adjoint exists for restriction along any map $f$ in the image of $\PSh(S)$ with representable target.
\end{rk}

\begin{lemma}[See \cite{CLL_Clefts}*{Lemma~4.11}]\label{lm:continuouity-clefts}
    A functor $F\colon\Cc\to\Dd$ of $S$-\hskip0ptcocomplete $T$-$\infty$-categories is $S$-cocontinuous if and only if it is fiberwise cocontinuous and the Beck--Chevalley map $s_!F\to Fs_!$ is invertible for every map $s$ in $S$.\qed
\end{lemma}

\subsubsection{Equivariant spaces vs.~global spaces}\label{subsubsec:equivariant-spaces-vs-gl}
The above already shows that the free equivariantly cocomplete global $\infty$-category $\ul\Spc_{\Orb\triangleright\Glo}$ is a full subcategory of $\ul\Spc_{\Glo}$. Thus, to compare it to $\ul\myS$, it will suffice to similarly embed $\ul\myS$ into $\ul\myS_\gl$, and then check that the two subcategories match up under the unique equivalence $\ul\Spc_{\Glo}\simeq\ul\myS_\gl$ from Theorem~\ref{thm:unstable-main}. To address the former, we will now show:

\begin{prop}\label{prop:triv-ff}
    Let $G$ be a compact Lie group. The functor
    \[
        \consto\colon\cat{$\bm G$-Top}\to\cat{$\bm G$-$\cat{L}$-Top}
    \] sending a topological $G$-space to the associated constant orthogonal $G$-space 
    is homotopical and descends to a fully faithful left adjoint $\myS_G\to\myS_\textup{$G$-gl}$.
\end{prop}

The proof will rely on the following model-categorical observation.

\begin{lemma}
    Let $V$ be a faithful $G$-representation. Then
    \begin{equation}\label{eq:const-made-left-Quillen}
        \cat{L}(V,-)\times\consto(-)\colon \cat{$\bm G$-Top}\to\cat{$\bm G$-${\cat L}$-Top}
    \end{equation}
    is left Quillen, and its left derived functor $\myS_G\to\myS_\textup{$G$-gl}$ is fully faithful.
    \begin{proof}
        It is clear that $\cat{L}(V,-)\times\consto(-)\dashv\ev_V$ is an adjunction with unit $\eta\colon X\to\cat{L}(V,V)\times X$ given by the inclusion of $\id_V\in\cat{L}(V,V)$. Combining Lemmas~\ref{lemma:orth-spc-ev-lrQ} and~\ref{lemma:restr-equiv-rQ} shows that $\ev_V$ is right Quillen for the $G$-global \emph{level} model structure; as the $G$-global model structure has the same cofibrations and more acyclic cofibrations (and hence the same acyclic fibrations and fewer fibrations) as the $G$-global level model structure, $\ev_V$ is then also right Quillen for the $G$-global model structure on $\cat{$\bm G$-$\cat{L}$-Top}$.
        
        For full faithfulness, it will suffice to show that the derived unit is invertible, for which we fix a fibrant replacement $i\colon\cat{L}(V,-)\to\mathcal L$. Then $\cat{L}(V,-)\times X\to\mathcal L\times X$ is a $G$-global weak equivalence for every $G$-space $X$ by Corollary~\ref{cor:cart-prod-homotopical}, and the target is fibrant by direct inspection; thus, the derived unit is given by the composite
        \[
            X\xrightarrow{\;x\mapsto(\id_V,x)\;}\cat{L}(V,V)\times X\xrightarrow{\;i(V)\times X\;}\mathcal L(V)\times X,
        \]
        and it suffices to show that the $G$-space $\mathcal L(V)$ (where $G$ acts diagonally) is $G$-equivariantly weakly contractible. For this, note that $\mathcal L\to 1$ is a $G$-global weak equivalence (by Lemma~\ref{lemma:resolve-pt}) and a $G$-global fibration (by choice of $\mathcal L$), hence also an acyclic fibration in the $G$-equivariant \emph{level} model structure. The claim therefore follows from the definition of the $G$-global level weak equivalences.
    \end{proof}
\end{lemma}

\begin{proof}[Proof of Proposition~\ref{prop:triv-ff}]
    Fix a faithful $G$-representation $V$. Combining Lemma~\ref{lemma:resolve-pt} with Corollary~\ref{cor:cart-prod-homotopical}, the projection $\cat L(V,-)\times X\to\consto(X)$ is a $G$-global weak equivalence for any $X\in\cat{$\bm G$-$\cat{L}$-Top}$. The claim therefore follows immediately from the previous lemma.
\end{proof}

\begin{rk}\label{rk:ev-U-ra-const}
    The above proof shows that the right adjoint $\myS_\text{$G$-gl}\to\myS_G$ of $\consto$ is given as the right derived functor of $\ev_V$ for any faithful $G$-representation $V$. As in the proof of Lemma~\ref{lemma:id-g-fixed-points-vs-ev} we may equivalently describe this as the \emph{left} derived functor of $\ev_{\Uu}$ for any complete $G$-universe $\Uu$.
\end{rk}

With this at hand, we can now show:

\begin{thm}\label{thm:unstable-equiv-main}
    The global $\infty$-category $\ul\myS$ is equivariantly presentable, and the unique equivariantly cocontinuous functor $\ul\Spc_{\Orb\triangleright\Glo}\to\ul\myS$ preserving the terminal object is an equivalence.
    \begin{proof}
        The topological functor $\consto$ induces a global functor $\smash{\Ntop(\cat{Top}^\dual)}\to\smash{\Ntop(\cat{$\cat{L}$-Top}^\dual)}$, which by Proposition~\ref{prop:triv-ff} descends to a fully faithful and levelwise left adjoint functor $\consto\colon\ul\myS\to\ul\myS_\gl$. To complete the proof, it will now suffice to show that the essential image of the composite
        \[
            \Psi\colon\ul\myS\xrightarrow{\;\consto\;}\ul\myS_\gl\xrightarrow[\smash{\raise3.5pt\hbox{$\scriptstyle\sim$}}]{\;\,\Phi\,\;}\ul\Spc_{\Glo}
        \]
        agrees with $\ul\Spc_{\Orb\triangleright\Glo}\subset\ul\Spc_{\Glo}$, where $\Phi$ denotes the unique equivalence.

        To see that $\Psi$ factors through $\ul\Spc_{\Orb\triangleright\Glo}$, it suffices by Elmendorf's Theorem to show that the left adjoint $\Psi_G\colon\myS_G\to\ul\Spc_{\Glo}(G)$ maps $G/H$ into $\ul\Spc_{\Orb}(G)$ for any compact Lie group $G$ and any closed subgroup $H\subset G$. To this end, let us write $i\colon H\hookrightarrow G$ for the inclusion, so that $\consto(G/H)\simeq\cat{L}i_!\consto(1_H)$ by Proposition~\ref{prop:quotient-we-global} and hence $\Psi_G(G/H)=\Phi_G\cat{L}i_!\consto(1_H)=i_!\Phi_H\consto(1)=i_!1\in\ul\Spc_{\Orb\triangleright\Glo}(G)$.

        For the converse, note that $\ul\Spc_{\Orb}(\BGcat{G})\simeq\PSh(\Orb_{/\BGcat{G}})$ is generated under ordinary colimits by the objects $(i\colon H\hookrightarrow G)=i_!1_H$ for closed subgroups $H\subset G$, so that also $\ul\Spc_{\Orb\triangleright\Glo}(\BGcat{G})$ is generated by the objects $i_!1_H$. We have just seen that $\Psi_G(G/H)\simeq i_!1_H$; thus, the essential image of the fully faithful left adjoint $\Psi_G$ contains a set of generators of $\ul\Spc_{\Orb\triangleright\Glo}(\BGcat{G})$, so it already has to be all of $\ul\Spc_{\Orb\triangleright\Glo}(\BGcat{G})$.
    \end{proof}
\end{thm}

Combining this with Proposition~\ref{prop:Spc-S-T-univ-prop} we now get:

\begin{cor}\label{cor:equiv-spaces-universal}
    For every equivariantly cocomplete global $\infty$-category $\Dd$, evaluation at the terminal object defines an equivalence
    \[
        \ul\Fun^{\textup{Orb-cc}}_{\Glo}(\ul\myS,\Dd)\iso\Dd.\qednow
    \]
\end{cor}

\subsection{Symmetric monoidal structures}\label{subsec:sym-mon}
As our final results in this section, we will show that $\ul\myS_\gl$ and $\ul\myS$ admit unique (levelwise) symmetric monoidal structures such the tensor product preserves global respectively equivariant colimits in each variable in a sense that will be made precise below, and that the resulting global symmetric monoidal categories again admit suitable universal properties.

\medskip
\subsubsection{Symmetric monoidal $T$-$\infty$-categories} We begin by introducing symmetric monoidal structures in the parametrized setting and discussing some examples.

\begin{definition}
    A \emph{symmetric monoidal $T$-$\infty$-category} is an object of the $\infty$-category $\CMon(\CAT_T)$. We refer to the morphisms of this category as \emph{symmetric monoidal $T$-functors}.
\end{definition}

\begin{notation}
    We will typically denote symmetric monoidal $T$-$\infty$-categories by $\Cc$, $\Dd$,\,\dots\ again. If we want to emphasize the difference between a symmetric monoidal $T$-$\infty$-category and its underlying $T$-$\infty$-category, we will decorate the former with a superscript `${}^\otimes$' (e.g. $\Cc^\otimes$, $\Dd^\otimes$,\,\dots).
\end{notation}

\begin{rk}
    By currying, a symmetric monoidal $T$-$\infty$-category can be equivalently described as a functor $T^\op\to\CMon(\CAT_\infty)$, and we will freely switch between these two perspectives.
\end{rk}

\begin{ex}\label{ex:cartesian-sym-mon}
    Let $\Cc$ be a $T$-$\infty$-category with fiberwise finite products, i.e.~such that $\Cc$ factors through the non-full subcategory $\CAT^{\Pi}_\infty\subset\CAT_\infty$ of $\infty$-categories with finite products and functors preserving them. Postcomposing with the functor $(-)^\times\colon\CAT^\Pi_\infty\to\CMon(\CAT_\infty)$ equipping an $\infty$-category with finite products with the cartesian symmetric monoidal structure, we may then upgrade $\Cc$ to a symmetric monoidal $T$-$\infty$-category $\Cc^\times\colon T^\op\to\CMon(\CAT_\infty)$. We will refer to this as the \emph{cartesian $T$-symmetric monoidal structure on $\Cc$}.
\end{ex}

\begin{rk}\label{rk:what-it-means-to-be-cartesian}
    As $(-)^\times\colon\CAT^\Pi_\infty\to\CMon(\CAT_\infty)$ is fully faithful, it is a \emph{property} for a symmetric monoidal $T$-$\infty$-category $\Cc$ to be cartesian, and this is moreover equivalent to each individual $\Cc(A)\in\CMon(\CAT_\infty)$ being cartesian. By \cite{CHLL_Bispans}*{Remark~2.4.2${}^\op$}, this is in turn equivalent to the unit being terminal and the tensor product $\Cc(A)\times\Cc(A)\to\Cc(A)$ admitting a left adjoint.
\end{rk}

Another family of examples comes from a refinement of the continuous Borel construction. We will begin by explaining how to associate a symmetric monoidal $\infty$-category to a \emph{topological} symmetric monoidal category (i.e.\ a symmetric monoidal category with a topological enrichment of the underlying category such that the tensor product is an enriched functor). While this could be done explicitly by adapting the constructions from \cite{HA}*{§2.1.1}, we prefer a more abstract approach.

\begin{constr}
    We write $\cat{topSymMonCat}$ for the $1$-category of small topological symmetric monoidal categories and enriched strong symmetric monoidal functors. For all $\cat{C}^\otimes,\cat{D}^\otimes\in\cat{topSymMonCat}$ we define $\cat{Fun}_\cont^\otimes(\cat{C}^\otimes,\cat{D}^\otimes)$ as the category of enriched strong symmetric monoidal functors and symmetric monoidal transformations, with each morphism space $\maps(F,G)$ topologized as a subspace of $\prod_{c\in\cat{C}}\maps(F(c),G(c))$ as before. We omit the straightforward verification that this makes $\cat{topSymMonCat}$ into a $\cat{topCat}$-enriched category, which is moreover cotensored over $\cat{topCat}$, with cotensoring given by sending $(\cat{K}\in\cat{topCat},\cat{C}^\otimes\in\cat{topSymMonCat})$ to $\FUN_\cont(\cat{K},\cat{C})$ with the pointwise symmetric monoidal structure, and with the evident (enriched) functoriality in $\cat{K}$ and $\cat{C}^\otimes$. We define the simplicially enriched category $\cat{topSymMonCat}_\Delta$ by transporting the enrichment along $\core\Ntop$ as before.
\end{constr}

\begin{lemma}
    The $\infty$-category $N_\Delta(\cat{topSymMonCat}_\Delta)$ is semiadditive, with products computed in the 1-category of topological categories.
    \begin{proof}
        It is clear that the projections induce even an isomorphism 
        \[
            \FUN_\cont^\otimes(\cat{C}^\otimes,\cat{D}^\otimes\times\cat{E}^\otimes)\xrightarrow{\;\,\smash{\cong}\,\;}\FUN_\cont^\otimes(\cat{C}^\otimes,\cat{D}^\otimes)\times\FUN_\cont^\otimes(\cat{C}^\otimes,\cat{E}^\otimes)
        \]
        for all topological symmetric monoidal categories $\cat{C}^\otimes,\cat{D}^\otimes,\cat{E}^\otimes$, and hence exhibit $\cat{D}^\otimes\times\cat{E}^\otimes$ as a product. Analogously, one sees that $1$ (with the unique symmetric monoidal structure) is terminal.
    
        One directly checks for any $\cat{C}^\otimes$, the map $\bbone\colon 1\to\cat{C}^\otimes$ classifying the symmetric monoidal unit admits a (unique) symmetric monoidal structure; any other symmetric monoidal functor $1\to\cat{C}^\otimes$ is isomorphic to $\bbone$ by the unitality condition for symmetric monoidal functors, and the unitality condition for symmetric monoidal \emph{transformations} implies that this isomorphism is unique. Thus, $\FUN_\cont^\otimes(\cat{1},\cat{C}^\otimes)$ is equivalent in the enriched sense (and hence in particular Dwyer--Kan equivalent) to the terminal category for any $\cat{C}^\otimes$, so that $1$ is also initial.

        Finally, we will show that for any $\cat{C}^\otimes,\cat{D}^\otimes,\cat{E}^\otimes$ the map
        \begin{equation}\label{eq:restrict-unit}
            \FUN_\cont^\otimes(\cat{C}^\otimes\times\cat{D}^\otimes,\cat{E}^\otimes)\to\FUN_\cont^\otimes(\cat{C}^\otimes,\cat{E}^\otimes)\times\FUN_\cont^\otimes(\cat{D}^\otimes,\cat{E}^\otimes)
        \end{equation}
        induced by the inclusions $i_1\colon\cat{C}^\otimes\to\cat{C}^\otimes\times\cat{D}^\otimes\gets\cat{D}^\otimes\noloc i_2$ of the units is an equivalence again, so that $i_1$ and $i_2$ exhibit $\cat{C}^\otimes\times\cat{D}^\otimes$ as the coproduct, which will then complete the proof of semiadditivity. For this, let us make the tensor product $\cat{E}\times\cat{E}\to\cat{E}$ into a strong symmetric monoidal functor $\mu$ via the associativitiy and symmetry isomorphisms of $\cat{E}^\otimes$, and consider the composite
        \begin{multline}\label{eq:multiply}
            \FUN_\cont^\otimes(\cat{C}^\otimes,\cat{E}^\otimes)\times\FUN_\cont^\otimes(\cat{D}^\otimes,\cat{E}^\otimes)\xrightarrow{\;{-}\times{-}\;}
            \FUN_\cont^\otimes(\cat{C}^\otimes\times\cat{D}^\otimes,\cat{E}^\otimes\times\cat{E}^\otimes)\\\xrightarrow{\;\mu\circ{-}\;}\FUN_\cont^\otimes(\cat{C}^\otimes\times\cat{D}^\otimes,\cat{E}^\otimes).
        \end{multline}
        For any $F\colon\cat{C}^\otimes\times\cat{D}^\otimes\to\cat{E}^\otimes$ and any $X\in\cat{C},Y\in\cat{D}$ we then have isomorphisms
        \[
            F(X,Y)\cong F\big((X,\bbone)\otimes(\bbone,Y)\big)\cong F(X,\bbone)\otimes F(\bbone,Y)
        \]
        via the unitality isomorphisms of $\cat{C}^\otimes$ and $\cat{D}^\otimes$ and the symmetric monoidal structure of $F$. It is then a routine, but somewhat lengthy argument to check that these assemble into a symmetric monoidal isomorphism $F\cong \mu\circ(Fi_1\times Fi_2)$ and that this isomorphism is natural in $F$, exhibiting $(\ref{eq:multiply})$ as left inverse to $(\ref{eq:restrict-unit})$. Similarly, one uses the unitality isomorphisms to exhibit $(\ref{eq:multiply})$ also as \emph{right} inverse to $(\ref{eq:restrict-unit})$.
    \end{proof}
\end{lemma}

Specializing \cite{Gepner-Groth-Nikolaus}*{Corollary~2.5} we therefore get:

\begin{cor}\label{cor:Ntop-sym-mon-lift}
    There exists a unique (dashed) lift in the diagram
    \[
        \begin{tikzcd}[anchor=south, baseline=.575em]
            N_\Delta(\cat{topSymMonCat}_\Delta)\arrow[r,dashed]\arrow[d,"\fgt"'] & \CMon(\Cat_\infty)\arrow[d,"\fgt"]\\
            N_\Delta(\cat{topCat}_\Delta)\arrow[r,"\Ntop"'] & \Cat_\infty\rlap.
        \end{tikzcd}
        \qednow
    \]
\end{cor}

\begin{rk}
    The explicit construction of $(\ref{eq:multiply})$ shows that the codiagonal $\cat{E}^\otimes\times\cat{E}^\otimes\to\cat{E}^\otimes$ is given for any $\cat{E}^\otimes\in\Ntop(\cat{topSymMonCat}_\Delta)$ by the tensor product. It follows that the tensor product on $\Ntop(\cat{E}^\otimes)$ is given by the composite
    \[
        \Ntop(\cat{E})\times\Ntop(\cat{E})\iso\Ntop(\cat{E}\times\cat{E})\xrightarrow{\;\Ntop(\otimes)\;}\Ntop(\cat{E}).
    \]
    In particular, we see that whenever $\cat{C}^\otimes$ is equipped with a collection $\Ww$ of weak equivalences stable under the tensor product, then the symmetric monoidal $\infty$-category $\Ntop(\cat{E}^\otimes)$ localizes to a symmetric monoidal enhancement of $\Ntop(\cat{E})[\Ww^{-1}]$.
\end{rk}

We are now ready to upgrade the continuous Borel construction to a functor $\Ntop(\cat{topSymMonCat}_\Delta)\to\CMon(\Cat_\Glo)$.

\begin{constr}\label{constr:sym-mon-Borel}
    Restricting the cotensoring of $\cat{topSymMonCat}$, we obtain an enriched functor $\cat{topSymMonCat}_\Delta\to\FUN(\cat{Glo},\cat{topSymMonCat}_\Delta)$, given on objects by sending $\cat{C}^\otimes$ to $\FUN_\cont(-,\cat{C}^\otimes)$. The composite
    \begin{multline*}
        N_\Delta(\cat{topSymMonCat}_\Delta)\to N_\Delta\big(\FUN(\cat{Glo},\cat{topSymMonCat}_\Delta)\big)\\
        \to\Fun\big(\Glo, N_\Delta(\cat{topSymMonCat}_\Delta)\big)\to
        \Fun\big(\Glo,\CMon(\Cat_\infty)\big)
    \end{multline*}
    with the natural comparison map and the lift of $\Ntop$ from Corollary~\ref{cor:Ntop-sym-mon-lift} then provides a lift of the continuous Borel construction. We will again denote this lift on objects by $\cat{C}^\otimes\mapsto\Ntop(\cat{C}^{\otimes,\dual})$.
\end{constr}

\subsubsection{The parametrized Lurie tensor product}\label{subsubsec:param-Lurie} Below, we will typically not be interested in general symmetric monoidal structures, but instead only with those where the tensor product preserves certain parametrized colimits separately in each variable. We begin by formally defining this notion:

\begin{definition}
    Let $S\subset T$ be a cleft and let $\Cc_1,\Cc_2,\Dd$ be $S$-cocomplete $T$-$\infty$-categories. A functor $\Cc_1\times\Cc_2\to\Dd$ is called \emph{$S$-cocontinuous in each variable} or \emph{$S$-bilinear} if for every $A\in T$ and $X_1\in\Cc_1(A), X_2\in\Cc_2(A)$ the two composites
    \begin{align*}
        \{X_1\}\times\pi_A^*\Cc_2&\lhook\joinrel\longrightarrow \pi_A^*\Cc_1\times\pi_A^*\Cc_2\xrightarrow{\;\pi_A^*F\;}\pi_A^*\Dd\\
        \pi_A^*\Cc_1\times\{X_2\}&\lhook\joinrel\longrightarrow\pi_A^*\Cc_1\times\pi_A^*\Cc_2\xrightarrow{\;\pi_A^*F\;}\pi_A^*\Dd
    \end{align*}
    preserve $(T_{/A}\times_TS)$-colimits in the usual sense. Analogously, we define what it means for $\Cc_1\times\cdots\times\Cc_n\to\Dd$ to be $S$-cocontinuous in each variable (\emph{$S$-multilinear}).
\end{definition}

\begin{rk}
    Applying Lemma~\ref{lemma:pointwise-criterion}, one easily sees that $F\colon\Cc_1\times\Cc_2\to\Dd$ is $S$-bilinear if and only if each individual functor $\Cc_1(A)\times\Cc_2(A)\to\Dd(A)$ preserves ordinary colimits in each variable and we moreover have \emph{projection formul\ae} in the sense that the Beck--Chevalley maps
    \[
        f_!F(f^*X_1,-)\to F(X_1,f_!(-))\qquad\text{and}\qquad
        f_!F(-,f^*X_2)\to F(f_!(-),X_2)
    \]
    associated to the naturality squares
    \[
        \begin{tikzcd}[column sep=4.5em]
            \Cc_2(B)\arrow[r, "{F(X_1,-)}"]\arrow[d,"f^*"'] & \Dd(B)\arrow[d,"f^*"] & \Cc_1(B)\arrow[r,"{F(-,X_2)}"]\arrow[d,"f^*"'] & \Dd(A)\arrow[d,"f^*"]\\
            \Cc_2(A)\arrow[r,"{F(f^*X_1,-)}"'] & \Dd(A) & \Cc_1(A)\arrow[r,"{F(-,f^*X_2)}"'] & \Dd(A)
        \end{tikzcd}
    \]
    are equivalences for every $f\colon A\to B$ in $S$, $X_1\in\Cc_1(B)$, and $X_2\in\Cc_2(A)$. 
\end{rk}

\begin{rk}
    Spelling out the definitions, the first Beck--Chevalley map from the previous remark is given as the composite
    \[
        f_!F(f^*X_1,X_2)\hskip0pt minus .75pt\xrightarrow{f_!F(f^*X_1,\eta)\hskip0pt plus .5pt}\hskip0pt minus .75pt
        f_!F(f^*X_1,f^*f_!X_2)\simeq
        f_!f^*F(X_1,f_!X_2)\hskip0pt minus .75pt\xrightarrow{\epsilon\hskip0pt plus .5pt}\hskip0pt minus .75ptF(X_1,f_!X_2),
    \]
    and in particular we see that it upgrades (in a preferred way) to a natural transformation of functors $\Cc_1(B)\times\Cc_2(A)\to\Dd(B)$ (whereas the original construction is a priori only natural in the second variable). Similarly, one shows that also the other Beck--Chevalley map `is' natural in both variables.
\end{rk}

\begin{rk}\label{rk:charact-bilinear}
    By \cite{martiniwolf2022presentable}*{Lemma 2.6.1.3}, the $T$-$\infty$-categories $\ul\Fun^\text{$S$-cc}_T(\Cc_i,\Dd)$ for $i=1,2$ are again $S$-cocomplete (with the inclusion into $\ul\Fun_T(\Cc_i,\Dd)$ being $S$-cocontinuous), and we can equivalently characterize $S$-bilinear maps $\Cc_1\times\Cc_2\to\Dd$ as those that adjoin to an $S$-cocontinuous map $\Cc_1\to\ul\Fun^\text{$S$-cc}(\Cc_2,\Dd)$ or equivalently $\Cc_2\to\ul\Fun^\text{$S$-cc}(\Cc_1,\Dd)$.
\end{rk}

Consider now the cartesian symmetric monoidal structure $\CAT_T^\times\to\text{Fin}_*$ on the (very large) $\infty$-category $\CAT_T$. We define a subcategory $\CAT_T^{\text{$S$-cc},\otimes}$ whose objects are the $S$-cocomplete $T$-$\infty$-categories, and such that the maps $(\Cc_1,\dots,\Cc_n)\to(\Dd_1,\dots,\Dd_m)$ over an $\alpha\colon\langle n\rangle\to\langle m\rangle$ in $\text{Fin}_*$ are those maps such that each component $\prod_{i\in\alpha^{-1}(j)}\Cc_i\to\Dd_j$ is $S$-cocontinuous in each variable.

\begin{prop}[See \cite{martiniwolf2022presentable}*{Corollary 2.6.2.5 and Remark~2.6.2.6}]\label{prop:lurie-tensor-exists}
    The composite $\CAT_T^{\textup{$S$-cc},\otimes}\hookrightarrow\CAT_T^\times\to\textup{Fin}_*$ is the cocartesian fibration for a symmetric monoidal structure on $\CAT_T^\textup{$S$-cc}$ with unit $\ul\Spc_{S\triangleright T}$.\qed
\end{prop}

We will refer to the tensor product of this symmetric monoidal structure as the \emph{parametrized Lurie tensor product}.

\begin{rk}\label{rk:CAlg-for-param-Lurie}
    By \cite{martiniwolf2022presentable}*{Remark 2.6.2.7}, $\CAT_T^{\textup{$S$-cc},\otimes}\hookrightarrow\CAT_T^\times$ is lax symmetric monoidal, and the induced functor $\CAlg(\CAT_T^{\textup{$S$-cc},\otimes})\hookrightarrow\CAlg(\CAT_T^\times)$ defines an equivalence onto the subcategory whose objects are the $S$-cocomplete symmetric monoidal $T$-$\infty$-categories such that the tensor product preserves $S$-colimits in each variable separately, and whose morphisms are the $S$-cocontinuous strong symmetric monoidal functors.
\end{rk}

\begin{rk}\label{rk:S-cc-closed}
    It follows from Remark~\ref{rk:charact-bilinear} that the symmetric monoidal structure afforded by the parametrized Lurie tensor product is closed, with internal homs given by $\ul\Fun^\text{$S$-cc}(-,-)$.
\end{rk}

Recall that the usual Lurie tensor product on $\CAT_\infty^\text{cc}$ restricts to a symmetric monoidal product on $\PrL$; the next goal will be to show the parametrized analogue of this result. Our proof will proceed by reducing to the case $T=S$ which was already handled by Martini and Wolf; we therefore begin by showing:

\begin{lemma}\label{lemma:restr-cleft-sym-mon}
    Let $i\colon S\hookrightarrow T$ be a cleft. Then the strong symmetric monoidal functor $i^*\colon\CAT_T^\times\to\CAT_S^\times$ restricts to a strong symmetric monoidal functor $\CAT_T^{\textup{$S$-cc},\otimes}\to\CAT_S^{\textup{$S$-cc},\otimes}$.
    \begin{proof}
        This proof is an elaboration on \cite{CLL_Clefts}*{Remark 4.23}.  Unravelling definitions, the lemma amounts to proving the following:
        \begin{enumerate}
            \item The unique $S$-cocontinuous $T$-functor $\ul\Spc_S\to i^*(\ul\Spc_{S\triangleright T})$ preserving the terminal object is an equivalence.
            \item For any $\Cc_1,\Cc_2\in\CAT_T^\text{$S$-cc}$, `the' universal $S$-bilinear map $\beta\colon\Cc_1\times\Cc_2\to\Dd$ forgets to the universal $S$-bilinear map $i^*\Cc_1\times i^*\Cc_2\to i^*\Dd$.
        \end{enumerate}
        The first statement is clear, while for the second statement it will suffice to show that for any $S$-cocomplete $S$-$\infty$-category $\Dd$ the induced map 
        \[
            -\circ i^*\beta\colon\Fun_S^\text{$S$-bilinear}(i^*\Cc_1\times i^*\Cc_2,\Ee)\to\Fun_S^\text{$S$-cc}(i^*\Dd,\Ee)
        \]
        is an equivalence. By \cite{CLL_Clefts}*{Construction~4.21 and Corollary~4.22}, the equivalence $\ul\Fun_T(\Dd,i_*\Ee)\iso i_*\ul\Fun_S(i^*\Dd,\Ee)$ obtained as the total mate of the obvious equivalence $i^*\Dd\times i^*(-)\cong i^*(\Dd\times-)$ restricts to $\ul\Fun_T^\text{$S$-cc}(\Dd,i_*\Ee)\iso i_*\ul\Fun_S^\text{$S$-cc}(i^*\Dd,\Ee)$, and evaluating at the terminal presheaf $1\in\PSh(T)$ yields a natural equivalence $\Fun_T(\Dd,i_*\Ee)\iso \Fun_S(i^*\Dd,\Ee)$ restricting to $\Fun_T^\text{$S$-cc}(\Dd,i_*\Ee)\iso \Fun_S(i^*\Dd,\Ee)$. On the other hand, applying this observation twice, Remark~\ref{rk:charact-bilinear} shows that also the equivalence $\Fun_S(i^*\Cc_1\times i^*\Cc_2,\Ee)\iso\Fun_T(\Cc_1\times\Cc_2,i_*\Ee)$ restricts to subcategories of $S$-bilinear functors. In summary, we obtain a commutative square
        \[
            \begin{tikzcd}
                 \Fun_S^\text{$S$-bilinear}(i^*\Cc_1\times i^*\Cc_2,\Ee)\arrow[r,"-\circ i^*\beta"]\arrow[d,"\sim"'] &[.5em] \Fun_S^\text{$S$-cc}(i^*\Dd,\Ee)\arrow[d,"\sim"]\\[-.05ex]
                 \Fun_T^\text{$S$-bilinear}(\Cc_1\times \Cc_2,i_*\Ee)\arrow[r,"-\circ\beta"'] & \Fun_T^\text{$S$-cc}(\Dd,i_*\Ee)
            \end{tikzcd}
        \]
        The bottom map is invertible by the defining property of $\beta$, whence so is the upper map by $2$-out-of-$3$, as desired.
    \end{proof}
\end{lemma}

\begin{prop}
    The symmetric monoidal structure on $\CAT_T^\textup{$S$-cc}$ given by the parametrized Lurie tensor product restricts to a symmetric monoidal structure on the full subcategory $\Pr^S_T$.
    \begin{proof}
        It is clear that $\Pr^S_T$ contains the unit. To see that it is closed under the tensor product, note that a fiberwise cocomplete $T$-$\infty$-category $\Ee$ is fiberwise presentable if and only if each $\Ee(A)$ is presentable; in particular, an $S$-cocomplete $T$-$\infty$-category is $S$-presentable if and only its underlying $S$-$\infty$-category is $S$-presentable. It will therefore suffice to show that for any $\Cc,\Dd\in\Pr^S_T$ the $S$-$\infty$-category $i^*(\Cc\otimes\Dd)$ is $S$-presentable. By the previous lemma, $i^*(\Cc\otimes\Dd)\simeq (i^*\Cc)\otimes (i^*\Dd)$, which is then indeed $S$-presentable by \cite{martiniwolf2022presentable}*{Proposition 2.6.2.9}.
    \end{proof}
\end{prop}

\subsubsection{Symmetric monoidal structures on equivariant and global spaces} As the unit of $\CAT_T^{\text{$S$-cc},\otimes}$, $\ul\Spc_{S\triangleright T}$ admits a unique symmetric monoidal structure with unit given by the terminal object. We will now explicitly describe this structure:

\begin{prop}\label{prop:spc-T-cart-pres-sym-mon}
    Let $S\subset T$ be a cleft.
    Then $\ul\Spc_{S\triangleright T}$ has fiberwise finite products, and the  inclusion $\ul{\Spc}_{S\triangleright T}\hookrightarrow \ul{\Spc}_T$ preserves fiberwise finite products. Moreover, the cartesian symmetric monoidal structure (Example~\ref{ex:cartesian-sym-mon}) on $\ul{\Spc}_{S\triangleright T}$ is $S$-presentably symmetric monoidal.

    \begin{proof}        
        It is clear that $\ul\Spc_T$ has fiberwise finite products. To see that the same holds more for generally for $\ul\Spc_{S\triangleright T}$, it will then suffice to show that the latter is closed under fiberwise finite products in $\ul\Spc_T$. For this we recall from \cite{CLL_Clefts}*{Lemma~3.16} that the functors $\smash{(i_!)_{/A}\colon\PSh(S)_{/A}\to\PSh(T)_{/A}}$ induced by left Kan extension along $i\colon S\hookrightarrow T$ assemble into an $S$-left adjoint $\smash{\ul\Spc_S\to i^*\ul\Spc_T}$. Clearly, each $\smash{(i_!)_{/A}}$ preserves the terminal object $\id_A$ (so that its image in particular is given by $\smash{\ul\Spc_{S\triangleright T}(A)}$), while \cite{CLL_Clefts}*{Lemma 3.14} shows that it preserves binary products. Thus, $\smash{\ul\Spc_{S\triangleright T}}$ has fiberwise finite products. 

        It only remains to prove that the tensor product preserves $S$-colimits in each variable separately. It is clear that each $\ul\Spc_{S\triangleright T}(A)\simeq\PSh(S)_{/A}$ is cartesian closed, so it remains to verify the projection formula $f_!(X\times f^*Y)\iso f_!X\times Y$ for all $f\colon A\to B$ in $S$, $X\in\ul\Spc_{S\triangleright T}(A)$, and $Y\in\ul\Spc_{S\triangleright T}(B)$. Since this only depends on the underlying $S$-$\infty$-category, we may assume without loss of generality that $T=S$, so that $\ul\Spc_{S\triangleright T}=\ul\Spc_S=\PSh(S)_{/{-}}$. Unravelling definitions, the Beck--Chevalley map is then simply given by the equivalence $X\times_A(A\times_BY)\iso X\times_BY$ from the pasting law for pullback squares.
    \end{proof}
\end{prop}

\cite{HA}*{Corollary 3.2.1.9} therefore specializes to the following:

\begin{cor}\label{cor:Spc-S-T-x-initial}
    The $S$-presentably symmetric monoidal $T$-$\infty$-category $\ul\Spc_{S\triangleright T}^\times$ is the initial object of $\CAlg(\CAT_T^{\textup{$S$-cc},\otimes})$.\qed
\end{cor}

As our last results in this section, we want to describe the cartesian symmetric monoidal structures on $\ul\myS_\gl$ and $\ul\myS$ in terms of the underlying model, which in light of the previous corollary will then in particular yield a corresponding description of the initial objects of $\CAlg(\CAT_\Glo^{\text{$\Glo$-cc},\otimes})$ and $\CAlg(\CAT_\Glo^{\text{$\Orb$-cc},\otimes})$, respectively.

\begin{constr}
    Applying Construction~\ref{constr:sym-mon-Borel} to the cartesian symmetric monoidal structures on $\cat{C}=\cat{Top}$ and localizing as before yields symmetric monoidal global $\infty$-categories $\ul\myS^\otimes$. Similarly, Corollary~\ref{cor:box-product-bifun} shows that the box product of orthogonal spaces induces a symmetric monoidal structure on $\ul\myS_\gl$.
\end{constr}

\begin{prop}
    Both $\ul\myS_\gl^\otimes$ and $\ul\myS^\otimes$ are cartesian.
    \begin{proof}
        It is clear from the definitions that the unit of $\ul\myS$ is terminal and that the derived functors $\myS_G\times\myS_G\to\myS_G$ of the product functors admit left adjoints (namely, the derived functor of the diagonal) for every $G$. Thus, $\ul\myS^\otimes$ is cartesian by Remark~\ref{rk:what-it-means-to-be-cartesian}. The statement for $\ul\myS_\gl^\otimes$ follows in the same way using Lemma~\ref{lemma:box-product-vs-product}.
    \end{proof}
\end{prop}

Combining Corollary~\ref{cor:Spc-S-T-x-initial} with Theorems~\ref{thm:unstable-main} and~\ref{thm:unstable-equiv-main} we therefore get:

\begin{cor}\label{cor:S-gl-times-initial}
    The globally presentably symmetric monoidal global $\infty$-category $\ul\myS_\gl^\otimes$ is the initial object of $\CAlg(\CAT_\Glo^{\textup{$\Glo$-cc},\otimes})$.\qed
\end{cor}

\begin{cor}\label{cor:S-equiv-times-initial}
    The equivariantly presentably symmetric monoidal global $\infty$-category $\ul\myS^\otimes$ is the initial object of $\CAlg(\CAT_\Glo^{\textup{$\Orb$-cc},\otimes})$.\qed
\end{cor}

\section{Basepoints}\label{sec:basepoints}
In this section we will establish universal properties for global $\infty$-categories of pointed global or equivariant spaces.

\subsection{Pointed \texorpdfstring{$\bm{T}$-$\bm{\infty}$}{T-∞}-categories} We begin by setting up the abstract theory.

\begin{definition}
    A $T$-$\infty$-category $\Cc$ is called \emph{pointed} if it factors through the non-full subcategory $\CAT_\infty^0\subset\CAT_\infty$ of $\infty$-categories with zero objects and functors preserving zero objects.
\end{definition}

Any $T$-$\infty$-category with a terminal object admits a universal map from a pointed $T$-$\infty$-category:

\begin{lemma}[See \cite{CLL_Global}*{Corollary~4.1.9}]
    Let $\Dd$ be a $T$-$\infty$-category with a terminal object. Then the forgetful map $\fgt\colon\Dd_*\to\Dd$ from the (levelwise) slice under the terminal object has the following universal property: for any pointed $T$-$\infty$-category $\Cc$ we have an equivalence
    \[
        \ul\Fun^0(\Cc,\fgt)\colon\ul\Fun^0(\Cc,\Dd_*)\iso\ul\Fun^0(\Cc,\Dd),
    \]
    where $\ul\Fun^0(\Cc,\Ee)\subset\ul\Fun(\Cc,\Ee)$ denotes the full subcategory spanned in degree $A\in T$ by those functors $\Cc\to\ul\Fun(\ul A,\Ee)$ preserving the terminal object.\qed
\end{lemma}

For $T$-cocomplete $\Cc$, there is also a universal left adjoint \emph{to} a pointed cocomplete $T$-$\infty$-category. We begin with the non-parametrized statement; while we expect this to be well-known, we weren't able to find a proof in this generality in the literature.

\begin{prop}\label{prop:univ-prop-*-cc}
    Let $\Cc$ be a cocomplete $\infty$-category with a terminal object. Then $\fgt\colon\Cc_*\to\Cc$ admits a left adjoint $(-)_+$, and this enjoys the following universal property: for any pointed cocomplete $\Dd$, restriction along $(-)_+$ defines an equivalence
    \[
        \Fun^\textup{cc}(\Cc_*,\Dd)\iso\Fun^\textup{cc}(\Cc,\Dd).
    \]
    \begin{proof}
        It is clear that the forgetful functor admits a left adjoint, given by taking coproduct with the terminal object. To complete the proof, it will then suffice to show that restriction along $(-)_+$ defines a bijection
        \begin{equation}\label{eq:restr-disj-bp}
            \mathop{\pi_0}\hom_{\CAT_\infty^\text{cc}}(\Cc_*,\Dd)\iso\mathop{\pi_0}\hom_{\CAT_\infty^\text{cc}}(\Cc,\Dd);
        \end{equation}
        the corresponding equivalence of functor categories will then follow by replacing $\Dd$ by $\Fun(\Kk,\Dd)$ for varying $\infty$-categories $\Kk$.

        To prove that $(\ref{eq:restr-disj-bp})$ is injective, let us define $\Sqcof(\Cc)\subset\Fun([1]\times[1],\Cc)$ as the full subcategory spanned by all pushout squares of the form
        \[
            \begin{tikzcd}
                A\arrow[r]\arrow[d]\arrow[dr,"\ulcorner"{very near end},phantom]& B\arrow[d]\\
                1 \arrow[r] & C
            \end{tikzcd}
        \]
        and similarly define $\Sqcof(\Dd)$. Note that the functors $\Sqcof(\Cc)\to\Cc$ and $\Sqcof(\Dd)\to\Dd$ restricting to the top row are equivalences. 

        If we view the counit $\epsilon$ of the adjunction $(-)_+\dashv\fgt$ as a functor $\Cc\to\Ar(\Cc)$, apply $\Ar$ on both sides, and restrict along $\Cc_*\hookrightarrow\Ar(\Cc)$, we obtain a functor $\Cc_*\to\Fun([1]\times[1],\Cc)$ sending $1\to X$ to the naturality square
        \begin{equation}\label{diag:nat-+-fgt}
            \begin{tikzcd}
                1_+\arrow[d,"\epsilon"']\arrow[r] & X_+\arrow[d,"\epsilon"]\\
                1\arrow[r] & X\rlap.
            \end{tikzcd}
        \end{equation}
        This functor factors through $\Sqcof(\Cc)$; postcomposing with $G$ then yields a functor $G'\colon\Cc_*\to\Sqcof(\Dd)$ such that ${\ev_{1,1}}\circ G'\simeq G$. On the other hand, the restriction of $G'$ to the top row only depends on $G\circ (-)_+$ by construction, so $(\ref{eq:restr-disj-bp})$ is injective.

        For surjectivity, let $F\colon\Cc\to\Dd$ arbitrary and define $G\colon\Cc_*\to\Dd$ as the composite
        \[
            \Cc_*=\Cc_{1/}\xrightarrow{\;F\;}\Dd_{F(1)/}\xrightarrow{\;p_!\;}\Dd_{0/}\iso\Dd,
        \]
        where $p_!$ is left adjoint to precomposition with $p\colon F(1)\to 0$, i.e.\ it is given by taking cofibers. This composite is then clearly cocontinuous. We now consider the diagram
        \[
            \begin{tikzcd}
                \Cc\arrow[r,"(-)_+"]\arrow[d,"F"'] & \Cc_{1/}\arrow[d,"F"]\\
                \Dd\arrow[r, "i_!"'] & \Dd_{F(1)/}\arrow[r,"p_!"'] & \Dd_{0/}
            \end{tikzcd}
        \]
        where the bottom left arrow is left adjoint to the forgetful functor. The left-hand square commutes as $F$ preserves pushouts (so that $\Ar(F)\colon\Ar(\Cc)\to\Ar(\Dd)$ preserves $\ev_0$-cocartesian edges), and the bottom composite is left adjoint, hence also inverse to the forgetful map $\fgt\colon\Dd_{0/}\iso\Dd$. We conclude that 
        \[
            G\circ(-)_+=\fgt^{-1}\circ p_!\circ F\circ(-)_+\simeq\fgt^{-1}\circ p_!\circ i_!\circ F\simeq F,
        \]
        so $G$ is the desired preimage.
    \end{proof}
\end{prop}

\begin{prop}\label{prop:univ-prop-*-cc-param}
    Let $S\subset T$ be a cleft, and let $\Cc$ be an $S$-cocomplete $T$-$\infty$-category with a terminal object. Then $\Cc_*$ is $S$-cocomplete and the forgetful functor $\fgt\colon\Cc_*\to\Cc$ admits a $T$-left adjoint $(-)_+$ enjoying the following universal property: for any pointed $S$-cocomplete $\Dd$, restriction along $(-)_+$ defines an equivalence
    \[
        \ul\Fun^\textup{$S$-cc}_T(\Cc_*,\Dd)\iso
        \ul\Fun^\textup{$S$-cc}_T(\Cc,\Dd).
    \]
    \begin{proof}
        It is clear that the forgetful functor admits a pointwise left adjoint (given by coproduct with the terminal object), and the Beck--Chevalley condition follows from fiberwise cocompleteness.

        To show that $\Cc_*$ is $S$-cocomplete, we apply Lemma~\ref{lemma:S-cocompleteness-repr}. First note that the composite
        \[
            \PSh(T)^\op\xrightarrow{\;\Cc\;}\CAT_\infty^1\xrightarrow{\;(-)_*\;}\CAT_\infty^*\subset\CAT_\infty,
        \]
        where $\CAT_\infty^1\subset\CAT_\infty$ denotes the subcategory of $\infty$-categories with a terminal object and functors preserving terminal objects, is limit-preserving (as $(-)_*$ is a right adjoint). Thus, this composite describes the limit extension of $\Cc_*$. Given any $Y\in\PSh(T)$ and $(f\colon X\to Y)\in\ul\Spc_{S\triangleright T}(Y)$, we may therefore identify $f^*\colon\Cc_*(Y)\to\Cc_*(X)$ with $f^*\colon\Cc(Y)_{1/}\to\Cc(X)_{1/}$ which has a left adjoint given by the composite
        \[
            \Cc(X)_{1/}\xrightarrow{\;f_!\;}\Cc(Y)_{f_!1/}\xrightarrow{\;\epsilon_!\;}\Cc(Y)_{1/}.
        \]
        To verify the Beck--Chevalley condition, consider a pullback
        \[
            \begin{tikzcd}
                A'\arrow[r, "t'"]\arrow[d,"s'"'] & B'\arrow[d,"s"]\\
                A\arrow[r,"t"'] & B
            \end{tikzcd}
        \]
        in $\PSh(T)$ with $s\in\ul\Spc_{S\triangleright T}(B)$. Then the Beck--Chevalley map $s'_!t'^*\to t^*s_!$ is a natural transformation of cocontinuous functors $\Cc(B')\to\Cc(A)$; we may therefore check invertibility after applying the equivalence from the previous proposition, i.e.\ after precomposing with $(-)_+\colon\Cc(B)\to\Cc(A)$. However, as Beck--Chevalley maps compose, this fits into a commutative diagram
        \[
            \begin{tikzcd}
                s'_!t^*(-)_+\arrow[r,"\sim"]\arrow[d,"\BC_!(-)_+"'] &[.5em] s'_!(-)_+t^*\arrow[r,"\BC_!t^*"] &[.5em] (-)_+s'_!t^*\arrow[d, "(-)_+\BC_!"]\\
                t^*s_!(-)_+\arrow[r,"t^*\BC_!"'] & t^*(-)_+s_!\arrow[r,"\sim"'] & (-)_+t^*s_!\rlap,
            \end{tikzcd}
        \]
        where the unlabelled equivalences come from $(-)_+$ being a $T$-functor. The two horizontal Beck--Chevalley maps are equivalences since $(-)_+$ is a left adjoint, as is the right-hand vertical map since $\Cc$ is $S$-cocomplete. Thus, also the left-hand vertical map is an equivalence by 2-out-of-3, finishing the proof that $\Cc_*$ is $S$-cocomplete.

        Finally, we will show that the restriction along $(-)_+$ defines an equivalence $\hom_{\CAT_T^\text{$S$-cc}}(\Cc_*,\Dd)\iso\hom_{\CAT_T^\text{$S$-cc}}(\Cc,\Dd)$; the corresponding claim for param\-etrized functor categories will then follow by replacing $\Dd$ by $\ul\Fun(\ul A\times[n],\Dd)$ for varying $A\in T$, $n\ge0$.

        Applying the previous proposition levelwise, we see that restriction defines an equivalence $\hom_{\CAT_T^\text{$\core T$-cc}}(\Cc_*,\Dd)\iso\hom_{\CAT_T^\text{$\core T$-cc}}(\Cc,\Dd)$; to complete the proof it will therefore suffice to show that a fiberwise cocontinuous $F\colon\Cc_*\to\Dd$ is $S$-cocontinuous if and only if the composite $\Cc\to\Cc_*\to\Dd$ is so. 
        This follows from Lemma~\ref{lm:continuouity-clefts}: Fix $B\in\PSh(T)$ and $(s\colon A\to B)\in\ul\Spc_{S\triangleright T}(B)$; arguing as before, we may check invertibility of the Beck--Chevalley map $\BC_!\colon s_!F\to Fs_!$ after precomposing with $(-)_+$, and we may identify $\BC_!(-)_+$ up to conjugation by equivalences with the Beck--Chevalley map for the composite $F\circ(-)_+\colon\Cc\to\Dd$. The claim follows immediately.
    \end{proof}
\end{prop}

\begin{cor}\label{cor:Spc-S-T-*-pres}
    The $T$-$\infty$-category $\ul\Spc_{S\triangleright T,*}$ is $S$-presentable, and for any $S$-cocomplete pointed $T$-$\infty$-category $\Cc$, evaluation at $1_+$ defines an equivalence
    \[
        \ul\Fun_T^\textup{$S$-cc}(\ul\Spc_{S\triangleright T,*},\Cc)\iso\Cc.
    \]
    \begin{proof}
        By the previous proposition, $\ul\Spc_{S\triangleright T,*}$ is $S$-cocomplete, and it is clear that each $\ul\Spc_{S\triangleright T,*}(A)\simeq\PSh(T_{/A})_{\id_A/}$ is presentable. Finally, combining the previous proposition with Corollary~\ref{cor:univ-prop-Spc-T} verifies the universal property.
    \end{proof}
\end{cor}

We can also give a different construction of the initial $S$-cocontinuous functor with pointed target using the Lurie tensor product:

\begin{prop}\label{prop:Spc*-levelwise-tensor}
    Let $\Cc$ be any $S$-cocomplete $T$-$\infty$-category. Then the levelwise tensor product ${\Spc_*}\otimes\Cc$ is again $S$-cocomplete, the map $\Cc\to{\Spc_*}\otimes\Cc$ induced by $(-)_+\colon\Spc\to\Spc_*$ is $S$-cocontinuous, and it moreover admits the following universal property: for every pointed $S$-cocomplete $T$-$\infty$-category $\Dd$, restriction defines an equivalence
    \[
        \ul\Fun^\textup{$S$-cc}_T({\Spc_*}\otimes\Cc,\Dd)\iso\ul\Fun^\textup{$S$-cc}_T(\Cc,\Dd).
    \]
    \begin{proof}
        In light of Proposition~\ref{prop:univ-prop-*-cc-param}, it will be enough to construct an equivalence ${\Spc_*}\otimes\Cc\iso\Cc_*$ under $\Cc$. For this, we consider the commutative diagram
        \[
            \begin{tikzcd}
                {\Spc}\otimes\Cc_{\hphantom{*}}\arrow[d]\arrow[r] & {\Spc_*}\otimes\Cc_{\hphantom{*}}\arrow[d]\\
                {\Spc}\otimes\Cc_*\arrow[r] & {\Spc_*}\otimes\Cc_*
            \end{tikzcd}
        \]
        with horizontal maps induced by $(-)_+\colon\Spc\to\Spc_*$ and vertical maps induced by $(-)_+\colon\Cc\to\Cc_*$. We claim that the right-hand vertical map and the bottom horizontal map are equivalences; this will then complete the proof of the corollary.

        To see that the bottom horizontal map is an equivalence, it suffices to show that ${\Spc}\otimes\Dd\to{\Spc_*}\otimes\Dd$ is invertible for every (non-parametrized) pointed cocomplete $\Dd$, or equivalently that the induced map between corepresented functors on $\CAT_\infty^\text{cc}$ is an equivalence. By the defining property of the tensor product, the latter amounts to saying that we have an equivalence
        \begin{equation}\label{eq:univ-prop-pted}
            {-}\circ(-)_+\colon\hom(\Spc_*,\Fun^\text{cc}(\Dd,\Ee))\iso\hom(\Spc,\Fun^\text{cc}(\Dd,\Ee))
        \end{equation}
        for any cocomplete $\Ee$. By \cite{CLL_Global}*{Corollary 4.1.10${}^\op$}, the full subcategory of $\Fun(\Dd,\Ee)$ spanned by the functors preserving the initial object is pointed, whence so is the smaller subcategory $\Fun^\text{cc}(\Dd,\Ee)$ as it still contains the initial object. Thus, applying the Proposition~\ref{prop:univ-prop-*-cc-param} for $S=T=1$ shows that $(\ref{eq:univ-prop-pted})$ is an equivalence.

        For the vertical map, we similarly compare corepresented functors on fiberwise cocomplete $\infty$-categories. Arguing as above, this reduces us to invertibility of
        \[
            \hom(\Cc_*,\Fun^\text{cc}(\Spc_*,-)\circ\Dd)\to
            \hom(\Cc,\Fun^\text{cc}(\Spc_*,-)\circ\Dd),
        \]
        and applying the above levelwise shows that $\Fun^\text{cc}(\Spc_*,-)\circ\Dd$ is pointed. The claim therefore follows from the special case $S=\core T$ of the proposition.
    \end{proof}
\end{prop}

We will now apply this to global $\infty$-categories of equivariant and global spaces.

\begin{constr}
    Localizing $\Ntop(\cat{$\cat{L}$-Top}_*^\dual)$ in each degree $G$ at the (underlying) $G$-global weak equivalences yields a global $\infty$-category $\ul\myS_{\gl,*}$. Analogously, we define $\ul\myS_{*}$ by localizing $\smash{\Ntop(\cat{Top}_*^\dual)}$ in each degree $G$ at the (underlying) $G$-equivariant weak equivalences.
\end{constr}

\begin{cor}\label{cor:gl-pointed}
    The global $\infty$-category $\ul\myS_{\gl,*}$ is pointed and globally presentable. For every pointed globally cocomplete $\Dd$, evaluation at $\consto(S^0)$ defines an equivalence
    \[
        \ul\Fun^\textup{$\Glo$-cc}_{\Glo}(\ul\myS_{\gl,*},\Dd)\iso\Dd.
    \]
    \begin{proof}
        It is clear that $\ul\myS_{\gl,*}$ has a zero object, so that there is a unique global functor $\ul\myS_{\gl,*}\to (\ul\myS_\gl)_*$ over $\ul\myS_\gl$, and applying Lemma~\ref{lemma:slice-over-fibrant}${}^\op$ levelwise shows that this map is an equivalence. Combining this with Theorem~\ref{thm:unstable-main}, we then get an equivalence $(\ul\Spc_\gl)_*\simeq\ul\myS_{\gl,*}$ sending $1_+$ to $\consto(S^0)$. The claim now follows from Proposition~\ref{prop:univ-prop-*-cc-param}.
    \end{proof}
\end{cor}

In the same way we deduce from Corollary~\ref{cor:equiv-spaces-universal}:

\begin{cor}\label{cor:equiv-pointed}
    The global $\infty$-category $\ul\myS_*$ is pointed and equivariantly presentable. For every pointed equivariantly cocomplete $\Dd$, evaluation at $S^0$ defines an equivalence
    \[
        \ul\Fun^\textup{$\Orb$-cc}_{\Glo}(\ul\myS_*,\Dd)\iso\Dd.\qednow
    \]
\end{cor}

\subsection{Model categorical properties} We devote this subsection to discussing some model categorical properties of based equivariant and global spaces, which will be used both in the next subsection to construct symmetric monoidal structures on $\ul\myS_*$ and $\ul\myS_{\gl,*}$ as well as in setting up the stable theory in Part II.

Let us begin by recalling a couple of generalities: for any model category $\Cc$ and any $A\in\Cc$, the slices $\Cc_{A/}$ and $\Cc_{/A}$ inherit model structures such that weak equivalences, fibrations, and cofibrations are created by the forgetful functor to $\Cc$. Moreover, $\Cc_{A/}$ and $\Cc_{/A}$ inherit properties like properness or cofibrant generation from $\Cc$ by \cite{hirschhorn-slice}. Similarly, if $\Cc$ is topological, then so are the slices $\Cc_{/A}$ and $\Cc_{A/}$ as the tensoring or cotensoring, respectively, are preserved by the forgetful functor. Let us spell this out for one of the model structures above, leaving the analogous results for other model categories implicit:

\begin{cor}
    The category $\cat{$\bm{G}$-$\cat{L}$-Top}_*$ of pointed orthogonal $G$-spaces admits a model structure in which a map is a weak equivalence, fibration, or cofibration if and only if it is so in $\cat{$\bm G$-$\cat{L}$-Top}$. This model structure is proper, topological, and cofibrantly generated.\qed
\end{cor}

\medskip
\subsubsection{Homotopy colimits} In §\ref{subsubsec:homotopy-po-equiv} we saw that pushouts along h-cofibrations compute homotopy pushouts in equivariant spaces. The analogous statement in the based setting is not quite true, even for the trivial group: for every based space $X$, the map $1\to X$ is an h-cofibration in $\cat{Top}_*$, but wedge sums of weak equivalences in $\cat{Top}_*$ are not always weak equivalences. We therefore introduce:

\begin{definition}
    Let $G$ be a compact Lie group. A map in $\cat{$\bm G$-Top}_*$ or $\cat{$\bm G$-$\cat{L}$-Top}$ is called an \emph{unbased h-cofibration} if it is a cofibration in $\cat{$\bm G$-Top}$ or $\cat{$\bm G$-$\cat{L}$-Top}$, respectively. We call an object $X$ \emph{well-based} if the unique map $1\to X$ is an unbased h-cofibration. 
\end{definition}

\begin{lemma}\label{lemma:h-cof-compute-po}
    Let $G$ be a compact Lie group, and let
    \begin{equation}\label{diag:pushout-h-cof}
        \begin{tikzcd}
            A\arrow[r,"i"]\arrow[d,"f"']\arrow[dr,phantom,"\ulcorner"{very near end}] & B\arrow[d]\\
            C\arrow[r] & D
        \end{tikzcd}
    \end{equation}
    be any pushout in $\cat{$\bm{G}$-$\cat{L}$-Top}_*$ such that $i$ is an unbased h-cofibration. Then $(\ref{diag:pushout-h-cof})$ is a homotopy pushout in the $G$-global model structure.

    The analogous statement holds for homotopy pushouts in the $\Ff$-model structure on $\cat{$\bm G$-Top}_*$ for any collection $\Ff$ of closed subgroups.
    \begin{proof}
        The second statement follows immediately from the corresponding unbased statement (Lemma~\ref{lemma:h-cof-compute-po-unbased}) as weak equivalences, cofibrations, and pushouts are created by the forgetful functor to $\cat{$\bm G$-Top}$.

        For the first statement, we then factor $f$ into a level cofibration $A\to C_0$ followed by a $G$-global level weak equivalence. Then $i(V)$ is an h-cofibration in $\cat{$\bm{(\O(V)\times G)}$-Top}$ by \cite{schwede2018global}*{Corollary~A.30(ii)}, so applying the second claim levelwise shows that the induced map $B\amalg_AC_0\to B\amalg_AC=D$ is even a $G$-global level weak equivalence, so that $(\ref{diag:pushout-h-cof})$ is a homotopy pushout as claimed.
    \end{proof}
\end{lemma}

\begin{lemma}\label{lemma:coprod-well-based}
    Let $G$ and $\Ff$ be as before. Let $(X_i)_{i\in I}$ be a small collection of based $G$-spaces, such that for every $i\in I$ the object $X_i$ is well-based. Then the maps $X_j\to\bigvee_{i\in I}X_i$ exhibit $\bigvee_{i\in I}X_i$ as a coproduct in the underlying $\infty$-category. 

    Moreover, the analogous statement for orthogonal $G$-spaces holds.
    \begin{proof}
        We will prove the equivariant statement, the proof of the global one being analogous. 
        
        Pick for every $i\in I$ a cofibrant replacement $X'_i\to X_i$. The claim amounts to saying that the induced map $\bigvee_{i\in I}X_i'\to\bigvee_{i\in I}X_i$ is an $\Ff$-weak equivalence, which can be checked after forgetting the basepoint. The map in question is then the induced map on pushouts
        \[
            \begin{tikzcd}
                1\arrow[d,"="'] & \arrow[l] \coprod_{i\in I}1\arrow[d,"="{description}]\arrow[r] & \coprod_{i\in I} X_i'\arrow[d]\\
                1 &\arrow[l]\coprod_{i\in I} 1\arrow[r]& \coprod_{i\in I} X_i\rlap.
            \end{tikzcd}
        \]       
        By assumption the bottom right horizontal map is a coproduct of h-cofibrations, hence itself an h-cofibration; thus, the previous lemma shows that the pushout of the bottom row computes the homotopy pushout. Similarly, the pushout of the top row computes the homotopy pushout as the top right map is a cofibration. The claim follows as $\coprod_{i\in I}X_i'\to\coprod_{i\in I}X_i$ is an $\Ff$-weak equivalence by direct inspection.
    \end{proof}
\end{lemma}

\subsubsection{Change of group} Next, we come to our usual discussion of functoriality as the group $G$ varies:

\begin{lemma}\label{lemma:res-pted-Quillen}
    Let $\alpha\colon G\to G'$ be a homomorphism of compact Lie groups. Then
    \[
        \alpha_!\colon\cat{$\bm G$-$\cat{L}$-Top}_*\rightleftarrows\cat{$\bm{G'}$-$\cat{L}$-Top}_*\noloc\alpha^*
    \]
    is a Quillen adjunction with homotopical right adjoint. If $\alpha$ is injective, then also
    \[
        \alpha^*\colon\cat{$\bm{G'}$-$\cat{L}$-Top}_*\rightleftarrows\cat{$\bm{G}$-$\cat{L}$-Top}_*\noloc \alpha_*
    \]
    is a Quillen adjunction.
    \begin{proof}
        As $\alpha^*$ commutes with the functors forgetting the basepoint, and since weak equivalences, cofibrations, and fibrations are created by the forgetful functor, this is immediate from the corresponding unbased statements recorded as Lemmas~\ref{lemma:restr-right-Quillen} and~\ref{lemma:restr-inj-left-Quillen}.
    \end{proof}
\end{lemma}

Similarly, Lemmas~\ref{lemma:restr-equiv-rQ} and~\ref{lemma:restr-equiv-lQ} have based analogues; again, we will not make them explicit, but instead simply refer to the unbased statements whenever we need the based versions.

\medskip
\subsubsection{Multiplicative properties} Finally, let us discuss how the smash product of based spaces interacts with equivariant model structures.

\begin{lemma}\label{lemma:smash-lQ}
    Let $G$ be a compact Lie group, and let $\Ff_1,\Ff_2$ be families of closed subgroups. Then
    \begin{equation}\label{eq:smash-equiv}
        {-}\smashp{-}\colon (\cat{$\bm G$-Top}_*)_{\Ff_1}\times(\cat{$\bm G$-Top}_*)_{\Ff_2}\to(\cat{$\bm G$-Top}_*)_{\Ff_1\cap\Ff_2}.
    \end{equation}
    is a left Quillen bifunctor for the model structures as indiciated.
    \begin{proof}
        Write $\Ff_1\times\Ff_2$ for the family of closed subgroups $K\subset G\times G$ such that there exist $H_1\in\Ff_1,H_2\in\Ff_2$ with $K\subset H_1\times H_2$. Then $(\ref{eq:smash-equiv})$ factors as
        \[
            (\cat{$\bm G$-Top}_*)_{\Ff_1}\times(\cat{$\bm G$-Top}_*)_{\Ff_2}
            \xrightarrow{\;{-}\smashp{-}\;}\big(\cat{$\bm{(G\times G)}$-Top}_*\big){}_{\Ff_1\times\Ff_2}
            \xrightarrow{\;\Delta^*\;}(\cat{$\bm G$-Top}_*)_{\Ff_1\cap\Ff_2},
        \]
        where the second functor is left Quillen by Lemma~\ref{lemma:restr-equiv-lQ}. We are therefore reduced to showing that the first functor is a left Quillen bifunctor.

        It suffices to verify the pushout product axiom for the generating (acyclic) cofibrations, which are of the form $(G_i/H_i\times f_i)_+$ for generating (acyclic) cofibrations $f$ of $\cat{Top}_*$ and $H_i\in\Ff_i$. The pushout product is then given by
        \[
            (G_1/H_1\times f_1)_+\ppo
            (G_2/H_2\times f_2)_+\cong
            \big((G_1\times G_2)/(H_1\times H_2) \times (f_1\ppo f_2)\big){}_+.
        \]
        As $\cat{Top}$ is a cartesian model category, $f_1\ppo f_2$ is a cofibration, acyclic if at least one of $f_1,f_2$ is so. As $(G_1\times G_2)/(H_1\times H_2)$ is $(\Ff_1\times\Ff_2)$-cofibrant by definition and since the $(\Ff_1\times\Ff_2)$-model structure is topological, the claim follows.
    \end{proof}
\end{lemma}

\begin{prop}\label{prop:smash-well-based}
    Let $G$ be a compact Lie group and let $\Ff_1,\Ff_2$ be two families of closed subgroups of $G$. If $X$ is $\Ff_1$-cofibrant, then $X\smashp-$ sends $(\Ff_1\cap\Ff_2)$-weak equivalences between well-based $G$-spaces to $\Ff_2$-weak equivalences.
    \begin{proof}
       Let $f\colon Y\to Y'$ be an $(\Ff_1\cap\Ff_2)$-weak equivalence of well-based $G$-spaces. We note that $(X\smashp f)^H$ agrees up to homeomorphism with $X^H\smashp f^H$ for any $H\subset G$. If $H\in\Ff_2\smallsetminus\Ff_1$, then $X^H=0$, so that $(X\smashp f)^H$ is a homeomorphism; on the other hand, if $H\in\Ff_1\cap\Ff_2$, then $X^H\smashp f^H$ is a weak equivalence as $f^H$ is a weak equivalence of well-based spaces and $X^H$ is well-based.
    \end{proof}
\end{prop}

The following seemingly random result will later turn out to be crucial in setting up the model structures for $G$-global stable homotopy theory:

\begin{cor}\label{cor:balanced-smash-product-oddly-specific}
    Let $G,H$ be compact Lie groups, and let $N\subset H$ be a closed normal subgroup. We write $\Ff_1$ for the family of closed subgroups of $G\times H$ that intersect $1\times N$ trivially and $\Ff_2$ for the family of graph subgroups for homomorphisms $H'\to G$ whose kernel is contained in $N$.
    
    If $X$ is cofibrant in the $\Ff_1$-model structure, then $X\smashp_H-$ sends $(\mathcal G_{H,G}\cap\mathcal G_{G,H})$-weak equivalences between $\Ff_2$-cofibrant based $(G\times H)$-spaces to $G$-equivariant weak equivalences.
\end{cor}

As the most important special case, for $N=H$ this says says that for any $\mathcal G_{H,G}$-cofibrant $X$, the functor $X\smashp_H-$ sends $(\mathcal G_{H,G}\cap\mathcal G_{G,H})$-weak equivalences between $\mathcal G_{G,H}$-cofibrant objects to $G$-weak equivalences. The extra flexibility afforded by the choice of normal subgroup $N$ will be useful for a reduction argument later.

\begin{proof}
    The family $\mathcal G_{H,G}\cap\mathcal G_{G,H}$ consists precisely of those groups that intersect $1\times H$ and $G\times 1$ trivially. If $H'\subset H$ and $\phi\colon H'\to G$ is a homomorphism, then $\Gamma_{H',\phi}\cap(G\times1)=1$; moreover, if $\ker(\phi)\subset N$, then $\Gamma_{H',\phi}\cap (1\times H)\subset 1\times N$, whence $\Gamma_{H',\phi}\cap (1\times H)=\Gamma_{H',\phi}\cap (1\times N)$. We conclude that $\Ff_1\cap\Ff_2\subset \mathcal G_{H,G}\cap\mathcal G_{G,H}$, so that any $(\mathcal G_{H,G}\cap\mathcal G_{G,H})$-weak equivalence is also an $(\Ff_1\cap\Ff_2)$-weak equivalence (in fact, it is not hard to show by similar arguments that $\Ff_1\cap\Ff_2=\mathcal G_{H,G}\cap\mathcal G_{G,H}$). Thus, if $X$ is $\Ff_1$-cofibrant and $f\colon Y\to Y'$ is a $(\mathcal G_{H,G}\cap\mathcal G_{G,H})$-weak equivalence of $\Ff_2$-cofibrant objects, then the previous proposition in particular shows that $X\smashp f$ is an $\Ff_2$-weak equivalence between $\Ff_2$-cofibrant based $(G\times H)$-spaces. Since $(-)/H\colon (\cat{$\bm{(G\times H)}$-Top}_*)_{\Ff_2}\to(\cat{$\bm G$-Top}_*)_{\All}$ is left Quillen by Lemma~\ref{lemma:restr-equiv-rQ}, Ken Brown's Lemma then shows that $X\smashp_Hf=(X\smashp f)/H$ is a $G$-weak equivalence as claimed.
\end{proof}

Now we turn to the analogue of the smash product in the global situation. 

\begin{constr}
    For every $X,Y\in\cat{$\cat{L}$-Top}$ we define $X\boxsmash Y$ as the cofiber of the pushout-product map $X\vee Y\cong(X\boxtimes 1)\amalg_{1\boxtimes 1}(1\boxtimes Y)\to (X\boxtimes Y)$. We omit the straight-forward but lengthy verification that this becomes an enriched functor preserving colimits in each variable via the universal property of the cofiber, and that the symmetry, unitality, and associativity isomorphisms for the box product on $\boxtimes$ then similarly induce isomorphisms exhibiting $\boxsmash$ as the symmetric monoidal product of a closed topologically enriched symmetric monoidal structure on $\cat{$\cat{L}$-Top}_*$ with unit $\consto(S^0)$. Moreover, the left adjoint $(-)_+\colon\cat{$\cat{L}$-Top}\to\cat{$\cat{L}$-Top}_*$ of the forgetful functor acquires the structure of a strong symmetric monoidal functor with unit isomorphism the unique one and with structure isomorphisms $(X_+)\boxsmash(Y_+)\to (X\boxtimes Y)_+$ inverse to the map induced by $X\boxtimes Y\to (X_+)\boxtimes(Y_+)\to (X_+)\boxsmash(Y_+)$.
\end{constr}

\begin{lemma}\label{lemma:boxsmash}
    For any compact Lie groups $G$ and $H$, 
    \[
        {-}\boxsmash{-}\colon\cat{$\bm G$-$\cat{L}$-Top}_*\times \cat{$\bm H$-$\cat{L}$-Top}_*\to\cat{$\bm{(G\times H)}$-$\cat{L}$-Top}_*
    \]
    is a left Quillen bifunctor.
\end{lemma}

Arguing as before, we then in particular see that it also defines a left Quillen bifunctor $\cat{$\bm G$-$\cat{L}$-Top}_*\times\cat{$\bm G$-$\cat{L}$-Top}_*\to\cat{$\bm G$-$\cat{L}$-Top}_*$ for any compact Lie group $G$.

\begin{proof}
    The pushout product axiom can be checked on generating (acyclic) cofibration. As these can be obtained from the generating (acyclic) cofibrations of the corresponding unbased model structure by adding disjoint basepoints, this then follows from the pushout-product axiom for the box product (Lemma~\ref{lemma:box-product-bifun-general}).
\end{proof}

\begin{definition}\label{defi:qcof}
    A (based) $G$-global space $X\in\cat{$\bm G$-$\cat{L}$-Top}_*$ is called \emph{quasi-cofibrant} if it is a restriction of a cofibrant (based) orthogonal $G'$-space along some $\alpha\colon G\to G'$. We write $\cat{$\bm G$-$\cat{L}$-Top}^\text{qcof}$ and $\cat{$\bm G$-$\cat{L}$-Top}^\text{qcof}_*$ for the full subcategory of the quasi-cofibrant (based) orthogonal $G$-spaces.
\end{definition}

\begin{prop}\label{prop:boxsmash-qcof}
    The functor $\boxsmash$ restricts to a homotopical functor
    \[
        \cat{$\bm G$-$\cat{L}$-Top}^\textup{qcof}_*
        \times
        \cat{$\bm G$-$\cat{L}$-Top}^\textup{qcof}_*\to
        \cat{$\bm G$-$\cat{L}$-Top}^\textup{qcof}_*
    \]
    \begin{proof}
        If $X$ and $Y$ are quasi-cofibrant, then we find $\alpha\colon G\to G_1$ and $\beta\colon G\to G_2$ together with cofibrant $X'\in\cat{$\bm{G_1}$-\cat{L}-Top}_*$ and $Y'\in\cat{$\bm{G_2}$-\cat{L}-Top}_*$ such that $X=\alpha^*X'$ and $Y=\beta^*Y'$. Then the based orthogonal $G$-space $X\boxsmash Y$ is the restriction of the based orthogonal $(G_1\times G_2)$-space $X'\boxsmash Y'$ along $(\alpha,\beta)\colon G\to G_1\times G_2$. As $X'\boxsmash Y'$ is cofibrant by Lemma~\ref{lemma:boxsmash}, we conclude that $X\boxsmash Y$ is quasi-cofibrant.

        Similarly, we deduce from Lemma~\ref{lemma:box-product-bifun-general} that for $X$ and $Y$ as above the top horizontal map in the pushout square
        \[
            \begin{tikzcd}
                X\vee Y\arrow[d]\arrow[r]\arrow[dr,"\ulcorner"{very near end},phantom] & X\boxtimes Y\arrow[d]\\
                1\arrow[r] & X\boxsmash Y
            \end{tikzcd}
        \]
        defining $X\boxsmash Y$ is a restriction of a $(G_1\times G_2)$-cofibration along $(\alpha,\beta)$, and hence in particular an h-cofibration. Thus, Lemma~\ref{lemma:h-cof-compute-po} shows that this square is also a homotopy pushout. Given now $G$-global weak equivalences $f\colon X_1\to X_2$ and $g\colon Y_1\to Y_2$ of quasi-cofibrant based $G$-spaces, the map $f\boxtimes g$ is a $G$-global weak equivalence by Corollary~\ref{cor:box-product-bifun}, as is $f\vee g$ by Lemma~\ref{lemma:coprod-well-based}. We conclude that also the induced map on homotopy pushouts $f\boxsmash g$ is a $G$-global weak equivalence.
    \end{proof}
\end{prop}

Finally, let us discuss the model-categorical properties of the (co)tensoring of $\cat{$\bm{G}$-$\cat{L}$-Top}_*$ over $\cat{$\bm G$-Top}_*$:

\begin{lemma}\label{lemma:G-Top*-tensored}
    The levelwise smash product $\cat{$\bm G$-Top}_*\times\cat{$\bm G$-$\cat{L}$-Top}_*\to \cat{$\bm G$-$\cat{L}$-Top}_*$ is a left Quillen bifunctor for the $\All$-model structure on $\cat{$\bm G$-Top}_*$ and the $G$-global model structure on $\cat{$\bm G$-$\cat{L}$-Top}_*$.
    \begin{proof}
        Adjoining over, we may equivalently show that the cotensoring over $\cat{$\bm G$-Top}$ is a right Quillen bifunctor, which can be checked after forgetting the basepoint. Adjoining back, it therefore suffices to prove that the cartesian product $\cat{$\bm G$-Top}\times\cat{$\bm G$-$\cat{L}$-Top}\to \cat{$\bm G$-$\cat{L}$-Top}$ is a left Quillen bifunctor. 
        
        For the pushout-product axiom for cofibrations, it then suffices to treat the case of the generating cofibrations, which by the same reasoning as before reduces to showing that $X\coloneqq G/K\times\cat{L}(V,-)\times_\phi G$ is a cofibrant orthogonal $G$-space for every closed subgroup $K\subset G$, any homomorphism $\phi\colon H\to G$, and any faithful $H$-representation $V$. However, if we write $i\colon K\to G$ for the inclusion, then $X$ is isomorphic to $i_!i^*(\cat{L}(V,-)\times_\phi G)$, so the claim follows as both $i_!$ and $i^*$ are left Quillen by Lemmas~\ref{lemma:restr-inj-left-Quillen} and~\ref{lemma:restr-right-Quillen}. With this established, the pushout product for acyclic cofibrations is then an instance of \cite{barrero2021}*{Corollary~A.10}.
    \end{proof}
\end{lemma}

\begin{cor}\label{cor:global-smash-well-based}
    If $X\in\cat{$\bm G$-Top}_*$ is $\All$-cofibrant, then $X\smashp{-}\colon\cat{$\bm G$-$\cat{L}$-Top}_*\to \cat{$\bm G$-$\cat{L}$-Top}_*$ preserves $G$-global weak equivalences between well-based orthogonal $G$-spaces.
    \begin{proof}
        By the previous lemma, $X\smashp{-}$ preserves $G$-global weak equivalences between \emph{cofibrant} objects. Picking functorial cofibrant replacements in the $G$-global level model structure, it will therefore suffice to show that $X\smashp{-}$ preserves $G$-global \emph{level} weak equivalence between well-based objects. This follows at once by applying Proposition~\ref{prop:smash-well-based} levelwise.
    \end{proof}
\end{cor}

\begin{lemma}\label{lemma:cotensoring-fully-homotopical}
    Let $X$ be a finite based $G$-CW-complex. Then the cotensoring 
    \begin{equation}\label{eq:L-Top*-cotensoring}
        \maps_*(X,-)\colon\cat{$\bm G$-$\cat{L}$-Top}_*\to\cat{$\bm G$-$\cat{L}$-Top}_*
    \end{equation}
    preserves $G$-global weak equivalences. 
    \begin{proof}
        Applying Lemma~\ref{lemma:smash-lQ} levelwise shows that $(\ref{eq:L-Top*-cotensoring})$ is right Quillen for the $G$-global \emph{level} model structure, hence it preserves $G$-global level weak equivalences by Ken Brown's Lemma. Taking functorial cofibrant replacements in the $G$-global level model structure, we are therefore reduced to proving that $(\ref{eq:L-Top*-cotensoring})$ preserves $G$-global weak equivalences between cofibrant pointed orthogonal $G$-spaces. Note that any cofibrant pointed orthogonal $G$-space $X$ is in particular closed: by Lemma~\ref{lemma:restr-inj-left-Quillen}, $1\to X$ is a cofibration in $\cat{$\cat{L}$-Top}$, and $1=\cat{L}(0,-)$ is cofibrant, so that $X$ is cofibrant in $\cat{$\cat{L}$-Top}$ and hence closed by Lemma~\ref{lemma:cofibrant-closed}. Thus, it suffices to prove that $\maps_*(X,-)$ preserves $G$-global weak equivalences between closed pointed orthogonal $G$-spaces.

        For this, we fix a compact Lie group $H$ and a complete $H$-universe $\Uu_H$. Then 
        \begin{equation}\label{eq:cotensoring-graph}
            \maps_*(X,-)\colon\cat{$\bm{(G\times H)}$-Top}_*\to\cat{$\bm{(G\times H)}$-Top}_*
        \end{equation}
        preserves sequential colimits along closed embeddings, so choosing a cofinal sequence in the poset of finite-dimensional subrepresentations of $\Uu_H$, we see that $\maps_*(X,Y(\Uu_G))\cong\maps_*(X,Y)(\Uu_G)$ natural in $Y\in\cat{$\bm G$-$\cat{L}$-Top}_*$. It therefore suffices that $(\ref{eq:cotensoring-graph})$ preserves $\mathcal G_{H,G}$-weak equivalences, which again follows from Lemma~\ref{lemma:smash-lQ} and Ken Brown's Lemma.
    \end{proof}
\end{lemma}

\subsection{Symmetric monoidal universal properties}\label{subsec:pointed-sym-mon} In this subsection, we will explain how one can use the smash products studied above to upgrade $\ul\myS_{\gl,*}$ and $\ul\myS_*$ to symmetric monoidal global $\infty$-categories, and we will moreover prove that these satisfy pointed analogues of the universal properties established in §\ref{subsec:sym-mon}. 

\medskip
\subsubsection{Idempotent algebras and smashing localizations} We begin by more generally studying symmetric monoidal structures on $\ul\Spc_{S\triangleright T,*}$ for any cleft $S\subset T$. Our goal will be to show that there is in fact a \emph{unique} $S$-presentably $T$-symmetric monoidal structure with unit $1_+$, and that this yields the initial \emph{pointed} symmetric monoidal $T$-$\infty$-category such that the tensor product is $S$-bilinear. The most convenient way to prove this result is via the framework of \emph{idempotent algebras} and \emph{smashing localizations}, and we begin by recalling the relevant general theory.

\begin{definition}
    Let $\Cc$ be a symmetric monoidal $\infty$-category with unit $\bbone$. An \emph{idempotent algebra} consists of an object $A\in\Cc$ together with a map $a\colon\bbone\to X$ such that the map $a\otimes A$ is an equivalence. 
\end{definition}

\begin{definition}
    Let $\Cc$ be a symmetric monoidal $\infty$-category, and let $a\colon\bbone\to A$ be idempotent. An object $B\in\Cc$ is called an \emph{$A$-module} if the map $a\otimes B\colon \bbone\otimes B\to A\otimes B$ is an equivalence.
\end{definition}

\begin{rk}\label{rk:idempotent-algs-are-algs}
    Given any idempotent algebra $a\colon \bbone\to A$, \cite{HA}*{Proposition 4.8.2.9} shows that $A$ admits a unique commutative algebra structure with unit given by $a$. Similarly, {Proposition 4.8.2.10} of \emph{loc.\ cit.}\ shows that the forgetful functor $\text{LMod}_A(\Cc)\to\Cc$ is fully faithful with essential image given by the $A$-modules in the sense of the previous definition, justifying the terminology.
\end{rk}

As a direct consequence of the preceeding remark we have, see also~\cite{HA}*{Proposition 4.8.2.4}:

\begin{cor}\label{cor:idempotent-la}
    Let $a\colon\bbone\to A$ be idempotent. Then the inclusion of the full subcategory spanned by the $A$-modules admits a left adjoint given by $A\otimes-$, with unit induced by $a$.\qed
\end{cor}

\begin{definition}
    We call a Bousfield localization $L\colon\Cc\rightleftarrows\Dd\noloc i$ a \emph{smashing localization} if there exists an idempotent algebra $(A,a)$ such that the essential image of the fully faithful right adjoint $i$ consists precisely of the $A$-modules.
\end{definition}

\begin{rk}
    By Corollary~\ref{cor:idempotent-la}, $L\simeq A\otimes{-}$ and $a\colon\bbone\to A$ can be recovered as the unit $\bbone\to iL(\bbone)$. In particular, the idempotent algebra $(A,a)$ is unique.
\end{rk}

Finally, we record:

\begin{lemma}[See \cite{CSY2021AmbiHeight}*{Proposition 5.1.7}]\label{lemma:unique-algebra-structure-initial}
    Let $(A,a)$ be an idempotent algebra, viewed as an object of $\CAlg(\Cc)$. If $B\in\CAlg(\Cc)$ is arbitrary, then $\hom_{\CAlg(\Cc)}(A,B)$ is empty or contractible, and there exists a map $A\to B$ of commutative algebras if and only if (the underlying object of) $B$ is an $A$-module. In particular, $A$ is the initial commutative algebra that is also an $A$-module.\qed
\end{lemma}

\subsubsection{Smashing localizations of $\Pr^S_T$}
We will now specialize the above discussion to the case where $\Cc=\CAT_T^\text{$S$-cc}$ or $\Cc=\Pr^S_T$ for a cleft $S\subset T$, both equipped with the parametrized Lurie tensor product from §\ref{subsubsec:param-Lurie}.

\begin{convention}
    By the universal property of the unit $\ul\Spc_{S\triangleright T}$, a morphism $\ul\Spc_{S\triangleright T}\to \Cc$ in one of these categories is described by the image $X$ of the terminal object; thus, we will typically denote idempotent algebras by $(\Cc,X)$. 
\end{convention}

The following result will often allow us to reduce to the case $T=S$:

\begin{lemma}\label{lemma:idempotent-restrict}
    Let $i\colon S\hookrightarrow T$ be a cleft, let $\Cc\in\CAT_T^\textup{$S$-cc}$, and let $X$ be a global section. Then $(\Cc,X)$ is an idempotent in $\CAT_T^\textup{$S$-cc}$ if and only if the restriction $(i^*\Cc,i^*X)$ is an idempotent in $\CAT^\textup{$S$-cc}_S$. Moreover, in this case a $\Dd\in\CAT_T^\textup{$S$-cc}$ is a $\Cc$-module if and only if $i^*\Dd$ is an $i^*\Cc$-module.
    \begin{proof}
        This follows at once from Lemma~\ref{lemma:restr-cleft-sym-mon} since the restriction functor $i^*\colon\CAT_T\to\CAT_S$ is conservative.
    \end{proof}
\end{lemma}

\begin{lemma}\label{lemma:smash-module-internal-hom}
    Let $(\Cc,X)$ be an idempotent in $\CAT_T^\textup{$S$-cc}$. Then $\Dd\in\CAT_T^\textup{$S$-cc}$ is a $\Cc$-module if and only if evaluation at $X$ defines an equivalence
    \begin{equation}\label{eq:cofree-counit}
        \ev_X\colon\ul\Fun^\textup{$S$-cc}_T(\Cc,\Dd)\iso\Dd.
    \end{equation}
    \begin{proof}
        By Remark~\ref{rk:S-cc-closed}, the symmetric monoidal structure on $\CAT_T^\textup{$S$-cc}$ afforded by the parametrized Lurie tensor product is closed, with internal homs given by $\ul\Fun^\textup{$S$-cc}$.
        Thus, the fully faithful functor $\fgt\colon\text{LMod}_\Cc(\CAT_T^\textup{$S$-cc})\to \CAT_T^\textup{$S$-cc}$ admits a right adjoint such that the counit is given for any $\Dd$ by the map $(\ref{eq:cofree-counit})$.
        The claim now follows since the essential image of any fully faithful left adjoint consists precisely of those objects for which the counit is invertible. 
    \end{proof}
\end{lemma}

Let us now specialize the above to the setting of pointed $T$-$\infty$-categories:

\begin{prop}
    The pair $\smash{(\ul\Spc_{S\triangleright T,*},1_+)}$ defines an idempotent algebra in $\Pr^S_T$, and hence also in $\CAT_T^\textup{$S$-cc}$. An $S$-cocomplete $T$-$\infty$-category $\Cc$ is a $\ul\Spc_{S\triangleright T,*}$-module if and only if it is pointed.
    \begin{proof}
        By Corollary~\ref{cor:Spc-S-T-*-pres}, $\ul\Spc_{S\triangleright T,*}$ is $S$-presentable. To see that it is idempotent and to characterize its modules, recall that $\smash{\ul\Spc_{S\triangleright T,*}}\simeq{{\Spc_*}\otimes\ul\Spc_{S\triangleright T}}$ by Proposition~\ref{prop:Spc*-levelwise-tensor}. By \cite{martiniwolf2022presentable}*{Proposition~2.6.3.7}, the functor ${-}\otimes\ul\Spc_S\colon\PrL\to\PrL_S$ can be upgraded to a strong symmetric monoidal functor; as $(\Spc_*,S^0)$ is an idempotent algebra \cite{HA}*{Proposition 4.8.2.11}, so is its image $(\ul\Spc_{S,*},S^0)$. Lemma~\ref{lemma:idempotent-restrict} then implies that also $\smash{(\ul\Spc_{S\triangleright T,*},S^0)}$ is idempotent.

        Combining Lemma~\ref{lemma:smash-module-internal-hom} with Proposition~\ref{prop:univ-prop-*-cc-param} already shows that every pointed $S$-cocomplete $T$-$\infty$-category is a $\ul\Spc_{S\triangleright T,*}$-module. For the converse, it will suffice to show that $\ul\Fun^\text{$S$-cc}_T(\ul\Spc_{S\triangleright T,*},\Dd)$ is pointed for any $S$-cocomplete $\Dd$. By \cite{CLL_Global}*{Corollary 4.1.10${}^\op$}, the full subcategory of $\ul\Fun_T(\ul\Spc_{S\triangleright T,*}^\op,\Dd)$ spanned by the functors preserving the initial object is pointed; as $\ul\Fun^\text{$S$-cc}_T(\ul\Spc_{S\triangleright T,*},\Dd)$ is a smaller subcategory still containing the initial object, it is then necessarily also pointed.
    \end{proof}
\end{prop}

Lemma~\ref{lemma:unique-algebra-structure-initial} therefore specializes to the following:

\begin{cor}\label{cor:smash-initial-general}
    There exists a unique $S$-presentably symmetric monoidal structure on $\ul\Spc_{S\triangleright T,*}$ with unit $1_+$. Moreover, this defines the initial object of ${\CAlg(\CAT_T^{\textup{$S$-cc},0})}$.\qed
\end{cor}

\subsubsection{Symmetric monoidal structures for global and equivariant based spaces}
We will now describe the symmetric monoidal structure from the corollary in model categorical terms for $T=\Glo$ and $S\in\{\Orb,\Glo\}$.

\begin{constr}\label{constr:smash-product-unstable}
    We write $\Ntop(\cat{$\cat{L}$-Top}_*^{\dual,\text{qcof}})\subset\Ntop(\cat{$\cat{L}$-Top}_*^{\dual})$ for the full subcategory given in degree $G$ by the quasi-cofibrant orthogonal $G$-spaces (Definition~\ref{defi:qcof}). By Proposition~\ref{prop:boxsmash-qcof}, the box smash product $\boxsmash$ then restricts to a symmetric monoidal structure on this subcategory and it is moreover homotopical. Thus, we may derive to a yield a symmetric monoidal global $\infty$-category $\ul\myS_{\gl,*}^{\text{qcof},\otimes}$. 
    
    Similarly, we define $\smash{\ul\myS_*^{\text{cof},\otimes}}$ by localizing the full subcategory of $\smash{\Ntop(\cat{Top}^\dual)}$ given in each degree by the $\All$-cofibrant spaces, equipped with the smash product.
\end{constr}

\begin{lemma}\label{lemma:(q)cof-resolution}
    The inclusions induce equivalences of global $\infty$-categories 
    \[
        \ul\myS_{\gl,*}^\textup{qcof}\iso\ul\myS_{\gl,*}
        \qquad\text{and}\qquad
        \ul\myS_*^\textup{cof}\iso\ul\myS_*.
    \]
    \begin{proof}
        If we fix a compact Lie group $G$, then an inverse is induced by cofibrant replacement in the $G$-global and $G$-equivariant model structure, respectively.
    \end{proof}
\end{lemma}

As an upshot, we can uniquely extend the symmetric monoidal structures from Construction~\ref{constr:smash-product-unstable}, yielding symmetric monoidal global $\infty$-categories $\ul\myS_{\gl,*}^\otimes$ and $\ul\myS_{*}^\otimes$.

\begin{thm}\label{thm:gl-smash-initial}
    The pointed symmetric monoidal global $\infty$-category $\ul\myS_{\gl,*}^\otimes$ is globally presentably symmetric monoidal, and it is initial in ${\CAlg(\CAT_{\smash{\Glo}}^{\textup{$\Glo$-cc},0})}$.
    \begin{proof}
        In light of Corollaries~\ref{cor:gl-pointed} and~\ref{cor:smash-initial-general} it only remains to show that this is a globally presentably symmetric monoidal $\infty$-category with unit $\consto(S^0)$.

        The claim about the unit is clear; moreover, we already know that the underlying global $\infty$-category $\ul\myS_{\gl,*}\simeq\ul\myS_{\gl,*}$ is globally presentable, and that for every $G$ \[{-}\boxsmash{-}\colon\myS_{\text{$G$-gl},*}\times\myS_{\text{$G$-gl},*}\to\myS_{\text{$G$-gl},*}\] preserves colimits in each variable as it can be identified with the left derived functor of the left Quillen bifunctor ${-}\boxsmash{-}\colon\cat{$\bm G$-$\cat{L}$-Top}\times\cat{$\bm G$-$\cat{L}$-Top}\to\cat{$\bm G$-$\cat{L}$-Top}$. Observe now that the symmetric monoidal functor $(-)_+\colon\cat{$\cat{L}$-Top}\to\cat{$\cat{L}$-Top}_*$ induces a homotopical symmetric monoidal global functor $\smash{\Ntop(\cat{$\cat{L}$-Top}^{\boxtimes,\dual,\text{qcof}})\to\Ntop(\cat{$\cat{L}$-Top}_*^{\smallboxsmash,\dual,\text{qcof}})}$, and hence localizes to a symmetric monoidal global functor $\ul\myS_\gl^\otimes\to\ul\myS_{\gl,*}^{\otimes}$. Thus, we in particular have a commutative square as depicted on the left:
        \[
            \begin{tikzcd}
                \ul\myS_\gl\times\ul\myS_\gl\arrow[d,"(-)_+\times(-)_+"']\arrow[r, "-\boxtimes-"] & \ul\myS_\gl\arrow[d,"(-)_+"]
                &[2em]\ul\myS_\gl\times\ul\myS_\gl\arrow[d,"(-)_+\times(-)_+"']\arrow[r, "-\boxtimes-"] & \ul\myS_\gl\arrow[d,"(-)_+"]
                \\
                \ul\myS_{\gl,*}\times\ul\myS_{\gl,*}\arrow[r,"{-}\boxsmash{-}"'] & \ul\myS_{\gl,*}&
                \ul\myS_{\gl,*}\times\ul\myS_{\gl,*}\arrow[r,"{-}\otimes{-}"',dashed] & \ul\myS_{\gl,*}\rlap.
            \end{tikzcd}
        \]
        On the other hand, by the universal property of $\ul\myS_\gl\to\ul\myS_{\gl,*}$ as the initial globally cocontinuous map with a pointed target and since $\ul\myS_\gl^\otimes$ is globally presentably symmetric monoidal (Corollary~\ref{cor:S-gl-times-initial}), there exists a unique bifunctor ${-}\otimes{-}$ preserving global colimits in each variable and making the diagram on the right commute. Using that the initial globally cocontinuous functor $\ul\myS_\gl\to\ul\myS_{\gl,*}$ is also the initial \emph{fiberwise} cocontinuous functor (see Proposition~\ref{prop:univ-prop-*-cc-param}), ${-}\otimes{-}$ is then also uniquely characterized by preserving \emph{fiberwise} colimits in each variable and making the diagram commute. Thus, ${-}\boxsmash{-}$ is equivalent to ${-}\otimes{-}$, and hence in particular preserves global colimits in each variable.
    \end{proof}
\end{thm}

Analogously one shows:

\begin{cor}\label{cor:equiv-smash-initial}
    The pointed symmetric monoidal global $\infty$-category $\ul\myS_*^\otimes$ is equivariantly presentably symmetric monoidal, and it defines an initial object of $\CAlg(\CAT_{\smash{\Glo}}^{\textup{$\Orb$-cc},0})$.\qed
\end{cor}

\chapter[The stable story]{\for{toc}{\phantom{I}}The stable story}
In this second part we will study equivariant and global spectra and show that they can be organized into global $\infty$-categories $\ul\mySp$ and $\ul\mySp_\gl$, respectively. As our main results, we will show that these global $\infty$-categories admit universal properties as the free equivariantly or globally cocomplete $\infty$-categories satisfying a genuine refinement of stability, that we call \emph{representation stability}.

\section[Representation stability]{Representation stability}\label{sec:rep-stable}
In this section we introduce the notion of {representation stability} for equivariantly or globally cocomplete pointed global $\infty$-categories. We will moreover show (Corollary~\ref{cor:stabilization-exists}) that any globally cocomplete $\Cc$ admits an initial globally cocontionus functor to a representation stable globally cocomplete global $\infty$-category, and similarly in the equivariantly cocomplete case. These two notions of stabilization will however differ in general, with the globally cocomplete case being more subtle. In particular, stabilization in the equivariant setting will later turn out to be a pointwise procedure, while the global version is a genuinely parametrized phenomenon.

\subsection{$\bm{\mathcal R}$-stability} The notion of representation stability is a special case of a more general concept, which we will introduce in this subsection. For this we will first need to talk about the canonical tensoring over pointed parametrized spaces:

\begin{lemma}
    Let $S\subset T$ be any cleft, and let $\Cc$ be any pointed $S$-cocomplete $T$-$\infty$-category. Then there exists a unique $S$-bilinear bifunctor 
    \begin{equation}\label{eq:canonical-tensoring}
        {-}\otimes{-}\colon\ul\Spc_{S\triangleright T,*}\times\Cc\to\Cc
    \end{equation} 
    such that $1_+\otimes{-}$ is the identity of $\Cc$.
    \begin{proof}
        Passing to adjuncts, the claim translates to saying that there is a unique $S$-cocontinuous $T$-functor $\theta\colon\Cc\to\ul\Fun_T^\text{$S$-cc}(\ul\Spc_{S\triangleright T,*},\Cc)$ such that postcomposing with $\smash{\ev_{1^+}\colon\Fun_T^\text{$S$-cc}(\ul\Spc_{S\triangleright T,*},\Cc)\to\Cc}$ recovers the identity. But $\ev_{1^+}$ is an equivalence by Proposition~\ref{prop:univ-prop-*-cc-param}, so this is obvious.
    \end{proof}
\end{lemma}

\begin{rk}
    The above uniqueness argument did not use $S$-cocontinuity of $\theta$. In other words, the tensor product is also uniquely characterized by preserving $S$-colimits \emph{in the first variable} and satisfying $1_+\otimes{-}=\id_\Cc$.
\end{rk}

\begin{rk}
    Using the theory of smashing localizations, one can upgrade the canonical tensoring to an $\smash{\ul\Spc_{S\triangleright T,*}^\otimes}$-module structure in a unique way; for our present purposes, the above less coherent version will however be sufficient.
\end{rk}

\begin{ex}\label{ex:based-global-spaces-tensored}
    If $S=T=\Glo$, then we have an equivalence $\smash{\ul\myS_{\gl,*}\simeq\ul\Spc_{\Glo,*}}$ sending $S^0$ to $1_+$. The canonical tensoring of a pointed globally cocomplete global $\infty$-category may therefore equivalently viewed as a bifunctor ${-}\otimes{-}\colon\ul\myS_{\gl,*}\times\Cc\to\Cc$ preserving global colimits in each variable separately and such that $S^0\otimes{-}=\id_\Cc$.
    The canonical tensoring of $\ul\myS_{\gl,*}$ over itself is then given by the symmetric monoidal product from Theorem~\ref{thm:gl-smash-initial}, i.e.\ by the derived box smash product.
\end{ex}

\begin{ex}
    Analogously to the previous example, the canonical tensoring of the pointed equivariantly presentable global $\infty$-category $\ul\myS_*$ may be identified with the derived smash product from Corollary~\ref{cor:equiv-smash-initial}.
\end{ex}

\begin{ex}\label{ex:tensoring-pointwise}
    Let $\Cc$ be $S$-cocomplete and pointed, and let $\Kk$ be any $T$-$\infty$-category. Then $\ul\Fun_T(\Kk,\Cc)$ is again $S$-cocomplete and pointed, and the composite
    \begin{multline*}
        \ul\Spc_{S\triangleright T,*}\times\ul\Fun_T(\Kk,\Cc)\xrightarrow{\;\const_{\Kk}\times\id\;}
        \ul\Fun_T(\Kk,\ul\Spc_{S\triangleright T,*})\times\ul\Fun_T(\Kk,\Cc)\\
        {}\simeq\ul\Fun_T(\Kk,\ul\Spc_{S\triangleright T,*}\times\Cc)\xrightarrow{\;\ul\Fun_T(\Kk,{-}\otimes{-})\;}\ul\Fun_T(\Kk,\Cc)
    \end{multline*}     
    is $S$-bilinear and restricts to the identity on $\{1_+\}\times\ul\Fun_T(\Kk,\Cc)$, so it agrees with the canonical tensoring of $\ul\Fun_T(\Kk,\Cc)$. Put more informally, the canonical tensoring on a functor category is pointwise.
\end{ex}

\begin{ex}\label{ex:tensoring-slice}
    Let $\Cc$ be a $T$-cocomplete and pointed $T$-$\infty$-category, let $A\in T$ arbitrary, and let $\fgt\colon T_{/A}\to T$ denote the forgetful functor. Then $\fgt^*\Cc$ and $\fgt^*\ul\Spc_{T,*}$ are $T_{/A}$-cocomplete pointed $T_{/A}$-$\infty$-categories by direct inspection. Moreover, \cite{martini2021yoneda}*{Remark~3.7.2} shows that the unique left adjoint $\ul\Spc_{T_{/A}}\to\fgt^*\ul\Spc_T$ preserving the terminal object is an equivalence, so passing to pointed objects levelwise provides an equivalence $\ul\Spc_{T_{/A},*}\simeq\fgt^*\ul\Spc_{T,*}$ preserving $1_+$. Thus, we may identify the canonical tensoring of $\fgt^*\Cc$ with
    \[
        \ul\Spc_{T_{/A},*}\times\fgt^*\Cc\simeq\fgt^*\big(\ul\Spc_{T,*}\times\Cc\big)\xrightarrow{\;\fgt^*({-}\otimes{-})\;}\fgt^*\Cc.
    \]
\end{ex}

For later use we note the following naturality property of the tensoring:

\begin{lemma}\label{lemma:tensoring-natural}
    Let $\Cc,\Dd$ be pointed $S$-cocomplete $T$-$\infty$-categories, and let $X$ be a global section of $\ul\Spc_{S\triangleright T,*}$. Then the two $T$-functors
    \[
        \ul\Fun^\textup{$S$-cc}(X\otimes{-},\Dd),
        \ul\Fun^\textup{$S$-cc}(\Cc,X\otimes{-})\colon\ul\Fun^\textup{$S$-cc}(\Cc,\Dd)\rightrightarrows\ul\Fun^\textup{$S$-cc}(\Cc,\Dd)
    \]
    are naturally equivalent.
\end{lemma}

In particular, passing to underlying $\infty$-categories this says that $F(X\otimes-)\simeq X\otimes F(-)$ naturally in $F\in\Fun^\text{$S$-cc}(\Cc,\Dd)$. 

\begin{proof}
    The adjunct of
    \[\hskip-15.16pt\hfuzz=15.17pt
        \ul\Fun^\text{$S$-cc}(\Cc,\Dd)\xrightarrow{\;\ul\Fun(-\otimes-,\Dd)\;}
        \ul\Fun^\text{$S$-bilinear}(\ul\Spc_{S\triangleright T,*}\times\Cc,\Dd)\simeq
        \ul\Fun^\text{$S$-cc}\big(\ul\Spc_{S\triangleright T,*},\ul\Fun^\text{$S$-cc}(\Cc,\Dd)\big)
    \]
    defines an $S$-bilinear functor $\Theta\colon\ul\Spc_{S\triangleright T,*}\times\ul\Fun^\text{$S$-cc}(\Cc,\Dd)\to \ul\Fun^{\text{$S$-cc}}(\Cc,\Dd)$. By a straight-forward naturality argument, the restriction of $\Theta$ to a global section $X$ is given by $\ul\Fun^\text{$S$-cc}(X\otimes-,\Dd)$; in particular, $\Theta(\id,-)$ is the identity. Thus, $\Theta$ necessarily agrees with the tensoring of $\ul\Fun^\text{$S$-cc}(\Cc,\Dd)$ from Example~\ref{ex:tensoring-pointwise}, and the claim follows by restricting both tensorings to $\{X\}\times\ul\Fun^\text{$S$-cc}(\Cc,\Dd)$.
\end{proof}

We are now ready to define the notion of $\Rr$-stability:

\begin{definition}
    Let $\Rr\subset\ul\Spc_{S\triangleright T,*}$ be a full subcategory. We say that a pointed $S$-cocomplete $T$-$\infty$-category $\Cc$ is \emph{$\Rr$-stable} if for every $A\in T$ and $X\in\Rr(A)$ the functor $X\otimes{-}\colon\Cc(A)\to\Cc(A)$ induced by the canonical tensoring is an equivalence.
\end{definition}

\begin{ex}
    Let $S\subset T$ be any cleft, and let $\Rr\subset\ul\Spc_{S\triangleright T,*}$ be the full subcategory spanned by the global section $S^0=1_+$. Then every $S$-cocomplete $T$-$\infty$-category is $\Rr$-stable.
\end{ex}

\begin{ex}
    Let $T$ be any small $\infty$-category, and let $S=\core T$ be the minimal cleft, so that $\ul\Spc_{S\triangleright T,*}$ is constant with value $\Spc_*$. We can then take $\Rr$ to be the full subcategory spanned in each degree by the circle $S^1$, in which case an $S$-cocomplete (i.e.,\ fiberwise cocomplete) $T$-$\infty$-category $\Cc$ is $\Rr$-stable if and only if it is fiberwise stable in the sense that each $\Cc(A)$ is a stable $\infty$-category (note that all restriction functors are automatically exact as they are cocontinuous).
\end{ex}

\begin{ex}
    Let $T=\Glo$ and $S=\Orb$. We take $\Rr=\textup{RepSph}$ to be the full subcategory of $\ul\Spc_{S\triangleright T,*}$ corresponding under the equivalence from Corollary~\ref{cor:equiv-pointed} to the full subcategory of $\ul\myS_*$ given in each degree $G$ by the \emph{representation spheres}, i.e.\ the 1-point compactifications of $G$-representations. We call a pointed equivariantly cocomplete global $\infty$-category $\Cc$ \emph{representation stable} if it is $\textup{RepSph}$-stable. Note that if $\Cc$ is actually \emph{globally} cocomplete, this is equivalent to being $\textup{RepSph}'$-stable, where $\textup{RepSph}'\subset\ul\Spc_{\Glo,*}$ denotes the image of $\textup{RepSph}$ under the unique symmetric monoidal equivariantly cocontinuous functor $\ul\Spc_{\Orb\triangleright\Glo,*}\to\ul\Spc_{\Glo,*}$ or, equivalently, the image of the representation spheres under the composite
    \[
        \ul\myS_*\xrightarrow{\;\consto\;}\ul\myS_{\gl,*}\iso\ul\Spc_{\Glo,*}
    \]
    of the functor sending a $G$-space to the constant orthogonal $G$-space and the equivalence from Corollary~\ref{cor:gl-pointed}.
\end{ex}

\begin{variant}
    For any non-empty collection $\Ff$ of compact Lie groups, we obtain an analogous notion of representation stability for (equivariantly or globally cocomplete) $\Ff$-global $\infty$-categories by restricting the previous example.
    For $\Ff=\Fin$, these notions have been studied by Linskens in \cite{Linskens2023globalization}. By the main result of \cite{kaif-sil-excisive}, this agrees with the notion of \emph{equivariant stability} for $\Glo_{\Ff\kern-.5pt in}$-$\infty$-categories from \cites{CLL_Global,CLL_Clefts} defined via `abstract Wirthmüller isomorphisms.'
\end{variant}

\begin{rk}
    Let $\Rr\subset\ul\Spc_{T,*}$ arbitrary. The notion of $\Rr$-stability is compatible with slicing the indexing category in the following sense: Write $\fgt\colon T_{/A}\to T$ for the forgetful functor and $\Rr_A$ for the image of $\fgt^*\Rr$ under the equivalence $\fgt^*\ul\Spc_{T,*}\simeq\ul\Spc_{T_{/A},*}$ induced by the unique equivalence $\fgt^*\ul\Spc_T\simeq\ul\Spc_{T_{/A}}$. If $\Cc$ is $T$-cocomplete and $\Rr$-stable, then Example~\ref{ex:tensoring-slice} shows that $\fgt^*\Cc$ is $\Rr_A$-stable.
\end{rk}

\begin{definition}
    We write $\Pr^{S,\Rr\text{-ex}}_T\subset\Pr^S_T$ and $\CAT_T^{\text{$S$-cc},\Rr\text{-ex}}\subset\CAT_T^\text{$S$-cc}$ for the full subcategories spanned by the $\mathcal R$-stable $T$-$\infty$-categories.
\end{definition}

\subsection{Colimits of cocomplete parametrized $\bm\infty$-categories} Fix a cleft $S\subset T$. As the main result of this section, we want to prove that the inclusion $\CAT_T^\text{$S$-cc,\,$\Rr$-ex}\hookrightarrow\CAT_T^\text{$S$-cc}$ admits a left adjoint (giving rise to \emph{$\Rr$-stabilizations}) for any small $\Rr\subset\ul\Spc_{S\triangleright T,*}$, and that this adjunction restricts to the presentable case. The key step in this proof will be the construction of the $\Rr$-stabilization of $\ul\Spc_{S\triangleright T}$, which will be as a certain colimit in $\CAlg(\Pr^S_T)$. To see that this has the correct universal property in $\CAT_T^\text{$S$-cc}$ requires understanding the relation between colimits in $\Pr^S_T$ and $\CAT_T^\text{$S$-cc}$, to which the present subsection is devoted.

We begin with the following more general statement about parametrized $\infty$-categories with all colimits of certain prescribed shapes, which for $T=1$ is a special case of \cite{HA}*{Lemma 4.8.4.2}:

\begin{lemma}\label{lemma:small-class-of-colimits-cc}
    Let $\cat{U}\subset\ul\Cat_T$ be small. Then the $\infty$-category $\Cat_T^\textup{$\cat{U}$-cc}$ of $\cat{U}$-cocomplete small $T$-$\infty$-categories and $\cat{U}$-cocontinuous functors has small colimits, and for any small diagram $\Cc_\bullet\colon I\to\Cat_T^\textup{$\cat{U}$-cc}$ and any $\Dd\in\Cat_T^\textup{$\cat{U}$-cc}$ the induced map
    \begin{equation}\label{eq:internal-colimit}
        \ul\Fun_T^\textup{$\cat{U}$-cc}(\colim\nolimits_I \Cc_\bullet,\Dd)\to\lim\nolimits_{I^\op}\ul\Fun_T^\textup{$\cat{U}$-cc}(\Cc_\bullet,\Dd)
    \end{equation}
    is invertible.
    \begin{proof}
        Fix any small diagram $\Cc_\bullet\colon I\to\Cat_T^\textup{$\cat{U}$-cc}$ and an extension to a colimit diagram $I^\triangleright\to\Cat_T$. Writing $\infty$ for the terminal object of $I^\triangleright$, it follows that for any $\cat{U}$-cocomplete $\Dd$, the induced equivalence
        \[
            \hom_{\Cat_T}(\Cc_\infty,\Dd)\iso\lim\nolimits_{I^\op}\hom_{\Cat_T}(\Cc_\bullet,\Dd)
        \]
        then restricts to an equivalence between the union of those path components of the left-hand side given by functors $F\colon\Cc\to\Dd$ such that each $\Cc_i\to\Cc\to\Dd$ is $\cat{U}$-cocontinuous and $\lim_{I^\op}\hom_{\Cat_T^\textup{$\cat{U}$-cc}}(\Cc_\bullet,\Dd)$. On the other hand, \cite{martiniwolf2022presentable}*{Proposition~A.4.2 and Remark~A.4.3} yield a map $\Cc_\infty\to\Cc_\infty'$ to a small $\cat{U}$-cocomplete $T$-$\infty$-category such that $\hom_{\Cat_T^\textup{$\cat{U}$-cc}}(\Cc_\infty',\Dd)\to\hom_{\Cat_T}(\Cc_\infty,\Dd)$ defines an equivalence of the union of the same set of path components. This provides the desired colimit.
        
        To see that $(\ref{eq:internal-colimit})$ is an equivalence, it suffices by the (non-parametrized) Yoneda lemma, that the induced map on $\hom(\Kk,-)$ is an equivalence for every $\Kk\in\Cat_T$. However, this induced map may be identified with the equivalence
        \[
            \hom_{\Cat_T^\text{$\cat{U}$-cc}}(\colim\nolimits_I\Cc_\bullet,\Fun(\Kk,\Dd))\iso
            \lim\nolimits_{I^\op}\hom_{\Cat_T^\text{$\cat{U}$-cc}}(\Cc_\bullet,\Fun(\Kk,\Dd)).
            \qedhere
        \]        
    \end{proof}
\end{lemma}

Applying this in a larger universe we get:

\begin{cor}
    The very large $\infty$-category $\CAT_T^\textup{$S$-cc}$ has all large colimits, and for any diagram $\Cc_\bullet\colon I\to\CAT_T^\textup{$S$-cc}$ and $\Dd\in\CAT_T^\textup{$S$-cc}$ the induced map
    \[
        \ul\Fun_T^\textup{$S$-cc}(\colim\nolimits_I \Cc_\bullet,\Dd)\to\lim\nolimits_{I^\op}\ul\Fun_T^\textup{$S$-cc}(\Cc_\bullet,\Dd)
    \]
    is invertible.\qed
\end{cor}

\begin{prop}\label{prop:PrS-T-vs-Cat-S-cc}
    The full subcategory $\Pr^S_T\subset\CAT_T^\textup{$S$-cc}$ is closed under \emph{small} colimits. In particular, $\Pr^S_T$ is cocomplete.
\end{prop}

The proof will require some preparations. The basic idea will be to factor any small diagram $I\to\Pr^S_T$ through a levelwise $\Ind_\kappa$-cocompletion functor from (small) $T$-$\infty$-categories with `$\kappa$-small $S$-colimits' to $\Pr^S_T$. For this to make sense, we have to put some assumptions on $\kappa$:

\begin{definition}
    A regular cardinal $\kappa$ is \emph{$(S,T)$-regular} if for every pullback
    \begin{equation}\label{diag:pb-regular-card}
        \begin{tikzcd}
            X\arrow[r, "f'"]\arrow[dr,pullback]\arrow[d, "s'"'] & A\arrow[d,"s"]\\
            C\arrow[r,"f"'] & B
        \end{tikzcd}
    \end{equation}
    in $\PSh(T)$ such that $A,B,C$ are representable and $s$ is a map in $S$ (so that $X\in\PSh(S)$ by the definition of a cleft), $X$ can be written as a $\kappa$-small colimit in $\PSh(S)$ of representables.
\end{definition}
Note that there is only a set's worth of pullback squares of the form $(\ref{diag:pb-regular-card})$ (up to isomorphism), so there always exists an $(S,T)$-regular $\kappa$. Moreover, if $\mu>\kappa$ is regular and $\kappa$ is $(S,T)$-regular, then so is $\mu$; thus, there exist arbitrarily large $(S,T)$-regular cardinals.
\begin{definition}
    Let $\kappa$ be $(S,T)$-regular. We say a $T$-$\infty$-category $\Cc$ has \emph{$\kappa$-small $S$-colimits} if it has $\cat{U}$-colimits for $\cat{U}$ the union of constant $\kappa$-small categories and the Yoneda images of all maps $s'$ as in $(\ref{diag:pb-regular-card})$.
\end{definition}

\begin{lemma}\label{lemma:fiberwise-kappa-compactly-gen}
    Let $\Cc\in\Pr^S_T$ and let $\kappa$ be $(S,T)$-regular such that each $\Cc(A)$ is $\kappa$-compactly generated and each restriction preserves $\kappa$-compact objects. 
    Denote by $\Cc^{\kappa}\subset \Cc$ the full $T$-subcategory given in degree $A\in T$ by the $\kappa$-compact objects. 
    Then $\Cc^\kappa$ has $\kappa$-small $S$-colimits, and for any $\Dd\in\CAT_T^\textup{$S$-cc}$, restriction along the inclusion defines an equivalence
    \[
        \ul\Fun^\textup{$S$-cc}(\Cc,\Dd)\iso\ul\Fun^\textup{$\kappa$-small $S$-cc}(\Cc^\kappa,\Dd). 
    \]
    \begin{proof}
        It follows directly from the assumptions that restriction along $\Cc^\kappa\hookrightarrow\Cc$ defines an equivalence $\ul\Fun_T^\text{fib $\kappa$-filt cc}(\Cc,\Dd)\iso\ul\Fun_T(\Cc^\kappa,\Dd)$. Moreover, $\Cc^\kappa$ has fiberwise $\kappa$-small colimits, and applying \cite{HTT}*{Example~5.3.6.8} pointwise, we see that a functor $\Cc\to\ul\Fun(\ul A,\Dd)$ preserves fiberwise colimits if and only if its restriction to $\Cc^\kappa$ preserves fiberwise $\kappa$-small colimits. 

        As each of the left adjoints $s_!\colon\Cc(A)\to\Cc(B)$ preserves $\kappa$-compact objects (its right adjoint admits another right adjoint, so it is in particular $\kappa$-accessible), we see that each $s^*\colon\Cc^\kappa(B)\to\Cc^\kappa(A)$ for $s\colon A\to B$ in $S$ admits a left adjoint. We claim that also $s_!'\colon\Cc(X)\to\Cc(C)$ for $s'$ as in $(\ref{diag:pb-regular-card})$ restricts to $\Cc^\kappa(X)\to\Cc^\kappa(C)$ (beware that the left-hand side is a priori different from $\Cc(X)^\kappa$); with this established the Beck--Chevalley condition for $\Cc^\kappa$ will follow from the one for $\Cc$, i.e.~$\Cc^\kappa$ will have $\kappa$-small $S$-colimits.

        To prove the claim, write $X\simeq\colim_iX_i$ as a $\kappa$-small colimit of representables in $\PSh(S)$, so that $s^*$ agrees up to equivalence with $(s_i^*)_{i\in I}\colon\Cc(C)\to\lim_{I^\op}\Cc(X_\bullet)$ for some cocone $(s_i)_{i\in I}$. By \cite{descent-lim}*{Theorem~B}, this then admits a left adjoint given on objects by sending an object $(y_i)_{i\in I}\in\lim_{I^\op}\Cc(X_\bullet)$ to the colimit of a certain diagram $i\mapsto s_{i!}(y_i)$. If each $y_i$ is $\kappa$-compact, then so is each $s_{i!}(y_i)$ by the above, whence so is their colimit as $I$ is $\kappa$-small. This completes the proof of the claim, and hence of $\Cc^\kappa$ having $\kappa$-small $S$-colimits.

        It only remains to show that a functor $F\colon\Cc\to\ul\Fun(\ul A,\Dd)$ for $S$-cocomplete $\Dd$ satisfies the Beck--Chevalley condition with respect to the maps $s\colon A\to B$ in $S$ if and only if this holds for the composite $Fi\colon\Cc^\kappa\to\ul\Fun(\ul A,\Dd)$. This follows at once by observing that the collection of $X\in\Cc(A)$ such that $s_!FX\iso Fs_!X$ is closed under ($\kappa$-filtered) colimits and contains $\Cc^\kappa(A)$, as $\Cc^\kappa\hookrightarrow\Cc$ satisfies the Beck--Chevalley condition by construction of the left adjoints for $\Cc^\kappa$.
    \end{proof}
\end{lemma}

\begin{lemma}\label{lemma:fiberwise-ind-kappa-completion}
    Let $\kappa$ be $(S,T)$-regular, and let $\Cc$ be a $T$-$\infty$-category with $\kappa$-small $S$-colimits. Then ${\Ind_\kappa}\circ\Cc$ is $S$-cocomplete, and the inclusion $\Cc\hookrightarrow {\Ind_\kappa}\circ\Cc$ preserves $\kappa$-small $S$-colimits.
    \begin{proof}
        It is clear that ${\Ind_\kappa}\circ\Cc$ is fiberwise cocomplete and that the inclusion preserves fiberwise $\kappa$-small colimits. Moreover, for every $s\colon A\to B$ in $S$, the restriction $\Ind_\kappa \Cc(B)\to\Ind_\kappa\Cc(A)$ admits a left adjoint, and the inclusion satisfies the Beck--Chevalley condition with respect to this adjoints.

        We now again apply \cite{descent-lim}*{Theorem~B} to see that for any $s'\colon X\to C$ in $\PSh(S)$ such that $X$ is a $\kappa$-small colimit $\colim_{I}X_\bullet$ of representables, the restriction $\Cc(B)\to ({\Ind_\kappa}\circ\Cc)(X)$ admits a left adjoint; beware that the target will be typically different from $\Ind_\kappa(\Cc(X))$. We now claim that the square
        \[
            \begin{tikzcd}
                \Cc(X)\arrow[r,hook]\arrow[from=d,"s'^*"] & ({\Ind_\kappa}\circ\Cc)(X)\arrow[from=d,"s'^*"']\\
                \Cc(C)\arrow[r,hook] & ({\Ind_\kappa}\circ\Cc)(C)
            \end{tikzcd}
        \]
        is vertically left adjointable. As the horizontal maps are fully faithful, this amounts to saying that $s_!'\circ\incl$ factors through the bottom inclusion. Applying the description of $s_!'$ from \cite{descent-lim} again shows that $s_!'\circ\incl$ sends a $y\in\Cc(X)$ with components $y_i\in\Cc(X_i)$ to the colimit of some $I$-indexed diagram $i\mapsto s_{i!}'y_i$. By the above, each $s_{i!}'y_i$ is contained the image of the bottom inclusion; the claim follows as the latter is closed under $\kappa$-small colimits.

        To complete the proof, it then only remains to show that the Beck--Chevalley map $s'_!f'^*\to f^*s_!$ of functors $\Ind_\kappa\Cc(A)\to\Ind_\kappa\Cc(C)$ is invertible for any pullback $(\ref{diag:pb-regular-card})$. As both source and target preserve ($\kappa$-filtered) colimits, it suffices to check this after precomposing with $\Cc(A)\hookrightarrow\Ind_\kappa\Cc(A)$. Using that Beck--Chevalley compose, this follows from the previous paragraph and the assumption that $\Cc$ have $\kappa$-small $S$-colimits.
    \end{proof}
\end{lemma}

\begin{proof}[Proof of Proposition~\ref{prop:PrS-T-vs-Cat-S-cc}]
    Let $\Cc_\bullet\colon I\to\Pr^S_T\subset\CAT_T^\text{$S$-cc}$ be a \emph{small} diagram, and pick an $(S,T)$-regular $\kappa>\aleph_0$ such that each $\Cc_i(A)$ is $\kappa$-compactly generated and all restrictions $\Cc_i(B)\to\Cc_i(A)$ as well as all transition maps $\Cc_i(A)\to\Cc_j(A)$ preserve $\kappa$-compact objects. We then obtain a subfunctor $\Cc^\kappa_\bullet\colon I\to\Cat_T^\text{$\kappa$-small cc}$, which admits a colimit $\Dd$ by Lemma~\ref{lemma:small-class-of-colimits-cc}. By the previous lemma, $\Ind_\kappa\circ\Dd$ is $S$-cocomplete, and as $\Dd$ is small, it is even fiberwise presentable, i.e.\ $S$-presentable. Moreover, $\Dd$ is fiberwise idempotent complete (as it has sequential colimits by assumption on $\kappa$), so $\Dd\hookrightarrow\Ind_\kappa\circ\Dd$ induces an equivalence onto $({\Ind_\kappa}\circ\Dd)^\kappa$. Applying Lemma~\ref{lemma:fiberwise-kappa-compactly-gen} twice then shows that ${\Ind_\kappa}\circ\Dd$ is a colimit of ${\Ind_\kappa}\circ\Cc_\bullet^\kappa\simeq\Cc_\bullet$, finishing the proof.
\end{proof}

\subsection{Existence of $\bm{\mathcal R}$-stabilizations}\label{subsec:invert-objects} In this subsection we will finally show that the subcategory of $\Rr$-stable categories is a smashing localization of the $\infty$-category of $S$-cocomplete, pointed $T$-$\infty$-categories. 
To this end, we will more generally study the process of inverting a collection of objects in an $S$-cocomplete symmetric monoidal $T$-$\infty$-category. We begin with the presheaf case:

\begin{rk}\label{rk:sym-mon-cocompletion}
    By \cite{martiniwolf2022presentable}*{Remarks~2.6.2.6 and~2.6.2.7}, the lax symmetric monoidal functor $(\CAT_T^\textup{$S$-cc})^\otimes\hookrightarrow\CAT_T^\times$ admits a strong symmetric monoidal left adjoint, which then lifts to a left adjoint $\mathcal P_{S\triangleright T}(-)^{\otimes}\colon\CAlg(\CAT_T)\to\CAlg(\CAT_T^\text{$S$-cc})$ of the forgetful functor. 
    By Remark~\ref{rk:PSh-ST}, the functor $\mathcal P_{S\triangleright T}(-)^{\otimes}$ restricts to $\CMon(\Cat_T)\to\CAlg(\Pr^S_T)$.
\end{rk}

\begin{lemma}
    Let $X\in\CMon(\SPC_T)$, and let $f\colon X\to X^{\grp}$ be its (pointwise) group completion. Then 
    \begin{equation}\label{eq:restr-Pf}
        \mathcal P(f)^*\colon\textup{LMod}_{\mathcal P_{S\triangleright T}(X^{\grp})^{\otimes}}(\CAT_T^\textup{$S$-cc})\to \textup{LMod}_{\mathcal P_{S\triangleright T}(X)^{\otimes}}(\CAT_T^\textup{$S$-cc})
    \end{equation}
    is fully faithful, and its essential image consists precisely of those $\mathcal P_{S\triangleright T}(X)^{\otimes}$-modules $\Mm$ such that for every $A\in T$ and $x\in X(A)$ the action map $x\cdot-\colon\Mm(A)\to\Mm(A)$ is invertible.
    \begin{proof}
        By definition, each $X(A)\to X^\grp(A)$ is an epimorphism in $\CMon(\SPC)$, so that $X\to X^\grp$ is an epimorphism in $\CMon(\SPC_T)\simeq\Fun\big(T^\op,\CMon(\SPC)\big)$. As left adjoints preserve epimorphisms, we conclude that also $\mathcal P(f)\coloneqq \mathcal P_{S\triangleright T}(f)^{\otimes}$ is an epimorphism, so that $(\ref{eq:restr-Pf})$ is fully faithful by \cite{Ro15}*{Proposition 2.4}.
        
        To characterize its essential image, it remains to show that a $\mathcal P_{S\triangleright T}(X)^{\otimes}$-module structure on a given $\Mm$ lifts to a $\mathcal P_{S\triangleright T}(X^{\grp})^{\otimes}$-module structure if and only if each $x\in X(A)$ acts invertibly. For this we note that since the symmetric monoidal structure on $\CAT_T^\textup{$S$-cc}$ is closed, it in particular admits endomorphism objects in the sense of \cite{HA}*{§4.7.1}. Thus, the composition on $\ul\Fun^\textup{$S$-cc}(\Mm,\Mm)$, i.e.~the map adjunct to
        \[
            \Mm\otimes\ul\Fun^\textup{$S$-cc}(\Mm,\Mm)\otimes\ul\Fun^\textup{$S$-cc}(\Mm,\Mm)\xrightarrow{\ev\otimes\id}\Mm\otimes\ul\Fun^\textup{$S$-cc}(\Mm,\Mm)\xrightarrow{\ev}\Mm,
        \]
        refines to an associative algebra structure, and the adjunct $\tilde\mu$ of any $\mu\colon\mathcal P_{S\triangleright T}(X)^{\otimes}\otimes\Mm\to\Mm$ that is part of a module structure refines to a map of associative algebras. By \cite{HA}*{Corollary 4.7.1.40}, lifting such a structure to $\mathcal P_{S\triangleright T}(X^{\grp})^{\otimes}$ then amounts to solving the extension problem in $\Alg(\CAT^\text{$S$-cc}_T)$ depicted on the left in the following diagram:
        \[
            \begin{tikzcd}
                \mathcal P_{S\triangleright T}(X)^{\otimes}\arrow[d,"\mathcal P(f)^{\otimes}"']\arrow[r,"\tilde\mu"] & \ul\Fun^\textup{$S$-cc}(\Mm,\Mm) &[2em] X\arrow[d,"f"']\arrow[r,"\tilde\mu|_X"] & \core\ul\Fun^\textup{$S$-cc}(\Mm,\Mm)\rlap. \\
                \mathcal P_{S\triangleright T}(X^{\grp})^{\otimes}\arrow[ur, bend right=15pt, dashed] && X^\grp\arrow[ur, bend right=15pt, dashed]
            \end{tikzcd}
        \]        
        By adjunction, this is equivalent to solving the extension problem on the right in $\Alg(\SPC_T)=\Fun(T^\op,\textup{Mon}(\SPC))$. As group completions in $\CMon(\SPC)$ and $\textup{Mon}(\SPC)$ agree (both are given by the loop space of the bar construction), the extension exists if and only if $X$ has image in the units with respect to the algebra structure. By direct inspection, the multiplication of the algebra structure is given in degree $A$ by sending elements corresponding to (suitably cocontinuous) endomorphisms of $\pi_A^*\Mm$ to their composite, i.e.~the units precisely correspond to \emph{auto}morphisms of $\pi_A^*\Mm$. Further unravelling definitions, $X\to\ul\Fun^\textup{$S$-cc}(\Mm,\Mm)$ is given in degree $A$ by sending $x\in X(A)$ to the map corresponding to
        \[
            \pi_A^*\Mm\simeq1\times\pi_A^*\Mm\xrightarrow{\;x\times\id\;}\pi_A^*X\times\pi_A^*\Mm\to\pi_A^*\Mm,
        \]
        where the last map comes from the module structure. This is invertible if and only if it is so after evaluating at each $f\colon B\to A$ in $T_{/A}$, i.e.~if and only $f^*x$ acts invertibly on $\Mm(B)$. The lemma follows.
    \end{proof}
\end{lemma}

\begin{thm}\label{thm:sym-mon-localization-for-modes}
    Let $\Cc\in\CAlg(\CAT_T^\textup{$S$-cc})$, let $X\in\CMon(\SPC_T)$, and let $I\colon X\to\Cc$ be a map in $\CMon(\CAT_T)$. 
    
    \begin{enumerate}
        \item There exists a map $F\colon\Cc\to\Dd$ in $\CAlg(\CAT^\textup{$S$-cc}_T)$ such that $FI\colon X\to\Dd$ lands in the tensor-invertible objects, and which is universal in the sense that for any $\Ee\in\CAlg(\CAT_T^\textup{$S$-cc})$ the induced map
        \[
            -\circ F\colon\hom_{\CAT_T^\textup{$S$-cc}}(\Dd,\Ee)\to\hom_{\CAT_T^\textup{$S$-cc}}(\Cc,\Ee)
        \]
        is the inclusion of the components consisting of all $G\colon\Cc\to\Ee$ such that $GI$ factors through the invertible objects.
        \item If $\Cc$ is $S$-presentable and $X$ is small, then $\Dd$ is $S$-presentable.
        \item If $\Cc$ is idempotent, then so is $\Dd$. Moreover, an $\Mm\in\CAlg(\CAT_T^\textup{$S$-cc})$ is a $\Dd$-module if and only if it is a $\Cc$-module and the map $I(x)\cdot-\colon\Mm(A)\to\Mm(A)$ given by the $\Cc$-module structure is an equivalence for every $A\in T$ and $x\in X(A)$.
    \end{enumerate}
    \begin{proof} Consider the pushout
        \begin{equation}\label{diag:pushout-inversion}
            \begin{tikzcd}
                \mathcal P_{S\triangleright T}(X)^{\otimes}\arrow[r,"\mathcal P(f)^{\otimes}"]\arrow[d,"\tilde I"']\arrow[dr,phantom,"\ulcorner"{very near end}] & \mathcal P_{S\triangleright T}(X^{\grp})^{\otimes}\arrow[d]\\
                \Cc\arrow[r] & \Dd
            \end{tikzcd}
        \end{equation}
        in $\CAlg(\CAT_T^\textup{$S$-cc})$, where the left-hand vertical map is adjunct to $S$ and the horizontal map is induced by a group completion. We may then identify for any $\Ee\in\CAlg(\CAT_T^\textup{$S$-cc})$ the map $\hom(\Dd,\Ee)\to\hom(\Cc,\Ee)$ with the projection 
        \[
            \hom(\mathcal P_{S\triangleright T}(X^{\grp})^{\otimes},\Ee)\times_{\hom(\mathcal P_{S\triangleright T}(X)^{\otimes},\Ee)}\hom(\Cc,\Ee)\to\hom(\Cc,\Ee),
        \]
        which by adjunction is the inclusion of those path components of symmetric monoidal $S$-cocontinuous $G\colon\Cc\to\Ee$ such that $GI$ factors through invertible objects.

        For the second claim we observe that the assumptions on $X$ and $\Cc$ guarantee that all objects in $(\ref{diag:pushout-inversion})$ except possibly the bottom right one are $S$-presentable. As $\CAlg(\Pr^S_T)\subset\CAlg(\CAT_T^\textup{$S$-cc})$ is closed under small colimits as a consequence of Proposition~\ref{prop:PrS-T-vs-Cat-S-cc}, we conclude that also $\Dd$ is $S$-presentable.

        For the final claim, we first observe that $\Dd$ is a $\Cc$-algebra, hence in particular a $\Cc$-module. Thus $\ul\Spc_{S\triangleright T}\otimes\Dd\to\Cc\otimes\Dd$ is an equivalence, and it will suffice for idempotency that 
        \begin{equation}\label{eq:check-idempotency}
            \hom_{\CAlg(\CAT_T^\textup{$S$-cc})}(\Dd\otimes\Dd,\Ee)\iso\hom_{\CAlg(\CAT_T^\textup{$S$-cc})}(\Cc\otimes\Dd,\Ee)
        \end{equation}
        for any $\Cc$-algebra $\Ee$. Since the symmetric monoidal product is the coproduct on $\CAlg(\CAT_T^\textup{$S$-cc})$, we may view $(\ref{eq:check-idempotency})$ as a map
        \[
            \hom(\Dd,\Ee)\times\hom(\Dd,\Ee)\to\hom(\Cc,\Ee)\times\hom(\Dd,\Ee).
        \]
        As $\Cc$ is idempotent and $\Ee$ was assumed to be a $\Cc$-algebra, $\hom(\Cc,\Ee)\simeq1$ by Lemma~\ref{lemma:unique-algebra-structure-initial}. By the universal property of $\Dd$, we then have $\hom(\Dd,\Ee)\simeq1$ if the composite $X\to\Cc\to\Ee$ lands in the invertible objects, and $\hom(\Dd,\Ee)\simeq\emptyset$ otherwise. Thus, either both sides $(\ref{eq:check-idempotency})$ are initial, or both sides are terminal; in any case, any map between them is an equivalence.

        It remains to prove the characterization of $\Dd$-modules. Clearly, any $\Dd$-module $\Mm$ is also a $\Cc$-module and its underlying $\mathcal P_{S\triangleright T}(X)^{\otimes}$-module lifts to $\mathcal P_{S\triangleright T}(X^{\grp})^{\otimes}$, i.e.~$X$ acts invertibly by the previous lemma. For the converse, we compute for a $\Cc$-module $\Mm$
        \[
            \Dd\otimes\Mm\simeq\Dd\otimes_\Cc\Mm\simeq(\mathcal P_{S\triangleright T}(X^{\grp})\otimes_{\mathcal P_{S\triangleright T}(X)} \Cc)\otimes_\Cc\Mm\simeq \mathcal P_{S\triangleright T}(X^{\grp})\otimes_{\mathcal P_{S\triangleright T}(X)}\Mm,
        \]
        where the first equivalence is idempotency of $\Cc$, the second one uses that $(\ref{diag:pushout-inversion})$ exhibits $\Dd$ as the coproduct in $\CAlg(\CAT_T^\textup{$S$-cc})_{\mathcal P_{S\triangleright T}(X)^{\otimes}/}\simeq\CAlg(\text{LMod}_{\mathcal P_{S\triangleright T}(X)^{\otimes}})$, while the last equivalence follows from \cite{HA}*{Propositions~4.4.3.14 and 4.4.3.16}. If $X$ acts invertibly on $\Mm$, then the counit $\mathcal P_{S\triangleright T}(X^{\grp})^{\otimes}\otimes_{\mathcal P_{S\triangleright T}(X)^{\otimes}}\Mm\to\Mm$ is an equivalence of $\mathcal P_{S\triangleright T}(X)^{\otimes}$-modules by another application of the previous lemma, so that $\Dd\otimes\Mm\simeq\Mm$ abstractly in $\CAT_T^\textup{$S$-cc}$. As the left-hand side is a $\Dd$-module, so is the right-hand side, finishing the proof.
    \end{proof}
\end{thm}

\begin{cor}\label{cor:stabilization-exists}
    Let $\Rr\subset\smash{\ul\Spc_{S\triangleright T,*}}$ small. Then $\smash{\CAT_T^{S\textup{-cc},\textup{$\Rr$-ex}}}\hookrightarrow\CAT_T^\textup{$S$-cc}$ admits a left adjoint restricting to $\Pr^S_T\to\Pr^{S,\textup{$\Rr$-ex}}_T$, and both of these adjunctions are smashing localizations. In particular, there exists for any pointed $S$-presentable $\Cc$ an \emph{$\bm\Rr$-stabilization}, i.e.\ an $S$-cocontinuous $\Sigma^\infty\colon\Cc\to\ul\Sp^\Rr_{S\triangleright T}(\Cc)$ to an $S$-presentable $\Rr$-stable $T$-$\infty$-category, such that for any $S$-cocomplete $\Rr$-stable $\Dd$
    \[
        \ul\Fun^\textup{$S$-cc}(\Sigma^\infty,\Dd)\colon\ul\Fun^\textup{$S$-cc}(\ul\Sp^\Rr_{S\triangleright T}(\Cc),\Dd)\iso\ul\Fun^\textup{$S$-cc}(\Cc,\Dd).
    \]
    \begin{proof}
        Note that if $X,Y\in\Rr(A)$ act invertibly, then so does $X\otimes Y$. We may therefore assume that $\Rr$ is a symmetric monoidal subcategory, and the first claim follows by applying Theorem~\ref{thm:sym-mon-localization-for-modes} to $\core\Rr\hookrightarrow\ul\Spc_{S\triangleright T,*}$. The second claim is then an instance of Lemma~\ref{lemma:smash-module-internal-hom}.
    \end{proof}
\end{cor}

Thanks to the abstract theory we obtain several further universal properties of this construction:

\begin{cor}\label{cor:SpR-initial}
    There exists a unique $S$-presentably symmetric monoidal structure on $\ul\Sp_{S\triangleright T}^\Rr\coloneqq\ul\Sp_{S\triangleright T}^\Rr(\ul\Spc_{S\triangleright T,*})$ with unit $\mathbb S\coloneqq\Sigma^\infty(1_+)$, and this makes it into the initial $\Rr$-stable $S$-presentably symmetric monoidal $T$-$\infty$-category. 
    \begin{proof}
        Existence and uniqueness of the symmetric monoidal refinement were recalled in Remark~\ref{rk:idempotent-algs-are-algs}, and this is initial by Lemma~\ref{lemma:unique-algebra-structure-initial}.
    \end{proof}
\end{cor}

\begin{cor}\label{cor:SpR-sym-mon-univ-prop}
    For any $S$-presentably symmetric monoidal $\Dd$, the restriction
    \[
        \hom(\Sigma^\infty,\Dd)\colon \hom_{\CAlg(\Pr^S_T)}(\ul\Sp_{S\triangleright T}^\Rr,\Dd)\to \hom_{\CAlg(\Pr^S_T)}(\ul\Spc_{S\triangleright T,*},\Dd)
    \]
    is the inclusion of the components consisting of those $S$-cocontinuous symmetric monoidal functors that send the objects of $\Rr$ to invertible objects.
    \begin{proof}
        As part of Theorem~\ref{thm:sym-mon-localization-for-modes} we constructed an $S$-presentably symmetric monoidal structure with the correct universal property. This then necessarily agrees with the structure from the previous corollary.
    \end{proof}
\end{cor}

\section[The universal property of equivariant spectra]{The universal property of equivariant spectra}\label{sec:equiv-spectra}
In this section we will recall the model of $G$-equivariant stable homotopy theory in terms of orthogonal spectra with $G$-action, and we will show that for varying $G$ the resulting $\infty$-categories $\mySp_G$ of \emph{(genuine) $G$-spectra} naturally assemble into the free representation stable and equivariantly cocomplete global $\infty$-category.

\subsection{A reminder on (equivariant) orthogonal spectra}\label{subsec:equiv-spectra}
We begin by recalling the coordinate-free description of orthogonal spectra:

\begin{constr}
    We write $\cat{O}$ for the following $\cat{Top}_*$-enriched category: the objects of $\cat{O}$ are the finite dimensional inner product spaces, and the pointed space of maps $V\to W$ is given as the $1$-point compactification of the total space of the $\R$-vector bundle
    \[
        \{(A, w)\in\cat{L}(V,W)\times W : \langle Av,w\rangle = 0\text{ for all $v\in V$}\}\xrightarrow{\;\pr\;}\cat{L}(V,W),
    \]
    where $\cat{L}(V,W)$ denotes the space of linear isometric embeddings as before; this receives a natural quotient map from the fiberwise 1-point compactification, so we will denote points in $\cat{O}(V,W)$ by $[A,w]$ where $w\in S^W=W\cup\{\infty\}$ and $\langle w,AV\rangle=0$ for $w\not=\infty$. The composition law is then given by
    \begin{align*}
        \cat{O}(V,W)\smashp\cat{O}(U,V)&\longrightarrow\cat{O}(U,W)\\
        [B,w]\smashp[A,v]&\longmapsto[BA, Bv+w]\rlap.
    \end{align*}
\end{constr}

\begin{definition}
    We write $\osp$ for the ($\cat{Top}_*$-enriched) category of $\cat{Top}_*$-enriched functors $\cat{O}\to\cat{Top}_*$ and natural transformations between them. An object of $\osp$ is called an \emph{orthogonal spectrum}.
\end{definition}

\begin{constr}
    If $i\colon U\to V$ is a linear isometric embedding and $X$ is any orthogonal spectrum, then we obtain a map $S^{V-i(U)}\smashp X(U)\to X(V)$ as the adjunct of the composite
    \begin{align*}
        S^{V-i(U)}&\longrightarrow\cat{O}(U,V)\xrightarrow{\;X\;}\maps(X(U),X(V))\\
        v&\longmapsto [i,v]\rlap;
    \end{align*}
    we will refer to this as the \emph{structure map} of $X$. We will be particularly interested in the case of a direct summand inclusion $V\to V\oplus W$, in which case we can write the structure map as $\sigma_{V,W}\colon S^W\smashp X(V)\to X(V\oplus W)$. It is easy to check that this map is $(\O(V)\times \O(W))$-equivariant with respect to the evident $\O(W)$-action on $S^W$ and the actions via functoriality of $X$.
\end{constr}

\begin{rk}
    The category $\cat{O}$ has a skeleton given by the vector spaces $\R^n$ equipped with the standard scalar products, so it is enough to specify an orthogonal spectrum on the full subcategory spanned by them. Similarly, it is not hard to show that any morphism $\R^m\to\R^n$ in $\cat{O}$ is given up to conjugation by automorphisms by a composite $\R^m\to\R^{m+1}\to\cdots\R^n$, where each $\R^k\to\R^{k+1}$ is of the form $[\incl,(0,\dots,0,x)]$ for some $x\in\R\cup\{\infty\}$. 
    
    It follows that we may equivalently describe an orthogonal spectrum as a sequence $X_m=X(\R^m)$ of based $\O(m)$-spaces for all $m\ge0$, together with maps $S^1\smashp X_m\to X_{m+1}$ satisfying some compatibility conditions. This is the classical `coordinatized' description of orthogonal spectra; as we will never need this description, we simply refer the reader to \cite{schwede2018global}*{Remark 3.1.6} for details.
\end{rk}

For later use, let us describe the morphism spaces in $\cat{O}$ in concrete terms:

\begin{lemma}\label{lemma:O-mapping-spaces}
    Let $V,W$ be finite-dimensional inner product spaces.
    \begin{enumerate}
        \item If $\dim(V)>\dim(W)$, then $\cat{O}(V,W)=0$.
        \item If $V=W$, then $\cat{O}(V,V)=\O(V)_+$.
        \item If $V\subset W$, then there exists a left $\O(W)$- and right $\O(V)$-equivariant homeomorphism $\cat{O}(V,W)\cong \O(W)_+\smashp_{\O(W-V)}S^{W-V}$, where on the right-hand side $\O(W)$ and $\O(V)$ act via left and right multiplication on $\O(W)$, while they act on the left-hand side via post- and precomposition, respectively.
    \end{enumerate}
    \begin{proof}
        The first two statements are clear from the definitions. For the third statement, one easily checks that 
        \begin{align*}
            \O(W)_+\smashp_{\O(W-V)}S^{W-V} &\longrightarrow \cat{O}(V,W)\\
            [A,w] &\longmapsto [A|_{V}, Aw]
        \end{align*}
        is well-defined and an equivariant homeomorphism, also cf.\ \cite{schwede2018global}*{Construction~3.1.1}.
    \end{proof}
\end{lemma}

\begin{prop}\label{prop:O-mapping-CW}
    Let $G,H$ be compact Lie, let $V$ be an $H$-representation, and let $U$ be a $G$-representation. Then $\cat{O}(V,U)$ admits the structure of a finite based $(G\times H)$-CW-complex. Moreover, if $K\subset H$ acts faithfully on $V$, then $K$ acts freely on $\cat{O}(V,U)$ outside the basepoint.
    \begin{proof}
        The claim is clear if $\cat{O}(V,U)=0$, so we may assume there exists an isometric embedding $V\to U$, which allows us to identify $V$ with a vector subspace of $U$, equipped with the restricted inner product. By the previous lemma, we may then identify $\cat{O}(V,U)\cong \O(U)_+\smashp_{\O(U-V)}S^{U-V}$ as $(G\times H)$-spaces. By Illman's triangulation theorem, $\O(U)$ is a $\big(G\times \O(U-V)\times H\big)$-CW-complex and $S^{U-V}$ is an $\O(U-V)$-CW-complex. In any such structure, $\infty$ is a $0$-cell (being an isolated $\O(U-V)$-fixed point), so that $\O(U)_+\smashp S^{U-V}$ inherits the structure of a based $(G\times \O(U-V)\times H)$-CW-complex, ultimately giving $\cat{O}(V,U)$ the structure of a based $(G\times H)$-CW-complex; moreover, this is finite by compactness.
    
        To see that $K$ acts freely outside the basepoint on $\O(U)_+\smashp_{\O(U-V)} S^{U-V}$, it will suffice to show that it acts freely on $\O(U)/{\O(U-V)}$. For this let $A\in \O(U)$ arbitrary and let $k\in K$. As $K$ acts faithfully on $V$ by assumption, there exists a $v\in V$ such that $k^{-1}.v\not=v$, and hence $(k.A)(v)=A(k^{-1}.v)\not=A(v)$ by injectivity of $A$. Acting with any element of $\O(U-V)$ does not change the values of $A$ on $V$, so that $k.A$ and $A$ define distinct classes in $\O(U)/{\O(U-V)}$, as claimed.
    \end{proof}
\end{prop}

\subsubsection{$G$-equivariant stable homotopy theory of orthogonal $G$-spectra} For any compact Lie group $G$, we write $\Gosp{G}$ for the ($\cat{Top}_*$-enriched) category of continuous $G$-objects in $\osp$, and we will refer to its objects as \emph{orthogonal $G$-spectra}. We will want to view them through the lens of a suitable notion of \emph{$G$-equivariant weak equivalences}, which we will introduce now.

\begin{constr}
    Let $H$ be any compact Lie group and let $V$ be an $H$-representation. If $\rho\colon H\to \O(V)$ is the corresponding homomorphism, then composing with the identification $\O(V)_+\cong\cat{O}(V,V)$ from Lemma~\ref{lemma:O-mapping-spaces} yields a natural $H$-action on $X(V)$ for any orthogonal spectrum $X$. In particular, if $X$ was already an orthogonal $G$-spectrum, then this $H$-action commutes with the $G$-action, making $X(V)$ into a $(G\times H)$-space. We will be particularly interested in the case where $H\subset G$ is a closed subgroup, in which case we will usually view $X(V)$ as an $H$-space via the \emph{diagonal} action.
\end{constr}

\begin{constr}
    Fix once and for all a complete $G$-universe $\Uu_G$ and write $s(\Uu_G)$ for the poset of finite dimensional subrepresentations of $\Uu_G$. We then define for every closed subgroup $H\subset G$, $n\ge 0$, and every $G$-orthogonal spectrum $X$ its $n$-th $H$-equivariant homotopy set via
    \[
        \pi_{n}^H(X)\coloneqq\colim_{V\in s(\Uu_{G})} [S^{V\oplus\R^n},X(V)]^H_*
    \]
    where $[-{,}-]_*^H$ denotes based $H$-equivariant homotopy classes of maps and the transition map for $V\subset W$ is given by the composite
    \[\hskip-26.05pt\hfuzz=26.05pt
        [S^{V\oplus\R^n},X(V)]_*^H\xrightarrow{\;\tilde\sigma\;}[S^{V\oplus\R^n},\Omega^{W-V}X(W)]_*^H\cong[S^{W-V}\smashp S^{V\oplus\R^n},X(W)]_*^H\cong[S^{W\oplus\R^n},X(W)]_*^H
    \]
    of the adjunct structure map, the adjunction isomorphism, and the isomorphism $S^{W-V}\smashp S^{V\oplus\R^n}\cong S^{W\oplus\R^n}$ given by compactifying the evident isomorphism of real vector spaces $(W-V)\oplus (V\oplus\R^n)\cong W\oplus\R^n$.

    If $n<0$, then we define $\pi_n^H(X)$ as the analogous colimit
    \[
        \colim_{V\in s(\Uu_G)}[S^V,X(V\oplus\R^{-n})].
    \]
    Both of these become functors in $X$ via postcomposition.
\end{constr}

\begin{rk}\label{rk:universe-subgroup}
    Note that $\pi_n^H$ is independent of the choice of universe up to isomorphism. In particular, since any complete $G$-universe $\Uu_G$ is also a complete $H$-universe for any $H\subset G$ by \cite{schwede2018global}*{Remark 1.1.13}, we could as well have defined $\pi_n^H$ as a colimit over a complete $H$-universe. In other words, we also get a natural isomorphism $\pi_n^H(X)\cong\pi_n^H(\res^G_HX)$ for any orthogonal $G$-spectrum $X$ with underlying orthogonal $H$-spectrum $\res^G_HX$.
\end{rk}

\begin{definition}
    A map $f\colon X\to Y$ in $\Gosp{G}$ is called a \emph{$G$-equivariant weak equivalence} if the induced map $\pi_n^H(X)\to\pi_n^H(Y)$ is bijective for every $n\in\Z$ and $H\subset G$. 
\end{definition}

\begin{definition}
    The $\infty$-category $\mySp_G$ of \emph{$G$-equivariant spectra} is the Dwyer--Kan localization of $\Gosp{G}$ at the $G$-equivariant weak equivalences.
\end{definition}

Our next goal is to show that $\mySp_G$ is a compactly generated stable $\infty$-category and that the smash product of orthogonal spectra derives to a closed symmetric monoidal structure. While these $\infty$-categorical statements are of course already recorded in the literature, we spend some time to review the model categorical underpinnings as several of the intermediate results as well as the general strategies will be needed later  to get the $G$-global theory off the ground.

We therefore begin by describing our reference model structure on orthogonal $G$-spectra.

\begin{definition}\label{def:level}
    A map $f\colon X\to Y$ in $\Gosp{G}$ is called an \emph{$G$-equivariant level fibration} if for every finite-dimensional inner product space $V$, $f(V)\colon X(V)\to Y(V)$ is a fibration in the $\mathcal G_{G,\O(V)}$-equivariant model structure, i.e.~if $f(V)^H$ is a Serre fibration for every closed subgroup $H\subset G$ and every continuous group homomorphism $H\to O(V)$. Analogously, we define the \emph{$G$-equivariant level weak equivalences} as those $f$ such that $f(V)$ is a $\mathcal G_{G,\O(V)}$-weak equivalence for all finite-dimensional inner product spaces $V$. 

    A map is called a \emph{$G$-equivariant cofibration} if it has the left lifting property against all maps that are simultaneously $G$-equivariant level fibrations and $G$-equivariant level weak equivalences.
\end{definition}

\begin{definition}
    A map $f\colon X\to Y$ is called a \emph{$G$-equivariant fibration} if it is a $G$-equivariant level fibration and for all closed subgroups $H\subset G$ and $H$-representations $V,W$, the naturality square
    \[
        \begin{tikzcd}
            X(V)\arrow[r,"\tilde\sigma"]\arrow[d,"f(V)"'] & \Omega^WX(V\oplus W)\arrow[d,"\Omega^Wf(V\oplus W)"]\\
            Y(V)\arrow[r,"\tilde\sigma"'] & \Omega^WY(V\oplus W)
        \end{tikzcd}
    \]
    for the adjunct structure maps is a homotopy pullback in $\cat{$\bm H$-Top}_*$.

    An orthogonal $G$-spectrum $X$ is called a \emph{$G$-$\Omega$-spectrum} if the unique map $X\to 0$ is a $G$-equivariant fibration, i.e.\ the adjunct structure map $X(V)\to \Omega^WX(V\oplus W)$ is an $H$-equivariant weak equivalence for each closed subgroup $H\subset G$ and all $H$-representations $V,W$.
\end{definition}

If we want the fibrations to be part of a cofibrantly generated model structure, we have to characterize them via a lifting property. For the level fibrations this is rather straightforward, while the following set of maps will turn out to detect the homotopy pullback condition:

\begin{constr}\label{constr:generating-equiv-acyclic}
    Let $H\subset G$ be any closed subgroup and let $V$ be any $H$-representation. Then it follows from the enriched Yoneda lemma that the $\cat{Top}_*$-enriched functor 
    \begin{align*} 
        \Gosp{G} &\longrightarrow \cat{Top}_*\\ 
        X&\longmapsto X(V)^H
    \end{align*} 
    is corepresented by $\cat{O}(V,-)\smashp_HG_+$ via evaluation at $[\id_V,0;1]\in\big(\cat{O}(V,-)\smashp_HG_+\big){}^H$. 
    
    If $W$ is any further $H$-representation, then the functor $X\mapsto\big(\Omega^WX(V\oplus W)\big){}^H$ is similarly corepresented by $\big(\cat{O}(V\oplus W,-)\smashp S^W\big)\smashp_H G_+$ via evaluation at $V\oplus W$ and restriction along
    \begin{align*}
        S^W&\longrightarrow\big(\cat{O}(V\oplus W,V\oplus W)\smashp S^W\big)\smashp_H G_+\\
        w&\longmapsto [\id_{V\oplus W},0;w;1].
    \end{align*}
    It follows from the enriched Yoneda lemma, that there exists a unique map
    \[
        \kappa_{V,W}^{H,G}\colon \big(\cat{O}(V\oplus W,-)\smashp S^W\big)\smashp_H G_+\to\cat{O}(V,-)\smashp_HG_+
    \] 
    such that $\maps(\kappa_{V,W}^{H,G},X)$ corresponds under the chosen identifications with the $H$-fixed points of the adjunct structure map $X(V)\to\Omega^WX(V\oplus W)$.
\end{constr}

Note that by construction an orthogonal $G$-spectrum $X$ is a $G$-$\Omega$-spectrum if and only if it is local (in the enriched sense) with respect to the maps $\kappa_{V,W}^{H,G}$. Detecting the pullback condition in the \emph{unenriched} sense will merely require a slight modification:

\begin{constr}
    We use the mapping cylinder to factor each $\kappa_{V,W}^{H,G}$ as
    \[
        (\cat{O}(V\oplus W,-)\smashp S^W\big)\smashp_HG_+\lhook\joinrel\xrightarrow{\;\lambda_{V,W}^{H,G}\;} Z \iso \cat{O}(V,-)\smashp_HG_+ 
    \]
    and we define the set $K_G$ to consist of all the pushout product maps
    \[
        \begin{tikzcd}
            (\cat{O}(V\oplus W,-)\smashp S^W\big)\smashp_HG_+\smashp\partial D^n_+\arrow[r,"\lambda_{V,W}^{H,G}\smashp\partial D^n_+"]\arrow[d,hook] &[2.2em] Z\smashp\partial D^n_+\arrow[d]\arrow[ddr,hook, bend left=15pt]\\
            (\cat{O}(V\oplus W,-)\smashp S^W\big)\smashp_HG_+\smashp D^n_+\arrow[r]\arrow[drr,bend right=15pt,"\lambda_{V,W}^{H,G}\smashp D^n_+"'] & \arrow[from=ul,phantom,"\ulcorner"{very near end}",xshift=1.3em,yshift=.2em]P\arrow[dr,dashed,"\lambda_{V,W}^{H,G}\ppo\incl"{description}]\\
            && Z\smashp D^n_+
        \end{tikzcd}
    \]
    where $n\ge 0$, $H$ is a closed subgroup of $G$, and $V,W$ run through (representatives of isomorphisms classes of) $H$-representations.
\end{constr}

\begin{thm}
    The $G$-equivariant weak equivalences, cofibrations, and fibrations are part of a stable, proper, and topological model structure on $\Gosp{G}$. This model structure is cofibrantly generated with generating cofibrations
    \[
        I=\{\cat{O}(V,-)\smashp_H G_+\smashp(\partial D^n\hookrightarrow D^n)_+ : \text{$H\subset G$,\,$V$ an $H$-representation,\,$n\ge0$}\}
    \]
    and generating acyclic cofibrations
    \[
        J\hskip0pt minus 1pt=\hskip 0pt minus 1pt\{\cat{O}(V,-)\smashp_H G_+\smashp(\partial D^n\hskip0pt minus 2pt\hookrightarrow\hskip0pt minus 2pt D^n)_+\hskip0pt minus 1pt :\hskip 0pt minus 1pt \text{$H\subset G$,\,$V$ an $H$-representation,\,$n\ge0$}\}\cup K_G.
    \]
    We call this the \emph{$\bm G$-equivariant model structure}.
    \begin{proof}
        All of these statements apart from the description of the generating (acyclic) cofibrations are part of \cite{DHLPS2019Proper}*{Theorem~1.37}. The description of the generating (acyclic) cofibrations is then noted in the beginning of part (iv) of the proof of \emph{loc.\ cit.}
    \end{proof}
\end{thm}

\begin{lemma}\label{lemma:Sigma-oo-left-Quillen}
    The adjunction 
    \[
        \Sigma^\infty\coloneqq\cat{O}(0,-)\smashp-\colon\cat{$\bm G$-Top}_*\rightleftarrows\Gosp{G}:\Omega^\infty\coloneqq\ev_0
    \]
    is a Quillen adjunction.
    \begin{proof}
        It is clear that the left adjoint sends generating (acyclic) cofibrations to generating (acyclic) cofibrations, so it is left Quillen.
    \end{proof}
\end{lemma}

\subsubsection{Homotopy pushouts} Similarly to various model structures considered before, the equivariant model structure on orthogonal $G$-spectra has rather few cofibrations, so that it is often hard to get concrete control over the factorization of a given map into a cofibration followed by a weak equivalence. In order to nevertheless be able to compute homotopy pushouts we introduce:

\begin{definition}
    A map $f\colon X\to Y$ in $\Gosp{G}$ is called an \emph{unbased level h-cofibration} if each $f(V)\colon X(V)\to Y(V)$ is an h-cofibration in $\cat{$\bm{(\O(V)\times G)}$-Top}$.
\end{definition}

\begin{lemma}\label{lemma:equiv-spectra-homotopy-po}
    Pushouts along unbased level h-cofibrations are homotopy pushouts in the $G$-equivariant model structure.
    \begin{proof}
        Consider a pushout as depicted on the left, where $i$ is an unbased level h-cofibration, and pick a factorization of $f$ as a cofibration $j$ followed by a \emph{level} weak equivalence $g$, leading to the iterated pushout depicted on the right:
        \[
            \begin{tikzcd}
                A\arrow[r, "f"]\arrow[dr,phantom,"\ulcorner"{very near end}]\arrow[d,"i"'] & B\arrow[d] &[3em] A\arrow[dr,phantom,"\ulcorner"{very near end}]\arrow[r, "j"]\arrow[d,"i"'] & X\arrow[dr,phantom,"\ulcorner"{very near end}]\arrow[d,"i'"]\arrow[r,"g","\sim"'] & B\arrow[d]\\
                C\arrow[r] & D & C\arrow[r]& Y\arrow[r,"h"'] & D\rlap.
            \end{tikzcd}
        \]
        If $H\subset G$ and $V$ is any $H$-representation, then $i(V)$ is in particular an h-cofibration in $\cat{$\bm H$-Top}$, whence so is $i'(V)$. As $g(V)$ is an $H$-equivariant weak equivalence, so is $h(V)$ by Lemma~\ref{lemma:h-cof-compute-po}, i.e.\ $h$ is again a $G$-equivariant level weak equivalence. This in particular shows that the original square is a homotopy pushout.
    \end{proof}
\end{lemma}

\subsubsection{Stability} For $G=1$, the above is a model of ordinary stable homotopy theory, and in particular the corresponding $\infty$-category $\mySp$ is stable in the sense that the suspension-loop adjunction is an equivalence. For general $G$, $\mySp_G$ satisfies the following genuine version of stability (which will later precisely correspond to representation stability of a certain global $\infty$-category $\ul\mySp$ of equivariant spectra):

\begin{prop}[See \cite{mandell-may}*{Proposition~III.3.9 and Theorem~III.3.11}]\label{prop:rep-spheres-invertible-equiv-model}
    Let $V$ be any $G$-representation. Then the endofunctors $S^V\smashp-$ and $\Omega^V$ of $\Gosp{G}$ preserve $G$-equivariant weak equivalences, and the induced endofunctors of $\mySp_G$ are mutually inverse equivalences.\qed
\end{prop}

\begin{cor}
    The $\infty$-category $\mySp_G$ is stable.
    \begin{proof}
        As the underlying $\infty$-category of a model category, $\mySp_G$ is complete and cocomplete, with the localization preserving initial and final objects. Thus, $\mySp_G$ is pointed. Observe now that the $\infty$-categorical suspension functor on $\mySp_G$ is induced by the suspension functor on the pointset level (e.g.~since the left Quillen bifunctor $\cat{Top}_*\times\Gosp{G}\to\Gosp{G}$ describing the tensoring descends to a functor $\Spc_*\times\mySp_G\to\mySp_G$ preserving colimits in each variable and restricting to the identity on $\{S^0\}\times\mySp_G$, and there is only one such functor). Thus, the claim follows by specializing the previous proposition to the 1-dimensional trivial $G$-representation.
    \end{proof}
\end{cor}

In particular, $\mySp_G$ is additive. Since each of the functors $\pi_n^H\colon\mySp_G\to\cat{Sets}$ preserves finite products by direct inspection, it therefore admits a unique lift to $\cat{Ab}$ (which would of course also be easy to describe directly).

\begin{prop}\label{prop:equiv-sp-compact-gen}
    The stable $\infty$-category $\mySp_G$ is compactly generated by the collection of suspension spectra $\Sigma^\infty G/H_+$ for all closed subgroups $H\subset G$. In particular, it is presentable.
    \begin{proof}
        By \cite{HA}*{Proposition 1.4.4.1 and Corollary 1.4.4.2} it only remains to show that the spectra $\Sigma^\infty G/H_+$ form a set of compact generators in the triangulated sense, which is an instance of \cite{DHLPS2019Proper}*{Corollary~1.59}.
    \end{proof}
\end{prop}

\subsubsection{Functoriality} In order to be able to later assemble the individual $\infty$-categories $\mySp_G$ into a global $\infty$-category $\ul\mySp$, we should once more discuss their relationship as the group $G$ varies.

\begin{prop}[See \cite{DHLPS2019Proper}*{Theorem~1.66(ii)}]\label{prop:restr-left-Quillen-G-Sp}
    Let $f\colon G\to G'$ be a homomorphism of compact Lie groups. Then
    \[
        f^*\colon\Gosp{G'}\rightleftarrows\Gosp{G}:f_*
    \]
    is a Quillen adjunction for the equivariant model structures, and hence in particular induces an adjunction $\cat{L}f^*\colon\mySp_{G'}\rightleftarrows\mySp_G\noloc\cat{R}f_*$.\qed
\end{prop}

\begin{prop}\label{prop:i!-left-Quillen-equiv-sp}
    Let $i\colon G\to G'$ be an \emph{injective} homomorphism of compact Lie groups. Then
    \[
        i_!\colon\Gosp{G}\rightleftarrows\Gosp{G'}\noloc i^*
    \]
    is a Quillen adjunction in which both adjoints are homotopical. In particular, they descend to an adjunction $i_!\colon\mySp_G\rightleftarrows\mySp_{G'}\noloc i^*$.
    \begin{proof}
        \cite{DHLPS2019Proper}*{Theorem~1.66(i)} shows that $i^*$ is homotopical right Quillen, while $i_!$ is homotopical by \cite{schwede2018global}*{Corollary~3.2.21}.
    \end{proof}
\end{prop}

\begin{warn}
    In contrast to all other models of equivariant or global homotopy theory considered in this paper, restriction along a non-injective homomorphism $f\colon G\to G'$ does \emph{not} send $G'$-equivariant weak equivalences to $G$-equivariant ones. Accordingly, some care will be needed in the next subsection to construct the global $\infty$-category $\ul\mySp$.
\end{warn}

\subsubsection{Smash product symmetric monoidal structure} Next, we recall the usual smash product of orthogonal spectra.

\begin{constr}
    The indexing category $\cat{O}$ for orthogonal spectra admits a $\cat{Top}_*$-enriched symmetric monoidal structure given on objects by direct sum and on morphisms by
    \begin{align*}
        \cat{O}(V_1,W_1)\smashp\cat{O}(V_2,W_2)&\longrightarrow\cat{O}(V_1\oplus V_2,W_1\oplus W_2)\\
        \big[[A_1,w_1],[A_2,w_2]]&\longmapsto [A_1\oplus A_2, (w_1,w_2)]
    \end{align*}
    (with the convention that $(\infty,w_2)=\infty$ and $(w_1,\infty)=\infty$). The unit is given by the zero representation, and the associativity, unitality, and symmetry isomorphisms are induced from the ones for the cocartesian symmetric monoidal structure on $\R$-vector spaces. By enriched Day convolution, this gives us a closed and $\cat{Top}_*$-enriched symmetric monoidal structure on $\osp$, which we call the \emph{smash product}. The internal homs are given by $F(X,Y)(V)=\maps(\sh^VX,Y)$ where $\sh^V$ is given by restriction along $V\oplus-\colon\cat{O}\to\cat{O}$, and the functoriality of $F(X,Y)$ in $V$ is via the functoriality of ${-}\oplus{-}$ in the first variable.
    
    Pulling through $G$-actions, we then more generally get the analogous structure on $\Gosp{G}$ for every compact Lie group $G$.
\end{constr}

\begin{prop}[See \cite{DHLPS2019Proper}*{Proposition~1.43}]\label{prop:equiv-smash-sp-left-Quillen-bifunctor}
    The smash product defines a left Quillen bifunctor $\Gosp{G}\times\Gosp{G}\to\Gosp{G}$ with respect to the equivariant model structures.\qed
\end{prop}

As a consequence, we see that smashing with a cofibrant orthogonal $G$-spectrum preserves weak equivalences between cofibrant orthogonal $G$-spectra. However, the smash product is actually homotopical under significantly weaker assumptions. In order to state this result, we will first need some terminology.

\begin{constr}
    Let $\cat{O}_{\le n}\subset\cat{O}$ denote the full subcategory spanned by the inner product spaces of dimension ${}\le n$, and write $\osp_{\le n}$ for the category of $\cat{Top}_*$-enriched functors $\cat{O}_{\le n}\to\cat{Top}_*$. Then the restriction $\osp\to\osp_{\le n}$ admits a left adjoint $i_!$ via enriched left Kan extension. We write $\sk^{n}\coloneqq i_!i^*$, which comes with a natural transformation $\epsilon\colon\sk^n\to\id_{\osp}$ given by the counit of the adjunction. We lift $\sk^n$ and $\epsilon$ to $\Gosp{G}$ by pulling through the $G$-actions as before.
\end{constr}

\begin{definition}\label{def:flat}
    Given an orthogonal $G$-spectrum $X$ and an $n\ge0$, the $n$-th \emph{latching map} $L_nX\to X(\R^n)$ is the map $\epsilon(\R^n)\colon\sk^{n-1}(X)(\R^n)\to X(\R^n)$. We say that $X$ is \emph{flat} if the $n$-th latching map is a cofibration in the $\All$-model structure on $\cat{$\bm{(\O(n)\times G)}$-Top}_*$ for every $n\ge0$.
\end{definition}

\begin{ex}\label{ex:equiv-cofibrant-flat}
    Every orthogonal $G$-spectrum that is cofibrant in the equivariant model structure is flat. This follows at once from the characterization of the $G$-equivariant cofibrations in \cite{DHLPS2019Proper}*{Definition~1.29}, see also Remark~1.38 of \emph{op.\ cit.}
\end{ex}

\begin{thm}[Stolz]\label{thm:equiv-flatness}
    If $X$ is a flat orthogonal $G$-spectrum, then the functor $X\smashp{-}\colon\Gosp{G}\to\Gosp{G}$ preserves $G$-equivariant weak equivalences. Moreover, for \emph{any} orthogonal $G$-spectrum $Y$ the functor ${-}\smashp Y$ preserves $G$-equivariant weak equivalences between flat orthogonal $G$-spectra.
    \begin{proof}
        For a published account of the first statement see \cite{schwede2018global}*{Theorem 3.5.10}. For the second statement, we let $X_1\to X_2$ be any $G$-equivariant weak equivalence between flat orthogonal $G$-spectra and we pick a cofibrant replacement $Y'\to Y$. Then all maps in the induced square 
        \[
            \begin{tikzcd}
                X_1\smashp Y'\arrow[r]\arrow[d] & X_2\smashp Y'\arrow[d]\\
                X_1\smashp Y\arrow[r] & X_2\smashp Y
            \end{tikzcd}
        \]
        except possibly the bottom one are $G$-equivariant weak equivalences by the first part; thus, also the bottom map is a $G$-equivariant weak equivalence by 2-out-of-3.
    \end{proof}
\end{thm}

\subsection{The universal property of equivariant spectra} As the main result of this section, we will now describe the free representation stable equivariantly presentable global $\infty$-category $\ul\Sp_{\Orb\triangleright\Glo}^\textup{RepSph}$ in terms of the $\infty$-categories $\mySp_G$ of equivariant spectra. We begin by explaining how the latter assemble into a global $\infty$-category:

\begin{constr}
    Consider the full subcategory of $\Ntop(\osp^\flat)$ given in degree $G$ by the cofibrant orthogonal $G$-spectra. By Proposition~\ref{prop:restr-left-Quillen-G-Sp}, this is indeed a global subcategory and all restriction functors preserve equivariant weak equivalences when restricted to this subcategory. We may therefore Dwyer--Kan localize, yielding a global $\infty$-category $\ul\mySp$.
\end{constr}

By Lemma~\ref{lemma:fgt-simpl-structure-cof}, we may equivalently define $\ul\mySp$ on $0$-cells and $1$-cells after forgetting the topological enrichment. Combining this with the existence of (functorial) cofibrant replacements, we therefore arrive at the following incoherent description of $\ul\mySp$: it sends an object $\BGcat{G}\in\Glo$ to the $\infty$-category $\mySp_G$ of $G$-equivariant spectra and a functor $\BGcat{G}\to\BGcat{G'}$ corresponding to a homomorphism $f\colon G\to G'$ to the left derived functor $\cat{L}f^*\colon\mySp_{G'}\to\mySp_G$.

\begin{lemma}
    The global $\infty$-category $\ul\mySp$ is fiberwise presentable.
    \begin{proof}
        By Proposition~\ref{prop:equiv-sp-compact-gen} each $\mySp_G$ is presentable, while Proposition~\ref{prop:restr-left-Quillen-G-Sp} shows that all restriction functors are left adjoints.
    \end{proof}
\end{lemma}

Below we will prove:

\begin{thm}\label{thm:equiv-main}
    The global $\infty$-category $\ul\mySp$ is the free equivariantly cocomplete and representation stable global $\infty$-category generated by the global section $\mathbb S$, i.e.\ for any other such $\Cc$ evaluation at $\mathbb S$ defines an equivalence
    \[
        \ul\Fun^\textup{$\Orb$-cc}_\Glo(\ul\mySp,\Cc)\iso\Cc.
    \]
\end{thm}

Our argument will proceed by instead establishing a universal property for a certain symmetric monoidal refinement $\ul\mySp$ and then appealing to the general theory established in the previous section. 

To construct the symmetric monoidal structure, recall first that by Proposition~\ref{prop:equiv-smash-sp-left-Quillen-bifunctor} the symmetric monoidal structure on $\Ntop(\osp^\flat)$ induced by the smash product restricts to the full subcategory of cofibrant objects, and this restricted tensor product preserves weak equivalences in each variable. This then upgrades $\ul\mySp$ to a symmetric monoidal global $\infty$-category $\ul\mySp^\otimes$ such that each individual $\mySp_G$ is presentably symmetric monoidal. We will now explain how to assemble the suspension spectrum functors $\cat{L}\Sigma^\infty\colon\myS_{G,*}\to\mySp_G$ into a global symmetric monoidal functor:

\begin{constr}
    The topologically enriched symmetric monoidal functor $\Sigma^\infty\colon\cat{Top}_*\to\osp$ induces a symmetric monoidal functor $\Ntop(\cat{Top}_*^\flat)\to\Ntop(\osp^\flat)$. By Lemma~\ref{lemma:Sigma-oo-left-Quillen}, this restricts to a homotopical functor between the corresponding subcategories of cofibrant objects and hence induces
    \[
        \cat{L}\Sigma^\infty\colon\ul\myS_{*}^\otimes\simeq(\ul\myS_*^\text{cof})^\otimes\to\ul\mySp^\otimes.
    \]
\end{constr}

As a direct consequence of work of Gepner--Meier, this construction enjoys the following `pointwise' universal property:

\begin{prop}\label{prop:equiv-sp-pw-univ-property}
    Let $f\colon T\to\Glo$ be any functor, and let $\Cc^\otimes$ be a fiberwise presentably symmetric monoidal $T$-$\infty$-category. Then the map
    \[
        \hom_{\Fun(T^\op,\CAlg(\PrL))}(f^*\ul\mySp^\otimes,\Cc^\otimes)\to
        \hom_{\Fun(T^\op,\CAlg(\PrL))}(f^*\ul\myS_*^\otimes,\Cc^\otimes) 
    \]
    induced by $\cat{L}\Sigma^\infty$ is an inclusion of path-components, and its image consists precisely of those functors that send representation spheres to invertible objects.
    \begin{proof}
        For any compact Lie group $G$, \cite{gepnermeier2020equivTMF}*{Proposition C.9} shows that $\cat{L}\Sigma^\infty\colon\myS_{G,*}^\otimes\to\mySp_G^\otimes$ is the initial example of a map in $\CAlg(\PrL)$ inverting representation spheres, i.e.~it exhibits $\smash{\mySp_G^\otimes}$ as a left adjoint object with respect to the functor $i\colon\PrL\to\Ar(\CAlg(\PrL))$ sending a symmetric monoidal $\infty$-category to the inclusion of the full subcategory of tensor invertible objects. Thus, $f^*\cat{L}\Sigma^\infty$ exhibits $f^*\ul\mySp^\otimes$ as a left adjoint object with respect to $\Fun(T^\op,i)$, which is exactly the statement of the proposition. 
    \end{proof} 
\end{prop}

It remains to understand how the above pointwise stabilization interacts with the notion of equivariant cocompleteness.

\begin{lemma}\label{lemma:Sigma-oo-equiv-left-adjointable}
    Let $i\colon G\to G'$ be an injective homomorphism of compact Lie groups. Then $i^*\colon\mySp_{G'}\to\mySp_G$ admits a left adjoint $i_!$, and the Beck--Chevalley map $i_!\cat{L}\Sigma^\infty\to\cat{L}\Sigma^\infty\cat{L}i_!$ of functors $\myS_{G,*}\to\mySp_{G'}$ is an equivalence.
    \begin{proof}
        Proposition~\ref{prop:i!-left-Quillen-equiv-sp} shows that we have a homotopical Quillen adjunction $i_!\colon\Gosp{G}\rightleftarrows\Gosp{G'}\noloc i^*$, proving the existence of the left adjoint $i_!$ of $i^*=\cat{L}i^*\colon\mySp_{G'}\to\mySp_G$. Since also $i_!\colon\cat{$\bm G$-Top}_*\to\cat{$\bm{G'}$-Top}$ and the suspension spectrum functors are left Quillen (see Lemmas~\ref{lemma:res-pted-Quillen} and~\ref{lemma:Sigma-oo-left-Quillen}, respectively), it then suffices to observe for the Beck--Chevalley condition that the \emph{underived} Beck--Chevalley map $i_!\Sigma^\infty\to\Sigma^\infty i_!$ is an isomorphism.
    \end{proof}
\end{lemma}

\begin{prop}\label{prop:equiv-sp-equiv-cc}
    The symmetric monoidal global $\infty$-category $\ul\mySp^\otimes$ is equivariantly presentably symmetric monoidal.
    \begin{proof}
        Combining Proposition~\ref{prop:equiv-sp-pw-univ-property} for $\Orb\hookrightarrow\Glo$ with \cite{twisted-ambidexterity}*{Proposition~2.46} shows that $\ul\mySp^\otimes|_{\Orb}$ is $\Orb$-presentably symmetric monoidal. Factoring a general map in $\Glo$ as a surjective homomorphism followed by a subgroup inclusion and appealing to Proposition~\ref{prop:pb-surj-glo}, it then only remains to verify that the Beck--Chevalley map $j_!q^*\to p^*i_!$ is an equivalence for any pullback
        \[
            \begin{tikzcd}
                H'\arrow[r,"j"]\arrow[d,"{q}"']\arrow[dr,pullback] & {G'}\arrow[d,"{p}"]\\
                {H}\arrow[r,hook,"{i}"'] & {G}
            \end{tikzcd}
        \]
        in the $1$-category of compact Lie groups, where $i$ is the inclusion of a subgroup and $p$ is surjective. As both source and target of the Beck--Chevalley map are left adjoint functors of presentable stable $\infty$-categories, it suffices to check the Beck--Chevalley condition on the compact generators from Proposition~\ref{prop:equiv-sp-compact-gen}, and hence in particular after precomposing with $\cat{L}\Sigma^\infty$. As Beck--Chevalley maps compose, we then have a commutative diagram of Beck--Chevalley maps
        \[
            \begin{tikzcd}
                j_!q^*\cat{L}\Sigma^\infty\arrow[r,"\sim"]\arrow[d,"\BC_!\circ\cat{L}\Sigma^\infty"'] & j_!\cat{L}\Sigma^\infty q^*\arrow[r,"\BC_!"] &  \cat{L}\Sigma^\infty j_!q^*\arrow[d,"\cat{L}\Sigma^\infty\BC_!"] \\
                p^*i_!\cat{L}\Sigma_!\arrow[r,"\BC_!"']&p^*\cat{L}\Sigma^\infty i_!\arrow[r,"\sim"'] & \cat{L}\Sigma^\infty p^*i_!\rlap,
            \end{tikzcd}
        \]
        where all unlabelled equivalences are given by naturality. We want to show that the left-hand certical map is an equivalence. By the previous lemma (and using that both $i$ and $j$ are maps in $\Glo$ as opposed to $\PSh(\Glo)$, by virtue of the explicit description of the pullback), the horizontal Beck--Chevalley maps are equivalences, and so is the right-hand vertical map by Corollary~\ref{cor:equiv-pointed}. The claim follows by 2-out-of-3.
    \end{proof}
\end{prop}

\begin{prop}\label{prop:sigma-oo-pres-refl-equiv-cc}
    The symmetric monoidal functor $\cat{L}\Sigma^\infty\colon\ul\myS_*^\otimes\to\ul\mySp^\otimes$ is equivariantly cocontinuous. Moreover, if $\Cc^\otimes$ is equivariantly presentably symmetric monoidal, then a fiberwise cocontinuous symmetric monoidal functor $F\colon\ul\mySp^\otimes\to\Cc^\otimes$ is equivariantly cocontinuous if and only if $F\circ\cat{L}\Sigma^\infty_+\colon\ul\myS_*^\otimes\to\Cc^\otimes$ is so.
    \begin{proof}
        To see that $\cat{L}\Sigma^\infty$ is equivariantly cocontinuous, we observe that each $\cat{L}\Sigma^\infty\colon\myS_{G,*}\to\mySp_G$ is cocontinuous (being the left derived functor of a left Quillen functor), while the Beck--Chevalley condition for left adjoints to restrictions along injective homomorphisms is the content of Lemma~\ref{lemma:Sigma-oo-equiv-left-adjointable}. This completes the proof of the first statement and the `only if' part of the second one.

        To prove the remaining direction, we left $F\colon\ul\mySp^\otimes\to\Cc^\otimes$ be a fiberwise cocontinuous symmetric monoidal functor such that $F\circ\cat{L}\Sigma^\infty$ is equivariantly cocontinuous; we have to show that  the Beck--Chevalley map $i_!F\to Fi_!$ is invertible for every injective $i\colon G\to G'$. Note that the presentably symmetric monoidal $\infty$-categories $\Cc(G)$ and $\Cc(G')$ are stable as they receive symmetric monoidal left adjoints from presentably symmetric monoidal stable $\infty$-categories (namely, $F_G$ and $F_G'$); thus, the Beck--Chevalley condition follows by exactly the same reduction argument as in the proof of the previous proposition.
    \end{proof}
\end{prop}
    
\begin{proof}[Proof of Theorem~\ref{thm:equiv-main}]
    Applying Proposition~\ref{prop:equiv-sp-pw-univ-property} for $f=\id$ shows that $\cat{L}\Sigma^\infty\colon\ul\myS_{*}^\otimes\to\ul\mySp^\otimes$ is the initial example of a map in $\Fun(\Glo^\op,\CAlg(\PrL))$ with source $\ul\myS_{*}^\otimes$ and inverting all representation spheres. Combining this with Propositions~\ref{prop:equiv-sp-equiv-cc} and~\ref{prop:sigma-oo-pres-refl-equiv-cc} shows that it is also the initial such map in $\CAlg(\Pr^{\Orb}_{\Glo})$. Thus, Corollary~\ref{cor:SpR-sym-mon-univ-prop} implies that there is an equivalence $\smash{\ul\Sp^\text{RepSph}_{\Orb\triangleright\Glo}\simeq\ul\mySp}$ under $\ul\Spc_{\Orb\triangleright\Glo,*}\simeq\myS_*$ (in fact, the corollary even provides a symmetric monoidal equivalence). The claim then follows from the defining universal property of $\smash{\ul\Sp^\text{RepSph}_{\Orb\triangleright\Glo}}$.
\end{proof}

In the same way one deduces from Corollary~\ref{cor:SpR-initial}:

\begin{cor}\label{cor:equiv-spectra-smash}
    The global $\infty$-category $\ul\mySp$ admits a unique equivariantly presentably symmetric monoidal structure with unit $\mathbb S$, and this is the symmetric monoidal structure induced by the smash product constructed above. Moreover, the resulting global symmetric monoidal $\infty$-category is the initial equivariantly presentably symmetric monoidal representation stable global $\infty$-category. \qed
\end{cor}

\section{The universal property of global spectra}\label{sec:global-spectra}
In this section, we will prove the universal property of global stable homotopy theory. In more detail, we explain in §\ref{subsec:gl-spectra} how looking at orthogonal $G$-spectra through a finer notion of weak equivalence than the $G$-equivariant weak equivalences considered in the previous section gives rise to a \emph{$G$-global model structure}, which for $G=1$ recovers Schwede's global model structure \cite{schwede2018global}. In §\ref{subsec:L-spectra} we introduce an alternative model of $G$-global stable homotopy theory using orthogonal spectrum objects in orthogonal $G$-spaces, which will be crucial for establishing the universal property; the proof that these two models are indeed equivalent occupies all of §\ref{subsec:comp-spectra}.
In §\ref{subsec:gl-oo-gl-sp} we show that the $\infty$-categories $\mySp_\textup{$G$-gl}$ of $G$-global spectra assemble into a representation stable globally presentable global $\infty$-category $\ul\mySp_\gl$, and in §\ref{subsec:gl-sp-univ} we finally prove that this is in fact the universal such category.

\subsection{The \texorpdfstring{$\bm G$}{G}-global model structure on orthogonal \texorpdfstring{$\bm G$}{G}-spectra}\label{subsec:gl-spectra}
In this subsection we will introduce our reference model of $G$-global stable homotopy theory in terms of a \emph{$G$-global model structure} on orthogonal $G$-spectra, generalizing Schwede's global model structure on orthogonal spectra \cite{schwede2018global}*{Theorem 4.3.17}. The $G$-global model structure we consider here was first constructed in the master's thesis of Vincent Grande \cite{grande-thesis}; as this thesis is not publicly available, we will give a new proof of the existence of this model structure for the reader's convenience.

\medskip
\subsubsection{The $G$-global level model structure.} We again start with a suitable level model structure:

\begin{prop}\label{prop:global-lvl-GSp}
    Let $G$ be a compact Lie group. There exists a (unique) model structure on $\Gosp{G}$ in which a map $f\colon X\to Y$ is a weak equivalence or fibration if and only if for every finite-dimensional inner product space $V$, $f(V)\colon X(V)\to Y(V)$ is a weak equivalence or fibration, respectively, in the $\mathcal G_{\O(V),G}$-equivariant model structure\footnote{We emphasize that this is different from the $\mathcal G_{G,\O(V)}$-model structure used to the define the \emph{$G$-equivariant} level weak equivalences.} on $\cat{$\bm{(\O(V)\times G)}$-Top}$. We call this the \emph{$\bm G$-global level model structure} and its weak equivalences the \emph{$\bm G$-global level weak equivalences}. If $f\colon X\to Y$ is any cofibration in the model structure, then $f(V)$ is in particular a $\mathcal G_{\O(V),G}$-cofibration for every inner product space $V$, and likewise for the acyclic cofibrations.
    
    Moroever, the $G$-global level model structure is proper, topological, and cofibrantly generated with generating cofibrations
    \[
        I=\{\cat{O}(V,-)\smashp_\phi G_+\smashp(\partial D^n\hookrightarrow D^n)_+ : \text{$V\in\cat{O}$, $H\subset \O(V)$, $\phi\colon H\to G$, $n\ge0$}\}
    \]
    and generating acyclic cofibrations
    \[
        J=\{\cat{O}(V,-)\smashp_\phi G_+\smashp(D^n\hookrightarrow D^n\times I)_+ : \text{$V\in\cat{O}$, $H\subset \O(V)$, $\phi\colon H\to G$, $n\ge0$}\}.
    \]
\end{prop}

For the proof of the proposition we will need:

\begin{lemma}\label{lemma:O-cobc-lQ}
    Let $V,W$ be finite dimensional inner product spaces. Then
    \[
    \cat{O}(V,W)\smashp_{\O(V)}{-}\colon\cat{$\bm{(\O(V)\times G)}$-Top}_*\to \cat{$\bm{(\O(W)\times G)}$-Top}_*
    \]
    is left Quillen for the $\mathcal G_{\O(V),G}$-model structure on the source and the $\mathcal G_{\O(W),G}$-model structure on the target.
    \begin{proof}
        It is clear that this functor is a left adjoint. By the explicit description of the generating (acyclic) cofibrations, it will then suffice to show that ${\cat{O}(V,W)\smashp_{\O(V)}(\O(V)\times G)_+/\Gamma_{H,\phi}}\cong \cat{O}(V,W)\smashp_\phi G_+$ is cofibrant in the $\mathcal G_{\O(W),G}$-model structure on $\cat{$\bm{(\O(W)\times G)}$-Top}_*$ for any $H\subset \O(V)$ and $\phi\colon H\to G$.

        By Proposition~\ref{prop:O-mapping-CW}, $\cat{O}(V,W)$ is an $(\O(W)\times H)$-CW complex with free $H$-action outside the basepoint, hence cofibrant in the $\mathcal G_{\O(W),H}$-model structure. To complete the proof it then suffices to observe that
        \[
            {-}\smashp_\phi G_+=\phi_!\colon\cat{$\bm{(\O(W)\times H)}$-Top}_*\to \cat{$\bm{(\O(W)\times G)}$-Top}_*
        \]
        is left Quillen for the graph model structures by Example~\ref{ex:graph-restr-rQ}.
    \end{proof}
\end{lemma}

\begin{proof}[Proof of Proposition~\ref{prop:global-lvl-GSp}]
    Observe first that if $f$ is a generating cofibration, then $f(V)$ is a $\mathcal G_{\O(V),G}$-cofibration by the previous lemma. As evaluation at $V$ is cocontinuous, we then see that the same holds more generally whenever $f$ is a retract of a relative $I$-cell complex; in particular, once we have set up the model structure, this will prove that cofibrations are levelwise cofibrations. In the same way, we see that retracts of relative $J$-cell complexes are levelwise acyclic cofibrations.

    For existence of the model structure, we want to apply the Crans--Kan transfer criterion \cite{hirschhorn-book}*{Theorem~11.3.2} for
    \begin{equation}\label{eq:Sp-transfer-RA}
        \Gosp{G}\xrightarrow{(\ev_V)_V}\prod_V\cat{$\bm{(\O(V)\times G)}$-Top}
    \end{equation}
    where the product runs over a set of isomorphism classes of finite-dimensional inner product spaces (e.g.\ just over the standard $\R^n$'s), and we equip the factor $\cat{$\bm{(\O(V)\times G)}$-Top}$ with the $\mathcal G_{\O(V),G}$-model structure. Note that $(\ref{eq:Sp-transfer-RA})$ is right adjoint to the functor sending a family of maps $f_V$ to the coproduct $\bigvee_{V}\cat{O}(V,-)\smashp_{\O(V)} (f_V)_+$. Moreover, the target of $(\ref{eq:Sp-transfer-RA})$ is cofibrantly generated, with a set of generating (acyclic) cofibrations given by those tuples $(f_V)_V$ where one $f_V$ is a generating (acyclic) cofibration while all other are the identity of the initial object. Thus, with respect to the standard choices of generating (acyclic) cofibrations for orthogonal $G$-spaces, the images of generating cofibrations or generating acyclic cofibrations of $\prod_V\cat{$\bm{(\O(V)\times G)}$-Top}$ are precisely the sets $I$ and $J$ from the theorem statement. It then remains to show:
    \begin{enumerate}
        \item The set $I$ permits the small object argument (i.e.\ the source of any map in $I$ is small with respect to transfinite compositions of relative $I$-cell complexes).
        \item The set $J$ permits the small object argument.
        \item Any relative $J$-cell complex is a $G$-global level weak equivalence.
    \end{enumerate}
    We begin by proving (1); the proof of (2) will then be completely analogous. By construction, any map in $I$ is of the form $\cat{O}(V,-)\smashp_{\O(V)}i$ for some inner product space $V$ and some generating cofibration $i$ of the $\mathcal G_{\O(V),G}$-equivariant model structure from Proposition~\ref{prop:equiv-model-structure}. By direct inspection, the source of any such $i$ is small with respect to transfinite compositions of closed embeddings. On the other hand, $\cat{O}(V,-)\smashp_{\O(V)}{-}$ is left adjoint to $\ev_V$, and the latter is cocontinuous and sends relative $I$-cell complexes to closed embeddings as a consequence of the previous lemma. Thus, the source of $\cat{O}(V,-)\smashp_{\O(V)}i$ is small with respect to transfinite composition of $I$-cell complexes as claimed.

    For the proof of $(3)$ on the other hand, it suffices to observe once more that evaluation at any $V$ sends maps in $J$ to $\mathcal G_{\O(V),G}$-acyclic cofibrations. This completes the proof that the model structure exists and is cofibrantly generated by $I$ and $J$; moreover, as announced in the beginning, we have also shown that cofibrations are in particular levelwise cofibrations. Properness is therefore immediate from properness of the $\mathcal G_{\O(V),G}$-model structure. Similarly, we conclude that the cotensoring over $\cat{Top}$ is a right Quillen bifunctor, i.e.\ the model structure is topological.
\end{proof}

For later use we record:

\begin{lemma}\label{lemma:cofibrant-G-gl-are-flat}
    Every orthogonal $G$-spectrum that is cofibrant in the $G$-global level model structure is flat in the sense of Definition~\ref{def:flat}.
    \begin{proof}
        Consider the class $\mathscr F$ of \emph{flat cofibrations}, i.e.\ maps $f\colon X\to Y$ such that the induced map $X(\R^n)\amalg_{L_nX}L_nY\to Y(\R^n)$ is an $(\O(n)\times G)$-equivariant cofibration for all $n\ge0$, where $L_nX, L_nY$ denote the nth latching objects of $X$ and $Y$ as defined in Definition~\ref{def:flat}. We will show that $\mathscr F$ contains all $G$-global cofibrations, which will then in particular imply the lemma.        
        
        For this we observe that $\mathscr F$ is closed under pushout, transfinite composition, and retracts: for example, using that both evaluation at $\R^n$ and the restriction functor $i^*$ from the previous construction also have \emph{right} adjoints, one can identify $\mathscr F$ with the class of maps having the right lifting property against a certain class of maps. It therefore suffices to show that the generating cofibrations are contained in $\mathscr F$.

        As $\cat{O}_{\le n-1}\hookrightarrow\cat{O}_{\le n}$ is fully faithful, the counit $\epsilon\colon\sk^{n-1}X\to X$ is an isomorphism if $X$ is left Kan extended from $\cat{O}_{\le n-1}$ and in particular for $X=\cat{O}(V,-)\smashp_\phi G_+\smashp K_+$ with $\dim(V)<n$. Thus, also the latching map is an isomorphism; we conclude that for the generating cofibration \[X=\cat{O}(V,-)\smashp_\phi G_+\smashp\partial D^m\hookrightarrow \cat{O}(V,-)\smashp_\phi G_+\smashp D^m=Y,\] the map $X(\R^n)\amalg_{L_nX} L_nY\to Y(\R^n)$ is an isomorphism, and hence in particular an equivariant cofibration, whenever $n>\dim(V)$. If $n\le\dim(V)$ on the other hand, then $\sk^{n-1}\cat{O}(V,-)=0$, and hence $L_n\big(\cat{O}(V,-)\smashp_HG_+\smashp K\big)=0$. For $X\to Y$ the generating cofibration from before, we can then identify $X(\R^n)\amalg_{L_nX}L_nY\to Y(\R^n)$ with $\cat{O}(V,\R^n)\smashp_\phi G_+\smashp(\partial D^m\hookrightarrow D^m)$. If $n<\dim(V)$, both sides are zero, so this is even an isomorphism; on the other hand, for $n=\dim(V)$, this can be further identified via Lemma~\ref{lemma:O-mapping-spaces} with $(\O(n)\times G)_+/\Gamma_{H,\phi}\smashp(\partial D^m\hookrightarrow D^m)_+$, which is one of the generating cofibrations.
    \end{proof}
\end{lemma}

\subsubsection{The $G$-global model structure} The model structure we are actually after will be obtained from the $G$-global level model structure by Bousfield localization. We can describe both the local objects and the weak equivalences of this localization explicitly:

\begin{definition}
    An orthogonal $G$-spectrum $X$ is called a \emph{$G$-global $\Omega$-spectrum} if for every compact Lie group $H$, and all $H$-representations $V,W$ such that $V$ is \emph{faithful}, the adjoint structure map $X(V)\to \Omega^WX(V\oplus W)$ is a $\mathcal G_{H,G}$-weak equivalence.
\end{definition}

\begin{definition}
    A map $f\colon X\to Y$ of orthogonal $G$-spectra is called a \emph{$G$-global weak equivalence} if $\phi^*f\colon\phi^*X\to\phi^*Y$ is an $H$-equivariant weak equivalence for every compact Lie group $H$ and every Lie group homomorphism $\phi\colon H\to G$.
\end{definition}

\begin{rk}\label{rk:projections-suffice}
    Factoring $\phi$ as the composite $\smash{H\xrightarrow{\;(\id,\phi)\;}H\times G\xrightarrow{\;\pr\;}G}$
    and using that equivariant weak equivalences are preserved by restriction along \emph{injective} homomorphisms (Proposition~\ref{prop:i!-left-Quillen-equiv-sp}), we see that a map is a $G$-global weak equivalence if and only if it is a $(G\times H)$-equivariant weak equivalence with respect to the trivial $H$-action for every compact Lie group $H$.
\end{rk}

\begin{lemma}
    Any $G$-global level weak equivalence is in particular a $G$-global weak equivalence.
    \begin{proof}
        Let $H$ be any compact Lie group. By \cite{broecker-tom-dieck}*{Theorem~4.1}, there exists a faithful $H$-representation $V$, which we may assume without loss of generality to be a subrepresentation of the chosen complete $H$-universe $\mathcal U_H$. Then the subposet of $s(\mathcal U_H)$ of all subrepresentations containing $V$ is cofinal, whence so is the (larger) subposet of all faithful subrepresentations. To complete the proof it then suffices to observe that if $f\colon X\to Y$ is any $G$-global level weak equivalence and $W$ a \emph{faithful} $H$-representation classified by some $\rho\colon H\to \O(W)$, then $(\phi^*f)(W)$ is an $H$-equivariant weak equivalence as $(\phi^*f)(W)^K=f(W)^{\phi\rho^{-1}|_{\rho(K)}}$ for any $K\subset H$.
    \end{proof}
\end{lemma}

\begin{thm}\label{thm:G-gl-Sp-model-struct}
    The $G$-global level model structure on $\Gosp{G}$ admits a Bousfield localization whose weak equivalences are precisely the $G$-global weak equivalences and whose fibrant objects are precisely the $G$-global $\Omega$-spectra. The resulting model category is again proper, topological, and cofibrantly generated with generating cofibrations $I$ and generating acyclic cofibrations $J\cup Z$, where $I$ and $J$ are the generating (acyclic) cofibrations for the $G$-global level model structure from Proposition~\ref{prop:global-lvl-GSp} and the set $Z$ is described in Construction~\ref{constr:new-gen-cof} below. Moreover, pushouts along unbased level h-cofibrations are homotopy pushouts.
\end{thm}

\begin{constr}\label{constr:new-gen-cof}
    Recall the map $\kappa_{V,W}^{H,H}$ from Construction~\ref{constr:generating-equiv-acyclic} together with its factorizations as a cofibration $\smash{\lambda_{V,W}^{H,H}}$ followed by a homotopy equivalence via the mapping cyclinder. We define $Z$ as the set given by the pushout products of the inclusions $\partial D^n\hookrightarrow D^n$ with the maps $\phi_!\lambda_{V,W}^{H,H}$ for all compact Lie groups $H$, $H$-representations $V,W$ where $V$ is \emph{faithful}, and all continuous group homomorphisms $\phi\colon H\to G$.
\end{constr}

The following will be the key non-formal ingredient for the proof of Theorem~\ref{thm:G-gl-Sp-model-struct}:

\begin{prop}\label{prop:kappa-gl}
    Let $\phi\colon H\to G$ be a homomorphism of compact Lie groups, and let $V,W$ be $H$-representations such that $V$ is faithful as a $\ker(\phi)$-representation. Then $\phi_!\kappa^{H,H}_{V,W}$ is a $G$-global weak equivalence.
\end{prop}

For setting up the model structure, we will only need the special case where all of $H$ acts faithfully; the other extreme where $\ker(\phi)=1$ (meaning that there is no condition on the representation $V$) will become relevant in §\ref{subsec:comp-spectra}.

\begin{proof}
    In light of Remark~\ref{rk:projections-suffice} we have to show that $\phi_!\kappa_{V,W}^{H,H}$ is a $(G\times K)$-weak equivalence for any compact Lie group $K$ (acting trivially on both sides). Viewing $V$ and $W$ as $(H\times K)$-representations via the trivial $K$-action, and replacing $\phi$ by $\phi\times\id_K\colon H\times K\to G\times K$, we may then assume that $K=1$; in other words, we only have to show that 
    $\phi_!\kappa_{V,W}^{H,H}\colon \phi_!\big(\cat{O}(V\oplus W,{-}\big)\smashp S^W)\to \phi_!\cat{O}(V,-)$
    is a $G$-\emph{equivariant} weak equivalence.
    
    For this pick a cofibrant replacement $p\colon\Ee\to (G\times H)/\Gamma_{H,\phi}$ in the ${(\mathcal G_{H,G}\cap\mathcal G_{G,H})}$-model structure on $\cat{$\bm{(G\times H)}$-Top}$, and consider the induced commutative square
    \[\hskip-20.56pt\hfuzz=20.57pt
        \begin{tikzcd}[cramped]
            (\cat{O}(V\oplus W,-)\smashp S^W) \smashp_H \Ee_+\arrow[d]\arrow[r, "\kappa_{V,W}^{H,H}\smashp_H\Ee_+"] &[6.66em] \cat{O}(V,-)\smashp_H \Ee_+\arrow[d]\\
            (\cat{O}(V\oplus W,-)\smashp S^W) \smashp_H (G\times H)_+/\Gamma_{H,\phi} \arrow[r, "\kappa_{V,W}^{H,H}\smashp_H(G\times H)_+/\Gamma_{H,\phi}"'] & \cat{O}(V,-)\smashp_H (G\times H)_+/\Gamma_{H,\phi}\rlap.
        \end{tikzcd}
    \]
    The bottom map agrees up to isomorphism with $\phi_!\kappa_{V,W}^{H,H}$; thus, it will suffice by 2-out-of-3 that the other three maps are $G$-equivariant weak equivalences. For the top map this is an instance of \cite{schwede2018global}*{Proposition 3.2.20} as $\Ee$ is both $G$-free and $H$-free by $(\mathcal G_{H,G}\cap\mathcal G_{G,H})$-cofibrancy while $\kappa_{V,W}^{H,H}$ is an $H$-equivariant weak equivalence as it agrees up to level weak equivalence with a generating acyclic cofibration.
    
    We will now show that the right-hand vertical map is a $G$-equivariant weak equivalence; an analogous argument will then show that the same holds for the left-hand vertical map, completing the proof. For this, we first note that for any $G'\subset G$, any $G'$-subrepresentation of $\Uu_G$ is contained in a $G$-subrepresentation by \cite{schwede_orbispaces_2020}*{Proposition A.7(i)}. It will therefore suffice that $\cat{O}(V,U)\smashp_H p_+$ is a $G$-equivariant weak equivalence for every $G$-representation $U$, for which we will verify the assumptions of Corollary~\ref{cor:balanced-smash-product-oddly-specific} with $N=\ker(\phi)$. By definition, $p_+$ is a $(\mathcal G_{H,G}\cap\mathcal G_{G,H})$-weak equivalence. Moreover, $\Gamma_{H,\phi}$ belongs to the family $\Ff_2$ of graph subgroups for homomorphisms with kernel contained in $N=\ker(\phi)$, so that $(G\times H)_+/\Gamma_{H,\phi}$ is $\Ff_2$-cofibrant, while $\Ee$ is even $(\mathcal G_{H,G}\cap\mathcal G_{G,H})$-cofibrant. Finally, Proposition~\ref{prop:O-mapping-CW} shows that $\cat{O}(V,U)$ is a based $(G\times H)$-CW-complex with free $N$-action outside the basepoint, i.e.\ it is $\Ff_1$-cofibrant, where $\Ff_1$ denotes the family of subgroups intersecting $N=\ker(\phi)$ trivially.
\end{proof}

\begin{prop}
    Every map in the set $Z$ from Construction~\ref{constr:new-gen-cof} is a $G$-global weak equivalence.
    \begin{proof}
        We will first show that each of the maps $\smash{\lambda_{V,W}^{(\phi)}\coloneqq\phi_!(\lambda_{V,W}^{H,H})}$ is a $G$-global weak equivalence. Recall that the map $\lambda_{V,W}^{H,H}$ was defined by factoring $\kappa_{V,W}^{H,H}$ into a cofibration followed by a homotopy equivalence using the mapping cylinder construction; write $z$ for the other map in this factorization. Then $\phi_!z$ is again a homotopy equivalence as $\phi_!$ is a topologically enriched left adjoint, while $\phi_!(\kappa_{V,W}^{H,H})$ is a $G$-global weak equivalence by the previous proposition; thus, $\smash{\lambda_{V,W}^{(\phi)}}$ is a $G$-global weak equivalence by 2-out-of-3.

        To complete the proof, it will suffice to show that the pushout product of $\smash{\lambda_{V,W}^{(\phi)}}$ with any generating cofibration of $\cat{Top}_*$ is a $G$-global weak equivalence. We will prove this more generally for any $G$-global weak equivalence $f\colon X\to Y$ in place of $\smash{\lambda_{V,W}^{(\phi)}}$. Restricting along any continuous group homomorphism $\psi\colon K\to G$ (and replacing $G$ by $K$), we may equivalently show the analogous statement for $G$-\emph{equivariant} weak equivalences. For this we consider the defining pushout square
        \[
            \begin{tikzcd}
                X\smashp \partial D^n_+\arrow[d,"f\smashp\partial D^n_+"']\arrow[dr, phantom, "\ulcorner"{very near end}]\arrow[r] & X\smashp D^n_+\arrow[ddr, "f\smashp D^n_+", bend left=20pt]\arrow[d,"g"]\\
                Y\smashp\partial D^n_+\arrow[r]\arrow[drr, bend right=20pt] & P\arrow[dr,dashed]\\
                && Y\smashp D^n_+
            \end{tikzcd}
        \]
        The left-hand vertical map $f\smashp\partial D^n_+$ is a $G$-equivariant weak equivalence by \cite{mandell-may}*{Theorem~III.3.11}. On the other hand, evaluated at $V\in\cat{O}$, the top horizontal map is a pushout in $\cat{$\bm{(\O(V)\times G)}$-Top}$ of the h-cofibration $X(V)\times\partial D^n\hookrightarrow X(V)\times D^n$; thus, the top map is an unbased level h-cofibration, and the pushout square is a homotopy pushout. We conclude that $g$ is a $G$-equivariant weak equivalence. On the other hand, also $f\smashp D^n_+$ is again a $G$-equivariant weak equivalence by   \cite{mandell-may}*{Theorem~III.3.11} (or simply because it agrees up to homotopy equivalence with $f$), hence so is the pushout-product map $P\to Y\smashp D^n_+$ by 2-out-of-3.
    \end{proof}
\end{prop}

\begin{proof}[Proof of Theorem~\ref{thm:G-gl-Sp-model-struct}]
    We consider the functor
    \begin{equation}\label{eq:brutal-right-adjoint}
        (\sh^V\phi^*)_{H,V,\phi}\colon\Gosp{G}\to\prod_{H,V,\phi}\Gosp{H},
    \end{equation}
    where the product runs over all compact Lie groups $H$, faithful $H$-representations $V$ (both up to isomorphism), and Lie group homomorphisms $\phi\colon H\to G$. We equip the product on the right with the product of the $H$-equivariant model structures, giving it the structure of a cofibrantly generated model category with generating (acyclic) cofibrations given by those tuples $(f_{H,V,\phi})$ such that all but one $f_{H,V,\phi}$ are the identity of the zero object, while the one non-trivial $f_{H,V,\phi}$ is one of the standard generating (acyclic) cofibration. We will denote the resulting sets of generating (acyclic) cofibrations by $I_{\Pi}$ and $J_{\Pi}$, respectively.

    If $V$ is any $H$-representation, then $\smash{\sh^V\simeq S^V\smashp{-}}$ by \cite{schwede2018global}*{Proposition 3.1.25}, so Proposition~\ref{prop:rep-spheres-invertible-equiv-model} shows that $\smash{\sh^V}$ preserves and reflects $H$-equivariant weak equivalences. Thus, a map of orthogonal $G$-spectra is a $G$-global weak equivalence if and only if $(\ref{eq:brutal-right-adjoint})$ sends it to a weak equivalence in the product model structure.
    
    On the other hand, a simple direct computation shows that the left adjoint $\bigvee_{H,V,\phi}\phi_!\big(\cat{O}(V,-)\smashp(-)\big)$ of $(\ref{eq:brutal-right-adjoint})$ sends $I_\Pi$ to the putative generating cofibrations $I$ of the $G$-global model structure, while $J_\Pi$ gets sent to the putative generating acyclic cofibrations $J\cup Z$. Thus, to prove that the desired Bousfield localization exists and is cofibrantly generated by the sets $I$ and $J\cup Z$, it will suffice to verify the assumptions of the Crans--Kan transfer criterion once more. The proof that the sets $I$ and $J\cup Z$ admit the small object argument is analogous to the case of the level model structure. To see that any relative $(J\cup Z)$-cell complex is a $G$-global weak equivalence, observe first that any map in $J$ is even a $G$-global level weak equivalence while any map in $Z$ is a $G$-global weak equivalence by the previous proposition. As any map in $J\cup Z$ is in particular an unbased level h-cofibration, Lemma~\ref{lemma:equiv-spectra-homotopy-po} shows that any pushout of a map in $J\cup Z$ is again a $G$-global weak equivalence. Since any such pushout is moreover in particular levelwise a closed embedding, any transfinite composition of them commutes with the formation of equivariant stable homotopy groups, and hence is again a $G$-global weak equivalence.
    
    This completes the proof that there exists a Bousfield localization with the prescribed weak equivalences and cofibrantly generated by $I$ and $J\cup Z$. Moreover, it follows directly from the construction as a transferred model structure that the fibrant objects of this model structure are precisely the $G$-global $\Omega$-spectra.

    Finally, the model structure is right proper and topological as this holds for the $H$-equivariant model structures, while Lemma~\ref{lemma:equiv-spectra-homotopy-po} together with \cite{g-global}*{Proposition~A.2.15} implies that the model structure is left proper and that pushouts along unbased level h-cofibrations are homotopy pushouts.
\end{proof}

\begin{rk}
    For $G=1$, the above model structure has the same generating (acyclic) cofibrations as---and hence agrees with---the $\All$-global model structure from \cite{schwede2018global}*{Theorem 4.3.17}.
\end{rk}

\begin{definition}
    We define the \emph{$\infty$-category of $G$-global spectra} $\mySp_\textup{$G$-gl}$ as the $\infty$-categorical localization of $\Gosp{G}$ at the $G$-global weak equivalences.
\end{definition}

\subsubsection{Functoriality} Once again, we will be interested in the functoriality of the above model categories as the group $G$ varies.

\begin{lemma}\label{lemma:restr-gl-sp-right-Quillen} 
    Let $\alpha\colon G\to G'$ be a homomorphism of compact Lie groups. Then
    \[
        \alpha_!\colon\Gosp{G}_\textup{$G$-global}\rightleftarrows\Gosp{G'}_\textup{$G'$-global} \noloc \alpha^*
    \]
    is a Quillen adjunction with homotopical right adjoint.
    \begin{proof}
        It is clear from the definitions that $\alpha^*$ preserves weak equivalences, while the construction as a transferred model structure in the proof of Theorem~\ref{thm:G-gl-Sp-model-struct} makes it clear that $\alpha^*$ also preserves fibrations.
    \end{proof}
\end{lemma}

\begin{lemma}\label{lemma:restr-inj-gl-sp-left-Quillen} 
    Let $\alpha\colon G\to G'$ be an \emph{injective} homomorphism of compact Lie groups. Then
    \[
        \alpha^*\colon\Gosp{G'}_\textup{$G'$-global}\rightleftarrows\Gosp{G}_\textup{$G$-global}\noloc \alpha_*
    \]
    is a Quillen adjunction.
    \begin{proof}
        Thanks to the previous lemma, it only remains to show that $\alpha^*$ preserves (level) cofibrations, or equivalently that $\alpha_*$ preserves level acyclic fibrations. This follows at once by applying Example~\ref{ex:graph-restr-lQ} levelwise.
    \end{proof}
\end{lemma}

\begin{lemma}\label{lemma:gl-sp-(co)ind-homotopical}
    Let $\alpha\colon G\to G'$ be an \emph{injective} homomorphism of compact Lie groups. Then $\alpha_!\colon\Gosp{G}_\textup{$G$-global}\to\Gosp{G'}_\textup{$G'$-global}$ is homotopical.
    \begin{proof}
        Let $f\colon X\to Y$ be a $G$-global weak equivalence, and let $H$ be any compact Lie group. Then $\triv_H \alpha_!f$ agrees up to conjugating by isomorphisms with $(\alpha\times H)_!\triv_Hf$, and the latter is a $(G'\times H)$-equivariant weak equivalence by Proposition~\ref{prop:i!-left-Quillen-equiv-sp}. Letting $H$ vary, this implies via Remark~\ref{rk:projections-suffice} that $\alpha_!f$ is a $G'$-global weak equivalence.
    \end{proof}
\end{lemma}

\subsubsection{Multiplicative properties} The same pointset level construction as in the $G$-equivariant setting gives us a smash product on orthogonal $G$-spectra. Again this interacts nicely with the model structure:

\begin{prop}\label{prop:gl-flatness}
    Let $X$ be any flat orthogonal $G$-spectrum. Then $X\smashp{-}$ preserves $G$-global weak equivalences. Moreover, for any $Y$, the functor ${-}\smashp Y$ preserves $G$-global weak equivalences between flat orthogonal $G$-spectra.
    \begin{proof}
        Let $f$ be a $G$-global weak equivalence. For the first statement, we have to show that $\phi^*(X\smashp f)=\phi^*(X)\smashp\phi^*(f)$ is an $H$-equivariant weak equivalence for every $\phi\colon H\to G$. As $\phi^*(X)$ is again flat by \cite{schwede2018global}*{Remark 3.5.8}, this follows at once from the first part of Theorem~\ref{thm:equiv-flatness}. The second part of the proposition follows in the same way from the second part of Theorem~\ref{thm:equiv-flatness}.
    \end{proof}
\end{prop}

\begin{prop}\label{prop:smash-gl-left-Quillen}
    The smash product of orthogonal spectra defines a left Quillen bifunctor 
    \[
        \Gosp{G}_\textup{$G$-global}\times\Gosp{H}_\textup{$H$-global}\to
        \Gosp{(G\times H)}_\textup{$(G\times H)$-global}.
    \]
    \begin{proof}
        To see that the pushout product axiom for cofibrations holds, it suffices to treat the case of generating cofibrations. Using that the enriched Yoneda embedding $\cat{O}^\op\to\osp$ admits a strong symmetric monoidal structure, we then compute that the pushout product of $\cat{O}(V,-)\smashp_{\phi} G_+\smashp (\partial D^n\hookrightarrow D^n)$ with $\cat{O}(W,-)\smashp_\psi G_+\smashp (\partial D^m\hookrightarrow D^m)$ (where $\phi\colon K\to G$ and $\psi\colon L\to H$ are homomorphisms) is given by
        \[
            \cat{O}(V\oplus W,-)\smashp_{\phi\times\psi}(G\times H)_+\smashp\big((\partial D^n\hookrightarrow D^n)\ppo(\partial D^m\hookrightarrow D^m)\big)_+,
        \]
        where we view $H\times K\subset \O(V)\times \O(W)$ as a subgroup of $\O(V\oplus W)$ via block sum. As $\cat{Top}$ is cartesian and $\Gosp{(G\times H)}$ is topological, this is then indeed a cofibration, completing the proof of the pushout product axiom for cofibrations.

        Next, we will show that $i\ppo f$ is a $(G\times H)$-global weak equivalence for any cofibration $i\colon A\to B$ between cofibrant objects and any $H$-global weak equivalence $f\colon X\to Y$ between cofibrant objects; this then in particular applies to the generating (acyclic) cofibrations, proving (up to exchanging the roles of $G$ and $H$) the remaining case for the pushout product axiom. For this we consider the defining pushout
        \[
            \begin{tikzcd}
                A\smashp X\arrow[r,"i\smashp X"]\arrow[d,"A\smashp f"']\arrow[dr,phantom,"\ulcorner"{very near end}] & B\smashp X\arrow[d]\arrow[ddr, bend left=15pt,"B\smashp f"]\\
                A\smashp Y\arrow[r]\arrow[drr, bend right=15pt,"i\smashp Y"'] & P\arrow[dr,dashed,"i\ppo f"{description}]\\
                && B\smashp Y\rlap.
            \end{tikzcd}
        \]
        By the previous proposition applied to $\triv_HA$ and $\triv_Gf$, the map $A\smashp f$ is a $(G\times H)$-global weak equivalence, and likewise for $B\smashp f$. Moreover, $i\smashp X$ is a cofibration by the previous step, so that $B\smashp X\to P$ is also a $(G\times H)$-global weak equivalence by left properness. We conclude from $2$-out-of-$3$ that also $P\to B\smashp Y$ is a $(G\times H)$-global weak equivalence, as claimed.
    \end{proof}
\end{prop}

\begin{cor}\label{cor:smash-gl-left-Quillen}
    The smash product defines a left Quillen bifunctor 
    \begin{equation}\label{eq:smash-product-G-gl}
        \Gosp{G}\times\Gosp{G}\to\Gosp{G}
    \end{equation} 
    with respect to the $G$-global model structures for any compact Lie group $G$.
    \begin{proof}
        We may factor $(\ref{eq:smash-product-G-gl})$ as
        \[
            \Gosp{G}\times\Gosp{G}\xrightarrow{\;{-}\smashp{-}\;}\Gosp{(G\times G)}\xrightarrow{\;\Delta^*\;}\Gosp{G}.
        \]
        The claim now follows from the previous proposition together with Lemma~\ref{lemma:restr-inj-gl-sp-left-Quillen}.
    \end{proof}
\end{cor}

\begin{lemma}\label{lemma:G-global-spectra-tensoring}
    The smash product $\cat{$\bm G$-Top}_*\times\Gosp{G}\to\Gosp{G}$ is a left Quillen bifunctor for the $\All$-model structure on $\cat{$\bm G$-Top}_*$ and the $G$-global model structure on $\Gosp{G}$. Moreover, if $X\in\cat{$\bm G$-Top}_*$ is $\All$-cofibrant, then $X\smashp{-}$ is fully homotopical, and if $Y\in\Gosp{G}$ is flat, so is ${-}\smashp Y\colon\cat{$\bm G$-Top}_*\to\Gosp{G}$.
    \begin{proof}
        The pushout product axiom for cofibrations follows from Lemmas~\ref{lemma:restr-gl-sp-right-Quillen} and~\ref{lemma:restr-inj-gl-sp-left-Quillen} as in the proof of Lemma~\ref{lemma:G-Top*-tensored}. Moreover, if $X\in\cat{$\bm G$-Top}_*$ is $\All$-cofibrant, then $X\smashp{-}\cong\Sigma^\infty X\smashp{-}$ is homotopical by Proposition~\ref{prop:gl-flatness}, as is ${-}\smashp Y$ for any flat $Y$. With these established, the pushout-product axiom for acyclic cofibrations then follows as in the proof of Proposition~\ref{prop:smash-gl-left-Quillen}.
    \end{proof}
\end{lemma}

\subsubsection{Suspension spectra} Next, we come to the analogue of the suspension spectrum functor in the global setting.

\begin{constr}
    We recall from \cite{schwede2018global}*{Construction~4.1.6} the construction of the (topologically enriched) functor $\Omega^\bullet\colon\osp\to\cat{$\cat{L}$-Top}_*$:
    
    The based orthogonal space $\Omega^\bullet X$ is given on objects by $(\Omega^\bullet X)(V) =\Omega^VX(V)$, and a linear isometric embedding $i\colon V\to W$ acts via the composite
    \[
        \Omega^V X(V) \xrightarrow{\;\Omega^V\tilde\sigma\;}\Omega^V\Omega^{W-i(V)} X(W)\cong\Omega^WX(W).
    \]
    The functoriality of $\Omega^\bullet$ in morphisms of orthogonal spectra is the obvious one.
\end{constr}

We will construct the enriched left adjoint $\Sigma^\bullet$ of $\Omega^\bullet$ in two steps:

\begin{constr}\label{constr:diag}
    Write $\cat{$\cat L$-Sp}$ for the enriched category of $\cat{Top}_*$-enriched functors $\cat{O}\to\cat{$\cat{L}$-Top}_*$, which we will conflate with the category of $\cat{Top}_*$-enriched functors $\cat{O}\smashp\cat{L}_+\to\cat{Top}_*$, where $\cat{L}_+$ is the $\cat{Top}_*$-enriched category obtained from $\cat{L}$ by adding disjoint basepoints to morphism spaces, and ${-}\smashp{-}$ denotes the tensor product of enriched categories, with objects all pairs $(V\in\cat{O},W\in\cat{L})$ and with morphism spaces given by the smash product of the morphism spaces in $\cat{O}$ and $\cat{L}_+$.

    We then define the $\cat{Top}_*$-enriched functor $\diag\colon\cat{$\cat{L}$-Sp}\to\osp$ as the restriction along the functor $\Delta\colon\cat{O}\to\cat{O}\smashp\cat{L}_+$ given on objects by sending a finite-dimensional inner-product space $V$ to $(V,V)$ and given on morphisms by sending $[A,w]\in\cat{O}(V,W)$ to $[A,v]\smashp A\in\cat{O}(V,W)\smashp\cat{L}(V,W)_+$. More concretely, this means that $\diag(X)(V)=X(V)(V)$ for any $X\colon\cat{L}\to\osp$, with the structure maps given by the composite of the spectral and orthogonal space structure maps; the functoriality in $X$ is the obvious one.
\end{constr}

\begin{constr}\label{constr:Sigma-bullet}
    Applying $\Sigma^\infty$ levelwise yields a $\cat{Top}_*$-enriched functor $\Sigma^\infty\colon\cat{$\cat{L}$-Top}_*\to\cat{$\cat{L}$-Sp}$, and we define $\Sigma^\bullet\coloneq\diag\circ\Sigma^\infty\colon\cat{\cat{L}-Top}_*\to\osp$. More concretely, this means that $(\Sigma^\bullet X)(V)=\Sigma^VX(V)\coloneqq S^V\smashp X(V)$, with the evident functoriality in $X$ and with structure maps induced by the structure maps of $X$. One readily verifies that this is an enriched left adjoint to the functor from the previous construction, with unit $X\to\Omega^\bullet\Sigma^\bullet X$ given after evaluation at each $V\in\cat{L}$ by the unit $X(V)\to\Omega^V\Sigma^VX(V)$ of the adjunction $\Omega^V\dashv\Sigma^V$.    
    
    As usual, we will write $\Sigma^\bullet_+\coloneq\Sigma^\bullet\circ(-)_+\colon\cat{$\cat{L}$-Top}\to\osp$; plugging in the definitions, the latter agrees with the functor of the same name from \cite{schwede2018global}*{Construction~4.1.7}.  
\end{constr}

\begin{prop}\label{prop:Sigma-bullet-homotopical-left-Quillen}
    Let $G$ be compact Lie. The adjunctions
    \begin{align*}
        \Sigma^\bullet_+\colon\cat{$\bm G$-$\cat{L}$-Top}&\rightleftarrows\Gosp{G}\noloc\Omega^\bullet\\
        \Sigma^\bullet\colon\cat{$\bm G$-$\cat{L}$-Top}_*&\rightleftarrows\Gosp{G}\noloc\Omega^\bullet
    \end{align*}
    are Quillen adjunctions, and $\Sigma^\bullet_+$ is fully homotopical.
    \begin{proof}
        It is clear from the definition that $\Omega^\bullet$ preserves fibrant objects, i.e.\ it sends $G$-global $\Omega$-spectra to static $G$-orthogonal spaces. To complete the proof that these are Quillen adjunctions, it will then suffice by \cite{HTT}*{Corollary A.3.7.2} that $\Sigma^\bullet$ preserves cofibrations or, equivalently, that $\Omega^\bullet$ preserves acyclic (level) fibrations as a functor to $\cat{$\bm G$-$\cat{L}$-Top}_*$. Plugging in the definitions, this amounts to saying that
        \[
            S^V\smashp{-}\colon\cat{$\bm{(\O(V)\times G)}$-Top}_*\rightleftarrows\cat{$\bm{(\O(V)\times G)}$-Top}_*\noloc\Omega^V
        \] 
        is a Quillen adjunction for the $\mathcal G_{\O(V),G}$-model structure. This is an instance of Lemma~\ref{lemma:smash-lQ}, using  that $S^V$ admits the structure of an $(\O(V)\times G)$-CW-complex.

        It remains to show that $\Sigma^\bullet_+$ is homotopical, for which it will suffice by factoring a general weak equivalence into an acyclic cofibration followed by an acyclic fibration that $\Sigma^\bullet_+$ preserves $G$-global level weak equivalences. Plugging in the definitions this amounts to saying that
        \[
            S^V\smashp (-)_+\colon\cat{$\bm{(\O(V)\times G)}$-Top}\to\cat{$\bm{(\O(V)\times G)}$-Top}_*
        \]
        preserves $\mathcal G_{\O(V),G}$-weak equivalences, which is an instance of Proposition~\ref{prop:smash-well-based}.
    \end{proof}
\end{prop}

\subsubsection{Stability} $G$-global spectra enjoy the same refined notion of stability as their equivariant cousins:

\begin{prop}\label{prop:resp-stability-G-gl}
    Let $G$ be any compact Lie group and let $V$ be a $G$-representation. Then both adjoints in
    \[
        \Sigma^V\colon\Gosp{G}_\textup{$G$-global}\rightleftarrows\Gosp{G}_\textup{$G$-global}\noloc\Omega^V
    \]
    are homotopical, and the induced adjunction $\mySp_\textup{$G$-gl}\rightleftarrows\mySp_\textup{$G$-gl}$ is an adjoint equivalence.
    \begin{proof}
        If $f$ is a $G$-global weak equivalence and $\phi\colon H\to G$ is any homomorphism, then $\phi^*(\Sigma^Vf)=\Sigma^{\phi^*V}\phi^*f$ is an $H$-equivariant weak equivalence by Proposition~\ref{prop:rep-spheres-invertible-equiv-model}. The same argument shows that $\Omega^V$ is homotopical.

        It remains to show that the unit $X\to\Omega^V\Sigma^VX$ and counit $\Sigma^V\Omega^VX\to X$ are invertible for every $X$. This follows again from Proposition~\ref{prop:rep-spheres-invertible-equiv-model} after restricting along each group Lie group homorphism $\phi\colon H\to G$.
    \end{proof}
\end{prop}

\begin{prop}\label{prop:gl-spectra-compact-generators}
    The $\infty$-category $\mySp_\textup{$G$-gl}$ is stable and compactly generated by the objects $\Sigma^\bullet_+\big(\cat{L}(V_H,-)\times_\phi G\big)$ where $\phi$ runs through homomorphisms $\phi\colon H\to G$ of compact Lie groups (up to isomorphism), and $V_H$ is a fixed faithful $H$-representation. In particular, $\mySp_\textup{$G$-gl}$ is presentable.
    \begin{proof}
        As in the equivariant case (Proposition~\ref{prop:equiv-sp-compact-gen}), stability follows from the previous proposition, and it will suffice to show that the homotopy category is compactly generated in the triangulated sense by the given objects. This is in turn the content of \cite{schwede-stiefel}*{Corollary A.19}.
    \end{proof}
\end{prop}

\subsubsection{Equivariant spectra vs.~global spectra} Finally, let us explain the relation between $G$-equivariant and $G$-global spectra.

\begin{lemma}\label{lemma:equiv-spectra-vs-gl-spectra}
    Any $G$-equivariant weak equivalence between equivariantly cofibrant orthogonal $G$-spectra is a $G$-global weak equivalence. Moreover, the inclusion induces a fully faithful functor $i\colon\mySp_G\simeq\mySp_G^\textup{cof}\hookrightarrow\mySp_\textup{$G$-gl}$, and this admits a right adjoint, which is a localization at the $G$-equivariant weak equivalences.
    \begin{proof}
        The first claim is immediate from Proposition~\ref{prop:restr-left-Quillen-G-Sp} together with Ken Brown's Lemma. 
        For the second claim, we fix a functorial cofibrant replacement $\epsilon\colon Q\to\id$ for the $G$-equivariant model structure on $\Gosp{G}$. By $2$-out-of-$3$, a map $f$ is a $G$-equivariant weak equivalence if and only if $Qf$ is a $G$-equivariant or equivalently (by the above) a $G$-global weak equivalence. In particular, the full subcategory $\Cc\subset\Gosp{G}$ of those objects for which $\epsilon_X\colon QX\to X$ is a \emph{$G$-global} weak equivalence is itself closed under $G$-global weak equivalences. Thus, \cite{g-global}*{Proposition~A.1.15} shows that the inclusion of $\Cc$ remains fully faithful after localizing both sides at the $G$-global weak equivalences. As on the other hand the inclusion of equivariantly cofibrant objects into $\Cc$ admits a homotopy inverse (namely, $Q$), this shows that the map $\mySp_G^\textup{cof}\to\mySp_\text{$G$-gl}$ induced by the inclusion is indeed fully faithful.

        We will now show that $Q$ defines a right adjoint to this map; as we have already seen that it precisely inverts the $G$-equivariant weak equivalences, this will then complete the proof of the lemma. For the claim it suffices to show that $\epsilon$ exhibits $Q$ as an idempotent comonad on $\mySp_\text{$G$-gl}$, i.e.\ that for every $X\in\Gosp{G}$ the two maps $Q\epsilon_X,\epsilon_{QX}\colon QX\rightrightarrows QQX$ are $G$-global weak equivalences. This follows at once from them being \emph{$G$-equivariant} weak equivalences between equivariantly cofibrant objects by construction.
    \end{proof}
\end{lemma}

\subsection{How to do global homotopy theory less efficiently}\label{subsec:L-spectra}
In this subsection, we will introduce another model of the $\infty$-category $\mySp_\text{$G$-gl}$, this time in terms of the category $\cat{$\bm G$-$\cat{L}$-Sp}$ of orthogonal spectrum objects in based $G$-orthogonal spaces from Construction~\ref{constr:diag}. At first, this might seem like a silly idea---after all, we just said that the same $\infty$-categorical information can be packaged in a much more efficient way using plain orthogonal $G$-spectra. However, this new bloated model will be central for our proof of the universal property of global spectra (Theorem~\ref{introthm:stable-main}). The reason for this is that it once again comes from Bousfield localizing a suitable \emph{$G$-global level model structure} on $\cat{$\bm G$-$\cat{L}$-Sp}$, and this level model structure will \emph{not} yet be equivalent to $G$-global level model structure on $\Gosp{G}$. For our purposes, the level model structure on $\cat L$-spectra will be better behaved: we will see later that the corresponding derived suspension spectrum functor $\Sigma^\infty$ is fully faithful and that its right adjoint $\Omega^\infty$ in turn admit another right adjoint; both of these formal properties will then be exploited crucially in the proof of Theorem~\ref{introthm:stable-main}.

\medskip
\subsubsection{The level model structure} As promised, we begin by constructing a suitable $G$-global level model structure.

\begin{definition}
    A map $f\colon X\to Y$ in $\cat{$\bm G$-$\cat{L}$-Sp}$ is called a \emph{$G$-global level weak equivalence} or \emph{$G$-global level fibration} if $f(V)\colon X(V)\to Y(V)$ is an $(\O(V)\times G)$-global weak equivalence or fibration, respectively, in $\cat{$\bm{(\O(V)\times G)}$-$\cat{L}$-Top}$ for every finite-dimensional inner product space $V$.
\end{definition}

\begin{prop}\label{prop:LSp-lvl}
    The $G$-global level weak equivalences and $G$-global level fibrations participate in a model structure on $\cat{$\bm G$-$\cat{L}$-Sp}$, which we call the \emph{$\bm G$-global level model structure}. This model structure is proper, topological, and cofibrantly generated with set $I$ of generating cofibrations given by the maps
    \[
        \phi_!(\cat{O}(V,-)\smashp\cat{L}(W,-)_+)\smashp(\partial D^n\hookrightarrow D^n)_+
    \]
    where we let $H$ run through all compact Lie groups, $V$ through all $H$-representations, $W$ through all \emph{faithful} $H$-representations, and $\phi$ through all homomorphisms $H\to G$.
\end{prop}

\begin{definition}
    We write $\myLSp_\textup{$G$-lvl}$ for the $\infty$-category associated to the $G$-global level model structure on $\cat{$\bm G$-$\cat{L}$-Sp}$.\footnote{The reader is invited to choose for themself whether the letter `$\mathfrak d$' stands for `double,' `diagonal,' or `dumb.'}
\end{definition}

The only (minor) difficulty in proving the proposition will be that while the weak equivalences are again defined levelwise, the evaluation functors are very far from preserving cofibrations this time. In order to understand the homotopical behaviour of pushouts along cofibrations, we will therefore again isolate a larger class of maps such at pushouts along them are homotopically well-behaved.

\begin{lemma}\label{lemma:fgt-unbased-h-cof}
    Let $G$ be a compact Lie group, let $W$ be a finite-dimensional inner product space, and write $K$ for the class of all maps in $\cat{$\bm{(\O(W)\times G)}$-$\cat{L}$-Top}_*$ of the form $\cat{O}(V,W)\smashp_{\O(V)} i$ where we let $V$ run through finite-dimensional inner product spaces and $i$ runs through cofibrations of $\cat{$\bm{(\O(V)\times G)}$-$\cat{L}$-Top}_*$. Then every retract of a relative $K$-cell complex is a cofibration in $\cat{$\bm G$-$\cat{L}$-Top}_*$ and an h-cofibration in $\cat{$\bm{(\O(W)\times G)}$-$\cat{L}$-Top}$.
\end{lemma}
\begin{proof}
    For the first claim it suffices to show that any map in $K$ is a cofibration in $\cat{$\bm G$-$\cat{L}$-Top}_*$ and we may then further restrict to the case that $i$ is a generating cofibration, i.e.~of the form $\big(\cat{L}(U,-)\times_{(\rho,\phi)}(\O(V)\times G)\times(\partial D^n\hookrightarrow D^n)\big){}_+$ for some homomorphisms $\rho\colon H\to\O(V)$ and $\phi\colon H\to G$, faithful $H$-representation $U$, and $n\ge0$. Plugging in the definitions, $\cat{O}(V,W)\smashp_{\O(V)} i$ can then be identified with the map $\phi_!\big(\cat{O}(V,W)\smashp\cat{L}(U,-)_+\smashp(\partial D^n\hookrightarrow D^n)_+\big)$. Lemma~\ref{lemma:G-Top*-tensored} shows that $\cat{O}(V,W)\smashp\cat{L}(U,-)_+\smashp(\partial D^n\hookrightarrow D^n)$ is an $H$-global cofibration, and so the claim follows from $\phi_!$ being left Quillen (Lemma~\ref{lemma:restr-right-Quillen}).
    
    For the second claim, we similarly use the closure properties from Remark~\ref{rk:h-cof-closure} to reduce to the case where $i$ is a generating cofibration, and hence in particular of the form $i=j_+$ for some cofibration of $\cat{$\bm{(\O(W)\times G)}$-$\cat{L}$-Top}$. Then $\cat{O}(V,W)\smashp_{\O(V)} j_+$ is a pushout in $\cat{$\bm{(\O(W)\times G)}$-$\cat{L}$-Top}$ of $\cat{O}(V,W)\times_{\O(V)} j$. Using Remark~\ref{rk:level-cof-are-h-cof} one easily checks that $\cat{O}(V,W)\times j$ is an (unbased) h-cofibration of orthogonal $(\O(V)\times\O(W)\times G)$-spaces, so that $\cat{O}(V,W)\times_{\O(V)} j$ is an unbased h-cofibration of orthogonal $(\O(W)\times G)$-spaces by \cite{schwede2018global}*{Corollary~A.30(ii)}. Using once more that h-cofibrations are closed under pushout, this yields the claim.
\end{proof}

\begin{definition}
    We call a map $f$ in $\cat{$\bm G$-$\cat{L}$-Sp}$ an \emph{unbased level h-cofibration} if $f(V)$ is an h-cofibration in $\cat{$\bm{(\O(V)\times G)}$-$\cat{L}$-Top}$ for every $V\in\cat{O}$. We call $X\in\cat{$\bm G$-$\cat{L}$-Sp}$ \emph{well-based} if the unique map $0\to X$ is an unbased level h-cofibration.
\end{definition}

\begin{proof}[Proof of Proposition~\ref{prop:LSp-lvl}]
    Similarly to the construction of the $G$-global level model structure on $\Gosp{G}$, we will apply the Crans--Kan transfer criterion for
    \begin{equation}\label{eq:LSp-transfer-RA}
        \cat{$\bm G$-$\cat{L}$-Sp}\xrightarrow{(\ev_V)_V}\prod_V\cat{$\bm{(G\times \O(V))}$-$\cat{L}$-Top}
    \end{equation}
    where the product runs over a set of representatives for the isomorphism classes of finite-dimensional inner product spaces. Assuming the conditions of the criterion are satisfied, one again directly checks that the resulting model structure would have the correct weak equivalences and fibrations. Moreover, the generating (acyclic) cofibrations would be the maps of the form $\cat{O}(V,-)\smashp_{\O(V)}i_+$ where $V$ runs through finite-dimensional inner product spaces and $i$ through generating (acyclic) cofibrations of $\cat{$\bm{(\O(V)\times G)}$-$\cat{L}$-Top}$. If we take the generating cofibrations of the latter to be the ones from Proposition~\ref{prop:orth-lvl}, this then simplifies to the claimed generating cofibrations.
    
    By Lemma~\ref{lemma:fgt-unbased-h-cof}, all maps in $I$ and $J$ are in particular levelwise closed embeddings, while \cite{barrero2021}*{Lemma A.16} shows that the sources of the standard generating (acyclic) cofibrations of $\cat{$\bm{(\O(V)\times G)}$-\cat{L}-Top}$ are small with respect to levelwise closed embeddings; thus, the same adjunction argument as in the proof of Proposition~\ref{prop:global-lvl-GSp} shows that the sets $I$ and $J$ permit the small object argument. It remains to prove that any relative $J$-cell complex is a $G$-global level weak equivalence.
    
    For this we first consider a generating acyclic cofibration $k\coloneqq\cat{O}(V,-)\smashp_{\O(V)}j_+$, where $j$ is a generating acyclic cofibration of $\cat{$\bm{(\O(V)\times G)}$-$\cat{L}$-Top}$; by Remark~\ref{rk:generating-ayclic-between-cof} we may assume that $j$ is a map between cofibrant orthogonal $(\O(V)\times G)$-spaces and hence in particular levelwise Hausdorff. We claim that $k$ is a $G$-global weak equivalence and an unbased levelwise h-cofibration. If we fix $W\in\cat{O}$, then $k(W)=\cat{O}(V,W)\smashp_{\O(V)}j_+$ is a pushout of $\cat{O}(V,W)\times_{\O(V)}j$, so it suffices to show that the latter is an unbased h-cofibration and an $(\O(W)\times G)$-global weak equivalence. The former was verified as part of the proof of Lemma~\ref{lemma:fgt-unbased-h-cof}. For the latter, we first observe that $\cat{O}(V,W)\times j$ is an $(\O(V)\times\O(W)\times G)$-global weak equivalence by Lemma~\ref{lemma:restr-right-Quillen} together with Corollary~\ref{cor:cart-prod-homotopical}. On the other hand, its source and target are levelwise Hausdorff by assumption on $j$, and they are $\O(V)$-free as $\O(V)$ acts freely on $\cat{O}(V,W)$; thus, Proposition~\ref{prop:quotient-we-global} shows that $\cat{O}(V,W)\times_{\O(V)}j$ is an $(\O(W)\times G)$-global weak equivalence. This completes the proof that all maps in $J$ are $G$-global level weak equivalences and levelwise unbased h-cofibrations.
    
    Lemma~\ref{lemma:h-cof-compute-po} then shows that any pushout of a map in $J$ is again a $G$-global level weak equivalence; moreover, it is also a levelwise unbased h-cofibration as the latter are closed under pushouts. Since a transfinite compositions in $\cat{$\bm{(\O(W)\times G)}$-$\cat{L}$-Top}$ of maps that are both $(\O(W)\times G)$-global weak equivalences and h-cofibrations are again $(\O(W)\times G)$-global weak equivalences by \cite{barrero2021}*{Corollary A.12}, we conclude that any relative $J$-cell complex is a $G$-global level weak equivalence. This completes the verification of the assumptions of the Crans--Kan criterion and hence finishes the proof that the model structure exists and is cofibrantly generated with the claimed generating cofibrations.
    
    The model structure is right proper and topological since it is transferred from a right proper topological model structure. Finally, left properness follows since by Lemma~\ref{lemma:fgt-unbased-h-cof} any cofibration is an unbased level h-cofibrations, and pushouts along the latter preserve $G$-global level weak equivalences by Lemma~\ref{lemma:h-cof-compute-po}.
\end{proof}

With the model structure established, Lemma~\ref{lemma:fgt-unbased-h-cof} has the following consequences, the first of which was already noted in the proof of the proposition:

\begin{cor}\label{cor:SpL-level-hcof}
    Any cofibration in $\cat{$\bm G$-${\cat L}$-Sp}$ is an unbased level h-cofibration. In particular, every cofibrant object is well-based.\qed
\end{cor}

\begin{cor}\label{cor:cofibrant-LSp-closed}
    If $X$ is cofibrant in $\cat{$\bm G$-$\cat{L}$-Sp}$, then the orthogonal space $X(V)$ is closed for every $V\in\cat{O}$.\qed
\end{cor}

\begin{rk}\label{rk:LSp-Lf!}
    If $\alpha\colon G\to G'$ is any homomorphism of compact Lie groups, then applying Lemma~\ref{lemma:restr-right-Quillen} levelwise shows that $\alpha^*\colon\cat{$\bm{G'}$-$\cat{L}$-Sp}\to\cat{$\bm G$-$\cat{L}$-Sp}$ is homotopical right Quillen, giving rise to a derived adjunction $\cat{L}\alpha_!\colon\myLSp_\text{$G$-gl}\rightleftarrows\myLSp_\text{$G'$-gl}\noloc \alpha^*$.
\end{rk}

\subsubsection{A presheaf description} We will now describe the $\infty$-category $\myLSp_\textup{$G$-lvl}$ in terms of reduced presheaves:

\begin{prop}\label{prop:mySplev-unpar-presh}
    Let $G$ be any compact Lie group, and write $\mySph_G$ for the full subcategory of $\myLSp_\textup{$G$-lvl}$ spanned by the zero object together with the objects $\phi_!(\cat{O}(V,-)\smashp\cat{L}(W,-)_+)$ for all homomorphisms $\phi\colon H\to G$, arbitrary $H$-representations $V$, and all faithful $H$-representations $W$. Then the inclusion extends to an equivalence \[\PSh^0(\mySph_G)\iso\myLSp_\textup{$G$-lvl},\] where the left-hand side denotes the full subcategory of $\PSh(\mySph_G)$ spanned by presheaves sending $0$ to the terminal object. In particular, $\myLSp_\textup{$G$-lvl}$ is presentable.
\end{prop}

To prove the proposition, we will use the following based variant of Lurie's detection result for presheaf $\infty$-categories \cite{HTT}*{Corollary~5.1.6.11}:

\begin{lemma}\label{lemma:detect-pointed-presheaves}
    Let $\Cc$ be a cocomplete, pointed, and locally small $\infty$-category. Let $\Ii\subset\Cc$ be a small full subcategory containing the zero object and satisfying the following properties:
    \begin{enumerate}
        \item For any $I\in\Ii$, the unique lift $\Cc\to\Spc_*$ of $\hom(I,-)$ is cocontinuous.
        \item The functors $\hom(I,-)\colon\Cc\to\Spc$ are jointly conservative.
    \end{enumerate}
    Then $\Ii\hookrightarrow\Cc$ extends to an equivalence $\PSh^0(\Ii)\iso\Cc$.
    \begin{proof}
        As $\Cc$ is cocomplete and locally small, the inclusion extends to a functor $\Lambda\colon\PSh(\Ii)\to\Cc$ left adjoint to the restricted Yoneda embedding $\Upsilon$. By what it means to be a zero object, $\hom(0,x)\simeq1$ for all $x\in\Cc$, i.e.\ $\Upsilon$ factors through $\PSh^0(\Ii)$, so we obtain a restricted adjunction $\PSh^0(\Ii)\rightleftarrows\Cc$. It remains to show that this is in fact an adjoint equivalence. As our second assumption guarantees that $\Upsilon$ is conservative, it will suffice for this that the unit $\eta\colon\id\to\Upsilon\Lambda$ is invertible. 

        Since $\Ii^\op$ has a zero object, postcomposing with the forgetful functor induces an equivalence $\Fun^0(\Ii^\op,\Spc_*)\iso\Fun^0(\Ii^\op,\Spc)$. Moreover, the full subcategory $\Fun^0(\Ii^\op,\Spc_*)\subset\Fun(\Ii^\op,\Spc_*)$ is closed under colimits (as any colimit of initial objects is initial), so that colimits in $\Fun^0(\Ii^\op,\Spc_*)$ are computed pointwise. Thus, our first assumption guarantees that the right adjoint $\Upsilon$ is cocontinuous, so that the class of all pointed presheaves $X$ for which $\eta_X$ is invertible, is closed under colimits. It will therefore suffice to show that the unit $\eta_X$ is invertible whenever $X=\Upsilon(I)$ is representable. If we now let $Z\in\Cc$ arbitrary, then the composite
        \[
            \hom(\Lambda(X),Z)\xrightarrow{\;\Upsilon\;}\hom(\Upsilon\Lambda(X),\Upsilon(Z))\xrightarrow{\;\hom(\eta,\Upsilon(Z))} \hom(X,\Upsilon(Z))
        \]
        is invertible as it is just the adjunction equivalence. By assumption, $\Lambda(X)\simeq I\in\Ii$; if also $Z\in I$, then the first map is an equivalence by the Yoneda lemma, so that also the second map is an equivalence by $2$-out-of-$3$. We conclude that $\hom(\eta_X,\Upsilon(Z))$ is an equivalence for every $Z\in I$. As $\eta_X\colon \Upsilon(I)=X\to \Upsilon\Lambda(X)\simeq\Upsilon(I)$ is a map in Yoneda image, this implies that $\eta_X$ is an equivalence.
    \end{proof}
\end{lemma}

Applying this criterion to the situation at hand will require some preparations.

\begin{lemma}\label{lemma:ev-RA}
    Let $V\in\cat{O}$ be arbitrary. Then we have Quillen adjunctions
    \begin{align*}
        \cat{O}(V,-)\smashp_{\O(V)}{(-)}\colon\cat{$\bm{(\O(V)\times G)}$-$\cat{L}$-Top}_*&\rightleftarrows\cat{$\bm G$-$\cat{L}$-Sp}\noloc\ev_V\\
        \cat{O}(V,-)\smashp_{\O(V)}{(-)_+}\colon\cat{$\bm{(\O(V)\times G)}$-$\cat{L}$-Top}&\rightleftarrows\cat{$\bm G$-$\cat{L}$-Sp}\noloc\ev_V
    \end{align*}
    with homotopical right adjoints.
    \begin{proof}
        It is clear from the definition that $\ev_V$ preserves weak equivalences and fibrations.
    \end{proof}
\end{lemma}

\begin{cor}\label{cor:ev_V-corep}
    Let $\phi\colon H\to G$ be a homomorphism of compact Lie groups, let $V$ be an $H$-representation with corresponding homomorphism $\rho\colon H\to\O(V)$, let $W$ be a faithful $H$-representation, and let $\Uu_H$ be any complete $H$-universe. Then $\phi_!(\cat{O}(V,-)\smashp\cat{L}(W,-)_+)$ corepresents the composite
    \begin{equation}\label{eq:semifree-corep}
        \myLSp_\textup{$G$-gl}\xrightarrow{\;\ev_V\;}\myS_{\textup{$(\O(V)\times G)$-gl},*}\iso\myS_{\textup{$(\O(V)\times G)$-gl},*}^\textup{closed}\xrightarrow{\;(-)(\Uu_H)^{(\rho,\phi)}\;}\Spc_*
    \end{equation}
    \begin{proof}
        Because of the isomorphism \[\phi_!(\cat{O}(V,-)\smashp\cat{L}(W,-)_+)\cong\cat{O}(V,-)\smashp_{\O(V)}\big(\cat{L}(W,-)\times_{(\rho,\phi)}(\O(V)\times G)\big)_+\] this follows at once from the previous lemma together with Lemma~\ref{lemma:fixed-point-corep}.
    \end{proof}
\end{cor}

\begin{lemma}\label{lemma:ev_V-cc}
    In the situation of the previous corollary, both the functor $\ev_V\colon\myLSp_\textup{$G$-lvl}\to\myS_{\textup{$(\O(V)\times G)$-gl},*}$ and the composite $(\ref{eq:semifree-corep})$ are cocontinuous.
    \begin{proof}
        In light of Lemma \ref{lemma:id-g-fixed-points-vs-ev}, it suffices to prove the first statement. For this it is in turn enough to verify preservation of coproducts and pushouts.

        For coproducts, let $(X_i)_{i\in I}$ be any family of cofibrant objects of $\cat{$\bm G$-$\cat{L}$-Sp}$. By Corollary~\ref{cor:SpL-level-hcof}, each $X_i(V)$ is well-based, so that the inclusions $X_i(V)\to\bigvee_{i\in I}X_i(V)=\big(\bigvee_{i\in I}X_i\big)(V)$ exhibit the latter as a coproduct in the underlying $\infty$-category of $\cat{$\bm{(\O(V)\times G)}$-$\cat{L}$-Top}_*$ by Lemma~\ref{lemma:coprod-well-based}, as desired.

        The statement about pushouts follows in the same way from Corollary~\ref{cor:SpL-level-hcof} together with Lemma~\ref{lemma:h-cof-compute-po}(2).
    \end{proof}
\end{lemma}

\begin{proof}[Proof of Proposition~\ref{prop:mySplev-unpar-presh}]
    We verify the assumptions of Lemma~\ref{lemma:detect-pointed-presheaves}:
    \begin{enumerate}
        \item[(0)] As the underlying $\infty$-category of a model category, $\myLSp_\textup{$G$-lvl}$ is locally small and cocomplete.
        \item Cocontinuity of corepresented functors follows by combining Corollary~\ref{cor:ev_V-corep} and Lemma~\ref{lemma:ev_V-cc}.
        \item Joint conservativity is clear from Corollary~\ref{cor:ev_V-corep} and the definition of the $G$-global level weak equivalences.\qedhere
    \end{enumerate}
\end{proof}

\subsubsection{$G$-global $\Omega$-spectra} 
To obtain a model of $G$-global stable homotopy theory, we will again have to pass to a suitable Bousfield localization. For simplicity, we will only do this on the level of $\infty$-categories (as opposed to performing a Bousfield localization on the level of model categories, which would entail verifying additional technical assumptions).

\begin{definition}
    We call $X\in \cat{$\bm G$-$\cat{L}$-Sp}$ a \emph{$G$-global $\Omega$-spectrum} if the adjoint structure map
    $\tilde\sigma\colon X(V)\to\Omega^WX(V\oplus W)$
    is a $(G\times H)$-global weak equivalence for all compact Lie groups $H$ and all $H$-representations $V,W$ (not necessarily faithful).
\end{definition}

\begin{lemma}\label{lemma:omega-fewer-quantifiers}
    Let $X\in\cat{$\bm G$-$\cat{L}$-Sp}$ be levelwise closed, i.e. for all inner product spaces $V\in \cat{O}$, $X(V)$ is a closed orthogonal $G$-space. Then $X$ is a $G$-global $\Omega$-spectrum if and only if $\tilde\sigma(\Uu_H)^\phi\colon X(V)(\Uu_H)^\phi\to(\Omega^WX(V\oplus W)(\mathcal U_H))^\phi$ is a weak homotopy equivalence for every homomorphism of compact Lie groups $\phi\colon H\to G$, some (hence any) complete $H$-universe $\Uu_H$, and all $H$-representations $V,W$.
    \begin{proof}
        The assumptions guarantee that both $X(V)$ and $\Omega^WX(V\oplus W)$ are closed orthogonal $(G\times K)$-spaces for all $K$-representations $V,W$. Thus, Lemma~\ref{lemma:we-between-closed} shows that $X$ is a $G$-global $\Omega$-spectrum if and only if \[\tilde\sigma(\Uu_H)\colon X(V)(\Uu_H)\to\Omega^WX(V\oplus W)(\Uu_H)\]
        is a $\mathcal G_{H,K\times G}$-equivariant weak equivalence for all compact Lie groups $H,K$ and all $K$-representations $V,W$. As $\Uu_H$ is also a complete universe for every closed subgroup of $H$ (Remark~\ref{rk:universe-subgroup}), this is equivalent via replacing $H$ by its subgroups to $\tilde\sigma(\Uu_H)$ inducing a weak homotopy equivalence on $(\psi,\phi)$-fixed points for all compact Lie groups $H$ and $K$, $K$-representations $V$ and $W$, and homomorphisms $\psi\colon H\to K$ and $\phi\colon H\to G$. Restricting $V$ and $W$ to $H$-representations along $\psi$, this simplifies to the condition from the statement of the lemma, finishing the proof.
    \end{proof}
\end{lemma}

As a consequence of Lemma~\ref{lemma:cotensoring-fully-homotopical}, $G$-global $\Omega$-spectra are stable under $G$-global level weak equivalences. Thus, it makes sense to ask whether an object of the $\infty$-categorical localization $\myLSp_\text{$G$-lvl}$ is a $G$-global $\Omega$-spectrum.

\begin{definition}
    We write $\myLSp_{G}^{\Omega}\subset\myLSp_\textup{$G$-lvl}$ for the full subcategory spanned by the $G$-global $\Omega$-spectra.
\end{definition}

Below we will prove:

\begin{prop}\label{prop:L-gl-Bousfield}
    The inclusion $\myLSp_{G}^{\Omega}\hookrightarrow\myLSp_\textup{$G$-lvl}$ admits a left adjoint $L$, giving rise to an accessible Bousfield localization.
\end{prop}

\begin{definition}
    We call a map in $\myLSp_\text{$G$-lvl}$ a \emph{$G$-global weak equivalence} if it is sent to an equivalence under the Bousfield localization $L$. A map in $\cat{$\bm G$-$\cat{L}$-Sp}$ will be called a $G$-global weak equivalence if and only if its image in $\myLSp_\textup{$G$-lvl}$ is so.\footnote{A concrete description of the $G$-global weak equivalences will be given in Proposition~\ref{prop:he-preserve-but-he-also-reflec}.} We denote the localization of $\cat{$\bm G$-$\cat{L}$-Sp}$ at the $G$-global weak equivalences (which is equivalent to $\myLSp_G^\Omega$ in a preferred way) by $\myLSp_\text{$G$-gl}$.
\end{definition}

To prove Proposition~\ref{prop:L-gl-Bousfield}, we will explicitly construct a set $S$ such that the $G$-global $\Omega$-spectra are precisely the $S$-local objects.

\begin{constr}\label{constr:S}
    Let $H$ be a compact Lie group, let $V,W$ be $H$-{\hskip0pt}representations, and recall the map $\kappa_{V,W}^{H,H}\colon \cat{O}(V\oplus W,-)\smashp S^W\to \cat{O}(V,-)$ of orthogonal $H$-spectra from Construction~\ref{constr:generating-equiv-acyclic} corepresenting the $H$-fixed points $\smash{X(V)^H\to\big(\Omega^WX(V\oplus W)\big){}^H}$ of the adjoint structure map.
    We define $S$ as the set consisting of the maps $\smash{\bar\kappa_{V,W}^{(\phi)}\coloneqq\cat{L}\phi_!(\consts\,\kappa_{V,W}^{H,H})}$ in $\myLSp_\text{$G$-lvl}$ for all homomorphisms $\phi\colon H\to G$ and representations $V,W$ (up to isomorphism).
\end{constr}

\begin{lemma}
    An $X\in\myLSp_\textup{$G$-lvl}$ is $S$-local precisely if it is a $G$-global $\Omega$-spectrum.
    \begin{proof}
        By adjunction, $X$ is $S$-local if and only if $\phi^*X$ is local with respect to the maps $\smash{\consts\,\kappa_{V,W}^{H,H}}$ for all $\phi\colon H\to G$ and all representations $V,W$. If we fix a faithful $H$-representation $U$, then the collapse map $\cat{L}(U,-)_+\smashp\cat{O}(V,-)\to\cat{O}(V,-)$ is $G$-global level weak equivalence by Proposition~\ref{prop:smash-well-based} applied levelwise, and likewise for $\cat{L}(U,-)_+\smashp\cat{O}(V\oplus W,-)\smashp S^W\to\cat{O}(V\oplus W,-)\smashp S^W$; thus, being local with respect to the map $\kappa_{V,W}^{H,H}$ is equivalent to being local with respect to $\cat{L}(U,-)_+\smashp\kappa_{V,W}^{H,H}$ (in the $\infty$-categorical sense). We now observe that $\cat{L}(U,-)_+\smashp\cat{O}(V,-)$ is cofibrant in $\cat{$\bm H$-$\cat{L}$-Sp}$ by the explicit description of the generating cofibrations, as is $\cat{L}(U,-)_+\smashp\cat{O}(V\oplus W,-)$. Since $S^W$ is an $H$-CW-complex, the endofunctor $\Omega^W$ of $\cat{$\bm H$-$\cat{L}$-Sp}$ is right Quillen by Lemma~\ref{lemma:G-Top*-tensored} applied levelwise, so that ${-}\smashp S^W$ is left Quillen; we conclude that ${\cat{L}(U,-)_+\smashp\kappa_{V,W}^{H,H}}$ is a map of cofibrant objects. Lemma~\ref{lemma:localization-on-homs} then shows that a fibrant object of $\cat{$\bm H$-$\cat{L}$-Sp}$ is local in the topologically enriched sense with respect to $\smash{\cat{L}(U,-)_+\smashp\kappa_{V,W}^{H,H}}$ if and only if its image in $\myLSp_\text{$H$-gl}$ is local in the $\infty$-categorical sense with respect to $\smash{\kappa_{V,W}^{H,H}}$. 
        
        Fixing $V$ and $W$, but letting $U$ vary and passing to colimits, we see that an object of $\myLSp_\text{$H$-lvl}$ represented by a fibrant and levelwise closed $H$-$\cat{L}$-spectrum $Y$ is local with respect to $\consts\,\kappa_{V,W}^{H,H}$ if and only if $Y(V)(\Uu_H)^H\to \Omega^WY(V\oplus W)(\Uu_H)^H$ is a weak homotopy equivalence for our favourite complete $H$-universe $\Uu_H$. As the latter condition is invariant under $H$-global level weak equivalences of levelwise closed $H$-$\cat{L}$-spectra, we conclude that this more generally holds without the fibrancy assumption. Varying $\phi\colon H\to G$ and the $H$-representations $V$ and $W$, Lemma~\ref{lemma:omega-fewer-quantifiers} then shows that a levelwise closed $X\in\cat{$\bm G$-$\cat{L}$-Sp}$ defines an $S$-local object of $\myLSp_\text{$G$-gl}$ if and only if it is a $G$-global $\Omega$-spectrum. As any $X\in\cat{$\bm G$-$\cat{L}$-Sp}$ is weakly equivalent to a levelwise closed $G$-$\cat{L}$-spectrum (Corollary~\ref{cor:cofibrant-LSp-closed}), this completes the proof. 
    \end{proof}
\end{lemma}

\begin{proof}[Proof of Proposition~\ref{prop:L-gl-Bousfield}]
    By Proposition~\ref{prop:mySplev-unpar-presh}, the $\infty$-category $\myLSp_\text{$G$-lvl}$ is presentable. Thus, the claim follows from the previous lemma via the general theory of Bousfield localizations \cite{HTT}*{Proposition 5.5.4.15}. 
\end{proof}

\subsection{Comparison of the approaches} \label{subsec:comp-spectra}
In this subsection we will prove:

\begin{thm}\label{thm:comp-gl-models}
    Let $G$ be any compact Lie group. Then the functors 
    \[
        \consts\colon\Gosp{G}\to\cat{$\bm G$-$\cat{L}$-Sp}
        \qquad\text{and}
        \qquad
        \diag\colon\cat{$\bm G$-$\cat{L}$-Sp}\to\Gosp{G}
    \]
    descend to mutually inverse equivalences $\mySp_\textup{$G$-gl}\simeq\myLSp_\textup{$G$-gl}$.
\end{thm}

The key step in the proof will be showing that $\diag$ preserves $G$-global weak equivalences. We begin with the case of $G$-global \emph{level} weak equivalences:

\begin{lemma}
    The functor $\diag\colon\cat{$\bm G$-$\cat{L}$-Sp}\to\Gosp{G}$ sends $G$-global level weak equivalences to $G$-global weak equivalences and hence descends to $\myLSp_\textup{$G$-lvl}\to{\mySp}_\textup{$G$-gl}$.
    \begin{proof}
        It is clear that $\diag$ sends acyclic fibrations (i.e.~maps $X\to Y$ such that each $X(V)(U)\to Y(V)(U)$ is an acyclic fibration in the $\mathcal G_{\O(U),\O(V)\times G}$-model structure) to $G$-global level weak equivalences. By Corollary~\ref{cor:cofibrant-LSp-closed} it will therefore suffice to show that $\diag$ sends any $G$-global level weak equivalence $f\colon X\to Y$ such that $X$ and $Y$ are levelwise closed to a $G$-global weak equivalence. 
        
        Restricting along any homomorphism of compact Lie groups $\phi\colon H\to G$ (and replacing $G$ by $H$), it suffices to show that $\pi_n^G(\diag\,f)$ is bijective for every $n\in\Z$. We will give the argument for $n=0$; the argument for general $n$ is analogous, but requires slightly more notation. For this we consider the functor $\Pi_X\colon s(\Uu_G)\times s(\Uu_G)\to\cat{Sets}$ sending $(V,W)$ to $[S^V, X(V)(W)]_*^G$ with the obvious functoriality in $V$ and $W$, and define $\Pi_Y$ analogously. If we write $\Delta\colon s(\Uu_G)\to s(\Uu_G)\times s(\Uu_G)$ for the diagonal embedding, then $\pi_0^G(\diag\,X)=\colim\Pi_X\circ\Delta$ and $\pi_0^G(\diag\,Y)=\colim\Pi_Y\circ\Delta$ by definition. As $\Delta$ is cofinal, it will therefore suffice to show that the map
        \[
            \colim_{(V,W)\in s(\Uu_G)\times s(\Uu_G)}[S^V,X(V)(W)]_*^G \longrightarrow 
            \colim_{(V,W)\in s(\Uu_G)\times s(\Uu_G)}[S^V,Y(V)(W)]_*^G
        \]
        induced by $f$ is bijective, for which it is enough by the Fubini Theorem for colimits to prove bijectivity of the analogous map
        \[
            \colim_{V\in s(\Uu_G)}\colim_{W\in s(\Uu_G)}[S^V, X(V)(W)]_*^G
            \longrightarrow
            \colim_{V\in s(\Uu_G)}\colim_{W\in s(\Uu_G)}[S^V, Y(V)(W)]_*^G.
        \]
        Arguing as in the proof of Lemma~\ref{lemma:cotensoring-fully-homotopical}, this may be identified with the map
        \[
            \colim_{V\in s(\Uu_G)}[S^V, X(V)(\Uu_G)]_*^G\longrightarrow
            \colim_{V\in s(\Uu_G)}[S^V, Y(V)(\Uu_G)]_*^G
        \]
        induced levelwise by $f(V)(\Uu_G)\colon X(V)(\Uu_G)\to Y(V)(\Uu_G)$. As the latter is a $G$-equivariant weak equivalence, the claim follows.
    \end{proof}
\end{lemma}

\begin{lemma}\label{lemma:delta*-cc}
    The functor $\diag\colon \myLSp_\textup{$G$-lvl}\to{\mySp}_\textup{$G$-gl}$ is cocontinuous.
    \begin{proof}
        It suffices that the pointset level functor $\cat{$\bm G$-$\cat{L}$-Sp}\to\Gosp{G}$ preserves homotopy pushouts and homotopy coproducts.

        Any cofibration $f$ in $\cat{$\bm G$-$\cat{L}$-Sp}$ is an unbased level h-cofibration by Corollary~\ref{cor:SpL-level-hcof}, and so $\diag\,f$ is an unbased level cofibration in $\Gosp{G}$. Pushouts along the latter are homotopy pushouts by Theorem~\ref{thm:G-gl-Sp-model-struct}; as $\diag\colon\cat{$\bm G$-$\cat{L}$-Sp}\to\Gosp{G}$ evidently preserves pointset level pushouts, this shows that it also preserves homotopy pushouts.
        
        On the other hand, \cite{schwede2018global}*{Corollary~3.1.37(i)} shows that coproducts in $\Gosp{G}$ preserve $G$-global weak equivalences in full generality, so $\diag$ sends any coproduct of cofibrant objects to a homotopy coproduct.
    \end{proof}
\end{lemma}

The previous lemma already reduces the question whether $\diag$ also descends to $\myLSp_\text{$G$-gl}\to\mySp_\text{$G$-gl}$ to understanding its effect on the set $S$ from Construction~\ref{constr:S}. Looking at their definition, we should then first understand the interaction of $\diag$ with the left adjoints to restriction (see Remark~\ref{rk:LSp-Lf!}).

\begin{lemma}\label{lemma:delta*-gl-cc}
    Let $\alpha\colon G\to G'$ be any homomorphism of compact Lie groups. Then the Beck--Chevalley map $\BC_!\colon\cat{L}\alpha_!\circ\diag\to\diag\circ\cat{L}\alpha_!$ is an equivalence.    
    \begin{proof}
        As a consequence of Lemma~\ref{lemma:delta*-cc}, both source and target of $\BC_!$ are cocontinuous functors $\myLSp_\text{$G$-lvl}\to\mySp_\text{$G'$-gl}$. In light of the presheaf description of $\myLSp_\text{$G$-lvl}$ from Proposition~\ref{prop:mySplev-unpar-presh}, it will therefore suffice to show that the Beck--Chevalley map is an equivalence on the objects of the subcategory $\mySph_G$. It is clear that the Beck--Chevalley map is an equivalence for the zero object; every other object of $\mySph_G$ is then of the form $X\coloneqq\cat{L}\phi_!\cat{O}(V,-)$ for some $\phi\colon H\to G$ and some $H$-representation $V$. As Beck--Chevalley maps compose, we see that the map $\BC_!\colon\cat{L}\alpha_!\,\diag(X)\to\diag\,\cat{L}\alpha_!(X)$ in question is an equivalence provided that both the Beck--Chevalley maps
            $\cat{L}(\phi\alpha)_!\,\diag(\cat{O}(V,-))\to\diag\,\cat{L}(\phi\alpha)_!\cat{O}(V,-)$
           and
            $\cat{L}\phi_!\,\diag(\cat{O}(V,-))\to\diag\,\cat{L}\phi_!\cat{O}(V,-)$
        are so. In other words, up to replacing $\alpha$ and $G$, it suffices to treat the case where $X=\cat{O}(V,-)$ for some $G$-representation $V$.
        
        Fix a faithful $G$-representation $U$, so that the collapse map $\cat{L}(U,-)_+\smashp\cat{O}(V,-)\to\cat{O}(V,-)$ is a cofibrant replacement of $X$. Plugging in the definitions, the Beck--Chevalley condition then amounts to saying that for some (hence any) cofibrant replacement $f\colon Y\to\diag(\cat{L}(U,-)_+\smashp\cat{O}(V,-))$ in $\Gosp{G}$ the induced map 
        \begin{equation}\label{eq:BC-unravelled}
            \alpha_!Y\to\alpha_!\,\diag\big(\cat{L}(U,-)_+\smashp\cat{O}(V,-)\big)=\diag\,\alpha_!\big(\cat{L}(U,-)_+\smashp\cat{O}(V,-)\big)
        \end{equation} 
        is a $G'$-global weak equivalence. We may choose the cofibrant replacement already in the $G$-global level model structure, and we will show that in this case $(\ref{eq:BC-unravelled})$ is even a $G'$-global \emph{level} weak equivalence, i.e.\ that for every inner product space $W$ the $\mathcal G_{\O(W),G}$-weak equivalence $f(W)\colon Y(W)\to\cat{L}(U,W)_+\smashp\cat{O}(V,W)$ is sent by $\alpha_!$ to a $\mathcal G_{\O(W),G'}$-weak equivalence. By Proposition~\ref{prop:global-lvl-GSp}, $Y(W)$ is $\mathcal G_{\O(W),G}$-cofibrant, as is $\cat{L}(U,W)_+\smashp\cat{O}(V,W)$ by Lemma~\ref{lemma:smash-lQ}. Thus, the claim follows from Example~\ref{ex:graph-restr-rQ}.
    \end{proof}
\end{lemma}

\begin{cor}\label{cor:diag-preserve}
    The functor $\diag\colon\myLSp_\textup{$G$-lvl}\to\mySp_\textup{$G$-gl}$ inverts $G$-global weak equivalences.
    \begin{proof}
        By cocontinuity (Lemma~\ref{lemma:delta*-cc}), it suffices to show that $\diag$ inverts the maps in the set $S$ from Construction~\ref{constr:S}. By definition of the latter together with the previous lemma, it then suffices (after possibly changing the compact Lie group $G$) to show that all maps of the form $\smash{\consts\,\kappa_{V,W}^{G,G}}$ for $G$-representations $V,W$ are inverted. As $\diag\circ\consts=\id_{\Gosp{G}}$, this is an instance of Proposition~\ref{prop:kappa-gl}.
    \end{proof}
\end{cor}

For the proof of the converse we will need:

\begin{lemma}\label{lemma:I'm-so-static}
    Let $X\in\cat{$\bm G$-$\cat{L}$-Sp}$ be level fibrant and a $G$-global $\Omega$-spectrum.
    \begin{enumerate}
        \item If $U,U'$ are $H$-representations and $V$ is any faithful $H$-representation, then $\tilde\sigma(V)\colon X(U)(V)\to \Omega^{U'}X(U\oplus U')(V)$ is $\mathcal G_{H,G}$-weak equivalence.
        \item If $U$ is any $H$-representation and $i\colon V\to W$ is an embedding of faithful $H$-representations, then $X(U)(i)\colon X(U)(V)\to X(U)(W)$ is a $\mathcal G_{H,G}$-equivariant weak equivalence.
    \end{enumerate}
    \begin{proof}
        For the first claim, note that $\tilde\sigma\colon X(U)\to X(U\oplus U')$ is a $(G\times H)$-global weak equivalence, and both its source and target are fibrant orthogonal $(G\times H)$-spaces by Lemma~\ref{lemma:restr-right-Quillen} and the definition of the $G$-global level fibrations. Thus, $\tilde\sigma$ is even a $(G\times H)$-global \emph{level} weak equivalence, so that $\tilde\sigma(V)$ is a $\mathcal G_{\O(V),G\times H}$-weak equivalence and hence in particular a $\mathcal G_{H,G}$-weak equivalence for the diagonal $H$-action.

        For the second claim, we similarly note that $X(U)(i)$ is a $\mathcal G_{H,G\times \O(U)}$-weak equivalence by fibrancy of $X(U)$, hence in particular a $\mathcal G_{H,G}$-weak equivalence.
    \end{proof}
\end{lemma}

\begin{cor}
    If $X\in\cat{$\bm G$-$\cat{L}$-Sp}$ is fibrant and a $G$-global $\Omega$-spectrum, then also $\diag\,X\in\Gosp{G}$ is a $G$-global $\Omega$-spectrum.
    \begin{proof}
        If $H$ is any compact Lie group and $V,W$ are $H$-representations such that $V$ is faithful, then the adjoint structure map $(\diag\,X)(V)\to\Omega^W(\diag\,X)(V\oplus W)$ factors as the composite
        \[
            X(V)(V)\xrightarrow{\;X(V)(i)\;} X(V)(V\oplus W)\xrightarrow{\;\tilde\sigma(V\oplus W)}\Omega^WX(V\oplus W)(V\oplus W)   
        \]
        where $i\colon V\hookrightarrow V\oplus W$ denotes the summand inclusion. By the previous lemma, both of these maps are $\mathcal G_{H,G}$-weak equivalences, so the claim follows by 2-out-of-3.
    \end{proof}
\end{cor}

\begin{prop}\label{prop:he-preserve-but-he-also-reflec}
    The functor $\diag\colon\cat{$\bm G$-$\cat{L}$-Sp}\to\Gosp{G}$ preserves and reflects $G$-global weak equivalences.
    \begin{proof}
        We may equivalently show that the induced functor $\myLSp_\text{$G$-lvl}\to\mySp_\text{$G$-gl}$ precisely inverts the $G$-global weak equivalences.
       
        That it inverts all $G$-global weak equivalences was established in Corollary~\ref{cor:diag-preserve}. For the converse, let $f\colon X\to Y$ be any $G$-global weak equivalence in $\myLSp_\text{$G$-lvl}$ and consider the naturality square
        \[
            \begin{tikzcd}
                X\arrow[d,"f"']\arrow[r,"\eta"] & iLX\arrow[d,"iLf"]\\
                Y\arrow[r,"\eta"'] & iLY
            \end{tikzcd}
        \]
        for the unit of the Bousfield localization $L\colon\myLSp_\text{$G$-lvl}\rightleftarrows\myLSp_\text{$G$}^\Omega\noloc i$. As we have just seen, $\diag$ sends the horizontal maps to equivalences in $\mySp_\text{$G$-gl}$. Thus, if $\diag\,f$ is a $G$-global weak equivalence, then so is $\diag\,iLf$ by 2-out-of-3; conversely, to show that $f$ is a $G$-global weak equivalence, it will suffice to show that $iLf$ is so. In other words, we may assume without loss of generality that $X$ and $Y$ are $G$-global $\Omega$-spectra. After cofibrant-fibrant replacement in the $G$-global level model structure, we may then assume that $f$ is represented by a map $X\to Y$ in $\cat{$\bm G$-$\cat{L}$-Sp}$ between cofibrant-fibrant $G$-global $\Omega$-spectra; we will denote this pointset level map by $f$ again, and we will similarly view $\diag\,f$ as a map in $\Gosp{G}$.

        By the previous corollary, both $\diag\,X$ and $\diag\,Y$ are $G$-global $\Omega$-spectra, so that $\diag\,f$ is a $G$-global \emph{level} weak equivalence. We will now complete the proof of the proposition by showing that $f(V)$ is an $(\O(V)\times G)$-global \emph{level} weak equivalence of orthogonal $(\O(V)\times G)$-spaces for any $V$, i.e.\ $f(V)(W)$ is a $\mathcal G_{H,G}$-equivariant weak equivalence for every compact Lie group acting arbitrarily on $V$ and faithfully on $W$. To prove this, we consider the commutative diagram
        \[
            \begin{tikzcd}
                X(V)(W)\arrow[d,"f(V)(W)"']\arrow[r,"X(V)(i)"] &[1.5em] X(V)(V\oplus W)\arrow[d,"f(V)(V\oplus W)"{description}]\arrow[r,"\tilde\sigma(V\oplus W)"]&[2em]\Omega^WX(V\oplus W)(V\oplus W)\arrow[d,"\Omega^Wf(V\oplus W)(V\oplus W)"]\\
                Y(V)(W)\arrow[r,"Y(V)(i)"'] & Y(V)(V\oplus W)\arrow[r,"\tilde\sigma(V\oplus W)"'] & \Omega^WY(V\oplus W)(V\oplus W)
            \end{tikzcd}
        \]
        By Lemma~\ref{lemma:I'm-so-static}, the horizontal maps are $\mathcal G_{H,G}$-weak equivalences. On the other hand, $f(V\oplus W)(V\oplus W)$ is a $\mathcal G_{H,G}$-weak equivalence by assumption on $\diag\,f$, whence so is $\Omega^Wf(V\oplus W)(V\oplus W)$ by Lemma~\ref{lemma:smash-lQ}. The claim follows by 2-out-of-3.
    \end{proof}
\end{prop}

\begin{cor}\label{cor:LSp-restr-homotopical}
    Let $\alpha\colon G\to G'$ be any homomorphism of compact Lie groups. Then $\alpha^*\colon\cat{$\bm{G'}$-$\cat{L}$-Sp}\to\cat{$\bm G$-$\cat{L}$-Sp}$ sends $G'$-global weak equivalences to $G$-global weak equivalences.\qed
\end{cor}

\begin{proof}[Proof of Theorem~\ref{thm:comp-gl-models}]
    As $\diag\circ\consts=\id_{\Gosp{G}}$, Proposition~\ref{prop:he-preserve-but-he-also-reflec} shows that these functors descend to the localizations at the $G$-global weak equivalences, with $\consts$ being right inverse to $\diag$. To prove that $\consts$ is also left inverse (thereby completing the proof of the theorem), we will now exhibit a pointset level zig-zag of natural transformations between the identity of $\cat{$\bm G$-$\cat{L}$-Sp}$ and the composite $\consts\circ\diag$. For this we define $F\colon\cat{$\bm G$-$\cat{L}$-Sp}\to\cat{$\bm G$-$\cat{L}$-Sp}$ by $F(X)(V)(W)\coloneqq X(V)(V\oplus W)$ with the obvious functoriality in each variable. The two summand inclusions $V\hookrightarrow V\oplus W\hookleftarrow W$ then induce natural transformations $i\colon\consts\,\diag\to F\gets\id\noloc j$. We will show that $j$ is a weak equivalence; the argument for $i$ is then analogous.

    By Proposition~\ref{prop:he-preserve-but-he-also-reflec}, it will suffice that $\diag\,j\colon\diag\,X\to\diag\,F(X)$ is a $G$-global weak equivalence for every $X$. If we write $\sh\colon\cat{$\cat{L}$-Top}_*\to\cat{$\cat{L}$-Top}_*$ for the functor given by precomposition with $V\mapsto V\oplus V$, then $\diag\,F(X)$ agrees with $\diag(\sh\circ X)$ by direct inspection, and the map $j$ is induced by the natural transformation $k\colon\id\to\sh$ given by the inclusion of the second summand. By \cite{barrero2021}*{Lemma 3.8}, $k\circ X\colon X\to\sh\circ X$ is a $G$-global level weak equivalence, so that $\diag\,j=\diag(k\circ X)$ is a $G$-global weak equivalence by another application of  Proposition~\ref{prop:he-preserve-but-he-also-reflec}.
\end{proof}

\subsection{The global \texorpdfstring{$\bm\infty$}{∞}-category of global spectra}\label{subsec:gl-oo-gl-sp}
We will now explain how the two pointset models of $G$-global spectra discussed above assemble into representation stable globally presentable global $\infty$-categories.

\begin{constr}
    Define $\ul\mySp_\gl$ as the levelwise localization of $\Ntop(\osp^\dual)$ at the $G$-global weak equivalences; this is well-defined as the restriction functors are homotopical by definition of the weak equivalences. We similarly define $\ulmySplev$ by localizing $\Ntop(\cat{$\cat{L}$-Sp}^\dual)$ at the $G$-global level weak equivalences and $\ul\myLSp_\gl$ by instead localizing at the $G$-global weak equivalences (which is well-defined by Corollary~\ref{cor:LSp-restr-homotopical}).
\end{constr}

Theorem~\ref{thm:comp-gl-models} immediately implies:

\begin{cor}\label{cor:LSp-vs-Sp-param}
    The topological functor $\consts\colon\osp\to\cat{$\cat{L}$-Sp}$ induces an equivalence of global $\infty$-categories $\ul\mySp_\gl\iso\ul\myLSp_\gl$.\qed
\end{cor}

\subsubsection{Presentability} As in the unstable situation, it is not clear how to verify the Beck--Chevalley condition for global cocompleteness by hand. Instead, our proof will proceed in two steps: we will first show that $\ulmySplev$ is globally presentable using presentability of $\ul\myS_{\gl,*}$, and then exhibit $\ul\mySp_\gl\simeq\ul\myLSp_\gl$ as a (parametrized) accessible Bousfield localization.

\begin{constr}
    Let $G$ be a compact Lie group and let $V$ be any $G$-representation. Evaluation at $V$ defines a topologically enriched homotopical functor $\cat{$\bm H$-$\cat{L}$-Sp}\to\cat{$\bm{(H\times G)}$-$\cat{L}$-Top}_*$ for every compact Lie group $H$. Applying $\Ntop\big((-)^\dual\big)$ and localizing, we therefore obtain a global functor 
    \[
        \ev_V\colon\ulmySplev\to\ul\myS_{\gl,*}(\BGcat{G}\times-).
    \]
\end{constr}

The definition of the $G$-global level weak equivalences immediately implies:

\begin{lemma}\label{lemma:ev-V-jointly-repr}
    The functors $\ev_V$ are jointly conservative when $G$ runs through all compact Lie groups and $V$ through all $G$-representations.\qed
\end{lemma}

\begin{lemma}\label{lemma:ev-V-param-la}
    The functor $\ev_V\colon\ulmySplev\to\ul\myS_{\gl,*}(\BGcat{G}\times-)$ is a global left adjoint.
    \begin{proof}
        By Lemma~\ref{lemma:ev_V-cc} and presentability of $\ulmySplev(\BGcat{H})=\myLSp_\text{$H$-lvl}$ (Proposition~\ref{prop:mySplev-unpar-presh}), each of the individual functors $\ev_V\colon\ulmySplev(\BGcat{H})\to\ul\myS_{\gl,*}(\BGcat{G}\times\BGcat{H})$ admits a right adjoint, so it will suffice to show that these right adjoints satisfy the Beck--Chevalley condition with respect to restriction along any $\alpha\colon H\to H'$ in $\Glo$ or, equivalently, that the Beck--Chevalley map $\cat{L}\alpha_!\ev_V\to\ev_V\cat{L}\alpha_!$ is invertible for any such $\alpha$, where $\cat{L}\alpha_!\colon\myLSp_\text{$H$-lvl}\to\myLSp_\text{$H'$-lvl}$ is as in Remark~\ref{rk:LSp-Lf!}.
        
        Note that both source and target of the latter Beck--Chevalley map are left adjoints, so the collection of objects for which it is invertible is closed under all colimits. By the presheaf description of $\myLSp_\text{$H$-lvl}$ from Proposition~\ref{prop:mySplev-unpar-presh}, it will therefore suffice to verify invertibility for the objects of the form $X\coloneqq\cat{L}\phi_!\cat{O}(W,-)$ for homomorphisms $\phi\colon K\to H$ and $K$-representations $W$. Arguing as in the proof of Lemma~\ref{lemma:delta*-gl-cc}, one reduces to the case $\phi=\id$, so that $X$ is modelled by the cofibrant object $\cat{L}(U,-)_+\smashp\cat{O}(W,-)$ where $U$ is our favourite faithful $H$-representation. If $f\colon Y\to \cat{L}(U,-)_+\smashp\cat{O}(W,V)$ is any cofibrant replacement in the $(G\times H)$-global level model structure on $\cat{$\bm{(G\times H)}$-$\cat{L}$-Top}_*$, then the Beck--Chevalley map in question is represented on the pointset level by $\alpha_!f$. We claim that the latter is a $(G\times H')$-global level weak equivalence, which we can reexpress as saying that for any inner product space $T$, any $K\subset\O(T)$, and any homomorphism $\psi\colon K\to G$ the map $\psi^*\alpha_!f(T)=\alpha_!\psi^*f(T)$ is a $\mathcal G_{K,H'}$-equivariant weak equivalence.
        By definition, $Y(T)$ is $\mathcal G_{K,G\times H}$-cofibrant, so that $\psi^*Y(T)$ is $\mathcal G_{K,H}$-cofibrant by Lemma~\ref{lemma:restr-equiv-lQ}. Moreover, also $\cat{L}(U,T)_+\smashp\cat{O}(W,\psi^*V)$ is $\mathcal G_{K,H}$-cofibrant by Lemma~\ref{lemma:smash-lQ} and faithfulness of the $H$-representation $U$. Thus, Example~\ref{ex:graph-restr-rQ} shows that $\alpha_!$ sends the $\mathcal G_{K,H}$-weak equivalence $\psi^*f(T)$ to a $\mathcal G_{K,H'}$-equivariant weak equivalence, as desired.
    \end{proof}
\end{lemma}

\begin{prop}\label{prop:Splev-cc}
    The global $\infty$-category $\ulmySplev$ is globally presentable.
    \begin{proof}
        We begin by proving that $\ulmySplev$ is fiberwise presentable. For this we note once more that if $G$ is any compact Lie group, then $\ulmySplev(\BGcat{G})$ is presentable by Proposition~\ref{prop:mySplev-unpar-presh}. Let now $V$ be any $H$-representation; by the previous lemma $\ev_V\colon\ulmySplev(\BGcat{G})\to\myS_{\textup{$(G\times H)$-gl},*}$ is cocontinuous. Thus, in the naturality square
        \[
            \begin{tikzcd}
                \ulmySplev(\BGcat{G'})\arrow[d,"\ev_V"']\arrow[r,"f^*"] & \ulmySplev(\BGcat{G})\arrow[d,"\ev_V"]\\
                \myS_{\text{$(G'\times H)$-gl},*}\arrow[r,"f^*"'] & \myS_{\text{$(G\times H)$-gl},*} 
            \end{tikzcd}
        \]
        all maps except possibly the top one preserve colimits. As for varying $V$ the right-hand vertical maps are jointly conservative, also $f^*\colon\ulmySplev(\BGcat{G'})\to\ulmySplev(\BGcat{G})$ is cocontinuous, hence a left adjoint by the Adjoint Functor Theorem.
        
        Next, we recall once more that each $f^*\colon\ulmySplev(\BGcat{G}')\to\ulmySplev(\BGcat{G})$ admits a left adjoint $\cat{L}f_!$, whence so does more generally any restriction along a map $f\colon X\to\BGcat{G}$ in $\PSh(\Glo)$ by Remark~\ref{rk:cleft-cc-addendum}. It remains to verify the Beck--Chevalley conditions from Lemma~\ref{lemma:S-cocompleteness-repr}, i.e.\ that for any pullback
        \[
            \begin{tikzcd}
                X'\arrow[r, "f'"]\arrow[d, "g'"']\arrow[dr,pullback] & \BGcat{H}\arrow[d, "g"]\\
                X\arrow[r, "f"'] & \BGcat{G}
            \end{tikzcd}
        \]
        in $\PSh(\Glo)$ the Beck--Chevalley map $\BC_!\colon\cat{L}f_!'g'^*\to g^*\cat{L}f_!
        $ is invertible. By Lemma~\ref{lemma:ev-V-jointly-repr} the functors $\ev_V$ are jointly conservative, so it will suffice to show that the Beck--Chevalley map is an equivalence after hitting it with $\ev_V$ for every representation $V$. But we have a commutative diagram of Beck--Chevalley maps
        \[
            \begin{tikzcd}
                \ev_V\cat{L}f_!'g'^*\arrow[r,"\BC_!"]\arrow[d,"\ev_V\BC_!"'] & \cat{L}f_!'\ev_Vg'^*\arrow[r,"\sim"] & \cat{L}f'_!g'^*\ev_V\arrow[d,"\BC_!"]\\
                \ev_Vg^*\cat{L}f_!\arrow[r,"\sim"'] & \ev_Vg^*\cat{L}f_!\arrow[r,"g^*\BC_!"'] & g^*\cat{L}f_!\ev_V
            \end{tikzcd} 
        \]
        in which the two horizontal Beck--Chevalley maps are equivalences as $\ev_V$ is a left adjoint, while the right-hand vertical map is an equivalence since $\ul\myS_{\gl,*}$ is globally cocomplete (Corollary~\ref{cor:gl-pointed}). We conclude by 2-out-of-3 that also the left-hand vertical map is invertible, finishing the proof of the proposition.
    \end{proof}
\end{prop}

\begin{lemma}
    The localization $\ulmySplev\to\ul\myLSp_\gl$ has a fully faithful right adjoint.
    \begin{proof}
        If $G$ is any compact Lie group, then Proposition~\ref{prop:L-gl-Bousfield} shows that the localization functor $\myLSp_\text{$G$-lvl}=\ulmySplev(\BGcat{G})\to\ul\myLSp_\gl(\BGcat{G})=\mySp_\text{$G$-gl}$ admits a right adjoint, which is fully faithful with essential image given by the $G$-global $\Omega$-spectra. It therefore only remains to verify the Beck--Chevalley condition. By full faithfulness of the right adjoints, this reduces to proving that the restriction functors preserve the essential images, which is clear from the definition of $G$-global $\Omega$-spectra. 
    \end{proof}
\end{lemma}

\begin{thm}\label{thm:gl-sp-pres}
    Both $\ul\myLSp_\gl$ and $\ul\mySp_\gl$ are globally presentable.
    \begin{proof}
        In light of Corollary~\ref{cor:LSp-vs-Sp-param}, it will suffice to prove the first statement. 

        By Proposition~\ref{prop:L-gl-Bousfield}, each individual category $\myLSp_\text{$G$-gl}$ is presentable, so it remains to show global cocompleteness. For this, let us write $R\colon\ul\myLSp_\gl\to\ulmySplev$ for the right adjoint from the previous lemma. If $X\in\PSh(\Glo)$ is any presheaf, then $R(X)\colon\ul\myLSp_\gl(X)\to\ulmySplev(X)$ is a fully faithful right adjoint (as a limit of such functors), so its left adjoint is a Bousfield localization. Moreover, if $f\colon X\to Y$ is any map in $\PSh(\Glo)$, then $f_!\colon\ulmySplev(X)\to\ulmySplev(Y)$ descends through the localization functor as $f^*$ preserves the essential image of the fully faithful right adjoint $R$ (the latter being part of a global functor); in particular, $f^*\colon\ul\mySp_\gl(Y)\to\ul\mySp_\gl(Y)$ has a left adjoint, and it only remains to verify the Beck--Chevalley condition for pullbacks in $\PSh(\Glo)$. As localizations are in particular essentially surjective, this may be done after precomposing with $\ulmySplev\to\ul\myLSp_\gl$. Using once more that Beck--Chevalley maps compose, this then follows at once from $\ulmySplev$ being globally cocomplete (Proposition~\ref{prop:Splev-cc}) and the localization being a left adjoint.
    \end{proof}
\end{thm}

Recall the homotopical functor $\Sigma^\bullet_+\colon\cat{$\bm G$-$\cat{L}$-Top}\to\Gosp{G}$ from Proposition~\ref{prop:Sigma-bullet-homotopical-left-Quillen}. Varying $G$, this yields a global functor $\Sigma^\bullet_+\colon\ul\myS_{\gl}\to\ul\mySp_\gl$; for later use, we record:

\begin{lemma}\label{lemma:sigma-bullet-left-adj}
    The global functor $\Sigma^\bullet_+\colon\ul\myS_\gl\to\ul\mySp_\gl$ is a left adjoint.
    \begin{proof}
        By the aforementioned Proposition \ref{prop:Sigma-bullet-homotopical-left-Quillen}, $\Sigma^\bullet_+$ is left Quillen, proving the existence of a pointwise right adjoint. To prove the Beck--Chevalley condition for these right adjoints it then suffices to observe that also the restriction functors $\Gosp{G'}\to\Gosp{G}$ are right Quillen (Lemma~\ref{lemma:restr-gl-sp-right-Quillen}) and that the Beck--Chevalley condition holds for the underived functors (in fact, for the standard choice of adjunctions and counits, the Beck--Chevalley map will even be the identity).
    \end{proof}
\end{lemma}

\subsubsection{Representation stability}
In order to prove representation stability, we will explicitly describe the tensoring of $\ul\mySp_\gl$ over pointed equivariant spaces. 

\begin{prop}
    The levelwise smash product
    \begin{equation}\label{eq:levelwise-smash-before-loc}
        \Ntop(\cat{Top}_*^\textup{cof,$\dual$})\times\Ntop(\osp^\dual)\to\Ntop(\osp^\dual)
    \end{equation}    
    descends to $-\otimes-\colon\ul\myS_*\times\ul\mySp_\gl\simeq \ul\myS_*^\textup{cof}\times\ul\mySp_\gl\to\ul\mySp_\gl$, and this is the unique bifunctor preserving $\Orb$-colimits in the first variable and such that $S^0\otimes-$ is the identity.
    \begin{proof}
        Lemma~\ref{lemma:G-global-spectra-tensoring} shows that the levelwise smash product descends accordingly and that each individual $\myS_{G,*}^\text{cof}\times\mySp_\text{$G$-gl}\to\mySp_\text{$G$-gl}$ preserves colimits in the first variable. It remains to show that for any injective $\alpha\colon G\to G'$, any cofibrant $G$-space $X$, and any $G'$-global spectrum $Y$ the Beck--Chevalley map $\cat{L}\alpha_!(X\otimes \alpha^*Y)\to (\cat{L}\alpha_!X)\otimes Y$ is invertible. In light of Lemma~\ref{lemma:gl-sp-(co)ind-homotopical}, we may equivalently show that the pointset level map $\alpha_!(X\smashp \alpha^*Y)\to(\alpha_!X)\smashp Y$ is a weak equivalence. A direct computation shows that the latter is even an isomorphism.
    \end{proof}
\end{prop}

In other words, the canonical tensoring of $\ul\mySp_\gl$ over $\ul\myS_*$ is given (up to cofibrant replacement) by the levelwise smash product. Proposition~\ref{prop:resp-stability-G-gl} therefore translates to the following:

\begin{cor}\label{cor:Spgl-rep-stable}
    The global $\infty$-category $\ul\mySp_\gl$ is representation stable.\qed
\end{cor}

\subsubsection{Symmetric monoidal structure} Once we know the universal property of $\ul\mySp_\gl$, the general theory will yield a globally presentably symmetric monoidal structure with unit $\mathbb S$, and shows that this unique. We will now construct such a symmetric monoidal structure by hand, using the smash product of orthogonal spectra; it will then later follow a posteriori that this agrees with the symmetric monoidal structure coming from the universal property.

\begin{constr}
    For any compact Lie group $G$, we write $\Gosp{G}^\text{flat}\subset\Gosp{G}$ for the full (topologically enriched) subcategory spanned by the flat $G$-orthogonal spectra in the sense of Definition~\ref{def:flat}. These assemble into a global subcategory $\Ntop(\osp^\text{flat,$\dual$})\subset\Ntop(\osp^\dual)$, and we write $\ul\mySp_\gl^\text{flat}$ for its levelwise localization at the $G$-global weak equivalences.
\end{constr}

\begin{lemma}\label{lemma:flat-resolution}
    The map $\ul\mySp_\gl^\textup{flat}\to\ul\mySp_\gl$ induced by the inclusion is an equivalence.
    \begin{proof}
        If $G$ is any compact Lie group, then Lemma~\ref{lemma:cofibrant-G-gl-are-flat} implies that functorial cofibrant replacement in the $G$-global model structure defines a homotopy inverse to the inclusion $\Gosp{G}^\text{flat}\hookrightarrow\Gosp{G}$.
    \end{proof}
\end{lemma}

\begin{constr}\label{constr:global-smash-product-param}
    By \cite{equiv-smash}*{Proposition 3.5.3}, the enriched symmetric monoidal structure on $\Gosp{G}$ given by the smash product restricts to $\Gosp{G}^\text{flat}$, and by Proposition~\ref{prop:gl-flatness} this descends (for varying $G$) to a symmetric monoidal structure on $\ul\mySp_\gl^\text{flat}\simeq\ul\mySp_\gl$. We will write $\smash{\ul\mySp_\gl^\otimes}$ for the resulting symmetric monoidal global $\infty$-category.
\end{constr}

We will now construct a symmetric monoidal refinement of $\Sigma^\bullet_+\colon\ul\myS_\gl\to\ul\mySp_\gl$. Note that if we already knew $\smash{\ul\mySp_\gl^\otimes}$ to be globally presentably symmetric monoidal, we would obtain a unique symmetric monoidal left adjoint from the general theory; we will conversely need the symmetric monoidal structure on $\Sigma^\bullet_+$ to establish that the smash product on $\ul\mySp_\gl$ preserves global colimits in each variable.

\begin{constr}\label{constr:sigma-bullet-sym-mon}
    Recall the universal property of the Day convolution product as the recipient of the universal bimorphism $\beta$. Given orthogonal spaces $X,Y$, we then have maps
    \begin{multline*}
        \hskip-5pt(\Sigma^\bullet_+X)(V)\smashp(\Sigma^\bullet_+Y)(W) = S^V\smashp X(V)_+\smashp S^W\smashp Y(W)_+
        \cong S^{V\oplus W}\smashp\big(X(V)\times Y(W)\big){}_+\\\xrightarrow{\;S^{V\oplus W}\smashp(\beta_{V,W})_+\;} S^{V\oplus W}\smashp(X\boxtimes Y)(V\oplus W)_+
        =(\Sigma^\bullet_+ X\boxtimes Y)(V\oplus W)
        \hskip-5pt
    \end{multline*}
    for all $V,W\in\cat{O}$; by \cite{schwede2018global}*{Proposition 4.1.18}, these define a bimorphism, the induced map $(\Sigma^\infty_+X)\smashp(\Sigma^\infty_+Y)\to \Sigma^\infty_+(X\boxtimes Y)$ is an isomorphism, and together with the unique isomorphism $\Sigma^\bullet_+1\cong\mathbb S$ this upgrades $\Sigma^\bullet_+\colon\cat{$\cat{L}$-Top}\to\osp$ to a symmetric monoidal functor.     
    Applying the symmetric monoidal continuous Borel construction, we obtain a symmetric monoidal global functor $\Ntop(\cat{$\cat{L}$-Top}^{\otimes,\dual})\to\Ntop(\osp^{\otimes,\dual})$ restricting to $\smash{\Ntop(\cat{$\cat{L}$-Top}^{\text{qcof},\otimes,\dual})\to\Ntop(\osp^{\text{flat},\otimes,\dual})}\vphantom{S^f}$, where $\Ntop(\cat{$\cat{L}$-Top}^{\text{qcof},\otimes,\dual})$ denotes the $\Glo$-subcategory on quasi-cofibrant orthogonal based $G$-spaces (Construction~\ref{constr:smash-product-unstable}). Localizing then yields the desired symmetric monoidal refinement
    \[
        \Sigma^\bullet_+\colon\ul\myS_\gl^\otimes\simeq\ul\myS_\gl^{\text{qcof},\otimes}\to
        \ul\mySp_\gl^{\text{flat},\otimes}\simeq\ul\mySp_\gl^\otimes.
    \]
\end{constr}

\begin{variant}\label{var:reduced-global-susp-spectrum}
    Analogously, one upgrades the reduced suspension spectrum functor $\Sigma^\bullet\colon\cat{$\cat{L}$-Top}_*\to\cat{OrthSp}$ to a symmetric monoidal functor, and hence obtains a symmetric monoidal global functor
    \[
        \Sigma^\bullet\colon\ul\myS_{\gl,*}^\otimes\simeq\ul\myS_{\gl,*}^{\text{qcof},\otimes}\to\ul\mySp_\gl^{\text{flat},\otimes}\simeq\ul\mySp_\gl^\otimes.
    \]
    Once the next proposition has been proven, the universal property of $\ul\myS_{\gl,*}^\otimes$ (Theorem~\ref{thm:gl-smash-initial}) will imply that $\Sigma^\bullet$ is the unique symmetric monoidal global left adjoint with the given source and target.
\end{variant}

\begin{prop}\label{prop:gl-spectra-pres-sym-mon}
    The symmetric monoidal global $\infty$-category $\ul\mySp_\gl^\otimes$ is globally presentably symmetric monoidal.
    \begin{proof}
        For any compact Lie group $G$, the symmetric monoidal product on $\mySp_\text{$G$-gl}$ is equivalently given as the left derived functor of the left Quillen bifunctor from Corollary~\ref{cor:smash-gl-left-Quillen}, and hence preserves colimits in each variable. As we already know that $\ul\mySp_\gl$ is globally presentable, it therefore only remains to show that the Beck--Chevalley map $\cat{L}\alpha_!(X\otimes \alpha^*Y)\to(\cat{L}\alpha_!X)\otimes Y$ is invertible for any $\alpha\colon G\to G'$ and $X\in\mySp_\text{$G$-gl},Y\in\mySp_\text{$G'$-gl}$. As both sides preserve ordinary colimits in both variable separately and since both $\mySp_\text{$G$-gl}$ and $\mySp_\text{$G'$-gl}$ are stable in the usual sense, we reduce to the case where $X$ and $Y$ are generators as in Proposition~\ref{prop:gl-spectra-compact-generators}, and hence in particular contained in the essential image of the symmetric monoidal left adjoint $\Sigma^\bullet_+$ from the previous construction. Using that Beck--Chevalley maps compose as in the proof of Proposition~\ref{prop:equiv-sp-equiv-cc}, the claim then follows from Corollary~\ref{cor:S-gl-times-initial}.
    \end{proof}
\end{prop}

\subsubsection{Equivariant vs. global spectra redux} Finally, we discuss the relation to the global $\infty$-category $\ul\mySp$ of \emph{equivariant spectra} from the previous section:

\begin{prop}\label{prop:mySp-vs-mySp-gl}
    The inclusion of equivariantly cofibrant orthogonal $G$-spectra induces a fully faithful global functor $\ul\mySp\to\ul\mySp_\gl$. Moreover, its underlying $\Orb$-functor admits a right adjoint which is levelwise a localization at the equivariant weak equivalences.
    \begin{proof}
        Lemma~\ref{lemma:equiv-spectra-vs-gl-spectra} already shows that $\ul\mySp\to\ul\mySp_\gl$ is fully faithful and admits a pointwise right adjoint which is a localization at the equivariant weak equivalences.
        To complete the proof, it remains to verify the Beck--Chevalley condition for the pointwise right adjoints with respect to \emph{injective} homomorphisms. As the right adjoints are localizations, this amounts to saying that $\alpha^*\colon\mySp_\text{$G'$-gl}\to\mySp_\text{$G$-gl}$ sends $G'$-equivariant weak equivalences to $G$-equivariant ones for any injective $\alpha\colon G\to G'$. This is however immediate from Proposition~\ref{prop:i!-left-Quillen-equiv-sp}.
    \end{proof}
\end{prop}

By the universal property of $\ul\mySp$, the first half can be rephrased as saying that the unique equivariantly cocontinuous functor $\ul\mySp\to\ul\mySp_\gl$ sending the sphere to the sphere is fully faithful. The additional information about the right adjoint provided by the proposition will become relevant later in the proof of Theorem~\ref{thm:Thom-vs-fgt}.

\subsection{Proof of the universal property}\label{subsec:gl-sp-univ} In this subsection we will finally prove the following strenghtening of Theorem~\ref{introthm:stable-main} from the introduction:

\begin{thm}\label{thm:stable-main}
    The pair $(\ul\mySp_\gl,\mathbb S)$ is an idempotent algebra in $\CAT_{\Glo}^\textup{$\Glo$-cc}$. The $\ul\mySp_\gl$-modules are precisely the representation stable globally cocomplete global $\infty$-categories. In particular, we have for any such $\Dd$ an equivalence
    \[
        \ev_{\mathbb S}\colon\ul\Fun^\textup{$\Glo$-cc}(\ul\mySp_\gl,\Dd)\iso\Dd.
    \]
\end{thm}

The proof will be given below after some preparations. The key step for this will be showing that $\ev_{\mathbb S}\colon\Fun^\textup{L}(\ul\myLSp_\gl,\Dd)\to\Dd$ admits a fully faithful right adjoint for any representation stable globally presentable $\Dd$. We begin by proving the analogous statement for $\ulmySplev$, which actually holds in larger generality.

\begin{lemma}\label{lemma:omega-oo-both-adjoints}
    The evaluation functor $\Omega^\infty\coloneqq\ev_0\colon\ulmySplev\to\ul\myS_{\gl,*}$ admits both a fully faithful left adjoint $\cat{L}\Sigma^\infty$ as well as a fully faithful right adjoint.
    \begin{proof}
        The existence of the right adjoint is a special case of Lemma~\ref{lemma:ev-V-param-la}. For the existence of the left adjoint, observe that $\ev_0\colon\cat{$\bm G$-$\cat{L}$-Sp}\to\cat{$\bm G$-$\cat{L}$-Top}_*$ is right Quillen by Lemma~\ref{lemma:ev-RA}, with left adjoint given by the usual suspension spectrum functor. It then only remains to verify the Beck--Chevalley condition. Unravelling definitions, this amounts to showing that for any $\alpha\colon G\to G'$, any cofibrant $X\in\cat{$\bm{G'}$-$\cat{ L}$-Top}_*$ and any cofibrant replacement $f\colon Y\to \alpha^*X$, the induced map $\Sigma^\infty f\colon \Sigma^\infty Y\to \Sigma^\infty \alpha^*X$ is a $G$-global level weak equivalence, i.e.\ for every $V$ the map $S^V\smashp f$ is an $(O(V)\times G)$-global weak equivalence. This is however simply an instance of Corollary~\ref{cor:global-smash-well-based}.

        Finally, it is clear that the ordinary unit $X\to\Omega^\infty \Sigma^\infty X$ is an isomorphism for every (cofibrant) $X$; thus, the derived unit $\id\to\Omega^\infty\,\cat{L}\Sigma^\infty$ is invertible, i.e.\ $\cat{L}\Sigma^\infty$ is fully faithful. 
        It follows formally (see e.g.~\cite{haine2025fullyfaithfulfunctorspushouts}*{Lemma 1.2}) that the right adjoint of $\Omega^\infty$ is also fully faithful.
    \end{proof}
\end{lemma}

\begin{prop}\label{prop:Sigma-oo-lev-ff}
    Let $\Cc$ be any pointed globally cocomplete global $\infty$-category. Then
    $\ev_{\mathbb S}\colon\Fun^\textup{L}(\ulmySplev,\Cc)\to\Cc(1)$ admits a both a fully faithful left adjoint $\mathbb S_!$ and a fully faithful right adjoint $\mathbb S_*$.
    \begin{proof}
        By the previous lemma and 2-functoriality of $\Fun^\text{L}(-,\Cc)$, the restriction functor $-\circ\cat{L}\Sigma^\infty\colon\Fun^\text{L}(\ulmySplev,\Cc)\to\Fun^\text{L}(\ul\myS_{\gl,*},\Cc)$ has a fully faithful left adjoint given by ${-}\circ\Omega^\infty$ (here we are using that $\Omega^\infty$ is also a \emph{left} adjoint). Combining this with the universal property of $\ul\myS_{\gl,*}$ from Corollary~\ref{cor:gl-pointed}, we conclude that $\ev_{\mathbb S}$ has a fully faithful left adjoint, and it only remains to show that it also admits a right adjoint, which will then automatically be fully faithful as before.

        For this we observe that $\Fun^\text{L}(\ulmySplev,\Cc)\hookrightarrow\Fun(\ulmySplev,\Cc)$ preserves colimits by \cite{martiniwolf2022presentable}*{Lemma 2.6.1.3}, so that $\ev_{\mathbb S}\colon\Fun^\text{L}(\ulmySplev,\Cc)\to\Cc(1)$ is cocontinuous. As the source is presentable (Proposition~\ref{prop:FunL-presentable}), the existence of a right adjoint then follows from the (non-parametrized) Adjoint Functor Theorem.
    \end{proof}
\end{prop}

An analogous argument already shows that $\ev_{\mathbb S}\colon\Fun^\text{L}(\ul\myLSp_\gl,\Cc)\to\Cc(1)$ admits a right adjoint for any pointed globally presentable global $\infty$-category $\Cc$; however, this adjoint will usually not be fully faithful, i.e.\ we might not have a Bousfield localization. For \emph{representation stable} $\Cc$, we will now show that in fact the Bousfield localization $\Fun^\text{L}(\ulmySplev,\Cc)\rightleftarrows\Cc(1)$ from the previous proposition restricts accordingly:

\begin{prop}\label{prop:RA-factor-through-LA}
    Let $\Cc$ be globally presentable and representation stable. Then the right adjoint $\mathbb S_*\colon\Cc(1)\to\Fun^\textup{L}(\ulmySplev,\Cc)$ of $\ev_{\mathbb S}$ lands in the full subcategory of those left adjoints factoring over the localization $\ulmySplev\to\ul\myLSp_\gl$.
\end{prop}

In Construction~\ref{constr:S}, we already explicitly described a set $S$ such that $\myLSp_\text{$G$-lvl}\to\myLSp_\text{$G$-gl}$ is a Bousfield localization at $S$, and varying $G$ this allows us to characterize the left adjoints $\ulmySplev\to\Cc$ factoring through $\ul\myLSp_\gl\to\Cc$ as those functors inverting an explicit set of maps. For the proof of the proposition, it will be convenient to recast this criterion in terms of the invertibility of certain Beck--Chevalley maps.

\begin{constr}\label{constr:sigma-as-BC}
    Let $G$ be a compact Lie group, and let $V,W$ be $G$-representations. For any $\cat L$-spectrum $X$, we get a natural map of based orthogonal $G$-spaces $\tilde\sigma\colon X(V)\to\Omega^WX(V\oplus W)$ as the adjoint of the structure map, yielding a natural transformation of topologically enriched functors $\cat{$\bm{\cat L}$-Sp}\rightrightarrows\cat{$\bm G$-${\cat L}$-Top}_*$; note moreover that both functors are homotopical, see Lemma~\ref{lemma:cotensoring-fully-homotopical} for the non-trivial case. Applying $\Ntop((-)^\dual)$ and passing to localizations we obtain a natural transformation
    \[
        \begin{tikzcd}
            \ulmySplev\arrow[from=d,"="]\arrow[dr,Rightarrow,"\tilde\sigma"{description},shorten=1em]\arrow[r,"\ev_V"] &[1em] \ul\myS_{\gl,*}(\BGcat{G}\times-)\arrow[from=d,"\Omega^W"']\\
            \ulmySplev\arrow[r,"\ev_{V\oplus W}"']&\ul\myS_{\gl,*}(\BGcat{G}\times-)\rlap.
        \end{tikzcd}
    \]
    We write $\sigma$ for the associated Beck--Chevalley map $\Sigma^W\ev_{V}\to\ev_{V\oplus W}$. Plugging in the definitions, this is given on any well-based $X\in\cat{$\bm G$-$\cat{L}$-Sp}$ by the usual structure map $\Sigma^WX(V)\to X(V\oplus W)$.
\end{constr}

\begin{lemma}
    Let $V$ be a $G$-representation. Then $\ev_V\colon\ulmySplev\to\ul\myS_{\gl,*}({\BGcat{G}\times{-}})$ admits a left adjoint $\cat{O}(V,-)\otimes_G{-}$.
    \begin{proof}
        As $\ulmySplev$ is globally cocomplete, the functor $\ulmySplev\to\ulmySplev(\BGcat{G}\times-)$ equipping an $H$-global spectrum with the trivial $G$-action admits a left adjoint. Thus, it will suffice to show that the functor $\ulmySplev(\BGcat{G}\times-)\to\ul\myS_{\gl,*}(\BGcat{G}\times-)$ induced by the functor $\ev_V\colon\cat{$\bm G$-$\cat{L}$-Sp}\rightarrow\cat{$\bm G$-$\cat{L}$-Top}$ sending $X$ to $X(V)$ with the diagonal $G$-action admits a left adjoint. For every compact Lie group $H$, the induced functor $\ev_V\colon\cat{$\bm{(G\times H)}$-$\cat{L}$-Sp}\to\cat{$\bm{(G\times H)}$-$\cat{L}$-Top}_*$ is right Quillen by Lemmas~\ref{lemma:ev-RA} and~\ref{lemma:restr-right-Quillen}, and its left adjoint $\cat{O}(V,-)\smashp{-}$ sends $(G\times H)$-global weak equivalences between well-based $(G\times H)$-orthogonal spaces to $(G\times H)$-global level weak equivalences by Corollary~\ref{cor:global-smash-well-based}. The same argument as in the proof Lemma~\ref{lemma:omega-oo-both-adjoints} (which considered the special case $G=1$ and $V=0$) therefore shows that the individual left adjoints $\myS_{\text{$(G\times H)$-gl},*}\to\myLSp_\text{$(G\times H)$-gl}$ satisfy the Beck--Chevalley condition and hence assemble into the desired global left adjoint.
    \end{proof}
\end{lemma}

\begin{constr}
    Let $V,W$ be $G$-representations. We write 
    \[
        \kappa\colon\cat{O}(V\oplus W,-)\otimes_G\Sigma^W(-)\to\cat{O}(V,-)\otimes_G{-}
    \]
    for the Beck--Chevalley transformation associated to the natural transformation $\sigma$ from Construction~\ref{constr:sigma-as-BC}. In other words, $\kappa$ is the total mate of the natural transformation $\tilde\sigma$.
\end{constr}

\begin{lemma}\label{lemma:kappa-bousfield}
    Let $\Cc$ be globally presentable. A left adjoint $F\colon\ulmySplev\to\Cc$ factors through $\ul\myLSp_\gl$ if and only if $F$ inverts the natural transformations $\kappa$ from the previous construction for all compact Lie groups $G$ and representations $V,W$.
    \begin{proof}
        The left adjoint $F$ factors through the localization if and only if its right adjoint $U\colon\Cc\to\ulmySplev$ factors through the full subcategory $\ul\myLSp^\Omega$. By definition, this precisely amounts to requiring that the pastings
        \[
            \begin{tikzcd}
                \Cc\arrow[from=d,"="]\arrow[r, "U"]&\ulmySplev\arrow[r,"\ev_V"]\arrow[from=d,"="{description}]\arrow[dr,Rightarrow,"\tilde\sigma"{description},shorten=1em] &[2em] \ul\myS_{\gl,*}(\BGcat{G}\times-)\\
                \Cc\arrow[r, "U"']& \ulmySplev\arrow[r,"\ev_{V\oplus W}"'] & \ul\myS_{\gl,*}(\BGcat{G}\times-)\arrow[u,"\Omega^W"']
            \end{tikzcd}
        \]
        are equivalences for all compact Lie groups $G$ and $G$-representations $V,W$. However, as a natural transformation is an equivalence if and only if its total mate is so and since passing to mates is compatible with pastings, this is in turn equivalent to the pastings
        \[
            \begin{tikzcd}
                \Cc\arrow[d,"="']\arrow[from=r, "F"']&\ulmySplev\arrow[from=r,"{\cat{O}(V,-)\otimes_G{-}}"']\arrow[d,"="{description}]\arrow[from=dr,Rightarrow,"\kappa"{description},shorten=1em] &[4em] \ul\myS_{\gl,*}(\BGcat{G}\times-)\\
                \Cc\arrow[from=r, "F"]& \ulmySplev\arrow[from=r,"\;{\cat{O}(V\oplus W,-)\otimes_G{-}}"] & \ul\myS_{\gl,*}(\BGcat{G}\times-)\arrow[from=u,"\Sigma^W"]
            \end{tikzcd}
        \]
        being equivalences, which precisely amounts to $F$ inverting the maps $\kappa$.
    \end{proof}
\end{lemma}

\begin{proof}[Proof of Proposition~\ref{prop:RA-factor-through-LA}]
    Let us fix a compact Lie group $G$ together with $G$-representations $V,W$, and consider the diagram
        \begin{equation}\label{diag:ess-im-S*}
            \begin{tikzcd}
                \Cc(1)&[1em]\arrow[l,"="'] \Cc(1)\\
                \Fun^\textup{L}(\ulmySplev,\Cc)\arrow[u,"\ev_{\mathbb S}"] & \arrow[l,"="'] \Fun^\textup{L}(\ulmySplev,\Cc)\arrow[from=dl,Rightarrow,"\sigma^{*}"{description}]\arrow[u,"\ev_{\mathbb S}"']\\
                \arrow[u,"-\circ\ev_V"] \Fun^\textup{L}(\ul\myS_{\gl,*}({\ul{\BGcat{G}}}\times{-}),\Cc)&\arrow[l,"-\circ\Sigma^W"] \Fun^\textup{L}(\ul\myS_{\gl,*}({\ul{\BGcat{G}}}\times{-}),\Cc)\arrow[u,"-\circ\ev_{V\oplus W}"']\rlap.
            \end{tikzcd}
        \end{equation}
        By (homotopy) 2-functoriality of $\Fun^\text{L}$, the two lower vertical maps have right adjoints, given by precomposition with the \emph{left} adjoints to evaluations, and the Beck--Chevalley map for the lower square is given on an $F\colon\ulmySplev\to\Cc$ by applying $F$ to the Beck--Chevalley map associated to $\sigma$, i.e.\ to $\kappa$. As Beck--Chevalley maps compose, we then conclude from the previous lemma that the right adjoint $\mathbb S_*$ factors through $\Fun^\text{L}(\ul\myLSp_\gl,\Cc)$ if and only if the Beck--Chevalley map associated to the pasting $(\ref{diag:ess-im-S*})$ is invertible for all compact Lie groups $G$ and all representations $V,W$.
        
        The pasting $(\ref{diag:ess-im-S*})$ is given on a left adjoint $F\colon\ul\myS_{\gl,*}(\BGcat{G}\times-)\to\Cc$ by applying $F$ to the map $\sigma_{\mathbb S}$ in $\ul\myS_{\gl,*}(\BGcat{G})$. By the explicit description from Construction~\ref{constr:sigma-as-BC}, $\sigma_{\mathbb S}$ is invertible, so that the pasting $(\ref{diag:ess-im-S*})$ is an equivalence. Thus, if we can show that the top and bottom horizontal maps are equivalences, then also the Beck--Chevalley map will be invertible for purely formal reasons.

        For the top map, there is nothing to prove. For the bottom map, it will suffice (after replacing $\Cc$ by $\Fun([n],\Cc)$ for varying $n\ge0$) to prove that the map is an equivalence on groupoid cores. We write $\fgt\colon\Glo_{/\BGcat{G}}\to\Glo$ for the forgetful functor, and we consider the bifunctor
        \[
            \Theta\colon\fgt^*\ul\myS_{\gl,*}\times\fgt^*\ul\myS_{\gl,*}=\fgt^*(\ul\myS_{\gl,*}\times\ul\myS_{\gl,*})\xrightarrow{\;\fgt^*(-{\otimes}-)}\fgt^*\ul\myS_{\gl,*}
        \]
        induced by the smash product, so that the effect of the bottom map in $(\ref{diag:ess-im-S*})$ on groupoid cores is given by 
        \[
            {-}\circ\fgt_*\Theta(S^W,-)\colon\hom_{\PrL_{\Glo}}(\fgt_*\fgt^*\ul\myS_{\gl,*},\Cc)\to \hom_{\PrL_{\Glo}}(\fgt_*\fgt^*\ul\myS_{\gl,*},\Cc).
        \]
        To see that this is an equivalence, it will suffice by Corollary~\ref{cor:PrLT-ambidextrous-fgt} that
        \begin{equation}\label{eq:precompose-Sigma}
            \hskip2pt{-}\circ\Theta(S^W,-)\colon \hom_{\PrL_{\Glo_{/\BGcat{G}}}}\hskip0pt minus 4pt(\fgt^*\ul\myS_{\gl,*}\hskip0pt minus 1pt,\fgt^*\Cc)\to \hom_{\PrL_{\Glo_{/\BGcat{G}}}}\hskip0pt minus 4pt(\fgt^*\ul\myS_{\gl,*}\hskip0pt minus 1pt,\fgt^*\Cc)\hskip2pt
        \end{equation}
        is so.  Combining Examples~\ref{ex:based-global-spaces-tensored} and~\ref{ex:tensoring-slice}, we see that $\Theta$ defines the canonical tensoring of the pointed cocomplete $\Glo_{/\BGcat{G}}$-$\infty$-category $\fgt^*\ul\myS_{\gl,*}$. Thus, Lemma~\ref{lemma:tensoring-natural} shows that $(\ref{eq:precompose-Sigma})$ is homotopic to
        \[
            (S^W\otimes{-})\circ{-}\colon \hom(\fgt^*\ul\myS_{\gl,*},\fgt^*\Cc)\to\hom(\fgt^*\ul\myS_{\gl,*},\fgt^*\Cc).
        \]
        This map is then indeed an equivalence by representation stability of $\Cc$.
\end{proof}

\begin{cor}
    Let $\Cc$ be globally presentable and representation stable. Then $\ev_{\mathbb S}\colon\Fun^\textup{L}(\ul\myLSp_\gl,\Cc)\to\Cc(1)$ is the left adjoint in a Bousfield localization.
    \begin{proof}
        As a consequence of the previous proposition, $\ev_{\mathbb S}$ has a right adjoint given by the composite
        \[
            \Cc(1)\xrightarrow{\;\mathbb S_*\;}\Fun'(\ulmySplev,\Cc)\iso\Fun^\text{L}(\ul\myLSp_\gl,\Cc),
        \]
        where $\Fun'(\ulmySplev,\Cc)$ denotes the full subcategory of left adjoints factoring through the localization $L\colon\ulmySplev\to\ul\myLSp_\gl$ and the unnamed equivalence is given by precomposing with a right adjoint to $L$. As $\mathbb S_*$ is fully faithful by Proposition~\ref{prop:Sigma-oo-lev-ff}, the claim follows.
    \end{proof}
\end{cor}

\begin{proof}[Proof of Theorem~\ref{thm:stable-main}]
    Since $\smash{\ul\mySp_\gl}$ is representation stable (Corollary~\ref{cor:Spgl-rep-stable}) and globally presentable (Theorem~\ref{thm:gl-sp-pres}), there is a unique left adjoint \[\ul\Sp_{\smash\Glo}^{\text{RepSph}}\to\ul\mySp_\gl\] sending $\mathbb S$ to $\mathbb S$, and we may equivalently show that this functor is an equivalence.

    Comparing corepresented functors on $\smash{\Pr^{\text{L},\text{RepSph-ex}}_\Glo}$, it will suffice to show that $\ev_{\mathbb S}\colon\Fun^\text{L}(\ul\mySp_\gl,\Cc)\to\Cc$ is an equivalence for any representation stable and presentable $\Cc$. By the previous corollary together with Corollary~\ref{cor:LSp-vs-Sp-param}, $\ev_{\mathbb S}$ is a Bousfield localization; we will complete the proof by showing it is conservative.
    
    For this let $\tau\colon F_1\to F_2$ be any morphism in $\Fun^\text{L}(\ul\mySp_\gl,\Cc)$ such that $\tau(\mathbb S)$ is an equivalence in $\Cc(1)$. By the universal property of $\ul\myS_\gl$, this guarantees that the natural transformation $\tau\circ\Sigma^\bullet_+\colon F_1\circ\Sigma^\bullet_+\to F_2\circ\Sigma^\bullet_+$ of left adjoint functors $\ul\myS_\gl\to\Cc$ is an equivalence. In particular, if $G$ is any compact Lie group, then the natural transformation $\tau_G$ between $F_1,F_2\colon\mySp_\text{$G$-gl}\rightrightarrows\Cc(\BGcat{G})$ is an equivalence on the compact generators from Proposition~\ref{prop:gl-spectra-compact-generators}. As $\Cc(\BGcat{G})$ is in particular stable in the usual non-parametrized sense, this implies that $\tau_G$ is an equivalence, as desired.
\end{proof}

\begin{cor}\label{cor:global-spectra-smash}
    The global $\infty$-category $\ul\mySp_\gl$ admits a unique globally presentably symmetric monoidal structure with unit $\mathbb S$, and this is the one induced by the smash product (Construction~\ref{constr:global-smash-product-param}). Moreover, this defines the initial representation stable globally presentably symmetric monoidal global $\infty$-category.
    \begin{proof}
        The symmetric monoidal structure from Construction~\ref{constr:global-smash-product-param} is globally presentably symmetric monoidal by Proposition~\ref{prop:gl-spectra-pres-sym-mon}. The remaining statements then follow from Corollary~\ref{cor:SpR-initial} together with Theorem~\ref{thm:stable-main}.
    \end{proof}
\end{cor} 

\begin{rk}\label{rk:not-naive-stab}
    In the previous section, we saw that the $\Orb$-cocomplete representation stabilization $\Sigma^\infty\colon\ul\myS^\otimes_*\to\ul\mySp^\otimes$ is a pointwise construction: for every compact Lie group $G$, the symmetric monoidal functor $\Sigma^\infty\colon\myS_{G,*}^\otimes\to\mySp_{G}^\otimes$ is given by universally inverting $G$-representation spheres in ${\CAlg(\PrL)}$. We will now show that the corresponding statement is \emph{not} true for the $\Glo$-cocomplete representation stabilization $\Sigma^\bullet\colon\ul\myS_{\gl,*}^\otimes\to\ul\mySp_\gl^\otimes$ by showing that $\myS_{\gl,*}^\otimes\to\mySp_\gl^\otimes$ is not simply given by inverting the non-equivariant spheres.

    By \cite{CHLL_NRings}*{Lemma 5.5.3}, we may equivalently show that $\Sigma^\bullet\colon\myS_{\gl,*}\to\mySp_\gl$ is not a stabilization in $\smash{\PrL}$. For this we will write $B_\gl(\Z/2)\coloneqq\cat{L}(\R[\Z/2],-)/(\Z/2)$. Specializing the computation from \cite{schwede2018global}*{Corollary 4.1.13}, we see that the group $\pi_0\hom(\Sigma^\bullet B_\gl(\Z/2)_+,\Sigma^\bullet B_\gl(\Z/2)_+)$ (with the group structure coming from additivity of $\mySp_\gl$) is free abelian of rank 3. On the other hand, $B_\gl(\Z/2)$ corresponds under our equivalence to $\PSh(\Glo)$ to the represented presheaf $\BGcat{(\Z/2)}$. If $L\colon\PSh(\Glo)_*\to\Cc$ denotes the stabilization, a simple computation shows that 
    \begin{equation}\label{eq:hom-naive}
        \pi_0\hom_{\Cc}(L\BGcat{(\Z/2)}_+,L\BGcat{(\Z/2)}_+)\cong \pi_0\big(\Sigma^{\infty}_+\hom_{\Glo}(\BGcat{(\mathbb Z/2)},\BGcat{(\mathbb Z/2)})\big).
    \end{equation}
    As $\hom_{\Glo}(\BGcat{(\mathbb Z/2)},\BGcat{(\mathbb Z/2)})$ has exactly two path components (corresponding to the identity homomorphism and the trivial homomorphism), the abelian group $(\ref{eq:hom-naive})$ has rank $2\not=3$.
\end{rk}

\begin{rk}\label{rk:gepner-nikolaus}
    In \cite{gepner-nikolaus}, Gepner and Nikolaus introduced a notion of genuine stability in the setting of $\myS_\gl$-modules in $\PrL$, which is closely related to representation stability: they declare an $\myS_\gl$-module $\Cc$ to be \emph{$\myS_\gl$-stable} if it is pointed and the relative tensor product ${-}\otimes_{\myS_\gl}\Cc$ inverts $\Sigma^V\colon\myS_{\text{$G$-gl},*}\to\myS_{\text{$G$-gl},*}$ for every compact Lie group $G$ and every $G$-representation $V$ (where we view $\myS_{\text{$G$-gl},*}$ as an $\myS_\gl$-module via the symmetric monoidal functor $\triv_G\circ(-)_+\colon\myS_\gl^\times\to\myS_{\gl,*}^\otimes\to\myS_{\text{$G$-gl},*}^\otimes$). They moreover announced that $\mySp_\gl$ is the free $\myS_\gl$-stable $\myS_\gl$-module on one generator; however, as they explained to us, their unpublished argument contained a gap. Let us compare this (hypothetical) universal property to the one established in Theorem~\ref{thm:stable-main} above:

    As $\smash{\ul\myS_\gl}$ is the unit for the parametrized Lurie tensor product, the lax symmetric monoidal evaluation functor $\Gamma\coloneqq\ev_1\colon\PrL_\Glo\to\PrL$ lifts to a functor $\Gamma^\text{lin}$ to $\myS_\gl$-modules, and by \cite{martiniwolf2022presentable}*{Propositions~2.6.3.2--2.6.3.7} this lift admits a fully faithful symmetric monoidal left adjoint given explicitly by $\ul\myS_\gl\otimes_{\myS_\gl}{-}$. A short computation using Proposition~\ref{prop:Spc*-levelwise-tensor} shows that $\ul\myS_{\gl,*}\simeq\ul\myS_\gl\otimes_{\myS_\gl}\myS_{\gl,*}$ and that an $\myS_\gl$-module $\Cc$ is $\myS_\gl$-stable if and only if the globally presentable global $\infty$-category $\ul\myS_\gl\otimes_{\myS_\gl}\Cc$ is representation stable in our sense. As we will now explain, Theorem~\ref{thm:stable-main} reduces the question whether $\mySp_\gl$ is the $\myS_\gl$-stabilization of $\myS_\gl$ in the sense of Gepner and Nikolaus to the question whether $\ul\mySp_\gl$ is contained in the essential image of the inclusion of $\myS_\gl$-modules. We do not know whether the latter condition holds, though we strongly suspect it does not.
    
    Assume first that $\ul\mySp_\gl$ lies in the essential image. Then $\ul\mySp_\gl\simeq\ul\myS_\gl\otimes_{\myS_\gl}\mySp_\gl$ via the counit, so $\mySp_\gl$ is $\myS_\gl$-stable and the universal property as an $\myS_\gl$-module just follows from the universal property of $\ul\mySp_\gl$ and full faithfulness of $\ul\myS_\gl\otimes_{\myS_\gl}{-}$. Conversely, assume that $\mySp_\gl=\Gamma^\text{lin}(\ul\mySp_\gl)$ is the free $\myS_\gl$-stable $\myS_\gl$-module. If $\Cc$ is any representation stable global stable $\infty$-category, then $\Gamma^\text{lin}(\Cc)$ is $\myS_\gl$-stable as it admits an $\mySp_\gl$-module structure induced by the unique $\ul\mySp_\gl$-module structure on $\Cc$. By adjunction, the representation stable globally presentable global $\infty$-category $\ul\myS_\gl\otimes_{\myS_\gl}\mySp_\gl$ therefore corepresents the same functor on $\Pr_{\Glo}^\text{L,RepSph-ex}$ as $\ul\mySp_\gl$, so that $\ul\mySp_\gl\simeq\ul\myS_\gl\otimes_{\myS_\gl}\mySp_\gl$ as claimed.
\end{rk}

\chapter{Global Thom spectra}
Reaping the rewards of our work in the previous sections, we will now give two constructions of a symmetric monoidal global \emph{Thom spectrum functor}: an $\infty$-categorical construction as a natural parametrized analogue of \cite{thom-oo} and a model-categorical construction as a global refinement of \cite{sagave-schlichtkrull-thom}. As our final main result (Theorem~\ref{thm:comparison-Thom}) we will then show that these two constructions are equivalent. 
\section{The parametrized global Thom spectrum functor}\label{sec:parametrized-thom}
In \cite{thom-oo}*{Definition~2.20}, Ando, Blumberg, Gepner, Hopkins, and Rezk introduced an $\infty$-categorical Thom spectrum functor as the unique left adjoint
\[
    \Th_{\mathbb S}\colon\Spc_{/\pic(\Sp)}\to\Sp
\]
extending the inclusion $i\colon\pic(\Sp)\hookrightarrow\Sp$; this admits a natural symmetric monoidal enhancement induced by the evident symmetric monoidal structure on the inclusion $i$.
Restricting this along the functor $\Spc_{/\BO\times\mathbb Z}\to\Spc_{/\pic(\Sp)}$ induced by the \emph{$J$-homomorphism} $J\colon\BO\times \mathbb Z\to \pic(\Sp)$ (which we will recall in Construction~\ref{constr:definitionjhomomorphism}), one recovers the classical Thom spectrum functor for virtual vector bundles, which can equivalently be described as the unique cocontinuous (symmetric monoidal) extension of $J$.

In this section, we introduce a globally parametrized version of the Thom spectrum functor. We begin in §\ref{subsec:slicecats} with some preliminaries on parametrized symmetric monoidal cocompletion. §\ref{subsec:J} is devoted to the construction of a global refinement of the $J$-homomorphism. In  §\ref{subsec:global-Thom}, we finally construct the globally parametrized Thom spectrum functor and equivariant variations thereof.

\subsection{Symmetric monoidal cocompletion}\label{subsec:slicecats}
If $\Cc^\otimes$ is a symmetric monoidal $T$-$\infty$-category and $\eta\colon\Cc\to\mathcal P(\Cc)$ a $T$-cocompletion of the underlying $T$-$\infty$-category of $\Cc^\otimes$, then we already saw in Remark~\ref{rk:sym-mon-cocompletion} that there is a unique pair of an enhancement of $\mathcal P(\Cc)$ to an object $\mathcal P(\Cc)^\otimes\in\CAlg(\CAT_T^\text{$T$-cc})$ together with an enhancement $\eta^\otimes$ of $\eta$ to a symmetric monoidal $T$-functor, and that $\eta^\otimes$ is then initial among symmetric monoidal $T$-functors into an object of $\CAlg(\CAT_T^\text{$T$-cc})$. 
Specializing this to the description of the $T$-cocompletion from Theorem~\ref{thm:univ-prop-PSh} we get:

\begin{cor}
    Let $\Cc^\otimes$ be a small symmetric monoidal $T$-$\infty$-category. Then there exists a unique way to equip $\ul\PSh_T(\Cc)$ with a $T$-presentably symmetric monoidal structure and to make the Yoneda embedding $y\colon\Cc\to\ul\PSh_T(\Cc)$ into a $T$-symmetric monoidal functor $y^\otimes$. Moreover, if $\Dd^\otimes$ is any $T$-presentably symmetric monoidal $T$-$\infty$-category, then restriction along $y^\otimes$ defines an equivalence
    \[
        \Fun^{\textup{L},\otimes}(\ul\PSh_T(\Cc^\otimes),\Dd^\otimes)\iso
        \Fun^\otimes(\Cc^\otimes,\Dd^\otimes).\qednow
    \]
\end{cor}

\subsubsection{Day convolution} It is immediate from the construction that the unit for the symmetric monoidal structure on $\ul\PSh_T(\Cc^\otimes)$ is given by the Yoneda image of the unit of $\Cc^\otimes$. We will now also describe the symmetric monoidal product on $\ul\PSh_T(\Cc)^\otimes$ fairly concretly in analogy with the usual formula for Day convolution. For this, let us define for any small $T$-$\infty$-categories $\Cc,\Dd$ the `external product' \[{-}\exttimes{-}\colon\ul\PSh_T(\Cc)\times\ul\PSh_T(\Dd)\to\ul\PSh_T(\Cc\times\Dd)\] as the composite
\[
    \ul\PSh_T(\Cc)\times\ul\PSh_T(\Dd)\xrightarrow{\;\pr_1^*\times\pr_2^*\;}
    \ul\PSh_T(\Cc\times\Dd)\times\ul\PSh_T(\Cc\times\Dd)\xrightarrow{\;{-}\times{-}\;}\ul\PSh_T(\Cc\times\Dd).
\]

\begin{prop}\label{prop:Day-convolution}
    Let $\Cc^\otimes$ be a small symmetric monoidal $T$-$\infty$-category. Then the symmetric monoidal product on $\ul\PSh_T(\Cc^\otimes)$ is equivalent to the composite
    \begin{equation}\label{eq:Day-convolution}
        \ul\PSh_T(\Cc)\times\ul\PSh_T(\Cc)\xrightarrow{\;{-}\exttimes{-}\;}\ul\PSh_T(\Cc\times\Cc)\xrightarrow{\;\otimes_!\;}\ul\PSh_T(\Cc),
    \end{equation}
    where the second functor is left Kan extension along the tensor product of $\Cc$.
    \begin{proof}
        The functors $\pr_1^*$ and $\pr_2^*$ are left adjoints, as is $\otimes_!$. Since the product on $\ul\PSh_T(\Cc\times\Cc)$ is $T$-bilinear as a consequence of Proposition~\ref{prop:spc-T-cart-pres-sym-mon}, we see that $(\ref{eq:Day-convolution})$ is $T$-bilinear. Moreover, in the diagram
        \[
            \begin{tikzcd}
                \Cc\times\Cc\arrow[d,"y_\Cc\times y_\Cc"']\arrow[r,equals] & \Cc\times\Cc\arrow[r,"{-}\otimes{-}"]\arrow[d,"y_{\Cc\times\Cc}"{description}]&\Cc\arrow[d,"y"]\\
                \ul\PSh_T(\Cc)\times\ul\PSh_T(\Cc)\arrow[r,"{-}\exttimes{-}"'] & \ul\PSh_T(\Cc\times\Cc)\arrow[r,"\otimes_!"'] & \ul\PSh_T(\Cc)
            \end{tikzcd}
        \]
        the left-hand square commutes by \cite{martini2022cocartesianfibrationsstraighteninginternal}*{Lemma~6.2.6} and the right-hand square does so by \cite{martiniwolf2021limits}*{Corollary~3.3.3}. Altogether we see that $(\ref{eq:Day-convolution})$ is a $T$-bilinear functor extending the tensor product on $\Cc$. As the same holds for the tensor product on $\ul\PSh_T(\Cc^\otimes)$ by definition, the claim now follows from the universal property of $\ul\PSh_T(\Cc)$ as the free (non-monoidal) $T$-cocompletion.
    \end{proof}
\end{prop}

If $f^\otimes\colon\Cc^\otimes\to\Dd^\otimes$ is a strong symmetric monoidal $T$-functor and $\Cc,\Dd$ are small, then it follows from the universal property that there is a unique symmetric monoidal left adjoint $f_!^\otimes$ fitting into a commutative square
\[
    \begin{tikzcd}
        \Cc^\otimes\arrow[d,"f^\otimes"']\arrow[r,"y^\otimes"] & \ul\PSh_T(\Cc^\otimes)\arrow[d,dashed,"f^\otimes_!"]\\
        \Dd^\otimes\arrow[r,"y^\otimes"'] & \ul\PSh_T(\Dd^\otimes)\rlap.
    \end{tikzcd}
\]
Its underlying $T$-functor is then simply the unique left adjoint completing the corresponding square of non-monoidal $T$-functors, and hence left adjoint to the restriction functor $f^*\colon\ul\PSh_T(\Dd)\to\ul\PSh_T(\Cc)$ by \cite{martiniwolf2021limits}*{Corollary~3.3.3}. We will now describe a sufficient criterion for when $f^*$ acquires a strong symmetric monoidal structure again, making $f_!^\otimes$ into a left adjoint in $\CMon(\CAT_T)$:

\begin{prop}\label{prop:Day-restriction}
    Let $\Jj$ be a small symmetric monoidal $T$-$\infty$-category, and let $\Ii\subset \Jj$ be a symmetric monoidal subcategory. Assume that for every $A\in T$ the inclusion $\Ii(A)\hookrightarrow\Jj(A)$ admits a right adjoint $R_A$ such that the Beck--Chevalley map $R_A(X)\otimes R_A(Y)\to R_A(X\otimes Y)$ is an equivalence for all $X,Y\in\Jj(A)$.

    Then the unique $T$-cocontinuous symmetric monoidal functor $i_!^\otimes\colon\ul\PSh_T(\Ii)\to\ul\PSh_T(\Jj)$ extending the inclusion has a strong symmetric monoidal right adjoint, i.e.\ a right adjoint in $(T\times\textup{Fin}_*^\op)$-$\infty$-categories.
\end{prop}

Note that if one were willing to assume that the right adjoints $R_A$ from the previous proposition assemble into a parametrized functor again, this would be a direct consequence of 2-functoriality. However, in the example of interest to us, this stronger condition will not be satisfied, and so we will need a different approach. This will rely on the notions of initial and final $T$-functors, which we now recall:

\begin{definition}
    A functor $f\colon\Ii\to\Jj$ of small $T$-$\infty$-categories is called \emph{final} if for every $T$-cocomplete $T$-$\infty$-category $\Cc$, the Beck--Chevalley map $\colim_{\Ii}\circ f^*\to\colim_\Jj$ is an equivalence of functors $\ul\Fun(\Jj,\Cc)\to \Cc$. The notion of an \emph{initial} functor is defined dually.
\end{definition}

\begin{rk}
    By \cite{martiniwolf2021limits}*{Proposition~4.6.1 and Remark~4.6.2} it suffices to verify the above assumption for $\Cc=\ul\Spc_T$ in the case of final functors or $\Cc=\ul\Spc_T^\op$ for initial functors. Moreover, \emph{loc.\ cit.}\ also shows that the above notion of initial and final functors is equivalent to the definition given in \cite{martini2021yoneda}*{Definition~4.3.1} in terms of lifting properties.
\end{rk}

\begin{lemma}[cf.~\cite{shah2021parametrized}*{Theorem~6.7}]
    A functor $f\colon\Ii\to\Jj$ is final if and only if $f(A)\colon\Ii(A)\to\Jj(A)$ is final in the non-parametrized sense for every $A\in T$. Dually, $f$ is initial if and only if each $f(A)$ is initial.
    \begin{proof}
        It suffices to prove the second claim. By the parametrized version of Quillen's Theorem~A from \cite{martini2021yoneda}*{Corollary~4.4.8}, $f$ is initial if and only if for every $T$-$\infty$-groupoid $X\in\PSh(T)$ and each functor $X\to\Jj$, the induced map
        $\Ar(\Ii)\times_{\Jj}\ul X\to\ul X$ (where the pullback is via the target map) becomes an equivalence after applying the groupoidification functor $\Cat_T\to\Spc_T$ left adjoint to the inclusion. As the latter is just given by applying groupoidification pointwise, we see that $f$ is initial if and only if the projection \[\Ar(\Ii(A))\times_{\Jj(A)}\hom(A,X)\to\hom(A,X)\] becomes an equivalence after groupoidification for every $A\in T$ and $\ul X\to\Jj$. Applying the statement of the lemma for $T=1$ (and the $\infty$-groupoid $\hom(A,X)$) then shows that this holds provided that $f(A)$ is initial for every $A\in T$.

        Conversely, assume that $f$ is initial. If $T$ has a terminal object $1$, then running the above argument backwards for $A=1$ shows that the $\infty$-category $\Ar(\Ii(1))\times_{\Jj(1)} X$ is weakly contractible for every $\infty$-groupoid $X$ and every $X\to\Jj(1)$, i.e.\ $f(1)$ is initial. The claim now follows since for any $A\in T$ the $T_{/A}$-functor obtained from $f$ by restriction along $T_{/A}\to T$ is initial by \cite{martini2021yoneda}*{Remark~4.4.9}.
    \end{proof}
\end{lemma}

\begin{proof}[Proof of Proposition~\ref{prop:Day-restriction}]
    By \cite{martiniwolf2021limits}*{Corollary~3.3.3}, the underlying $T$-functor of $i_!^\otimes$ is left adjoint to the restriction functor $i^*\colon\ul\PSh_T(\Jj)\to\ul\PSh_T(\Ii)$ along $i\colon\Ii\hookrightarrow\Jj$. Thus, if we view the strong symmetric monoidal $T$-functor $i_!^\otimes$ as a $(T\times\textup{Fin}_*^\op)$-functor, then it has a levelwise right adjoint, and this satisfies the Beck--Chevalley condition with respect to maps in $T$. Factoring a general map in $T\times\textup{Fin}_*^\op$ and using that Beck--Chevalley maps compose, it will therefore suffice to show that the right adjoints also satisfy the Beck--Chevalley maps with respect to $\textup{Fin}_*^\op$. This precisely amounts to saying that for any fixed $A\in T$ the natural \emph{lax} symmetric monoidal structure on the functor $i^*\colon\ul\PSh_{T}(\Jj)(A)\to\ul\PSh_{T}(\Ii)(A)$ is strong, which is in turn equivalent to the Beck--Chevalley maps for the tensor product and the inclusion of the symmetric monoidal unit being invertible.

    As the $T$-functor $i_!$ is fully faithful by Theorem~\ref{thm:Kan-extension}, the Beck--Chevalley condition for the unit is automatic, and the Beck--Chevalley condition for the tensor product is equivalent to the existence of \emph{some} equivalence filling
    \[
        \begin{tikzcd}
            \ul\PSh(\Ii)(A)\times\ul\PSh(\Ii)(A)\arrow[r,"{-}\otimes{-}"] &[1em] \ul\PSh(\Ii)(A)\\
            \ul\PSh(\Jj)(A)\times\ul\PSh(\Jj)(A)\arrow[r,"{-}\otimes{-}"']\arrow[u,"i^*\times i^*"] & \ul\PSh(\Jj)(A)\arrow[u,"i^*"']
        \end{tikzcd}
    \]
    for each $A\in T$. Plugging in the description of the symmetric monoidal product from Proposition~\ref{prop:Day-convolution}, it will suffice to show that $i^*$ commutes with both the external product $\exttimes$ and with the left Kan extension $\otimes_!$ separately. The former is immediate from the definitions (without any assumptions on $i$). For the latter, we will show that the Beck--Chevalley map associated to the naturality equivalence
    \[
        \begin{tikzcd}
            \ul\PSh(\Jj)(A)\arrow[d,"i^*"']\arrow[r, "\otimes^*"] &\ul\PSh(\Jj\times\Jj)(A)\arrow[d,"(i\times i)^*"]\\
            \ul\PSh(\Ii)(A)\arrow[r,"\otimes^*"'] & \ul\PSh(\Ii\times\Ii)(A)
        \end{tikzcd}
    \]
    is an equivalence. Replacing $T$ by $T_{/A}$, we may assume that $A$ is terminal. For any $B\in T$ and $Y\in\Jj(B)$, we write $Y\downarrow{\otimes_\Jj}$ for the comma $T_{/B}$-category, i.e.~the pullback
    $\{Y\}\times_{\fgt^*\Jj}\Ar(\fgt^*\Jj)\times_{\fgt^*\Jj}\fgt^*(\Jj\times\Jj)$ where $\fgt\colon T_{/B}\to T$ denotes the forgetful functor and the pullback is via the tensor product; analogously we define $X\downarrow{\otimes_\Ii}$ for $X\in\Ii(B)$. By the pointwise formula for parametrized Kan extensions \cite{martiniwolf2021limits}*{Remark~6.3.6}, see also \cite{LLP}*{Remark~3.36}, the Beck--Chevalley condition then translates to saying that for every $B\in T$, $X\in\Ii(B)$, and $F\colon\Jj^\op\times\Jj^\op\to\ul\Spc_T$ the Beck--Chevalley map
    \[
        \colim_{(X\downarrow\otimes_\Ii)^\op} (\fgt^*F)\circ\pr
        \to
        \colim_{(X\downarrow\otimes_\Jj)^\op} (\fgt^*F)\circ\pr
    \]
    (i.e. the map induced by the inclusion of indexing categories) is invertible. For this it will in turn suffice to show that $X\downarrow\otimes_\Ii\hookrightarrow X\downarrow\otimes_\Jj$ is initial (so that the result of applying $(-)^\op$ is \emph{final}). By the previous lemma, this can be checked levelwise, i.e.\ we are altogether reduced to showing that for every $C\in T$, $Y\in\Ii(C)$ the inclusion
    \[
        \Ii(C)_{Y/}\times_{\Ii(C)}\big(\Ii(C)\times\Ii(C)\big)\hookrightarrow       
        \Jj(C)_{Y/}\times_{\Jj(C)}\big(\Jj(C)\times\Jj(C)\big)
    \]
    is an initial functor of (non-parametrized) $\infty$-categories. However, our assumptions guarantee that this functor admits a right adjoint, induced by $R_C$, and every left adjoint functor is in particular initial.
\end{proof}

\subsubsection{Slice symmetric monoidal structures} We will now give an even more concrete description of the symmetric monoidal cocompletion of a symmetric monoidal $T$-space. For this let us first recall from Remark~\ref{rk:Yoneda-image} that for any $T$-space $X$ the Yoneda embedding $X\to\ul\PSh_T(X)$ corresponds under our usual equivalence $\ul\PSh_T(X)\simeq\ul\Spc_T(X\times{-})$ to the functor 
\begin{equation}\label{eq:yoneda-concrete}
    y\colon\ul X\to\ul\Spc_T(X\times {-})
\end{equation} 
classifying the object $\Delta\colon X\to X\times X$ of $\ul\Spc_T(X\times X)=\PSh(T)_{/X\times X}$; in particular, $(\ref{eq:yoneda-concrete})$ exhibits $\ul\Spc_T(X\times{-})$ as the parametrized cocompletion of $X$ by Theorem~\ref{thm:univ-prop-PSh}. Let us recast this in terms of slice $T$-$\infty$-categories:

\begin{rk}\label{rk:slice-vs-shift}
    As before (Construction~\ref{constr:slice-cat}), we define $\smash{\ul\Spc_{T/X}\coloneqq(\ul\Spc_{T})_{/X}}$ as the pullback of $\ev_1\colon\Ar(\ul\Spc_T)\to\ul\Spc_T$ along the map $1\to\ul\Spc_T$ classifying the global section $X$. 
    
    Recall that the cartesian unstraightening of $\ul\Spc_T\in \Cat_T\subset \Fun(\PSh(T)^{\op},\Cat)$ over
    $\Spc_T=\PSh(T)$ is simply given by $\ev_1\colon\Ar(\Spc_T)\to\Spc_T$. As explained in the proof of Theorem~\ref{thm:criterion-spc-T}, it follows that the cartesian unstraightening of $\Ar(\ul\Spc_T)$ is given by
    the full subcategory
    $\Ar^\text{fw}(\Ar(\Spc_T))\subset\Ar(\Ar(\Spc_T))$ spanned by all
    maps in $\Ar(\Spc_T)$ inverted by $\ev_1$. Passing to cartesian unstraightenings over $\Spc_T=\PSh(T)$, we therefore obtain a pullback square as depicted on the left in the following diagram
    \[
        \begin{tikzcd}
            \Un^\ct(\ul\Spc_{T/X})\arrow[d]\arrow[r]\arrow[dr,pullback] & \Ar^\text{fw}(\Ar(\Spc_T))\arrow[d,"\ev_1"']\arrow[r,"\Ar^\text{fw}(\ev_0)"] &[1em] \Ar(\Spc_T)\arrow[d,"\ev_1"]\\
            \Spc_T\arrow[r,"f"'] & \Ar(\Spc_T)\arrow[r,"\ev_0"'] & \Spc_T\rlap,        \end{tikzcd}
    \] 
    where $f$ is the unique map of cartesian fibrations sending $1$ to $X$, i.e.~it is the functor $Y\mapsto(\pr\colon X\times Y\to Y)$. By direct inspection, also the right-hand square is a pullback, hence so is the total rectangle. As $\ev_0\circ f$ is simply the functor $X\times{-}$, naturality of straightening therefore shows that $\Un^\ct(\ul\Spc_{T/X})\simeq\Un^\ct(\ul\Spc_T(X\times{-}))$, and hence $\ul\Spc_{T/X}\simeq\ul\Spc_T(X\times{-})$. Plugging in the definitions, this equivalence is given levelwise by the evident equivalence $(\PSh(T)_{/Y})_{/\pr\colon X\times Y\to Y}\simeq\PSh(T)_{/X\times Y}$. We will denote the composite of the Yoneda embedding with this equivalence again by $y\colon\ul X\to\ul\Spc_{T/X}$ and refer to it as the \emph{Yoneda embedding} again; it in particular exhibits $\ul\Spc_{T/X}$ as the $T$-cocompletion of $X$. Chasing the identification of the Yoneda image from Remark~\ref{rk:Yoneda-image} through the above equivalence, we see that the image of $y\colon\ul X\to\ul\Spc_{T/X}$ agrees with the preimage of the terminal object under the forgetful functor $\ul\Spc_{T/X}\to\ul\Spc_T$.
\end{rk}

Suppose that $\mathcal C^{\otimes}\in \CMon(\CAT_T)\simeq \Fun(T^{\op},\CMon(\CAT_\infty))$ is a symmetric monoidal $T$-$\infty$-category and $X\in\CAlg(\mathcal C^{\otimes}(1))$ is a commutative algebra in its global sections. We now describe a slice symmetric monoidal structure on the slice $T$-$\infty$-category $\mathcal C_{/X}$ from Construction~\ref{constr:slice-cat}.
Let us start by recalling the pointwise symmetric monoidal structure on the arrow category:  

\begin{constr}
    Postcomposing with the limit-preserving functor  $\Ar=\Fun([1],-)\colon \CAT_\infty\to \CAT_\infty$ defines a functor 
    \begin{equation}\label{def:ar-otimes}
        \Ar\colon \CMon(\CAT_\infty)\to\CMon(\CAT_\infty).
    \end{equation}
    Denote by ${\Ar}\colon \CMon(\CAT_T)\to \CMon(\CAT_T)$ postcomposition with (\ref{def:ar-otimes}) and by $\ev_i\colon {\Ar}\to \id_{\CMon(\CAT_{T})}$, $i=0,1$ the functor induced by evaluation at $i\in[1]$. 
\end{constr}

\begin{rk}
    Note that if $\Cc^{\times}$ is a cartesian symmetric monoidal $T$-$\infty$-category in the sense of Remark~\ref{rk:what-it-means-to-be-cartesian}, then the symmetric monoidal $T$-$\infty$-category $\Ar(\mathcal C^{\times})$ is cartesian as well. 
\end{rk}

\begin{constr}
    Let $\mathcal C^{\otimes}\in\CMon(\CAT_T)\simeq\Fun(T^{\op},\CMon(\CAT_\infty))$ be a symmetric monoidal $T$-$\infty$-category. 
    Given any commutative algebra $X^\otimes$ in $\mathcal C^{\otimes}(1)$, we can view it as a lax symmetric monoidal functor $1\to\mathcal C^{\otimes}(1)$ and hence a map $1=\const\,1\to\mathcal C^{\otimes}$ in $\Fun(T^\op,\text{OP}_\infty)$, where $\text{OP}_\infty$ denotes the very large $\infty$-category of $\infty$-operads (containing $\CMon(\CAT_\infty)$ as a non-full subcategory).
    We then define $\mathcal C^\otimes_{/X}$ as the pullback 
    \begin{center}
    \begin{tikzcd}
    \mathcal C^{\otimes}_{/X}\arrow[dr,pullback]\arrow[r]\arrow[d] & \Ar(\mathcal C^{\otimes})\arrow[d,"\ev_1"] \\ 
    1 \arrow[r,"X^\otimes"']& \mathcal C^{\otimes }
    \end{tikzcd}
    \end{center} 
    in $\Fun(T^\op,\text{OP}_\infty)$. Since the forgetful functor $\text{OP}_\infty\to\CAT_\infty$ preserves limits, this refines the $T$-$\infty$-category $\Cc_{/X}$ from Construction~\ref{constr:slice-cat}. Moreover, in each degree $A\in T$, the above is just the usual slice symmetric monoidal structure \cite{HA}*{Theorem~2.2.2.4 and Remark~2.2.2.5}; in particular, each $\smash{\Cc^\otimes_{/X}(A)}$ is a symmetric monoidal $\infty$-category (and not just an operad), and the composite 
    \begin{equation}\label{eq:refine-fgt}
        \Cc^\otimes_{/X}\xrightarrow{\;\;\;}\Ar(\Cc^\otimes)\xrightarrow{\;\ev_0\;}\Cc^\otimes
    \end{equation} 
    is given in each degree by a strong symmetric monoidal functor (and not just a lax one). As this composite is moreover conservative, it follows that for every $f\colon A\to B$ in $T$ the lax symmetric monoidal restriction functor $\smash{f^*\colon\Cc_{/X}^\otimes(B)\to\Cc_{/X}^\otimes(A)}$ is strong symmetric monoidal, i.e.\ $\smash{\Cc_{/X}^\otimes}$ is a symmetric monoidal $T$-$\infty$-category. We will refer to this as the \emph{slice symmetric monoidal structure} on $\Cc_{/X}$. With this established, what we said above shows that $(\ref{eq:refine-fgt})$ is a refinement of the forgetful functor $\fgt\colon\Cc_{/X}\to\Cc$ to a strong symmetric monoidal $T$-functor.
\end{constr}

\begin{lemma}\label{lm:tensor-bil}
    Let $S\subset T$ be a cleft and let $\Cc^\otimes$ be an $S$-cocomplete symmetric monoidal $T$-$\infty$-category such that the tensor product is $S$-bilinear. Then also $\smash{\Cc_{/X}^\otimes}$ is $S$-cocomplete, and the tensor product is $S$-bilinear.
    \begin{proof}
        By \cite{martiniwolf2022presentable}*{Corollary~A.1.4}, the underlying $T$-$\infty$-category of $\smash{\Cc_{/X}^\otimes}$ is $S$-cocomplete, and the forgetful functor is $S$-cocontinuous. Using that the latter is also conservative and symmetric monoidal, we immediately see that the tensor product on $\smash{\Cc_{/X}^\otimes}$ preserves $S$-colimits in each variable as this holds in $\Cc^\otimes$.
    \end{proof}
\end{lemma}

\begin{prop}\label{prop:slice-sym-mon-is-universal}
    Let $X^\otimes\in\CMon(\Spc_T)$. Then the Yoneda embedding $\smash{y\colon \ul X\to\ul\Spc_{T/X}}$ refines to a strong symmetric monoidal functor $\smash{\ul X^\otimes\to\ul\Spc_{T/X}^\otimes}$ for the slice symmetric monoidal structure, thus exhibiting $\smash{\ul\Spc_{T/X}^\otimes}$ as the symmetric monoidal $T$-cocompletion of $X^\otimes$.
    \begin{proof}
        Observe first that any lax symmetric monoidal refinement of $y$ is already strong: namely, we can check this after postcomposing with the conservative strong symmetric monoidal functor $\fgt\colon\ul\Spc_{T/X}^\otimes\to\ul\Spc_T^\times$, and as $\fgt\circ y\colon X\to\ul\Spc_T$ takes values in the full subcategory spanned by the terminal object, it admits a unique lax symmetric monoidal structure, which is automatically strong. 

        Since the tensor product on $\ul\Spc_{T/X}$ is $T$-bilinear (Lemma~\ref{lm:tensor-bil}), any symmetric monoidal enhancement $y^{\otimes}$ of the Yoneda embedding extends uniquely to a symmetric monoidal, $T$-cocontinuous functor $\ul{\PSh}_T(X)^{\otimes}\to \ul\Spc_{T/X}^{\otimes}$ from the free symmetric monoidal $T$-cocompletion of $X$. As explained in Remark~\ref{rk:slice-vs-shift}, the underlying functor $\ul\PSh_T(X)\to {\ul{\Spc}_T}_{/X}$ is an equivalence, which then implies that $y^{\otimes}\colon \ul X\to \ul\Spc_{T/X}^{\otimes}$ is a free symmetric monoidal $T$-cocompletion. 

        We are therefore reduced to constructing a lax symmetric monoidal structure on the Yoneda embedding $y$, which by the definition of $\ul\Spc_{T/X}^{\otimes}$ as a pullback amounts to constructing a commutative square in $\Fun(T^{\op},\text{OP}_\infty)$
        \[
            \begin{tikzcd}
                \hom(-,X)^\otimes\arrow[r]\arrow[d] & \Ar(\ul\Spc_T)^\times\arrow[d,"\ev_1"]\\
                1^\otimes\arrow[r,"X"'] & \ul\Spc_T^\times
            \end{tikzcd}
        \]
        such that $\id_X\in\hom(X,X)$ gets sent to its Yoneda image, viewed as an element in the fiber over $(\pr_2\colon X\times X\to X)\in\ul\Spc_T(X)$.
        By the universal property of cartesian symmetric monoidal structures \cite{HA}*{Proposition 2.4.1.7}, this is equivalent to constructing a commutative square of $T$-$\infty$-categories
        \begin{equation}\label{diag:lax-yoneda-before-un}
            \begin{tikzcd}
                \Un^\cc\circ \hom(-,X)^\otimes\arrow[r]\arrow[d] & \Ar(\ul\Spc_T)\arrow[d,"\ev_1"]\\
                \consts\,\text{Fin}_*\arrow[r,"X"'] & \ul\Spc_T
            \end{tikzcd}
        \end{equation}
        sending $\id_X$ in the fiber over $1\in\text{Fin}_*$ to the Yoneda image of $\id_X$ and such that $\Un^\cc(\hom(Y,X)^\otimes)\to\Ar(\ul\Spc_T(Y))$ satisfies the Segal condition for every $Y\in\PSh(T)$. Note moreover that if the top arrow factors through $\ul\Spc_{T,*}\subset\Ar(\ul\Spc_T)$, then the Segal condition is automatic as the forgetful functor $\ul\Spc_{T,*}\to\ul\Spc_T$ is conservative and preserves fiberwise limits.
        
        Let us abbreviate $\Bb\coloneqq\PSh(T)$. By \cite{HHLNb}*{Example~2.5.10}, the composite of the Yoneda embedding $\Bb^\op\to\Fun(\Bb,\Spc)$ with straightening for left fibrations is given by $X\mapsto\Bb_{/X}$, i.e.~as the \emph{cartesian} straightening of $\ev_1\colon\Ar(\Bb)\to\Bb$. Thus, a diagram of the form $(\ref{diag:lax-yoneda-before-un})$ is equivalent (via cartesian unstraightening) to
        \begin{equation}\label{diag:lax-yoneda-after-un}
            \begin{tikzcd}
                \Ar(\Bb)\times_\Bb\text{Fin}_*\arrow[d]\arrow[r,"\alpha"] &[6em] \Ar^\text{fw}(\Ar(\Bb))\arrow[d,"\ev_1"]\\
                \Bb\times\text{Fin}_*\arrow[r,"{(Y,\langle n\rangle)\mapsto (\pr\colon Y\times X(\langle n\rangle)\to Y)}\,"'] & \Ar(\Bb)\rlap,
            \end{tikzcd}
        \end{equation}
        where the pullback is via the target map $\Ar(\Bb)\to\Bb$, and the left-hand vertical map is induced by the source map.        
        Plugging in the definitions, the map $\Un^\cc\circ\hom(-,X)^\otimes\to\Ar(\ul\Spc_T)$ corresponding to $\alpha$ factors through $\ul\Spc_{T,*}$ if and only if $\alpha$ takes values in the full subcategory of those squares in $\Bb$ such that the diagonal composite is an equivalence. Moreover, the resulting lax symmetric monoidal functor will be a refinement of the Yoneda embedding if and only if $\alpha$ further sends the object $(\id_X,1)$ in the fiber over $(X,1)$ to the object
        \[
            \begin{tikzcd}
                X\arrow[d,equals]\arrow[r,"\Delta"] & X\times X\arrow[d,"\pr"]\\
                X\arrow[r,equals] & X
            \end{tikzcd}
        \]
        in the fiber over $X\times X\to X$. We now simply note that the pasting
        \[
            \begin{tikzcd}
                \arrow[d]\Ar(\Bb)\times_\Bb\text{Fin}_*\arrow[r,"\pr"] & \Ar(\Bb)\arrow[d,"{(\ev_0,\ev_1)}"{description}]\arrow[r,"\beta"] & \Ar(\Bb)\arrow[d,"\ev_1"]\\
                \Bb\times\text{Fin}_*\arrow[r,"\id\times X^\otimes"']\arrow[r] &\Bb\times\Bb\arrow[r,"{-}\times{-}"'] & \Bb\rlap,
            \end{tikzcd}
        \]
        where $\beta$ sends $f\colon Y\to Z$ to 
        \[
            \begin{tikzcd}
                Y\arrow[d,equals]\arrow[r,"{(\id,f)}"] & Y\times Z\arrow[d,"\pr"]\\
                Y\arrow[r,equals] & Y
            \end{tikzcd}
        \]
        provides the desired square $(\ref{diag:lax-yoneda-after-un})$.
    \end{proof}
\end{prop}

Finally, let us generalize the above discussion to the setting of a cleft $S\subset T$:

\begin{prop}\label{prop:S-T-slice-univ-prop}
    Let $S\subset T$ be a cleft, and let $X\in\Spc_T$. The Yoneda embedding $y\colon\ul X\to\ul\Spc_{T/X}$ factors uniquely through $\ul\Spc_{S\triangleright T/X}\coloneqq\smash{\ul\Spc_{T/X}\times_{\ul\Spc_T}\ul\Spc_{S\triangleright T}}$, exhibiting the latter as $S$-cocompletion. Moreover, given any commutative monoid structure $X^\otimes$ on $X$, $\ul\Spc_{S\triangleright T/X}$ is a symmetric monoidal subcategory with respect to the slice symmetric monoidal structure on $\smash{\ul\Spc_{T/X}}$, and $\smash{y^\otimes\colon\ul X^\otimes\to\ul\Spc_{S\triangleright T/X}^\otimes}$ is a symmetric monoidal $S$-cocompletion.
    \begin{proof}
        Proposition~\ref{prop:spc-T-cart-pres-sym-mon} shows that $\smash{\ul\Spc_{S\triangleright T/X}^\otimes}$ is a symmetric monoidal subcategory for any commutative monoid $X^\otimes$. In light of the previous proposition, it will therefore once more suffice to prove the non-monoidal universal property. By \cite{martiniwolf2021limits}*{Theorem 7.1.13} this reduces to showing that $\smash{\ul\Spc_{S\triangleright T/X}}\subset\ul\Spc_{T/X}$ is the smallest subcategory containing the Yoneda image and closed under $S$-colimits.

        By design, $\ul\Spc_{S\triangleright T}\subset\ul\Spc_T$ contains the terminal object, and so the description of the Yoneda image from Remark~\ref{rk:slice-vs-shift} shows that the Yoneda embedding factors through $\ul\Spc_{S\triangleright T/X}$. Moreover, as $\ul\Spc_{S\triangleright T}\subset\ul\Spc_T$ is closed under $S$-colimits, so is $\ul\Spc_{S\triangleright T/X}\subset\ul\Spc_T$. Conversely, let $A\in T$ be arbitrary, let $a\colon A\to 1$ be the unique map in $\PSh(T)$, and let $Y\to a^*X$ define an element of $\ul\Spc_{S\triangleright T/X}(A)\subset(\PSh(T)_{/A})_{/X\times A}$; we want to show that $Y$ is an $S$-colimit of objects in the Yoneda image, which translates to saying that every $(f,g)\colon Z\to X\times A$ with $Z\in\Spc_S\subset\Spc_T$ is an $S$-colimit in $\ul\Spc_{T}(X\times{-})$ of objects in the Yoneda image. But it follows straight from the definitions that $(f,g)$ is the image under $f_!$ of $(\id,g)\colon X\to X\times A$, and the latter belongs to the Yoneda image by Remark~\ref{rk:Yoneda-image}.
    \end{proof}
\end{prop}

\begin{ex}\label{ex:post-composition-slice}
Suppose $f\colon X\to Y$ is a map of $T$-spaces and denote by 
    \[f_{!}\colon \ul{\Spc}_{S\triangleright T/X}\to \ul{\Spc}_{S\triangleright T/Y}\] the unique $S$-cocontinuous extension of $\smash{X\xrightarrow{\;f\;}Y\xhookrightarrow{\;y\;} {\ul\Spc}_{T/Y}}$.
    By construction of the Yoneda embedding (Remark~\ref{rk:slice-vs-shift}), $f_{!}$ factors as \[{\ul\Spc_{S\triangleright T}}_{/X}\simeq \ul\Spc_{S\triangleright T}(X\times -)\xrightarrow{(f\times -)_!}\ul\Spc_{S\triangleright T}(Y\times -)\simeq {\ul\Spc_{S\triangleright T}}_{/Y},\] where the equivalences are the ones from Remark~\ref{rk:slice-vs-shift}. 
    In particular, $f_!$ is given in degree $A\in T$ by the functor  $\PSh(T_{/A})_{/\pi_A^*f}\colon\PSh(T_{/A})_{/\pi_A^*X}\to \PSh(T_{/A})_{/\pi_A^*Y}$ postcomposing with $\pi_A^*f$. 
\end{ex} 

Suppose now that $S\subset T$ is a cleft, $X\in\Spc_T$, and $i\colon X\to \mathcal C$ is a functor into an $S$-cocomplete $T$-$\infty$-category $\Cc$.
Denote by \[{T}_i\colon \ul\Spc_{S\triangleright T/X}\to \Cc\]
the unique $S$-cocontinuous extension of $i$.
This has the following pointwise description on global sections:

\begin{lemma}\label{lm:pointwise-par-colimit}
For any object $(f\colon Y\to X)\in \Spc_{S\triangleright T/X}$, we have an equivalence
\[ T_i(f)\simeq\colim_{Y}(i\circ f)\]
in $\Cc(1)$, where $\colim_{Y}\colon \ul\Fun(\ul Y,\mathcal C)\to {\mathcal C}$ denotes the parametrized colimit-functor (left adjoint to the inclusion of constant diagrams).
\begin{proof}
    As $Y\in\Spc_{S\triangleright T}$ and $T_i$ is an $S$-cocontinuous extension of $i$, it suffices to show that $f\in\Spc_{S\triangleright T/X}$ is the parametrized colimit of the composite $y\circ f\colon\ul Y\to\ul X\to\ul\Spc_{T/X}$.

    If we write $f_!\colon\ul\Spc_{S\triangleright T/Y}\to\ul\Spc_{S\triangleright T/X}$ for the unique $S$-cocontinuous extension of $yf$ as in the previous example, then we have seen in said example that $f_!(\id_Y)\simeq f$. On the other hand, $yf=f_!y$ by construction. As $f_!$ is $S$-cocontinuous, this reduces us to considering the case $f=\id_Y$, i.e.~we want to show that the colimit of $y\colon \ul Y\to\ul\Spc_{S\triangleright T/Y}$ is the terminal object $\id_Y$. By \cite{martiniwolf2021limits}*{Proposition 6.1.3}, this holds in $\ul\Spc_{T/Y}$. As $\ul\Spc_{S\triangleright T/Y}\subset\ul\Spc_{T/Y}$ is a full subcategory containing the terminal object, this then also holds in $\ul\Spc_{S\triangleright T/Y}$.
\end{proof}
\end{lemma}

\subsection{The global \texorpdfstring{$\bm J$}{J}-homomorphism}\label{subsec:J}
We begin by recalling the classical $J$-homorphism $\BO\times\mathbb Z\to \pic(\Sp)$: 

\begin{constr}\label{constr:definitionjhomomorphism}
    Consider the topologically enriched symmetric monoidal functor $\cat{j}^\otimes\colon\cat{L}^\oplus\to\cat{Top}_*^\smashp$ sending an inner product space to its 1-point compactification $S^V$ (with basepoint the point at $\infty$) and sending a linear isometric embedding $i\colon V\to W$ to its unique basepoint-preserving extension $S^V\to S^W$; the unit isomorphism is the unique one, while the compatibility isomorphism for the tensor product is given by the isomorphisms $S^{V\oplus W}\to S^V\smashp S^W, (v,w)\mapsto v\smashp w$.
    Applying $\Ntop$ and postcomposing with the localisation $\Ntop(\cat{Top}_*^\smashp)\to \Spc_{*}^{\otimes}$ yields a functor 
    $j^{\otimes}\colon \Ntop(\cat{L}^\oplus)\to \Spc_{*}^{\otimes}$. 
    The composite $\Sigma^{\infty}\circ j^{\otimes}$ restricts to a map of commutative monoids $\coprod_{n\ge0}B\O(n)\simeq\Ntop(\core\cat{L}^\oplus)\to\pic(\Sp)$, and using the universal property of group completion this extends uniquely to a map $J^\otimes\colon B\O\times\Z\to\pic(\Sp)$ of commutative monoids, which we call the \emph{$J$-homomorphism}. 
\end{constr} 

Using our concrete pointset models $\smash{\ul\myS_{\gl,*}^\otimes}$ of the initial pointed globally presentably symmetric monoidal $\infty$-category (Theorem~\ref{thm:gl-smash-initial}) and $\smash{\ul\mySp_\gl^\otimes}$ of the initial representation stable globally presentably symmetric monoidal $\infty$-category (Corollary~\ref{cor:global-spectra-smash}) in terms of pointed orthogonal spaces and orthogonal spectra, respectively, we can easily write down a global version of the previous construction:

\begin{constr}\label{constr:global-sph}
    By applying the symmetric monoidal continuous Borel construction (Construction~\ref{constr:sym-mon-Borel}) to the topological symmetric monoidal functor $\cat{j}^\otimes$ from above and localizing at the equivariant weak equivalences, we obtain a symmetric monoidal global functor $\smash{\ul{\mathfrak j}^\otimes\colon\Ntop(\cat{L}^{\oplus,\dual})\to\ul\myS^\otimes_*}$ sending a $G$-representation $V$ to the $G$-space $S^V$. We will write $\ul{\mathfrak j}_\gl^\otimes$ for the composite 
    \[
        \Ntop(\cat{L}^{\oplus,\dual})\xrightarrow{\;\ul{\mathfrak j}^\otimes\;}\ul\myS_*^\otimes\xhookrightarrow{\;\consto\;}\ul\myS_{\gl,*}^{\otimes}
    \]
    with the unique symmetric monoidal equivariantly cocontinuous functor (modelled on the pointset level by sending a topological space to the constant orthogonal space).
\end{constr}

\begin{rk}\label{rk:equivariant-j}
    By definition, $\ul{\mathfrak j}^\otimes(1)\colon\Ntop(\cat{L}^\otimes)\to\Spc_*^\otimes$ recovers the symmetric monoidal functor $\smash{j^\otimes}$ from Construction~\ref{constr:definitionjhomomorphism}.
\end{rk}

\begin{constr}\label{constr:global-J}
    We write $\ul\Vect^\oplus\coloneqq\core\Ntop(\cat{L}^{\oplus,\flat})$, and we fix once and for all a (pointwise) group completion $\ul\Rep^\oplus\to\ul\Rep^{\oplus,\text{grp}}\eqqcolon\ul\VRep^\oplus$. We may think of objects of $\ul\VRep(\BGcat{G})$ as \emph{virtual $G$-representations}, i.e.\ formal differences $V\ominus W$ of $G$-representations.

    For any $G$-representation $V$, the $G$-global spectrum $\Sigma^\infty S^V$ is invertible with respect to the smash product as a consequence of Proposition~\ref{prop:resp-stability-G-gl}. Thus, the composite
    \[
        \ul\Vect^\oplus\xrightarrow{\;\ul{\mathfrak j}_\gl^{\smash\otimes}\;}\ul\myS_{\gl,*}^{\otimes}\xrightarrow{\;\Sigma^\bullet\;}\ul\mySp_\gl^\otimes
    \]
    with the unique symmetric monoidal global left adjoint (see  Variant~\ref{var:reduced-global-susp-spectrum}) factors uniquely through a symmetric monoidal global functor $\smash{\ul{\mathfrak J}_\gl^\otimes\colon\ul\VRep^\oplus\to\ul\mySp_\gl^\otimes}$. 
    We call $\ul{\mathfrak J}_\gl^\otimes$ the \emph{global $J$-homomorphism}.
\end{constr}

The above construction of the global $J$-homomorphism in terms of the pointset level models has the advantage that it is very explicit, and this description will be particularly convenient for the comparison of global Thom spectrum constructions in the next section. Nevertheless, in more categorical contexts it can be preferable to have a model-independent construction. While there is probably no way to avoid using at some point that topological spaces model $\infty$-groupoids, we will now give an alternative construction of $\ul{\mathfrak j}_\gl^\otimes$ that only requires the symmetric monoidal functor $\smash{j^\otimes\colon\Ntop(\cat{L}^\oplus)\to\Spc_*^\otimes}$ as geometric input. To emphasize the model-independence, we phrase this in terms of the construction $\ul\Spc_{\Glo,*}^\otimes$ of the initial pointed globally presentably symmetric monoidal $\infty$-category obtained from the abstract theory developed in §\ref{subsec:pointed-sym-mon} (which admits a \emph{unique} equivalence $\Phi\colon\ul\myS_{\gl,*}^\otimes\iso\ul\Spc_{\Glo,*}^\otimes$).

\begin{thm}\label{thm:who-is-afraid-of-the-J-homomorphism}
    There exists an initial example of an extension of 
    \begin{equation}\label{eq:constj}
        \Ntop(\cat{L}^\oplus)\xrightarrow{\;j^\otimes\;}\Spc_*^\otimes\xrightarrow{\;\const\;}\Spc_{\Glo,*}^\otimes
    \end{equation}
    (where $\const$ denotes the unique symmetric monoidal left adjoint) to a symmetric monoidal global functor $\ul{j}_\gl^\otimes\colon\Ntop(\cat{L}^{\oplus,\flat})\to\ul\Spc_{\Glo,*}^\otimes$, and the diagram
    \begin{equation}\label{diag:comparison-of-js}
        \begin{tikzcd}[row sep=1ex]
            & \ul\myS_{\gl,*}^{\otimes}\arrow[dd,"\Phi"',"\sim"]\\
            \Ntop(\cat{L}^{\oplus,\dual})\arrow[ur, bend left=15pt,"\ul{\mathfrak j}_\gl^\otimes"]\arrow[dr, bend right=15pt," \ul{j}_\gl^\otimes"']\\
            & \ul\Spc_{\Glo,*}^\otimes
        \end{tikzcd}
    \end{equation}
    commutes. More precisely, any identification of the resulting symmetric monoidal functors on global sections $\smash{\Ntop(\cat{L}^\oplus)\rightrightarrows\Spc_{\Glo,*}^\otimes}$  extends uniquely to an equivalence filling $(\ref{diag:comparison-of-js})$.
\end{thm}

Let us begin by observing that $(\ref{eq:constj})$ corresponds by adjunction to a symmetric monoidal global functor $\smash{\kern1.33pt\hat{\kern-1.33pt\textit{\j}}^\otimes\colon\const\,\Ntop(\cat{L}^\oplus)\to\ul\Spc_{\Glo,*}^\otimes}$. If we ignore monoidal structures for a second, then extending $(\ref{eq:constj})$ to a global functor is equivalent to extending $\kern1.33pt\hat{\kern-1.33pt\textit{\j}}$ along the (fully faithful) adjunction counit $\smash{\const\,\Ntop(\cat{L})\to\Ntop(\cat{L}^\dual)}$. The initial such extension is given by the left Kan extension of $\kern1.33pt\hat{\kern-1.33pt\textit{\j}}$, so we have to understand when a functor $\smash{\Ntop(\cat{L}^\dual)\to\ul\Spc_{\Glo,*}}$ is left Kan extended from $\const\,\Ntop(\cat{L})$.

\begin{lemma}
    Let $T$ be a small $\infty$-category with a terminal object $1$, and let $\Cc$ be a small $T$-$\infty$-category such that for every $p\colon A\to1$ the functor $p^*\colon\Cc(1)\to\Cc(A)$ is fully faithful and admits a right adjoint $p_*$. Moreover, let $\Ee$ be any presentable (non-parametrized) $\infty$-category and denote by $\ul\Ee_T$ the category of $T$-objects in $\Ee$ (Example~\ref{ex:T-objects}). Then a $T$-functor $F\colon\Cc\to\ul\Ee_T$ is left Kan extended from $\const\,\Cc(1)$ if and only if for every $A\in T$ and $X\in\Cc(A)$ the map $F(\epsilon\colon p^*p_*X\to X)$ in $\ul\Ee_T(A)=\Fun(T_{/A}^\op,\Ee)$ induced by the counit becomes an equivalence after evaluation at $\id_A$.
    \begin{proof}
        Write $\pi\colon\const\,\Ee\to\ul\Ee_T$ for the unique $T$-functor extending $\const\colon\Ee\to\ul\Ee_T(1)=\Fun(T^\op,\Ee)$. By \cite{LLP}*{Theorem~6.9} (for $M=\core T$), the individual right adjoints $\ul\Ee_T(A)\to\Ee$ (evaluating at $\id_A\in T_{/A}^\op$) assemble into a right adjoint $\rho\colon\Un^\cc(\ul\Ee_T)\to\Un^\cc(\const\,\Ee)=\Ee\times T^\op$ over $T^\op$ such that the composite
        \[
            \Fun_T(\Dd,\ul\Ee_T)\xrightarrow{\;\Un^\cc\;}
            \Fun_{/T^\op}(\Un^\cc(\Dd),\Un^\cc(\ul\Ee_T))\xrightarrow{\;\rho\circ{-}\;}
            \Fun_{/T^\op}(\Un^\cc(\Dd),\Ee\times T^\op)
        \]        
        is an equivalence for every $T$-$\infty$-category $\Dd$. Applying this for $\Dd=\Cc$ and $\Dd=\const\,\Cc(1)$ yields a commutative square
        \begin{equation}\label{diag:laxify}
            \begin{tikzcd}[cramped]
                \Fun_T(\Cc,\ul\Ee_T)\arrow[r,"\sim"]\arrow[d,"\res"'] & \Fun_{/T^\op}(\Un^\cc(\Cc),\Ee\times T^\op)\arrow[d,"\res"]\\
                \Fun_T(\const\,\Cc(1),\ul\Ee_T)\arrow[r,"\sim"'] & \Fun_{/T^\op}(\Cc(1)\times T^\op,\Ee\times T^\op)\rlap.
            \end{tikzcd}
        \end{equation}
        By \cite{HA}*{Proposition 7.3.2.6${}^\op$}, the individual right adjoints $p_*\colon\Cc(A)\to\Cc(1)$ assemble into a right adjoint $\gamma\colon\Un^\cc(\Cc)\to\Cc(1)\times T^\op$ over $T^\op$, so that the right-hand vertical arrow in $(\ref{diag:laxify})$ admits a \emph{left} adjoint given by precomposing with $\gamma$. On the other hand, the left-hand vertical arrow admits a left adjoint given by parametrized left Kan extension, and it follows formally that also the square where we pass to left adjoints vertically commutes.

        As we assumed each $p^*\colon\Cc(1)\to\Cc(A)$ to be fully faithful, the right adjoint $\gamma$ is a localization at the images of the counits $\epsilon\colon p^*p_*X\to X$ for all $A\in T$ and $X\in\Cc(A)$. Thus, $F\colon\Cc\to\ul\Ee_T$ is in the image of the left adjoint of the left-hand vertical map (i.e., $F$ is left Kan extended) if and only if $\rho\circ\Un^\cc(F)$ inverts the above counit maps. By the explicit fiberwise description of $\rho$, this is precisely demanding that each $F(\epsilon)$ become an equivalence after evaluation at $\id_A\in T_{/A}^\op$.
    \end{proof}
\end{lemma}

\begin{cor}\label{cor:J-LKE}
    The composite
    $\smash{\Ntop(\cat{L}^\dual)\xrightarrow{\;\ul{\mathfrak j}_\gl\;}\ul\myS_{\gl,*}\xrightarrow[\raise3.5pt\hbox{$\scriptstyle\smash{\sim}$}]{\;\;\Phi\;\;}\ul\Spc_{\Glo,*}}$
    is left Kan extended from $\const\,\Ntop(\cat{L})$.
    \begin{proof}
        It is clear that $\triv_G\colon\Ntop(\cat{L}^\dual)(1)\to\Ntop(\cat{L}^\dual)(\BGcat{G})$ is fully faithful for every compact Lie group $G$, with right adjoint given by taking $G$-fixed points. Applying the previous lemma for $\Ee=\Spc_*$ and using that the composite
        \[
            \myS_{G,*}\xrightarrow{\;\const\;}\myS_\text{$G$-gl}\xrightarrow[\raise3.5pt\hbox{$\scriptstyle\smash{\sim}$}]{\;\;\Phi\;\;}\ul\Spc_{\Glo,*}(\BGcat{G})\xrightarrow{\;\ev_{\id_{\BGcat{G}}}}\Spc_*
        \]
        is given by taking $G$-fixed points (see Lemma~\ref{lemma:id-g-fixed-points-vs-ev}), we only have to show that for every $G$-representation $V$ the inclusion $\smash{S^{(V^G)}}\hookrightarrow S^V$ of $G$-spaces induces an equivalence on $G$-fixed points. This is obvious.
    \end{proof}
\end{cor}

In order to lift the above results to statements about \emph{symmetric monoidal} global functors, we will use the following proposition:

\begin{prop}\label{prop:symmetric-monoidal-LKE}
    Let $\Cc^\otimes$ be globally presentably symmetric monoidal. Then the restriction 
    $i^{*}_\otimes\colon\Fun^\otimes(\Ntop(\cat{L}^{\oplus,\dual}),\Cc^\otimes)\rightarrow\Fun^\otimes(\const\,\Ntop(\cat{L}^\oplus),\Cc^\otimes)$
    admits a left adjoint $\smash{i_!^\otimes}$, and the Beck--Chevalley map $\smash{i_!\circ\fgt\to\fgt\circ i_!^\otimes}$ associated to
    \[
        \begin{tikzcd}
            \Fun^\otimes(\Ntop(\cat{L}^{\oplus,\dual}),\Cc^\otimes)\arrow[r,"i^*_\otimes"]\arrow[d,"\fgt"'] & \Fun^\otimes(\const\,\Ntop(\cat{L}^\oplus),\Cc^\otimes)\arrow[d,"\fgt"]\\
            \Fun(\Ntop(\cat{L}^\dual),\Cc)\arrow[r,"i^*"'] & \Fun(\const\,\Ntop(\cat{L}),\Cc)
        \end{tikzcd}
    \]
    is an equivalence.
    \begin{proof}
        We may identify the diagram in question with
        \begin{equation}\label{diag:too-many-forgetful-functors}
            \begin{tikzcd}
                \Fun^{\text{L},\otimes}(\ul\PSh(\Ntop(\cat{L})^\flat)^\otimes,\Cc^\otimes)\arrow[r]\arrow[d,"\fgt"'] & \Fun^{\text{L},\otimes}(\ul\PSh(\const\,\Ntop(\cat{L}))^\otimes,\Cc^\otimes)\arrow[d,"\fgt"]\\
                \Fun^\text{L}(\ul\PSh(\Ntop(\cat{L})^\flat),\Cc)\arrow[r] & \Fun^\text{L}(\ul\PSh(\const\,\Ntop(\cat{L})),\Cc)\rlap,
            \end{tikzcd}
        \end{equation}
        where the top and bottom horizontal arrow are given by precomposing with the symmetric monoidal global left adjoint $L^\otimes\colon\ul\PSh(\const\,\Ntop(\cat{L}))^\otimes\to\ul\PSh(\Ntop(\cat{L}^\dual))^\otimes$ induced by the inclusion and its underlying global functor $L\colon\smash{\ul\PSh(\const\,\Ntop(\cat{L}))}\to\smash{\ul\PSh(\Ntop(\cat{L}^\dual))}$, respectively. The latter admits a right adjoint $R$ (namely, the restriction functor) which is again globally cocontinuous, and the left adjoint to the bottom arrow in the above diagram is then given by precomposing with $R$, with unit and counit again given by precomposition. If we can show that $L\dashv R$ upgrades to a symmetric monoidal adjunction $L^\otimes\dashv R^\otimes$, then a left adjoint to the top horizontal arrow in $(\ref{diag:too-many-forgetful-functors})$ will be similarly given by precomposing with $R^\otimes$ and with unit and counit given analogously, which then clearly satisfies the Beck--Chevalley condition.

        For this we will verify the assumptions of Proposition~\ref{prop:Day-restriction}. As observed in the proof of the previous corollary, each $\triv_G\colon\Ntop(\cat{L}^\dual)(1)\to\Ntop(\cat{L}^\dual)(\BGcat{G})$ is fully faithful and has a right adjoint induced by taking $G$-fixed points. Unravelling the definitions, the Beck--Chevalley condition amounts to saying that $V^G\oplus W^G\hookrightarrow V\oplus W$ defines an isomorphism onto $(V\oplus W)^G$, which is obvious.
    \end{proof}
\end{prop}

\begin{proof}[Proof of Theorem~\ref{thm:who-is-afraid-of-the-J-homomorphism}]
    As noted above, the symmetric monoidal functor $(\ref{eq:constj})$ corresponds to a symmetric monoidal global functor $\kern1.33pt\hat{\kern-1.33pt\textit{\j}}^\otimes\colon\const\,\Ntop(\cat{L}^\oplus)\to\ul\Spc_{\Glo,*}^\otimes$ and Proposition~\ref{prop:symmetric-monoidal-LKE} provides the desired initial extension $\smash{\ul j_\gl^\otimes}$. The universal property then shows that any identification $\smash{\ul j_\gl^\otimes(1)=j^\otimes\simeq(\Phi\circ\ul{\mathfrak j}_\gl)(1)}$ (e.g. the one obtained from Remark~\ref{rk:equivariant-j}) extends uniquely to a symmetric monoidal natural transformation $\ul j_\gl^\otimes\to \Phi\circ\ul{\mathfrak j}_\gl^\otimes$ and it only remains to show that this transformation is an equivalence. This can be checked after forgetting the monoidal structure, where this follows immediately from the fact that both maps sides are left Kan extended from the full subcategory $\const\,\Ntop(\cat{L})$ by Corollary~\ref{cor:J-LKE} and Proposition~\ref{prop:symmetric-monoidal-LKE}.
\end{proof}

\begin{rk}\label{rk:J-with-fewer-models}
    One can similarly define a `model independent' version of the global $J$-homomorphism as the unique extension of the composite of $\ul j_\gl^\otimes\colon\ul\Vect^\oplus\to\ul\Spc_{\Glo,*}^\otimes$ with the initial symmetric monoidal left adjoint $\smash{\ul\Spc_{\Glo,*}^\otimes\to\ul\Sp_{\Glo}^{\text{RepSph},\otimes}}$ inverting the image of $\smash{\ul j_\gl^\otimes}$, as provided by §\ref{subsec:invert-objects}. The previous theorem then implies that the result corresponds under the unique equivalence $\smash{\ul\mySp_\gl^\otimes\iso\ul\Sp_{\Glo}^{\text{RepSph},\otimes}}$ to the global $J$-homomorphism $\ul{\mathfrak J}_\gl^\otimes$ constructed above via the pointset model.
\end{rk}

\subsection{Parametrized Thom spectrum functors}\label{subsec:global-Thom} 
Applying Proposition~\ref{prop:slice-sym-mon-is-universal}, we can now  construct the global Thom spectrum functor and its unstable analogue:
\begin{definition}
    We define the \emph{global Thom space functor}
    \[
        \ul\th_\gl^\otimes\colon \ul\Spc_{\Glo/\ul\Vect^\oplus}^\otimes\to\ul\myS_{\gl,*}^\otimes
    \]
    as the unique globally cocontinuous extension of $\ul{\mathfrak j}_\gl^\otimes\colon\ul\Vect^\oplus\to\ul\myS_{\gl,*}^\otimes$ (Construction~\ref{constr:global-sph}) along the Yoneda embedding, and we define the \emph{global Thom spectrum functor} as the unique globally cocontinuous extension
    \[
        \ul\Th_\gl^\otimes\colon\ul\Spc_{\Glo/\ul\VRep^{\oplus}}^\otimes\to\ul\mySp_\gl^\otimes
    \]
    of the global $J$-homomorphism $\ul{\mathfrak J}_\gl^\otimes\colon\ul\VRep^{\oplus}\to\ul\mySp_\gl^\otimes$ from Construction~\ref{constr:global-J}.
\end{definition}
\begin{rk}\label{rk:incoherent-description-Thom-gl}
    In particular, we get for any compact Lie group $G$ symmetric monoidal left adjoints
    \begin{align*}
        \th_\text{$G$-gl}^\otimes\colon\Spc_{\text{$G$-gl}/\triv_G\ul\Vect^\oplus}^\otimes\coloneqq\PSh(\Glo_{/\BGcat{G}})_{/\triv_G\ul\Vect^\oplus}^\otimes&\to \myS_{\text{$G$-gl},*}^{\otimes}\\
        \Th_\text{$G$-gl}^\otimes\colon\Spc_{\text{$G$-gl}/\triv_G\ul\VRep^{\oplus}}^\otimes\coloneqq\PSh(\Glo_{/\BGcat{G}})_{/\triv_G\ul\VRep^{\oplus}}^\otimes&\to\mySp_\text{$G$-gl}^\otimes
    \end{align*}
    where the sources carry the slice symmetric monoidal structures and the targets carry the symmetric monoidal structure induced by the respective pointset level smash product. The fact that we started with global functors provides (coherent) symmetric monoidal equivalences ${\th_\text{$G$-gl}^\otimes\circ f^*\simeq f^*\circ\th_\text{$G'$-gl}^\otimes}$ and $\smash{\Th_\text{$G$-gl}^\otimes\circ f^*\simeq f^*\circ\Th_\text{$G'$-gl}^\otimes}$ for any homomorphism $f\colon G\to G'$ of compact Lie groups, while the fact that the global functors $\ul\th_\gl$ and $\ul\Th_\gl$ are \emph{parametrized} left adjoints implies that the corresponding Beck--Chevalley maps provide equivalences $f_!\circ\th_\text{$G$-gl}\simeq\th_\text{$G'$-gl}\circ f_!$ and $f_!\circ\Th_\text{$G$-gl}\simeq\Th_\text{$G'$-gl}\circ f_!$ (of non-monoidal functors).
\end{rk}

Corollary~\ref{cor:SpR-sym-mon-univ-prop} shows that the unique globally cocontinuous symmetric monoidal functor $\Sigma^{\bullet}\colon \ul\myS_{\gl,*}^{\otimes}\to \ul\mySp_{\gl}^{\otimes}$ (also see Variant~\ref{var:reduced-global-susp-spectrum} for a pointset level description) is the initial map in $\smash{\CAlg(\PrL_\Glo)}$ that inverts representation spheres. By construction of the global Thom space and Thom spectrum functors, this fits into a commutative square
\begin{equation}\label{diag:thom-space-spectrum}
    \begin{tikzcd}
    \ul{\Spc}_{\Glo/\ul\Vect^{\oplus}}\arrow[r,"\ul\th_\gl^\otimes"]\arrow[d] & \ul{\myS}_{\gl,*}^\otimes\arrow[d,"\Sigma^{\bullet}"]\\
    \ul{\Spc}_{\Glo/\ul\VRep^{\oplus}}\arrow[r,"\ul\Th_\gl^\otimes"'] & \ul{\mySp}_\gl^\otimes
    \end{tikzcd}
\end{equation}
of presentably symmetric monoidal global $\infty$-categories, where the left vertical map is the unique globally cocontinuous symmetric monoidal extension of the composition 
\[
    \ul\Vect^{\oplus}\to \ul\VRep^{\oplus}\xrightarrow{y^{\otimes}}{\ul\Spc_{\Glo}}_{/\ul\VRep^{\oplus}}
\] 
of the group completion map with the Yoneda embedding (Proposition~\ref{prop:slice-sym-mon-is-universal}).
We owe the following corollary of the symmetric monoidal universal property of global spectra (Theorem~\ref{thm:stable-main} and Corollary~\ref{cor:global-spectra-smash}) to David Gepner, Thomas Nikolaus, and Stefan Schwede: 
    \begin{cor}\label{cor:Thom-space-po}
    The diagram $(\ref{diag:thom-space-spectrum})$ is a pushout of presentably symmetric monoidal global $\infty$-categories. 
\begin{proof} 
    By Corollaries~\ref{cor:global-spectra-smash} and \ref{cor:gl-pointed} and the construction of the various maps, $(\ref{diag:thom-space-spectrum})$ is a commutative square in $\CAlg(\Pr^\Glo_\Glo)$.
    It therefore suffices to show that for every presentably symmetric monoidal global $\infty$-category $\Cc$, 
    the square $(\ref{diag:thom-space-spectrum})$ induces an equivalence \[\hom(\ul\mySp_{\gl}^{\otimes},\Cc)\iso  \hom(\ul\myS_{\gl,*}^{\otimes},\Cc)\times_{\hom({\ul\Spc_{\Glo/\ul\Vect^{\oplus}}},\Cc)}\hom(\ul\Spc_{\Glo/\ul\VRep^{\oplus}},\Cc),\] where $\hom({-},{-})\coloneqq \hom_{\CAlg(\Pr^{\Glo}_{\Glo})}({-},{-})$ denotes the $\infty$-groupoid of symmetric monoidal left adjoints.  
  
    Fix a presentably symmetric monoidal global $\infty$-category $\Cc$.
    It follows from Proposition~\ref{prop:slice-sym-mon-is-universal} and the universal property of group completion that restriction along $\ul{\Spc}_{\Glo/\ul\Vect^{\oplus}}\to \ul\Spc_{\Glo/\ul\VRep^{\oplus}}$ yields an inclusion of path components \[ \hom(\ul{\Spc}_{\Glo/\ul\VRep^{\oplus}},\Cc)\hookrightarrow \hom(\ul{\Spc}_{\Glo/\ul\Vect^{\oplus}},\Cc)\] with image those globally cocontinuous symmetric monoidal functors $\ul{\Spc}_{\Glo/\ul\Vect^{\oplus}}\to \Cc$ which send elements in the Yoneda image of $\smash{\ul\Vect^{\oplus}}$ to invertible objects.   
    By construction of the Thom space functor $\ul\th_\gl^\otimes$, this implies that the projection \begin{align*}\hom(\ul\myS_{\gl,*}^{\otimes},\Cc)\times_{\hom(\ul\Spc_{\Glo/\ul\Vect^{\oplus}},\Cc)}\hom(\ul\Spc_{\Glo/\ul\VRep^{\oplus}},\Cc) \to \hom(\ul\myS_{\gl,*}^{\otimes},\Cc)\end{align*} is the inclusion of the components consisting of those globally cocontinuous symmetric monoidal functors which in every degree $G$ send all representation spheres to invertible objects. The map $\hom(\Sigma^\bullet,\Cc)\colon\hom(\ul\mySp_\gl^\otimes,\Cc)\to\hom(\ul\myS_{\gl,*},\Cc)$ admits the same description by the universal property of global spectra (Corollary~\ref{cor:global-spectra-smash}) and Corollary~\ref{cor:SpR-sym-mon-univ-prop},  finishing the proof.
\end{proof}
\end{cor}

\begin{variant}
    People allergic to pointset models might prefer defining the global Thom spectrum functor as the unique left adjoint $\ul\Spc_{\Glo/\ul\VRep^{\oplus}}^\otimes\to\ul\Sp_{\Glo}^{\text{RepSph},\otimes}$ extending the `model-independent' global $J$-homomorphism from Remark~\ref{rk:J-with-fewer-models}. By said remark, the result will agree with the above construction up to postcomposition with the unique symmetric monoidal equivalence $\ul\Sp_{\Glo}^{\text{RepSph},\otimes}\simeq\ul\mySp_\gl^\otimes$. Similarly, we may of course change the model of the symmetric monoidal cocompletion of $\smash{\ul\VRep^{\oplus}}$, and the two resulting global Thom spectrum functors will agree up to precomposition with the unique equivalence under $\ul\VRep^{\oplus}$.
\end{variant}

We now explain how to adapt the global construction from the previous subsection to the equivariant setting. For simplicity, we will restrict to the case of equivariant Thom \emph{spectra}, leaving the analogous case of Thom \emph{spaces} to the interested reader. 

\begin{lemma}
    The global $J$-homomorphism $\ul{\mathfrak J}_\gl^\otimes\colon\ul\VRep^{\oplus}\to\ul\mySp_\gl^\otimes$ factors uniquely through a symmetric monoidal $\Glo$-functor $\smash{\ul{\mathfrak J}^\otimes\colon\ul\VRep^{\oplus}\to\ul\mySp^\otimes}$.
    \begin{proof}
        By Proposition~\ref{prop:mySp-vs-mySp-gl}, the unique equivariantly cocontinuous symmetric monoidal global functor $i\colon \ul\mySp^\otimes\to\ul\mySp_\gl^\otimes$ is fully faithful and at every compact Lie group $G$, it admits a right adjoint $r_G\colon\mySp_{\text{$G$-gl}}^{\otimes}\to \mySp_{G}^{\otimes}$ which is a localization at the equivariant weak equivalences.
        As equivariant weak equivalences are preserved by smashing with flat orthogonal $G$-spectra (Theorem~\ref{thm:equiv-flatness}), it follows that the induced lax symmetric monoidal structure on the right adjoint $r_G$ is strong. 
        By \cite{puetzstueck-new}*{Corollary 6.3}, this implies that $\smash{i_G\colon \Sp_G^{\otimes}\rightleftarrows \Sp_{\text{$G$-gl}}^{\otimes}\noloc r_G}$ induce mutually inverse equivalences 
        $\smash{\pic(\Sp_G^{\otimes})\simeq \pic(\Sp_\text{$G$-gl}^{\otimes})}$ for all compact Lie groups $G$.  
        In particular, $i$ restricts to an equivalence $\smash{\pic(\ul\mySp^{\otimes})\iso \pic(\ul\mySp_{\gl}^{\otimes})}$, so that $\ul{\mathfrak J}_\gl^\otimes$ factors uniquely through $\ul\mySp^\otimes$ as claimed.
    \end{proof}
\end{lemma}

\begin{definition}
    We define the \emph{equivariant Thom spectrum functor}
    \begin{equation}\label{eq:equiv-thom}
        \ul\Th^\otimes\colon\ul\Spc^\otimes_{\Orb\triangleright\Glo/\ul\VRep^{\oplus}}\to\ul\mySp^\otimes
    \end{equation}
    as the unique $\Orb$-cocontinuous symmetric monoidal extension of $\ul{\mathfrak J}^\otimes$. 
\end{definition}

\begin{rk}
    Writing $\Orb_\text{fin}\subset\Glo$ for the subcategory spanned by {finite} groups and injective homomorphisms, the underlying symmetric monoidal $\Orb_\text{fin}$-functor of $\ul\Th^\otimes$ was first constructed in \cite{juran-thesis}.
\end{rk}
\begin{rk}\label{rk:incoherent-description-Thom-equiv}
Let us fix a compact Lie group $G$, and write $\Orb_G$ for the nerve of the topological category of transitive $G$-spaces as in Example~\ref{ex:G-oo}. \cite{LNP}*{Lemma 6.12} decribes a functor $\fgt\colon\Orb_G\to\Orb$ lifting to an equivalence $\Orb_G\to\Orb_{/\BGcat{G}}$ sending $G/H$ to $\BGcat{H}\hookrightarrow\BGcat{G}$. If we abbreviate $\smash{\ul\VRep_G^{\oplus}\coloneqq\fgt^*\ul\VRep^{\oplus}}$, then we may view $(\ref{eq:equiv-thom})$ in degree $G$ as a symmetric monoidal functor 
\[
    \Th_G^\otimes\colon\Spc_{G/\ul\VRep^{\oplus}_G}^\otimes\coloneqq\PSh(\Orb_G)_{/\ul\VRep^{\oplus}_G}^\otimes\to\mySp_G^\otimes
\]
which we call the \emph{$G$-equivariant Thom spectrum functor}.
The fact that we started with global functors provides (coherent) symmetric monoidal equivalences \[\smash{\Th_\text{$G$}^\otimes\circ f^*\simeq f^*\circ\Th_{G'}^\otimes}\] for any homomorphism $f\colon G\to G'$ of compact Lie groups, while the fact that the global functor $\ul\Th$ is an \emph{Orb-parametrized} left adjoint implies that $\Th_G$ is a left adjoint for every compact Lie group $G$ and that for every \emph{injective} group homomorphism $f\colon G\to G'$, the corresponding Beck--Chevalley map provides a (non-monoidal) equivalence $f_!\circ\Th_\text{$G$}\simeq\Th_\text{$G'$}\circ f_!$.
\end{rk}

\begin{prop}\label{prop:equiv-thom-from-global}
    The square
    \begin{equation}\label{diag:equiv-Thom-from-global}
        \begin{tikzcd}
            \ul\Spc^\otimes_{\Orb\triangleright\Glo/\ul\VRep^{\oplus}}\arrow[d,"L"']\arrow[r,"\ul\Th^\otimes"] &[.5em] \ul\mySp^\otimes\arrow[d,hook]\\
            \ul\Spc_{\Glo/\ul\VRep^{\oplus}}^\otimes\arrow[r,"\ul\Th_\gl^\otimes"'] & \ul\mySp_\gl^\otimes
        \end{tikzcd}
    \end{equation}
    of symmetric monoidal global functors commutes up to canonical equivalence, where the map on the left is the unique symmetric monoidal left adjoint compatible with the Yoneda embeddings and the map on the right is the unique symmetric monoidal left adjoint.
    \begin{proof}
        By construction, both paths through the diagram are symmetric monoidal $\Orb$-cocontinuous functors extending $\ul{\mathfrak J}_\gl^\otimes$, so the claim follows from the universal property of $\Orb$-cocompletion.
    \end{proof}
\end{prop}

While the vertical functors in $(\ref{diag:equiv-Thom-from-global})$ do not have \emph{global} right adjoints (they only preserve $\Orb$-colimits, not all $\Glo$-colimits), their underlying $\Orb$-functors do admit right adjoints. We then have:

\begin{thm}\label{thm:Thom-vs-fgt}
    The induced lax symmetric monoidal structures on the right adjoints of the vertical maps in $(\ref{diag:equiv-Thom-from-global})$ are in fact strong, and the Beck--Chevalley map
    \begin{equation}\label{diag:Thom-fgt-BC}
        \begin{tikzcd}
            \ul\Spc_{\Glo/\ul\VRep^{\oplus}}^\otimes|_{\Orb}\arrow[from=d]\arrow[r,"\ul\Th^\otimes"] \arrow[dr,Rightarrow,shorten=15pt]&[1.5em] \ul\mySp^\otimes|_{\Orb}\arrow[from=d]\\
            \ul\Spc^\otimes_{\Orb\triangleright\Glo/\ul\VRep^{\oplus}}|_{\Orb}\arrow[r,"\ul\Th_\gl^\otimes"'] & \ul\mySp_\gl^\otimes|_{\Orb}
        \end{tikzcd}
    \end{equation}
    is an equivalence.
\end{thm}

The proof will require some preparations.

\begin{lemma}
    The lax symmetrix monoidal $\Orb$-right adjoint $\ul\mySp_\gl^\otimes|_\Orb\to\ul\mySp^\otimes|_\Orb$ of the inclusion is strong symmetric monoidal. Moreover, it is $\Orb$-cocontinuous.
    \begin{proof}
        Let us begin by showing that the right adjoint $r$ is strong symmetric monoidal. By Lemma~\ref{lemma:equiv-spectra-vs-gl-spectra}, $r$ is given in each degree by a localization at the $G$-equivariant weak equivalences. Thus, the Beck--Chevalley condition for the unit is automatic, while the Beck--Chevalley condition for the tensor product amounts to saying that a tensor product of $G$-equivariant weak equivalences is a $G$-equivariant weak equivalence. Plugging in the definition of the tensor product on $\ul\mySp_\gl^\otimes$ as the derived smash product of flat $G$-orthogonal spectra, this is then a direct consequence of Theorem~\ref{thm:equiv-flatness}.
        
        It remains to show $\Orb$-cocontinuity. By \cite{schwede2018global}*{Corollary 3.1.37} equivariant weak equivalences are stable under wedge sum, so the same argument as above shows that $r$ preserves fiberwise coproducts. As it moreover preserves fiberwise finite colimits by stability, this shows that $f$ is fiberwise cocontinuous. Finally, Proposition~\ref{prop:i!-left-Quillen-equiv-sp} implies by the same kind of argument that $r$ satisfies the Beck--Chevalley condition for the left adjoints to restrictions along injective homomorphisms, so that $r$ is $\Orb$-cocontinuous by Lemma~\ref{lm:continuouity-clefts}.
    \end{proof}
\end{lemma}

\begin{lemma}\label{lemma:right-adjoint-right-Bousfield}
    Let $S\subset T$ be a cleft and let $X\in\Spc_T$. Then the $S$-parametrized right adjoint $R$ of the inclusion $\ul\Spc_{S\triangleright T/X}\hookrightarrow\ul\Spc_{T/X}$ is $S$-cocontinuous. Moreover, if $X$ is equipped with a commutative monoid structure, then the lax symmetric monoidal structure on $R$ (with respect to the slice symmetric monoidal structures) is strong.
    \begin{proof}
        By \cite{CLL_Clefts}*{Lemmas~3.16 and~3.17}, we have a right Bousfield localization $\incl\colon\ul\Spc_{S\triangleright T}|_S\rightleftarrows\ul\Spc_T|_S\noloc r$. It follows formally that the Beck--Chevalley map
        \[
            \begin{tikzcd}
                (\ul\Spc_{T/X})|_S\arrow[r,"R"]\arrow[d,"\fgt"'] & (\ul\Spc_{S\triangleright T/X})|_S\arrow[d,"\fgt"]\arrow[dl,Rightarrow,shorten=10pt]\\
                (\ul\Spc_T)|_S\arrow[r,"r"'] & \ul\Spc_{S\triangleright T}|_S
            \end{tikzcd}
        \]
        is an equivalence. By \cite{CLL_Clefts}*{Proposition 4.28}, the bottom map $r$ is $S$-cocontinuous; as the vertical maps are conservative and $S$-cocontinuous, we conclude that also $R$ is $S$-cocontinuous. Similarly, $r$ is strong symmetric monoidal with respect to the cartesian symmetric monoidal structures (being a right adjoint), so that $R$ is strong symmetric monoidal with respect to the the slice symmetric monoidal structures. 
    \end{proof}
\end{lemma}

\begin{proof}[Proof of Theorem~\ref{thm:Thom-vs-fgt}]
    That the lax symmetric monoidal structures on the right adjoints in $(\ref{diag:Thom-fgt-BC})$ are strong was verified as part of the previous two lemmas. Moreover, the lemmas also show that all functors in $(\ref{diag:Thom-fgt-BC})$ are $\Orb$-cocontinuous, so that it suffices to verify the Beck--Chevalley condition after restricting to the Yoneda image.

    If we write $R$ for the left-hand vertical map in $(\ref{diag:Thom-fgt-BC})$ and $r$ for the right-hand one, then the Beck--Chevalley map is given by definition on an object $X\in\ul\Spc_{\Orb\triangleright\Glo/\ul\VRep}|_{\Orb}(\BGcat{G})$ by the composite
    \begin{multline*}
        (\Th_G\circ R)(X)\xrightarrow{\;\eta\;}(r\circ\incl\circ\Th_G\circ R)(X)\\\iso(r\circ \Th_\text{$G$-gl}\circ \incl\circ R)(X)\xrightarrow{\;\epsilon\;}(r\circ\Th_\text{$G$-gl})(X).
    \end{multline*}
    The map $\eta$ is an equivalence by full faithfulness of $\incl\colon\ul\mySp\hookrightarrow\ul\mySp_\gl$; moreover, if $X$ is contained in the Yoneda image (and hence in particular in the essential image of $\incl\colon \ul\Spc_{\Orb\triangleright\Glo/\ul\VRep}\hookrightarrow\ul\Spc_{\Glo/\ul\VRep}$), then so is $\epsilon\colon(\incl\circ R)(X)\to X$ by the same argument. Thus, the whole composite is an equivalence for any such $X$, finishing the proof.
\end{proof}

Evaluating at $\BGcat{G}$ for any compact Lie group $G$, this shows:

\begin{cor}
    We have a commutative diagram
    \[
        \begin{tikzcd}
            \Spc_{\textup{$G$-gl}/\triv_G\ul\VRep^{\oplus}}^\otimes\arrow[d]\arrow[r,"\Th_\textup{$G$-gl}"] &[1em] \mySp_\textup{$G$-gl}^\otimes\arrow[d]\\
            \Spc_{G/\ul\VRep^{\oplus}_G}^\otimes\arrow[r,"\Th_G"'] & \mySp_G^\otimes\rlap,
        \end{tikzcd}
    \]
    where the vertical map on the left is given by restriction along the unique lift $\Orb_G\to\Glo_{/\BGcat{G}}$ of $\fgt\colon\Orb_G\to\Glo$, while the map on the right is the localization at the equivariant equivalences (which is induced by the identity of the 1-category of orthogonal $G$-spectra).\qed
\end{cor}
Combining this with the global functoriality of $\ul\Th_{\gl}$, we in particular obtain for every compact Lie group $G$ a commutative diagram
    \[
        \begin{tikzcd}
            \Spc_{\textup{gl}/\ul\VRep^{\oplus}}^\otimes\arrow[d,"\res_G"']\arrow[r,"\Th_\textup{gl}"] &[1em] \mySp_\textup{gl}^\otimes\arrow[d,"\res_G"]\\
            \Spc_{G/\ul\VRep^{\oplus}_G}^\otimes\arrow[r,"\Th_G"'] & \mySp_G^\otimes\rlap.
        \end{tikzcd}
    \]
\begin{rk}\label{rk:Thom-pic}
    Analogously to the non-equivariant setting discussed before, one can define more general versions of the global and equivariant Thom spectrum functors by instead extending the inclusions $\pic(\ul\mySp_\gl^\otimes)\hookrightarrow\mySp_\gl^\otimes$ and $\pic(\ul\mySp^\otimes)\hookrightarrow\ul\mySp^\otimes$, see also~\cite{horev-klang-zou} for the corresponding definition in the setting of $G$-$\infty$-categories for a fixed finite group $G$. For us, however, the approach where the $J$-homomorphism is already built into the definitions will be more convenient. In particular, this is closer to the geometrical intuition of virtual equivariant vector bundles, and the examples of global (homotopical) bordism spectra we are interested in arise much more naturally in this setting.\footnote{In fact, in order to realize these examples over $\pic$ we would have to define the global $J$-homomorphism anyhow.} Moreover, while we will see in the next section that the above global Thom spectrum construction can be equivalently described in terms of a concrete pointset level construction, we do not even have a candidate for a pointset level model of the more general Thom spectrum construction for global spaces over $\pic(\ul\mySp_\gl^\otimes)$.
\end{rk}

\section{Comparison with classical global Thom spectra}\label{sec:model-thom}
In this section we will describe a model-categorical construction of a global Thom spectrum in the form of an enriched symmetric monoidal functor from the category of orthogonal spaces over a certain \emph{periodic global Grassmannian} $\cat{BOP}$ to the category of orthogonal spectra. On the level of ordinary categories, this will merely be a minor (and inconsequential) variation of a non-equivariant construction due to Sagave and Schlichtkrull \cite{sagave-schlichtkrull-thom}. What is new here is the observation, due to Schwede, that this construction has excellent model-categorical properties, from which we will deduce that the induced functor on continuous Borel categories  $\Ntop(\cat{$\cat{L}$-Top}_{/\cat{BOP}}^{\boxtimes,\dual})\to\Ntop(\osp^{\smashp,\dual})$ descends to a global left adjoint
\[
    \Ntop(\cat{$\cat{L}$-Top}_{/\cat{BOP}}^{\boxtimes,\dual})[\text{$G$-gl.~weak equiv.}^{-1}]\to\ul\mySp_\gl^\otimes.
\]
The main goal of this section (Theorem~\ref{thm:comparison-Thom}) will then be to identify this with the parametrized global Thom spectrum functor from the previous section.

\subsection{Global Thom spectra via model categories}\label{subsec:model-thom}
In this subsection we will describe our pointset model of the ($G$-)global Thom spectrum functor. We will restrict ourselves throughout to giving the necessary definitions, and most of the proofs of the good homotopical properties of the global Thom spectrum functor are relegated to Appendix~\ref{app:thom-Stefan-functor} by Stefan Schwede.

\medskip
\subsubsection{Thom spaces of (equivariant) vector bundles}
We begin with a brief review of Thom spaces of (equivariant) Euclidean vector bundles:

\begin{constr}
    Let $\zeta\colon E\to B$ be a Euclidean vector bundle. Its \emph{Thom space} $T_\zeta$ is the topological space with underlying set $E\cup\{\infty\}$, topologized in such a way that
    \begin{equation}\label{eq:thom-space-tolopogify}
        \begin{aligned}
            S(E)\times[0,1]&\longrightarrow T_\zeta\\
            (e,\lambda)&\longmapsto \begin{cases}
                {\lambda\over1-\lambda}\cdot e & \lambda\not=1\\
                \infty & \lambda=1
            \end{cases}
        \end{aligned}
    \end{equation}
    is a quotient map, where $S(E)\coloneqq\{e\in E: \|e\|=1\}$ denotes the sphere bundle. We view $T_\zeta$ as a based space, with basepoint $\infty$. If $G$ is a compact Lie group and $\zeta$ is a Euclidean $G$-vector bundle (i.e. $G$ acts continuously on $E$ and $X$ by Euclidean bundle maps), then we view $T_\zeta$ as a based $G$-space via the unique extension of the given action on $E$. Similarly, given a map of Euclidean $G$-bundles, i.e.\ a pullback square
    \begin{equation}\label{diag:map-of-G-bun}
        \begin{tikzcd}
            E_1\arrow[d,"\zeta_1"']\arrow[r,"\bar f"]\arrow[dr,pullback] & E_2\arrow[d,"\zeta_2"]\\
            B_1\arrow[r,"f"'] & B_2
        \end{tikzcd}
    \end{equation} of topological $G$-spaces
    such that $\bar f$ is fiberwise linear and isometric, we define $T_{f}$ as the unique based extension of $\bar f$.
\end{constr}

\begin{rk}
    For any $b\in B$, the map $(\ref{eq:thom-space-tolopogify})$ identifies all points $\zeta|_{S(E)}^{-1}(b)\times\{0\}$ with each other, and it moreover identifies all points in $S(E)\times\{1\}$. Thus, $(\ref{eq:thom-space-tolopogify})$ descends to a homeomorphism from the (unreduced) mapping cone of $\zeta|_{S(E)}\colon S(E)\to B$. Similarly, we can also describe the topology on $T_{\zeta}$ as the quotient topology with respect to 
    \begin{align*}
        D(E)&\longrightarrow T_{\zeta}\\
        e&\longmapsto 
        \begin{cases}
            {1\over 1-\|e\|}\cdot e & \|e\|\not=1\\
            \infty & \|e\|=1,\\
        \end{cases}
    \end{align*}
    where $D(E)=\{e\in E: \|e\|\le 1\}$ denotes the \emph{disk bundle} of $E$, i.e.~$T_\zeta$ is homeomorphic to $D(E)/S(E)$.
\end{rk}

\begin{ex}\label{ex:thom-space-1pt}
    If $B$ is compact, then $(\ref{eq:thom-space-tolopogify})$ is a continuous surjection when we equip the target with the topology of the 1-point compactification. As the latter is Hausdorff, this shows that $T_\zeta$ is just the 1-point compactification in this case. In particular, if $V$ and $W$ are inner product spaces, then the mapping space $\cat{O}(V,W)$ in the indexing category of orthogonal spectra is \emph{equal} to the Thom space of the bundle $\omega\coloneqq\pr\colon\cat{O}_0(V,W)\coloneqq\{(i,w)\in\cat{L}(V,W)\times W: w\perp i(V)\}\to\cat{L}(V,W)$.
\end{ex}

\begin{rk}\label{rk:classical-thom-sym-mon}
    If $\zeta_i\colon E_i\to B_i$, $i=1,2$ are two Euclidean $G$-bundles, then we obtain a map $T_{\zeta_1}\smashp T_{\zeta_2}\to T_{\zeta_1\times\zeta_2}$ by extending $e_1\smashp e_2\mapsto (e_1,e_2)$. By \cite{husemoller}*{proof of Proposition~16.1.5}, this is a homeomorphism. 
\end{rk}

\begin{constr}\label{constr:thom-over-Gr}
    Given any $G$-representation $V$ and any map of $G$-spaces $f\colon A\to\cat{Gr}(V)$, we write $T(f)$ for the Thom space $T_{f^*\zeta}$ of the pullback of the tautological bundle $\zeta\coloneqq\pr\colon \{(U,u) : U\in\cat{Gr}(V),u\in U\}\to\cat{Gr}(V)$. This then becomes a topological functor $T\colon\cat{$\bm G$-Top}_{/\cat{Gr}(V)}\to\cat{$\bm G$-Top}_*$ in the obvious way.
\end{constr}

The following result will be proven as Lemma~\ref{app-lemma:thom-spaces-homotopical} in the appendix:

\begin{lemma}\label{lemma:classical-thom-homotopical}
    Let $G$ be a compact Lie group, let $V$ be a $G$-representation, and let $\Ff$ be a family of closed subgroups. Then $T\colon\cat{$\bm G$-Top}_{/\cat{Gr}(V)}\to\cat{$\bm G$-Top}_*$ preserves $\Ff$-weak equivalences.\qed
\end{lemma}

\subsubsection{Global Thom spaces} Now we come to the construction of the \emph{global} Thom space functor. We begin by introducing what will later turn out to be the model analogues of the symmetric monoidal global space $\ul\Vect^\oplus$ and its group completion $\ul\VRep^{\oplus}$. These naturally come to us as highly structured objects:

\begin{definition}
    An \emph{ultra-commutative monoid} is a (strictly) commutative algebra in $\cat{$\cat{L}$-Top}$ with respect to the box product.
\end{definition}

Recall that the box product on $\cat{$\cat{L}$-Top}$ was constructed as enriched Day convolution, and so an ultra-commutative monoid structure on an orthogonal space $X$ can be equivalently described as a collection of \emph{multiplication maps} $\mu_{V,W}\colon X(V)\times X(W)\to X(V\oplus W)$ for all $V,W\in\cat{L}$ together with a \emph{unit map} $\eta\colon 1\to X(0)$ satisfying the expected naturality, unitality, associativity, and commutativity relations.

\begin{ex}[See \cite{schwede2018global}*{Example~2.3.12}]
    For any $n\ge0$ we define $\cat{Gr}_n$ as the orthogonal space $\cat{L}(\R^n,-)/\O(n)$ where $\O(n)$ acts via the tautological action on $\R^n$. Note that this means that two embeddings $\R^n\rightrightarrows V$ are identified in $\cat{Gr}_n(V)$ if and only if they have the same image, and we will therefore typically denote points in $\cat{Gr}_n(V)$ as subspaces of $V$. In this notation, the structure map $\cat{Gr}_n(V)\to\cat{Gr}_n(W)$ for a linear isometric embedding $f\colon V\to W$ is given by $U\mapsto f(U)$.

    The \emph{global Grassmannian} $\cat{Gr}$ is then the orthogonal space $\cat{Gr}\coloneqq\coprod_{n\ge0}\cat{Gr}_n$. This becomes an ultra-commutative monoid whose unit is the unique point of $\cat{Gr}(0)$ and with multiplications given by the coproduct of the maps
    \begin{align*}
        \cat{Gr}_m(V)\times\cat{Gr}_n(W)&\longrightarrow\cat{Gr}_{m+n}(V\oplus W)\\
        (U_1,U_2)&\longmapsto U_1\oplus U_2.
    \end{align*}
\end{ex}

\begin{ex}[See \cite{schwede2018global}*{Example~2.4.1}]\label{ex:BOP}
    The \emph{periodic global Grassmannian} $\cat{BOP}$ is given on an inner product space $V$ by $\cat{BOP}(V)=\cat{Gr}(V\oplus V)$, but with the structure map associated to a linear isometry $f\colon V\to W$ defined as
    \begin{align*}
        \cat{Gr}(V\oplus V)&\longrightarrow\cat{Gr}(W\oplus W)\\
        U&\longmapsto (f\oplus f)(U) + (\im(f)^\perp\oplus 0)\rlap.
    \end{align*}
    One can check that this becomes an ultra-commutative monoid with unit the unique point of $\cat{BOP}(0)$ and with multiplication maps $\cat{Gr}_m(V\oplus V)\times\cat{Gr}_n(W\oplus W)\to\cat{Gr}_{m+n}((V\oplus W)\oplus (V\oplus W))$ given by sending $(U_1,U_2)$ to the image of $U_1\oplus U_2$ under the isomorphism $(V\oplus V)\oplus (W\oplus W)\cong (V\oplus W)\oplus (V\oplus W), (v_1,v_2;w_1,w_2)\mapsto (v_1,w_1;v_2,w_2)$. We then have a map of ultra-commutative monoids $i\colon\cat{Gr}\to\cat{BOP}$ given at an inner product space $V$ by $\cat{Gr}(V)\to\cat{Gr}(V\oplus V),U\mapsto V\oplus U$.
\end{ex}

Our model-categorical global Thom space functor will be defined as a certain enriched symmetric monoidal functor $\cat{th}^\otimes\colon\cat{$\cat{L}$-Top}_{/\cat{Gr}}\to\cat{$\cat{L}$-Top}_*$, where we equip the source with the 1-categorical analogue of the \emph{slice symmetric monoidal structure}, i.e.\ the tensor product is given on objects by sending $(f\colon A\to \cat{Gr}, g\colon B\to \cat{Gr})$ to the composite
\[
    A\boxtimes B\xrightarrow{\;f\boxtimes g\;} \cat{Gr}\boxtimes \cat{Gr}\xrightarrow{\;\mu\;} \cat{Gr},
\]
while the symmetric monoidal product on morphisms and the unitality, associativity, and symmetry isomorphisms are defined in the unique way for which the forgetful functor to $\cat{$\cat{L}$-Top}$ is \emph{strict} symmetric monoidal. Note that since the mapping space between $(f\colon A\to \cat{Gr})$ and $(g\colon B\to \cat{Gr})$ in $\cat{$\cat{L}$-Top}_{/\cat{Gr}}$ is a subspace of the mapping space $\maps_{\cat{$\cat{L}$-Top}}(A,B)$, this symmetric monoidal structure is again topological. 

\begin{constr}\label{constr:Thom-space-functor-gl}
    For any $V\in\cat{L}$, we define 
    \[
        \cat{EGr}(V)\coloneqq\{(v,U)\in V\times\cat{Gr}(V) : v \in U\},
    \]
    and we write $\zeta_V\coloneqq\pr_2\colon\cat{EGr}(V)\to\cat{Gr}(V)$ for the tautological bundle. If we make $\cat{EGr}$ into an orthogonal space by letting $u\colon V\to W$ act via $(v,U)\mapsto (u(v),u(U))$, then $\zeta$ is a map of orthogonal spaces. 

    Given now any object $f\colon B\to\cat{Gr}$ of $\cat{$\cat{L}$-Top}_{/\cat{Gr}}$, we define $\cat{th}(f)$ as the based orthogonal space sending $V\in\cat{L}$ to the Thom space of the pulled back bundle $f_V^*\zeta_V\colon\cat{EGr}(V)\times_{\cat{Gr}} B(V)\to B(V)$. A linear isometric embedding $u\colon V\to W$ acts via the map $T_{f_U^*\zeta_U}\to T_{f^*_V\zeta_V}$ induced by the fiber product of the structure maps
    \[
        \cat{EGr}(U)\times_{\cat{Gr}(U)}{B}(U)\to\cat{EGr}(V)\times_{\cat{Gr(V)}} B(V);
    \]
    in other words, $\cat{th}(f)$ is obtained by levelwise applying the Thom space functor to the pulled back map $f^*\zeta$ in $\cat{OrthSpc}$. In the notation of Construction~\ref{constr:thom-over-Gr}, we then have $\cat{th}(f)(V)=T(f_V\colon B(V)\to\cat{Gr}(V))$ for any $V\in\cat{L}$.
    Given a commutative diagram
    \[
        \begin{tikzcd}[column sep=small]
            A\arrow[rr,"\alpha"]\arrow[dr, bend right=15pt,"f"'] && B\arrow[dl,"g", bend left=15pt]\\
            &\cat{Gr}
        \end{tikzcd}
    \]
    in $\cat{$\cat{L}$-Top}$ we similarly define $\cat{th}(\alpha)\colon\cat{th}(f)\to\cat{th}(g)$ as the map given in degree $V$ by $T(\alpha_V\colon f_V\to g_V)$.
\end{constr}

\begin{ex}\label{ex:th-on-corep}
    Let $V\in\cat{L}$ and let $W\subset V$ be any subspace. If we write $w\colon\cat{L}(V,-)\to\cat{Gr}$ for the map corresponding to $W\in\cat{Gr}(V)$, then the pulled back bundle $w^*\zeta_U\colon\cat{EGr}(U)\times_{\cat{Gr}(U)}\cat{L}(V,U)\to\cat{L}(V,U)$ admits a trivialization for any $U\in\cat{L}$ by sending a point $(i,x)\in\cat{L}(V,U)\times W$ to $(i(x),i(W);i)$. It follows that  we have an isomorphism $\cat{L}(V,-)_+\smashp S^W\cong\cat{th}(w)$ given in degree $U$ by sending $i\smashp x$ to $\infty$ if $x=\infty$ and to $(i(x),i(W);i)\in\cat{EGr}(U)\times_\cat{Gr(U)}\cat{L}(V,U)$ if $x\in W$.
\end{ex}

The following result will again be proven in the appendix:

\begin{prop}[See Theorem~\ref{thm:Thom-marvelous-appendix}(1)]\label{prop:thom-space-marvelous}
    For every compact Lie group $G$,
    \[
        \cat{th}\colon\cat{$\bm G$-$\cat{L}$-Top}_{/\cat{Gr}}\to\cat{$\bm G$-$\cat{L}$-Top}_*
    \]
    is homotopical and left Quillen.\qed
\end{prop}

We now explain in what sense the isomorphisms from Example~\ref{ex:th-on-corep} are natural:

\begin{constr}\label{constr:th-on-corep}
    Write $\cat{y}\colon\cat{L}^\op\to\cat{$\cat{L}$-Top}$ for the enriched Yoneda embedding and $\cat{y}_{/\cat{Gr}}$ for the comma category as usual, so that objects of $\cat{y}_{/\cat{Gr}}$ are given by maps $w\colon\cat{L}(V,-)\to\cat{Gr}$ while a morphism $w\to w'$ is given by an isometric embedding $\alpha\colon V\to V'$ making the diagram    
    \begin{equation}\label{diag:morphism-y-Gr}
        \begin{tikzcd}[column sep=small]
            \cat{L}(V',-)\arrow[rr,"-\circ\alpha"]\arrow[dr, bend right=15pt,"w'"'] && \cat{L}(V,-)\arrow[dl, bend left=15pt,"w"]\\
            & \cat{Gr}
        \end{tikzcd}
    \end{equation}
    commute (note the change in variance).

    We define an enriched functor $(\cat{y}_{/\cat{Gr}})^\op\to\cat{L}$ by sending $w\colon\cat{L}(V,-)\to\cat{Gr}$ to $w(\id_V)\in\cat{Gr}(V)$ and by sending a map $\alpha\colon w\to w'$ to the \emph{inverse} of the restriction of $\alpha\colon V\to V'$ to $w(\id_{V})\to w'(\id_{V'})$; note that this is indeed well-defined since $\alpha(w(\id_V))=w(\alpha)=w'(\id_{V'})$ by commutativity of $(\ref{diag:morphism-y-Gr})$. Postcomposing with $\cat{j}\colon\cat{L}\to\cat{Top}_*$ then gives a map $\cat{k}\colon(\cat{y}_{/\cat{Gr}})^\op\to\cat{Top}_*$ sending $w\colon\cat{L}(V,-)\to\cat{Gr}$ to the sphere $S^{w(\id_V)}$. One directly checks that the isomorphisms from the previous example define a natural isomorphism between the restriction of $\cat{th}$ and the composite
    \begin{equation}\label{diag:thom-space-on-yoneda}
        (\cat{y}_{/\cat{Gr}})^\op\xrightarrow{\;(\cat{k},(-)_+\circ\cat{y}\circ\fgt)\;}
        \cat{Top}_*\times\cat{$\cat{L}$-Top}_*\xrightarrow{-\smashp-}
        \cat{$\cat{L}$-Top}_*.
    \end{equation}
\end{constr}

If we equip $\cat{y}_{/\cat{Gr}}$ with the symmetric monoidal structure inherited from $\cat{L}$, then 
\[
    \big((w\colon\hskip0pt minus 1pt\cat{L}(V,-)\hskip0pt minus 1pt\to\hskip0pt minus 1pt\cat{Gr})\boxtimes (w'\hskip0pt minus .5pt\colon\hskip0pt minus 1pt\cat{L}(V',-)\hskip0pt minus 1pt\to\hskip0pt minus 1pt\cat{Gr})\big)(\id_{V\oplus V'})\hskip0pt minus .5pt=\hskip0pt minus .5ptw(\id_V)+w'(\id_{V'})\hskip0pt minus .5pt\subset\hskip0pt minus .5pt V\oplus V'\rlap,
\]
and one directly checks that $\cat{k}$ becomes a symmetric monoidal functor via the evident isomorphism $S^{w(\id_V)}\smashp S^{w'(\id_{V'})}\cong S^{w(\id_V)\oplus w'(\id_{V'})}$ and the unique isomorphism $S^0\cong S^{\{0\}}$. As also all the functors in $(\ref{diag:thom-space-on-yoneda})$ have natural symmetric monoidal structures, this upgrades the whole composite to a strong symmetric monoidal functor. We will now explain how one can make $\cat{th}$ into a strong symmetric monoidal functor $\cat{th}^\otimes$ such that the natural isomorphism from the previous construction is a symmetric monoidal isomorphism.

\begin{lemma}\label{lemma:extend-nat-trafo}
    Let $n\ge0$, $\cat{I}_1,\dots,\cat{I}_n$ be topologically enriched categories, and let $X_i\in\cat{Fun}_\cont(\cat{I}_i,\cat{Top})$ for $i=1,\dots,n$ be arbitrary. Moreover, let $\cat{D}$ be any cocomplete category enriched and tensored over $\cat{Top}$, and let
    \[F,G\colon\prod_{i=1}^n\cat{Fun}_\cont(\cat{I}_i,\cat{Top})_{/X_i}\rightrightarrows\cat{D}\]
    be topological functors that preserve colimits and tensors in each variable. Then any natural transformation $\sigma$ between the restrictions $\prod_{i=1}^n(\cat{y}_{/X_i})^\op\rightrightarrows\cat{D}$ of $F$ and $G$ admits a unique extension to a natural transformation $F\to G$.
    \begin{proof}
        We will give the proof for $n=1$; the general case then follows inductively by adjoining over one of the factors.

        Given $Y\colon\cat{I}\to\cat{Top}$, we have for all $i,j\in\cat{I}$ maps 
        \begin{align*}
            \alpha_{i,j}\colon\maps(j,-)\times\maps(i,j)\times Y(i)&\to\maps(i,-)\times Y(i)\\
            \beta_{i,j}\colon\maps(j,-)\times\maps(i,j)\times Y(i)&\to\maps(j,-)\times Y(j)
        \end{align*}
        via precomposition and the functoriality of $Y$, respectively. Moreover, we have for every $i\in \cat{I}$ a unique map $\gamma_i\colon\maps(i,-)\times Y(i)\to Y$ such that $\gamma_i(\id_i,-)=\id$. These then assemble into a diagram
        \begin{equation}\label{eq:coyoneda}
            \begin{tikzcd}[column sep=1.5em]
                \coprod\limits_{i,j\in I} \maps(j,-)\times\maps(i,j)\times Y(i)\arrow[r, yshift=2pt,"\alpha"]\arrow[r,yshift=-2pt,"\beta"']& \coprod\limits_{i\in I} \maps(i,-)\times Y(i) \arrow[r,"\gamma"] & Y\rlap,
            \end{tikzcd}
        \end{equation}
        in $\Fun_\cont(\cat{I},\cat{Top})$ and the enriched Yoneda lemma says that this is a coequalizer diagram. Given any map $f\colon Y\to X$, we can then uniquely lift $(\ref{eq:coyoneda})$ to a diagram in $\Fun_\cont(\cat{I},\cat{Top})_{/X}$ expressing $f$ as a coequalizer. As $F$ and $G$ are assumed to preserve tensors, this then similarly expresses $F(Y)$ and $G(Y)$ as coequalizers. As any  transformation $F\to G$ is automatically enriched and hence preserves tensors, this shows that in any extension of $\sigma$ to a natural transformation $\bar\sigma\colon F\to G$, the map $\bar \sigma_Y\colon F(Y\to X)\to G(Y\to X)$ must be the map induced on coequalizers by
        \[
            \begin{tikzcd}[column sep=1.1em,cramped]
                \arrow[d,"{\coprod_{i,j}\eta\otimes (\maps(i,j)\times Y(i))}"']\coprod\limits_{i,j\in I} F(\maps(j,-)\to X)\otimes (\maps(i,j)\times Y(i))\arrow[r, yshift=2pt]\arrow[r,yshift=-2pt]& \coprod\limits_{i\in I} F(\maps(i,-)\to X)\otimes Y(i)\arrow[d,"\coprod_{i}\eta\otimes Y(i)"]\\
                \coprod\limits_{i,j\in I} G(\maps(j,-)\to X)\otimes (\maps(i,j)\times Y(i))\arrow[r, yshift=2pt]\arrow[r,yshift=-2pt]& \coprod\limits_{i\in I} G(\maps(i,-)\to X)\otimes Y(i)\rlap.
            \end{tikzcd}
        \]
        In summary, we have shown that there exists for every $Y$ a unique map $\bar\sigma$ making each of the diagrams
        \begin{equation}\label{diag:sigma-bar-compat}
            \begin{tikzcd}
                F(\maps(i,-)\to X)\otimes Y(i)\arrow[d,"\sigma"']\arrow[r]& F(Y)\arrow[d,"\bar\sigma"]\\
                G(\maps(i,-)\to X)\otimes Y(i)\arrow[r]& G(Y)
            \end{tikzcd}
        \end{equation}
        commute. It only remains to show that these maps are indeed natural, for which we consider any map $\alpha\colon (f\colon Y\to X)\to (g\colon Z\to X)$ in $\Fun_\cont(\cat{I},\cat{Top})_{/X}$. We want to show that $G(\alpha)\circ\bar\sigma_f=\bar\sigma_g\circ F(\alpha)$, which can be checked after restricting along the epimorphism $\coprod_{i\in I} F(\maps(i,-)\to X)\otimes Y(i)\to F(Y)$ and hence after restricting along each $F(\maps(i,-)\to X)\otimes Y(i)$. Using the compatibility relation $(\ref{diag:sigma-bar-compat})$ for both $Y$ and $Z$ the claim then translates to commutativity of
        \[
            \begin{tikzcd}
                F(\maps(i,-)\to X)\otimes Y(i)\arrow[d, "\sigma\otimes\id"']\arrow[r,"\id\otimes\alpha(i)"] &[1em] F(\maps(i,-)\to X)\otimes Z(i)\arrow[d,"\sigma\otimes\id"]\\
                G(\maps(i,-)\to X)\otimes Y(i)\arrow[r,"\id\otimes\alpha(i)"'] & G(\maps(i,-)\to X)\otimes Z(i)\rlap,
            \end{tikzcd}
        \]
        which is clear.        
    \end{proof}
\end{lemma}

\begin{rk}
    It seems tempting to try to prove the theorem by first identifying $\cat{Fun}_\cont(\cat{I},\cat{Top})_{/X}\simeq\cat{Fun}_\cont(\cat{y}_{/X},\cat{Top})$ and then appealing to the result for terminal $X$ from \cite{kelly-enriched}*{Theorem~4.51}. However, while such an identification holds in the unenriched (or $\infty$-categorical) setting, it fails here for pointset topological reasons: by \cite{kelly-enriched}*{Theorem 5.26} this would entail that the corepresented functor of any $x\colon\maps(i,-)\to X$ be cocontinuous---this corepresented functor sends $Y\to X$ to the fiber of $Y(i)$ over $\id_i$, and taking fibers is \emph{not} in general cocontinuous in $\cat{Top}$. 
\end{rk}

\begin{prop}\label{prop:th-on-corep-sym-mon}
    There exists a unique strong symmetric monoidal structure on $\cat{th}$ such that the isomorphism from Construction~\ref{constr:th-on-corep} becomes an isomorphism of symmetric monoidal functors.
    \begin{proof}
        As $\cat{y}_{/\cat{Gr}}$ contains the monoidal unit $0\colon\cat{L}(0,-)\to\cat{Gr}$, there is a unique isomorphism $\smash{S^0\cong\cat{th}(0)}$ compatible with the identification from Construction~\ref{constr:th-on-corep}; similarly, there is a unique compatible isomorphism $\cat{th}(-)\boxsmash\cat{th}(-)\to\cat{th}(-\boxtimes-)$ on $(\cat{y}_{/\cat{Gr}})^\op$. By Lemma~\ref{lemma:extend-nat-trafo} the latter uniquely extends to a transformation on all of $\cat{$\cat{L}$-Top}_{/\cat{Gr}}$, which is then necessarily an isomorphism. Moreover, appealing to uniqueness multiple times shows that this satisfies all the coherence conditions necessary for the structure map of a strong symmetric monoidal functor.
    \end{proof}
\end{prop}

Having established all the necessary model-categorical properties, we can now explain how $\cat{th}^\otimes$ induces a symmetric monoidal left adjoint of global $\infty$-categories.

\begin{constr}
    For any orthogonal space $X$ we write $\ul\myS_\gl[X]$ for the levelwise localization of the continuous Borel category $\Ntop(\cat{$\cat{L}$-Top}_{/X}^{\smash\dual})$ at the underlying ($G$-)global weak equivalences. If $X$ is an ultra-commutative monoid, then we write $\ul\myS_\gl[X]^\otimes$ for the analogous localization of the \emph{symmetric monoidal} continuous Borel construction applied to the slice symmetric monoidal structure on $\cat{$\cat{L}$-Top}_{/X}$. These two constructions define functors $\cat{$\cat{L}$-Top}\to\CAT_{\Glo}$ and $\CAlg(\smash{\cat{$\cat{L}$-Top}^\boxtimes})\to\CMon(\Cat_{\Glo})$ via the strict functoriality of the 1-categorical slice and the functoriality of the continuous Borel construction.
\end{constr}

\begin{rk}
    We have already seen in Proposition~\ref{prop:model-categorical-model} that $\ul\myS_\gl[X]$ is equivalent to the slice $\smash{(\ul\myS_\gl)_{/X}}$. The reason to introduce new notation here is that $\ul\myS_\gl[X]$ naturally comes to us as a functor in $X$, and we did not show that the aforementioned identification is natural in $X$; in fact, even defining the functoriality of the $\infty$-categorical slice would need tools like the parametrized unstraightening from \cite{martini2022cocartesianfibrationsstraighteninginternal} that we did not introduce here.    
    We have therefore decided to introduce the above ad-hoc notation in order to avoid any kind of ambiguity in the comparison to the constructions from the previous section.
\end{rk}

As $\smash{\cat{th}\colon\cat{$\bm G$-$\cat{L}$-Top}_{/\cat{Gr}}\to\cat{$\bm G$-$\cat{L}$-Top}_*}$ is homotopical for every $G$ (see Proposition~\ref{prop:thom-space-marvelous}), the continuous Borel construction $\smash{\Ntop(\cat{th}^\dual)}$ descends to a global functor $\ul\myth_\gl\colon\ul\myS_\gl[\cat{Gr}]\to\ul\myS_{\gl,*}$. For the symmetric monoidal version we have to work only slightly harder:

\begin{constr}\label{constr:gl-Thom-space-sym-mon}
    As $\cat{th}$ is left Quillen for every $G$, it preserves cofibrant objects, and hence also quasi-cofibrant ones (as introduced in Definition~\ref{defi:qcof}). As a functor $\Ntop(\cat{$\cat{L}$-Top}_{/\cat{Gr}}^{\text{qcof},\dual})\to\Ntop(\cat{$\cat{L}$-Top}_*^{\text{qcof},\dual})$ it is of course still homotopical; moreover, also the symmetric monoidal product on both sides is now homotopical (see Proposition~\ref{prop:boxsmash-qcof}), and hence the restriction of $\Ntop(\cat{th}^{\otimes,\dual})$ descends to a symmetric monoidal global functor
    \[
        \ul\myth_\gl^\otimes\colon\ul\myS_\gl[\cat{Gr}]^\otimes\simeq\ul\myS_\gl[\cat{Gr}]^{\text{qcof},\otimes}\to\ul\myS_{\gl,*}^{\text{qcof},\otimes}\simeq\ul\myS_{\gl,*}^\otimes
    \]
    refining $\ul\myth_\gl$.
\end{constr}

\begin{prop}\label{prop:gl-Thom-space-la}
    The global functor $\ul\myth_\gl$ is a left adjoint.
    \begin{proof}
        Proposition~\ref{prop:thom-space-marvelous} immediately implies the existence of a pointwise right adjoint. To see that these pointwise right adjoints satisfy the Beck--Chevalley condition, we may then equivalently show that $\ul\myth_\gl$ satisfies the Beck--Chevalley condition with respect to the left adjoints along restriction. This follows as usual by observing that the Beck--Chevalley condition holds on the pointset level and that all functors in sight are left Quillen.
    \end{proof}
\end{prop}

\subsubsection{Global Thom spectra} We will now construct a global Thom spectrum functor $\smash{\cat{$\cat{L}$-Top}_{/\cat{BOP}}^\boxtimes\to\osp^\smashp}$, where $\cat{BOP}$ is the ultra-commutative monoid from Example~\ref{ex:BOP}.

\begin{constr}
    Recall from Example~\ref{ex:thom-space-1pt} that the mapping space $\cat{O}(V,W)$ in the indexing category for orthogonal spectra agrees with the Thom space of the bundle $\omega=\pr\colon\cat{O}_0(V,W)=\{(i,w):w\perp\im(i)\}\to\cat{L}(V,W)$. As before we morever write 
    $\zeta=\pr\colon\cat{EGr}(V\oplus V)=\{(U,u) : u\in U\subset V\}\to\cat{Gr}(V\oplus V)=\cat{BOP}(V)$ for the tautological bundle. Then we have a commutative diagram
    \begin{equation*}
        \begin{tikzcd}
            \cat{O}_0(V,W)\times \cat{EGr}(V\oplus V)\arrow[d,"\omega\times \zeta"']\arrow[r] & \cat{EGr}(W\oplus W)\arrow[d,"\zeta"]\\
            \cat{L}(V,W)\times\cat{BOP}(V)\arrow[r,"\text{act}"'] & \cat{BOP}(W)
        \end{tikzcd}
    \end{equation*}
    where the bottom map is the action map and the top map sends a tuple $(i,w;U,u)$ to $((i\oplus i)(U)+(\im(i)^\perp\oplus 0),i(u)+(w,0))$. One directly checks that this is continuous and fiberwise linear bijective, so that the above square is in fact a pullback.
    
    By the pasting law we then get for any map $f\colon X\to\cat{BOP}$ of orthogonal spaces a pullback square
    \begin{equation}\label{diag:thomify-this!}
        \begin{tikzcd}
            \cat{O}_0(V,W)\times f_V^*\cat{EGr}(V\oplus V)\arrow[dr,pullback]\arrow[d,"\omega\times f_V^*\zeta"']\arrow[r] & f_W^*\cat{EGr}(W\oplus W)\arrow[d,"f_W^*\zeta"]\\
            \cat{L}(V,W)\times X(V)\arrow[r,"\text{act}"'] & X(W)
        \end{tikzcd}
    \end{equation}
    where the top map sends a tuple consisting of $i\colon V\to W$, $w\in\im(i)^\perp$, $x\in X(V)$, and $u\in f_V(x)\subset V\oplus V$ to $\big(X(i)(x)\in X(W),(w,0)+i(u) \in\cat{BOP}(i)(f_V(x))\big)$. 
\end{constr}

\begin{constr}\label{constr:Thom-spectrum-functor}
    Given a map $f\colon X\to\cat{BOP}$ of orthogonal spaces, we define the orthogonal spectrum $\cat{Th}(f)$ as follows: for every inner product space $V$, we set $\cat{Th}(f)(V)=T_{f_V^*\zeta_{V\oplus V}}=T(f_V\colon X(V)\to\cat{Gr}(V\oplus V))$; given another inner product space $W$, we define the structure map $\cat{O}(V,W)\smashp\cat{Th}(f)(V)\to \cat{Th}(f)(W)$ as the composite
    \[
        T_\omega\smashp T_{f_V^*\zeta}\cong T_{\omega\times f_V^*\zeta}\to T_{f_W^*\zeta}
    \]
    of the isomorphism from the symmetric monoidal structure on $T$ and the map induced by $(\ref{diag:thomify-this!})$. More concretely this says that for any $i\colon V\to W$, the structure map $S^{W-i(V)}\smashp\cat{Th}(f)(V)\to\cat{Th}(f)(W)$ is the based map sending $w\smashp [x,u]$ for $w\in W-i(V)$, $x\in X(V)$, and $u\in f_V(X)\subset V\oplus V$ to $[X(i)(x), (i\oplus i)(u)+(w,0)]$. We omit the straightforward computations that this indeed defines an orthogonal spectrum and that $\cat{Th}$ becomes a topological functor $\cat{$\cat{L}$-Top}_{/\cat{BOP}}\to\osp$ by sending a map $\alpha\colon f\to g$ in $\cat{$\cat{L}$-Sp}_{/\cat{BOP}}$ to the map given levelwise by \[T(\alpha_V)\colon\cat{Th}(f)(V)\to\cat{Th}(g)(V), \, [x,u]\mapsto [\alpha(x),u].\]
\end{constr}

\begin{ex}\label{ex:MO}
    The orthogonal spectrum $\cat{Th}(\id_{\cat{BOP}})$ is literally equal to the spectrum denoted $\cat{MOP}$ in \cite{schwede2018global}*{Example~6.1.7}. Similarly, if we write $\cat{BOP}^{[d]}\subset\cat{BOP}$ for the subspace given in degree $V$ by $\smash{\cat{Gr}_{d+\dim(V)}(V\oplus V)}$, then $\cat{Th}(\cat{BOP}^{[d]}\hookrightarrow\cat{BOP})$ is the spectrum denoted $\cat{MOP}^{[d]}$ in \emph{loc.\ cit.} In particular, taking $d=0$ this says that the \emph{global homotopical real bordism spectrum} $\cat{MO}$ is obtained by applying the above construction to  $\cat{BO}\coloneqq\smash{\cat{BOP}^{[0]}}\hookrightarrow\cat{BOP}$.
\end{ex}

\begin{ex}\label{ex:Th-on-corep}
    Let $V\in\cat{L}$ and let $w\colon\cat{L}(V,-)\to\cat{BOP}$ classify the subspace $W\coloneqq w(\id_V)\subset V\oplus V$. Analogously to Example~\ref{ex:th-on-corep} we have an isomorphism $\cat{O}(V,-)\smashp S^W\cong\cat{Th}(w)$ given in degree $U$ by sending $[i,u]\smashp x$ for $i\colon V\to U$, $u\in\im(i)^\perp$, and $x\in W$ to $[i,(i\oplus i)(x)+(u,0)]\in\cat{Th}(w)(U)$.
\end{ex}

Just like its unstable cousin, the global Thom space functor has excellent model-categorical properties:

\begin{thm}[See Theorem~\ref{thm:Thom-marvelous-appendix}(2)]\label{thm:thom-spectrum-marvelous}
    For any compact Lie group $G$,
    \[
        \cat{Th}\colon\cat{$\bm G$-$\cat{L}$-Top}_{/\cat{BOP}}\to\Gosp{G}
    \]
    is homotopical and left Quillen for the $G$-global model structures.\qed
\end{thm}

Similarly to our strategy for the global Thom space functor, we will now upgrade the computation from Example~\ref{ex:Th-on-corep} to a functorial comparison and use this to make $\cat{Th}$ into a symmetric monoidal functor.

\begin{constr}\label{constr:Th-on-corep}
    We define $\cat{f}\colon(\cat{y}_{/\cat{BOP}})^\op\to\osp$ as follows: for $V\in\cat{L}$ and $w\colon\cat{L}(V,-)\to\cat{BOP}$, we set $W\coloneqq w(\id_V)\subset V\oplus V$ and $\cat{f}(w)\coloneqq\cat{O}(V,-)\smashp S^{W}$. Given another map $w'\colon\cat{L}(V',-)\to\cat{BOP}$ and an isometric embedding $\alpha\colon V\to V'$ satisfying $w\circ\cat{L}(\alpha,-)=w'$ (i.e.\ $(\alpha\oplus\alpha)(W)+(\im(\alpha)^\perp\oplus0)=W'$), we define $\cat{f}(\alpha)$ as the adjunct of the map $S^{W'}\to\cat{O}(V,V')\smashp S^W$ sending an $x\in W'$ with components $x_1\in(\alpha\oplus\alpha)(W)$ and $x_2\in\im(\alpha)^\perp\subset V'\oplus 0$ to $[\alpha,x_2]\smashp (\alpha\oplus\alpha)^{-1}(x_1)$. We omit the straightforward verification that the diagrams
    \[
        \begin{tikzcd}
            \cat{O}(V',-)\smashp S^{W'}\arrow[r,"\cat{f}(\alpha)"]\arrow[d,"\cong","\text{Ex.~\ref{ex:Th-on-corep}}"'] &[2em] \cat{O}(V,-)\smashp S^W\arrow[d,"\cong"',"\text{Ex.~\ref{ex:Th-on-corep}}"]\\
            \cat{Th}(w')\arrow[r,"{\cat{Th}(\cat{L}(\alpha,-))}"'] & \cat{Th}(w)
        \end{tikzcd}
    \]
    commute, so that $\cat{f}$ is a topological functor isomorphic to the restriction of $\cat{Th}$.
\end{constr}

\begin{constr}
    Given $w\colon\cat{L}(V,-)\to\cat{BOP}$ and $w'\colon\cat{L}(V',-)\to\cat{BOP}$ with $W\coloneqq w(\id_V)$ and $W'\coloneqq w'(\id_{V'})$, then $(w\oplus w')(\id_{V\oplus V'})$ is the image of $W\oplus W'\subset (V\oplus V)\oplus (V'\oplus V')$ under the isomorphism $\kappa\colon V\oplus V\oplus V'\oplus V'\to V\oplus V'\oplus V\oplus V'$ swapping the two terms in the middle. We omit the easy verification that the composites
    \begin{multline*}
        \cat{f}(w)\smashp \cat{f}(w')=\big(\cat{O}(V,-)\smashp S^W\big)\smashp\big(\cat{O}(V',-)\smashp S^{W'}\big)\xrightarrow{\;\;\cong\;\;}\\\cat{O}(V\oplus V',-)\smashp S^{W\oplus W'}
        \xrightarrow[\raise3.5pt\hbox{$\scriptstyle\cong$}]{\;\;\kappa\;\;}\cat{O}(V\oplus V',-)\smashp S^{\kappa(W\oplus W')}=\cat{f}(w\boxtimes w')
    \end{multline*}
    together with the unique isomorphism $\cat{f}(0\colon\cat{L}(0,-)\to\cat{BOP})\cong\mathbb S$ make $\cat{f}$ into a strong symmetric monoidal functor.
\end{constr}

As in Proposition~\ref{prop:th-on-corep-sym-mon} we deduce:

\begin{cor}
    There exists a unique symmetric monoidal structure on $\cat{Th}$ such that the isomorphism from Construction~\ref{constr:Th-on-corep} is symmetric monoidal.\qed
\end{cor}

As in Construction~\ref{constr:gl-Thom-space-sym-mon} we then get a symmetric monoidal global functor
\[
    \ul\myTh_\gl^\otimes\colon\ul\myS_\gl[\cat{BOP}]^\otimes\simeq\ul\myS_\gl[\cat{BOP}]^{\text{qcof},\otimes}\to\mySp_\gl^{\text{flat},\otimes}\simeq\mySp_\gl^\otimes
\]
by applying the symmetric monoidal continuous Borel construction to $\cat{Th}^\otimes$. Moreover, we deduce from Theorem~\ref{thm:thom-spectrum-marvelous} by the same argument as before:

\begin{cor}\label{cor:gl-Thom-spectrum-la}
    The global functor $\ul\myTh_\gl$ is a left adjoint.\qed
\end{cor}

Finally, let us explain how the global Thom space and Thom spectrum functors relate to each other:

\begin{prop}
    The diagram
    \[
        \begin{tikzcd}
            \cat{$\cat{L}$-Top}_{/\cat{Gr}}\arrow[d," \cat{$\cat{L}$-Top}_{/i}"']\arrow[r,"\cat{th}"] & \cat{$\cat{L}$-Top}_*\arrow[d,"\Sigma^\bullet"]\\
            \cat{$\cat{L}$-Top}_{/\cat{BOP}}\arrow[r,"\cat{Th}"']&\cat{Sp}
        \end{tikzcd}
    \]
    of symmetric monoidal (topological) functors commutes up to symmetric monoidal isomorphism.
    \begin{proof}
        As before, it suffices to prove this after restriction to $(\cat{y}_{/\cat{Gr}})^\op$.
        
        By definition, the bottom composite becomes $\cat{f}\circ\cat{y}_{/i}^\op$. Recall that $i\colon\cat{Gr}\to\cat{BOP}$ was defined in degree $V$ by $U\mapsto V\oplus U$, so that this composite is given on objects by sending $w\colon\cat{L}(V,-)\to\cat{Gr}$ with $W\coloneqq w(\id_V)$ to $\cat{O}(V,-)\smashp S^{V\oplus W}$. We have an isomorphism $\cat{O}(V,-)\smashp S^V\cong\Sigma^\bullet\cat{L}(V,-)_+$ given in degree $U\in\cat{O}$ by sending $(i\colon V\to U,u\in\im(i)^\perp)\smashp x$ to $(i(x)+u)\smashp i$, which then induces an isomorphism $\cat{O}(V,-)\smashp S^{V\oplus W}\cong \Sigma^\bullet(\cat{L}(V,-)_+\smashp S^W)$. We omit the straightforward but lengthy computation this is natural and compatible with the symmetric monoidal structures and hence identifies $\cat{f}\circ\cat{y}_{/i}$ with the description of $\Sigma^\bullet\circ\cat{th}|_{(\cat{y}_{/\cat{Gr}})}$ from Proposition~\ref{prop:th-on-corep-sym-mon}.
    \end{proof}
\end{prop}

\begin{cor}\label{cor:thom-space-vs-thom-spectra} 
    We have a commutative square
    \[
        \begin{tikzcd}[anchor=south,baseline=8.4pt]
            \ul\myS_\gl[\cat{Gr}]^\otimes\arrow[d,"{\ul\myS_\gl^\otimes[i]}"']\arrow[r,"\ul\myth_\gl^\otimes"] & \ul\myS_{\gl,*}^\otimes\arrow[d,"\Sigma^\bullet"]\\
            \ul\myS_\gl[\cat{BOP}]^\otimes\arrow[r,"\ul\myTh_\gl^\otimes"'] & \ul\mySp_\gl^\otimes\rlap.
        \end{tikzcd}\qednow
    \]
\end{cor}

\subsection{Rigidity of the global Grassmannian} The two global Thom spectrum constructions expect inputs from two different global $\infty$-categories, so as part of our comparison we in particular have to produce an equivalence between the latter. If there were several such equivalences, we would of course be under pressure to justify why our equivalence is the right choice. In order to settle this question once and for all, we will establish a uniqueness property for our comparison; this subsection is devoted to proving the following key ingredient for this:

\begin{thm}\label{thm:Aut-Vect-oplus}
    The identity is the unique automorphism of the symmetric monoidal global $\infty$-category $\ul\Vect^\oplus$ which sends the object $\R\in\ul\Vect(1)$ to itself; more precisely, the fiber of $\ev_{\R}\colon\Aut(\ul\Vect^\oplus)\to\ul\Vect(1)$ over the point $\R$ is contractible.
\end{thm}
\subsubsection{Homomorphisms into orthogonal groups}\label{subsubsec:reptheo} Because of the equivalence $\ul\Vect\simeq\coprod_{n\ge0}\BGcat{\O(n)}$, (non-monoidal) endomorphisms of $\ul\Vect$ can be described in terms of conjugacy classes of homomorphisms between orthogonal groups. We therefore begin with some elementary representation-theoretic considerations.

\begin{lemma}\label{lemma:conjugate-to-On}
    Let $G$ be a compact Lie group, let $n\ge0$, and let $f_1,f_2\colon G\rightrightarrows \O(n)$ be homomorphisms. Then $f_1$ and $f_2$ are conjugate (i.e.\ homotopic in $\Glo$) if and only if for every $g\in G$ the two matrices $f_1(g)$ and $f_2(g)$ are conjugate in $\GL(n,\C)$.
    \begin{proof}
        The condition is clearly necessary. To see that it is also sufficient, we first observe that $f_1$ and $f_2$ are conjugate if and only if they are conjugate as homomorphisms to $\GL(n,\R)$ \cite{broecker-tom-dieck}*{Exercise II.6.10.5}, i.e.\ if and only if the two representations $f_1^*\R^n$ and $f_2^*\R^n$ are isomorphic \emph{without} regarding the inner product. We will show that these representations become isomorphic after scalar extension to $\C$; forgetting the complex structure will then in particular give an isomorphism $f_1^*\R^n\oplus f_1^*\R^n\cong f_2^*\R^n\oplus f_2^*\R^n$ and counting the multiplicities of irreducible representations on each side yields the claim.

        To prove that $f_1^*\C^n$ and $f_2^*\C^n$ are isomorphic, it will suffice to show that their characters $\chi_i\colon g\mapsto\text{Tr}(f_i(g))$ agree \cite{broecker-tom-dieck}*{Theorem II.4.12}. This follows at once from cyclic invariance of the trace and the assumption that $f_1(g)$ be conjugate to $f_2(g)$ for each fixed $g\in G$.
    \end{proof}
\end{lemma}

\begin{cor}\label{cor:uniqueness-maps-to-O(n)}
    Let $f_1,f_2\colon \O(m)\to \O(n)$ be homomorphisms. Assume that for every block diagonal matrix
    \begin{equation}\label{eq:orth-normal-form}
        A=\begin{pmatrix}
            A_1\\
            &\ddots\\
            &&A_r
        \end{pmatrix}
    \end{equation}
    where $A_i\in \O(m_i)$, $1\le m_i\le 2$, $m_1+\cdots+m_r=m$, the two matrices $f_1(A)$ and $f_2(A)$ are conjugate in $\GL(n,\C)$. Then $f_1$ and $f_2$ are homotopic in $\Glo$.
    \begin{proof}
        By elementary linear algebra, every orthogonal matrix is conjugate to a block matrix of the form $(\ref{eq:orth-normal-form})$. As any homomorphism sends conjugate elements to conjugate elements, the claim then follows from the previous lemma.
    \end{proof}
\end{cor}

\begin{lemma}
    Any automorphism of $\O(2)$ is an inner automorphism.
    \begin{proof}
        As $\SO(2)$ is the identity component of $\O(2)$, $f$ restricts to an automorphism of $\SO(2)$. By \cite{broecker-tom-dieck}*{Proposition~II.8.1}, we therefore need to have $f(A)=A$ or $f(A)=A^{-1}$ for all $A\in\SO(2)$. As $A$ and $A^{-1}$ become conjugate in $\O(2)$ via the reflection $R\coloneqq\begin{psmallmatrix}-1&0\\0&1\end{psmallmatrix}$, we see that $f\colon\O(2)\to\O(2)$ preserves conjugacy classes of elements in $\SO(2)$. 
        
        On the other hand, any matrix in $\O(2)\smallsetminus\SO(2)$ has $-1$ and $1$ as eigenvalues and hence is conjugate to $R$. As $f$ also has to preserve $\O(2)\smallsetminus\SO(2)$, we conclude that $f$ preserves also the conjugacy class of elements in $\O(2)\smallsetminus\SO(2)$, so that it is conjugate to the identity by Lemma~\ref{lemma:conjugate-to-On}.
    \end{proof}
\end{lemma}

Using this, we can now prove the following incoherent precursor to Theorem~\ref{thm:Aut-Vect-oplus}:

\begin{prop}\label{prop:Aut-h-Vect}
    The only automorphism of $\ul\Vect^\oplus$ in $\CMon(\textup{Ho}(\Spc_\Glo))$ is the identity.
    \begin{proof}
        If $f$ is any such automorphism, then $f(\R)\cong\R$ since this is the only object of $\ul\Vect^\oplus(1)$ generating it under the symmetric monoidal product. As $f$ is symmetric monoidal up to incoherent equivalence, we conclude that $f(\R^n)\cong\R^n$, and hence $\dim f(V)= \dim V$ for any representation $V$ by naturality. Thus, $f$ preserves the decomposition $\ul\Vect\simeq\coprod_{n\ge0}\BGcat{\O(n)}$, and we can therefore identify it with a family of conjugacy classes of isomorphisms $f_n\colon\O(n)\iso\O(n)$ such that the diagram
        \[
            \begin{tikzcd}
                \O(m)\times\O(n)\arrow[r,"f_m\times f_n"]\arrow[d,hook] &[1em] \O(m)\times\O(n)\arrow[d,hook]\\
                \O(m+n)\arrow[r,"f_{m+n}"'] & \O(m+n)
            \end{tikzcd}
        \]
        commutes up to conjugation for every $m,n\ge0$. We want to show that all $f_n$ are inner automorphisms, which is trivial for $n=0,1$, while for $n=2$ this is an instance of the previous lemma. For $n\ge3$, we then see by applying the compatibility condition inductively that $f_n$ sends any block matrix of the form $(\ref{eq:orth-normal-form})$ to an element conjugate to itself, so it is an inner automorphism by Corollary~\ref{cor:uniqueness-maps-to-O(n)}.
    \end{proof}
\end{prop}

We now turn our attention to describing the set of homotopies between suitable homomorphisms to $\O(n)$, which amounts to understanding the centralizers of certain subgroups of $\O(n)$.

\begin{convention}
    Throughout, we identify the wreath product $\Sigma_n\wr C_2$ with a subgroup of $\O(n)$ via the homomorphism that sends a permutation to the corresponding permutation matrix and that sends an element of $C_2^n=\{\pm1\}^n$ to the corresponding diagonal matrix.

    For $n_1,\dots,n_r\ge0$ we will abbreviate $\O(n_1|\cdots|n_r)\coloneqq\O(n_1)\times\cdots\times\O(n_r)$ and $\Sig(n_1|\cdots|n_r)\coloneqq(\Sigma_{n_1}\wr C_2)\times\cdots\times(\Sigma_{n_r}\wr C_2)$. We will view $\O(n_1|\cdots|n_r)$ as a subgroup of $\O(n_1+\cdots+n_r)$ via block sum, and we will similarly view $\Sig(n_1|\cdots|n_r)$ as a subgroup of $\Sig(n_1+\cdots+n_r)=\Sigma_{n_1+\cdots+n_r}\wr C_2$. Finally, we will view $\Sig(n_1|\cdots|n_r)$ as a subgroup of $\O(n_1|\cdots|n_r)$ via the product of the above embeddings.
\end{convention}

\begin{lemma}\label{lemma:centralizer-block-permutations}
    Let $n_1,\dots,n_r\ge0$ and write $n\coloneqq n_1+\cdots+n_r$. Then the centralizer of $\Sig(n_1|\cdots|n_r)$ in $\O(n)$ consists precisely of the matrices of the form
    \begin{equation}\label{eq:block-centralizer}
        \begin{pmatrix}
            \pm\id_{\R^{n_1}}\\
            &\ddots\\
            &&\pm\id_{\R^{n_r}}
        \end{pmatrix}.
    \end{equation}    
    \begin{proof}
        It is clear that the matrices of the above form are contained in the centralizer. To prove the converse, we view $\R^n$ as a representation of $G\coloneqq \Sig(n_1|\cdots|n_r)$, so that any element $A$ of the centralizer defines an isometric $G$-equivariant automorphism of $\R^n$. For every $1\le k\le r$ and every $n_1+\cdots+n_{k-1}<i\le n_1+\cdots+n_k$, the standard basis vector $e_i$ is fixed by the subgroup 
        \[
            G_{\lnot k}\coloneqq(\Sigma_{n_1}\wr C_2)\times\cdots(\Sigma_{n_{k-1}}\wr C_2)\times1\times(\Sigma_{n_{k+1}}\wr C_2)\times\cdots\times(\Sigma_{n_r}\wr C_2).
        \]
        Conversely, if $x\in\R$ is not contained in the span of these $e_i$'s, then it is not fixed by $\big((\id;-1,\dots,-1),(\id;-1,\dots,-1),\dots\big)\in G_{\lnot k}$. Thus, the $G_{\lnot k}$-fixed points are precisely given by the span $V_k$ of the $e_i$'s for $n_1+\cdots+n_{k-1}<i\le n_1+\cdots+n_k$; in particular, $A$ has to restrict to an automorphism of $V_k$. As $\R^n=\bigoplus_{k=1}^r V_k$, this expresses $A$ as an orthogonal direct sum of automorphisms of $V_k$ (centralizing $\Sigma_{n_k}\wr C_2$), and we may therefore restrict to the case $r=1$.
        
        For this, let us identify $A$ with an element of $\text{U}(n)$. Then $A$ has a complex eigenvalue $\lambda$ with corresponding eigenvector $v=(v_1,\dots,v_n)\not=0$. As $A$ centralizes $\Sigma_n\wr C_2$, we have $ABv=BAv=B\lambda v=\lambda Bv$ for any $B\in\Sigma_n\wr C_2$, i.e. $Bv$ is again an eigenvector with eigenvalue $\lambda$. If we fix $1\le i\le n$ such that $v_i\not=0$ and let $B\in\Sigma_n\wr C_2$ be the matrix sending the standard basis vector $e_j$ to $-e_j$ for $j\not=i$ and sending $e_i\mapsto e_i$, this implies that $e_i=(2v_i)^{-1}\cdot(v+Bv)$ is an eigenvector with eigenvalue $\lambda$. Similarly, for any $1\le j\le n$ we find a $B\in\Sigma_n\wr C_2$ with $e_i\mapsto e_j$, so that every standard basis vector is an eigenvector with eigenvalue $\lambda$, i.e.~$A=\diag(\lambda)$. As $A$ is an orthogonal matrix, we then have $\lambda=\pm1$, finishing the proof.
    \end{proof}
\end{lemma}

\begin{cor}\label{cor:so-many-centralizers}
    Let $n_1,\dots,n_r\ge0$ and $n\coloneqq n_1+\cdots+n_r$. Then the inclusions
    \[
        C_{\Sig(n)}\Sig(n_1|\cdots|n_r)\hookrightarrow
        C_{\O(n)}\Sig(n_1|\cdots|n_r)\hookleftarrow
        C_{\O(n)}\O(n_1|\cdots|n_r)
    \]
    of centralizers are actually identities.
    \begin{proof}
        In light of the previous lemma, it suffices to observe for the first inclusion that each element of the form $(\ref{eq:block-centralizer})$ belongs to $\Sig(n)$, while it suffices for the second inclusion to note that each such element centralizes all of $\O(n_1|\cdots|n_r)$.
    \end{proof}
\end{cor}

The techniques used in the proof of Proposition~\ref{prop:Aut-h-Vect} are of course woefully inadequate when one wants to prove a homotopy coherent statement. On the other hand, the definition of $\ul\Vect^\oplus$ involves the topology of orthogonal groups, and without a purely categorical definition of it at hand, it is not clear how we should get control over its automorphisms by abstract nonsense.

As our solution to this, we will approximate $\ul\Vect^\oplus$ by a categorically defined gadget $\ul\Vect_\delta^\amalg$, serving as a discrete replacement of real representations. To motivate this, note that as a first-order approximation, the difference between an automorphism $f$ of $\ul\Vect^\oplus$ in $\CMon(\text{Ho}(\Spc_\Glo))$ and one in $\CMon(\Spc_\Glo)$ is that in the former we do not record the choice of the homotopies filling the squares
\begin{equation}\label{diag:the-sad-tale-of-the-forgotten-2-cell}
    \begin{tikzcd}
        \BGcat{\O(n_1|\cdots|n_r)}\arrow[d]\arrow[r,"f_1\times\cdots\times f_r"] &[2em] \BGcat{\O(n_1|\cdots|n_r)}\arrow[d]\\
        \BGcat{\O(n)}\arrow[r,"f_n"'] & \BGcat{\O(n)}.
    \end{tikzcd}
\end{equation}
As the top map is assumed to be an automorphism, the space of such homotopies is a torsor over the automorphisms of the right-hand vertical map, i.e.~over the centralizer of ${\O(n_1|\cdots|n_r)}$ in ${\O(n)}$. Corollary~\ref{cor:so-many-centralizers} says that this agrees with the centralizer of $\Sig(n_1|\cdots|n_r)$ in either $\Sig(n)=\Sigma_n\wr C_2$ or $\O(n)$. Assuming we can make $\coprod_{n\ge0}\BGcat{\Sig(n)}$ into a symmetric monoidal global $\infty$-category and the inclusion to $\coprod_{n\ge0}\BGcat{\O(n)}=\ul\Vect^\oplus$ into a symmetric monoidal global functor,  to remember the above 2-cells it is therefore enough to know the restriction of $f$ to $\coprod_{n\ge0}\BGcat{\Sig(n)}$. Similarly, the aforementioned corollary also tells us that the 2-cell $(\ref{diag:the-sad-tale-of-the-forgotten-2-cell})$ actually comes from $\BGcat{\Sig(n)}$. One might therefore try to understand automorphisms of $\ul\Vect^\oplus$ in terms of automorphisms of (a symmetric monoidal refinement of) $\coprod_{n\ge0}\BGcat{\Sig(n)}$. This is exactly what we will do now: in §\ref{subsubsec:Vect-delta} we will introduce this refinement and compute its space of automorphisms, and in §\ref{subsubsec:comparison} we will show that this space is equivalent to $\Aut(\ul\Vect^\oplus)$, which will then easily yield the proof of Theorem~\ref{thm:Aut-Vect-oplus}.

\medskip
\subsubsection{Discretizing real representations}\label{subsubsec:Vect-delta}
Just like the symmetric monoidal refinement $\ul\Vect^\oplus$ of $\coprod_{n\ge0} \BGcat{\O(n)}$, our refinement $\ul\Vect^\amalg_\delta$ of $\coprod_{n\ge0}\BGcat{\Sig(n)}$ will arise via the continuous Borel construction.

\begin{constr}
    We write $\Ff_{C_2}^\amalg$ for the symmetric monoidal $1$-groupoid of finite free $C_2$-sets under disjoint union, and we write $\ul\Vect_\delta^{\amalg}\coloneqq\smash{\Ntop(\Ff_{C_2}^{\amalg,\flat})}$ for the corresponding continuous Borel category. More concretely, this is the symmetric monoidal global 1-groupoid sending a compact Lie group $G$ to the groupoid of finite $(C_2\times\pi_0(G))$-sets such that the $C_2$-action is free, with the evident restrictions.
\end{constr}

\begin{constr}
    Define a functor from the category of finite $C_2$-sets to the category of finite-dimensional inner product spaces and $\R$-linear maps as follows:

    We send a finite free $C_2$-set $X$ to the vector space $|X|$ of $C_2$-equivariant functions $\chi\colon X\to\R$, where the generator of $C_2$ acts on the target via multiplication by $-1$; the inner product is given by $\langle \chi,\xi\rangle\coloneqq{1\over2}\sum_{x\in X} \chi(x)\xi(x)$. For a $C_2$-equivariant map $f\colon X\to Y$, we define $|f|$ by $|f|(\chi)(y)\coloneqq\sum_{f(x)=y}\chi(x)$ for every $\chi\colon X\to\R$.

    One then directly checks that $|{\cdot}|$ preserves coproducts and sends \emph{injective} $C_2$-equivariant maps to isometric embeddings, thus yielding a symmetric monoidal functor $\smash{|{\cdot}|\colon\Ff_{C_2}^\amalg\to\core\cat{L}^\oplus}$. Passing to continuous Borel categories, we obtain a symmetric monoidal global functor $|{\cdot}|\colon\ul\Vect_\delta^\amalg\to\ul\Vect^\oplus$.
\end{constr}

\begin{rk}\label{rk:BSig}
    Note that $\coprod_{n\ge 0}\BGcat{\Sig(n)}\simeq\Ff_{C_2}$ as (discrete) topological categories, with the equivalence given on the $n$-th coproduct summand by the map classifying the object $\{1,\dots,n\}\times C_2$ with its evident action of $\Sig(n)=(\Sigma_n\wr C_2)$. 

    If $X$ is a finite free $C_2$-set and $\{x_1,\dots,x_n\}$ a system of orbit representatives, then $|X|$ can be identified as an inner product space with $\R^n$ with its standard scalar product via evaluation at $x_1,\dots,x_n$. One then directly checks that as a non-symmetric monoidal functor $|{\cdot}|$ agrees under the chosen equivalences with the map $\coprod_{n\ge0}\BGcat{\Sig(n)}\to\coprod_{n\ge0}\BGcat \O(n)$ induced by our usual embeddings $\Sig(n)\hookrightarrow \O(n)$. Similarly, one computes that the monoidal product on $\coprod_{n\ge0}\BGcat{\Sig(n)}$ is induced by the block sum embeddings $\Sig(m|n)\to\Sig(m+n)$.
\end{rk}

\begin{prop}\label{prop:Vect-delta-RKE}   
    Evaluation at $1\in\Glo$ defines an equivalence
    \[
        \End_{\CMon(\Spc_\Glo)}(\ul\Vect^\amalg_\delta)\iso\End_{\CMon(\Spc)}(\Ff_{C_2}^\amalg).
    \]
    \begin{proof}
        We factor the evaluation functor $\CMon(\Spc_\Glo)\to\CMon(\Spc)$ through the $\infty$-category $\CMon(\Spc_{\Glo_{\Ff in}})$ of commutative monoids in $\Fin$-global spaces, and we will show that each of these individual restrictions defines an equivalence on endomorphism spaces.

        For the first restriction we observe that $i\colon \Glo_{\Ff in}\hookrightarrow\Glo$ admits a left adjoint $\text{Ho}\colon\hskip 0pt minus 1pt\Glo\hskip0pt minus 1pt\to\hskip 0pt minus 1pt\Glo_{\Ff in}$ induced by the enriched functor $\BGcat{G}\mapsto\BGcat\pi_0(G)$. As $i$ is fully faithful, $\text{Ho}$ is a (Bousfield) localization at those group homomorphisms that induce bijections of path components. It follows that $\text{Ho}^*\coloneqq\Fun(\text{Ho},\CMon(\Spc))$ is fully faithful with essential image those functors that invert these maps in $\Glo$. This condition is evidently satisfied for $\ul\Vect^\amalg_\delta$; as $i^*$ is left inverse to $\text{Ho}^*$ this shows that it induces an equivalence
        \[
            \End_{\CMon(\Spc_\Glo)}(\ul\Vect^\amalg_\delta)\iso\End_{\CMon(\Spc_{\Glo_{\Ff in}})}(i^*\ul\Vect^\amalg_\delta).
        \]

        To show that also restriction along $\{1\}\hookrightarrow\Glo_{\Ff in}$ induces an equivalence of endomorphism spaces, it will suffice that $i^*\ul\Vect^\amalg_\delta\colon\Glo_{\Ff in}^\op\to\CMon(\Spc)$ is right Kan extended from $\{1\}$. As $\CMon(\Spc)\to\Spc$ is a conservative right adjoint, it suffices to show that the underlying $\Fin$-global space $i^*\ul\Vect_\delta\colon\Glo_{\Ff in}^\op\to\Spc$ is right Kan extended from the point. Plugging in the definitions, this can be written as
        \[
            \Glo_{\Ff in}^\op\hookrightarrow\Cat_\infty\xrightarrow{\;\Fun(-,\Ff_{C_2})\;}\Cat_\infty,
        \]
        which is indeed right Kan extended by \cite{swan}*{Example~2.8}.
    \end{proof}
\end{prop}

\begin{lemma}\label{lemma:FC2-free}
    Let $\Cc^\otimes$ be any symmetric monoidal $\infty$-category. Then evaluation at the free $C_2$-set $C_2$ defines an equivalence $\Fun^\otimes(\Ff_{C_2}^\amalg,\Cc^\otimes)\iso\Fun(BC_2,\Cc)$.
    \begin{proof}
        Write $\mathbb P\colon\Cat_\infty\to\CMon(\Cat_\infty)$ for the left adjoint to the forgetful functor. Then $\smash{\mathbb P(\Xx)\simeq\coprod_{n\ge 0}(\Xx^{\times n})_{h\Sigma_n}}$ is a $1$-groupoid whenever $\Xx$ is so; it will therefore suffice to prove the universal property of $\smash{\Ff_{C_2}^\amalg}$ with respect to symmetric monoidal 1-categories, which is an instance of \cite{global-mackey}*{Lemma~2.21}.
    \end{proof}
\end{lemma}

\begin{prop}\label{prop:Aut-Vect-delta}
    The space $\Aut(\ul\Vect^\amalg_\delta)$ is an Eilenberg--MacLane space of type $K(\Z/2,1)$. The unique non-trivial automorphism of the identity is given by acting with the generator of $C_2$.
    \begin{proof}
        For the first half, note that in the diagram
        \begin{equation}\label{eq:Aut-Vect-delta-vs-Aut-BC2}
            \Aut(\ul\Vect^\amalg_\delta)\xrightarrow{\;\ev\;}\End(\Ff_{C_2}^\amalg)\xrightarrow{\;\text{res}\;}\hom(BC_2,\Ff_{C_2})
        \end{equation}
        the first map is a monomorphism (i.e., an inclusion of path components) by Proposition~\ref{prop:Vect-delta-RKE}, while the second one is so by Lemma~\ref{lemma:FC2-free}. We now observe that any automorphism of $\Ff_{C_2}$ has to send the object $C_2$ to itself (as this is the only object with precisely two endomorphisms), so that $(\ref{eq:Aut-Vect-delta-vs-Aut-BC2})$ factors through $\Aut(BC_2)$. The latter is connected (and $\Aut(\ul\Vect^\amalg_\delta)$ is non-empty), so the inclusion of path components $\Aut(\ul\Vect^\amalg_\delta)\to\Aut(BC_2)$ has to be an equivalence. It is clear that $\Aut(BC_2)$ is a $K(\Z/2,1)$, with the unique non-trivial automorphism of the identity of $BC_2$ given by the natural transformation $\tau$ corresponding to the generator of $C_2$, and that $(\ref{eq:Aut-Vect-delta-vs-Aut-BC2})$ sends the natural transformation given by acting with the generator to $\tau$, which then completes the proof of the proposition. 
    \end{proof}
\end{prop}

\subsubsection{Automorphisms of $\ul\Vect^\oplus$}\label{subsubsec:comparison}
As promised, we will now prove the following comparison between automorphisms of $\ul\Vect^\oplus$ and $\ul\Vect^\amalg_\delta$:

\begin{prop}\label{prop:comparing-Aut}
    The maps
    \begin{equation}\label{eq:discretivatization}
        \Aut(\ul\Vect^\oplus)\to\hom(\ul\Vect^\amalg_\delta,\ul\Vect^\oplus)\gets
        \Aut(\ul\Vect^\amalg_\delta)
    \end{equation}
    induced by $|{\cdot}|$ are monomorphisms, and they have the same image.
    \begin{proof}
        View both $\ul\Vect^\oplus$ and $\ul\Vect^\amalg_\delta$ as functors $\text{Fin}_*\to\mathbb\Spc_{\Glo}$ satisfying the Segal condition. Throughout the proof, we will make use of the usual formula \cite{GHN17}*{Proposition~5.1} for mapping spaces in a functor category as the end
        \begin{multline*}
            \hom(F,G)\simeq\int_{\langle r\rangle\in\text{Fin}_*}\hom(F(\langle r\rangle),G(\langle r\rangle))\\{}\coloneqq
            \lim\big({\Tw(\text{Fin}_*)\xrightarrow{\;\ev\;} \text{Fin}_*\!\!\!{}^\op\times\text{Fin}_*\xrightarrow{\hom(F(-),G(-))\;}\Spc}\big),
        \end{multline*}
        where $\Tw$ denotes the twisted arrow category.

        Consider the subfunctor $F_\oplus\colon\Tw(\text{Fin}_*)\to\Spc$ of $\hom(\ul\Vect^\oplus(-),\ul\Vect^\oplus(-))$ given at $\alpha\colon\langle r\rangle \to\langle s\rangle$ by the component of $\ul\Vect^\oplus(\alpha)\colon\ul\Vect^\oplus(\langle r\rangle)\to\ul\Vect^\oplus(\langle s\rangle)$. Then the map $\lim F_\oplus\to\lim\hom\big(\Vect^\oplus(-),\Vect^\oplus(-)\big)\simeq\End(\ul\Vect^\oplus)$ induced by the inclusion is a monomorphism (as a limit of monomorphisms), and its image consists precisely of those natural transformations whose underlying functor is the identity. By Proposition~\ref{prop:Aut-h-Vect}, this image agrees with $\Aut(\ul\Vect^\oplus)$. In the same way, Proposition~\ref{prop:Aut-Vect-delta} identifies $\Aut(\ul\Vect^\amalg_\delta)$ with the limit of the analogously defined subfunctor $F_\delta\subset\hom(\ul\Vect^\amalg_\delta(-),\ul\Vect^\amalg_\delta(-))$, reducing us to proving that the maps
        \[
            \lim_{\alpha\colon\langle r\rangle\to\langle s\rangle}\hskip0pt minus 15ptF_\oplus(\alpha)\to
            \lim_{\alpha\colon\langle r\rangle\to\langle s\rangle}\hskip0pt minus 15pt\hom(\ul\Vect^\amalg_\delta(\langle r\rangle),\ul\Vect^\oplus(\langle s\rangle))\gets
            \lim_{\alpha\colon\langle r\rangle\to\langle s\rangle}\hskip0pt minus 15pt F_\delta(\alpha)
        \]
        induced by $|{\cdot}|$ are monomorphisms with the same image. Using once more that limits preserve monomorphisms, it will then suffice to show that for each $\alpha$ the maps
        \begin{equation}\label{eq:monos-without-ends}
            F_\oplus(\alpha)\to
            \hom\big(\ul\Vect^\amalg_\delta(\langle r\rangle),
            \ul\Vect^\oplus(\langle s\rangle)\big)\gets
            F_\delta(\alpha)
        \end{equation}
        are monomorphisms with the same image. 

        By construction, the image of the first map is the component of $\ul\Vect^\oplus(\alpha)\circ|{\cdot}|$, while the image of the second one is the component of $|{\cdot}|\circ\ul\Vect^\amalg_\delta(\alpha)$, and these maps are indeed homotopic. To prove that the maps $(\ref{eq:monos-without-ends})$ are monomorphisms, we first note that the Segal conditions for $\ul\Vect^\amalg_\delta$ and $\ul\Vect^\oplus$ imply that we have equivalences
        \begin{align*}
            F_\oplus(\alpha)&\iso\prod_{i=1}^s F_\oplus(\chi_i\circ\alpha)\\
            F_\delta(\alpha)&\iso\prod_{i=1}^s F_\delta(\chi_i\circ\alpha)\\
            \hom\big(\ul\Vect^\amalg_\delta(\langle r\rangle), \ul\Vect^\amalg_\delta(\langle s\rangle)\big)&\iso
            \prod_{i=1}^s\hom\big(\ul\Vect^\amalg_\delta(\langle r\rangle), \ul\Vect^\amalg_\delta(\langle 1\rangle)\big)
        \end{align*}
        induced by the maps $\chi_i\colon\langle s\rangle\to\langle1\rangle$ with $\chi_i^{-1}(1)=\{i\}$. As monomorphisms are closed under products, we may therefore assume that $s=1$. Using moreover that the class of all $\alpha\colon\langle r\rangle\to\langle1\rangle$ for which the above maps are monomorphisms is closed under equivalences in $\Tw(\text{Fin}_*)$, we may then further assume for (notational) simplicity that $\alpha$ is order preserving. We set $a\coloneqq\min(\alpha^{-1}(1)\cup\{r+1\})$.
        
        Now it is time to plug in the definitions. We may identify $\ul\Vect^\oplus(\alpha)$ with the map $\coprod_{m_1,\dots,m_r}\BGcat{\O(m_1|\cdots|m_r)}\to\coprod_{n}\BGcat{\O(n)}$ given on the summand $(m_1,\dots,m_r)$ by the composite 
        \[
            \BGcat{\O(m_1|\cdots|m_r)}\xrightarrow{\;\pr\;}
            \BGcat{\O(m_a|\cdots|m_r)}\lhook\joinrel\longrightarrow
            \BGcat{\O(m_a+\cdots+m_r)}.
        \]
        Writing $\hom^*\big(\BGcat{\O(m_1|\dots|m_r)},\BGcat{\O(m_a+\cdots+m_r)}\big)$ for the component of this map, we therefore have
        \[
            F_\oplus(\alpha)\simeq\prod_{m_1,\dots,m_r}\hom^*\big(\BGcat{\O(m_1|\dots|m_r)},\BGcat{\O(m_a+\cdots+m_r)}\big).
        \]
        Analogously, $F_\delta(\alpha)\simeq\prod_{m_1,\dots,m_r}\hom^*\big(\BGcat{\Sig(m_1|\cdots|m_r)},\BGcat{\Sig(m_a+\cdots+m_r)}\big)$, while we may identify the relevant path component of $\hom\big(\ul\Vect^\amalg_\delta(\langle r\rangle),\ul\Vect^\oplus(\langle 1\rangle)\big)$ with $\prod_{m_1,\dots,m_r}\hom^*\big(\BGcat{\Sig(m_1|\cdots|m_r)},\BGcat{\O(m_a+\cdots+m_r)}\big)$, where $\hom^*$ in each case denotes the component of the composite of the projection to $\BGcat{\Sig(m_a|\cdots|m_r)}$ followed by the evident embedding. In this description, the maps from $(\ref{eq:monos-without-ends})$ are given by postcomposition with the usual embedding $\BGcat{\Sig(m_a+\cdots+m_r)}\to\BGcat{\O(m_a+\cdots+m_r)}$ and precomposition with $\BGcat{\Sig(m_1|\cdots|m_r)}\to\BGcat{\O(m_1|\cdots|m_r)}$, respectively. Plugging in the explicit description of mapping spaces in $\Glo$, we have therefore altogether reduced to showing that the maps
        \begin{align*}
            B\big(C_{\Sig(m_a+\cdots+m_r)}\Sig(m_a|\cdots|m_r)\big)&\to
            B\big(C_{\O(m_a+\cdots+m_r)}\Sig(m_a|\cdots|m_r)\big)\\
            B\big(C_{\O(m_a+\cdots+m_r)}\O(m_a|\cdots|m_r)\big)&\to
            B\big(C_{\O(m_a+\cdots+m_r)}\Sig(m_a|\cdots|m_r)\big)
        \end{align*}
        of classifying spaces of centralizers induced by pre- and postcomposition are equivalences. This is precisely the content of Corollary~\ref{cor:so-many-centralizers}.
    \end{proof}
\end{prop}

\begin{proof}[Proof of Theorem~\ref{thm:Aut-Vect-oplus}]
    Remember that we want to show that the fiber of $\ev_{\R}\colon\Aut_{\CMon(\Spc_{\Glo})}(\ul\Vect^\oplus)\to\ul\Vect(1)$ over $\R$ is contractible. Replacing the target by the full subcategory $B\Aut(\R)$ spanned by $\R$ does not change the fiber, so we may equivalently show that the evaluation map $\ev_{\R}\colon\Aut(\ul\Vect^\oplus)\to B\Aut(\R)$ is an equivalence. The source is an Eilenberg--MacLane space of type $K(\Z/2,1)$ by Proposition~\ref{prop:comparing-Aut} combined with Proposition~\ref{prop:Aut-Vect-delta}, as is the target since $\Aut(\R)=\{\pm1\}$. Moreover, the automorphism of the identity of $\ul\Vect^\oplus$ given on each representation by scalar multiplication with $-1$ evaluates to the non-trivial automorphism of $\R$, so that the evaluation map is non-zero and hence an equivalence as claimed.
\end{proof}

\subsection{Comparing the inputs} In this section we will prove the promised comparison between the inputs of our two global Thom spectrum constructions:

\begin{thm}\label{thm:comparing-inputs}
    \begin{enumerate}
        \item There exists a unique equivalence 
        \[
            \psi^\otimes\colon\ul\Spc_{\Glo/\ul\Vect^\oplus}^\otimes\iso
            \ul\myS_\gl[\cat{Gr}]^\otimes
        \] of symmetric monoidal global $\infty$-categories 
        sending the object $y(\R)=(\R\colon 1\to\ul\Vect)$ of $\Spc_{\Glo/\ul\Vect}$ to the map $i_{\R}\colon\cat{L}(\R,-)\to\cat{Gr}$ classifying the unique point $\R\in\cat{Gr}_1(\R)$. 
        \item There exists a unique equivalence $\Psi^\otimes$ of symmetric monoidal global $\infty$-categories completing the commutative square
        \begin{equation}\label{diag:defining-Psi-otimes}
            \begin{tikzcd}
                \ul\Spc_{\Glo/\ul\Vect^\oplus}^\otimes\arrow[r,"\psi^\otimes","\sim"']\arrow[d] & \ul\myS_{\gl}[\cat{Gr}]^\otimes\arrow[d,"{\ul\myS_\gl[i]}"]\\
                \ul\Spc_{\Glo/\ul\VRep^{\oplus}}^\otimes\arrow[r,dashed,"\Psi^\otimes","\sim"'] & \ul\myS_\gl[\cat{BOP}]^\otimes\rlap,
            \end{tikzcd}
        \end{equation}
        where the left-hand vertical map is induced via the universal property of symmetric monoidal cocompletion by the chosen group completion.
    \end{enumerate}
\end{thm}

Before we give the formal proof, let us explain the main idea. By Corollary~\ref{cor:S-gl-times-initial} there exists a unique symmetric monoidal equivalence $\Phi^\otimes\colon\ul\myS_\gl^{\otimes}\iso\smash{\ul\Spc_{\Glo}^\times}$, inducing for every fixed orthogonal space $X$ an equivalence
\begin{equation}\label{eq:slice-comparison}
    \ul\myS_\gl[X]\xrightarrow[\raise3.5pt\hbox{$\scriptstyle\smash{\sim}$}]{\,\textup{Prop.~\ref{prop:model-categorical-model}}\,}(\ul\myS_{\gl})_{/X} \xrightarrow[\raise3.5pt\hbox{$\scriptstyle\smash{\sim}$}]{\,\;\Phi\;\,}\ul\Spc_{\Glo/\Phi(X)}.
\end{equation}
If we could upgrade these maps for varying ultra-commutative $X$ into a {natural} transformation of {symmetric monoidal} global $\infty$-categories, this would reduce the proof of existence of $\psi^\otimes$ to identifying $\Phi^\otimes(\cat{Gr})$ as a commutative monoid in $\Spc_{\Glo}$ with $\ul\Vect^\oplus$, while the second part would similarly reduce to proving that $\Phi^\otimes(\cat{Gr})\to\Phi^\otimes(\cat{BOP})$ is a group completion. While this is the general idea the reader should keep in mind, both establishing naturality of $(\ref{eq:slice-comparison})$ as well as identifying $\Phi^\otimes(\cat{Gr})$ as a commutative monoid are rather non-trivial, and we will therefore use a couple of tricks to get away with significantly less: we will only need to explicitly identify $\Phi(\cat{Gr})$ as a global space, and we will basically only establish a natural version of $(\ref{eq:slice-comparison})$ after evaluating at any fixed compact Lie group $G$.
 
\medskip
\subsubsection{The non-group-completed comparison} Naturally, we want to construct $\psi^\otimes$ using the universal property of $\ul\Spc_{\Glo/\ul\Vect^\oplus}^\otimes$ as symmetric monoidal cocompletion, so we should first describe a suitable functor $\smash{\ul\Vect^\oplus\to\ul\myS_\gl[\cat{Gr}]^\otimes}$.

\begin{constr}[cf.~\cite{juran-thesis}*{Construction~5.5}]
    Consider the enriched functor $\cat{y}\colon\core\cat{L}\to\cat{$\cat{L}$-Top}_{/\cat{Gr}}$ sending an inner product space $V$ to the map $i_V\colon\cat{L}(V,-)\to\cat{Gr}$ classifying the point $V\in\cat{Gr}(V)$, and sending a linear isometric isomorphism $f\colon V\to W$ to the commutative triangle
    \[
        \begin{tikzcd}[column sep=small]
            \cat{L}(V,-)\arrow[rr,"{\cat{L}(f^{-1},-)}"]\arrow[dr, bend right=15pt,"i_V"']&&\cat{L}(W,-)\rlap.\arrow[dl,bend left=15pt,"i_W"]\\
            & \cat{Gr}
        \end{tikzcd}
    \]
    This upgrades to a (enriched) symmetric monoidal functor $\cat{y}^\otimes$ using the symmetric monoidal structure on the enriched Yoneda embedding $\cat{L}^\op\to\cat{$\cat{L}$-Top}$. We will write $y_\top^\otimes$ for the resulting functor symmetric monoidal global functor 
    \begin{equation}\label{eq:ytop}
        \ul\Vect^\oplus=\Ntop(\core\cat{L}^{\oplus,\dual})\to\Ntop(\cat{$\cat{L}$-Top}_{/\cat{Gr}}^{\boxtimes,\dual})\to\ul\myS_\gl[\cat{Gr}]^\otimes.
    \end{equation}
\end{constr}

\begin{prop}\label{prop:construction-of-psi-otimes}
    The symmetric monoidal global functor $(\ref{eq:ytop})$ extends to an equivalence $\psi^\otimes\colon\ul\Spc_{\Glo/\Vect^\oplus}^\otimes\iso\ul\myS_\gl[\cat{Gr}]^\otimes$.
    \begin{proof}Proposition~\ref{prop:model-categorical-model} implies that $\ul\myS_{\gl}[\Gr]$ is globally presentable. 
    The forgetful functor $\ul\myS_{\gl}[\Gr]\to \ul\myS_{\gl}$ enhances naturally to a symmetric monoidal functor and is easily seen to preserve global colimits. As the forgetful functor is conservative and $\ul\myS_{\gl}^\otimes$ is globally presentably symmetric monoidal, it follows that also $\ul\myS_{\gl}[\Gr]^{\otimes}$ is globally presentably symmetric monoidal. 
    We now want to show that the unique symmetric monoidal left adjoint $\psi^\otimes$ extending $y_\top^\otimes$ is an equivalence, which can be checked after forgetting the symmetric monoidal structure. Moreover, we may just as well work with any other model of the global cocompletion, reducing us to constructing some \emph{non-monoidal} equivalence fitting into a commutative triangle as depicted on the left in the following diagram:
        \[
            \begin{tikzcd}[column sep=small,cramped]
                &[-1.4em] \ul\Vect\arrow[dl, bend right=15pt,"y"']\arrow[dr,bend left=15pt,"y_\top"] &[.6em]&[0em]&[-3em] \coprod\limits_{\smash{n\ge0}}\BGcat{\O(n)}\arrow[dl,bend right=15pt,"y"']\arrow[dr,bend left=15pt,"y_\top\circ\tau"]\\
                \ul\Spc_{\Glo}(\ul\Vect\times{-})\arrow[rr,"\psi","\sim"'] && \ul\myS_\gl[\cat{Gr}] &
                \ul\Spc_{\Glo}\big(\kern-1pt{\coprod\limits_{\smash{n\ge0}}\BGcat{\O(n)}\times{-}}\big) \arrow[rr,"\tilde\psi","\sim"'] && \ul\myS_\gl[\cat{Gr}]
            \end{tikzcd}
        \]
        In light of the equivalence $\tau\colon\coprod_{n\ge0}\BGcat{\O(n)}\to\ul\Vect$ given on the $n$-th coproduct summand by the map classifying the tautological $\O(n)$-representation $\R^n$, we may equivalently construct an equivalence $\tilde\psi$ making the triangle on the right commute. 
        
        Note that $y_\top\circ\tau$ is given on the $n$-th summand by the map classifying the projection \[p_n\colon\cat{L}(\R^n,-)\to\cat{L}(\R^n,-)/\O(n)=\triv_{\O(n)}\cat{Gr}_n\hookrightarrow\triv_{\O(n)}\cat{Gr},\] viewed as an object of $\ul\myS_\gl[\cat{Gr}](\BGcat{\O(n)})$. We will now construct an equivalence $\alpha_n\colon\Phi(\cat{Gr}_n)\iso{\BGcat{\O(n)}}$ fitting into a commutative diagram
        \begin{equation}\label{diag:alpha-compatible}
            \begin{tikzcd}[cramped]
                \Phi(\cat{L}(\R^n,-))\arrow[d,"\Phi(p_n)"']\arrow[rr,"\sim"] &[-.67em]&[2em] \id_{\BGcat{\O(n)}}\arrow[d,"\Delta"]\\
                \Phi(\triv_{\O(n)}\cat{Gr}_n)\arrow[r,"\sim"] & \triv_{\O(n)}\Phi(\cat{Gr}_n)\arrow[r,"\triv_{\O(n)}\alpha_n"] & \big(\pr\colon\BGcat{\O(n)}^{\times2}\to\BGcat{\O(n)}\big)\rlap.
            \end{tikzcd}
        \end{equation}
        For this we observe that as the $\O(n)$-orthogonal space $\cat{L}(\R^n,-)$ is cofibrant, $p_n$ exhibits $\cat{Gr}_n$ as a left adjoint object under $\triv_{\O(n)}\colon\myS_\gl\to\myS_\text{$\O(n)$-gl}$. As the diagonal embedding similarly exhibits $\BGcat{\O(n)}$ as a left adjoint object for $\triv_{\O(n)}$, we may simply take $\alpha_n$ to be the corresponding Beck--Chevalley map for the naturality equivalence $\Phi\circ\triv_{\O(n)}\simeq\triv_{\O(n)}\circ\Phi$. 
        
        This completes the construction of the $\alpha_n$'s.  
        Using in addition that $\Phi$ preserves coproducts, we then get an equivalence $\alpha\colon\Phi(\cat{Gr})\to\coprod_{n\ge0}\BGcat{\O(n)}$ restricting to the $\alpha_n$'s. Commutativity of $(\ref{diag:alpha-compatible})$ then says that
        \begin{equation}\label{eq:psi-tilde-inverse}
            \ul\myS_\gl[\cat{Gr}]\iso\ul\Spc_{\Glo}\big(\Phi(\cat{\Gr})\times{-}\big)\xrightarrow[\raise3pt\hbox{$\scriptstyle\sim$}]{\;\alpha_!\;}\ul\Spc_{\Glo}\big(\textstyle\coprod\limits_{\smash{n\ge0}}\BGcat{\O(n)}\times{-}\big)
        \end{equation}
        sends $p_n\in\ul\myS_\gl[\cat{Gr}](\BGcat{\O(n)})$ to the object
        \[
            \BGcat{\O(n)}\xrightarrow{\;\Delta\;}
            \BGcat{\O(n)}\times\BGcat{\O(n)}\hookrightarrow
            \big(\textstyle\coprod\limits_{\smash{n\ge0}}\BGcat{\O(n)}\big)\times\BGcat{\O(n)}.
        \]
        As the latter is the Yoneda image of $\BGcat{\O(n)}\hookrightarrow\coprod_{n\ge0}\BGcat{\O(n)}$ (see Remark~\ref{rk:Yoneda-image}), this shows that the inverse of $(\ref{eq:psi-tilde-inverse})$ is the desired $\tilde{\psi}$.
    \end{proof}
\end{prop}

To establish uniqueness of $\psi^\otimes$ we will need the following basic observation:

\begin{lemma}
    Let $T$ be a small $\infty$-category, and let $X$ be a small $T$-$\infty$-groupoid. Then every automorphism of $\ul\Spc_{T/X}$ sends the Yoneda image to itself.
    \begin{proof}
        Again, we may equivalently use the model of the global completion as $\ul\Spc_{T}(X\times{-})$.
        By global completeness of $\ul\Spc_T$, the restriction $\pi^*\colon\ul\Spc_T\to\ul\Spc_{T}(X\times{-})$ along $\pi\colon X\to 1$ is a left adjoint. Moreover, $\pi^*$ preserves the terminal object. If $F$ is any automorphism of $\ul\Spc_{T}(X\times{-})$, then $F\pi^*$ is again a left adjoint preserving the terminal object and hence is equivalent to $\pi^*$ by the universal property of $\ul\Spc_T$. Passing to mates, we obtain an equivalence $\pi_!F\simeq\pi_!$ of functors $\ul\Spc_{T}(\ul X\times{-})\to\ul\Spc_T$. As the Yoneda image is precisely the preimage of the terminal object under $\pi_!$ (see Remark~\ref{rk:Yoneda-image}), the claim follows.  
    \end{proof}
\end{lemma}

\begin{proof}[Proof of Theorem~\ref{thm:comparing-inputs}(1)]
    As $y_\top(\R)=i_\R$ by definition, the existence of such an equivalence follows from Proposition~\ref{prop:construction-of-psi-otimes}. For uniqueness, we may then equivalently show that there is only one symmetric monoidal automorphism of $\ul\Spc_{\Glo/\ul\Vect^\oplus}^\otimes$ sending $y(\R)$ to itself. By the universal property, restriction along the Yoneda embedding defines an equivalence 
    \[
        \smash{\End_{\CAlg(\PrL_{\Glo})}\big(\ul\Spc_{\Glo/\ul\Vect^\oplus}^\otimes\big)\iso\hom_{\CAlg(\CAT_\Glo)}\big(\ul\Vect^\oplus,\ul\Spc_{\Glo/\ul\Vect^\oplus}^\otimes\big)},
    \]
    and by the previous lemma this restricts to a map $\Aut(\ul\Spc_{\Glo/\ul\Vect^\oplus}^\otimes)\to\Aut(\ul\Vect^\oplus)$, which is then necessarily a summand inclusion. We conclude that the fiber of $\ev_\R\colon\Aut(\ul\Spc_{\Glo/\ul\Vect^\oplus}^\otimes)\to\ul\Spc_{\Glo}(\ul\Vect)$ over $y(\R)$ is a summand of the fiber of $\ev_\R\colon\Aut(\ul\Vect^\oplus)\to\ul\Vect(1)$ over $\R$, and so the claim follows from Theorem~\ref{thm:Aut-Vect-oplus}.
\end{proof}

\subsubsection{And now with group completions}
By the universal property of the group completion $\ul\Vect^\oplus\to\ul\VRep^{\oplus}$ and the universal property of symmetric monoidal cocompletion, the left-hand vertical map $\ul\Spc_{\Glo/\ul\Vect^\oplus}^\otimes\to\ul\Spc_{\Glo/\ul\VRep^{\oplus}}^\otimes$ in $(\ref{diag:defining-Psi-otimes})$ is initial among maps in $\smash{\CAlg(\PrL_\Glo)}$ sending the Yoneda image to invertible objects. Our strategy to produce the equivalence $\Psi^\otimes$ and show that it is unique will thus be to establish the same universal property for the right-hand vertical map $\ul\myS_\gl[i]^\otimes\colon\ul\myS_\gl[\cat{Gr}]^\otimes\to\ul\myS_\gl[\cat{BOP}]^\otimes$.

\begin{definition}
    We call a map $f\colon X\to Y$ of ultra-commutative monoids a \emph{global group completion} if it is sent to a (levelwise) group completion under the functor
    \[
        \CAlg(\cat{$\cat{L}$-Top}^\boxtimes)\xrightarrow{\;\CAlg(\text{loc})\;}\CAlg(\myS_\gl^\otimes)\xrightarrow[\raise3.5pt\hbox{$\smash{\scriptstyle\sim}$}]{\;\CAlg(\Phi^\otimes)\;}\CMon(\Spc_{\Glo})
    \]
\end{definition}

\begin{thm}\label{thm:Gr-oup-completion}
    The map $i\colon\cat{Gr}\to\cat{BOP}$ is a global group completion.
    \begin{proof}
        For $f\colon V\to W$ in $\cat{L}$ and $n\ge0$, the map $\cat{Gr}_n(f)\colon\cat{Gr}_n(V)\to\cat{Gr}_n(W)$ is a closed embedding as it is a continuous injection between compact spaces; in the same way, one shows that ${-}\oplus W\colon\cat{Gr}_n(V)\to\cat{Gr}_{n+\dim(W)}(V\oplus W)$ is a closed embedding. With this established, it follows straight from the definitions that the orthogonal spaces $\cat{Gr}$ and $\cat{BOP}$ are closed. In light of the description of the equivalence $\Phi$ from Lemma~\ref{lemma:id-g-fixed-points-vs-ev}, the theorem therefore amounts to saying that for every complete $G$-universe $\Uu_G$ the induced map $i(\Uu_G)^G\colon\cat{Gr}(\Uu_G)^G\to\cat{BOP}(\Uu_G)^G$ is a group completion in $\CMon(\Spc)$ (where we make $(-)({\Uu_G})^G$ into a symmetric monoidal functor $\myS_\gl^{\text{closed},\otimes}\to\Spc^\times$ in the unique way).

        As the key step in the proof of \cite{schwede2018global}*{Theorem~2.5.33}, it is shown that the induced map $H_*(i(\Uu_G)^G;\Z)$ on integral homology is a localization (in the ring-theoretic sense) at $\pi_0(\cat{Gr}(\Uu_G)^G)\subset H_0(\cat{Gr}(\Uu_G)^G;\Z)$. In particular, $\cat{BOP}(\Uu_G)^G$ is grouplike, and so $i(\Uu_G)^G$ factors through a map $f\colon X\to\cat{BOP}(\Uu_G)^G$ from the group completion. By the McDuff--Segal group completion theorem (see \cite{nikolaus-group-completion}*{Theorem~1} for a proof in this language in generality), also the group completion $\cat{Gr}(\Uu_G)^G\to X$ induces a localization on homology at the same classes, so that $f$ is a homology isomorphism. As both its source and target are grouplike h-spaces (and hence their fundamental groups act trivially on higher homotopy groups), $f$ is therefore already an equivalence by the homological Whitehead theorem.
    \end{proof}
\end{thm}

The desired universal property of $\ul\myS_\gl[i]$ will thus be an instance of the following more general proposition:

\begin{prop}\label{prop:group-complete-on-slices}
    Let $f\colon X\to Y$ be a global group completion. Then the induced map $\ul\myS_\gl[X]^\otimes\to\ul\myS_\gl[Y]^\otimes$ is initial among symmetric monoidal left adjoints inverting every object in the full subcategory $\smash{\ul\myS_\gl^\circ[X]}$ spanned in each degree ${G}$ by those maps $A\to X$ such that the $G$-orthogonal space $A$ is weakly contractible.
\end{prop}

This will require some preparations.

\begin{lemma}
    Let $G$ be any compact Lie group. Then there exists a natural equivalence $\myS_\textup{$G$-gl}^{\circ}[-]\iso\Phi_G(-)(\id_G)$ of functors $\cat{$\bm G$-$\cat{L}$-Top}\to\Spc$, where $\Phi_G\colon\myS_\textup{$G$-gl}\to\PSh(\Glo_{/\BGcat{G}})$ is obtained from the unique equivalence $\ul\myS_\gl\iso\ul\Spc_\Glo$ by evaluating at $\BGcat{G}$.
    \begin{proof}
        If we continue to view $\cat{$\bm G$-$\cat{L}$-Top}$ as an \emph{unenriched} category, the cocartesian unstraightening of the 1-functor $\cat{$\bm G$-$\cat{L}$-Top}\to\cat{Cat}\hookrightarrow\Cat_\infty,X\mapsto \cat{$\bm G$-$\cat{L}$-Top}_{/X}$ is given by $\Ar(\cat{$\bm G$-$\cat{L}$-Top})$. Postcomposition with $\Phi_G$ and localizing therefore yields a natural transformation $\myS_\text{$G$-gl}[-]\to\PSh(\Glo_{/\BGcat{G}})_{/\Phi_G(-)}$ given pointwise by the composite equivalence
        \[
            \myS_\text{$G$-gl}[X]\xrightarrow[\raise3.5pt\hbox{$\smash{\scriptstyle\sim}$}]{\,\textup{Lemma~\ref{lemma:slice-over-fibrant}}\,}(\myS_\text{$G$-gl})_{/X}\xrightarrow[\raise3.5pt\hbox{$\smash{\scriptstyle\sim}$}]{\;\,\Phi_G\,\;}\PSh(\Glo_{/\BGcat{G}})_{/\Phi_G(X)},
        \]
        which clearly restricts to $\myS_\text{$G$-gl}^\circ[X]\iso\PSh(\Glo_{/\BGcat{G}})_{/\Phi_G(X)}^\circ$. We will now construct for every small $\infty$-category $T$ with a terminal object $1$ a natural equivalence $\PSh(T)_{/-}^\circ\iso\ev_1$ of functors $\PSh(T)\to\Spc$; applying this to $T=\Glo_{/\BGcat{G}}$ (with terminal object $\id_{\BGcat{G}}$) will then complete the proof of the theorem.
        
        To prove the claim, note that the cocartesian unstraightening of $\PSh(T)_{/-}$ is given by the target map $\Ar(\PSh(T))\to\PSh(T)$, so that the cocartesian unstraightening of $\PSh(T)_{/-}^\circ$ is given by the forgetful map $\PSh(T)_{1/}\to\PSh(T)$. The latter is also the cocartesian unstraightening of $\hom(1,-)$, which agrees with $\ev_1\colon\PSh(T)\to\Spc$ by the Yoneda lemma.
    \end{proof}
\end{lemma}

\begin{lemma}
    For any compact Lie group $G$, we have a commutative diagram
    \begin{equation}\label{diag:circ-vs-phi-sym-mon}
        \begin{tikzcd}
            \CAlg(\cat{$\cat{L}$-Top}^\boxtimes)\arrow[rr,"{\myS_\textup{$G$-gl}^\circ[\triv_G(-)]^\otimes}"]\arrow[dr,"\CAlg(\Phi^\otimes\circ\textup{loc})"', bend right=15pt] && \CMon(\Spc)\rlap.
            \\
            &\CMon(\Spc_{\Glo})\arrow[ur,"\ev_{\BGcat{G}}"', bend right=15pt]
        \end{tikzcd}
    \end{equation}
    \begin{proof}
        By \cite{schwede2018global}*{Theorem 2.1.15}, $\CAlg(\cat{$\cat{L}$-Top})$ admits a model structure with weak equivalences and fibrations defined after forgetting the commutative algebra structure; in particular, the localization $\mathfrak U_\gl$ of $\CAlg(\cat{$\cat{L}$-Top})$ has finite coproducts and products with products created by the forgetful functor to $\myS_\gl$. Moreover, as discussed right before Definition 2.5.15 of \emph{op.\ cit.}, $\mathfrak U_\gl$ is semiadditive.

        It is clear that the bottom composite through the above diagram descends to a finite product-preserving functor $\mathfrak U_\gl\to\CMon(\Spc)$. By \cite{Gepner-Groth-Nikolaus}*{Corollary~2.5}, it will therefore suffice to prove that $(\ref{diag:circ-vs-phi-sym-mon})$ commutes after postcomposition with the forgetful functor $\CMon(\Spc)\to\Spc$. This however follows immediately from the previous lemma.
    \end{proof}
\end{lemma}

\begin{proof}[Proof of Proposition~\ref{prop:group-complete-on-slices}]
    Applying the previous lemma levelwise shows that the left-hand vertical map in the commutative square
    \[
        \begin{tikzcd}
            \ul\myS_\gl^{\circ}[X]^\otimes\arrow[r,hook]\arrow[d,"{\ul\myS_\gl^{\circ}[f]}"'] & \ul\myS_\gl[X]^\otimes\arrow[d,"{\ul\myS_\gl[f]}"]\\
            \ul\myS_\gl^{\circ}[Y]^\otimes\arrow[r,hook] & \ul\myS_\gl[Y]^\otimes\hskip-1pt\rlap.\hskip1pt
        \end{tikzcd}
    \]
    is a group completion in $\CMon(\Spc_\Glo)$. To complete the proof it will then suffice to show that the horizontal inclusions are symmetric monoidal global cocompletions.
    
    We begin with the non-monoidal statement. It is clear from the construction that the equivalence $\smash{\ul\myS_\gl[Z]\iso\ul\Spc_{\Glo/\Phi(Z)}}$ from $(\ref{eq:slice-comparison})$ restricts to an equivalence between $\ul\myS_\gl^{\circ}[Z]$ and the full subcategory of $\ul\Spc_{\Glo/\Phi(Z)}$ spanned by the objects sent under the forgetful functor to the terminal object. By Remark~\ref{rk:slice-vs-shift} the latter precisely agrees with the Yoneda image, so the claim follows from the universal property of $\ul\Spc_{\Glo/\Phi(Z)}$.

    For the symmetric monoidal version it will then suffice by Remark~\ref{rk:sym-mon-cocompletion} that $\ul\myS_\gl[X]^\otimes$ and $\ul\myS_\gl[Y]^\otimes$ are globally presentably symmetric monoidal. Since the forgetful functors to $\smash{\ul\myS_\gl^\otimes}$ are conservative symmetric monoidal left adjoints, this follows at once from Corollary~\ref{cor:S-gl-times-initial}.
\end{proof}

We can now prove the remaining half of our comparison:

\begin{proof}[Proof of Theorem~\ref{thm:comparing-inputs}(2)]
    Recall from Proposition~\ref{prop:construction-of-psi-otimes} that the equivalence $\psi^\otimes$ satisfies $\psi^\otimes\circ y^\otimes\simeq y_\top^\otimes\colon\ul\Vect^\oplus\to\ul\myS_\gl[\cat{Gr}]^\otimes\vphantom{y_\top}$; in particular $\myS_\gl[i]\circ\psi^\otimes\circ y^\otimes$ lands in $\vphantom{S^\otimes}\smash{\ul\myS_\gl^{\circ}[\cat{BOP}]^\otimes}\subset\ul\myS_\gl[\cat{BOP}]^\otimes$. As every object in this subcategory is invertible by Proposition~\ref{prop:group-complete-on-slices}, the universal property shows that there exists a unique symmetric monoidal left adjoint $\smash{\Psi^\otimes\colon\ul\Spc_{\Glo/\ul\VRep^{\oplus}}^\otimes\to\ul\myS_\gl[\cat{BOP}]^\otimes}$ together with an equivalence making the square $(\ref{diag:defining-Psi-otimes})$ commute. It only remains to show that $\Psi^\otimes$ is an equivalence, for which we may equivalently show that the composite $\smash{\ul\Spc_{\Glo/\ul\Vect^\oplus}^\otimes\to\ul\myS_\gl[\cat{BOP}]^\otimes}$ is again initial among all symmetric monoidal functors inverting the Yoneda image of $\ul\Vect^\oplus$. As $\psi^\otimes$ is an equivalence and the image of $\smash{\psi^\otimes\circ y^\otimes\simeq y_\top^\otimes}$ is precisely $\smash{\ul\myS_\gl^{\circ}[\cat{Gr}]^\otimes}$, this follows from the previous proposition together with Theorem~\ref{thm:Gr-oup-completion}.
\end{proof}

\begin{rk}\label{rk:explicit-description-psi}
    We now give a concrete pointwise description of the equivalences $\psi$ and $\Psi$ underlying the equivalences $\psi^{\otimes}$ and $\Psi^{\otimes}$ from Theorem \ref{thm:comparing-inputs}. 
    Denote by $\Phi_{\Gr}\colon \ul\myS_{\gl}[\Gr]\iso {\ul\Spc_{\Glo}}_{/\Phi(\Gr)}$ the equivalence from Proposition~\ref{prop:model-categorical-model}.
    As $\Phi_{\Gr}$ sends $\ul\myS_{\gl}^{\circ}[\Gr]$ to the Yoneda image of $\Phi(\Gr)\hookrightarrow {\ul\Spc_{\Glo}}_{/\Phi(\Gr)}$, the composite $\Phi_{\Gr}\circ \psi$ restricts to a map $\alpha\colon \ul\Vect\to \Phi(\Gr)$ as depicted in the diagram on the left:
    \[
    \begin{tikzcd}
    \ul\Vect\arrow[r,"\alpha", "\exists !"', dashed] \arrow[d, hookrightarrow]&[1em]  \Phi(\Gr)\arrow[d, hookrightarrow ] & \ul\VRep\arrow[r,"A","\exists !"', dashed] \arrow[d, hookrightarrow]&[1em]  \Phi(\BOP)\arrow[d, hookrightarrow]\\ 
    \ul\Spc_{\Glo/\ul\Vect}\arrow[r,"\Phi_{\Gr}\psi", "\sim"'] & {\ul\Spc_{\Glo}}_{/\Phi(\Gr)} &     {\ul\Spc_{\Glo}}_{/\ul\VRep}\arrow[r,"\Phi_{\BOP}\Psi","\sim"'] & \ul\Spc_{\Glo/\Phi(\BOP)}.
    \end{tikzcd}    \]
    As $\Phi_{\Gr}\circ \psi$ is an equivalence, $\alpha$ has to be an equivalence. 
    It follows from Example~\ref{ex:post-composition-slice} that  for every compact Lie group $G$, $\psi^{-1}$ sends an object of $\ul\myS_{\gl}[\Gr](\BGcat{G})$ given by a map of orthogonal $G$-spaces 
    $\xi \colon X\to \triv_G\,\Gr$ to the object of $\Spc_{\text{$G$-gl}/\triv_G\,\ul\Vect}$ given by the composite    
    \[\Phi(X)\xrightarrow{\;\Phi(\xi)\;} \Phi(\triv_G\,\Gr)\iso \triv_G\,\Phi(\Gr) \xrightarrow[\raise3.5pt\hbox{$\scriptstyle\smash{\sim}$}]{\;\triv_G\,\alpha^{-1}\;}\triv_G\,\ul\Vect,\] where $\Phi$ denotes the unique equivalence $\ul\myS_{\gl}\iso \ul\Spc_{\Glo}$ and the unlabelled equivalence comes from naturality of $\Phi$.

    We can analogously describe $\Psi$ on the level of objects: 
    The equivalence $\Phi_{\BOP}\colon \ul\myS_{\gl}[\BOP]\iso {\ul\Spc_{\Glo}}_{/\Phi(\BOP)}$ from Proposition~\ref{prop:model-categorical-model} sends $\ul\myS_{\gl}^{\circ}[\BOP]$ to the Yoneda image of $\Phi(\BOP)$. 
    Arguing as before, $\Phi_{\BOP}\circ \Psi$ restricts to an equivalence $A\colon \ul\VRep\to \Phi(\BOP)$ as depicted in the diagram on the right above, and we see that $\Psi^{-1}$ sends a map of orthogonal $G$-spaces $\xi\colon X\to \triv_G\,\BOP$ to the object of  $\Spc_{\text{$G$-gl}/\triv_G\,\ul\VRep}$ given by the composite\[\Phi(X)\xrightarrow{\;\Phi(\xi)\;}\Phi(\triv_G\,\BOP)\iso \triv_G\,\Phi(\BOP)\xrightarrow[\raise3.5pt\hbox{$\scriptstyle\smash{\sim}$}]{\;\triv_GA^{-1}\;}\triv_G\,\ul\VRep.\] 

        With slightly more work one could show that for any ultra-commutative monoid $M$, there is a symmetric monoidal equivalence $\ul\myS_{\gl}[M]^{\otimes}\iso \ul{\Spc}_{\Glo/M}^{\otimes}$ enhancing the equivalence $\Phi_M$ from Proposition~\ref{prop:model-categorical-model}. This then implies that $\alpha$ and $A$ enhance to symmetric monoidal equivalences \[\alpha^{\otimes}\colon \ul\Vect^{\oplus}\iso\Phi^{\otimes}(\Gr)\qquad\text{and}\qquad A^{\otimes}\colon \ul\VRep^{\oplus}\iso \Phi^{\otimes}(\BOP).\]   
\end{rk}

\subsection{Comparing the outputs} We can now state our comparison:

\begin{thm}\label{thm:comparison-Thom}
    The diagrams
    \[
        \begin{tikzcd}[row sep=tiny]
            \ul\Spc_{\Glo/\ul\Vect^\oplus}^\otimes\arrow[dd,"\psi^\otimes"',"\sim"]\arrow[dr, bend left=15pt, "\ul\th_\gl^\otimes"] &&[1.5em]
            \ul\Spc_{\Glo/\ul\VRep^{\oplus}}^\otimes\arrow[dr, bend left=15pt,"\ul\Th_\gl^\otimes"]\arrow[dd,"\Psi^\otimes"',"\sim"]\\
            & \ul\myS_{\gl,*}^\otimes && \ul\mySp_\gl^\otimes\\
            \ul\myS_\gl[\cat{Gr}]^\otimes\arrow[ur,"\ul\myth_\gl^\otimes"', bend right=15pt]&&\ul\myS_\gl[\cat{BOP}]^\otimes\arrow[ur,"\ul\myTh_\gl^\otimes"', bend right=15pt]
        \end{tikzcd}
    \]
    of symmetric monoidal global functors commute up to equivalence, where $\psi^\otimes$ and $\Psi^\otimes$ are as in Theorem~\ref{thm:comparing-inputs}.
    \begin{proof}
        We begin with the comparison of Thom spaces. We defined $\smash{\ul\th_\gl^\otimes}$ to be a global left adjoint, while $\smash{\ul\myth_\gl^\otimes}$ is so by Proposition~\ref{prop:gl-Thom-space-la}; thus, it will suffice to prove the comparison after restricting to $\ul\Rep^\oplus$. 
        
        By definition, the upper path through the diagram is simply $\ul{\mathfrak j}_\gl^\otimes\colon\ul\Vect^\oplus\to\smash{\ul\myS_{\gl,*}^\otimes}$, which is obtained from the continuous Borel construction of \[\cat{const}^\otimes\circ\cat{j}^\otimes\colon\core\cat{L}^\oplus\to\cat{$\cat{L}$-Top}_*^{\smallboxsmash}\] through localization. On the other hand, the equivalence $\psi^\otimes$ was constructed as extension of $y_\top^\otimes$, so that the lower composite can be identified via Proposition~\ref{prop:th-on-corep-sym-mon} with the symmetric monoidal global functor induced by
        \[
            \core\cat{L}^\oplus\xrightarrow{\;(\cat{y}_+^\otimes,\cat{j}^\otimes)\;}\cat{$\cat{L}$-Top}^{\smallboxsmash}_*\times\cat{Top}^{\smashp}_*\xrightarrow{\;{-}\smashp{-}\;}\cat{$\cat{L}$-Top}_*^{\smallboxsmash}
        \]
        If we replaced $\smash{\cat{y}_+^\otimes\colon\core\cat{L}^\oplus\to\cat{$\cat{L}$-Top}_*^{\smallboxsmash}}$ by the functor constant at the unit, this composite would be isomorphic to $\cat{const}^\otimes\circ\cat{j}^\otimes$ as topological symmetric monoidal functors via the unit isomorphism. By functoriality of $\Ntop((-)^\dual)$, it will therefore suffice to construct a natural equivalence $\fgt^\otimes\circ y_\top^\otimes\simeq(\const\,1)^\otimes$ of symmetric monoidal global functors $\smash{\ul\Vect^\oplus\to\ul\myS_\gl^\otimes}$. However, both of these functors send each object to a terminal object, so they are equivalent in a unique way. This completes the comparison of Thom spaces.

        For the comparison of Thom spectra, it similarly suffices (by additionally invoking the universal property of group completion) to construct the equivalence after restriction to $\ul\Vect^\oplus$. Again by design, the top path through the diagram then becomes $\Sigma^\bullet\circ\ul{\mathfrak j}_\gl^\otimes$. On the other hand, commutativity of $(\ref{diag:defining-Psi-otimes})$ shows that the bottom path is given by $\smash{\ul\myTh_\gl^\otimes\circ\ul\myS_\gl[i]^\otimes\circ y_\top^\otimes}$, which Corollary~\ref{cor:thom-space-vs-thom-spectra} further identifies with $\smash{\Sigma^\bullet\circ\ul\myth_\gl^\otimes\circ y_\top^\otimes}$. Thus, the claim follows from the above unstable comparison.
    \end{proof}
\end{thm}

As a concrete consequence of this, we will now express the global homotopical real bordism spectrum $\cat{MO}$ (Example~\ref{ex:MO}) as a global parametrized colimit. We begin with the analogous statement for $\cat{MOP}$.

\begin{constr}\label{constr:mop-cocone}
    Let $\bar y_\top^\otimes\colon\ul\VRep^{\oplus}\to\ul\myS_\gl[\cat{BOP}]^\otimes$ be the unique extension of $\smash{\ul\myS_\gl[i]^\otimes\circ y_\top^\otimes}$. As $\id_{\cat{BOP}}$ is terminal, there is then a unique natural transformation $\tilde\sigma$ from $\bar y_\top$ to the functor constant at $\id_{\cat{BOP}}$. Applying $\ul\myTh_\gl$ and appealing to the identification from the previous theorem, we then obtain a cocone
    \[
        \sigma\colon\ul{\mathfrak J}_\gl\simeq \ul\myTh_\gl\circ\bar y_\top\xrightarrow{\;\ul\myTh_\gl(\tilde\sigma)} \ul\myTh_\gl(\id_{\cat{BOP}})=\cat{MOP}.       
    \]
\end{constr}

\begin{rk}\label{rk:univ-cocone-thom-classes}
    Unravelling the identifications from the proof of the above theorem, the restriction of the cocone $\sigma$ along $\ul\Vect\to\ul\VRep$ is given incoherently as follows: if $G$ is any compact Lie group, and $V$ is a $G$-representation, then $\sigma_V$ is given by the zig-zag $S^V\xleftarrow{\;\smash{\lower2pt\hbox{$\scriptstyle\sim$}}\;}\cat{O}(V,-)\smashp S^{V\oplus V}\to\cat{MOP}$ where the first map is adjunct to the evident isomorphism $S^V\smashp S^V\cong S^{V\oplus V}$ and the second one is adjunct to $S^{V\oplus V}\to\cat{MOP}(V),x\mapsto(V\oplus V,x)$. Thus, viewing $\sigma_V$ as an element of $\pi_0^G(\Omega^V\cat{MOP})$, this is precisely the \emph{Thom class} $\sigma_{G,V}$ associated to the $G$-representation $V$ from \cite{schwede2018global}*{Construction~6.1.15}, motivating our notation. Using Theorem~6.1.17 of \emph{op.\ cit.}, one can then similarly describe the effect of $\sigma$ on a virtual representation of the form $0\ominus V$ in terms of the \emph{inverse Thom class} from (6.1.12) of \emph{op.\ cit.} In this sense, the cocone $\sigma$ is a coherent arrangement of all the (inverse) Thom classes of $\cat{MOP}$, as well as their products.
\end{rk}

\begin{prop}\label{prop:MOP-as-colim} 
    The cocone $\sigma$ expresses $\cat{MOP}$ as the colimit of the global functor $\ul{\mathfrak J}_\gl\colon\ul\VRep\to\ul\mySp_\gl$, i.e.\ as a left adjoint object with respect to $\const\colon\mySp_\gl\to\Fun(\ul\VRep,\ul\mySp_\gl)$.
    \begin{proof}
        As $\ul\myTh_\gl$ is a global left adjoint, it suffices to prove that the terminal object of $\ul\myS_\gl[\cat{BOP}]$ is equivalent to the global colimit of $\bar y_\top$ (in a necessarily unique way). This on the other hand may be checked after postcomposing with the equivalence $\Psi^{-1}\colon\ul\myS_\gl[\cat{BOP}]\iso\ul\Spc_{\Glo/\ul\VRep^{\oplus}}$. We claim that the resulting composite $\ul\VRep^{\oplus}\to\ul\Spc_{\Glo/\ul\VRep^{\oplus}}^\otimes$ is equivalent, even as a \emph{symmetric monoidal} global functor, to the parametrized Yoneda embedding; the proposition will then follow since the colimit of the latter is indeed terminal by \cite{martiniwolf2021limits}*{Proposition 6.1.3}.
        To prove the claim, we may again restrict to $\ul\Vect^\oplus$, i.e.~it will suffice to identify $(\psi^\otimes)^{-1}\circ y_\top^\otimes$ with the Yoneda embedding of $\ul\Vect^\oplus$. This however holds by construction of $\psi^\otimes$.
    \end{proof}
\end{prop}

\begin{constr}
    Sending an inner product space to its dimension defines an enriched symmetric monoidal functor $\core\cat{L}^\oplus\to(\mathbb N,+)$. Applying the continuous Borel construction and group completing we therefore obtain $\ul\VRep^{\oplus}\to\const(\mathbb Z,+)$. For every $d\ge0$, we define $\smash{\ul\VRep_{[d]}}$ as the preimage of $d\in\Z$; if we think of equivalence classes of objects of $\ul\VRep$ as virtual representations $V\ominus W$, then $\smash{\ul\VRep_{[d]}}$ consists precisely of those virtual representations with $\dim(V)-\dim(W)=d$.
\end{constr}

\begin{lemma}
    The restriction of the cocone $\tilde\sigma$ to $\ul\VRep_{[d]}$ factors (necessarily uniquely) through $\smash{\cat{BOP}^{[d]}\hookrightarrow\cat{BOP}}$.
    \begin{proof}
        For any map of $G$-orthogonal spaces $f\colon X\to\cat{BOP}$, $X$ splits up to isomorphism as $\coprod_{d\in\Z} f^{-1}(\cat{BOP}^{[d]})$. If the underlying space of $X$ is path-connected (e.g.\ when $X$ is weakly contractible), exactly one of these summands is non-empty, and we define $d(f)$ to be the corresponding index; equivalently, this is the unique $d$ such that $f$ factors through $\cat{BOP}^{[d]}$ (on the pointset level). It is clear that if $f\to g$ is any map over $\cat{BOP}$ and $g$ factors through $\cat{BOP}^{[d]}$, then so does $f$; in particular, if $f$ and $g$ are moreover in $\myS_\text{$G$-gl}^\circ[\cat{BOP}]$, then $d(f)=d(g)$.
        
        With this established, it is easy to check that $f\mapsto d(f)$  defines a symmetric monoidal global functor $\smash{\ul\myS_\gl^\circ[\cat{BOP}]^\otimes\to\const(\Z,+)}$. Moreover, its restriction along $\myS_\gl[i]^\otimes\circ y_\top^\otimes$ agrees with $\dim\colon\ul\Vect^\oplus\to\const(\Z,+)$ by direct inspection, so $d\circ\bar y_\top=\dim$ by the universal property of group completion. Thus, $\bar y_\top$ maps $\smash{\ul\VRep_{[d]}}$ into the subcategory of $\ul\myS_\gl[\cat{BOP}]$ of maps factoring through ${\cat{BOP}^{[d]}}$, which is in turn equivalent to $\tilde\sigma$ factoring accordingly.
    \end{proof}
\end{lemma}

Applying $\ul\myTh_\gl$ we therefore get cocones
$\smash{\sigma_{[d]}\colon\mathfrak J_\gl|_{\ul\VRep_{[d]}}\to\smash{\cat{MOP}^{[d]}}}$ for $d\in\Z$.

\begin{thm}\label{thm:MO-as-colim}
    The cocone $\sigma_{[0]}$ exhibits $\cat{MO}=\smash{\cat{MOP}^{[0]}}$ as the global colimit of $\smash{\ul{\mathfrak J}_\gl|_{\ul\VRep_{[0]}}}$.
    \begin{proof}
        We will more generally prove this for any $d\in\Z$, and similarly to before it will suffice to show that $\tilde\sigma_{[d]}$ exhibits $\cat{BOP}^{[d]}\hookrightarrow\cat{BOP}$ as a colimit.

        As $\ul\VRep=\coprod_{d\in\Z}\smash{\ul\VRep_{[d]}}$, we get an equivalence
        \[
            \colim_{\ul\VRep}\bar y_\top \simeq \coprod_{d\in\Z}\,\colim_{\ul\VRep_{[d]}}\,\bar y_\top|_{\ul\VRep_{[d]}}
        \]
        in $\myS_\gl[\cat{BOP}]$. As the left-hand side is terminal by Proposition~\ref{prop:MOP-as-colim}, so is the coproduct on the right-hand side. Thus, the coproduct of the maps 
        \[
            \colim\bar y_\top|_{\ul\VRep_{[d]}}\to\big(\cat{BOP}^{[d]}\hookrightarrow\cat{BOP}\big)
        \]
        induced by $\tilde\sigma_{[d]}$ for all $d\in\Z$ is a map between terminal objects and hence an equivalence in $\myS_\gl[\cat{BOP}]$. As a coproduct of maps in $\myS_\gl[\cat{BOP}]\simeq(\myS_{\gl})_{/\cat{BOP}}$ (or any other $\infty$-topos) is an equivalence only if each summand is so, this completes the proof of the theorem.
    \end{proof}
\end{thm}

\begin{rk}More generally, Remark~\ref{rk:explicit-description-psi} and Lemma~\ref{lm:pointwise-par-colimit} allow us to describe the $G$-global Thom spectrum of any map $\xi\colon X\to \triv_G\,\BOP$ of orthogonal $G$-spaces as a parametrized colimit: we have an equivalence
\[ \myTh_\text{$G$-gl}(\xi)\simeq\colim_{\Phi(X)}\big(\triv_G(\ul{\mathfrak{J}}_{\gl})\circ \triv_G(A)\circ \Phi(\xi)\big)\]
in $\mySp_\text{$G$-gl}$, where $A\colon \Phi(\BOP)\iso \ul\VRep$ denotes the equivalence constructed in Remark~\ref{rk:explicit-description-psi}, $\Phi\colon \ul\myS_{\gl}\iso \ul\Spc_{\Glo}$ is the unique equivalence as before, and $\colim_{\Phi(X)}$ denotes the $\Glo_{/\BGcat{G}}$-parametrized colimit-functor.
\end{rk}

Let us write $\cat{MOP}_G$ for the underlying $G$-spectrum of the global spectrum $\cat{MOP}$, i.e.\ the image under the composite $\mySp_\gl\to\mySp_\text{$G$-gl}\to\mySp_G$, and we similarly define $\cat{MO}_G$ as the underlying $G$-spectrum of $\cat{MO}$. Put differently, if we start with the underlying $\Orb$-section of the global section $\cat{MOP}$ of $\ul\mySp_\gl$ and apply the right adjoint $\ul\mySp_\gl|_\Orb\to\ul\mySp|_\Orb$, then the resulting $\Orb$-section $\cat{MOP}_\bullet$ of $\ul\mySp|_{\Orb}$ has components $\cat{MOP}_G$, and analogously for $\cat{MO}$.

Analogously to Construction~\ref{constr:mop-cocone} we then get an ($\Orb$-parametrized) cocone $\sigma$ from $\ul{\mathfrak J}|_{\Orb}\colon\ul\VRep|_{\Orb}\to\ul\mySp|_{\Orb}$ to $\cat{MOP}_\bullet$ as the image of the unique cocone $\bar y_\top|_{\Orb}\to\id_{\cat{BOP}}$ in $\ul\mySp_\gl[\cat{BOP}]|_{\Orb}$, and as before this restricts to a cocone from $\ul{\mathfrak J}|_{\ul\VRep_{[0]}|_{\Orb}}$ to $\cat{MO}_\bullet$.

\begin{thm}\label{thm:MOP-Orb-colim}
    The induced maps
    \[
        \colim_{\ul\VRep|_{\Orb}}\ul{\mathfrak J}\to\cat{MOP}_\bullet\qquad\text{and}\qquad
        \colim_{\ul\VRep_{[0]}|_{\Orb}}\ul{\mathfrak J}\to\cat{MO}_\bullet
    \]
    are equivalences in $\ul\mySp|_{\Orb}$.
\end{thm}

The proof will require some preparations. Let us call a map in $\ul\Spc_{\Glo/\ul\VRep}$ an \emph{equivariant weak equivalence} if it is inverted by the right adjoint $R$ of the inclusion $(\ul\Spc_{\Orb\triangleright\Glo/\ul\VRep})|_\Orb\hookrightarrow(\ul\Spc_{\Glo/\ul\VRep})|_\Orb$; by the proof of Lemma~\ref{lemma:right-adjoint-right-Bousfield} this is equivalent to being inverted by the composite
\[
    (\ul\Spc_{\Glo/\ul\VRep})|_\Orb\xrightarrow{\,\fgt\,}(\ul\Spc_{\Glo})|_\Orb\xrightarrow{\;r\;}(\ul\Spc_{\Orb\triangleright\Glo})|_{\Orb}
\]
of the forgetful functor with the right adjoint to the inclusion. We moreover define a map in $\ul\myS_\gl[\cat{BOP}]$ to be an \emph{equivariant weak equivalence} if it is inverted by the analogous composite
\[
    \ul\myS_\gl[\cat{BOP}]|_\Orb\xrightarrow{\,\fgt\,}\ul\myS_\gl|_\Orb\xrightarrow{\;r\;}\ul\myS|_\Orb.
\] 
The following will be one of the key ingredients for proving the theorem:

\begin{lemma}
    The equivalence $\Psi\colon\ul\Spc_{\Glo/\ul\VRep}\iso\ul\myS_\gl[\cat{BOP}]$ preserves and reflects equivariant weak equivalences.
    \begin{proof}
        We claim that the square
        \[
            \begin{tikzcd}
                \ul\Spc_{\Glo/\ul\VRep^\oplus}^\otimes\arrow[d,"\fgt"']\arrow[r,"\Psi"] & \ul\myS_\gl[\cat{BOP}]^\otimes\arrow[d,"\fgt"]\\
                \ul\Spc_{\Glo}^\times\arrow[r,"\sim"'] & \ul\myS_\gl^\times
            \end{tikzcd}
        \]
        of symmetric monoidal global functors commutes up to equivalence. As all functors in question are symmetric monoidal left adjoints, this can be checked after restricting to $\ul\Rep$, where it follows immediately from the definition of the equivalence $\Psi$ that both composites send every object to the terminal object, so that the square even commutes up to unique equivalence.

        Plugging in the definitions, we are then reduced to constructing a commutative square as depicted on the left in the following diagram, where the vertical maps are the underlying $\Orb$-functors of the unique equivalences:
        \[
            \begin{tikzcd}
                \ul\Spc_{\Orb\triangleright\Glo}|_\Orb\arrow[from=r,"r"']\arrow[d,"\sim"'] & \ul\Spc_{\Glo}|_\Orb\arrow[d,"\sim"] &[1em] \ul\Spc_{\Orb\triangleright\Glo}\arrow[r,hook]\arrow[d,"\sim"'] & \ul\Spc_{\Glo}\arrow[d,"\sim"]\\
                \ul\myS|_\Orb\arrow[from=r,"r"] &\ul\myS_\gl|_\Orb & \ul\myS\arrow[r,"c"'] & \ul\myS_\gl
            \end{tikzcd}
        \]
        The square on the right commutes up to unique equivalence as both paths are $\Orb$-left adjoints preserving the terminal object. Passing to underlying $\Orb$-functors and then to mates therefore provides the desired square on the left.
    \end{proof} 
\end{lemma}

\begin{proof}[Proof of Theorem~\ref{thm:MOP-Orb-colim}]
    Combining the previous lemma with Theorems~\ref{thm:Thom-vs-fgt} and~\ref{thm:comparison-Thom}, the left adjoint $\ul\Th_\gl\colon\ul\myS_\gl[\cat{BOP}]|_{\Orb}\to\ul\mySp_\gl|_{\Orb}$ descends to the localizations at the equivariant weak equivalences, and so does its right adjoint as a consequence of Proposition~\ref{prop:equiv-thom-from-global}. We can rephrase this as saying that the composite $\ul\myS_\gl[\cat{BOP}]|_{\Orb}\to\ul\mySp_\gl|_{\Orb}\to\ul\mySp|_{\Orb}$ descends through the localization $\Xx$ of the source at the equivariant weak equivalences, and that the resulting $\Orb$-functor is again a left adjoint. To prove the first statement, it will therefore suffice (analogously to the global setting) that the colimit of $\bar y_\top\colon\ul\VRep|_{\Orb}\to\Xx$ is terminal. Translating this through the equivalence $\Psi$, this is in turn equivalent to the $\Orb$-colimit of 
    \begin{equation}\label{eq:should-have-terminal-colimit}
        \ul\VRep|_{\Orb}\xrightarrow{\;y\;}\ul\Spc_{\Glo/\ul\VRep}|_{\Orb}\xrightarrow{\;R\;}\ul\Spc_{\Orb\triangleright\Glo/\ul\VRep}|_{\Orb}
    \end{equation}
    being terminal. The left adjoint $L$ of $R$ was defined so that $y=L\circ y$. Together with full faithfulness of $L$, this implies that $(\ref{eq:should-have-terminal-colimit})$ agrees up to postcomposition with an equivalence with the ($\Orb$-parametrized) Yoneda embedding $\smash{\ul\VRep}|_{\Orb}\to\ul\PSh_{\Orb}(\ul\VRep)$. Thus, its colimit is terminal by \cite{martiniwolf2021limits}*{Proposition~6.1.3}, finishing the proof of the first claim.

    The second claim will then follow formally as in the global setting once we show that the localization $\myS_\text{$G$-gl}[\cat{BOP}]\to\Xx(\BGcat{G})$ preserves coproducts for every $G$, or equivalently that equivariant weak equivalences are stable under arbitrary small coproducts. As before, we may prove the corresponding statement for $\ul\Spc_{\Glo/\ul\VRep}$ instead, which in turn follows immediately from the localization being an $\Orb$-left adjoint by Lemma~\ref{lemma:right-adjoint-right-Bousfield}.
\end{proof}

Fix a compact Lie group $G$ and recall once more the functor $\fgt\colon\Orb_{G}\to\Glo$ lifting to an equivalence $\Orb_G\iso\Orb_{/\BGcat{G}}$ sending $G/H$ to $\BGcat{H}\hookrightarrow\BGcat{G}$. By \cite{martiniwolf2021limits}*{Remark 4.1.9}, $\fgt^*\colon\Fun(\Orb,\CAT_\infty)\to\Fun(\Orb_G,\CAT_\infty)$ preserves colimiting cocones; thus, writing $\ul{\mathfrak J}_G\coloneqq\fgt^*\ul{\mathfrak J}$ and $\ul\VRep_G\coloneqq\fgt^*\ul\VRep$, we obtain the following equivalent reformulation of the previous theorem:

\begin{cor}\label{cor:MOG-as-colim}
    The maps
    \[
        \colim_{\ul\VRep_G}\ul{\mathfrak J}_G\to\cat{MOP}_G\qquad\text{and}\qquad
        \colim_{\ul\VRep_{G,[0]}}\ul{\mathfrak J}_G\to\cat{MO}_G
    \]
    induced by the $\Orb_G$-cocone $\fgt^*\sigma$ are equivalences in $\mySp_G$.\qed
\end{cor}

\appendix
\setcounter{section}{0}

\renewcommand{\thesection}{\Alph{section}}
\chapter*{\for{toc}{\hspace{-1em}}\except{toc}{Appendix}}
\section[A Thom spectrum functor with marvelous properties]{\for{toc}{\hspace{-.22em}}A Thom spectrum functor with marvelous pointset properties\except{toc}{\mdseries\\[.5ex]\itshape{by} \textsc{Stefan Schwede}}}
\addtocontents{toc}{\hskip1.2em{\small--- \textit{by} \textsc{Stefan Schwede}}}\label{app:thom-Stefan-functor}

Constructions~\ref{constr:Thom-space-functor-gl} and~\ref{constr:Thom-spectrum-functor} define a global Thom space and global Thom spectrum functor, respectively, as (symmetric monoidal) topologically enriched functors out of certain slices of the topological category $\cat{$\cat{L}$-Top}$ of orthogonal spaces. Taking $G$-objects on both sides for any compact Lie group $G$ then yields the $G$-global Thom space and Thom spectrum functor. The goal of this appendix is to prove the following result: 
  \begin{thm}\label{thm:Thom-marvelous-appendix} Let $G$ be a compact Lie group. 
\begin{enumerate}
\item The Thom space functor $\Gospc{G}_{/\Gr}\to \cat{$\bm G$-$\cat{L}$-Top}_*$ is left Quillen and fully homotopical with respect to the $G$-global model structures. 
\item The Thom spectrum functor $\Gospc{G}_{/\BOP}\to \Gosp{G}$ is left Quillen and fully homotopical with respect to the $G$-global model structures.
\end{enumerate}
  \end{thm}
Recall that the $G$-global model structures on orthgonal $G$-spaces and orthogonal $G$-spectra were constructed by localizing the respective $G$-level model structures. To prove the above theorem, we will first show the analogous statements for the $G$-level model structures on source and target, and then deduce from this the above statement. 

We begin by describing right adjoints of the Thom space and Thom spectrum functor. For that, we first show that for an inner product space $V$, the Thom space functor 
\[ T\colon\GTop{G}_{/\Gr(V)}\to \GTop{G}_*\] from Construction~\ref{constr:thom-over-Gr} admits a right adjoint $q_V$. 
\begin{constr}\label{constr:thom-level-right-adjt}Fix an inner product space $V$.
For a based topological space $X$, 
we define \[\Omega_V(X)\coloneqq \coprod_{d\ge0}\cat L(\R^d,V)\times_{\O(d)}\Omega^dX,\] where $\O(d)$ acts on $\Omega^dX=\maps_*(S^d,X)$ via precomposition.
The reference map to the Grassmannian is the projection \[\coprod_{d\ge0}\cat L(\R^d,V)\times_{\O(d)}\Omega^dX\to \coprod_{d\ge0}\cat L(\R^d,V)/\O(d)=\Gr(V).\] 

 This map is a locally trivial fiber bundle whose fiber over $L\in \Gr(V)$ `is'
  the loop space $\Omega^L=\maps_*(S^L,Z)$. This then defines a topological functor $q_V\colon\cat{Top}_*\to\cat{Top}_{/\cat{Gr}(V)}$ with functoriality via postcomposition. If $G$ is any compact Lie group, then pulling through the $G$-actions lifts this to $\cat{$\bm G$-Top}_*\to\cat{$\bm G$-Top}_{/\cat{Gr}(V)}$.

  For a continuous map $f\colon A\to\Gr(V)$, we define the adjunction unit
  \[
     \eta(f)\colon  A \to \Omega_V(T(f)) \text{\qquad by\qquad} \eta(f)(a) = [\psi,\{x\mapsto (\psi(x),a)\}]; 
  \]
  here $\psi\in\cat{L}(\mathbb R^d,V)$ is any linear isometric embedding with $\psi(\mathbb R^d)=f(a)$.
  We omit the easy verification that the map $\eta(f)$ is well-defined and lies over $\Gr(V)$, i.e., $q_{T(f)}\circ\eta(f)=f$.
  For a based space $Z$, we define the adjunction counit
  \[ 
    \epsilon_Z\colon  T(q_V(Z)) \to Z \text{\qquad by\qquad}
    \epsilon_Z(w,[\psi,f]) =  f(\psi^{-1}(w));
  \]
  here $\psi\colon\mathbb R^d\to V$ is a linear isometric embedding, $f\colon S^d\to Z$ a based continuous map,
  and $w\in \psi(\mathbb R^d)$. We omit the verifications that the maps $\eta(f)$ and $\epsilon_Z$
  are continuous and satisfy the triangle equality of an adjunction.
  We have thus extended the two topologically enriched functors
  \begin{align}\label{eq:Thom-level-rightadjt} T\colon \cat{$\bm G$-Top}_{/\Gr(V)}\rightleftarrows \cat{$\bm G$-Top}_*\colon q_V\end{align}
  to an adjoint pair.
\end{constr}

As we will now explain, the right adjoints $q_V$ for varying $V\in\cat L$ assemble into a right adjoint $\cat{q}\colon \Gospc{G}_*\to \Gospc{G}_{/\Gr}$ of the Thom space functor $\cat{th}$ from Construction~\ref{constr:Thom-space-functor-gl}:

\begin{constr} 
  For a based orthogonal space $X$, define 
  \[ \bm\Omega(X)= \coprod_{d\geq 0}\cat L(\R^d,-)\times_{\O(d)}\Omega^dX ;\]
  the reference map $\cat{q}(X)\colon\bm\Omega(X)\to\cat{Gr}$ is induced by the {projection} as before, so that at an inner product space $V$, $q(V)$ is the image $q_V\colon\Omega_V(X(V))\to\Gr(V)$ of $X(V)\in \cat{Top}_*$ under the right adjoint $q_V$ from the previous construction. This becomes a functor via the functoriality of the individual $q_V$, and this functor is again topological as mapping spaces in functor categories are topologized as subspaces of the product, and since each $q_V$ was topological. For a map $f\colon A\to \Gr$ of orthogonal spaces the unit transformations $\eta_V(f(V)), V\in\cat L$ from \ref{constr:thom-level-right-adjt} assemble into a map of orthogonal spaces $\eta(f)\colon \cat{q}(\cat{th}(f))\to X$ with $\cat{q}_{\Th(f)}\circ \eta(f)=f$. 

  For a pointed orthogonal space $X$, the counit natural transformations \[\epsilon_V(X(V))\colon \cat{th}(\bm\Omega(X)(V))=T(\Omega_V(X(V)))\to X(V), V\in\cat L\] define a map of pointed orthogonal spaces $\epsilon_X\colon \cat{th}(\bm\Omega(X))\to X$.
  The maps $\epsilon(X)$ for $X\in \cat{L-Top}_*$ and $\eta(f)$ for ${f\in \cat{L-Top}_{/\Gr}}$ define natural transformations $\id\to \textbf{q}\circ \cat{th}$ and $\cat{th}\circ \cat{q}\to\id$ which exhibit $\cat{q}$ as left adjoint to $\cat{th}$. Pulling through the actions as before, we then obtain for every compact Lie group $G$ an adjoint pair
  \[\cat{th}\colon  \Gospc{G}_{/\Gr}\rightleftarrows \Gospc{G}_*\colon \textbf{q}.\]
\end{constr}

Similarly to the above, we can construct a topological right adjoint 
  \[ \textbf{Q}\colon \Gosp{G}\to \Gospc{G}_{/\BOP} \]
  to the Thom spectrum functor  $\cat{Th}\colon \Gospc{G}_{/\BOP}\to \Gosp{G}$ from Construction~\ref{constr:Thom-spectrum-functor} by assembling the right adjoints to the individual Thom space functors.
  In essence, the functor $\textbf{Q}$ `bundles' all loop spaces of all values of an orthogonal spectrum
  into a single object, organized as an orthogonal space over $\BOP$: 
\begin{constr}\label{con:bq}
  We let $X$ be  an orthogonal spectrum, and we let $V$ be an inner product space. We set
   \[ (\bm\Omega X)(V) \coloneqq \Omega_{V\oplus V}(X(V))= \coprod_{d\geq 0}\cat L(\R^d,V\oplus V)\times_{\O(d)} \Omega^dX(V) \ , \]
  the total space of the locally trivial bundle over $\BOP(V)=\Gr(V\oplus V)$ whose fiber over
  a linear subspace $L$ of $V\oplus V$ `is' the $L$-loop space $\Omega^LX(V)$.
  The reference map 
  \[ \cat{Q}(X)(V)\coloneqq q_{V\oplus V}(X(V))\colon (\bm\Omega X)(V) \to \BOP(V)\]
  is the bundle projection sending $[\psi,f]$ to $\psi(\mathbb R^d)$.

  In order to define the structure maps of the orthogonal space $\bm\Omega(X)$, we consider the continuous map
  \[ j\colon \cat L(\R^d,V\oplus V)\to\cat L(\R^{d+e},(V\oplus\mathbb R^e)^2),\quad
    j(\psi)(x,y) = \kappa_{V,\mathbb R^e}( \psi(x), (y,0)), \]
  where $x\in\mathbb R^d$, $y\in\mathbb R^e$, and  $\kappa_{V,\mathbb R^e}\colon (V\oplus V)\oplus(\R^e\oplus\R^e)\to (V\oplus\R^e)^2$
  is the shuffling isomorphism.
  We define the structure map
  \[ (\bm\Omega X)(V,W)\colon \cat L(V,W)\times (\bm\Omega X)(V)\to (\bm\Omega X)(W)\]
  as the disjoint union, over $d\geq 0$, of the maps
  \begin{equation}\label{eq:d_summand}\begin{aligned}
      \cat L(V,W)\times \left( \cat L(\R^d,V^2)\times_{\O(d)} \Omega^dX(V)\right)
      &\longrightarrow\cat L(\R^{d+e},W\oplus W)\times_{\O(d+e)} \Omega^{d+e}X(W)\\
      (\varphi,[\psi,f]) &\longmapsto[ (A\oplus A)\circ j(\psi), X(A)\circ\sigma\circ (f\smashp S^e)],
  \end{aligned}
  \end{equation}
  where $e=\dim(W)-\dim(V)$, $A\colon V\oplus \mathbb R^e\iso W$ is any linear isometry such that $A|_V=\varphi$, and $\sigma\colon X(V)\smashp S^e\to X(V\oplus\R^e)$ and $X(A)\colon X(V\oplus \R^e)\to X(W)$ are the structure maps of $X$. We omit the easy verification that this is independent of the choice of $A$ and compatible with the composition of $\cat{L}$ as well as with the projections to $\cat{BOP}$.
  \end{constr}
While clearly expected, it might not be entirely obvious that the resulting action map is
continuous; so we provide an argument.

\begin{prop}\label{prop:structure_continuous}
  For every orthogonal spectrum $X$ and all inner product spaces $V$ and $W$, the map $(\textup{\ref{eq:d_summand}})$ is continuous. 
\end{prop}
\begin{proof}
  Product with $\cat L(V,W)$ preserves disjoint unions, so we may show that the restricted map
  \eqref{eq:d_summand} is continuous for every $d\geq 0$, where $e=\dim(W)-\dim(V)$.
  We may assume without loss of generality that $W=V\oplus\mathbb R^e$.
  The map
  \[ j\colon\cat L(\mathbb R^d,V^2)\to \cat L(\mathbb R^{d+e},(V\oplus\mathbb R^e)^2),\quad
    j(\psi)(x,y)= \kappa_{V,\mathbb R^e}( \psi(x), (y,0)) \]
  is continuous, where $x\in\mathbb R^d$, $y\in\mathbb R^e$, and $\psi\in\cat L(\mathbb R^d,V^2)$. Moreover, the homomorphism \begin{align*}\O(V\oplus \mathbb R)&\longrightarrow\O((V\oplus \mathbb R)^{2})\\ A&\longmapsto \begin{pmatrix}
  A & 0 \\ 
  0 & A 
  \end{pmatrix}\end{align*} defines a continuous $\O(V\oplus\mathbb R)$-action on the target of $j$, so the map
  \begin{align*}
    \O(V\oplus\mathbb R^e)\times \cat L(\mathbb R^d,V^2)\times\Omega^dX(V)
    &\longrightarrow  \cat L(\mathbb R^{d+e},(V\oplus\mathbb R^e)^2)\times\Omega^{d+e}X(V\oplus\mathbb R^e)\\
    (A,\psi,f)&\longmapsto ((A\oplus A)\circ j(\psi),\sigma\circ(f\smashp S^e)) 
  \end{align*}
  is continuous, too. This map participates in the commutative square:
  \[\hskip-26.5pt\hfuzz=26.55pt
  \begin{tikzcd}[cramped]
      \O(V\oplus\mathbb R^e)\times \cat L(\mathbb R^d,V^2)\times\Omega^dX(V)
      \arrow[r]\arrow[d,"{(A,\psi,f)\mapsto (A|_V,[\psi.f])}"{description}]&
      \cat L(\mathbb R^{d+e},(V\oplus\mathbb R^e)^2)\times\Omega^{d+e}X(V\oplus\mathbb R^e)\arrow[d,"{\text{proj}}"]\\
      \cat L(V,V\oplus\mathbb R^e)\times \left( \cat L(\mathbb R^d,V^2)\times_{\O(d)}\Omega^dX(V)\right) 
      \arrow[r,"(\text{\ref{eq:d_summand}})"'] &
      \cat L(\mathbb R^{d+e},(V\oplus\mathbb R^e)^2)\times_{\O(d+e)}\Omega^{d+e}X(V\oplus\mathbb R^e)
     \end{tikzcd}
  \]
  Since the composite through the upper right corner
  is continuous and the left vertical map is a quotient projection,
  the lower horizontal map is continuous, too.
  This completes the proof.
\end{proof}
  This concludes the definition of the topological functor  $\textbf{Q}$.
  The unit and counit of the adjunction between the functors
  \[ \cat{Th}\colon \Gospc{G}_{/\BOP}\rightleftarrows \Gosp{G}\colon\textbf{Q}\]
  are levelwise the units and counits of the adjunction from Construction~\ref{constr:thom-level-right-adjt}. 

  We will need the following observation to show that the Thom space and Thom spectrum functor are fully homotopical:

\begin{lemma}\label{app-lemma:thom-spaces-homotopical}
 Let $G$ be a compact Lie group with a collection $\Ff$ of closed subgroups, and let $V$ be a $G$-representation. Then for any commutative triangle
 \[
    \begin{tikzcd}[column sep=small]
        B_1\arrow[dr,"fh"', bend right=15pt]\arrow[rr,"h"] && B_2\arrow[dl, "f", bend left=15pt]\\
        & \cat{Gr}(V)
    \end{tikzcd}
 \]
 in $\cat{$\bm G$-Top}$ such that $h$ is an $\Ff$-weak equivalence, the induced map $T(fh)\to T(f)$ is again an $\Ff$-weak equivalence.
\end{lemma}
\begin{proof}
  We start with the non-equivariant case, i.e. when the group $G$ is trivial. If $E_i,i=1,2$ denote the two pullbacks of the tautological bundle over $\Gr(V)$, then the square
  \[
  \begin{tikzcd}
  S(E_1)\arrow[r] \arrow[d]& B_1\arrow[d,"h"]  \\ 
  S(E_2)\arrow[r] & B_2
  \end{tikzcd}
  \]
  is cartesian and both horizontal maps are locally trivial fiber bundles, hence Serre fibrations. Thus, this is already a homotopy pullback. Since $h$ is a weak equivalence, the left vertical map in the square is a weak equivalence, too.
  The Thom spaces  $T(f h)$ and $T(f)$ are the unreduced mapping cones of the two horizontal maps
  in the pullback square.  
  The formation of unreduced mapping cones is homotopical for weak equivalences,
  for example by \cite{dugger-isaksen}*{Lemma~A.1} for $n=1$;
  so the induced map $T(fh)\to T(f)$ is a weak equivalence.

  For the general case, we have to understand fixed points of Thom spaces of equivariant bundles. If $g\colon X\to\cat{Gr}(V)$ is an equivariant continuous map and $H\subset G$, then one directly checks that the $H$-fixed points of the pullback $g^*\zeta$ of the tautological bundle agree as a set with the pullback of the tautological bundle over $\cat{Gr}(V^H)$ along the composite
  \[
    X^H\to\cat{Gr}(V)^H\xrightarrow{(-)^H}\cat{Gr}(V^H),
  \]
  where $(-)^H$ denotes the map sending an $H$-invariant subspace to its $H$-fixed points. The latter map is continuous: by \cite{schwede2018global}*{Proposition~B.17}, $\Gr_d(V)^H=(\cat{L}(\R^d,V)/\O(d))^H$ splits as the coproduct over conjugacy classes of homomorphisms $\phi\colon H\to\O(d)$ of the spaces $\cat{L}(\phi^*(\R^d),V)^H/C(\phi)$, where $C(\phi)\subset\O(d)$ denotes the centralizer of the image of $\phi$. If $f\colon\phi^*(\R^d)\to V$ is $H$-equivariant and isometric, then $f(\R^d)^H=f((\phi^*(\R^d))^{H})$, so that the restriction of $(-)^H$ to $\cat{L}(\phi^*(\R^d),V)^H$ factors as the map $\cat{L}(\phi^*(\R^d),V)^H\to\cat{L}((\phi^*(\R^d))^H,V)$ induced by the inclusion followed by the quotient map. This completes the proof that $(-)^H$ is continuous. As both the $H$-fixed points of $g^*\zeta_V$ as well as the pullback of $\zeta_{V^H}$ along the continuous map $(-)^H\circ g^H$ are closed subspaces of $g^*\zeta$, this then shows that they agree as topological spaces.
  
  The identification of the fixed point bundle immediately implies the analogous description of the sphere bundle of the fixed points; as formation of unreduced mapping cones commutes with fixed points by \cite{schwede2018global}*{Proposition~B.1}, this gives us a natural homeomorphism $T(g)^H\cong T((-)^H\circ g^H)$. Applying this to $f$ and $fh$ and varying $H\in\Ff$, the general case of the lemma therefore follows from the above non-equivariant special case.
\end{proof}

\begin{lemma}\label{lm:Thom-left-Quillen-wrt-level}
    \begin{enumerate}
\item The Thom space functor $\cat{th}\colon\cat{$\bm G$-$\cat{L}$-Top}_{/\cat{Gr}}\to\cat{$\bm G$-$\cat{L}$-Top}_*$ is left Quillen and fully homotopical with respect to the $G$-global \emph{level} model structures. 
\item The Thom spectrum functor $\cat{th}\colon\cat{$\bm G$-$\cat{L}$-Top}_{/\cat{BOP}}\to\Gosp{G}$ is left Quillen and fully homotopical with respect to the $G$-global \emph{level} model structures. 
    \end{enumerate}
\end{lemma}
\begin{proof}
    By Lemma~\ref{app-lemma:thom-spaces-homotopical}, the Thom space and spectrum functor preserve $G$-global level weak equivalences. 
    We have already seen that the Thom space and Thom spectrum functors are left adjoints, so it only remains to show that they preserve cofibrations, then it follows that they are left Quillen. 

    By Proposition~\ref{prop:orth-lvl} together with \cite{hirschhorn-book}*{Theorem 1.5}, the maps 
    \begin{equation}\label{diag:generic-cof-over-Gr}
        \begin{tikzcd}[column sep=small]
            \cat L(V,-)\times_{\phi}G\times \partial D^n\arrow[dr,bend right=15pt]\arrow[rr,hook] && \cat L(V,-)\times_{\phi}G\times D^n\arrow[dl,bend left=15pt,"f^\dual"]\\
            & \cat{Gr}\rlap,
        \end{tikzcd}
    \end{equation}
    where $V$ is a faithful $H$-representation, $H$ is a compact Lie group, $\phi\colon H\to G$ is a continuous group homomorphism, $n\ge0$, and $f^\dual$ is arbitrary, form a set of generating cofibrations for  $\Gospc{G}_{/\cat{Gr}}$. By adjunction, $f^\dual$ corresponds to a map $f\colon D^l\to\cat{Gr}(V)^\phi$. Denote by \[\pi\colon E_f\coloneqq \eGr(V)\times_{\Gr(V)}D^l\to D^l\] the bundle classified by $f$ and by $j\colon E_f\to \eGr(V)$ the projection.
  For an inner product space $W$, the map
  \[G\times_{\phi}\cat L(V,W)\times E_f\to \eGr(W),\, (g,\psi,x)\mapsto \psi(j(x))\] exhibits  
  \[ (G\times_{\phi}\cat L(V,W))\times E_f\xrightarrow{G\times_{\phi}\cat L(V,W)\times \pi}G\times_{\phi}\cat L(V,W)\times D^l\] as the bundle classified by $f^{\dual}$, i.e. 
  \[
  \begin{tikzcd}
  G\times_{\phi}\cat L(V,W)\times E_f \arrow[r]\arrow[d]\arrow[dr, pullback]& \eGr(W)\arrow[d]\\ 
  G\times_{\phi}\cat L(V,W)\times D^l\arrow[r,"f^{\dual}(W)"']& \Gr(W)
  \end{tikzcd}
  \] is a pullback of topological $G$-spaces. 
  It follows that $\cat{th}$ sends $(\text{\ref{diag:generic-cof-over-Gr}})$ up to isomorphism to $(G\times_{\phi}\cat L(V,-))_+\smashp(T(f|_{\partial D^l})\hookrightarrow T(f))$.

  By Lemma~\ref{lemma:G-Top*-tensored}, $(G\times_{\phi}\cat L(V,-))_{+}\smashp {-}\colon \cat{$\bm G$-Top}_*\to \Gospc{G}_*$ is left Quillen.
  We are therefore reduced to showing the following claim: for every inner product space $V$ and every continuous map $f\colon D^l\to\Gr(V)$, the associated map $i_l$ of Thom spaces
  over $\partial D^l$ and $D^l$ is a cofibration of based $G$-spaces.
  Because $D^l$ is contractible, the pulled back bundle is trivial, so
  after choosing a trivialization, we can identify the map $i_l$ with the cofibration
  \[ \partial D^l_{+}\smashp S^V\hookrightarrow D^l_{+}\smashp S^V ,\]
  for some representation $V$, completing the proof of the case of Thom spaces.

  In the case of Thom spectra, we similarly let $V$ be a faithful representation for a compact Lie group $H$, $\phi\colon H\to G$ a homomorphism, and $h\colon D^l\to \BOP(V)^H$ a continuous map. Denote by $\pi\colon E_h\coloneqq \eGr(V\oplus V)\times_{\BOP(V)}D^l\to D^l$ the associated bundle and by $i\colon E_h\to \BOP(V)$ the canonical map. 
  For an inner product space $W$, denote by $e^{\perp}(V,W)=\{ (A,w)\in\cat L(V,W)\times W : w\in A(V)^{\perp}\}\to \cat L(V,W)$ the orthogonal complement bundle. 
  Then \[(G\times_{\phi}e^{\perp}(V,W))\times E_h\to \eGr(W\oplus W), \, (g,(A,w),e)\mapsto (w,0)+(A\oplus A)(i(e))\] 
  (where $g\in G$, $A\in\cat L(V,W)$, $w\in A(V)^{\perp}$, and $e\in E_h$)
  exhibits \[G\times_{\phi}e^{\perp}(V,W)\times \pi\colon (G\times_{\phi}e^{\perp}(V,W))\times E_h\to G\times_{\phi}\cat L(V,W)\times D^l\] as the $G$-equivariant bundle classified by $h^{\dual}(W)$. We may therefore identify the image under $\cat{Th}$ of a generic generating cofibration
  \[
    \begin{tikzcd}[column sep=small]
            \cat L(V,-)\times_{\phi}G\times \partial D^n\arrow[dr,bend right=15pt]\arrow[rr,hook] && \cat L(V,-)\times_{\phi}G\times D^n\arrow[dl,bend left=15pt,"h^\dual"]\\
            & \cat{BOP}\rlap,
        \end{tikzcd}
  \]
  with $G_{+}\smashp_{\phi}\cat O(V,-)\smashp (T(h|_{\partial D^l})\hookrightarrow T(h))$. As above, the second smash factor is a cofibration of based $G$-spaces, so the whole expression is a $G$-global cofibration of orthogonal $G$-spectra by Lemma~\ref{lemma:G-global-spectra-tensoring}, proving that the Thom spectrum functor preserves cofibrations. 
\end{proof}

Theorem~\ref{thm:Thom-marvelous-appendix} is now a consequence of the above and the following observation:  
\begin{lemma}\label{lm:fibrationsbetweenfibrantobjectscartesian}
    Suppose that $\phi\colon V\to W$ is an equivariant linear isometric embedding of faithful $G$-representations. 
    \begin{enumerate}
    \item If $X\in\Gospc{G}_{*}$ is fibrant in the $G$-global model structure, the square 
\[
\begin{tikzcd}
\bm\Omega(X)(V)^G\arrow[d,"\bm\Omega(X)(\phi)^G"']\arrow[r,"\cat q(X)"]& \Gr(V)^G\arrow[d,"\cat{Gr}(\phi)^G"] \\ 
\bm\Omega(X)(W)^G\arrow[r,"\cat q(X)"']& \Gr(W)^G
\end{tikzcd}
\] is a homotopy pullback. 
\item Suppose that $X\in \Gosp{G}$ is fibrant in the $G$-global model structure. Then the square 
\[
\begin{tikzcd}
\bm\Omega(X)(V)^G\arrow[d,"\bm\Omega(X)(\phi)^G"']\arrow[r,"\cat Q(X)"]& \BOP(V)^G\arrow[d,"\cat{BOP}(\phi)^G"] \\ 
\bm\Omega(X)(W)^G\arrow[r,"\cat Q(X)"']& \BOP(W)^G
\end{tikzcd}
\]
 is a homotopy pullback.
\end{enumerate}
\end{lemma}
    
\begin{proof}Let $Z$ be a based $G$-space. We first describe the fixed-point space $\Omega_R(Z)^G$ for a $G$-representation $R$. 
  The group $\O(d)$ acts freely on $\cat L(\mathbb R^d,R)\times \Omega^dZ$; so the $G$-fixed point space $(\Omega_T X)^G$
  decomposes as the disjoint union, indexed by all $d\geq 0$ and all conjugacy classes of
  continuous homomorphisms $\alpha\colon G\to\O(d)$, of the spaces
  \[( \cat L(\mathbb R^d,R)\times \Omega^dZ)^{\alpha}/C(\alpha)\]
  by another application of \cite{schwede2018global}*{Proposition~B.17}; here $C(\alpha)$ again denotes the centralizer in $\O(d)$ of the image
  of $\alpha$. If we write $U$ for the $d$-dimensional $G$-representation $\alpha^*(\R^d)$, then $C(\alpha)$ 
  is the compact Lie group $\O^G(U)$ of $G$-equivariant linear self-isometries of $U$.
  We can thus rewrite the fixed point space $(\Omega_T(Z))^G$
  as the disjoint union, indexed by representatives $U$ of all isomorphism classes
  of $G$-representations, of the spaces
  \[\cat L^G(U,V)\times_{\O^G(U)} (\Omega^{U}Z)^G.\]

  Suppose now that $X$ is a pointed orthogonal $G$-space. 
  By the above, the commutative square in question is the disjoint union, indexed over representatives of isomorphism classes $U$ of $G$-representations 
  \[
  \begin{tikzcd}  \cat L^G(U,V)\times_{\O^G(U)} (\Omega^UX(V))^G\arrow[d,"{[\psi,f]\mapsto [\phi\circ \psi,X(\phi)\circ f]}"'] \arrow[r]&  \cat L^G(U,V)/\O^G(U) \arrow[d,"{[\psi]\mapsto [\phi\circ \psi]}"] \\
        \cat L^G(U,W)\times_{\O^G(U)} \mathcal (\Omega^UX(W))^G \arrow[r]& \cat L^G(U,W)/\O^G(U).
  \end{tikzcd}
\]
The two horizontal maps arise from the projections to the first factors, and they are both locally trivial fiber bundles, and in particular Serre fibrations.   The generic fiber of the top horizontal map is $(\Omega^UX(V))^G$, which maps to the corresponding fiber $(\Omega^UX(W))^G$ of the bottom horizontal map via $(\Omega^UX(\phi))^G$. 
  As $V$ is a faithful $G$-representation and $X$ is a fibrant orthogonal $G$-space, $X(\phi)$ is an equivariant weak equivalence of pointed $G$-spaces. By Lemma~\ref{lemma:smash-lQ}, $\Omega^U\colon\cat{$\bm G$-Top}_*\to\cat{$\bm G$-Top}_*$ is right Quillen, and hence homotopical by Ken Brown's Lemma. Thus, $(\Omega^UX(\phi))^G$ is a weak homotopy equivalence, proving the statement for Thom spaces. 

Now we prove the statement for Thom spectra.  By replacing $W$ by the orthogonal complement of the image of $\varphi\colon V\to W$,
  we can assume without loss of generality that $\varphi\colon V\to V\oplus W$ is
  the embedding of the first summand of a direct sum.
  Then the commutative square in question becomes the disjoint union, indexed over 
  representatives $U$ of isomorphism classes of $G$-representations, of commutative squares of the form 
  \[\hskip-23.56pt\hfuzz=23.6pt
  \begin{tikzcd}[cramped]
  \cat L^G(U,V^2)\times_{\O^G(U)} (\Omega^UX(V))^G\arrow[d]\arrow[r] & 
  \cat L^G(U,V^2)/\O^G(U)\arrow[d]\\
  \cat L^G(U\oplus W,V^2\oplus W^2)\times_{\O^G(U\oplus W)}(\Omega^UX(V\oplus W))^G\arrow[r]& \cat L^G(U\oplus W, V^2\oplus W^2)/\O^G(U\oplus W)\rlap.
  \end{tikzcd}
\]
  Here the vertical maps are induced by the map $\cat{L}(U,V^2)\to\cat{L}(U\oplus W,V^2\oplus W^2)$ sending $\phi$ to $(v,w)\mapsto (\phi(v),(w,0))$ and by the map $(\Omega^UX(V))^G\to(\Omega^UX(V\oplus W))^G$ induced by the adjunct structure map of $X$. The two horizontal maps again arise from the projections to the first factors, and they are both locally trivial 
  fiber bundles, and hence Serre fibrations. 

   The generic fiber of the top horizontal map is $(\Omega^UX(V))^G$, which maps to the corresponding fiber $(\Omega^UX(V\oplus W))^G$ of the bottom map via the map induced by the adjunct structure map. As the adjunct structure map is a $G$-weak equivalence by fibrancy of $X$, the claim follows as before.
\end{proof}

We can now finally prove Theorem~\ref{thm:Thom-marvelous-appendix} stated at the beginning of this section. 
\begin{proof}[Proof of Theorem~\ref{thm:Thom-marvelous-appendix}]
    As the $G$-global model structure on $\cat{$\bm G$-$\cat{L}$-Top}$ was constructed as a Bousfield localization of the $G$-global level model structure, Lemma~\ref{lm:Thom-left-Quillen-wrt-level} implies that that the Thom space functor preserves cofibrations (with respect to the $G$-global model structures on both sides). 
This implies that its right adjoint $\textbf{q}$ preserves acyclic $G$-global fibrations.
To show that $\textbf{q}$ also preserves $G$-global fibrations, we use a criterion of Dugger \cite{dugger:replacing}*{Corollary~A.2}: it suffices to show that the right adjoint preserves $G$-global fibrations $\alpha\colon X\to Y$
\emph{between $G$-globally fibrant objects}. 
Suppose $\alpha\colon X\to Y$ is a $G$-global fibration between fibrant objects in $\Gospc{G}_*$. 
As $\cat{q}$ is right Quillen with respect to the $G$-level structure on $\Gospc{G}_{/\Gr}$, $\cat{q}(\alpha)$ is a $G$-global level fibration. If $\phi\colon H\to G$ is any homomorphism, and $V\hookrightarrow W$ is an embedding of faithful $H$-representations, it then only remains to show by the explicit description of $G$-global fibrations (see Definition~\ref{def:g-gl-fib-orth-spc}) that the top square in the following commutative diagram is a homotopy pullback:
\[
  \begin{tikzcd}
    \bm\Omega(X)(V)^{\phi}\arrow[r]\arrow[d]& \bm\Omega(X)(W)^{\phi}\arrow[d]\\ 
    \bm\Omega(Y)(V)^{\phi}\arrow[r]\arrow[d]& \bm\Omega(Y)(W)^{\phi}\arrow[d]\\ 
    \Gr(V)^{H}\arrow[r]& \Gr(W)^{H},
  \end{tikzcd}
\]
Applying the previous lemma to the $H$-orthogonal spaces $\phi^*X$ and $\phi^*Y$ (which are again fibrant by Lemma~\ref{lemma:restr-right-Quillen}) shows that the bottom and the total square are homotopy pullbacks, whence so is the top square by the pasting law. 
This proves that the Thom space functor is left Quillen for the $G$-global model structure.

That the Thom spectrum functor is left Quillen can be seen completely analogously.

  A $G$-global weak equivalence of orthogonal $G$-spaces factors as a composite of an acyclic $G$-global cofibration followed by an acyclic $G$-global fibration.
  As a left Quillen functor, the Thom space/spectrum functor takes acyclic $G$-global cofibrations to $G$-global weak equivalences. Since acyclic fibrations in the $G$-global model structure
  are $G$-global level weak equivalences,
  the Thom space/Thom spectrum functor takes these to $G$-global level weak equivalences by
  Lemma~\ref{lm:Thom-left-Quillen-wrt-level}, and hence in particular to $G$-global weak equivalences. Thus, both the Thom space and Thom spectrum functor preserve $G$-global weak equivalences.
\end{proof}

\section[Background on model categories]{\for{toc}{\hspace{-.22em}}Background on model categories}\label{app:model-cat}
In this section we collect some general facts about model categories and their associated $\infty$-categories for easy reference. We begin with the following observation relating $1$-categorical and $\infty$-categorical slices:

\begin{lemma}\label{lemma:slice-over-fibrant}
    Let $\cat{C}$ be a model category with collection of weak equivalences $\Ww$, and let $Y\in\cat{C}$ fibrant. Then the ($\infty$-categorical) localization $\cat{C}\to\cat{C}[\Ww^{-1}]$ induces a localization $\cat{C}_{/Y}\to\cat{C}[\Ww^{-1}]_{/Y}$. If $\cat{C}$ is right proper, the claim holds without the fibrancy assumption on $Y$.
    \begin{proof}
        The first claim is a special case of \cite{cisinski-book}*{Corollary 7.6.13}. For the second claim, it will then suffice that for right proper $\cat{C}$ and any weak equivalence $Y\to Y'$ to a fibrant $Y'$, the induced map $\cat{C}_{/Y}\to\cat{C}_{/Y'}$ becomes an equivalence after localization, which is in turn a direct consequence of \cite{rezk-proper}*{Proposition~2.5}.
    \end{proof}
\end{lemma}

The model categories we encounter throughout our paper are simplicial (in fact, even topological), and several arguments exploit the simplicial structure to say something about the associated $\infty$-category, although the latter a priori only depends on the underlying unenriched model category. We therefore have to explain how the simplicial and homotopical structures relate to each other.

\begin{lemma}[See \cite{hirschhorn-book}*{Proposition 9.5.16}]\label{lemma:hom-equiv-we}
    Let $\cat{C}$ be a simplicial model category. Then any homotopy equivalence (with respect to the simplicial structure) is a weak equivalence (with respect to the model structure).\qed
\end{lemma}

\begin{cor}
    Let $\cat{C}$ be a simplicial model category and let $X\in\cat{C}$ be arbitrary. Then the map $X\to \Delta^1\otimes X$ induced by $d_1\colon\Delta^0\to\Delta^1$ is a weak equivalence.
    \begin{proof}
        One easily checks that this is a homotopy equivalence, also see \cite{hirschhorn-book}*{Proposition 9.5.20}.
    \end{proof}
\end{cor}

Note that crucially there are no cofibrancy conditions on $X$ in either of these results. Using this we can now show:

\begin{lemma}\label{lemma:fgt-simpl-structure}
    Let $\cat{C}$ be a simplicial model category, and let $\cat{D}\subset\cat{C}$ be a full simplicial subcategory closed under weak equivalences. Assume moreover that $\cat{D}$ is locally fibrant. Then the inclusion $N(\cat{D}_0)\hookrightarrow N_\Delta(\cat{D})$ of the nerve of the underlying unenriched category induces an equivalence after localizing at the weak equivalences.
    \begin{proof}
        By the previous corollary, $\cat{D}$ is closed under $\Delta^1\otimes-$, while Lemma~\ref{lemma:hom-equiv-we} shows that any homotopy equivalence is a weak equivalence. The claim is now an instance of \cite{HA}*{Proposition 1.3.4.7}.
    \end{proof}
\end{lemma}

\begin{lemma}\label{lemma:fgt-simpl-structure-cof}
    Let $\cat{C}$ be a simplicial model category, and assume that the full simplicial subcategory $\cat{C}^c\subset\cat{C}$ of cofibrant objects is locally fibrant. Then $N(\cat{C}^c_0)\hookrightarrow N_\Delta(\cat{C}^c)$ induces an equivalence on $\infty$-categorical localizations.
    \begin{proof}
        The left Quillen functor $\Delta^1\otimes-$ preserves $\cat{C}^c$, so this follows by the same argument as in the proof of the previous lemma.
    \end{proof}
\end{lemma}

Finally, the simplicial structure allows us to compute mapping spaces in the localization:

\begin{lemma}\label{lemma:localization-on-homs}
    Let $\cat{C}$ be a locally fibrant simplicial model category (e.g.\ the underlying simplicial model category of topological model category), and let $X,Y\in\cat{C}$ with arbitrarily chosen cofibrant replacement $\alpha\colon X'\iso X$ and fibrant replacement $\beta\colon Y\iso Y'$. Then the map
    \[
        \maps_{\cat{C}}(X,Y)\simeq\hom_{N_\Delta(\cat{C})}(X,Y)\to\hom_{N_\Delta(\cat{C}')[\Ww^{-1}]}(X,Y)
    \]
    induced by the localization $N_\Delta(\cat{C})\to N_\Delta(\cat{C})[\Ww^{-1}]$ agrees up to postcomposition with an equivalence with $\beta\circ{-}\circ\alpha\colon\maps(X,Y)\to\maps(X',Y')$. 
    \begin{proof}
        As the map $\hom_{\kern-.5ptN_\Delta(\cat{C})[\Ww^{-1}]\kern-.5pt}(X,Y)\to\hom_{\kern-.5ptN_\Delta(\cat{C})[\Ww^{-1}]\kern-.5pt}(X',Y')$ induced by $\alpha$ and $\beta$ is an equivalence, it will suffice to consider the case $\alpha=\id_X,\beta=\id_Y$, so that $X$ is cofibrant and $Y$ is fibrant. Using that $\maps\colon\cat{C}^\op\times\cat{C}\to\cat{SSet}$ is a right Quillen bifunctor, we are further reduced to showing that the localization is fully faithful when restricted to the full subcategory of cofibrant-fibrant objects. But this restriction is even an equivalence by \cite{dk-modern}*{Corollary~1.4.3}.
    \end{proof}
\end{lemma}

\begin{cor}
    Let $\cat{C}$ be a locally fibrant simplicial model category, and let $G\colon\cat{C}\to\cat{SSet}$ be a functor that is corepresented (in the enriched sense) by a cofibrant $X\in\cat{C}$. Then $G$ is right Quillen and
    \begin{equation}\label{eq:right-derived-rep}
        \cat{C}[\Ww^{-1}]\xrightarrow{\;\cat{R}G\;}\cat{SSet}[\Ww_\textup{Quillen}^{-1}]\simeq\Spc
    \end{equation}
    is corepresented (in the $\infty$-categorical sense) by the same $X$; here the unlabelled equivalence in $(\textup{\ref{eq:right-derived-rep}})$ is the unique one.
    \begin{proof}
        As the mapping space functor for $\cat{C}$ is a right Quillen bifunctor, $G$ is indeed right Quillen. It then suffices to show that the two functors $\cat{C}[\Ww^{-1}]$ agree after restriction to the full subcategory of $\cat{C}$ spanned by the cofibrant-fibrant objects, where this follows by the same argument as in the previous lemma.
    \end{proof}
\end{cor}

Applying this to the underlying simplicial model category of a topological model category, we obtain:

\begin{cor}\label{cor:corep-Top}
    Let $\cat{C}$ be a topological model category, and let $G\colon\cat{C}\to\cat{Top}$ be corepresented (in the enriched sense) by a cofibrant $X\in\cat{C}$. Then $G$ is right Quillen and 
    \[
        \cat{C}[\Ww^{-1}]\xrightarrow{\;\cat{R}G\;}\cat{Top}[\Ww_\textup{Quillen}^{-1}]\simeq\Spc
    \]
    is corepresented (in the $\infty$-categorical sense) by the same $X$.\qed
\end{cor}

\section[Compact Lie groups acting on Hausdorff spaces]{\for{toc}{\hspace{-.22em}}Compact Lie groups acting on Hausdorff spaces}\label{app:here-be-pointsets}
In this appendix we collect some pointset topological facts about free actions of compact Lie groups on Hausdorff spaces that are used throughout the main text.

\begin{lemma}[See \cite{koerschgen}*{Theorem~A.9}]\label{lemma:quotient-map-is-fib}
    Let $X$ be a Hausdorff space with a free action of a compact Lie group $G$. Then the quotient map $X\to X/G$ is a (non-equivariant) Serre fibration.\qed
\end{lemma}

\begin{lemma}\label{lemma:quotient-pb}
    Let $G$ be a compact Lie group and let $f\colon X\to Y$ be a map of free $G$-spaces, such that $X$ is Hausdorff. Then
    \[
        \begin{tikzcd}
            X\arrow[r,"\textup{proj}"]\arrow[d,"f"'] & X/G\arrow[d,"f/G"]\\
            Y\arrow[r,"\textup{proj}"'] & Y/G
        \end{tikzcd}
    \]
    is a pullback square.
    \begin{proof}
        By \cite{schwede2018global}*{Proposition~B.13(i)}, the quotients $X/G$ and $Y/G$ are preserved by the forgetful functor to sets. Thus, to show that the comparison map $c\colon X\to X/G\times_{Y/G}Y$ is bijective, we may instead prove the analogous statement in the category of sets. This can be checked fiberwise over $X/G$; the fiber over some class $[x]$ is then given by the orbit $Gx$ (for some fixed representative), while the fiber of the projection $X/G\times_{Y/G}Y$ over the same $[x]$ is given by $\{[x]\}\times Gf(x)$. By assumption, both of these are free $G$-orbits, and any equivariant self-map of $G$ is bijective, proving the claim.
    
        It remains to show that the induced map 
        \[
            \hom(K,c)\colon\hom_{\cat{Top}}(K,X)\to\hom_{\cat{Top}}(K,X/G\times_{Y/G}Y)
        \]
        is bijective for every compact space $K$. As we already know that $c$ is bijective, it will be enough to show surjectivity, for which we consider any $(\alpha,\beta)\colon K\to X/G\times_{Y/G}Y$. The space $X/G$ is again Hausdorff by \cite{tom-Dieck-trafo}*{Proposition~3.1(v)}, so $\alpha(K)\subset X/G$ is compact, and Proposition~3.6(ii) of \emph{op.\ cit.} shows that also the preimage $X'$ of $\alpha(K)$ under $X\to X/G$ is compact.         
        By \cite{schwede2018global}*{Proposition B.13(iii)}, the map $X'/G\times_{Y/G}Y\to X/G\times_{Y/G}Y$ induced by the inclusion $X'\hookrightarrow X$ is an embedding, so $(\alpha,\beta)$ factors through a continuous map $K\to X'/G\times_{Y/G}Y$, and it will be enough to show that this admits a lift to $X'$. Replacing $X$ by $X'$, we may therefore assume without loss of generality that $X$ is compact. But then $X\to X/G\times_{Y/G}Y$ is a continuous bijection from a compact space to a weak Hausdorff space, hence a homeomorphism.
    \end{proof}
\end{lemma}

\begin{lemma}[cf.\ \cite{koerschgen}*{proof of Theorem~A.7}]\label{lemma:quotient-we}
    Let $G$ be a compact Lie group, and let $f\colon X\to Y$ be a $G$-equivariant map of $G$-free Hausdorff spaces. Then $f$ is a (non-equivariant) weak homotopy equivalence if and only if the induced map $f/G\colon X/G\to Y/G$ on orbits is so.
    \begin{proof}
        By Lemmas~\ref{lemma:quotient-map-is-fib} and~\ref{lemma:quotient-pb}, the commutative square
        \[
            \begin{tikzcd}
                X\arrow[d,"f"']\arrow[dr,pullback]\arrow[r] & X/G\arrow[d,"f/G"]\\
                Y\arrow[r] & Y/G
            \end{tikzcd}
        \]
        is a pullback in $\cat{Top}$ and the horizontal maps are Serre fibrations. Thus, this square defines a homotopy pullback in the $\infty$-topos $\Spc$. As $Y\to Y/G$ is surjective on $\pi_0$, i.e.\ an effective epimorphism, the claim follows from \cite{HTT}*{Lemma 6.2.3.16}.
    \end{proof}
\end{lemma}

\backmatter
\bibliography{reference}
\end{document}